\documentclass[11pt,a4paper,oneside]{book}

\usepackage{amsmath,amssymb,amsthm,mathtools}
\usepackage[utf8]{inputenc}
\usepackage[T1]{fontenc}
\usepackage{geometry}
\usepackage{hyperref}
\usepackage{xcolor}
\usepackage{booktabs}
\usepackage{array}
\usepackage{longtable}
\usepackage{graphicx}
\usepackage{cleveref}
\usepackage{enumitem}

\hypersetup{colorlinks=true, linkcolor=blue, citecolor=red, bookmarksnumbered=true,
  pdftitle={Regular Arithmetic Functions, Volume I. Theory, Applications, Examples},
  pdfauthor={Benoit Cloitre},
  pdfsubject={Regular arithmetic functions and functions of good variation, a branch of
    Tauberian theory},
  pdfkeywords={regular arithmetic function, function of good variation, regularity index,
    transparency threshold, Mellin transform, Tauberian theory, discrete Volterra
    equation, Ingham kernel, Riemann hypothesis}}

\newtheorem{theorem}{Theorem}[chapter]
\newtheorem{lemma}[theorem]{Lemma}
\newtheorem{proposition}[theorem]{Proposition}
\newtheorem{corollary}[theorem]{Corollary}
\theoremstyle{definition}
\newtheorem{definition}[theorem]{Definition}
\newtheorem{conjecture}[theorem]{Conjecture}
\newtheorem{openproblem}[theorem]{Open Problem}
\newtheorem{conditionaltheorem}[theorem]{Conditional Theorem}
\crefname{conditionaltheorem}{Conditional Theorem}{Conditional Theorems}
\newtheorem{example}[theorem]{Example}
\theoremstyle{remark}
\newtheorem{remark}[theorem]{Remark}

\newtheorem{numobs}[theorem]{Numerical Observation}
\newtheorem*{remarkx}{Remark}

\newcommand{\eps}{\varepsilon}
\newcommand{\R}{\mathbb{R}}
\newcommand{\C}{\mathbb{C}}
\newcommand{\N}{\mathbb{N}}
\newcommand{\Z}{\mathbb{Z}}
\newcommand{\Q}{\mathbb{Q}}
\newcommand{\gstar}{g^{*}}
\providecommand{\og}{\guillemotleft}\providecommand{\fg}{\guillemotright}
\newcommand{\rafepigraph}[3]{%
  \begin{flushright}
  \begin{minipage}{0.80\linewidth}
    \small\itshape #1\par
    \vspace{0.4\baselineskip}
    \normalfont\footnotesize\raggedleft #3\par
  \end{minipage}
  \end{flushright}
  \vspace{0.6\baselineskip}}
\newcommand{\nm}[2]{#1\index[names]{#1, #2}}
\newcommand{\Gfun}{\mathcal{G}}
\newcommand{\gcirc}{g^{\circ}}
\DeclareMathOperator{\Res}{Res}
\renewcommand{\Re}{\operatorname{Re}}
\renewcommand{\Im}{\operatorname{Im}}

\makeatletter
\renewcommand*\l@chapter[2]{%
  \ifnum \c@tocdepth >\m@ne
    \addpenalty{-\@highpenalty}%
    \vskip 1.0em \@plus\p@
    \setlength\@tempdima{2.6em}%
    \begingroup
      \parindent \z@ \rightskip \@pnumwidth
      \parfillskip -\@pnumwidth
      \leavevmode \bfseries
      \advance\leftskip\@tempdima
      \hskip -\leftskip
      #1\nobreak\hfil
      \nobreak\hb@xt@\@pnumwidth{\hss #2\kern-\p@\kern\p@}\par
      \penalty\@highpenalty
    \endgroup
  \fi}
\makeatother

\newcounter{dossier}
\newcommand{\dossierchapter}[2]{%
  \cleardoublepage
  \setcounter{chapter}{0}%
  \stepcounter{dossier}%
  \renewcommand{\thechapter}{D\thedossier}%
  \renewcommand{\theHchapter}{dossier\thedossier}%
  \renewcommand{\appendixname}{Research dossier}%
  \renewcommand{\chaptername}{Research dossier}%
  \chapter{#2}%
  \label{app:dossier_#1}%
}

\newcommand{\galleryentry}[7]{%
  \chapter[#2]{}%
  \label{app:#1}%
  \par\nobreak\vspace{-1.2em}
  \noindent\colorbox{black!4}{%
    \parbox{\dimexpr\linewidth-2\fboxsep}{%
      \vspace{3pt}
      \begin{tabular}{@{}p{0.155\linewidth}@{}p{0.815\linewidth}@{}}
        \textsc{kernel}    & #3 \\[2pt]
        \textsc{class}     & #4 \\[2pt]
        \textsc{transform} & #5 \\[2pt]
        \textsc{index}     & #6 \\[2pt]
        \textsc{status}    & #7 \\
      \end{tabular}
      \vspace{3pt}}}%
  \par\medskip
}

\newenvironment{proofstatus}[1]{%
  \par\medskip
  \noindent\colorbox{yellow!30}{\parbox{\dimexpr\linewidth-2\fboxsep}{%
    \textbf{Proof status note.} #1}}%
  \par\medskip
}{}

\graphicspath{{./}{figures/}{Figures/}{fig/}{figs/}{img/}{images/}%
  {manuscript/figures/}{../figures/}{../manuscript/figures/}}
\newcommand{\rafshowfig}[1]{\includegraphics[width=0.72\textwidth]{#1}}
\newcommand{\rafgalleryfig}[3]{%
  \begin{figure}[htbp]
  \centering
  \rafwithfig{#1}{\rafshowfig}%
    {\fbox{\parbox[c][1.75in][c]{0.72\textwidth}{\centering\small
       Figure missing. Place the \texttt{figures} folder in the same directory as
       this source file, and run the compiler from that directory. The log records
       every lookup on a line beginning \texttt{RAF-FIG}. The box has the height of
       the figure it replaces, so the pagination is unaffected.}}}
  \caption{#2}
  \label{#3}
  \end{figure}}
\newcommand{\raffound}[3]{\typeout{RAF-FIG: found #1 as #2}#3{#2}}
\newcommand{\rafwithfig}[3]{%
  \IfFileExists{#1}{\raffound{#1}{#1}{#2}}{%
  \IfFileExists{figures/#1}{\raffound{#1}{figures/#1}{#2}}{%
  \IfFileExists{Figures/#1}{\raffound{#1}{Figures/#1}{#2}}{%
  \IfFileExists{fig/#1}{\raffound{#1}{fig/#1}{#2}}{%
  \IfFileExists{figs/#1}{\raffound{#1}{figs/#1}{#2}}{%
  \IfFileExists{img/#1}{\raffound{#1}{img/#1}{#2}}{%
  \IfFileExists{images/#1}{\raffound{#1}{images/#1}{#2}}{%
  \IfFileExists{manuscript/figures/#1}{\raffound{#1}{manuscript/figures/#1}{#2}}{%
  \IfFileExists{../figures/#1}{\raffound{#1}{../figures/#1}{#2}}{%
  \IfFileExists{../manuscript/figures/#1}{\raffound{#1}{../manuscript/figures/#1}{#2}}{%
  \typeout{RAF-FIG: MISSING #1}#3}}}}}}}}}}}

\title{\Huge \textbf{Regular Arithmetic Functions}\\ \vspace{1.4em}
  \huge \textsc{Volume I}\\ \vspace{0.7em}
  \Large Theory, Applications, Examples}
\author{\Large Benoît Cloitre}
\date{}

\usepackage{imakeidx}
\makeindex[name=names,title=Index of Names,columns=2,intoc]
\makeindex[name=terms,title=Index of Terms,columns=2,intoc]

\begin{document}

\frontmatter
\hypersetup{pageanchor=false}
\maketitle
\cleardoublepage

\thispagestyle{empty}
\vspace*{0.32\textheight}
\begin{center}
{\large\itshape \`A Anne, Matthieu et Thomas}
\end{center}
\cleardoublepage
\hypersetup{pageanchor=true}

\setcounter{tocdepth}{1}%
\makeatletter
\renewcommand*\l@section{\@dottedtocline{1}{1.5em}{2.9em}}
\makeatother
\tableofcontents

\chapter*{Preface to the two volumes}

\rafepigraph{Heureusement pour les chercheurs, à mesure que les brouillards se dissipent sur un point, c'est pour se reformer sur un autre.}{Fortunately for the researcher, as the fog lifts at one point, it is only to gather again at another.}{André Weil\index[names]{Weil, A.}, \emph{De la métaphysique aux mathématiques} (1960)~\cite{WeilMetaphysique1960}}

\addcontentsline{toc}{chapter}{Preface to the two volumes}
\markboth{PREFACE TO THE TWO VOLUMES}{}

This book studies a critical exponent arising in a class of difference equations. These equations
are triangular summation systems in which a fixed arithmetic kernel acts on an unknown sequence
to produce a prescribed rate of decay. Written out, the system is

\[
\sum_{k=1}^{n} a_k\,G(n,k)=n^{-\beta},\qquad n\ge 1,
\]

where the kernel $G$ is fixed and each exponent $\beta$ determines the sequence $(a_k)$ by forward
substitution. The unknown appears only under the sum, so the system is a discrete Volterra
equation\index[terms]{Volterra equation} of the first kind. The kernels studied here have a
nonvanishing diagonal, which isolates $a_n$ and turns the system into one of the second kind, the
form the substitution uses. What the theory watches is the partial sum $A(n)=\sum_{k\le n}a_k$,
whose dependence on $G$ and on $\beta$ is left implicit throughout.

For small $\beta$ the partial sums reproduce the rate they are given, and the equation transmits.
At and above a certain exponent, the kernel controls their rate of decay, and the equation
absorbs. The exponent separating the two behaviors is the regularity index\index[terms]{regularity index} of the kernel. The
equation itself is classical. The threshold is not, and the first volume is built around it and
the transition it marks.

The Riemann hypothesis\index[terms]{Riemann hypothesis} enters the theory through the Ingham\index[names]{Ingham, A. E.} kernel

\[
G(n,k)=\Phi(k/n)=\frac{k}{n}\left\lfloor\frac{n}{k}\right\rfloor,
\qquad \Phi(x)=x\lfloor 1/x\rfloor.
\]

This kernel underlies Ingham's 1945 summation method, used in his proof of the prime number
theorem. Its regularity index equals one half if and only if the Riemann hypothesis holds.
Chapter~\ref{chap:equivalence} proves this equivalence. Before the rate comes the limit.
Read on the boundary of its transform, a kernel returns limits and nothing more, and
Chapter~\ref{chap:boundary} shows that for the Ingham kernel this return is the prime number
theorem, Ingham's own theorem being proved there by Wiener's method.

In the cases treated here, the transmitted power comes with a coefficient given by the reciprocal
of the Mellin transform\index[terms]{Mellin transform} of the kernel, evaluated at the prescribed exponent. The connection with
the zeta function passes through that transform. Its
zeros provide a natural candidate for the regularity threshold, but the examples show that the
two need not agree. Seven of the worked kernels have transforms with no zeros at all, and six of
them nevertheless carry a proved regularity index. The kernel broken on the lattice of the square
root of two shows a separation of another kind. All the zeros of its transform lie on the line
$\Re z=1$, and the discrete equation already loses transparency\index[terms]{transparency} at $\beta=1/2$. Its exact
regularity index remains open. What is proved is the separation of the transparency frontier from
the zeros, and that is enough to show that the transform does not determine the arithmetic
response.

The regularity index therefore cannot always be read from the zeros of the transform. It is an
arithmetic quantity defined by the discrete equation, independently of those zeros, although in
special cases the transform carries the same information. Known equivalents of the Riemann
hypothesis are often organized into arithmetic and analytic families, as in Broughan\index[names]{Broughan, K.}~\cite{Broughan2017I,Broughan2017II}. The
equivalence proved here cuts across that division. The Ingham function is the case where the two
readings meet, and the Riemann hypothesis is precisely the statement that they do so at one half.
The hypothesis is thus recast as a compatibility between the arithmetic behavior of a discrete
equation and the analytic distribution of the zeros of its transform.

The framework also reaches a problem of a different kind. Applied to the orthorecursive expansion
of unity studied by Kalmynin\index[names]{Kalmynin, A. B.} and Kosenko\index[names]{Kosenko, P. R.}, it improves the known decay of the partial sums from
$O(N^{-1/2})$ to $O_\eps(N^{-\alpha_1+\eps})$, sharpens the pointwise bound to
$O(N^{-2})$, and develops the associated spectral picture. Chapter~\ref{chap:ortho} thus shows the theory at
work outside the setting of the Riemann zeta function. Full optimality is not reached, the
nonvanishing of the leading amplitude being left as an open problem.

The proof combines the Mellin transform with a Volterra resolvent constructed by a Neumann
series. This makes it possible to shift the contour on the smooth resolvent rather than on the
step function of partial sums, and to transfer the information carried by the zeros back to the
original asymptotic problem.

I have chosen a deliberately bottom-up approach. Observations lead to examples, examples to first
notions, and the first notions, once some case breaks them, to corrections and then to general
classes. Definitions arrive when the examples have made them necessary. The regularity index was
not posited and then tested. It appeared in the behavior of particular sequences, it was measured
before it was defined, and it took its present form only after several attempts at a definition
had been broken by an example. What appears here is the framework that survived that process and
now provides a solid foundation for a theory that may develop in many directions.

The gallery of kernels is deliberately varied. Some kernels are purpose-built rather than drawn
from a classical problem, because each isolates a different behavior or a different obstacle to
proving regularity. Their proofs use exact recurrences, Gamma products, discrete resolvents and
the homogeneous solutions of the underlying difference equations. Many are constructive, in that
even when they integrate over a contour or use special functions they show how the particular
operator acts. No single method covers all the examples. Even the appealing idea that the
regularity index should be read from the homogeneous solution fails for some kernels. The variety
of the examples is therefore part of the theory, and it shows both what the present framework
captures and where a more general method is still missing.

The examples are not added to illustrate a finished theory. They are the material from which the
theory is built.

It is not a finished theory, and a theory built from examples inherits the incompleteness of its
examples. That is the accepted cost of the method. Statements that a settled subject would carry
as theorems appear here as conjectures, as conditional theorems, or as open problems, and they
are marked as such wherever they occur. The register at the end of the volume gathers them with
their page numbers. The classes of functions on which the index lives are not delimited, and the
general membership question is open. Where the transform has a zero, whether the leftmost one
bounds the index read off the equation, in either direction, is not known. The volume exhibits a
kernel where the equation gives strictly less than the first zero of the transform, together with
a kernel where the reverse inequality is the one that is proved.

The last chapters deform the kernel by a gauge, a change of the coordinates in which the kernel
is read. One gauge turns the transform of the Ingham kernel\index[terms]{Ingham kernel} into a rational function of the kind
a curve over a finite field carries, for the admissible prime powers. That the gauged transform
has all its zeros on the critical line is proved and unconditional. The gauged arithmetic index
is also proved to be one half, by the companion result quoted in
Theorem~\ref{numobs:gauge_critical}. Equilibrium with the ordinary Ingham index remains
conjectural, and it is equivalent to the Riemann hypothesis.

Volume II explores lines of attack on the Riemann hypothesis and on its generalized form, built
on the framework of the first volume and on the gauge theory of its last chapters. The first
volume can be read independently of the second.

There is more than one way into the book, and the table of contents is built so that each can be
followed. Five parts carry the theory. Part I builds the notion of a function of good variation
out of Ingham's weight, recovers the prime number theorem from the boundary of its transform and proves the
equivalence, in five chapters that need nothing from the rest. Part II widens
the kernel to a function of two separate arguments, the regular arithmetic function that names the
volume, and inverts the operator twice, on the half line for the Ingham kernel and at the discrete level for every kernel.
Part III shows that the index belongs to the equation and not to the transform, and that the word
regular is earned, by reading the equation as a Mercerian problem next to Karamata's classes. Part IV changes
the coordinates in which a kernel is read and asks what the index does under such a change. Part V
reads the transform beyond its first zero, as a density on its half plane, through the fractional
part sums where a gauge reaches the divisor problem, and as a sum over all its zeros, the trace
formulas, before settling the orthorecursive expansion of unity. The epilogue names the threads the
second volume takes up. The gallery of appendices works the kernels one by one, from A to Q, each
entry stating its kernel, its transform, its indices and the status of their proof before the
proofs themselves, and a research dossier on the kernel broken at the square root of two closes it.
Each part opens with a card of what it does, each chapter with the order of its sections, and the
chapter on notations fixes the conventions.

A reader who wants the equivalence with the Riemann hypothesis and nothing else may read
Chapters~\ref{chap:ingham} to~\ref{chap:equivalence}, then Chapter~\ref{chap:volterra}, where the same equivalence is proved a second time through the
Volterra resolvent, Chapter~\ref{chap:diophantine}, which says why the equivalence has content by exhibiting kernels
where the arithmetic and the analytic readings of the index come apart, and Chapter~\ref{chap:gauge_ingham} with the
epilogue, where a change of coordinates gives the transform of the Ingham function the shape of a
zeta function over a finite field.

A reader arriving from Tauberian theory or from regular variation may start with Chapter~\ref{chap:ingham}, which
places Ingham's method among the summability procedures and follows what the literature made of
it, take Chapter~\ref{chap:boundary}, where Wiener's theorem is written in the coordinates of the
equation, then go to Chapter~\ref{chap:regularity}, to Chapter~\ref{chap:gauge}, whose gauges are a device borrowed from that theory, and
to Chapter~\ref{chap:abelian_densities}, where the Abelian direction is treated for its own sake. Appendices B and F work the
two logarithmic profiles, unbounded at the origin.

A reader who prefers exact machinery, difference equations, resolvents and spectra, may begin with
Chapter~\ref{chap:raf} and its scalar reduction, take Chapter~\ref{chap:principles} for the kernels solved in closed form and
Chapter~\ref{chap:discrete_volterra} for the discrete resolvent, and go on to Chapters~\ref{chap:trace_poly} and~\ref{chap:ortho}, where the resolvent is
expanded over the zeros of the transform. The gallery is that reader's natural companion.

A reader who, like me, would rather meet examples before definitions may read Chapter~\ref{chap:examples} and the
last section of Chapter~\ref{chap:raf}, then the gallery, and return to Chapters~\ref{chap:fgv} and~\ref{chap:raf} for the definitions
once the examples have made them necessary.

Corrections, questions and solutions are welcome, at \href{mailto:benoit.cloitre+raf@proton.me}{\texttt{benoit.cloitre+raf@proton.me}}. I expect some
of the open problems to be within reach, and I would rather learn that from a reader than not
learn it at all.

\bigskip

\chapter*{Acknowledgements}
\addcontentsline{toc}{chapter}{Acknowledgements}
\markboth{ACKNOWLEDGEMENTS}{}

This work owes a great deal to conversations held over fifteen years. The people named
individually below, in alphabetical order, each left a mark on the theory, on its proofs, or on the
road that led to it.

Jean-Paul Allouche has been a trusted friend and mathematical interlocutor for many years, from
our early exchanges through the community of integer sequences to our joint paper, with Vladimir
Shevelev, \emph{Beyond odious and evil}. Over the years he has read and discussed many of my ideas and
projects, offered practical advice, and given me a confidence that has mattered throughout the
development of this theory.

Brian Conrey showed me why an early attempt to perturb the Ingham function by an external
multiplicative factor could not work, directing me to Titchmarsh\index[names]{Titchmarsh, E. C.}'s classical theorem on the value
distribution of the zeta function. The obstruction left that route no room. The factor later moved
inside, changing the coordinates in which the kernel is read rather than multiplying the kernel
itself, and that became the gauge theory of Chapter~\ref{chap:gauge}.

Michel Lapidus, whom I had contacted after encountering his inverse-problem perspective on the
Riemann hypothesis, kindly met me during a visit to Paris. As we walked through the Luxembourg
Gardens, our exploratory conversation helped me see that my search for the arithmetic question
behind the critical-line statement might have an abstract form. That thought remained with me and
became one of the origins of Volume II.

Ricardo P\'erez-Marco, over a drink at Odéon in Paris, listened to an early attempt of mine to
read the Riemann hypothesis topologically and showed me that the homotopy argument I had proposed
did not justify treating the regularity index as an invariant in that setting. He encouraged me to
set that programme aside until its foundations were clearer, while urging me to write up the
equivalence at the heart of this volume on its own, and those questions return, in a different
form, in Volume II.

Neil J. A. Sloane created the On-Line Encyclopedia of Integer Sequences, which became my
mathematical workshop for many years. Submitting sequences, adding comments and formulas to
existing entries, and searching the database trained me in experimental pattern recognition. At the
beginning of these investigations, when I entered the values of the floor transform of the
Liouville function, the OEIS identified them as $\lfloor\sqrt n\rfloor$. That recognition was the
observation from which this book began.

G\'erald Tenenbaum\index[names]{Tenenbaum, G.} was the first to tell me that the question I was asking was a Tauberian one.
He met my first formulations with counterexamples that broke them, and he sent me to Korevaar\index[names]{Korevaar, J.},
where I found Ingham's theorem and with it the kernel this book is built on. Years later he
constructed a lacunary profile whose transform carries a natural boundary. That example opened a
question I had not thought to ask, how such an analytic barrier can coexist with arithmetic
regularity, and it became the base of the family of functions of good variation of
\S\ref{sec:ex-lacunary}, where its monotone members answer that question and his own profile
stays open.

Mark Daniel Ward gave me great support during the early years and opened the literature on
Ingham summation to me. He later read substantial versions of the project, notably in 2014, and has
encouraged me since to seek expert reading and a wider audience.

Herbert Wilf discussed binomial recurrences with me, on a question that had nothing to do with
the ideas that became this theory. I sent him a draft at the time and asked him in jest to prove
the Riemann hypothesis. He laughed, and then said seriously that I should look at whether the
divisor problem could be tied to the properties of Ingham's weight. Nothing came of it then. It is
only recently, returning to his remark from the Abelian side, that the divisor problem entered the
volume, in the Abelian sums and the fractional-part sums of Part V.

Doron Zeilberger took an interest in the early, still experimental form of the theory and
invited me in 2016 to present it at the Rutgers Experimental Mathematics Seminar, my first
mathematical talk at a university. His confidence at that early stage, and his continued
invitations and encouragement since, have mattered greatly.

This theory also grew within a wider community, to which I owe a different kind of debt.
Jean-Paul Delahaye, Boris Gour\'evitch and Steven R. Finch made me want to explore $\pi$, the
prime numbers and the mathematical constants, and so shaped in me a mathematical culture far wider
than the one I had started with. Paul D. Hanna and the contributors to the OEIS and to SeqFan made
that exploration collaborative. PARI/GP was the computational workshop in which much of this
theory was first explored, and I thank its authors and maintainers at the University of Bordeaux,
and the PARI Group, for making such a system for number theory freely available. I owe this
circle the curiosity and the experimental practice from which the work eventually grew. I also
thank Olivier Bordell\`es and my other coauthors, with whom I worked on several papers in number
theory before the present theory took shape, and Badih Ghusayni, whose invitation led to the
publication of the 2016 article.

\section*{Note on method and on the use of AI assistants}

This monograph is the result of research begun around 2010 and publicly documented from 2011
onward, more than a decade before contemporary generative AI. Roughly four fifths of its content
existed before any assistant was involved. The volume takes up, verifies and develops the results
and the ideas of my earlier papers, the article of 2016~\cite{Cloitre2016}, the arXiv
preprint~\cite{CloitreFloor}, and the papers on the orthorecursive expansion of
unity~\cite{CloitreOrtho} and on the fractional parts of Fibonacci
quotients~\cite{CloitreFibonacci}, together with dozens of unpublished notes,
technical reports of thousands of computations carried out with PARI/GP over the years, and
dozens of worked examples, of which those printed here are the most significant and the most
illustrative. The first service of these tools was therefore integration, and what they changed
was the pace at which that material could be brought into one manuscript.

The theory is mine. So are the definitions it rests on, the questions it asks, the organizing ideas
that give it its shape, and the results that structure the volume. So is the strategy of
presentation, the decision to move from observations to examples and only then to general
statements, the decision to treat the gallery of worked kernels as the material of the theory
rather than as illustrations of it, and the decision to state the connection with the Riemann
hypothesis as an equivalence. The wording and the tone are mine as
well, and so are the notation and the names of the mathematical objects, both kept to a minimum.
Nothing is named for the sake of naming, and an object receives a name only when understanding
requires one and the object recurs often enough to earn it.

Beyond integration, the assistants helped me experiment. They ran and checked computations, which
let me uncover new patterns and confirm ones I already trusted. They searched and organized the
literature, reread proofs both old and new, helped with the English, the typesetting and the
composition of the book, and kept a manuscript of this length consistent with itself. They also
helped me carry out lengthy technical calculations, including the truncated Perron inversion and
contour estimates in Appendix~\ref{app:perron}. I subsequently reviewed and verified these
calculations by hand. For a few technically difficult examples, an assistant proposed the first draft of an intermediate argument,
which I then verified, reconstructed and rewrote.

The responsibility is mine, for every statement printed here and, above all, for its status.
Whether a result is proved, conditional, numerical, conjectural or open is a claim I make myself,
those not proved are identified as such in the text, and the register at the end collects them.
What an assistant produces is a proposal to be examined, never a mathematical authority.

The second volume, in preparation at the time of writing, uses these tools differently. There
they help me develop the specialized theory around the Riemann hypothesis and the Ingham function,
first by recording the results obtained on the gauged Ingham function, which appeared
in~\cite{CloitreFloor}, then by taking me into areas of mathematics more abstract than those of
the present volume, category theory, algebraic geometry and the history of the Weil conjectures.
That is what oriented the work now under way toward the Riemann hypothesis, its generalized form
and intermediate questions such as the existence of Siegel zeros. What comes of it belongs to the
second volume.

\bigskip

\chapter*{Prologue}
\addcontentsline{toc}{chapter}{Prologue}
\markboth{PROLOGUE}{}

\rafepigraph{Rien n'est plus fécond, tous les mathématiciens le savent, que ces obscures analogies, ces troubles reflets d'une théorie à une autre, ces furtives caresses, ces brouilleries inexplicables.}{Nothing is more fruitful, as every mathematician knows, than those obscure analogies, those hazy reflections of one theory in another, those furtive caresses, those inexplicable rifts.}{André Weil\index[names]{Weil, A.}, \emph{De la métaphysique aux mathématiques} (1960)~\cite{WeilMetaphysique1960}}

I did not set out to build a theory. What became one began around 2010, with an elementary
identity about the Liouville function\index[terms]{Liouville function}. Write $\lambda(k)$ for the value $+1$ or $-1$ according to
the parity of the number of prime factors of $k$, counted with multiplicity. The identity is

\[
\sum_{k=1}^{n}\lambda(k)\left\lfloor\frac{n}{k}\right\rfloor=\lfloor\sqrt n\rfloor,
\]

and it is easy to prove and easy to check by hand.

I immediately noticed the same exponent in two different places. On the left the terms
$\lambda(k)$ are weighted by $\lfloor n/k\rfloor$, and the sum is known exactly, of size the
square root of $n$. Remove the weight and take the partial sums $\lambda(1)+\cdots+\lambda(n)$ of
the same terms. The assertion that these partial sums are $O(n^{1/2+\eps})$ for every
$\eps>0$ is equivalent to the Riemann hypothesis\index[terms]{Riemann hypothesis}, and no one knows whether that holds. The
exponent one half sat on both sides. On one side it was an identity anyone could verify. On the
other it was a question that has stood open for more than a century.

The two quantities differ by a single thing, the weight. The weighted sum carries the exact rate.
The plain one carries the hypothesis. So I asked what the weight knew that the plain sum did not,
and whether the exact control of the one could be carried over to the other. To recover a plain
sum from a weighted one is a Tauberian\index[terms]{Tauberian} problem, and I began to look for a Tauberian theorem that
used this particular weight.

I found the framework in Ingham\index[names]{Ingham, A. E.}'s 1945 paper. What Ingham averages against is not
$\lfloor n/k\rfloor$ itself but its normalization $\Phi(k/n)=(k/n)\lfloor n/k\rfloor$. After
division by $n$, and after the factor $1/k$ is absorbed into the sequence, the same identity takes
the form of an equation with Ingham's kernel. Ingham uses that kernel as an element of a summation
method. The theory built here fixes it instead as the kernel of a family of equations in which the
prescribed power varies, and reading it that way changed the shape of the question. The exponent
one half was no longer an isolated bound in the complex plane. It was a number attached to an
operator, the point at which the operator stops passing a signal through and begins to absorb it.

Once the weight was read as a kernel, the question stopped depending on the Liouville function
that had raised it. It could be put to any kernel, and its answer became a property of the kernel
itself. That is where this book begins.

\bigskip

\chapter*{Notations and Conventions}
\addcontentsline{toc}{chapter}{Notations and Conventions}
\markboth{NOTATIONS AND CONVENTIONS}{}

The following pages collect the notation used throughout the monograph. They are divided into five parts, namely standard mathematical symbols and asymptotics, reserved letters, classical arithmetic functions, a glossary of the principal objects of the theory, and the acronyms. The glossary entries give precise references to the chapters where each object is first defined.

Three conventions hold throughout. A sequence indexed by the integers carries a subscript, as in $a_n$ or $u_n$. A function carries its argument, as in $f(x)$ or $\mu(n)$. The partial sums of a sequence form a summatory function, $A(x)=\sum_{n\le x}a_n$, and $A(n)$ is that function at an integer argument.

\section*{I. Standard Symbols and Asymptotics}

\begin{itemize}
    \item $\N$ (natural numbers, starting from $1$), $\Z$ (integers), $\R$ (reals), $\C$ (complex numbers).
    \item $\lfloor x \rfloor$: greatest integer $\le x$ (floor function\index[terms]{floor function}). $\{x\} := x - \lfloor x \rfloor$: fractional part\index[terms]{fractional part}.
    \item $x^{\underline{m}} := x(x-1)\cdots(x-m+1)$: falling factorial ($x^{\underline{0}} := 1$).
    \item $\Gamma(z)$: Euler Gamma function\index[terms]{Gamma function}. $\psi(z) := \Gamma'(z)/\Gamma(z)$: digamma\index[terms]{digamma function} function.
    \item $f(x) \sim g(x)$: $\lim_{x\to\infty} f(x)/g(x) = 1$.
    \item $f(x) = \mathcal{O}(g(x))$, equivalently $f(x) \ll g(x)$: there exists $C > 0$ with $|f(x)| \le C|g(x)|$ for all large $x$. The subscript $\ll_\eps$ means the constant $C$ depends on $\eps > 0$.
    \item $f(x) = o(g(x))$: $f(x)/g(x) \to 0$.
    \item $a \star b$: Dirichlet convolution, $(a \star b)(n) = \sum_{d | n} a(d)b(n/d)$.
    \item $(F \star H)(x) := \int_1^x F(x/y)\,H(y)\,dy/y$: multiplicative convolution on $[1,\infty)$.
    \item $(u*v)(x) := \int_{\R}u(x-y)\,v(y)\,dy$: additive convolution on $\R$, Chapters~\ref{chap:boundary} and~\ref{chap:ortho}. $\widehat u(\xi) := \int_{\R}u(v)e^{-i\xi v}\,dv$: Fourier transform, Chapter~\ref{chap:boundary} and Appendix~\ref{app:boundary}.
\end{itemize}

\section*{II. Reserved Letters}

The following conventions are maintained throughout.

\begin{center}
\begin{tabular}{l l p{9cm}}
\toprule
Letter(s) & Type & Role \\
\midrule
$n, k, m, d, j$ & integers & summation and index variables \\
$p$ & prime & prime number \\
$x, y, t$ & reals & continuous variables \\
$z = \sigma + it$ & complex & argument of Arithmetic Mellin Transforms $G^*(z)$ \\
$s = \sigma + it$ & complex & argument of classical Dirichlet series and $\zeta(s)$ \\
$\beta$ & real & exponent in the defining equation $A_G(n) = n^{-\beta}$ \\
$\lambda$ & kernel parameter & $\lambda\in(0,1)$ for the affine kernel, $\lambda>1$ for the broken harmonic family $g_\lambda$. The Liouville function is written $\lambda(n)$, always with an argument \\
$\eps$ & real $> 0$ & arbitrarily small, implicit in $\mathcal{O}_\eps$ and $\ll_\eps$ \\
$f$ & gauge & an increasing unbounded function through which a kernel is re-read, Chapter~\ref{chap:gauge} \\
$q$ & integer $\ge2$ & base of the exponential gauge $f(x)=q^{x}+1$ of Chapter~\ref{chap:gauge_ingham} and of the epilogue, where it is also the cardinality of the finite field $\mathbb F_q$ when $q$ is a prime power \\
$Q$ & integer $\ge2$ & modulus of the Dirichlet characters of Chapter~\ref{chap:boundary} and Appendix~\ref{app:boundary}, fixed in each statement, distinct from the gauge base $q$ \\
\bottomrule
\end{tabular}
\end{center}

\section*{III. Classical Arithmetic Functions}

\begin{itemize}
    \item $\mu(n)$: Möbius function. $\mu(1)=1$; $\mu(n) = (-1)^k$ if $n$ is a product of $k$ distinct primes; $\mu(n)=0$ otherwise.
    \item $\Lambda(n)$: von Mangoldt function\index[terms]{von Mangoldt function}. $\Lambda(n) = \log p$ if $n = p^k$, $0$ otherwise.
    \item $\varphi(n)$: Euler totient\index[terms]{Euler totient}. $\tau(n) = \sum_{d|n}1$: number of divisors.
    \item $\lambda(n) = (-1)^{\Omega(n)}$: Liouville function\index[terms]{Liouville function}, where $\Omega(n)$ counts prime factors with multiplicity.
    \item $\chi_4$: non-principal Dirichlet character\index[terms]{Dirichlet character} modulo 4: $\chi_4 = (1, 0, -1, 0, 1, 0, -1, \ldots)$.
    \item $M(x) := \sum_{n \le x}\mu(n)$: Mertens function\index[terms]{Mertens function}\index[names]{Mertens, F.}. $M_{-1}(x) := \sum_{n\le x}\mu(n)/n$: weighted Mertens function.
    \item $\zeta(s) := \sum_{n\ge 1} n^{-s}$: Riemann zeta function ($\Re(s) > 1$, meromorphically continued).
    \item $L(s,\chi) := \sum_{n\ge 1}\chi(n)n^{-s}$: Dirichlet $L$-function. $\beta(s) := L(s,\chi_4)$: Dirichlet beta function.
    \item $F_{n}$: Fibonacci\index[terms]{Fibonacci sequence}\index[names]{Fibonacci} numbers ($F_{1}=F_{2}=1$). $L_{n}$: Lucas numbers.
\end{itemize}

\section*{IV. Glossary of Principal Objects}

The objects below are the central protagonists of the theory. Each entry gives the symbol, a brief description, and the chapter of first definition.

\medskip
\noindent\textbf{The defining equation.}
\begin{equation*}
\sum_{k=1}^n a_k\, G(n,k) = n^{-\beta}, \qquad n \ge 1, \quad \beta \in \R.
\end{equation*}
The sequence $(a_n)_{n\ge 1}$ is the unknown, fixed by the equations through $a_1 = 1/G(1,1)$, equal to $1$ when $G(1,1)=1$. The function $G : \N^* \times \N^* \to \R$ is the kernel. The exponent $\beta$ parametrizes the family of problems. It is fixed by the context and left implicit in $a_n$ and $A(n)$. The subscript, as in $A_\beta(n)$, appears only where several forcing exponents are compared at once, and each such passage says so.

\medskip
\renewcommand{\arraystretch}{1.4}
\begin{longtable}{p{2.8cm}p{8.5cm}p{2.2cm}}
\toprule
\textbf{Symbol} & \textbf{Description} & \textbf{Defined in} \\
\midrule
\endfirsthead
\toprule
\textbf{Symbol} & \textbf{Description} & \textbf{Defined in} \\
\midrule
\endhead
$(a_n)$, $A(x)$ & Unknown sequence; partial sums $A(x) := \sum_{n\le x}a_n$, written $A(n)$ at integer arguments & \eqref{eq:defining_relation} \\
$A_1(x)$ & Weighted partial sums $A_1(x) := \sum_{n\le x}na_n$ & Ch.~\ref{chap:equivalence} \\
$G(n,k)$ & \textbf{Regular Arithmetic Function (RAF)}: a kernel of two integer variables satisfying the RAF asymptotic dichotomy, a finite index, with transparency below it and absorption at and above it & Def.~\ref{def:reg_index} \\
$g(x)$ & \textbf{Function of Good Variation (FGV)}: a single-variable profile on $(0,1]$ inducing the RAF $G(n,k) = g(k/n)$ & Def.~\ref{def:reg_index_fgv} \\
$\tau(G)$ & \textbf{Transparency threshold}: $\tau(G)=\sup\{a:\text{every }\beta<a\text{ is transparent}\}$, defined for every kernel and asserting nothing about absorption & Def.~\ref{def:reg_index_fgv} \\
$\alpha(G)$ & \textbf{Regularity Index} (arithmetic index): equal to $\tau(G)$ when the kernel is a RAF, that is when $\tau(G)$ is finite and the partial sums are absorbed above it ($A(n)=\mathcal{O}(n^{-\alpha+\eps})$). Transparency below the threshold reads $A(n)=\Xi_G(\beta)\,n^{-\beta}+o(n^{-\beta})$, with $\Xi_G=1/G^*$ wherever that identification is proved, zero at the poles of $G^*$, see Open Problem~\ref{op:xi_reciprocal}. Writing $\alpha(G)$ asserts that $G$ is a RAF & Def.~\ref{def:reg_index_fgv}, \ref{def:reg_index} \\
$g^*(z)$ & \textbf{Arithmetic Mellin Transform of a profile}: $g^*(z) = -z\int_0^1 g(t)t^{-z-1}\,dt$, initially for $\Re(z)<0$, extended meromorphically. See normalization note below. & Def.~\ref{def:mellin} \\
$G^*(z)$ & \textbf{Arithmetic Mellin Transform of a kernel}: the limit of the finite probes, $\lim_n(-z/n)\sum_{k\le n}G(n,k)(k/n)^{-z-1}$, which need not be given by an integral & Def.~\ref{def:reg_index} \\
$\eta(G)$ & \textbf{Analytic Index}: $\eta(G) = \inf\{\Re(\rho) : G^*(\rho)=0\}$, undefined when $G^*$ has no zero & Def.~\ref{def:mellin} \\
$\tau(g)$, $\alpha(g)$, $\eta(g)$ & The same three quantities read on the profile, used whenever the kernel is induced by one, $G(n,k)=g(k/n)$. The letter records which object carries them & Def.~\ref{def:reg_index} \\
$\Phi(x)$ & \textbf{Ingham function}: $\Phi(x) := x\lfloor 1/x\rfloor$. Its regularity index satisfies $\alpha(\Phi)=1/2 \iff \mathrm{RH}$ & Ch.~\ref{chap:fgv}, Thm.~\ref{thm:tauberian_rh} \\
$\nu$, $\jmath$ & \textbf{Jump measure and diagonal jump} of a profile: $\nu$ is the Lebesgue--Stieltjes measure of $g$ on $(0,1)$, with an absolutely continuous part and point masses at the jumps of $g$, and $\jmath=g(1)-g(1^{-})$ is the jump at the diagonal, which the transform does not see & Thm.~\ref{thm:xi_reciprocal} \\
$S_G(n)$, $E_g(n)$ & \textbf{Kernel sum and Abelian defect}: $S_G(n):=\sum_{k\le n}G(n,k)$, of density $G^{*}(-1)$, and the second order term $E_g(n)$ of the profile sum & Def.~\ref{def:kernel_sums}, Cor.~\ref{cor:density_defect} \\
$A_G(n)$ & Averaged sums, the discrete arithmetic operator $A_G(n) := \sum_{k=1}^n a_k G(n,k)$. Written $A_g(n)$ for a univariate profile and $A_\Phi(n)$ for the Ingham weight & Ch.~\ref{chap:fgv} \\
HLR criterion & Hardy-Littlewood-Ramanujan criterion: $A_G(n)=n^{-\beta}\Rightarrow na_n=o(n^{\eps})\ \forall\eps>0$; strong form $na_n=\mathcal{O}(1)$ & Def.~\ref{def:HLR} \\
$w_g$ & \textbf{Wiener companion} of a profile normalized at the origin, the integrable function $e_++h_g-e_+*h_g$ on $[0,\infty)$, with $e_+(v)=e^{-v}$ and $h_g(v)=g(e^{-v})-1$, whose Fourier transform is $g^{*}(-i\xi)/(1+i\xi)$ & \eqref{eq:wiener_companion} \\
$\nu_\Phi$, $R_\Phi$ & Perturbation measure and resolvent measure of the Ingham operator in the Volterra-Stieltjes decomposition & Def.~\ref{def:nu_Phi} \\
$G(N,k)$, $K(N,k)$ & Primitive triangular kernel and its Abel difference $K(N,k)=G(N,k)-G(N,k+1)$ in the discrete setting & Def.~\ref{def:discrete_resolvent} \\
$\mathcal{G}_N(z)$ & Finite arithmetic Mellin probe at level $N$; converges to $G^*(z)$ as $N\to\infty$ & Def.~\ref{def:probes} \\
$\alpha_f(G)$ & \textbf{Gauged regularity index}: regularity index of $G$ under the gauge $f$, defined by the deformed equation $A_G(f(n),f(\cdot)) = f(n)^{-\beta}$ & Def.~\ref{def:stability_equilibrium} \\
$G^*_f(z)$ & Arithmetic Mellin Transform of $G$ with respect to the gauge $f$ & Def.~\ref{def:gauge_probe} \\
$\tau_{\mathrm{tr}}(G)$ & \textbf{Transparency frontier}: supremum of the exponents $c$ such that every $\beta<c$ is transparent & Ch.~\ref{chap:diophantine} \\
$\Xi_G(\beta)$ & \textbf{Additive transparency coefficient}: $\Xi_G=1/G^*$, the constant of the transparent regime & Def.~\ref{def:transparency_fgv} \\
\bottomrule
\end{longtable}
\renewcommand{\arraystretch}{1}

\section*{V. Acronyms}

\begin{center}
\begin{tabular}{lll}
\toprule
Acronym & Full form & First used \\
\midrule
RAF & Regular Arithmetic Function & Def.~\ref{def:reg_index} \\
FGV & Function of Good Variation & Def.~\ref{def:reg_index_fgv} \\
HLR & Hardy-Littlewood-Ramanujan (criterion) & Def.~\ref{def:HLR} \\
RH & Riemann hypothesis & Ch.~1 \\
\bottomrule
\end{tabular}
\end{center}

\section*{Note on the normalization of $g^*(z)$}

The definition $g^*(z) = -z\int_0^1 g(t)t^{-z-1}\,dt$ differs from the classical Mellin transform\index[terms]{Mellin transform} $g^{\sharp}(z) = \int_0^1 g(t)t^{-z-1}\,dt$ used in the earlier paper \cite{Cloitre2016}. The $-z$ factor has three justifications. First, on the geometric side, $g^*(z) = \int_0^1 g(t)\,d(t^{-z})$ is a Stieltjes integral. Second, it removes the pole at $z=0$ present in $g^{\sharp}$ when $g(0^+)\neq 0$, giving $g^*(0) = g(0^+)$. Third, it produces the natural asymptotic constant $A(n)\sim n^{-\beta}/g^*(\beta)$ in the transparent case, without spurious signs or factors of $\beta$. The zeros of $g^*$ and $g^{\sharp}$ coincide for $z\neq 0$, so all spectral results from \cite{Cloitre2016} carry over unchanged. The star $\star$ is reserved throughout for convolution. The convention of \cite{CloitreOrtho} and of this monograph is $g^*(z) = -z\int_0^1 g(t)t^{-z-1}\,dt$ throughout.

\bigskip

\mainmatter

\part{From Ingham's Tauberian theorem to functions of good variation}
\label{part:one}

\rafepigraph{Là où Lagrange voyait des analogies, nous voyons des théorèmes.}{Where Lagrange saw analogies, we see theorems.}{André Weil, \emph{De la métaphysique aux mathématiques} (1960)~\cite{WeilMetaphysique1960}}

This part lays the foundations of the theory. It opens with Ingham's Tauberian\index[terms]{Tauberian} theorem and
the summation method it produced, builds from a single arithmetic coincidence the concept of
a function of good variation\index[terms]{function of good variation}, proves the equivalence between the Riemann hypothesis\index[terms]{Riemann hypothesis} and the
value one half of the index of the Ingham function\index[terms]{Ingham function} $\Phi(x)=x\lfloor 1/x\rfloor$, and closes
with a first gallery of examples. Between the definition and the equivalence it reads the
weight at the boundary of its transform, where it returns limits and the prime number theorem,
so that the rate is asked for once the limit is understood. The coincidence is the one recounted
in the prologue, and understanding it is what the five chapters are for.

\chapter{Ingham's Tauberian theorem and summation method}
\label{chap:ingham}

Ingham\index[names]{Ingham, A. E.} proved in 1945 a Tauberian\index[terms]{Tauberian} theorem and drew from it a proof of the prime number
theorem. Read for that conclusion alone, the result stands in the line of Littlewood\index[names]{Littlewood, J. E.}, one more
approach to a theorem already established. Beneath the method lies a fixed weight, the averaging
of a sequence against a single function of the ratio $k/n$, and the theory of this volume begins
by reading that weight as the kernel of an equation rather than as a device for summing a series.
This chapter sets out what the method was, how it was studied in the decades after 1945, and why
that body of work, Segal\index[names]{Segal, S. L.}'s approach to the Riemann hypothesis\index[terms]{Riemann hypothesis}\index[terms]{Riemann hypothesis} included, did not take the step the
theory takes. Following the weight past the value the averages approach and into the rate at which
they fall away is the work of this part, and with it the foundation the volume rests on.

\section{Before Ingham: Eratosthenian averages}\label{sec:eratosthenian}

The averaging that carries Ingham's theorem\index[terms]{Ingham theorem} was not new in 1945. Its arithmetical form goes back
to Wintner\index[names]{Wintner, A.}, who studied mean values of arithmetical functions through a divisor-sum averaging and
named it Eratosthenian, after the sieve that the divisor sum carries \cite{Wintner1943}. Given a
sequence $(a_n)$, the average in question is
\begin{equation}\label{eq:ingham_averages}
I(x)=\frac1x\sum_{n\le x}\sum_{d\mid n}d\,a_d
   =\frac1x\sum_{k\le x}k\,a_k\Big\lfloor\frac xk\Big\rfloor,
\end{equation}
the two forms equal because each $k\le x$ divides $\lfloor x/k\rfloor$ of the integers up to $x$.
The inner divisor sum is the sieve made arithmetic, and it is this object, the weight
$k\lfloor x/k\rfloor$ read as a function of the ratio $k/x$, that the later theory will isolate.

Wintner\index[names]{Wintner, A.}'s subject was the mean value of an arithmetical function, the limit of
$\frac1x\sum_{n\le x}f(n)$, and the Eratosthenian average is the divisor-sum form of that mean,
in which the plain average is replaced by the average of $\sum_{d\mid n}d\,a_d$. The name records
the sieve of Eratosthenes, whose crossing out of multiples is the divisor structure the inner sum
carries. The averaging existed, then, as a chapter in the theory of mean values, before it became
the vehicle of a Tauberian\index[terms]{Tauberian} theorem.

Two things are worth separating at the outset. One is the averaging device, the Eratosthenian
mean, an object of the theory of arithmetical functions with a history of its own. The other is
the Tauberian\index[terms]{Tauberian} theorem that Ingham attached to it, the return from the average to the series, which
gave the device a reach it did not have as an averaging alone. Hardy, ordering the summability
procedures of analysis, attached Ingham's name to the method and fixed the notation that became
usual \cite[App.~IV]{Hardy1949}. The name that survives, Ingham summability, records the theorem
rather than the average, and the account that follows keeps the two apart, the mean from Wintner\index[names]{Wintner, A.}
and the Tauberian\index[terms]{Tauberian} return from Ingham.

\section{Ingham's Tauberian theorem}\label{sec:ingham_theorem}

In Ingham's notation the series $\sum a_n$ is $I$-summable to $A$ when $I(x)\to A$ as
$x\to\infty$. A method of this kind assigns a value to a series through weighted averages of its
terms, and the substance of Ingham's\index[names]{Ingham, A. E.} 1945 paper is the converse, a return from the averages to the
series under one side condition \cite{Ingham1945}.

\begin{theorem}[Ingham, 1945]\label{thm:ingham}
Let $(a_n)$ satisfy the one-sided bound $na_n\ge -C$ for some constant $C\ge0$. If
$I(x)\to A$ as $x\to\infty$, then the series converges to the same value,
\[
\sum_{n=1}^{\infty}a_n=A .
\]
\end{theorem}
The statement is that of \nm{Korevaar}{J.}~\cite[p.~110]{Korevaar2004} and of the original paper. The one-sided
bound is the Tauberian side condition, the single restraint that turns summability back into
convergence, and it cannot be dropped, since without it the averages can settle while the series
swings. The theorem belongs to the tradition of Tauber\index[names]{Tauber, A.} and Littlewood, where a summability
statement is converted to convergence at the cost of a growth restriction on the terms, and what
is Ingham's own is the weight, the arithmetic divisor sum in place of the analytic kernels of the
classical converse theorems.

The difficulty of the return is arithmetic. The averaged data $I(x)$ carry the divisor sum
$\sum_{d\mid n}d\,a_d$, and the individual terms are recovered from that data by M\"obius
inversion, since $\mu\star1=e$ inverts the divisor convolution. The depth of the theorem is
therefore the size of the M\"obius sum\index[terms]{M\"obius function}, and the analytic weight of turning a statement about
averages back into a statement about the series is carried entirely by cancellation in
$\sum_{n\le x}\mu(n)$. A theorem about summing a series meets the zeta function at exactly this
point, because the size of that sum is a statement about the zeros.

Ingham placed the method among the Ces\`aro\index[terms]{Ces\`aro mean}\index[names]{Ces\`aro, E.} means. He proved that $I$-summability implies
Ces\`aro summability $(C,\delta)$ for every $\delta>0$, and that $(C,-\delta)$-summability
implies $I$-summability for $0<\delta<1$, while $I$-summability is not comparable with ordinary
convergence. The one implication he published in full, $I\Rightarrow(C,1)$, already leans on the
size of the M\"obius function, on the bound $\sum_{n\le x}\mu(n)/n=O\big((\log x)^{-k}\big)$, and
the stronger implications call on it for every $k$. That a method built from divisor sums should
ask for the M\"obius function is the first sign that it reaches toward the zeta function, whose
zeros are what the size of $\sum_{n\le x}\mu(n)$ records.

\section{The prime number theorem}\label{sec:ingham_pnt}

The application Ingham drew from the theorem is the prime number theorem, and it fixes what the
method reaches at a constant limit. The natural sequence is the one whose Dirichlet series\index[terms]{Dirichlet series} is the
logarithmic derivative of the zeta function. Take $a_n=(\Lambda(n)-1)/n$, with $\Lambda$ the
von Mangoldt function\index[terms]{von Mangoldt function},
\[
-\frac{\zeta'(s)}{\zeta(s)}=\sum_{n\ge1}\frac{\Lambda(n)}{n^{s}},
\]
the subtracted $1$ removing the pole of $\zeta$ at $s=1$ so that the averaged series can settle at
a finite value. The averages converge,
\[
I(x)=\frac1x\sum_{k\le x}\bigl(\Lambda(k)-1\bigr)\Big\lfloor\frac xk\Big\rfloor
\ \longrightarrow\ -2\gamma\qquad(x\to\infty),
\]
with $\gamma$ the Euler constant\index[terms]{Euler constant}. The limit is elementary. Summing $\Lambda$ over
divisors gives $\sum_{k\le x}\Lambda(k)\lfloor x/k\rfloor=\sum_{m\le x}\log m=\log\lfloor x\rfloor!$,
while $\sum_{k\le x}\lfloor x/k\rfloor=\sum_{m\le x}\tau(m)$, so that
\[
x\,I(x)=\log\lfloor x\rfloor!-\sum_{m\le x}\tau(m).
\]
Stirling\index[terms]{Stirling's formula}\index[names]{Stirling, J.}'s formula gives
$\log\lfloor x\rfloor!=x\log x-x+\mathcal{O}(\log x)$, Dirichlet's formula\index[terms]{Dirichlet's formula} for the
divisor sum gives $\sum_{m\le x}\tau(m)=x\log x+(2\gamma-1)x+\mathcal{O}(\sqrt x)$, and the two
principal terms cancel, leaving $x\,I(x)=-2\gamma x+\mathcal{O}(\sqrt x)$. The one-sided bound
holds since $\Lambda(n)\ge0$ gives $na_n=\Lambda(n)-1\ge-1$. Theorem~\ref{thm:ingham} returns this value to
the series,
\[
\sum_{n\ge1}\frac{\Lambda(n)-1}{n}=-2\gamma,
\]
which by partial summation\index[terms]{partial summation} is the prime number theorem in the form $\psi(x)\sim x$, see \nm{Tenenbaum}{G.}~\cite[II.4]{Tenenbaum2015}, where $\psi(x)=\sum_{n\le x}\Lambda(n)$ is Chebyshev\index[terms]{Chebyshev function}\index[names]{Chebyshev, P. L.}'s function. The
one analytic ingredient is the non-vanishing of the zeta function on the line $\Re s=1$, a
statement equivalent to $\sum_{n\le x}\mu(n)=o(x)$. It enters through the Tauberian return of
Theorem~\ref{thm:ingham} and not through the averages, whose convergence the computation above
makes elementary. Chapter~\ref{chap:boundary} proves that return and shows where the
non-vanishing is used. At
a constant limit the method reaches the prime number theorem and no more, and the finer behavior
it can carry appears only when the average is allowed to decay.

\section{The method among other summability procedures}\label{sec:ingham_posterity}

After 1945 the method was studied as a summability procedure, and the questions asked of it were
the questions one asks of such a procedure, inclusion between methods, comparison with neighbouring
means, and the side conditions under which summability returns convergence. Through all of it the
weight stays fixed and the target stays a constant.

Pennington\index[names]{Pennington, W. B.} compared Ingham summability with summability by Lambert series\index[terms]{Lambert series}\index[names]{Lambert, J. H.}, relating the oscillation
of the Lambert transform\index[terms]{Lambert transform} $L_\lambda(t)$ as $t\to0^{+}$ to that of the Ingham transform\index[terms]{Ingham transform} as
$x\to\infty$ \cite{Pennington1955}. Rajagopal\index[names]{Rajagopal, C. T.} supplied the Tauberian counterpart, an Abelian
theorem carrying the oscillation of the Lambert transform to the Ingham transform and Tauberian
converses recovering convergence under monotonicity and moment conditions on the kernel, built on
Pennington\index[names]{Pennington, W. B.}'s lemmas \cite{Rajagopal1955}. The comparison places the Lambert and Ingham methods
beside the Abel and Riesz families to which they answer. Segal\index[names]{Segal, S. L.} opened his own study of the method
in the same setting \cite{Segal1966}, and Jukes\index[names]{Jukes, K. A.} set the Ingham method against the $(D,h(n))$
methods \cite{Jukes1971}, a further comparison within the same frame.

The sharpest statement of the method as a summability procedure came later. Zacharovas\index[names]{Zacharovas, V.}
characterized $I$-summability exactly. The average $\frac1n\sum_{k\le n}a_k\lfloor n/k\rfloor$
converges to $A$ if and only if $\sum_{k\le n}a_k\lfloor n/k\rfloor\log k=o(n\log n)$ and the
series is Abel summable to the same $A$, and he located the method next to the $(A,\log n)$ mean,
as Ces\`aro summability sits next to $(A,n)$ \cite{Zacharovas2011}. The condition is a weighted
logarithmic bound on the same divisor-averaged sum, and it fixes how tightly the method sits
against Abel summability, the two agreeing exactly when that bound holds. This is the method
understood on its own terms, a necessary and sufficient condition for the averages to reach a
constant.

A unified account of the whole line, from Wintner\index[names]{Wintner, A.}'s averages through Ingham's theorem to Segal\index[names]{Segal, S. L.}'s
smoothed weights, was later given by Kanemitsu\index[names]{Kanemitsu, S.}, Kuzumaki\index[names]{Kuzumaki, T.} and Tanigawa\index[names]{Tanigawa, Y.} in the language of Stieltjes\index[names]{Stieltjes, T. J.}
integration, where partial summation and integration by parts become one process and the divisor
sum is handled as a Stieltjes measure \cite{Kanemitsu2018}. It is the most complete modern
treatment of the method, and it treats it throughout as a summability procedure. What none of
these studies asks, from Wintner\index[names]{Wintner, A.} to the Stieltjes synthesis, is what happens when the constant
target is replaced by a prescribed rate of decay, the weight held exactly where Ingham left it.

\section{Segal's route toward the Riemann hypothesis}\label{sec:segal_rh}

\index[terms]{Segal's route}\index[names]{Segal, S. L.}One line in this literature reached past the primes toward the zeros. Segal asked in 1975 whether
Ingham summability registers the zeros of the zeta function as well as its primes, and set out to
link the method to the Riemann hypothesis\index[terms]{Riemann hypothesis} \cite{Segal1975}. The question is worth recording with
care, because a reader may object that the method had already been tied to the hypothesis. It had,
but by another route and toward another invariant, and that difference is what the present volume
turns on.

Segal\index[names]{Segal, S. L.}'s device was to replace the arithmetic weight $1-v\{1/v\}$ by a smooth Riesz weight
$(1-v)^{\delta}$, passing from the divisor sum to a family of continuous means, and to tie the
strength of the resulting implication to the error term in the prime number theorem, and so to the
size of $\sum_{n\le x}\mu(n)$ and to the zero-free region\index[terms]{zero-free region} of the zeta function. The reconstruction
of Kanemitsu\index[names]{Kanemitsu, S.}, Kuzumaki\index[names]{Kuzumaki, T.} and Tanigawa\index[names]{Tanigawa, Y.}, in the language of Stieltjes\index[names]{Stieltjes, T. J.} integration, makes the mechanism
explicit. The transfer runs from $I$-summability to summability by the smooth weight, valid for
every weight in a class of twice-differentiable decreasing functions, and the strength of the
conclusion is set by the size of the M\"obius sum in two regimes, as in \nm{Kanemitsu}{S.}, \nm{Kuzumaki}{T.} and \nm{Tanigawa}{Y.}~\cite{Kanemitsu2018}.
Unconditionally the reducing factor is the one the classical zero-free region provides, of the
shape $\exp(-c(\log v)^{3/5-\eps})$. Under a bound $\sum_{n\le x}\mu(n)=O(x^{b})$ with $b<1$, a
weak form of the hypothesis, the transfer sharpens to a power saving. What the route reaches, in
the second regime, is the existence of a zero-free strip $\Re s>a$ with $a<1$, controlled by the
M\"obius sum, and not the Riemann hypothesis\index[terms]{Riemann hypothesis} itself. A simplification Segal\index[names]{Segal, S. L.} added was withdrawn the following year, leaving the
implication $I\Rightarrow(C,1)$ as it stood \cite{Segal1976}.

The same line marked the edge of the method from the other side, through limitation theorems, the
necessary conditions that the terms of any $I$-summable series must satisfy. Erd\H{o}s\index[names]{Erd\H{o}s, P.} and Segal\index[names]{Segal, S. L.}
showed that the bound of order $\log\log n$ on the terms is best possible, constructing an
$I$-summable series that approaches it arbitrarily slowly \cite{ErdosSegal1978}, and Segal\index[names]{Segal, S. L.} showed
shortly after that a companion bound on the partial sums, $\sum_{n\le x}a_n=o(\log x)$, is best
possible as well \cite{Segal1979}. These results fix how far a summable series can be stretched
while the method still returns it, and they leave untouched the object the method acts on. Segal\index[names]{Segal, S. L.}'s
route smooths the weight, and the smoothing is what ties the outcome to a remainder term rather
than to an exact threshold. Keeping the weight exactly as $x\lfloor 1/x\rfloor$ and asking not for
a limit but for a prescribed rate of decay is the step that turns the weak link into an
equivalence, and it is taken in Chapter~\ref{chap:equivalence}. The limitation of Erd\H{o}s\index[names]{Erd\H{o}s, P.} and
Segal\index[names]{Segal, S. L.} returns at the criterion of \S\ref{sec:boundary_problem}, where it marks the boundary of what the
method controls.

\section{What the literature did not ask}\label{sec:not_asked}

The record of the previous sections is of a method studied thoroughly and from many sides. It was
compared with Lambert series and with Riesz and Ces\`aro means, its Tauberian side conditions were
found and shown to be best possible, and it was carried, in Segal\index[names]{Segal, S. L.}'s hands, to the edge of the
Riemann hypothesis\index[terms]{Riemann hypothesis}. A summability method has a matrix, a weight, a kernel, and this one was handled
by authors who knew its weight. The change of viewpoint this volume takes is therefore not that the
weight went unnoticed.

The change is more precise. No treatment in this literature takes the fixed weight as the primary
object of a family of equations forced by prescribed powers, studies the ordinary partial sums of
their solutions as the forcing exponent varies, or attaches to the resulting operator an index
marking the passage between transparency\index[terms]{transparency} and absorption\index[terms]{absorption}. The method was read to answer one
question, whether it sums a series to its value. The theory that follows reads the same weight to
answer another, how the response of the equation varies when the imposed rate changes. The four
displacements are set side by side below.

\begin{center}
\begin{tabular}{ll}
\toprule
Classical literature & The theory that follows \\
\midrule
a summation method & a kernel taken as object \\
a given limit & a family of forcings \\
comparison of methods & the partial sums of the solution \\
Tauberian side conditions & an index of transparency and absorption \\
\bottomrule
\end{tabular}
\end{center}

\noindent The left column is where the literature of this chapter stands. The right column is where
the next chapters go. The passage from one to the other keeps the weight of Wintner\index[names]{Wintner, A.} and Ingham
exactly, and changes only the question asked of it.

\section{From a summation method to a kernel}\label{sec:method_to_kernel}

The step out of this chapter is short to state and carries the whole of what follows. The weight
$x\lfloor 1/x\rfloor$ that Wintner\index[names]{Wintner, A.} averaged and Ingham summed is a single function of one variable.
Read it not as the matrix of a summation but as a fixed function against which a sequence is
summed, and let the value the averages are asked to approach be replaced by a prescribed power. The
function stays the Ingham weight. The constant target becomes a free parameter. What was a method
for assigning a value to a series becomes an equation whose solution has a rate of decay to study.

Written out, the weight is the Ingham function\index[terms]{Ingham function}
\begin{equation}\label{eq:ingham_function_ch1}
\Phi(x)=x\Big\lfloor\frac1x\Big\rfloor,\qquad 0<x\le1,
\end{equation}
and it is worth pausing on it, since the rest of the volume is about this one function. On each
interval between two consecutive reciprocals of integers it is a straight segment through the
origin. At every $1/m$ it returns exactly to the value one and drops at once after, so the graph
is the fan of segments of Figure~\ref{fig:ingham}, and it never falls below one half. The return
to one at $x=1$ is not a detail. It is what makes the diagonal of the triangular array below
nonzero, and hence what lets the equation determine its solution with no hypothesis of any kind.
Why this weight and not another is what \S\ref{sec:two_identities} answers, with three classical
identities.

The function is classical. $\Phi(1/t)=\lfloor t\rfloor/t$ is the kernel \nm{Bingham}{N. H.} and \nm{Inoue}{A.} call the
P\'olya\index[names]{P\'olya, G.} kernel and trace to P\'olya's 1917 paper
\cite{Polya1917}. See also \cite[p.~189]{BinghamInoue2000b} and
\cite[\S6.4.4]{Bingham1989}. Its classical Mellin transform is
$\zeta(z+1)/(z+1)$ on $\Re z>0$ \cite[p.~189]{BinghamInoue2000b}. Moreover,
$\Phi(k/n)=(k/n)\lfloor n/k\rfloor$, and the Eratosthenian average
\eqref{eq:ingham_averages} of a sequence $(a_k)$ at $x=n$ is exactly
$\sum_{k\le n}a_k\Phi(k/n)$. Replacing the constant target by a prescribed power gives
\begin{equation}\label{eq:defining_relation_ch1}
\sum_{k=1}^{n}a_k\,\Phi\Big(\frac kn\Big)=n^{-\beta},\qquad n\ge1,
\end{equation}
one equation for each real $\beta$. The array is triangular and its diagonal entry is
$\Phi(1)=1$, so \eqref{eq:defining_relation_ch1} determines $(a_k)$ rank by rank, with no
hypothesis of any kind and no convergence to establish. What is left to study is the rate at
which the partial sums $A(n)=\sum_{k\le n}a_k$ fall away as $\beta$ varies. Ingham's theorem
answers with a limit, the value one at $\beta=0$ and zero at every $\beta>0$ once the side
condition is checked, as Chapter~\ref{chap:boundary} shows, and it says nothing about the rate.
The theory that follows asks for the exponent. The weight itself is drawn in Figure~\ref{fig:ingham}.

\rafgalleryfig{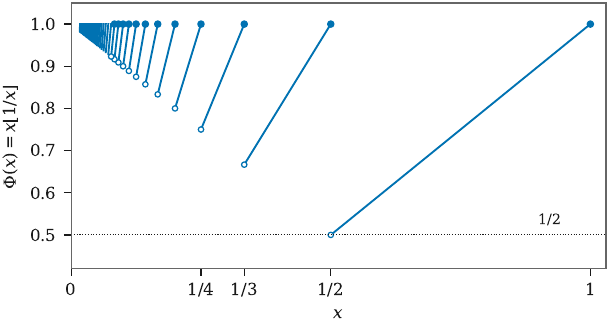}{The Ingham function $\Phi(x)=x\lfloor 1/x\rfloor$ on $(0,1]$. On each interval $\left(\tfrac1{m+1},\tfrac1m\right]$ the value of $\lfloor 1/x\rfloor$ is $m$, so $\Phi(x)=mx$, and the graph is a fan of segments, each reaching the height one at $x=1/m$ and dropping to $m/(m+1)$ just to the right of it. The dotted line marks the height one half, the value the volume will show to be the index of this weight if and only if the Riemann hypothesis holds.}{fig:ingham}

The next chapter makes this precise. It marks out the class of weights for which the construction
is stable, the functions of good variation, brings in the transform and the index that measure the
decay, and states the recurrence the weight generates. None of that machinery is needed here. What
this chapter leaves in place is the object itself, the Ingham weight, and the decision to read it as
a kernel.

\chapter{Functions of good variation}
\label{chap:fgv}

Ingham's theorem\index[terms]{Ingham theorem} is a statement about one target and one conclusion. Its weight, however, is a
function, and a function can be varied. This chapter takes the weight as the object and the
target as a parameter, and asks what survives.

Three exact identities open the chapter and settle whether there is anything to see. Against
the same Ingham weight, three classical sequences give averages that fall away at three
different rates, a constant for the sequence that yields the prime number theorem\index[terms]{prime number theorem}, $n^{-1/2}$
for Liouville\index[names]{Liouville, J.} and $n^{-1}$ for M\"obius\index[names]{M\"obius, A. F.}. The rate belongs neither to the weight alone nor to
the sequence alone. It is what the equation returns when a rate is imposed on it, and reading
it in that direction is what the rest of the volume does.

The apparatus follows in four steps. The defining recurrence\index[terms]{defining equation} attaches a sequence to a weight
and an exponent, with no hypothesis at all, its array being triangular with a nonzero diagonal.
The Mellin transform\index[terms]{Mellin transform} is attached to the weight and computed from it. The regularity
index\index[terms]{regularity index} is the threshold at which the partial sums stop following the imposed rate and settle on
a rate of their own, and the weights for which that threshold exists are the functions of good
variation\index[terms]{function of good variation}. The chapter closes with an existence theorem for the class and with the question of
how far the definition reaches when the weight is unbounded at the origin.

\section{From Ingham's Tauberian theorem to functions of good variation}\label{sec:from_ingham_to_fgv}

Ingham's theorem was applied at one target, a constant, and gave the prime number theorem.
The same average, formed from other sequences, behaves differently, and two exact identities
show how. Reading the three cases together raises the question the rest of the chapter takes
up, and carries the subject from a single theorem to a family, and then to a class of
weights.

\subsection{The first insight}\label{sec:two_identities}

In \S\ref{sec:ingham_pnt} the average of one sequence tended to a constant and gave the
prime number theorem. Two more identities do the same for two other sequences.

The Liouville function\index[terms]{Liouville function} $\lambda(n)=(-1)^{\Omega(n)}$, where $\Omega(n)$ counts the prime
factors of $n$ with multiplicity, satisfies
\[
\sum_{k=1}^{n}\lambda(k)\Big\lfloor\frac nk\Big\rfloor=\lfloor\sqrt n\rfloor,\qquad n\ge1.
\]
The left side is a sum of the kind that appears in an Ingham average. To read it as one,
divide by $n$ and gather the terms against a single weight. Recall the Ingham function\index[terms]{Ingham function} \eqref{eq:ingham_function_ch1} of
\S\ref{sec:method_to_kernel},
\[
\Phi(x)=x\Big\lfloor\frac1x\Big\rfloor,\qquad 0<x\le1,
\]
so that $\Phi(k/n)=(k/n)\lfloor n/k\rfloor$, and for a sequence $(a_k)$ set
\[
A_\Phi(n):=\sum_{k=1}^{n}a_k\,\Phi\!\Big(\frac kn\Big)
=\frac1n\sum_{k=1}^{n}k\,a_k\Big\lfloor\frac nk\Big\rfloor.
\]
This is the Ingham average of $(a_k)$ against the Ingham function\index[terms]{Ingham function} of
Figure~\ref{fig:ingham}, and its convergence
$A_\Phi(n)\to A$ is the summability of Theorem~\ref{thm:ingham}. Taking $a_k=\lambda(k)/k$,
the identity gives
\[
A_\Phi(n)=\frac{\lfloor\sqrt n\rfloor}{n}=n^{-1/2}+O(n^{-1}),
\]
so the average no longer tends to a constant. It goes to zero like $n^{-1/2}$.

A second identity, due to Meissel\index[names]{Meissel, E.}, points the same way,
\[
\sum_{k=1}^{n}\mu(k)\Big\lfloor\frac nk\Big\rfloor=1,\qquad n\ge1,
\]
so that with $a_k=\mu(k)/k$ the average is $A_\Phi(n)=1/n$, going to zero like $n^{-1}$.
Both identities are classical. More generally, for every arithmetic function $f$,
\[
 \sum_{k=1}^{n}f(k)\Big\lfloor\frac nk\Big\rfloor
 =\sum_{m=1}^{n}\sum_{d\mid m}f(d),
\]
because each $k$ on the left is counted once for every multiple $m\le n$. The two formulas
follow on using $\sum_{d\mid m}\lambda(d)=\mathbf 1[m\text{ is a perfect square}]$ and
$\sum_{d\mid m}\mu(d)=\mathbf 1[m=1]$, respectively.

Three sequences now stand side by side. The prime number sequence holds the average at a
constant, the Liouville sequence\index[terms]{Liouville function} sends it to zero like $n^{-1/2}$, and the M\"obius sequence
like $n^{-1}$. Writing the common shape as $n^{-\beta}$ names the exponent $\beta$, equal to
$0$, $\tfrac12$ and $1$ in the three cases. These are not three isolated facts. At $\beta=0$
the plain sum $A(n)=\sum_{k\le n}a_k$ recovers the prime number theorem, and at
$\beta=\tfrac12$, where the Liouville average sits, the plain sum $\sum_{k\le n}\lambda(k)/k$
stays below $n^{-1/2+\eps}$ for every $\eps>0$ exactly when the Riemann
hypothesis holds. Read as one, the three suggest a single law indexed by $\beta$ that would
carry the plain sum from the prime number theorem at $\beta=0$ to the Riemann hypothesis\index[terms]{Riemann hypothesis} at
$\beta=\tfrac12$. Whether such a law exists, and what fixes the passage, is the question the
rest of the chapter takes up.

\begin{figure}[t]
\centering
\rafwithfig{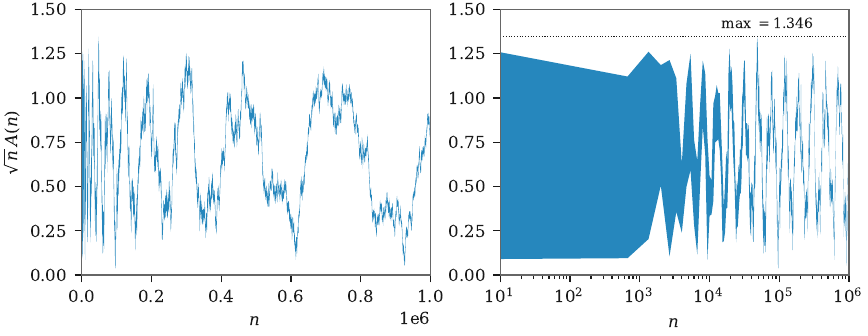}{\rafshowfig}%
  {\fbox{\parbox[c][2.70in][c]{0.72\textwidth}{\centering\small
     Figure~2.1 is missing. Place \texttt{liouville\_normalized.pdf} next to the source file or in a
     subfolder named \texttt{figures}, \texttt{fig} or \texttt{img}. The box has the height of the figure it
     replaces, so the pagination of the volume is unaffected.}}}
\caption{The normalized sum $\sqrt n\,A(n)$, with $A(n)=\sum_{k\le n}\lambda(k)/k$, for
$n\le10^{6}$, the shaded region being the range of the sum over each abscissa. Its rate of
decay is the Riemann hypothesis\index[terms]{Riemann hypothesis}, since $A(n)\ll n^{-1/2+\eps}$ holds if and only if
the hypothesis is true. Over six decades the quantity stays between $0.037$ and $1.346$, the
largest value being reached before $10^{5}$ and not exceeded afterwards. The apparent
boundedness is deceptive. Ingham\index[names]{Ingham, A. E.} showed the matching sum to be unbounded above and below
whenever the imaginary parts of the zeros carry no rational linear relation or only finitely
many, so a bound $\mathcal{O}(n^{-1/2})$ would force infinitely many such relations. Nothing
rules that out and nobody expects it, which is why the volume carries the $\eps$ from the
start.}
\label{fig:liouville}
\end{figure}

A first look at that plain sum calls for a warning. Figure~\ref{fig:liouville} shows the
normalized sum $\sqrt n\,A(n)$ plotted against $n$. Over this interval it stays inside a
narrow band and looks bounded, and the guess $\sqrt n\,A(n)=\mathcal{O}(1)$ is a natural
one. The guess is almost certainly wrong. Studying the matching normalized sum of
the M\"obius function\index[terms]{M\"obius function}, \nm{Ingham}{A. E.} showed that under the hypothesis, with the zeros
simple, and provided the imaginary parts of the nontrivial zeros carry no rational linear
relation or only finitely many, the normalized sum is unbounded, its limits superior and
inferior both infinite~\cite{Ingham1942}. Each of those conditions is expected to hold and none
of them is proved, so the argument gives a reason to doubt the clean bound and not a disproof of
it. The same mechanism drives the Liouville sum. This is the reason the
equivalent of the Riemann hypothesis is the bound $A(n)\ll n^{-1/2+\eps}$, with the
exponent relaxed by an arbitrarily small $\eps$, and not the cleaner
$\mathcal{O}(n^{-1/2})$. The two conjectures of that paper's title are the two these sums
invite, Mertens\index[names]{Mertens, F.} on the M\"obius side and
P\'olya\index[names]{P\'olya, G.} on the Liouville side, each asserting a clean bound that
would have given the hypothesis. Both had been checked far, and the argument above is what made
them doubtful before either was refuted. The natural guess made here is of the same kind. The
theory that follows carries that $\eps$ with care from the very start.

\subsection{A general equation}\label{sec:general_equation}

Read the other way, the three cases are solutions of one equation. For a real parameter
$\beta$, ask for a sequence whose Ingham average is $n^{-\beta}$ up to a smaller error,
\begin{equation}\label{eq:forcing}
\sum_{k=1}^{n}a_k\,\Phi\!\Big(\frac kn\Big)=n^{-\beta}+O(n^{-\delta}),\qquad \delta>\beta.
\end{equation}
At each of the three exponents a solution is already in hand, and in each the plain partial
sum $A(n)=\sum_{k\le n}a_k$ meets a classical statement of prime number theory. At $\beta=0$
the right side is a constant, the endpoint of \S\ref{sec:ingham_pnt}. The prime number
sequence scaled to that constant, $a_k=(1-\Lambda(k))/(2\gamma k)$, has average tending to
$1$ and plain sum recovering the prime number theorem. This endpoint carries no rate, only
the one-sided Tauberian bound of Ingham's\index[names]{Ingham, A. E.} theorem, satisfied here by $-a_k$, so $a_1$ is fixed by the sequence and not
chosen, and \eqref{eq:forcing} is read here as a Tauberian constraint, not the defining
recurrence of the next section. At $\beta=\tfrac12$ the Liouville sequence $a_k=\lambda(k)/k$
solves \eqref{eq:forcing}, its plain sum staying below $n^{-1/2+\eps}$ exactly under
the Riemann hypothesis, and at $\beta=1$ the M\"obius sequence $a_k=\mu(k)/k$ does the same.
The equivalence at $\beta=1$ is the criterion of \nm{Littlewood}{J. E.}~\cite{Littlewood1912},
recorded as Lemma~\ref{lem:littlewood} below, and both cases are set out through the floor
function\index[terms]{floor function} in \cite{CloitreFloor}.

Reading $\beta$ as free to move and asking for the behavior of $A(n)$ at each value, a
threshold appears. Below it the plain sum follows the right side, $A(n)\sim n^{-\beta}$, and
above it the plain sum settles at a floor set by the weight alone. This threshold is a
single number attached to the Ingham function\index[terms]{Ingham function}, and the sections that follow define it,
compute it, and show that for the Ingham function\index[terms]{Ingham function} it is the point $\beta=\tfrac12$ where the
plain sum meets the Riemann hypothesis. A second reading keeps the exponent at $\beta=1$, where
the solution is $\mu(n)/n$, and moves instead the hypothesis on the zeros of $\zeta$, from the
line $\Re s=1$ into the critical strip. Chapter~\ref{chap:boundary} begins that reading at the
line itself, and \S\ref{sec:boundary_to_strip} sets the two side by side.

The threshold and the floor are already visible below zero, where the average can be
inverted by hand. For $\beta<0$ the series $\sum_{j\ge1}\mu(j)\,j^{-(1-\beta)}$ converges
absolutely, and inverting the average gives, with no complex analysis,
\[
A(n)\sim\frac{1}{\Phi^{*}(\beta)}\,n^{-\beta},\qquad
\Phi^{*}(z)=\frac{z}{z-1}\,\zeta(1-z).
\]
The transform $\Phi^{*}$ enters here for the first time, as the constant of that inversion,
and its analytic shape is carried by the zeta function. Below the threshold the plain sum is
transparent to the forcing, following $n^{-\beta}$ with the constant $1/\Phi^{*}(\beta)$, and
the natural guess is that this transparency\index[terms]{transparency} holds as $\beta$ climbs toward $\tfrac12$, the
value at which the first zeros of $\Phi^{*}$ line up under the Riemann hypothesis. That the
transparency\index[terms]{transparency} reaches all the way to $\tfrac12$ is the content of the Riemann hypothesis,
made precise in \S\ref{sec:master_equiv}. The next sections define $\Phi^{*}$ and the
threshold, and turn this guess into a theorem.
\section{The defining recurrence}\label{sec:defining_recurrence}

A recurrence-admissible profile\index[terms]{recurrence-admissible profile} is a real-valued function
$g:(0,1]\to\R$, finite at every point, with $g(1)\neq0$. This is exactly the hypothesis needed for the finite
triangular recurrence, and it says nothing about boundedness, integrability, or the behavior of the
solution. An admissible profile\index[terms]{admissible profile} is a recurrence-admissible profile that is bounded and
Riemann integrable\index[terms]{Riemann integrability}. The general membership conjecture below concerns positive profiles
that are either admissible in this bounded sense or slowly varying at the origin. The wider
recurrence-admissible class accommodates, in particular, the two logarithmic profiles treated by exact
recurrences in Appendices~\ref{app:B} and~\ref{app:F}. Fix a recurrence-admissible profile $g$
and a function $r:\mathbb{N}^*\to\R$. The defining equation\index[terms]{defining equation}
\begin{equation}\label{eq:defining_relation}
\sum_{k=1}^n a_k\,g\!\Big(\frac kn\Big)=r(n)\qquad(n\ge1)
\end{equation}
determines the sequence $(a_n)_{n\ge1}$ uniquely. The relation at $n$ isolates the last term,
$a_n\,g(1)+\sum_{k<n}a_k\,g(k/n)=r(n)$, so that, $g(1)$ being nonzero,
\[
a_n=\frac{1}{g(1)}\Big(r(n)-\sum_{k=1}^{n-1}a_k\,g\!\Big(\frac kn\Big)\Big),
\]
which fixes $a_n$ from $a_1,\dots,a_{n-1}$ and from the single value $r(n)$. Starting from
$a_1=r(1)/g(1)$, induction on $n$ gives existence and uniqueness of the whole sequence. No
analytic hypothesis on $r$ enters, only that its value at each integer is prescribed. The
forcing $r(n)=n^{-\beta}$, with $\beta$ fixed and the sequence $(a_n)$ the unknown, is the case
carried through Part~I, and the quantity under study is the ordinary partial sum
$A(n):=\sum_{k=1}^n a_k$.

Two variants met later are not defining equations\index[terms]{defining equation} in this sense, and the distinction matters.
When the right-hand side is a Tauberian input $r(n)=n^{-\beta}+\mathcal{O}(n^{-c})$, the value
$r(n)$ is not prescribed but only constrained, so the recurrence no longer singles out one
sequence. The statements that assume such an input hold for every sequence compatible with it,
which is what makes them Tauberian rather than definitional. A constraint of the form
$\sum_{k\le x}a_k\,g(k/x)=h(x)$ imposed for every real $x$ is a different object again,
overdetermined relative to the integer recurrence, since $g(k/x)$ varies with $x$ between the
integers while the partial sum changes only at them. Only \eqref{eq:defining_relation}, read at
the integers with a prescribed right-hand side, defines a sequence.

It is worth pausing on the equation itself, because it is not new. With $g(1)\neq 0$
the relation $\sum_{k\le n}a_k\,g(k/n)=n^{-\beta}$ is a discrete Volterra\index[terms]{discrete Volterra} equation of the
first kind. Isolating the diagonal term turns it into the triangular recursion
\[
a_n=\frac{1}{g(1)}\Big(n^{-\beta}-\sum_{k=1}^{n-1}a_k\,g(k/n)\Big),
\]
solved uniquely by forward substitution, and the passage from the second member to the
sequence is a summation against a resolvent\index[terms]{resolvent} kernel, the discrete face of the variation of
constants formula and of its Green's matrix\index[terms]{Green's matrix} \cite{Agarwal2000}. The continuous counterpart, the Volterra integral equation of the first kind and its resolvent\index[terms]{resolvent}, is treated at book length by Brunner\index[names]{Brunner, H.} \cite{Brunner2017}. The classical theory of
such equations is long and complete. It supplies existence and uniqueness, the Green's
matrix with its cocycle law, closed forms through generating functions, and a full
account of boundedness, periodicity, and stability.

And yet one quantity is absent from that theory. Its results sort a
solution by growth, by boundedness, by periodicity, by stability, and none of them
names the threshold exponent on which the transfer law\index[terms]{transfer law} turns. The regularity index\index[terms]{regularity index} is
that exponent, and it is an arithmetic invariant of the kernel. The reason it is missing
is not an oversight but a matter of substance. The coefficients of a classical difference
equation are analytic in nature, whereas the kernel $g(k/n)$ carries arithmetic in the
ratio $k/n$, and it is that arithmetic the index reads. No invariant playing the same
role has been found in the classical literature, and none appears in a reference as
complete as \nm{Agarwal}{R. P.}~\cite{Agarwal2000}. The
classical theory offers the frame, and the arithmetic of the kernel supplies what fills
it.

Reading the asymptotics of a sum $\sum_n a_n\,g(n/x)$ off a Mellin transform is a developed
technique. \nm{Flajolet}{P.}, \nm{Gourdon}{X.} and \nm{Dumas}{P.} organize it for harmonic
sums\index[terms]{harmonic sum} \cite{FlajoletGourdonDumas1995}, factoring the transform of the sum into the
transform of the profile times a Dirichlet series\index[terms]{Dirichlet series} and then moving the contour to the left,
each pole crossed contributing one term of the expansion. The device is available whenever
the factored transform can be controlled on the line one wants to reach, and where it applies
it is the shortest route to an asymptotic expansion. The factorization itself is an identity on
the half plane of convergence and not an asymptotic statement. A condition on the edge of that
half plane, together with positivity or a bound on the terms, turns it into a limit by a
Tauberian argument, as Chapter~\ref{chap:boundary} shows, and a rate asks for estimates beyond
the edge.

The theory below meets a point where it is not available. The factored transform is under
control on the half plane where its integral converges, and the contour has to be moved past
the first zero of $g^{*}$, which lies outside that half plane. Nothing in the factorization
supplies the bound needed there. The route taken in \S\ref{sec:absorption} avoids the point
rather than forcing it. The recurrence is rewritten as a Volterra equation\index[terms]{Volterra equation}\index[names]{Volterra, V.}, and the contour is
moved on its resolvent\index[terms]{resolvent}, which is a continuous function whose decay is proved
directly from the zero free region. The step function never carries a contour. Both routes
read the same zeros, and only the second is available at the exponent this theory needs.

\section{The Mellin transform}\label{sec:mellin_transform}

The spectral object attached to the recurrence is a limit of arithmetic averages, and for a
kernel of the form $g(k/n)$ that limit is a Riemann sum. It coincides with a classical object.

\begin{proposition}\label{prop:mellin_coincidence}
Let $g:(0,1]\to\R$ be bounded and Riemann integrable. For $\Re z<0$ the limit
\[
g^{*}(z)=\lim_{n\to\infty}\frac{-z}{n}\sum_{k=1}^{n}g\!\Big(\frac kn\Big)\Big(\frac kn\Big)^{-z-1}
\]
exists and equals
\[
g^{*}(z)=-z\int_{0}^{1}g(t)\,t^{-z-1}\,dt,
\]
the convergence being uniform on every compact subset of $\{\Re z<0\}$. On that half plane
$g^{*}$ is holomorphic. Its continuation to $\Re z\ge0$ depends on the analytic nature of $g$,
and where $g^{*}$ has a zero the analytic index\index[terms]{analytic index} is
$\eta(g)=\inf\{\Re\rho:g^{*}(\rho)=0\}$, a quantity that is undefined when there is none.
\end{proposition}

\begin{proof}
Fix a compact set $K\subset\{\Re z<0\}$, and put $\sigma_{0}=\max_{z\in K}\Re z<0$ and
$M=\sup_{(0,1]}|g|$. For $z\in K$ write $\sigma=\Re z\le\sigma_{0}$.

Because $\sigma<0$, the bound $|g(t)t^{-z-1}|\le M\,t^{-\sigma-1}$ together with
$\int_{0}^{1}t^{-\sigma-1}\,dt=1/(-\sigma)$ shows that the integral
$I(z)=\int_{0}^{1}g(t)t^{-z-1}\,dt$ converges absolutely, with $|I(z)|\le M/(-\sigma_{0})$. Let
$S_{n}(z)=\frac1n\sum_{k=1}^{n}g(k/n)(k/n)^{-z-1}$ be the right-endpoint Riemann sum, so that
the average in the statement is $-z\,S_{n}(z)$. The claim reduces to $S_{n}\to I$ uniformly on
$K$.

Fix $\delta\in(0,1)$ and split both $S_{n}$ and $I$ at $\delta$. For the integral near $0$,
\[
\Big|\int_{0}^{\delta}g(t)t^{-z-1}\,dt\Big|\le M\int_{0}^{\delta}t^{-\sigma-1}\,dt
=\frac{M\,\delta^{-\sigma}}{-\sigma}\le\frac{M\,\delta^{-\sigma_{0}}}{-\sigma_{0}}.
\]
For the sum near $0$, the terms with $k\le\delta n$ contribute
\[
\Big|\frac1n\sum_{k\le\delta n}g(k/n)(k/n)^{-z-1}\Big|\le M\,n^{\sigma}\sum_{k\le\delta n}k^{-\sigma-1}.
\]
The exponent $-\sigma-1$ exceeds $-1$, so $\sum_{k\le m}k^{-\sigma-1}\le C\,m^{-\sigma}$ with a
constant $C=C(\sigma_{0})$, and with $m=\lfloor\delta n\rfloor$ the right-hand side is at most
$MC\,\delta^{-\sigma}\le MC\,\delta^{-\sigma_{0}}$. Both tails are of order $\delta^{-\sigma_{0}}$
uniformly in $z\in K$, hence tend to $0$ as $\delta\to0$.

On $[\delta,1]$ the factor $t^{-z-1}$ is continuous and bounded by $\delta^{-\sigma_{0}-1}$, and
the family $\{\,t\mapsto t^{-z-1}:z\in K\,\}$ is equicontinuous there. For a bounded
Riemann-integrable $g$ multiplied by such a factor, the right-endpoint Riemann sum over
$(\delta,1]$ converges to $\int_{\delta}^{1}g(t)t^{-z-1}\,dt$, uniformly in $z\in K$.

Given $\eps>0$, choose $\delta$ so that both tails stay below $\eps/3$ for every $z\in K$, then
$N$ so that the difference on $[\delta,1]$ stays below $\eps/3$ for $n\ge N$ and every $z\in K$.
Then $|S_{n}(z)-I(z)|<\eps$ for $n\ge N$ and $z\in K$, so the convergence is uniform on $K$.
Each $S_{n}$ is holomorphic on $\{\Re z<0\}$, and a locally uniform limit of holomorphic
functions is holomorphic, so $g^{*}=-z\,I$ is holomorphic there.
\end{proof}

The rate of convergence is set by the boundary term $g(1/n)\,n^{z}$, of size $n^{\Re z}$, which
decays slowly as $\Re z\to0^{-}$. The convergence is uniform on compact subsets of the open
half plane but not up to the imaginary axis, and it is faster for a $g$ that vanishes at $0$
than for one that does not. The passage across $\Re z=0$, where the integral no longer
converges, is the region that carries the zeros of $g^{*}$, and it rests on the analytic form
of the particular kernel. For the Ingham function the continuation is elementary and closes in
the zeta function.

The transform carries a number of its own, and it is one of the two the volume is about. Where
$g^{*}$ vanishes, the leftmost real part of its zeros is the analytic index $\eta(g)$ named
above, read off the transform alone. The index of \S\ref{sec:reg_index_def} is read off the
equation instead. The two are distinct quantities, and their agreement, wherever it holds, is a
theorem about the particular kernel and never a consequence of the definitions.
Chapter~\ref{chap:diophantine} exhibits a kernel whose transparency threshold is proved to lie
strictly below the leftmost zero of its transform, and sets it beside the rational kernel of
\S\ref{sec:rational_kernel}, whose transform has no zero at all while its index exists.

Before the computation, the kernel has to be placed in the class for which the previous
proposition is stated.

\begin{lemma}\label{lem:ingham_class}
The Ingham function $\Phi(x)=x\lfloor1/x\rfloor$ satisfies $\tfrac12<\Phi(x)\le1$ on $(0,1]$
and $\Phi(1)=1$. Its discontinuities are contained in the countable set
$\{1/m:m\ge2\}$, and it is Riemann integrable on $(0,1]$.
\end{lemma}

\begin{proof}
On $(1/(m+1),1/m]$ one has $\lfloor1/x\rfloor=m$, so $\Phi(x)=mx$, an affine function
increasing from the excluded value $m/(m+1)$ to the value $1$ attained at $x=1/m$. These
intervals partition $(0,1]$, so $\Phi$ is bounded between $\inf_{m\ge1}m/(m+1)=\tfrac12$,
which is not attained, and $1$. The case $m=1$ gives $\Phi(1)=1$, so the diagonal of the
kernel does not vanish and \eqref{eq:defining_relation} determines the sequence. Being affine
on each interval of the partition, $\Phi$ can be discontinuous only at the endpoints $1/m$
with $m\ge2$, a countable set and therefore of Lebesgue\index[names]{Lebesgue, H.} measure zero. A bounded function
whose set of discontinuities is null is Riemann integrable, by Lebesgue's\index[names]{Lebesgue, H.} criterion.
\end{proof}

The transform of the Ingham function is computed once here, and it is the place where the zeta
function enters the theory.

\begin{proposition}\label{prop:ingham_mellin}
For $\Re z<0$ the transform of the Ingham function $\Phi(x)=x\lfloor 1/x\rfloor$ is
\[
\Phi^{*}(z)=\frac{z}{z-1}\,\zeta(1-z).
\]
The right-hand side is meromorphic on $\C$, holomorphic apart from a simple pole at $z=1$ of
residue $-\tfrac12$, and it continues $\Phi^{*}$ across the line $\Re z=0$.
\end{proposition}

\begin{proof}
Lemma~\ref{lem:ingham_class} puts $\Phi$ in the class of
Proposition~\ref{prop:mellin_coincidence}, so the arithmetic averages converge and, for
$\Re z<0$,
\[
\Phi^{*}(z)=-z\int_{0}^{1}\Phi(t)\,t^{-z-1}\,dt
           =-z\int_{0}^{1}\Big\lfloor\frac1t\Big\rfloor t^{-z}\,dt,
\]
the second equality because $\Phi(t)=t\lfloor 1/t\rfloor$. The substitution $t=1/x$ carries the
interval $(0,1]$ to $[1,\infty)$ and gives
\[
\int_{0}^{1}\Big\lfloor\frac1t\Big\rfloor t^{-z}\,dt
=\int_{1}^{\infty}\lfloor x\rfloor\,x^{z-2}\,dx
=\sum_{m\ge1}m\int_{m}^{m+1}x^{z-2}\,dx
=\frac{1}{z-1}\sum_{m\ge1}m\big((m+1)^{z-1}-m^{z-1}\big),
\]
the middle equality because $\lfloor x\rfloor=m$ on $[m,m+1)$. With $c_{m}=m^{z-1}$, summation
by parts gives
\[
\sum_{m=1}^{N}m\,(c_{m+1}-c_{m})=N\,c_{N+1}-\sum_{m=1}^{N}c_{m}
=N(N+1)^{z-1}-\sum_{m=1}^{N}m^{z-1}.
\]
For $\Re z<0$ the boundary term obeys $|N(N+1)^{z-1}|\le N^{\Re z}\to0$, while
$\sum_{m\le N}m^{z-1}\to\zeta(1-z)$, the series converging absolutely because $\Re(z-1)<-1$.
Hence $\sum_{m\ge1}m\big((m+1)^{z-1}-m^{z-1}\big)=-\zeta(1-z)$, and the three displays combine to
\[
\Phi^{*}(z)=\frac{-z}{z-1}\big(-\zeta(1-z)\big)=\frac{z}{z-1}\,\zeta(1-z).
\]
The factor $\zeta(1-z)$ is holomorphic on $\C$ apart from a simple pole at $z=0$, where the
factor $z$ vanishes, so $z\,\zeta(1-z)$ is entire and the only pole of $\Phi^{*}$ is the simple
pole at $z=1$ from $1/(z-1)$, of residue $\zeta(0)=-\tfrac12$. The right-hand side is meromorphic
on $\C$ and agrees with $\Phi^{*}$ on $\Re z<0$, so it continues $\Phi^{*}$ across $\Re z=0$.
\end{proof}

The same computation is carried out through the floor function\index[terms]{floor function} in \cite[\S2.1]{CloitreFloor}.
The zeros of $\Phi^{*}$ are governed by those of the zeta function, and the regularity index\index[terms]{regularity index} of
the next section is defined without reference to them, the two readings being compared only in
Chapter~\ref{chap:equivalence}.

\section{The regularity index}\label{sec:reg_index_def}

For the forcing $r(n)=n^{-\beta}$ the regularity index\index[terms]{regularity index} compares the decay of the partial sum
$A(n)$ with the exponent $\beta$. Two notions are needed, and the transform enters neither.

\begin{definition}[Transparency]\label{def:transparency_fgv}
The exponent $\beta$ is transparent for the recurrence-admissible profile\index[terms]{recurrence-admissible profile} $g$ when there is a constant
$\Xi_g(\beta)$ with
\[A(n)=\Xi_g(\beta)\,n^{-\beta}+o(n^{-\beta}).
\]
Such a constant is unique when it exists, since two of them would differ by a $\delta$ with
$\delta\,n^{-\beta}=o(n^{-\beta})$, forcing $\delta=0$. It is called the transparent coefficient\index[terms]{transparent coefficient}
at $\beta$.
\end{definition}

\begin{definition}[The transparency threshold and functions of good variation]\label{def:reg_index_fgv}
For a recurrence-admissible profile\index[terms]{recurrence-admissible profile} $g$ put
\[
\tau(g)=\sup\big\{\,a\in\R:\ \text{every }\beta<a\ \text{is transparent\index[terms]{transparency} for }g\,\big\},
\]
the transparency threshold\index[terms]{transparency threshold}, with the convention $\tau(g)=-\infty$ when no exponent is transparent.
It is defined for every recurrence-admissible profile and asserts nothing beyond the transparent range. Every
$\beta<\tau(g)$ is transparent, and for every $\delta>0$ some $\beta\in[\tau(g),\tau(g)+\delta)$ is
not transparent. That last property is the sharpness\index[terms]{sharpness} carried by the supremum, and it is weaker
than opacity of the whole range above the threshold, which no result of this volume asserts.

The profile is a function of good variation\index[terms]{function of good variation} when $\tau(g)$ is finite and the partial sums are
absorbed\index[terms]{absorption} above it,
\[
A(n)=\mathcal{O}\!\left(n^{-\tau(g)+\eps}\right)\qquad(\beta\ge\tau(g),\ \eps>0),
\]
and its regularity index\index[terms]{regularity index} is then $\alpha(g)=\tau(g)$. The notation $\alpha(g)$ is reserved for that
case throughout the volume, so that writing $\alpha(g)$ asserts that $g$ is a function of good
variation, while $\tau(g)$ is used whenever the transparent range is discussed and membership is
not in question.
\end{definition}

At $\beta=\alpha(g)$ the absorption clause applies and may carry a logarithmic correction. The
constant profiles are the degenerate case. The kernel $g\equiv c$ follows the forcing exactly at
every real $\beta$, the defining equation reading $c\,A(n)=n^{-\beta}$, so every exponent is
transparent and $\tau(g)=+\infty$, while its transform is the nonzero constant $c$ and carries no
zero at all. The arithmetic reading of these profiles therefore sits at infinity, and no analytic
index is defined for them. The class is defined by a finite threshold, and the constants are set aside for that reason
rather than for a failure of their own.

The transform defines neither transparency nor the index. What it supplies is the model
coefficient on the initial half plane, and the law whose persistence the theory studies.
Wherever the identification is used in this volume it is proved, first for the Ingham
function in
Chapter~\ref{chap:equivalence} and then for each entry of the two galleries, and it reads
\begin{equation}\label{eq:xi_reciprocal}
\Xi_g(\beta)=\frac{1}{g^{*}(\beta)},
\end{equation}
with the convention $\Xi_g(\beta)=0$ at a pole of $g^{*}$, where transparency\index[terms]{transparency} asserts
$A(n)=o(n^{-\beta})$. Wherever this identification holds, a zero of $g^{*}$ cannot be a
transparent exponent, since the reciprocal is not finite.

Abel summation\index[terms]{Abel summation} against the defining equation\index[terms]{defining equation} turns the identification into the convergence of the
Riemann--Stieltjes\index[terms]{Riemann--Stieltjes integral} sums $\sum_{k<n}A(k)\big(g(\tfrac{k+1}{n})-g(\tfrac kn)\big)$ toward
$\Xi_g(\beta)n^{-\beta}\big(g(1)-g^{*}(\beta)\big)$, an interchange that needs control of the
variation of $g$ and not only its Riemann integrability. A weighted mass condition supplies that
control, and it holds at every transparent exponent of every profile of bounded variation in the
two galleries.

\begin{theorem}\label{thm:xi_reciprocal}
Let $g$ be an admissible profile, left continuous on $(0,1)$ and of bounded variation on
$[\delta,1]$ for every $\delta>0$. Let $\nu$ be its Lebesgue--Stieltjes measure on $(0,1)$,
in the left-continuous convention
$\nu([x,y))=g(y)-g(x)$ for $0<x<y<1$, and
$\jmath=g(1)-g(1^{-})$ its diagonal jump\index[terms]{diagonal jump}. Let $\beta$ satisfy the
weighted mass condition
\begin{equation}\label{eq:xi_mass}
\mathcal V_\beta=\int_{(0,1)}t^{-\beta}\,d|\nu|(t)<\infty .
\end{equation}
If $\beta$ is transparent for $g$, then \eqref{eq:xi_reciprocal} holds there, with
\[
g^{*}(\beta)=g(1^{-})-\int_{(0,1)}t^{-\beta}\,d\nu(t).
\]
In particular $g^{*}(\beta)\neq0$, so no zero of $g^{*}$ meeting \eqref{eq:xi_mass} is a
transparent exponent.
\end{theorem}

\begin{proof}
For $\Re z<0$, Stieltjes integration by parts in the left-continuous convention gives
\[
 \int_{[\delta,1)}t^{-z}\,d\nu(t)
 +\int_{\delta}^{1}g(t)\,d(t^{-z})
 =g(1^{-})-\delta^{-z}g(\delta).
\]
As $\delta\downarrow0$, the second integral tends to
$g^{*}(z)=-z\int_0^1g(t)t^{-z-1}\,dt$, and the boundary term tends to zero because $g$ is
bounded. Hence
\[
 g^{*}(z)=g(1^{-})-\int_{(0,1)}t^{-z}\,d\nu(t)\qquad(\Re z<0).
\]
Under \eqref{eq:xi_mass}, the right-hand side is holomorphic on $\Re z<\beta$ and extends
continuously to $z=\beta$, by dominated convergence. It agrees with the original transform on
their common half plane and therefore supplies its continuation and boundary value there. If
$\beta<0$, this is simply the value of the original convergent integral. The transform is thus
blind to the diagonal jump, as it must be, since changing $g(1)$ alone changes no integral.

Abel summation turns the defining equation into
$n^{-\beta}=A(n)g(1)-\sum_{k<n}A(k)\Delta_{n,k}$ with
$\Delta_{n,k}=g(\tfrac{k+1}n)-g(\tfrac kn)$. Left continuity gives
$\Delta_{n,k}=\nu\bigl([\tfrac kn,\tfrac{k+1}n)\bigr)$ for $k\le n-2$, while the last
increment carries in addition the jump $\jmath$. Write $\Xi=\Xi_g(\beta)$ and insert
$A(k)=\Xi k^{-\beta}+\eps_k$
with $\eps_k=o(k^{-\beta})$. Let $\varphi_n$ be the step function equal to
$(k/n)^{-\beta}$ on $[\tfrac kn,\tfrac{k+1}n)$ and to zero on $(0,\tfrac1n)$, so that
$\varphi_n\to t^{-\beta}$ pointwise on $(0,1)$ and $|\varphi_n|\le2^{|\beta|}t^{-\beta}$.
Dominated convergence under \eqref{eq:xi_mass} gives
\[
\sum_{k<n}\Bigl(\frac kn\Bigr)^{-\beta}\Delta_{n,k}
=\int_{(0,1)}\varphi_n\,d\nu+\Bigl(\frac{n-1}n\Bigr)^{-\beta}\jmath
\longrightarrow\int_{(0,1)}t^{-\beta}\,d\nu+\jmath .
\]
For the error part, fix $\delta>0$ and split at the rank $\delta n$. Above it the factor
$\sup_{k>\delta n}|\eps_k|k^{\beta}$ tends to zero and multiplies a quantity at most
$2^{|\beta|}n^{-\beta}(\mathcal V_\beta+|\jmath|)$. Below it $|\eps_k|\le Ck^{-\beta}$
and the increments involved lie in $(0,\delta+\tfrac1n)$, so the contribution is at most
$C2^{|\beta|}n^{-\beta}\int_{(0,\delta+1/n)}t^{-\beta}d|\nu|$, which is small with $\delta$
by \eqref{eq:xi_mass}. Dividing the Abel identity by $n^{-\beta}$, letting $n$ grow and then
$\delta$ decrease leaves
$1=\Xi\bigl[g(1)-\jmath-\int_{(0,1)}t^{-\beta}d\nu\bigr]=\Xi\,g^{*}(\beta)$.
\end{proof}

\begin{remark}\label{rem:xi_reciprocal_scope}
The diagonal jump cancels between the diagonal weight $A(n)g(1)$ and the last increment, which
is why the identification survives for profiles carrying one, as in Appendices~\ref{app:O}
and~\ref{app:P}. The threshold guarantees transparency only for $\beta<\tau(g)$ and is silent
at $\beta=\tau(g)$. Thus the theorem gives \eqref{eq:xi_reciprocal} throughout that guaranteed
range wherever \eqref{eq:xi_mass} holds, and at the endpoint only if transparency and the mass
condition are separately established. The binomial harmonic\index[terms]{binomial harmonic kernel}
kernel of Appendix~\ref{app:Q} illustrates why no endpoint value follows from the threshold
alone. It is a bivariate kernel induced by no profile and hence lies outside the scope of this
theorem. It satisfies the analogous identity $\Xi_G=1/G^{*}$ at every $\beta<1$, while at the endpoint
$\beta=1$ its partial sums still follow an exact power, with the different constant $1$.
Under \eqref{eq:xi_mass} the Stieltjes expression for $g^{*}(\beta)$ is finite, so $g^{*}$ has no
pole there. The condition fails for the Ingham function at every positive exponent, its jumps
at the points $1/m$ summing to a divergent series. That is the case
Chapter~\ref{chap:equivalence} settles by the M\"obius route.
\end{remark}

If the regularity hypotheses of Theorem~\ref{thm:xi_reciprocal} are dropped, the identification
can fail, and a two-valued profile is enough to show it.

\begin{proposition}\label{prop:xi_mass_needed}
For $\lambda\notin\{0,1\}$ define
\[
g_\lambda(x)=
\begin{cases}
1,&x=1/m\ \text{for an integer }m\ge1,\\
\lambda,&\text{otherwise.}
\end{cases}
\]
Then $g_\lambda$ is an admissible profile, its transform is the constant $g_\lambda^{*}\equiv\lambda$,
which has no zero, and the exponent $\beta=0$ is transparent with
\[
\Xi_{g_\lambda}(0)=1\neq\frac1\lambda=\frac{1}{g_\lambda^{*}(0)} .
\]
\end{proposition}

\begin{proof}
The profile is bounded, its discontinuities form the countable set $\{1/m\}$, so it is Riemann
integrable, and $g_\lambda(1)=1\neq0$. It equals $\lambda$ outside a set of measure zero, so
Proposition~\ref{prop:mellin_coincidence} gives, for $\Re z<0$,
\[
g_\lambda^{*}(z)=-z\int_{0}^{1}\lambda\,t^{-z-1}\,dt=\lambda .
\]
On the grid, $k/n$ is the reciprocal of an integer exactly when $k$ divides $n$, so
\[
g_\lambda\Big(\frac kn\Big)=\lambda+(1-\lambda)\mathbf 1_{k\mid n},
\]
and the defining equation \eqref{eq:defining_relation} reads
\[
\lambda A(n)+(1-\lambda)\sum_{k\mid n}a_k=n^{-\beta}.
\]
At $\beta=0$ the sequence $a_1=1$ and $a_n=0$ for $n\ge2$ satisfies it at every rank, both
$A(n)$ and $\sum_{k\mid n}a_k$ being equal to $1$. The diagonal coefficient $g_\lambda(1)=1$ does
not vanish, so the triangular recurrence has this solution and no other, and $A(n)=1$ for every
$n$. The exponent $0$ is therefore transparent with transparent coefficient $1$, while
$1/g_\lambda^{*}(0)=1/\lambda$.
\end{proof}

At every interior point $1/m$ both one-sided limits of $g_\lambda$ equal $\lambda$ while its value
is $1$. These are isolated values away from the limit and not jumps between two distinct
one-sided limits, so $g_\lambda$ is not left continuous, and the increment convention
$\nu([x,y))=g_\lambda(y)-g_\lambda(x)$ of Theorem~\ref{thm:xi_reciprocal} defines no measure for
it. The sets $[x,1/m)$ shrink to the empty set as $x\uparrow1/m$ while their assigned value stays
$1-\lambda$. The hypotheses of that theorem therefore fail, and \eqref{eq:xi_mass} has no meaning
for this profile. What the example shows is that a transform with a continuation and a transparent
exponent do not by themselves force the reciprocal identity in the full admissible class. It does
not show that \eqref{eq:xi_mass} is necessary inside the left continuous class of the theorem,
and the question is what has to be added.

\begin{openproblem}\label{op:xi_reciprocal}
Which conditions on $g$, weaker than the weighted mass condition \eqref{eq:xi_mass}, still force
\eqref{eq:xi_reciprocal} at a transparent exponent? Theorem~\ref{thm:xi_reciprocal} settles the
question under that condition, Chapter~\ref{chap:equivalence} settles it for the Ingham profile
despite its failure, the exact recurrences of Appendices~\ref{app:B} and~\ref{app:F} settle it
for the two unbounded logarithmic profiles, and Theorem~\ref{thm:fgv_power_singularity} settles
it for the power-singular family on its whole transparent interval. Every profile in that list is
left continuous, which $g_\lambda$ is not, and left continuity at the points $1/m$ is one
candidate for the missing hypothesis.
\end{openproblem}

The existence of $\alpha(g)$ is a property of the discrete
equation, whereas the equality $\alpha(g)=\eta(g)$ with the analytic index\index[terms]{analytic index} is a theorem proved
where it holds. It is proved for the affine and logarithmic profiles of this part, and for the
Ingham function it is a theorem under the Riemann hypothesis\index[terms]{Riemann hypothesis} and not otherwise, so that
$\alpha(g)=\eta(g)$ is nowhere a convention.

The threshold at which the transfer law changes has a meaning older than the problem, borrowed
from the equidimensional equation of Euler\index[names]{Euler, L.} and Cauchy\index[names]{Cauchy, A.-L.}, in which the exponents of the
solutions are the roots of an indicial polynomial\index[terms]{indicial polynomial} and the transform takes the part of that
polynomial. The analogy runs through three degrees, a differential equation of the first order,
one of the second, and an integral equation of Volterra type that covers the rest, and it has
limits that decide what the continuous picture can and cannot supply.
Section~\ref{sec:euler_cauchy} sets out the three together.

The arithmetic side is another matter, and here caution is due. The equality
$\alpha(g)=\eta(g)$ is proven for the affine kernel\index[terms]{affine kernel}. It is not a consequence of
continuity. A kernel continuous on $[0,1]$ may behave otherwise than one merely
continuous on $\left]0,1\right]$, where the conduct near the origin comes into play. The
conditions under which the two indices coincide are not settled, bounded variation\index[terms]{bounded variation} being the
natural candidate for the control that is missing, and this book states the equality only where a
proof stands behind it.

\section{The space of FGV and an existence theorem}

That such functions exist beyond trivial cases is not obvious from the definitions,
and the general membership question is open.

\begin{theorem}\label{thm:fgv_existence}
Every affine function $g(x)=c_1x+c_0$ with $c_0,c_1>0$ is a function of good
variation, of index $\alpha(g)=c_0/(c_0+c_1)$.
\end{theorem}

\begin{proof}
This is Theorem~\ref{thm:affine_body} below. The defining equation\index[terms]{defining equation} is solved exactly
through a first order recurrence, the full computation is carried out in
Appendix~\ref{app:A}, and the result was first stated in
\cite[Theorem~1.1]{Cloitre2016}.
\end{proof}

Membership is not automatic, and the next profile shows how it can fail. It is the first
elementary example in this volume of an admissible profile that is not a function of good
variation.

\begin{example}[The diagonal profile]\label{ex:diagonal_profile}
Let $g_\Delta(x)=0$ for $0<x<1$ and $g_\Delta(1)=1$. The profile is bounded, Riemann integrable,
nonconstant, and $g_\Delta(1)\neq0$, so it is admissible. Every term of the defining
equation\index[terms]{defining equation} with $k<n$ carries the factor $g_\Delta(k/n)=0$, leaving
\[
a_n\,g_\Delta(1)=n^{-\beta},\qquad\text{so}\qquad a_n=n^{-\beta},\qquad
A(n)=\sum_{k\le n}k^{-\beta}.
\]
No exponent is transparent. For $\beta<1$ the partial sums satisfy
$A(n)\sim n^{1-\beta}/(1-\beta)$, so $n^{\beta}A(n)\sim n/(1-\beta)$ and the ratio is unbounded.
For $\beta=1$ they satisfy $A(n)\sim\log n$, and $n^{\beta}A(n)\sim n\log n$. For $\beta>1$ they
converge to $\zeta(\beta)$, and $n^{\beta}A(n)\sim\zeta(\beta)n^{\beta}$. In each regime
$n^{\beta}A(n)$ tends to infinity, so no constant $\Xi_{g_\Delta}(\beta)$ can satisfy
Definition~\ref{def:transparency_fgv}. Hence $\tau(g_\Delta)=-\infty$ and $g_\Delta$ is not a
function of good variation\index[terms]{function of good variation}.
\end{example}

The failure is a failure of membership and not an index at $-\infty$, the letter $\alpha$ being
reserved by Definition~\ref{def:reg_index_fgv} for profiles that belong to the class. The example
also shows what the membership conjecture must exclude. The profile vanishes on all of $(0,1)$,
and it is that degeneracy, not any irregularity, that destroys transparency.

\begin{conjecture}\label{conj:fgv_membership}
Let $g:(0,1]\to(0,\infty)$ be a nonconstant recurrence-admissible profile. Suppose that
either
\begin{enumerate}[label=\textup{(\roman*)}]
\item $g$ is bounded and Riemann integrable on $(0,1]$, or
\item $g$ is Riemann integrable on $[\delta,1]$ for every $\delta>0$ and is slowly varying at
the origin, in the sense that $u\mapsto g(1/u)$ is slowly varying at infinity.
\end{enumerate}
Then $g$ is a function of good variation\index[terms]{function of good variation}: the threshold
$\tau(g)$ of Definition~\ref{def:reg_index_fgv} is finite and the partial sums are absorbed
above it.
\end{conjecture}
\begin{proofstatus}{The conjecture is not proved in the present volume, and the general membership
question is open. The evidence recorded here consists of the affine kernels of
Theorem~\ref{thm:fgv_existence}, the bounded strictly positive profiles in the gallery whose
membership has a complete proof, and the two unbounded slowly varying logarithmic profiles of
Appendices~\ref{app:B} and~\ref{app:F}. The power-singular family of
Theorem~\ref{thm:fgv_power_singularity} lies outside the two hypotheses and therefore does not
count as evidence for the conjecture itself. Genuinely bivariate RAFs lie outside the statement
as well. Positivity is a hypothesis of the
conjecture and not a convenience. Without it the statement is false,
Example~\ref{ex:diagonal_profile} exhibiting an admissible profile, nonconstant and vanishing on
$(0,1)$, for which no exponent is transparent. No theorem in this volume assumes the
conjecture.}
\end{proofstatus}

\begin{remark}[A stronger form, and what it would need]
The weaker hypothesis $\int_0^1g(t)\,dt\neq0$ would cover profiles that change sign, and the
question of whether it suffices is open. It is recorded here as a possible strengthening and it is
not part of the conjecture above, nor of the definition of an admissible profile. Nothing in this
volume treats a sign-changing profile whose integral vanishes.
\end{remark}

\begin{remark}
The bounded branch of the statement is formulated in \cite{Cloitre2016} for continuous profiles
and is extended here to the Riemann integrable setting. The slowly varying branch is the natural
working enlargement suggested by Appendices~\ref{app:B} and~\ref{app:F}. In the original
formulation the class was named by its analytic side alone, and the conjecture
would have been empty. Separating admissibility from the arithmetic property, as
\S\ref{sec:reg_index_def} does, is what gives it content. What this volume proves is membership for the affine family and for the
kernels of the gallery whose appendices carry complete arguments, the tables of the two
galleries recording which entries those are. For the Ingham function membership is not
proved here. A zero free half plane $\Re s>\sigma$ with $\sigma<1$ would give the positive
threshold $\tau(\Phi)\ge1-\sigma$ by Corollary~\ref{cor:zerofree}, and would still leave
membership open, but no such half plane is known, the supremum of the
real parts of the zeros of $\zeta$ not being known to be less than $1$. Conversely, membership
with a positive index would already give such a half plane, by the second direction of that
corollary. The exact value $\alpha(\Phi)=\tfrac12$, and not membership by itself, is what
Theorem~\ref{thm:tauberian_rh} makes equivalent to the Riemann
hypothesis\index[terms]{Riemann hypothesis}.
The content of the theory lies in identifying the regularity index precisely, and
the chapters that follow are organized around that question. Membership for the positive
absolutely continuous class, the first result of this kind beyond explicit families, is
proved in Chapter~\ref{chap:raf}.
\end{remark}

\subsection{Unbounded kernels beyond the bounded scale}
\label{sec:fgv_slowly_varying}

Boundedness enters the definition as a convenience and not as a necessity. Two
entries of the gallery are unbounded at the origin, the logarithmic kernel
$1-\lambda\log x$ of Appendix~\ref{app:B} and the kernel $x-\log x$ of
Appendix~\ref{app:F}, and both carry a complete determination of their
indices. Both blow up like $-\log x$, and that scale is the one under which the whole
apparatus of Chapter~\ref{chap:fgv} survives untouched. Read at $u=1/x$, a
logarithmic singularity is a slowly varying\index[terms]{slowly varying} function in the sense of Karamata\index[names]{Karamata, J.} \cite{Seneta1976}, the
notion around which the vocabulary of the theory is organized in
\S\ref{sec:rv}.

\begin{proposition}
\label{prop:fgv_slowly_varying}
Let $g:(0,1]\to\R$ be Riemann integrable on $[\delta,1]$ for every $\delta>0$, with
$g(1)\neq0$, and suppose that $g$ is of subpower growth\index[terms]{subpower growth} at the origin,
\[
\text{for every }\eps>0\text{ there is }C_\eps\text{ with }|g(t)|\le C_\eps\,t^{-\eps}
\text{ for all }t\in(0,1].
\]
Then
for every $\sigma<0$ the integral $\int_{0}^{1}|g(t)|\,t^{-\sigma-1}\,dt$ converges,
the transform
\[
g^{*}(z)=-z\int_{0}^{1}g(t)\,t^{-z-1}\,dt
\]
is holomorphic on $\{\Re z<0\}$, and the arithmetic Riemann sums of
Proposition~\ref{prop:mellin_coincidence} converge to it uniformly on compact subsets
of that half plane. The defining equation \eqref{eq:defining_relation} determines
$(a_n)$ as before, so transparency and the threshold $\tau(g)$ are defined for such a kernel,
and the regularity index of Definition~\ref{def:reg_index_fgv} follows once absorption is
proved.
\end{proposition}

\begin{proof}
Fix $\sigma<0$ and take $\eps<-\sigma$ in the subpower bound. Then
\[
\int_{0}^{\delta}|g(t)|\,t^{-\sigma-1}\,dt
\le C_\eps\int_{0}^{\delta}t^{-\eps-\sigma-1}\,dt
=\frac{C_\eps\,\delta^{-\eps-\sigma}}{-\eps-\sigma},
\]
finite because $-\eps-\sigma>0$, and the remaining range $[\delta,1]$ carries a
bounded Riemann integrable function. Holomorphy follows from local uniform
domination by the same bound. For the sums, the proof of
Proposition~\ref{prop:mellin_coincidence} splits at $\delta$ and needs only that the
part near the origin be uniformly small. Its arithmetic counterpart is
\[
\frac1n\sum_{k\le\delta n}\Big|g\Big(\frac kn\Big)\Big|\Big(\frac kn\Big)^{-\sigma-1}
\le\frac{C_\eps}{n}\sum_{k\le\delta n}\Big(\frac kn\Big)^{-\eps-\sigma-1}
\le\frac{C_\eps\,\delta^{-\eps-\sigma}}{-\eps-\sigma}+o(1),
\]
by comparison of the sum with its integral, the exponent $-\eps-\sigma-1$ being
larger than $-1$. The right side tends to $0$ with $\delta$, uniformly in $n$, and
the rest of the argument is unchanged. The convergence of the sums is uniform on
compact subsets of $\{\Re z<0\}$ and not merely pointwise, the constant $C_\eps$ and the
threshold $\delta$ depending on the compact through $\sigma$ alone.
\end{proof}

Slow variation is the case of the proposition that the theory meets, and it enters through Potter's
bound rather than as a hypothesis of its own.

\begin{corollary}\label{cor:fgv_slowly_varying}
Let $g:(0,1]\to\R$ be Riemann integrable on $[\delta,1]$ for every $\delta>0$, with $g(1)\neq0$,
and suppose that $u\mapsto g(1/u)$ is slowly varying\index[terms]{slowly varying} at infinity. Then
$g$ satisfies the hypothesis of Proposition~\ref{prop:fgv_slowly_varying}, and its conclusions
hold.
\end{corollary}

\begin{proof}
Potter's bound\index[terms]{Potter bound}, \cite[Theorem~1.5.6]{Bingham1989}, gives for
every $\eps>0$ constants $C_\eps$ and $U_\eps$ with $|g(1/u)|\le C_\eps u^{\eps}$ for $u\ge
U_\eps$, that is $|g(t)|\le C_\eps t^{-\eps}$ for $t\le1/U_\eps$. On the remaining interval
$[1/U_\eps,1]$ the function is bounded, so enlarging $C_\eps$ extends the bound to all of $(0,1]$.
\end{proof}

The hypothesis is one of size and carries no sign condition. A profile of subpower growth may
change sign near the origin, and nothing above uses positivity.

The analytic scale is sharp in the following sense. A kernel with a power singularity,
$g(t)\asymp t^{-a}$ with $a>0$, has $\int_{0}^{1}|g(t)|t^{-\sigma-1}\,dt$ convergent
only for $\sigma<-a$, and the term of index $k=1$ in the arithmetic sum equals
$n^{a+z}$ up to a constant, which fails to tend to $0$ as soon as $\Re z\ge-a$. Both
the integral and the arithmetic sums then live on the smaller half plane
$\{\Re z<-a\}$. Thus subpower growth at the origin is what the analytic construction of this chapter needs in
order to reach the whole half plane $\Re z<0$ without change, and slow variation is the classical
scale that supplies it. Neither is the boundary of the recurrence-defined class.

\begin{theorem}[Power singularities beyond slow variation]
\label{thm:fgv_power_singularity}
For every $\lambda>0$, the positive profile
\[
g_\lambda(x)=x^{-\lambda},\qquad 0<x\le1,
\]
is a function of good variation with
\[
\alpha(g_\lambda)=0.
\]
More precisely, if $A_\beta(n)$ denotes the partial sum determined by the defining equation at
the exponent $\beta$, then
\begin{align*}
A_\beta(n)&=\frac{\beta+\lambda}{\beta}\,n^{-\beta}+o(n^{-\beta})
&& (\beta<0),\\
A_0(n)&=-\lambda\log n+\mathcal O_\lambda(1),\\
A_\beta(n)&=\mathcal O_{\beta,\lambda}(1)
&& (\beta>0).
\end{align*}
The transform, initially defined on $\Re z<-\lambda$, has the meromorphic continuation
\[
g_\lambda^{*}(z)=\frac{z}{z+\lambda}.
\]
Consequently $\eta(g_\lambda)=0$ and the two indices coincide throughout the power-singular
family: $\alpha(g_\lambda)=\eta(g_\lambda)$.
\end{theorem}

\begin{proof}
Fix $\beta$ and write $A(n)=\sum_{k\le n}a_k$. Since
$g_\lambda(k/n)=n^\lambda k^{-\lambda}$, the defining equation is equivalent to
\[
\sum_{k\le n}a_k k^{-\lambda}=n^{-(\beta+\lambda)}.
\]
It follows by taking consecutive differences that $a_1=1$ and, for $n\ge2$,
\begin{equation}\label{eq:power-singular-coefficients}
a_n=n^{-\beta}\left[1-\left(1-\frac1n\right)^{-(\beta+\lambda)}\right]
=-(\beta+\lambda)n^{-\beta-1}
 +\mathcal O_{\beta,\lambda}(n^{-\beta-2}).
\end{equation}
If $\beta<0$, the standard power-sum estimate gives
\[
\sum_{k=2}^{n}k^{-\beta-1}=\frac{n^{-\beta}}{-\beta}+o(n^{-\beta}),
\qquad
\sum_{k=2}^{n}k^{-\beta-2}=o(n^{-\beta}),
\]
where the second relation includes the logarithmic case $\beta=-1$. Summing
\eqref{eq:power-singular-coefficients} therefore yields
\[
A(n)=\frac{\beta+\lambda}{\beta}\,n^{-\beta}+o(n^{-\beta}).
\]
At the apparent exceptional value $\beta=-\lambda$ the weighted partial sum above is identically
one, so $a_1=1$, $a_n=0$ for $n\ge2$, and the same formula reads
$A(n)=o(n^\lambda)$ with transparent coefficient zero.

At $\beta=0$, equation \eqref{eq:power-singular-coefficients} gives
$a_n=-\lambda/n+\mathcal O_\lambda(n^{-2})$, hence
$A(n)=-\lambda\log n+\mathcal O_\lambda(1)$, so the exponent zero is not transparent. If
$\beta>0$, the same equation gives $a_n=\mathcal O_{\beta,\lambda}(n^{-\beta-1})$, so the
partial sums are bounded. Every $\beta<0$ is thus transparent, zero is not, and every
$\beta\ge0$ satisfies the absorption bound $A(n)=\mathcal O_{\beta,\lambda,\eps}(n^\eps)$.
Definition~\ref{def:reg_index_fgv} gives $\alpha(g_\lambda)=0$.

Finally, for $\Re z<-\lambda$,
\[
-z\int_0^1t^{-\lambda-z-1}\,dt=\frac{z}{z+\lambda}.
\]
The right side continues the transform meromorphically and has its only zero at $z=0$, proving
$\eta(g_\lambda)=0$. Its reciprocal $(\beta+\lambda)/\beta$ is exactly the transparent
coefficient found above.
\end{proof}

\begin{remark}
The theorem answers affirmatively the former open problem about unbounded FGVs beyond slow
variation. Indeed $g_\lambda(1/u)=u^\lambda$ is regularly varying with positive index and is not
slowly varying, since $g_\lambda(1/(cu))/g_\lambda(1/u)=c^\lambda$ for $c\ne1$. It also separates
the two roles of the origin. The Mellin integral lives only on $\Re z<-\lambda$, while the discrete
recurrence still has a complete threshold theory and its transform reaches the threshold by
meromorphic continuation.

There is also a direct contribution to Open Problem~\ref{op:xi_reciprocal}. The weighted
variation condition of Theorem~\ref{thm:xi_reciprocal} holds here exactly for
$\beta<-\lambda$, whereas the reciprocal formula remains valid throughout
$-\lambda\le\beta<0$, including the pole $\beta=-\lambda$ under the zero-reciprocal convention.
Thus the power family supplies a new infinite-variation test case, without deciding the general
question.
\end{remark}
The chapter has turned a summation method into a class. A weight, an exponent and one nonzero
value at the diagonal produce a sequence with no hypothesis at all. The transform of the weight
and the threshold of the equation give two numbers attached to it, and the weights whose
threshold is finite and carries absorption above it are the functions of good variation. What is
not settled is which weights those are. Conjecture~\ref{conj:fgv_membership} states the
membership question for nonconstant positive profiles, bounded or slowly varying at the origin,
and is not proved here, and Open
Problem~\ref{op:xi_reciprocal} asks what forces the transparent coefficient to be the reciprocal
of the transform. The next chapter leaves the general question aside and asks what a weight
returns at the edge of its transform, a limit, before Chapter~\ref{chap:equivalence} turns to the
one weight whose index is not elementary to compute.

\chapter{Boundary transfer and the prime number theorem}
\label{chap:boundary}

Ingham's theorem returns a limit. The averages of a sequence against the Ingham function
settle at a value, a one-sided condition holds on the terms, and the partial sums settle at the
same value. The next chapter asks for more, a rate of decay and the exponent at which that rate
changes nature, and the statement that this exponent is one half is equivalent to the Riemann
hypothesis\index[terms]{Riemann hypothesis}. This chapter stays at the limit and asks what a weight alone can return there. At the limit the answer is a
characterization, and seeing it first is what makes the rate the right next question.

The return from an average to the partial sums is a convolution problem in the logarithmic
variable, and Wiener\index[names]{Wiener, N.}'s theorem on the translates of an integrable function decides it. The decision
is taken on the imaginary axis, the edge of the half plane where the transform of the weight
converges, and the condition is that $g^{*}$ have no zero there. For the Ingham function that
condition is the non-vanishing of $\zeta$ on the line $\Re s=1$, the classical input of the
prime number theorem\index[terms]{prime number theorem}, and the chapter proves that theorem three times, from the three
sequences that opened Chapter~\ref{chap:fgv}, and then in every arithmetic progression of a fixed
modulus.

Three mechanisms are kept apart throughout. A two-sided bound on the terms,
$na_n=\mathcal O(1)$, is what the strong Hardy--Littlewood--Ramanujan criterion supplies for
the defining equation, and it carries the main theorem. A one-sided bound together with the
positivity of an auxiliary function carries Ingham's theorem itself. Positivity of the terms,
with no bound on them, carries the Wiener--Ikehara theorem through the Mellin transform of the
averaged sum. Each returns a limit and none returns a rate. A line free of zeros is not a strip
free of zeros, and the last section says what changes when the boundary is moved into the
critical strip. None of the results of this chapter bears on the existence of the regularity
index, and none identifies it with a zero of a transform.

\section{The boundary problem and the HLR criterion}\label{sec:boundary_problem}

For a profile $g$ on $(0,1]$ and a sequence $(a_n)$, the averaged sum and the partial sum are
\[
A_g(x)=\sum_{n\le x}a_n\,g\Big(\frac nx\Big),\qquad A(x)=\sum_{n\le x}a_n,
\]
read at the integers unless real arguments are needed. The boundary problem\index[terms]{boundary problem} asks when
\[
A_g(n)\longrightarrow L\qquad\Longrightarrow\qquad A(n)\longrightarrow L .
\]
Two conditions enter and they are of different kinds. One bears on the profile and is spectral.
The other bears on the terms and excludes oscillations too fast for an average to see. The
defining equation enters through the second, and the Hardy-Littlewood-Ramanujan
criterion\index[terms]{Hardy--Littlewood--Ramanujan criterion} is the form it takes.

\begin{definition}\label{def:HLR}
A function $g$ satisfies the HLR criterion if for every $\beta\ge0$,
$\sum_{k=1}^n a_k g(k/n)=n^{-\beta}$ implies $na_n=o(n^{\eps})$ for
every $\eps>0$. It satisfies the strong HLR criterion if under the
same hypothesis $na_n=\mathcal{O}(1)$. The strong form implies the HLR
criterion. The name records the Hardy-Littlewood condition of the first
Tauberian theorems and the Ramanujan bound\index[terms]{Ramanujan bound}\index[names]{Ramanujan, S.} on the coefficients of
Dirichlet series\index[terms]{Dirichlet series} in the Selberg\index[names]{Selberg, A.} class\index[terms]{Selberg class}.
\end{definition}

The criterion is a property of the profile, tested on the solutions of the defining equation and
on nothing else. A theorem that needs a bound on $na_n$ for some other sequence has to check
that bound for that sequence, and the applications below do so one by one. The defining
equation is what makes the strong bound available at all. Under an input that only converges,
$A_\Phi(n)=\ell+o(1)$ with no rate, the best individual bound is $a_n=o(\log\log n)$, and it
cannot be improved \cite{ErdosSegal1978}.

The Ingham function meets the strong criterion, and the proof uses no information on the zeros
of $\zeta$.

\begin{proposition}\label{prop:hlr_ingham}
The Ingham function satisfies the strong HLR criterion, hence the HLR
criterion.
\end{proposition}

\begin{proof}
At $\beta=0$ the solution of the defining equation is $a=(1,0,0,\dots)$, since
$\Phi(1/n)=1$, and $na_n$ takes the values $1$ and $0$. Let $\beta>0$ and put $b_n=na_n$.
Multiplying the defining equation by $n$ and counting the integers $i$ with $ik\le n$ gives
$\sum_{k\le n}b_k\lfloor n/k\rfloor=n^{1-\beta}$, that is
$\sum_{m\le n}(b\star\mathbf 1)(m)=n^{1-\beta}$. Differencing at the ranks $n$ and $n-1$,
\[
(b\star\mathbf 1)(n)=c(n),\qquad c(1)=1,\quad c(n)=n^{1-\beta}-(n-1)^{1-\beta}\ \ (n\ge2),
\]
and M\"obius inversion gives $b=\mu\star c$. Write $c(n)=(1-\beta)n^{-\beta}+u(n)$. Then
$u(1)=\beta$, and for $n\ge2$ the mean value theorem, applied to $t\mapsto t^{1-\beta}$ and then
to $t\mapsto t^{-\beta}$, gives $u(n)=\mathcal O_\beta(n^{-\beta-1})$, so that $u$ is summable.
The function
\[
J_{-\beta}(n):=\sum_{d\mid n}d^{-\beta}\mu\Big(\frac nd\Big)
=\prod_{p^{a}\|n}\big(p^{-a\beta}-p^{-(a-1)\beta}\big)
\]
is multiplicative, and each factor lies in $[-1,1]$. Since $b=(1-\beta)J_{-\beta}+\mu\star u$,
\[
|na_n|\le|1-\beta|+\sum_{d\ge1}|u(d)|,
\]
and $na_n=\mathcal O_\beta(1)$ for each fixed $\beta$.
\end{proof}

The strong form is sharp. At a prime the sum over divisors reduces to two terms, so
\[
  pa_p=\big(p^{1-\beta}-(p-1)^{1-\beta}\big)-1\longrightarrow-1\qquad(p\to\infty),
\]
and the bound $na_n=\mathcal{O}(1)$ cannot be improved to $na_n=o(1)$. The value $-1$ is the
contribution of the term $\mu(p)$ in the M\"obius inversion, and it persists for every $\beta>0$.

The name of the criterion records two classical strands. The strong form $na_n=\mathcal{O}(1)$
is the condition under which Littlewood\index[names]{Littlewood, J. E.} established the converse
of Abel's\index[names]{Abel, N. H.} theorem, that Abel summability of $\sum a_n$ together with
$a_n=\mathcal{O}(1/n)$ forces convergence \cite{Littlewood1911}. Littlewood's theorem relaxes
the earlier condition $na_n=o(1)$ of Tauber\index[names]{Tauber, A.} \cite{Tauber1897}, and Hardy
and Littlewood later reduced it to the one-sided form $na_n\ge-C$ \cite{HardyLittlewood1914}. The weak form
$na_n=\mathcal{O}(n^{\eps})$ is the Ramanujan bound imposed on the Dirichlet coefficients of the
members of the Selberg\index[names]{Selberg, A.} class \cite{Selberg1992}. The survey of
Korevaar\index[names]{Korevaar, J.} \cite{Korevaar2004} covers the Tauberian side. Both forms
hold for the Ingham function, which is the object of Proposition~\ref{prop:hlr_ingham}, and they
separate for its arithmetic generalizations, which \S\ref{sec:gen_ingham} studies. Here the
strong form is used as a Tauberian condition\index[terms]{Tauberian theorem}. The next two sections show that it suffices for a
boundary theorem, and that for an arbitrary sequence the weak form does not.

\section{Wiener's theorem in the logarithmic variable}\label{sec:wiener_transfer}

The profiles of this section are normalized at the origin. The normalization is local to the
chapter. It fixes the value $g^{*}(0)=g(0^{+})$ of the note on the normalization of $g^{*}$ to
be one, and it asks nothing of $g(1)$.

\begin{definition}\label{def:normalized_origin}
A profile $g:(0,1]\to\C$ is normalized at the origin\index[terms]{normalized at the origin} when it
is bounded, Riemann integrable on every interval $[\eps,1]$ with $\eps>0$, and satisfies, for
some $\delta>0$,
\begin{equation}\label{eq:normalized_origin}
g(t)=1+\mathcal O(t^{\delta})\qquad(t\to0^{+}).
\end{equation}
\end{definition}

For such a profile the transform takes the form
\begin{equation}\label{eq:normalized_transform}
g^{*}(z)=1-z\int_0^1\big(g(t)-1\big)\,t^{-z-1}\,dt ,
\end{equation}
obtained by subtracting the constant part, $-z\int_0^1t^{-z-1}\,dt=1$ for $\Re z<0$. The integral
converges for $\Re z<\delta$, so $g^{*}$ is holomorphic there and $g^{*}(0)=1$. Zeros and poles
of a continuation further to the right do not enter what follows. The Ingham function is
normalized at the origin with $\delta=1$, since $\Phi(t)=1-t\{1/t\}$ gives $|\Phi(t)-1|\le t$.

\begin{theorem}[Boundary transfer]\label{thm:boundary_transfer}
\index[terms]{boundary transfer}Let $g$ be a real profile normalized at the origin. The following
are equivalent.
\begin{enumerate}
\item[(i)] $g^{*}(i\xi)\neq0$ for every real $\xi$.
\item[(ii)] For every real sequence with $\sup_{n}n|a_n|<\infty$, $A_g(n)\to0$ implies $A(n)\to0$.
\end{enumerate}
When they hold, every complex sequence with $na_n=\mathcal O(1)$ satisfies
\begin{equation}\label{eq:boundary_transfer}
A_g(n)\longrightarrow L\qquad\Longrightarrow\qquad A(n)\longrightarrow L .
\end{equation}
Condition (i) implies \eqref{eq:boundary_transfer} for a complex profile normalized at the origin
as well.
\end{theorem}

No positivity or monotonicity of $g$ is assumed, and no bounded variation. The spectral condition
bears on the axis alone and not on the left half plane. The Fourier transform is
$\widehat f(\xi)=\int_{\R}f(v)e^{-i\xi v}\,dv$, the star $u*v$ denotes additive convolution on $\R$,
as in Chapter~\ref{chap:ortho}, and $\star$ keeps its meaning of Dirichlet convolution. The step that Wiener's theorem supplies in
general is isolated in the next lemma, in the form the present class requires.

\begin{lemma}\label{lem:wiener_step}
Let $w$ be integrable on $\R$ with $\int_{\R}|v|^{m}|w(v)|\,dv<\infty$ for every $m\ge0$,
$\int_{\R}w=1$ and $\widehat w(\xi)\neq0$ for every real $\xi$. Let $F$ be bounded and measurable
on $\R$ with $(w*F)(u)\to L$ as $u\to+\infty$. Then $(j*F)(u)\to L\int_{\R}j$ for every integrable
$j$.
\end{lemma}

\begin{proof}
The moments make $\widehat w$ smooth, each derivative being the transform of an integrable
function. Let $\varphi$ be a Schwartz function whose transform is smooth with compact support.
The quotient $\widehat\varphi/\widehat w$ is then smooth with compact support. Its inverse
transform $\psi_0$ is a Schwartz function, and the Fourier inversion theorem
\cite[\S7.2, (7.16), pp.~218--219]{Folland1992} gives $\widehat{\psi_0}=\widehat\varphi/\widehat w$,
so that $\varphi=w*\psi_0$ and
\[
(\varphi*F)(u)=\big(\psi_0*(w*F)\big)(u)\longrightarrow L\int_{\R}\psi_0
=L\,\frac{\widehat\varphi(0)}{\widehat w(0)}=L\int_{\R}\varphi
\]
by dominated convergence, $w*F$ being bounded and tending to $L$. The integrable $j$ with
$(j*F)(u)\to L\int j$ form a subspace, closed in $L^{1}(\R)$ because
$\|(j-j')*F\|_{\infty}\le\|j-j'\|_{1}\|F\|_{\infty}$. It contains every such $\varphi$, and these
are dense. Indeed, if $\vartheta$ is smooth with compact support and $\vartheta(0)=1$, a Schwartz
function $\varphi$ is the limit in $L^{1}(\R)$, as $R\to\infty$, of the functions with transforms
$\widehat\varphi(\xi)\vartheta(\xi/R)$, which are the convolutions of $\varphi$ with an approximate
identity, and Schwartz functions are dense in $L^{1}(\R)$. The subspace is therefore all of
$L^{1}(\R)$.
\end{proof}

In general $\widehat w$ is only continuous, and the division by $\widehat w$ is Wiener's
theorem\index[terms]{Wiener's theorem} on the translates of an integrable function whose transform
has no real zero \cite{Wiener1932}. The moments available below make it elementary.

\begin{proof}[Proof of Theorem~\ref{thm:boundary_transfer}, sufficiency]
It is enough to treat $L=0$. Put $M=\sup_n n|a_n|$. On the logarithmic scale set
\[
b(u)=e^{u}a_{\lfloor e^{u}\rfloor}\ \ (u\ge0),\qquad b(u)=0\ \ (u<0),\qquad
F(u)=\int_{-\infty}^{u}b(v)\,dv .
\]
Then $|b|\le2M$, so $F$ is Lipschitz with constant $2M$, and
\begin{equation}\label{eq:F_partial}
F(\log n)=\int_1^{n}a_{\lfloor y\rfloor}\,dy=A(n-1).
\end{equation}
Let
\[
K_g(v)=g(e^{-v})\,\mathbf 1_{[0,\infty)}(v),\qquad
h_g(v)=\big(g(e^{-v})-1\big)\mathbf 1_{[0,\infty)}(v).
\]
By \eqref{eq:normalized_origin} the function $h_g$ is $\mathcal O(e^{-\delta v})$ and has
moments of every order. The convolution $U=K_g*b$ is defined pointwise, and
\begin{equation}\label{eq:U_sampled}
U=F+h_g*b,\qquad U(\log n)=\int_1^{n}a_{\lfloor y\rfloor}\,g\Big(\frac yn\Big)\,dy .
\end{equation}

The first estimate compares the integral with the sum. The difference $U(\log n)-A_g(n)$ is at
most
\[
|a_n|\,|g(1)|+\sum_{k=1}^{n-1}|a_k|\int_k^{k+1}\Big|g\Big(\frac yn\Big)-g\Big(\frac kn\Big)\Big|\,dy .
\]
Fix $0<\eps<\tfrac14$. Enlarging the constant in \eqref{eq:normalized_origin} so that it holds on
all of $(0,1]$, the terms with $k\le\eps n$ contribute
\[
\ll n^{-\delta}\sum_{k\le\eps n}\frac{(k+1)^{\delta}}{k}=\mathcal O(\eps^{\delta})+o(1).
\]
The terms with $k>\eps n$ contribute at most
$(M/\eps n)\sum_{\eps n<k<n}\operatorname{osc}\big(g;[k/n,(k+1)/n]\big)$, which tends to zero by
the Riemann integrability of $g$ on $[\eps,1]$. Letting $n\to\infty$ and then $\eps\to0$,
\begin{equation}\label{eq:sampling}
U(\log n)-A_g(n)\longrightarrow0 .
\end{equation}

Consequently $U(\log n)\to0$. The convolution $h_g*b$ is uniformly continuous, since
$\|(h_g*b)(\cdot+s)-h_g*b\|_{\infty}\le\|b\|_{\infty}\|h_g(\cdot+s)-h_g\|_{1}\to0$, and $F$ is
Lipschitz, so $U$ is uniformly continuous. As $\log(n+1)-\log n\to0$, the limit $U(u)\to0$
holds along all real $u\to+\infty$. In particular $U$ is bounded, and \eqref{eq:U_sampled} then
bounds $F$. The boundedness of the partial sums is a conclusion here, not a hypothesis.

Put
\begin{equation}\label{eq:wiener_companion}
e_+(v)=e^{-v}\mathbf 1_{[0,\infty)}(v),\qquad w_g=e_++h_g-e_+*h_g .
\end{equation}
The function $w_g$ is the Wiener companion\index[terms]{Wiener companion} of $g$. It is integrable with moments of every
order, and $\int w_g=1$. Since $F'=b$ almost everywhere and $F$ vanishes on $(-\infty,0]$,
integration by parts gives $e_+*b=F-e_+*F$, and therefore
\begin{equation}\label{eq:companion_identity}
e_+*U=e_+*F+h_g*(e_+*b)=w_g*F .
\end{equation}
Hence $(w_g*F)(u)\to0$. On the Fourier side, $\widehat{e_+}(\xi)=1/(1+i\xi)$, and the substitution
$t=e^{-v}$ in \eqref{eq:normalized_transform} gives $g^{*}(-i\xi)=1+i\xi\,\widehat{h_g}(\xi)$, so
\begin{equation}\label{eq:companion_transform}
\widehat{w_g}(\xi)=\frac{1+i\xi\,\widehat{h_g}(\xi)}{1+i\xi}=\frac{g^{*}(-i\xi)}{1+i\xi}.
\end{equation}
Condition (i) is exactly the non-vanishing of $\widehat{w_g}$ on the real line, and
Lemma~\ref{lem:wiener_step} applies with $j=\ell^{-1}\mathbf 1_{[0,\ell]}$ for each fixed
$\ell>0$. The backward average $\ell^{-1}\int_0^{\ell}F(u-v)\,dv$ tends to zero, and it differs
from $F(u)$ by at most $M\ell$. Taking the upper limit as $u\to\infty$ and then letting
$\ell\to0$ gives $F(u)\to0$, and \eqref{eq:F_partial} with $a_n\to0$ gives $A(n)\to0$.

For a general limit $L$, replace $a_1$ by $a_1-L$. The averaged sum changes by
$-L\,g(1/n)\to-L$ and the partial sum by $-L$. For a complex profile or a complex sequence the
argument is unchanged, every estimate above bearing on absolute values.
\end{proof}

\begin{proof}[Proof of Theorem~\ref{thm:boundary_transfer}, necessity]
Suppose $g^{*}(-i\xi)=0$ for some real $\xi$, necessarily $\xi\neq0$ since $g^{*}(0)=1$. Set
$a_1=1$ and $a_n=n^{i\xi}-(n-1)^{i\xi}$ for $n\ge2$, so that $A(n)=n^{i\xi}$ and
$a_k=i\xi k^{i\xi-1}+\mathcal O_\xi(k^{-2})$ for $k\ge2$, whence $na_n=\mathcal O_\xi(1)$.
Writing $g=1+(g-1)$,
\[
A_g(n)=n^{i\xi}+\sum_{k\le n}a_k\Big(g\Big(\frac kn\Big)-1\Big)
=n^{i\xi}\Big(1+i\xi\int_0^1\big(g(t)-1\big)t^{i\xi-1}\,dt\Big)+o(1)
=n^{i\xi}g^{*}(-i\xi)+o(1).
\]
The term $k=1$ is $\mathcal O(n^{-\delta})$. The error terms $\mathcal O(k^{-2})$, multiplied by
$g(k/n)-1$, give $o(1)$ by dominated convergence. The principal sum is an improper Riemann sum
whose part near the origin is $\mathcal O(\eps^{\delta})$ uniformly, as in the sampling estimate.
Hence $A_g(n)\to0$ while $A(n)=n^{i\xi}$ does not converge. For real $g$ the real parts give a
real sequence with $na_n=\mathcal O(1)$, $A_g(n)\to0$ and $A(n)=\cos(\xi\log n)$, and the zero
at $-i\xi$ is the conjugate of a zero at $i\xi$.
\end{proof}

Three examples mark the limits of the statement. The profile equal to $0$ at the points $1/m$,
$m\ge2$, and to $1$ elsewhere has $g^{*}\equiv1$, entire and without zero, yet the sequence
$a=(1,0,0,\dots)$ gives $A_g(n)=g(1/n)=0$ for $n\ge2$ and $A(n)=1$. Integral data do not see the
values of $g$ on the sampling grid, and the condition at the origin cannot be removed. The profile
$g=\tfrac12\big(1+\mathbf 1_{(0,1/2]}\big)$ is normalized at the origin and has
$g^{*}(z)=\tfrac12(1+2^{z})$, with a zero at $i\pi/\log2$. Here
$A_g(n)=\tfrac12\big(A(n)+A(\lfloor n/2\rfloor)\big)$, and $A(n)=\cos(\pi\log n/\log2)$ makes the
average tend to zero while the partial sums oscillate. The profile $g=2\cdot\mathbf 1_{(0,1/2]}-1$
has $g^{*}(z)=2^{z+1}-1$, whose zeros all have real part $-1$. None lies on the axis, and the
theorem applies. A zero-free left half plane would be the wrong hypothesis.

The companion $w_g$ is not the Volterra kernel $k_g$ of \S\ref{sec:absorption}, which acts by
multiplicative convolution and serves the inversion on the half plane. The relation between the
two is recorded in Remark~\ref{rem:volterra_wiener}.

\section{Ingham's theorem and the prime number theorem}\label{sec:ingham_proved}

The criterion of \S\ref{sec:boundary_problem} and the theorem of \S\ref{sec:wiener_transfer}
combine at once on the defining equation.

\begin{corollary}\label{cor:boundary_absorption}
Let $g$ be a real profile normalized at the origin, with $g(1)\neq0$, satisfying condition (i)
of Theorem~\ref{thm:boundary_transfer} and the strong HLR criterion. For every $\beta>0$ the
solution of the defining equation \eqref{eq:defining_relation} with $r(n)=n^{-\beta}$ satisfies
$A(n)\to0$, and for $\beta=0$ it satisfies $A(n)\to1$.
\end{corollary}

\begin{proof}
Apply \eqref{eq:boundary_transfer} with the limits $0$ and $1$.
\end{proof}

This is a statement at the boundary. It gives no rate $A(n)=\mathcal O(n^{-c})$ with $c>0$,
it does not produce a regularity index, and it does not identify a threshold with a zero of the
transform. In Theorem~\ref{thm:boundary_transfer} the bound $na_n=\mathcal O(1)$ cannot be relaxed
to $na_n=o(n^{\eps})$. The profile $g(t)=1-t$ on $(0,1)$ with $g(1)=1$ is normalized at the origin and has $g^{*}(z)=1/(1-z)$, which has no zero at all.
Taking $A(n)=\sin\big((\log n)^{2}\big)$ and $a_n=A(n)-A(n-1)$ gives $na_n=\mathcal O(\log n)$,
which is $o(n^{\eps})$ for every $\eps>0$, while summation by parts gives
$A_g(n)=n^{-1}\sum_{j<n}A(j)+a_n$. Comparing that mean with
$\int_1^{n}\sin\big((\log t)^{2}\big)\,dt=\mathcal O(n/\log n)$, obtained by one integration by
parts in the variable $\log t$, shows $A_g(n)\to0$, and the partial sums do not converge. The
example does not decide whether the weak criterion suffices in
Corollary~\ref{cor:boundary_absorption}, where the forcing is a pure power.

Ingham's own condition is one-sided, and it is met by a sequence the two-sided theorem cannot
reach. The next theorem trades the two-sided bound for the positivity of the companion, and it
asks for the average at every real argument. For a discontinuous profile the
real and the integer premises are different statements, and for the Ingham function they turn out
to coincide.

\begin{theorem}[One-sided transfer]\label{thm:onesided_transfer}
Let $g$ be a real profile normalized at the origin with $g^{*}(i\xi)\neq0$ for every real $\xi$
and $w_g\ge0$ almost everywhere. For a real sequence with
$na_n\ge-C$,
\[
A_g(x)\longrightarrow L\quad(x\to\infty\text{ through real values})
\qquad\Longrightarrow\qquad A(n)\longrightarrow L .
\]
\end{theorem}

\begin{proof}
This time let $F$ be the step function
\[
F(u)=\sum_{\log n\le u}a_n\ \ (u\ge0),\qquad F(u)=0\ \ (u<0),
\]
and put $U(u)=A_g(e^{u})$ for $u\ge0$ and $U(u)=0$ otherwise. Then $U$ is bounded and tends
to $L$. The lower bound on the terms gives, for a constant $C_1$
depending on $C$,
\begin{equation}\label{eq:slow_global}
F(u-s)\le F(u)+Cs+C_1\qquad(u,s\ge0),
\end{equation}
and, comparing harmonic sums with $\log(e^{u+s}/e^{u})=s$, for each fixed $\ell$,
\begin{equation}\label{eq:slow_local}
F(u+s)-F(u)\ge-Cs-\mathcal O(e^{-u})\qquad(0\le s\le\ell).
\end{equation}

Let $J=e_+*K_g$. Then $J(0)=0$ and $J'=K_g-J=w_g$ almost everywhere, by the definition
\eqref{eq:wiener_companion}. Summing over the finitely many terms and integrating,
\[
(w_g*F)(u)=\sum_{\log n\le u}a_n\int_0^{u-\log n}w_g(v)\,dv
=\sum_{\log n\le u}a_n\,J(u-\log n)=(e_+*U)(u),
\]
so $w_g*F$ is bounded and tends to $L$, both functions vanishing on the negative half line.

The function $F$ is bounded. Integrating \eqref{eq:slow_global} against the non-negative $w_g$
gives $(w_g*F)(u)\le F(u)+C\int_0^{\infty}v\,w_g(v)\,dv+C_1$, so $F\ge-B$ for some $B$. Choose
$\ell_0>0$ with $m_0=\int_0^{\ell_0}w_g>0$. For $0\le v\le\ell_0$, \eqref{eq:slow_global} gives
$F(u+\ell_0-v)\ge F(u)-C\ell_0-C_1$, and the remaining values of $F$ are at least $-B$, so
\[
(w_g*F)(u+\ell_0)\ge m_0\big(F(u)-C\ell_0-C_1\big)-(1-m_0)B ,
\]
an upper bound for $F$ since $w_g*F$ is bounded.

Lemma~\ref{lem:wiener_step} now applies to $w_g$, whose transform \eqref{eq:companion_transform}
has no real zero. The forward and the backward averages of $F$ over an interval of fixed length
$\ell>0$ both tend to $L$. By \eqref{eq:slow_local}, $F(u)$ is at most its forward average plus
$C\ell+o(1)$ and at least its backward average minus $C\ell+o(1)$. Letting $u\to\infty$ and then
$\ell\to0$ gives $F(u)\to L$, and $F(\log n)=A(n)$.
\end{proof}

For the Ingham function the positivity of the companion can be read off directly. With $x=e^{v}$ and
$n=\lfloor x\rfloor$,
\[
K_\Phi(v)=\frac nx,\qquad (e_+*K_\Phi)(v)=\frac1x\int_1^{x}\frac{\lfloor y\rfloor}{y}\,dy
=\frac{n\log x-\log(n!)}{x},
\]
and therefore
\begin{equation}\label{eq:companion_ingham}
w_\Phi(v)=\frac{n-n\log x+\log(n!)}{x}>0 .
\end{equation}
Indeed $\log(n!)\ge n\log n-n+1$, by comparison with $\int_1^{n}\log t\,dt$, and $x<n+1$, so the
numerator is at least $1-n\log(x/n)>1-n\log(1+1/n)>0$. Finally, for every sequence and $m=\lfloor x\rfloor$,
\begin{equation}\label{eq:real_integer_ingham}
A_\Phi(x)=\frac1x\sum_{k\le m}k\,a_k\Big\lfloor\frac xk\Big\rfloor=\frac mx\,A_\Phi(m),
\end{equation}
since $\lfloor x/k\rfloor=\lfloor m/k\rfloor$ for an integer $k$, so convergence at the integers
gives convergence along the reals.

What remains is the spectral condition, $\Phi^{*}(i\xi)\neq0$, which by
Proposition~\ref{prop:ingham_mellin} is $\zeta(1-i\xi)\neq0$ for $\xi\neq0$. That is the
classical theorem of Hadamard\index[names]{Hadamard, J.} and de la Vall\'ee Poussin\index[names]{Vall\'ee Poussin, C.-J. de la@de la Vall\'ee Poussin, C.-J.} \cite[Ch.~III]{Titchmarsh1986}, and its
short proof uses no Tauberian argument.

\begin{proposition}\label{prop:zeta_line}
$\zeta(1+it)\neq0$ for every real $t\neq0$. Consequently $\Phi^{*}(i\xi)\neq0$ for every
real $\xi$.
\end{proposition}

\begin{proof}
The identity $\zeta(s)=s/(s-1)-s\int_1^{\infty}\{x\}x^{-s-1}\,dx$ continues $\zeta$ to $\Re s>0$
with a simple pole of residue $1$ at $s=1$. For $\sigma>1$ the Euler product gives
$\log|\zeta(s)|=\Re\sum_{p}\sum_{k\ge1}k^{-1}p^{-ks}$, and the inequality
$3+4\cos\theta+\cos2\theta=2(1+\cos\theta)^{2}\ge0$, applied with $\theta=kt\log p$, gives
\[
\zeta(\sigma)^{3}\,|\zeta(\sigma+it)|^{4}\,|\zeta(\sigma+2it)|\ge1 .
\]
If $t\neq0$ and $\zeta(1+it)=0$, the left side tends to zero as $\sigma\to1^{+}$, the pole
contributing an order $-3$, the zero an order at least $4$, and the last factor staying bounded.
This contradicts the inequality. At $\xi=0$, $\Phi^{*}(0)=1$.
\end{proof}

\begin{corollary}[Ingham's theorem]\label{cor:ingham_proved}
Let $(a_n)$ be real with $na_n\ge-C$. If $A_\Phi(n)\to L$, then $A(n)\to L$. In particular
Theorem~\ref{thm:ingham} holds.
\end{corollary}

\begin{proof}
Proposition~\ref{prop:zeta_line}, the positivity \eqref{eq:companion_ingham} and the passage \eqref{eq:real_integer_ingham} from integer to real
arguments verify the hypotheses of Theorem~\ref{thm:onesided_transfer}. The average $I(x)$ of
\eqref{eq:ingham_averages} is $A_\Phi(x)$.
\end{proof}

The theorem quoted in Chapter~\ref{chap:ingham} is thereby proved from the elementary case of
Wiener's theorem, Lemma~\ref{lem:wiener_step}, with an explicit companion, and the zeros of $\zeta$ enter only through the line $\Re s=1$. Combined with
Proposition~\ref{prop:hlr_ingham}, Corollary~\ref{cor:boundary_absorption} gives for the Ingham
function, with no hypothesis, $A(n)\to0$ at every forcing exponent $\beta>0$. That limit is all the
boundary yields. Its rate is the subject of Chapter~\ref{chap:equivalence}.

The three sequences of \S\ref{sec:two_identities} now give the prime number theorem, each
under the condition on its terms that it actually satisfies.

\begin{proposition}\label{prop:pnt_boundary}
\index[terms]{prime number theorem}The partial sums of the three sequences of
\S\ref{sec:two_identities} converge, with
\[
\sum_{n\le x}\frac{\mu(n)}{n}\longrightarrow0,\qquad
\sum_{n\le x}\frac{\lambda(n)}{n}\longrightarrow0,\qquad
\sum_{n\le x}\frac{\Lambda(n)-1}{n}\longrightarrow-2\gamma .
\]
Consequently $\psi(x)\sim x$ and $\pi(x)\sim x/\log x$.
\end{proposition}

\begin{proof}
The M\"obius sequence $a_n=\mu(n)/n$ is the solution of the defining equation at $\beta=1$, by
the identity of Meissel, and $|na_n|\le1$. Corollary~\ref{cor:boundary_absorption}, or directly
Theorem~\ref{thm:boundary_transfer} with $A_\Phi(n)=1/n$, gives $M_{-1}(x)\to0$, and Abel
summation gives $M(x)=xM_{-1}(x)-\int_1^{x}M_{-1}(t)\,dt=o(x)$. The Liouville sequence
$a_n=\lambda(n)/n$ has $|na_n|=1$ and $A_\Phi(n)=\lfloor\sqrt n\rfloor/n\to0$, a source that is not
a pure power, so the bound on the terms is checked here directly. Theorem~\ref{thm:boundary_transfer}
gives the second limit, and the same Abel summation gives $\sum_{n\le x}\lambda(n)=o(x)$. The
identity $\mu(n)=\sum_{d^{2}\mid n}\mu(d)\lambda(n/d^{2})$, which is $1/\zeta(s)=\zeta(2s)^{-1}\cdot
\zeta(2s)/\zeta(s)$ read on coefficients, then gives $M(x)=\sum_{d\le\sqrt x}\mu(d)\sum_{m\le
x/d^{2}}\lambda(m)=o(x)$, each fixed $d$ contributing $o(x)$ and the tail at most
$x\sum_{d>D}d^{-2}$. The statement $M(x)=o(x)$ is one of the classical equivalents of the prime
number theorem \cite[\S~II.4]{Tenenbaum2015}, and Proposition~\ref{prop:mobius_bridge} of
Appendix~\ref{app:boundary} turns it into $\psi(x)\sim x$ by elementary means.

For the von Mangoldt sequence\index[terms]{von Mangoldt function} $a_n=(\Lambda(n)-1)/n$ the
product $na_n=\Lambda(n)-1$ is not bounded, and the two-sided theorem does not apply to it. It satisfies the
one-sided bound $na_n\ge-1$, and its averages converge to $-2\gamma$ by the elementary computation
of \S\ref{sec:ingham_pnt}, through Stirling's formula and Dirichlet's divisor estimate.
Corollary~\ref{cor:ingham_proved} gives the third limit. With
$S(x)=\sum_{n\le x}(\Lambda(n)-1)/n$, partial summation gives
$\sum_{n\le x}(\Lambda(n)-1)=xS(x)-\int_1^{x}S(t)\,dt=o(x)$, that is $\psi(x)=x+o(x)$.

Finally the prime powers of exponent at least two contribute $\mathcal O(\sqrt x\,\log^{2}x)$ to
$\psi(x)$, so $\theta(x)=\sum_{p\le x}\log p\sim x$, and splitting the primes at $x/\log^{2}x$ gives
\[
\frac{\theta(x)}{\log x}\le\pi(x)\le\frac{x}{\log^{2}x}+\frac{\theta(x)}{\log x-2\log\log x},
\]
whence $\pi(x)\sim x/\log x$.
\end{proof}

Each route checks its own condition on the terms. The M\"obius and Liouville sequences are
bounded after multiplication by $n$ and use the two-sided theorem. The von Mangoldt sequence is
bounded below only, and it needs the one-sided theorem with the positivity of $w_\Phi$. The
analytic input is the same in all three, the non-vanishing of $\zeta$ on $\Re s=1$, and it enters
once, through the transform of the weight on the boundary of its half plane.

\section{The Mellin transform of the averaged sum}\label{sec:averaged_mellin}

The transform of \S\ref{sec:mellin_transform} belongs to the weight. The averaged sum has a Mellin
transform of its own, and it separates the sequence from the weight. This is the separation of
harmonic sums\index[terms]{harmonic sum} organized by \nm{Flajolet}{P.}, \nm{Gourdon}{X.} and
\nm{Dumas}{P.}~\cite{FlajoletGourdonDumas1995}, recalled in \S\ref{sec:defining_recurrence} as a
device for asymptotic expansions. Here it serves another purpose. It is an identity on a half
plane, and recovering a limit from it on the edge of that half plane is a Tauberian\index[terms]{Tauberian} step, which
positivity of the terms makes available without any bound on their size.

For a bounded Borel measurable profile $g$ on $(0,1]$, not necessarily normalized at the origin, put
\[
m_g(s)=\int_0^1g(t)\,t^{s-1}\,dt=\frac{g^{*}(-s)}{s}\qquad(\Re s>0).
\]
If the Dirichlet series $D_a(s)=\sum_{n\ge1}a_nn^{-s}$ converges absolutely at $s=\sigma>0$, then
\begin{equation}\label{eq:fgd_factorization}
\int_1^{\infty}A_g(x)\,x^{-s-1}\,dx=m_g(s)\,D_a(s)=\frac{g^{*}(-s)}{s}\,D_a(s)
\qquad(\Re s=\sigma).
\end{equation}
Indeed the term of $a_n$ contributes $\int_n^{\infty}g(n/x)x^{-s-1}\,dx=n^{-s}m_g(s)$, and the
sum of the absolute values is at most $\|g\|_{\infty}\sigma^{-1}\sum_n|a_n|n^{-\sigma}$. The same
computation applies to a locally finite positive Borel measure $\nu$ on $[1,\infty)$, with
\begin{gather}
D_\nu(s)=\int_{[1,\infty)}y^{-s}\,d\nu(y),\qquad V_g(x)=\int_{[1,x]}g\Big(\frac yx\Big)\,d\nu(y),\notag\\
\int_1^{\infty}V_g(x)\,x^{-s-1}\,dx=m_g(s)\,D_\nu(s)\qquad(\Re s=\sigma,\ D_\nu(\sigma)<\infty).\label{eq:fgd_measure}
\end{gather}
For $a_n=b_n/n$ one has $D_a(s)=D_b(1+s)$, and this shift is what separates the axis of
Theorem~\ref{thm:boundary_transfer} from the line $\Re s=1$ of the problem below. The
factorization is an identity and not an asymptotic statement. The boundary step that turns it
into one rests on the following lemma, proved in Appendix~\ref{app:boundary}.

\begin{lemma}\label{lem:positive_boundary}
Let $\kappa>0$ and let $\nu$ be a locally finite positive Borel measure on $[1,\infty)$ with
$D_\nu(s)$ finite for real $s>\kappa$. If, for some $R\ge0$, the function
$D_\nu(s)-R/(s-\kappa)$ has a continuous extension from $\Re s>\kappa$ to the closed half plane
$\Re s\ge\kappa$, then
\[
\nu([1,x])=\frac{R}{\kappa}\,x^{\kappa}+o(x^{\kappa}).
\]
\end{lemma}

The continuity is local, on compact subsets of the closed half plane, and no growth estimate is
assumed on the vertical line.

\begin{theorem}[Positive transfer]\label{thm:positive_transfer}
Let $\kappa>0$, let $\nu\ge0$ be as in Lemma~\ref{lem:positive_boundary} with $D_\nu(s)$
convergent for $\Re s>\kappa$, and let $g$ be bounded and Borel measurable on $(0,1]$ with
$m_g(\kappa+i\xi)\neq0$ for every real $\xi$. If, for some $C\ge0$, the function
\begin{equation}\label{eq:positive_pole}
\int_1^{\infty}V_g(x)\,x^{-s-1}\,dx-\frac{C\,g^{*}(-\kappa)}{s-\kappa}
\end{equation}
has a continuous extension to $\Re s\ge\kappa$, then $\nu([1,x])=Cx^{\kappa}+o(x^{\kappa})$.
\end{theorem}

\begin{proof}
Write $E_g(s)$ for \eqref{eq:positive_pole} and note $g^{*}(-\kappa)=\kappa\,m_g(\kappa)$. The
function $m_g$ is holomorphic on $\Re s>0$ and has no zero on the line $\Re s=\kappa$, hence none
in a neighbourhood of each of its points. There, by \eqref{eq:fgd_measure},
\[
D_\nu(s)-\frac{C\kappa}{s-\kappa}
=\frac{E_g(s)}{m_g(s)}+C\kappa\,\frac{m_g(\kappa)-m_g(s)}{(s-\kappa)\,m_g(s)},
\]
and the second quotient has a removable singularity at $s=\kappa$. Both terms extend continuously
to the line, and at the points of $\Re s>\kappa$ the series $D_\nu$ is already holomorphic, zeros
of $m_g$ there included. Lemma~\ref{lem:positive_boundary} with $R=C\kappa$ concludes.
\end{proof}

For $g\equiv1$ one has $m_g(s)=1/s$, $g^{*}\equiv1$ and $V_g(x)=\nu([1,x])$, and the hypothesis
reduces to the continuity of $D_\nu(s)-C\kappa/(s-\kappa)$ on $\Re s\ge\kappa$. The theorem then
contains the continuous-boundary form of the theorem of Wiener\index[names]{Wiener, N.} and
Ikehara\index[names]{Ikehara, S.}\index[terms]{Wiener--Ikehara theorem} \cite{Ikehara1931,Wiener1932}, for
arbitrary positive measures and hence for generalized Dirichlet series. With
$\nu=\sum_n\Lambda(n)\delta_n$ and $\kappa=1$, $D_\nu=-\zeta'/\zeta$, and
Proposition~\ref{prop:zeta_line} together with the simple pole of $\zeta$ at $1$ makes
$-\zeta'(s)/\zeta(s)-1/(s-1)$ holomorphic near every point of $\Re s\ge1$. The theorem gives
$\psi(x)=x+o(x)$ a fourth time, now with positivity in place of any condition on the size of the
terms.

The two branches test different lines. The two-sided and one-sided theorems test $g^{*}(i\xi)$,
on the axis. The positive theorem at the growth exponent $\kappa$ tests $m_g(\kappa+i\xi)$, that
is $g^{*}(-\kappa-i\xi)$, on another line. Neither theorem contains the other, and in neither
does a meromorphic continuation replace the Tauberian step. The positive branch is stated with
the classical continuous boundary hypothesis, and its refinements with weaker boundary behaviour
are not used here.

\section{Characters at the boundary}\label{sec:boundary_characters}

The sufficiency half of the two-sided theorem holds for complex profiles, and a character
supplies one whose transform is a Dirichlet $L$-function. Let $Q\ge2$ be fixed and let $\chi$ be a nonprincipal Dirichlet
character\index[terms]{Dirichlet character} modulo $Q$, extended by zero to the integers not prime
to $Q$. Its harmonic character profile\index[terms]{harmonic character profile} is
\[
g_\chi(t)=\sum_{m\le1/t}\frac{\chi(m)}{m}\qquad(0<t\le1).
\]
The partial sums of $\chi$ are bounded, by periodicity and zero mean, so summation by parts gives
$g_\chi(t)=L(1,\chi)+\mathcal O_Q(t)$. The value $L(1,\chi)$ is not zero by
Lemma~\ref{lem:L_line}, and $\widetilde g_\chi=g_\chi/L(1,\chi)$ is normalized at the origin with
$\delta=1$. Termwise integration, justified by absolute convergence on $\Re z<0$, gives
\[
g_\chi^{*}(z)=L(1-z,\chi),\qquad \widetilde g_\chi^{\,*}(z)=\frac{L(1-z,\chi)}{L(1,\chi)},
\]
and the bound $\mathcal O_Q(t)$ continues the identity to $\Re z<1$. The boundary condition of
Theorem~\ref{thm:boundary_transfer} for $\widetilde g_\chi$ is therefore
$L(1+it,\chi)\neq0$ for every real $t$, which Lemma~\ref{lem:L_line} supplies without any
appeal to primes in progressions.

The profile $g_\chi$ is not the generalized Ingham function of \S\ref{sec:gen_ingham},
\[
\Phi_\chi(t)=\sum_{m\le1/t}\frac{\chi(m)}{m}\,\Phi(mt),
\]
whose transform $\Phi^{*}(z)L(1-z,\chi)$ carries $\zeta$ and the $L$-function together. The harmonic character profile isolates the
$L$-function on the boundary. For the principal character it is not available, since
$g_{\chi_0}(t)=(\varphi(Q)/Q)\log(1/t)+\mathcal O_Q(1)$ has no finite limit at the origin, and the
principal mode is carried by the Ingham function.

\begin{proposition}\label{prop:twisted_mobius}
For every nonprincipal character $\chi$ modulo $Q$,
\[
\sum_{n\le x}\frac{\mu(n)\chi(n)}{n}\longrightarrow\frac{1}{L(1,\chi)},\qquad
\sum_{n\le x}\mu(n)\chi(n)=o(x).
\]
\end{proposition}

\begin{proof}
Take $a_n=\mu(n)\chi(n)/n$, so that $n|a_n|\le1$. Complete multiplicativity of $\chi$ and the
divisor identity for $\mu$ give, at every real $x\ge1$, the exact identity
\[
A_{g_\chi}(x)=\sum_{km\le x}\frac{\mu(k)\chi(k)\chi(m)}{km}
=\sum_{r\le x}\frac{\chi(r)}{r}\sum_{k\mid r}\mu(k)=1 .
\]
The averaged sum against $\widetilde g_\chi$ is therefore the constant $1/L(1,\chi)$, and
Theorem~\ref{thm:boundary_transfer}, in its form for complex profiles, gives the first limit. With
$M_{-1}^{\chi}(N)=\sum_{n\le N}\mu(n)\chi(n)/n$, Abel summation gives
$\sum_{n\le N}\mu(n)\chi(n)=NM_{-1}^{\chi}(N)-\sum_{j<N}M_{-1}^{\chi}(j)=o(N)$, the limit
cancelling.
\end{proof}

The bound on the terms has been checked for this sequence alone, and no strong HLR property of
$\widetilde g_\chi$ for other forcings is inferred. In the factorization
\eqref{eq:fgd_factorization} the identity reads
$1/s=\big(L(1+s,\chi)/s\big)\cdot\big(1/L(1+s,\chi)\big)$ on $\Re s>0$, and the passage to
ordinary convergence remains the Tauberian part of the argument. The same profile applied to
$\lambda(n)\chi(n)/n$ gives the exact square-counting source
$\sum_{j\le\sqrt x}\chi(j)^{2}/j^{2}\to L(2,\chi^{2})$, and the transfer gives
$\sum_{n\le x}\lambda(n)\chi(n)/n\to L(2,\chi^{2})/L(1,\chi)$.

\begin{proposition}\label{prop:pnt_progressions}
\index[terms]{prime in an arithmetic progression}Let $Q\ge2$ be fixed and let $a$ be prime to $Q$.
Then
\[
\psi(x;Q,a)=\sum_{\substack{n\le x\\ n\equiv a\bmod Q}}\Lambda(n)=\frac{x}{\varphi(Q)}+o(x),\qquad
\pi(x;Q,a)\sim\frac{x}{\varphi(Q)\log x}.
\]
\end{proposition}

\begin{proof}
Let $\chi$ be nonprincipal and put $R_\chi(y)=\sum_{m\le y}\chi(m)\log m$, which is
$\mathcal O_Q(\log2y)$ by summation by parts. The identity $(\mu\chi)\star(\chi\log)=\chi\Lambda$,
read on the Dirichlet series as $(1/L)\cdot(-L')=-L'/L$, gives
\[
\psi_\chi(x):=\sum_{n\le x}\chi(n)\Lambda(n)=\sum_{d\le x}\mu(d)\chi(d)\,R_\chi\Big(\frac xd\Big).
\]
Lemma~\ref{lem:weighted_cancellation}, applied with $b_d=\mu(d)\chi(d)$ and $f(t)=R_\chi(1/t)$,
gives $\psi_\chi(x)=o(x)$ from Proposition~\ref{prop:twisted_mobius}. For the principal character,
$\psi_{\chi_0}(x)=\psi(x)-\sum_{p\mid Q}\sum_{p^{j}\le x}\log p=x+o(x)$ by
Proposition~\ref{prop:pnt_boundary}. The orthogonality of the characters modulo $Q$ gives
\[
\psi(x;Q,a)=\frac{1}{\varphi(Q)}\sum_{\chi\bmod Q}\overline{\chi(a)}\,\psi_\chi(x)
=\frac{x}{\varphi(Q)}+o(x).
\]
The prime powers of exponent at least two contribute $o(x)$, and the comparison between
$\theta$ and $\pi$ in the proof of Proposition~\ref{prop:pnt_boundary}, restricted to the class
of $a$, gives the second statement.
\end{proof}

The modulus is fixed and no uniformity in $Q$ is asserted. For the von Mangoldt sequence the
twisted analogue of the source of \S\ref{sec:ingham_pnt} exists, but $|na_n|=|\chi(n)|\Lambda(n)$
is not bounded, so the two-sided theorem does not apply. Nor does the one-sided theorem, which is a
statement about real terms. For a complex character the terms are not real, and for a real
character a bound on $\chi(n)\Lambda(n)$ from either side would leave only finitely many primes $p$
with $\chi(p)=-1$, or with $\chi(p)=1$, which the proposition excludes. The route therefore goes
through the M\"obius function instead. The
positive branch does apply, with $g\equiv1$, to the measure carried by the prime powers of the
class of $a$, whose Dirichlet series differs from $\varphi(Q)^{-1}(s-1)^{-1}$ by a function
holomorphic near $\Re s\ge1$, by Lemma~\ref{lem:L_line} and Proposition~\ref{prop:zeta_line},
and Theorem~\ref{thm:positive_transfer} gives the proposition again. Nothing in the proof
uses primes in progressions, Lemma~\ref{lem:L_line} resting on the Euler product and on a theorem
of Landau\index[names]{Landau, E.} on series with non-negative coefficients. The equidistribution
modulo four used in \S\ref{sec:gen_ingham}, in Lemma~\ref{lem:gen_chi4_split} and
Theorem~\ref{thm:gen_unbounded}, is the case $Q=4$ of the proposition. The converse order is
excluded. The results on $\Phi_{\chi_4}$ use that equidistribution and are never used to obtain it.

\section{From the boundary to a strip}\label{sec:boundary_to_strip}

Return to the Ingham function and to the M\"obius sequence $a_n=\mu(n)/n$, the solution at
$\beta=1$, whose average is exactly $1/n$. On the boundary the transfer returns $M_{-1}(x)\to0$
and nothing more, and the condition it uses is that $\Phi^{*}$ have no zero on the imaginary
axis, that is $\zeta$ none on the line $\Re s=1$. In the variable $z=1-s$ of the transform,
asking for more decay means moving that condition to the right, into the critical strip. If
$\zeta$ has no zero in a half plane $\Re s>\sigma$ with $\tfrac12<\sigma<1$, the proof of
Corollary~\ref{cor:zerofree}, which rests on the contour estimates of Appendix~\ref{app:perron},
gives at $\beta=1$
\[
M_{-1}(x)=\mathcal O_{\sigma,\eps}\big(x^{-(1-\sigma)+\eps}\big).
\]
At the critical depth the bound $M_{-1}(x)\ll_\eps x^{-1/2+\eps}$ is equivalent to the Riemann
hypothesis, by Abel summation, Lemma~\ref{lem:littlewood} and the value
$\sum_n\mu(n)/n=0$ of Proposition~\ref{prop:pnt_boundary}. The three levels are set side by side below and drawn in
Figure~\ref{fig:boundary_planes}.

\begin{center}
\begin{tabular}{lll}
\toprule
No zero of $\zeta$ & What follows & Status of the hypothesis \\
\midrule
on the line $\Re s=1$ & $M_{-1}(x)\to0$ & proved, Proposition~\ref{prop:zeta_line} \\
in $\Re s>\sigma$, $\tfrac12<\sigma<1$ & $M_{-1}(x)\ll x^{-(1-\sigma)+\eps}$, $\tau(\Phi)\ge1-\sigma$ & open for every such $\sigma$ \\
in $\Re s>\tfrac12$ & $M_{-1}(x)\ll x^{-1/2+\eps}$, $\alpha(\Phi)=\tfrac12$ & the Riemann hypothesis \\
\bottomrule
\end{tabular}
\end{center}

\rafgalleryfig{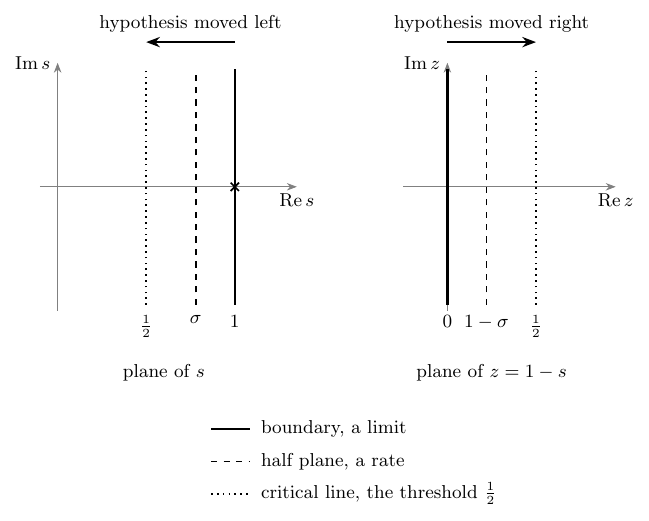}{Levels of hypothesis, not regions known to be free of zeros. In the plane of $s$, left, the three vertical lines are the boundary $\Re s=1$, the line $\Re s=\sigma$ bounding a half plane and the critical line $\Re s=\tfrac12$. In the plane of the transform, right, they are the axis $\Re z=0$, the line $\Re z=1-\sigma$ and the line $\Re z=\tfrac12$, the two planes being exchanged by $z=1-s$. A condition on the boundary returns a limit, a zero-free half plane returns a rate, and at the critical line the threshold one half is equivalent to the Riemann hypothesis. Moving the hypothesis left in $s$ moves it right in $z$. The cross marks the pole of $\zeta$ at $s=1$. No strip of positive width is known to be free of zeros, and none is shaded.}{fig:boundary_planes}

The passage from the first level to the second is not a continuity argument. Non-vanishing at each
point of the line $\Re s=1$ gives no strip of positive width, the classical zero-free regions
narrowing toward that line as the height grows, as the remark following
Corollary~\ref{cor:zerofree} records. A power rate $x^{-c}$ needs a quantitative input, a
zero-free half plane $\Re s>1-c$ and the contour estimates that go with it, the classical
zero-free regions giving only rates slower than every power, and the boundary theorem supplies
neither. In the same
way Corollary~\ref{cor:boundary_absorption} gives $A(n)\to0$ for the Ingham function at every
$\beta>0$ and says nothing about the transparency threshold, which asks for the rate $n^{-\beta}$
itself.

The chapter has shown what a weight returns at the edge of its transform. The return is a limit,
it is decided by the zeros of the transform on one line, and for the Ingham function it is the
prime number theorem. The next chapter lets the forcing exponent move and asks at which exponent
the partial sums stop following it. That is a question about a rate, and the statement that its
answer is one half is equivalent to the Riemann hypothesis.

\chapter{The equivalence with the Riemann hypothesis}
\label{chap:equivalence}

The regularity index\index[terms]{regularity index} of Chapter~\ref{chap:fgv} is defined by a discrete equation, and for the
weights met so far it is elementary to compute. For the Ingham weight it is not, and this
chapter proves that its value is equivalent to the Riemann hypothesis\index[terms]{Riemann hypothesis}. Chapter~\ref{chap:boundary}
obtained a limit at every non-negative forcing exponent from the boundary of the transform, with
no hypothesis. This chapter asks for the rate, and the rate depends on how far into the critical
strip the zeros of $\zeta$ are kept away. A half plane $\Re s>\sigma$ free of zeros, with
$\tfrac12<\sigma<1$, gives transparency below $1-\sigma$ and absorption at the rate $1-\sigma$,
Corollary~\ref{cor:zerofree}. The critical depth $\sigma=\tfrac12$ is the equivalence itself.

The equivalence has two halves and they are not symmetric. Under the hypothesis, transparency\index[terms]{transparency}
below one half and absorption\index[terms]{absorption} at and above it follow from the transfer law\index[terms]{transfer law} of
\S\ref{sec:master_equiv}, whose analytic input is the bound of Littlewood\index[names]{Littlewood, J. E.} on the summatory
M\"obius function\index[terms]{M\"obius function}. That the index cannot exceed one half needs no hypothesis at all. It rests on the existence, for
each fixed exponent in the relevant range, of a simple zero of $\zeta$ on the critical line that
is not a zero of the Dirichlet series attached to the forcing.
Proposition~\ref{prop:C_noncommon_zero} produces such a zero by comparing two counts, a Jensen
bound on the zeros of that series against the positive density of simple critical zeros of
$\zeta$. No zero is located and none is evaluated. In the other direction the exponent that carries the argument is $\beta=1$, where
an identity of Meissel\index[names]{Meissel, E.} makes the solution of the defining equation\index[terms]{defining equation} the M\"obius function divided
by the coordinate, so that the index bound becomes the classical bound on its summatory
function.

The mechanism is worth separating from the statement. The zeros of $\zeta$ are needed on both
sides and they are analytic objects, so nothing here turns the hypothesis into something other
than what it is. What the equivalence displays is where the passage is made, and it is made
through the divisor structure of the weight and through M\"obius inversion\index[terms]{M\"obius inversion}, which are
arithmetic.

The chapter opens on the case of negative exponents, where a proof is available that uses no
complex analysis and gives a first check of the theory before the general one. Three sections
follow the equivalence. The Hardy-Littlewood-Ramanujan\index[terms]{Hardy--Littlewood--Ramanujan criterion}\index[names]{Ramanujan, S.} criterion,
defined in Chapter~\ref{chap:boundary} where the Ingham weight was shown to satisfy it in its
strong form, is followed to the generalized Ingham functions, where its two forms separate. A
character twisted weight then carries the
same equivalence for two zero sets at once, and that one is proved here as well. A further section turns the same machinery on counting problems,
where the weight is read against an arithmetic count rather than against a zero set. The chapter
closes on the intermediate depths, where a half plane free of zeros\index[terms]{zero-free region} gives a lower bound on the
transparency threshold, and on
what a positive index would imply, the prime number theorem\index[terms]{prime number theorem} among other things.

The equivalence proved in this chapter has a well known neighbour, and the two are best set side
by side before either is used. The criterion of
\nm{Nyman}{B.}~\cite{Nyman1950} and \nm{Beurling}{A.}~\cite{Beurling1955} states that the Riemann
hypothesis holds if and only if the constant function belongs to the closed linear span, in
$L^{2}(0,1)$, of the functions $x\mapsto\{\theta/x\}-\theta\{1/x\}$ for $0<\theta\le1$, and the
refinement of \nm{B\'aez-Duarte}{L.}~\cite{BaezDuarte2003} restricts the dilations $\theta$ to the
reciprocals of the integers. Both belong to the closure problems collected, with their neighbours,
in \nm{Broughan}{K.}~\cite[Ch.~3 and Ch.~8]{Broughan2017II}.

The kinship with the present work is real, and it is read on the transform. The Nyman\index[names]{Nyman, B.} criterion
rests on the classical Mellin transform of the fractional part\index[terms]{fractional part},
\[
\int_0^{\infty}\Big\{\frac1t\Big\}\,t^{z-1}\,dt=-\frac{\zeta(z)}{z},
\qquad 0<\Re z<1,
\]
whereas the theory developed here rests on the Mellin transform\index[terms]{Mellin transform} of the Ingham
kernel\index[terms]{Ingham kernel},
\[
\Phi^{*}(z)=-z\int_0^{1}\Phi(t)\,t^{-z-1}\,dt=\frac{z}{z-1}\,\zeta(1-z),
\qquad \Re z<0,
\]
continued to the plane. The two kernels are one object seen twice, since
$\Phi(x)=x\lfloor1/x\rfloor=1-x\{1/x\}$, and the functional equation of $\zeta$ makes the two
transforms carry the same information on the location of the nontrivial zeros\index[terms]{nontrivial zero} in the critical
strip\index[terms]{critical strip}.

What separates the closure criteria from the construction of this volume is the question put to
the transform. A closure criterion asks for a distance in a Hilbert space to vanish. The
construction here asks at which exponent the response of a summation equation changes nature, and
the arithmetic enters through the M\"obius inversion\index[terms]{M\"obius inversion} that the
floor function carries, which is what places Perron's formula\index[terms]{Perron's formula} and
the divisor estimates of \S\ref{sec:master_equiv} on the direction from the hypothesis to the
index.
That direction is the
one \nm{Levinson}{N.}~\cite{Levinson1956} found to be of trivial character in the closure setting,
and his reservation about the prospects of closure theorems was addressed to it.

The two constructions descend from two theorems of Wiener\index[names]{Wiener, N.} rather than from
one, and the distance between them is the distance between those two theorems. The Tauberian
theorem for $L^{1}$, which concludes when the transform has no zero on the line where the problem
is posed, is the one that carries Ingham's\index[names]{Ingham, A. E.} method and the transfer used
in this volume, and Chapter~\ref{chap:boundary} writes it in the coordinates of the defining
equation. The closure theorem for $L^{2}$, which asks only that the transform be nonzero
almost everywhere for the translates of a function to span the space, is the one that carries the
criterion of Nyman\index[names]{Nyman, B.} and Beurling\index[names]{Beurling, A.}. Both are in \cite{Wiener1932}. What separates the two families is
therefore not their subject but which of the two theorems is in use, and with it how much a zero
of a transform is allowed to decide.

\section{A constructive proof for $\beta<0$}\label{sec:beta_negative}

The transfer law of the next section rests at $\beta=0$ on Ingham's theorem, proved in
Chapter~\ref{chap:boundary}, and for $\beta>0$, under the Riemann hypothesis, on a contour argument. For $\beta<0$ the same conclusion follows with no
complex analysis at all, and the computation is worth seeing before the general one. It shows the
whole mechanism at an exponent where nothing is conditional, the M\"obius inversion that carries
the arithmetic, the exact constant $1/\Phi^{*}(\beta)$, and the Abel summation that returns the
partial sums.

\begin{lemma}\label{lem:ingham_beta_neg}
Let $\beta < 0$ and $\delta < -\beta$. Suppose the sequence $(a_n)$ satisfies
\[
\sum_{k=1}^n a_k \Phi\!\left(\frac{k}{n}\right) = n^{-\beta} + \mathcal{O}(n^\delta).
\]
Then $A(x) := \sum_{n \le x} a_n$ satisfies
\[
A(x) \sim \frac{1}{\Phi^*(\beta)}\, x^{-\beta} \qquad (x \to \infty).
\]
\end{lemma}

\begin{proof}
Write $b_k=ka_k$ and $T(N)=\sum_{k\le N}b_k\lfloor N/k\rfloor$. Multiplying the hypothesis
by $N$ gives $T(N)=N^{1-\beta}+\mathcal{O}(N^{1+\delta})$. Since $T(N)=\sum_{m\le N}(b\star\mathbf 1)(m)$ by the interchange of summation,
this gives
$\sum_{j\le N}\mu(j)T(\lfloor N/j\rfloor)=\sum_{jm\le N}\mu(j)(b\star\mathbf 1)(m)
=\sum_{t\le N}(\mu\star b\star\mathbf 1)(t)=\sum_{t\le N}b_t$, that is
\[
\sum_{k\le N}ka_k=\sum_{j\le N}\mu(j)\,T\Big(\Big\lfloor\frac Nj\Big\rfloor\Big),
\]
and $1-\beta>1$ makes $\sum_{j\ge1}\mu(j)\,j^{-(1-\beta)}=1/\zeta(1-\beta)$ absolutely
convergent. Inserting, with $u=N/j\ge1$,
\[
 T(\lfloor u\rfloor)=u^{1-\beta}
 +\mathcal{O}_\beta\big(u^{-\beta}\big)
 +\mathcal{O}\big(u^{1+\delta}\big),
\]
and summing, the floor error contributes $\mathcal{O}(N\log N+N^{-\beta})$ and the remainder $\mathcal{O}((N^{1+\delta}+N)\log N)$, all $o(N^{1-\beta})$ since $1+\delta<1-\beta$, so
\[
\sum_{k\le N}ka_k=\frac{N^{1-\beta}}{\zeta(1-\beta)}+o(N^{1-\beta}).
\]
Abel summation then gives
$A(x)\sim\dfrac{\beta-1}{\beta\,\zeta(1-\beta)}\,x^{-\beta}=\dfrac{1}{\Phi^*(\beta)}\,x^{-\beta}$.
\end{proof}

\begin{remark}
The computation above establishes case (I) of the transfer law, Theorem~\ref{thm:transfer_law}
below, without any appeal to the Riemann hypothesis or to complex contour integration. This exact
M\"obius inversion is special to the Ingham kernel. What it exhibits, and what the gallery of
appendices meets again, is the mechanism of initial transparency, the discrete equation
reproducing the imposed power with the Mellin coefficient in the domain of absolute convergence.
Later kernels obey analogous laws, and their proofs may require entirely different tools.
\end{remark}

\section{The master equivalence}
\label{sec:master_equiv}

The Ingham function $\Phi(x):=x\lfloor x^{-1}\rfloor$ has Mellin transform\index[terms]{Mellin transform}
$\Phi^*(z)=\frac{z}{z-1}\zeta(1-z)$. The equivalence established below joins a long list,
collected with its sources in \nm{Borwein}{P.}~\cite{Borwein2009}. Throughout the section the averages are written
\begin{equation}\label{eq:strong_input}
A_\Phi(n):=\frac1n\sum_{k=1}^n k\,a_k\Big\lfloor\frac nk\Big\rfloor
=n^{-\beta}+r(n),
\end{equation}
so that the defining equation corresponds to $r\equiv0$. The admissible size of the
error $r(n)$ depends on $\beta$. The next theorem gives, for each case, an input
condition calibrated to the level at which the M\"obius inversion is performed, and
the remark following the proof explains why the tolerance changes at $\beta=0$.

Two elementary facts are used repeatedly in this section. First, for $t>1$ the sum
of $\mu$ over the divisors of $t$ runs over the squarefree divisors and equals
$(1-1)^{\omega(t)}=0$, where $\omega(t)$ counts the distinct prime factors, while it
equals $1$ for $t=1$. In convolution form $\mu\star\mathbf 1=e$, the unit of
Dirichlet convolution, so that $f=g\star\mathbf 1$ is equivalent to $g=\mu\star f$,
and, multiplying the absolutely convergent series for $\Re s>1$,
\[
\zeta(s)\sum_{n\ge1}\frac{\mu(n)}{n^{s}}
=\sum_{t\ge1}\frac{(\mathbf 1\star\mu)(t)}{t^{s}}=1,
\qquad\text{hence}\qquad
\sum_{n\ge1}\frac{\mu(n)}{n^{s}}=\frac{1}{\zeta(s)}.
\]
Second, for arithmetic functions $u$ and $v$,
\[
\sum_{n\le x}(u\star v)(n)=\sum_{jk\le x}u(k)\,v(j)
=\sum_{k\le x}u(k)\sum_{j\le x/k}v(j),
\]
a plain interchange of the order of summation. The following bound, which goes back
to Littlewood \cite{Littlewood1912}, is the analytic input for the values $\beta>0$.

\begin{lemma}\label{lem:littlewood}
The Riemann hypothesis is equivalent to the bound
$M(x):=\sum_{n\le x}\mu(n)\ll_\eps x^{1/2+\eps}$ for every $\eps>0$.
\end{lemma}

\begin{proof}
This is Theorem~14.25(C) of \nm{Titchmarsh}{E. C.}~\cite[Theorem~14.25(C)]{Titchmarsh1986}.
For completeness, the reverse implication is short. If the displayed bound holds, partial summation
makes $\sum_{n\ge1}\mu(n)n^{-s}$ locally uniformly convergent and holomorphic on
$\Re s>\tfrac12$. It agrees with $1/\zeta(s)$ on $\Re s>1$ and hence throughout that connected
half plane, so $\zeta$ has no zero there. The functional equation and the symmetry of the
nontrivial zeros then put every such zero on the critical line. For the direct implication, Corollary~\ref{cor:perron-mobius-halfplane} with
$\theta=1/2$ supplies the complete finite contour argument. The vertical bounds it uses are
proved in Lemma~\ref{lem:perron-zerofree-growth}. These arguments recover the cited
classical result without leaving its Perron step implicit.
\end{proof}

The transfer law is the analytic core of the equivalence, and it separates the forcing exponents
into cases.

\begin{theorem}\label{thm:transfer_law}
Let $A(x)=\sum_{n\le x}a_n$ and let $(a_n)$ satisfy \eqref{eq:strong_input}.
\begin{enumerate}
\item[(I)] If $\beta<0$ and $r(n)=\mathcal{O}(n^{\delta})$ for some $\delta<-\beta$, then
\[
A(x)\sim\frac{1}{\Phi^*(\beta)}\,x^{-\beta}\qquad(x\to\infty).
\]
\item[(II)] If $\beta=0$, $r(n)=o(1)$ and $na_n\ge-C$ for some constant $C>0$, then
\[
A(x)\longrightarrow\frac{1}{\Phi^*(0)}=1\qquad(x\to\infty).
\]
\item[(III)] If $0<\beta<\tfrac12$ and $r(n)=\mathcal{O}(n^{-\beta-1-\eta})$ for some
$\eta>0$, then under the Riemann hypothesis
\[
A(x)\sim\frac{1}{\Phi^*(\beta)}\,x^{-\beta}\qquad(x\to\infty).
\]
\item[(IV)] If $\beta\ge\tfrac12$ and $r(n)=\mathcal{O}(n^{-\beta-1-\eta})$ for some
$\eta>0$, then under the Riemann hypothesis
\[
A(x)\ll_\eps x^{-1/2+\eps}.
\]
\end{enumerate}
\end{theorem}
\begin{proofstatus}{Cases (I), (III), and (IV) are established by the arguments below and feed
Theorem~\ref{thm:tauberian_rh}. Case (II) is Ingham's theorem, proved in
Corollary~\ref{cor:ingham_proved}. It is included for completeness and is used in no later
proof. The central equivalence instead treats $\beta=0$ through the explicit solution
$a=(1,0,0,\ldots)$. The case $r\equiv0$ of the coefficient decomposition below reproduces the proof of
Proposition~\ref{prop:hlr_ingham}, and parts of it are reused in
Proposition~\ref{prop:gen_unbounded_leading}, Theorem~\ref{thm:gen_equiv}, and
Corollary~\ref{cor:zerofree}.}
\end{proofstatus}
\begin{proof}
Write $A_1(x)=\sum_{n\le x}na_n$ throughout.

(I). This case is Lemma~\ref{lem:ingham_beta_neg}, proved in the preceding section by a
self-contained argument that uses no complex analysis and no information on the zeros
of the zeta function.

(II). This is Ingham's theorem\index[terms]{Ingham theorem} \cite{Ingham1945}, recalled in \S\ref{sec:ingham_theorem} and proved in
Corollary~\ref{cor:ingham_proved}. At $\beta=0$ the input
\eqref{eq:strong_input} reads $A_\Phi(n)=1+o(1)$, the one-sided bound
$na_n\ge-C$ is the Tauberian hypothesis of that theorem, and the conclusion is the
convergence $\sum a_k=1$. No rate on $r$ is required. The value agrees with
$\Phi^*(0)^{-1}$ since $\Phi^*(z)\to1$ as $z\to0$.

(III). Differencing \eqref{eq:strong_input} at $n$ and $n-1$ gives, for $n\ge2$,
\[
\sum_{d\mid n}d\,a_d=c(n)+w(n),\qquad
c(n)=n^{1-\beta}-(n-1)^{1-\beta},\qquad
w(n)=n\,r(n)-(n-1)\,r(n-1),
\]
with $w(n)=\mathcal{O}(n^{-\beta-\eta})$. The two sequences are set separately at $n=1$, where the
displayed formulas do not apply, the term $(n-1)^{1-\beta}$ being undefined at $\beta\ge1$ and
$r(0)$ carrying no meaning. Taking \eqref{eq:strong_input} itself at $n=1$ gives $a_1=1+r(1)$, so
\[
c(1)=1,\qquad w(1)=r(1),
\]
which is the pair that keeps $\sum_{d\mid 1}d\,a_d=c(1)+w(1)$ true. The left side being the
convolution
of $n\mapsto na_n$ with $\mathbf 1$, M\"obius inversion turns this into
$na_n=(\mu\star c)(n)+(\mu\star w)(n)$. The expansion
$c(n)=(1-\beta)n^{-\beta}+u(n)$ with $u(n)=\mathcal{O}(n^{-\beta-1})$, valid for $n\ge2$ and
completed by $u(1)=\beta$, then splits the
first term, so that
\[
na_n=(1-\beta)\,J_{-\beta}(n)+(\mu\star u)(n)+(\mu\star w)(n),
\qquad J_{-\beta}(n):=\sum_{d\mid n}d^{-\beta}\mu\Big(\frac nd\Big),
\]
and the three partial sums are estimated in turn.

Since $J_{-\beta}$ is the convolution of $n\mapsto n^{-\beta}$ with $\mu$,
multiplying the absolutely convergent series for $\Re s>1$ gives
$\sum_{n\ge1}J_{-\beta}(n)n^{-s}=\zeta(s+\beta)/\zeta(s)$, the right side continuing
the series meromorphically. The divisors of $n$ factor uniquely along its prime
powers and $\mu$ is multiplicative, so that
\[
J_{-\beta}(n)=\prod_{p^{a}\|n}\big(p^{-a\beta}-p^{-(a-1)\beta}\big),
\]
every factor lying in $[-1,1]$, hence $|J_{-\beta}(n)|\le1$ for every $n$.
This is the coefficient bound needed for finite Perron inversion.
The full proof of the inversion and of the three contour estimates is in
Appendix~\ref{app:perron}, so only their arithmetic data are fixed here.
For a series $F(s)=\sum_{n\ge1}f(n)n^{-s}$ with $|f(n)|\le A$, set
$y=N+1/2\ge5/2$ and $c=5/4$. Lemma~\ref{lem:perron-half-integer},
with $\kappa=1$, gives for $T\ge2$
\begin{equation}\label{eq:perron_effective}
 \sum_{n\le y}f(n)=\frac1{2\pi i}\int_{c-iT}^{c+iT}
                  F(s)\frac{y^s}{s}\,ds
       +\mathcal O\!\left(A\frac{y^{5/4}+y\log y}{T}\right).
\end{equation}
The half integer is essential for this particular error bound, since it
keeps every $|y-n|$ at least $1/2$.

Apply the formula to $F(s)=\zeta(s+\beta)/\zeta(s)$. Fix
$0<\eps<\min(1/4,1/2-\beta,\eta)$ and choose
\[
 \delta=\min\left(\frac\eps4,\frac{1/2-\beta}{2}\right),\qquad
 b=\frac12+\delta<1-\beta,\qquad T=y^2,\qquad \nu=\frac\eps8.
\]
Under RH, Lemma~\ref{lem:perron-zerofree-growth} with $\theta=1/2$
applies to $\zeta(s)$ and $\zeta(s+\beta)$ on this strip. Applying it
with exponent $\nu/2$ to each factor gives
$F(\sigma+it)\ll_{\beta,\eps}(1+|t|)^\nu$ at large heights.
The pole of $\zeta(s)$ at $1$ gives a removable zero of its reciprocal.
Thus the only pole of $F(s)/s$ between $b$ and $c$ is $s=1-\beta$,
with residue
\[
 \Res_{s=1-\beta}\left(F(s)\frac{y^s}{s}\right)
       =\frac{y^{1-\beta}}{(1-\beta)\zeta(1-\beta)}.
\]
The estimates proved in Proposition~\ref{prop:perron-rectangle} now give
\[
 E_{\rm vert}\ll y^{1/2+\delta+2\nu},\qquad
 E_{\rm hor}\ll\frac{y^{-3/4+2\nu}}{\log y},\qquad
 E_{\rm trunc}\ll y^{-3/4}+y^{-1}\log y.
\]
Here $\delta+2\nu\le\eps/2$, so every error is
$\mathcal O_{\beta,\eps}(y^{1/2+\eps})$. All gaps and exponents are
fixed before $y$ tends to infinity. Consequently
\[
 (1-\beta)\sum_{n\le y}J_{-\beta}(n)
 =\frac{y^{1-\beta}}{\zeta(1-\beta)}
       +\mathcal O_{\beta,\eps}(y^{1/2+\eps}).
\]
For any real $x\ge2$ take $y=\lfloor x\rfloor+1/2$.
The partial sums agree and
$y^{1-\beta}-x^{1-\beta}=\mathcal O_\beta(x^{-\beta})$, so the
same formula holds with $x$ in place of $y$.
The two remaining sums are handled by the
identity $\sum_{n\le x}(f\star\mu)(n)=\sum_{k\le x}f(k)M(x/k)$ together with
$M(y)\ll y^{1/2+\eps}$, supplied under the hypothesis by Lemma~\ref{lem:littlewood}, and by
the elementary partial sum estimate
\begin{equation}\label{eq:harmonic_two_terms}
\sum_{k\le x}k^{-s}\ \ll_{s,\eps}\ 1+x^{1-s+\eps}\qquad(s\in\R,\ x\ge2),
\end{equation}
whose two terms record the two regimes of the series, convergence when $s>1$ and growth of
order $x^{1-s}$ when $s<1$, the $\eps$ covering the logarithm at $s=1$. For the second sum,
$u(k)=\mathcal{O}(k^{-\beta-1})$ gives
\[
\sum_{n\le x}(u\star\mu)(n)\ll x^{1/2+\eps}\sum_{k\le x}k^{-\beta-3/2-\eps}
\ll x^{1/2+\eps},
\]
where only the first term of \eqref{eq:harmonic_two_terms} survives, the exponent
$\beta+\tfrac32$ exceeding $1$ throughout the present case. For the third,
$w(k)=\mathcal{O}(k^{-\beta-\eta})$ produces the exponent $\beta+\eta+\tfrac12$, which can
lie on either side of $1$ because $\eta$ is only known to be positive, so both terms of
\eqref{eq:harmonic_two_terms} are kept and the factor $x^{1/2+\eps}$ multiplies each of
them,
\[
\sum_{n\le x}(w\star\mu)(n)\ll x^{1/2+\eps}\sum_{k\le x}k^{-\beta-\eta-1/2-\eps}
\ll x^{1/2+\eps}\big(1+x^{1/2-\beta-\eta}\big)
\ll x^{1/2+\eps}+x^{1-\beta-\eta+\eps}.
\]
Collecting the three estimates gives
\[
A_1(x)=\frac{x^{1-\beta}}{\zeta(1-\beta)}
+\mathcal{O}\big(x^{1/2+\eps}+x^{1-\beta-\eta+\eps}\big),
\]
and Abel summation\index[terms]{Abel summation} yields
\[
A(x)=\frac{A_1(x)}{x}+\int_1^x\frac{A_1(t)}{t^2}\,dt
=\Big(1-\frac1\beta\Big)\frac{x^{-\beta}}{\zeta(1-\beta)}+\kappa
+\mathcal{O}\big(x^{-1/2+\eps}+x^{-\beta-\eta+\eps}\big)
\]
for some constant $\kappa$, the integral of the error term over $[1,\infty)$
converging absolutely. Finally the constant vanishes. The last display shows $A(x)\to\kappa$. Termwise
$|c(d)|+|w(d)|\ll d^{-\beta}\le1$, so $|na_n|\le\sum_{d\mid n}(|c(d)|+|w(d)|)\ll d(n)$,
and $\sum_{n\ge1}d(n)n^{-1-s}=\zeta(1+s)^{2}$ for $s>0$, by squaring the series of
$\zeta(1+s)$, so the identity $\sum_{d\mid n}da_d=c(n)+w(n)$ multiplies out to
\[
\zeta(1+s)\sum_{n\ge1}\frac{a_n}{n^{s}}=\sum_{t\ge1}\frac{c(t)+w(t)}{t^{1+s}}
\qquad(s>0),
\]
all series converging absolutely. On the other hand partial summation\index[terms]{partial summation} gives
$\sum_{n\ge1}a_nn^{-s}=s\int_1^\infty A(t)t^{-s-1}\,dt$ for $s>0$, the function $A$
being bounded, and, since $s\int_1^\infty t^{-s-1}\,dt=1$,
\[
\Big|\sum_{n\ge1}\frac{a_n}{n^{s}}-\kappa\Big|
\le s\int_1^\infty\big|A(t)-\kappa\big|\,t^{-s-1}\,dt\longrightarrow0
\qquad(s\to0^{+}),
\]
the integral over $[1,T]$ being $\mathcal{O}(s\log T)$ and the tail at most
$\sup_{t\ge T}|A(t)-\kappa|$. Hence
\[
\kappa=\lim_{s\to0^{+}}\sum_{n\ge1}\frac{a_n}{n^{s}}
=\lim_{s\to0^{+}}\frac{1}{\zeta(1+s)}\sum_{t\ge1}\frac{c(t)+w(t)}{t^{1+s}}=0,
\]
the sum on the right tending to the finite value $\sum_t(c(t)+w(t))t^{-1}$ while
$\zeta(1+s)\ge\int_1^\infty t^{-1-s}\,dt=1/s$. Therefore
\[
A(x)=\Big(1-\frac1\beta\Big)\frac{x^{-\beta}}{\zeta(1-\beta)}
+\mathcal{O}\big(x^{-1/2+\eps}+x^{-\beta-\eta+\eps}\big)
=\frac{1}{\Phi^*(\beta)}\,x^{-\beta}+o\big(x^{-\beta}\big),
\]
the pole of the zeta function at $1$ absorbing the constant term.

(IV). Take $0<\eps<1/2$. The same decomposition applies, and the swap identity
now bounds all three sums at once. Bounds for larger $\eps$ follow by weakening
one of these estimates. Hence
\[
\sum_{n\le x}J_{-\beta}(n)=\sum_{k\le x}k^{-\beta}M\Big(\frac xk\Big)
\ll x^{1/2+\eps}\sum_{k\le x}k^{-\beta-1/2-\eps}\ll x^{1/2+\eps},
\]
the series converging since $\beta\ge\tfrac12$ and the bound on $M$ coming again
from Lemma~\ref{lem:littlewood}. The $u$ and $w$ parts obey the same bound, so that
$A_1(x)\ll x^{1/2+\eps}$. Abel summation and the vanishing of
the constant, exactly as in (III), give $A(x)\ll x^{-1/2+\eps}$.
\end{proof}

\begin{remark}\label{rem:input_tolerance}
The tolerances of Theorem~\ref{thm:transfer_law} are calibrated to the level at which
the M\"obius inversion operates. For $\beta<0$ the inversion acts on the partial sums
$\sum_{k\le N}ka_k\lfloor N/k\rfloor$, whose main term is $N^{1-\beta}$, and any error
a polynomial margin below it is absorbed. For $\beta>0$ the argument descends to the
divisor level, where the differenced input $\sum_{d\mid n}da_d$ carries the main term
$c(n)\asymp n^{-\beta}$, and the differenced error must in turn stay a polynomial
margin below, which is the condition $r(n)=\mathcal{O}(n^{-\beta-1-\eta})$. The point
$\beta=0$ is the only one where no rate is required, because there the inversion
reduces to the exact identity $\sum_{j\le N}\mu(j)\lfloor N/j\rfloor=1$,
established in the proof of Theorem~\ref{thm:tauberian_rh}, and
the polynomial margin is replaced by the one-sided Tauberian condition of Ingham's\index[names]{Ingham, A. E.}
theorem.

Two conditions on the individual terms enter, and they are not the same condition. At $\beta=0$
the hypothesis is the one-sided Tauberian condition $na_n\ge-C$ of Hardy\index[names]{Hardy, G. H.}
and Littlewood\index[names]{Littlewood, J. E.}, which is what Ingham's theorem asks. The HLR
criterion\index[terms]{Hardy--Littlewood--Ramanujan criterion} of Definition~\ref{def:HLR} is
two-sided, and its strong form $na_n=\mathcal{O}(1)$ implies the one-sided condition without being
implied by it. For $\beta>0$ the argument uses neither as a hypothesis. What it uses is the
multiplicative
mechanism behind the criterion, the product form of $J_{-\beta}(n)$ displayed in the proof, which bounds
the coefficients and validates the truncation of the Perron formula, so the fundamental theorem of
arithmetic enters the argument exactly there. No arithmetic condition on $(a_n)$ is
needed for $\beta<0$. The polynomial margin for $\beta>0$ is not an artifact of the
method. Under the rateless input $\ell+o(1)$ the best possible individual bound is
$a_n=o(\log\log n)$ \cite{ErdosSegal1978}, so already the bound
$na_n=\mathcal{O}(1)$ of the strong HLR criterion, established for the defining equation in
Proposition~\ref{prop:hlr_ingham}, requires the defining equation or a polynomially small error.
\end{remark}

The bound rests on a single nontrivial zero, entering through a series attached to the forcing
exponent.

\begin{lemma}
\label{lem:phi_sharpness}
Let $\tfrac12<\beta_0<1$ and put
\begin{equation}\label{eq:C_beta_series}
C_{\beta_0}(s)=\sum_{t\ge1}\bigl(t^{1-\beta_0}-(t-1)^{1-\beta_0}\bigr)t^{-s},
\end{equation}
a series that converges absolutely on $\Re s>1-\beta_0$. If $C_{\beta_0}(\rho)\ne0$ for
at least one zero $\rho$ of $\zeta$ with $\Re\rho>1-\beta_0$, then the exponent
$\beta_0$ is not transparent for the Ingham function.
\end{lemma}

\begin{proof}
Write $b_n=na_n$ for the solution at the exponent $\beta_0$. Multiplying the defining
equation\index[terms]{defining equation} by $n$ and counting the integers $i$ with $ik\le n$ gives
$\sum_{k\le n}b_k\lfloor n/k\rfloor=n^{1-\beta_0}$, that is
$\sum_{m\le n}(b\star\mathbf 1)(m)=n^{1-\beta_0}$, so
\[
(b\star\mathbf 1)(n)=n^{1-\beta_0}-(n-1)^{1-\beta_0}=:c_{\beta_0}(n),
\qquad b=\mu\star c_{\beta_0}
\]
by M\"obius inversion. Since $c_{\beta_0}(t)\ll t^{-\beta_0}$, the series
\eqref{eq:C_beta_series} converges absolutely on $\Re s>1-\beta_0$, and on $\Re s>1$ the
convolution identity reads
\begin{equation}\label{eq:D_over_zeta}
D(s):=\sum_{n\ge1}b_n n^{-s}=\frac{C_{\beta_0}(s)}{\zeta(s)} .
\end{equation}

Suppose $\beta_0$ transparent, so that $A(n)=\Xi n^{-\beta_0}+o(n^{-\beta_0})$ for some
constant $\Xi$, the value $\Xi=0$ included. Partial summation gives
\[
\sum_{n\le x}b_n=xA(x)-\sum_{n<x}A(n)
=\Xi\Bigl(1-\frac1{1-\beta_0}\Bigr)x^{1-\beta_0}+o\bigl(x^{1-\beta_0}\bigr)+\mathcal O(1),
\]
so the summatory function of $(b_n)$ is $\mathcal O(x^{1-\beta_0})$ and $D$ converges
and is holomorphic on $\Re s>1-\beta_0$. Both sides of \eqref{eq:D_over_zeta} are
meromorphic there and agree on $\Re s>1$, hence they agree on the whole half plane. At a
zero $\rho$ of $\zeta$ with $\Re\rho>1-\beta_0$ the quotient has a pole unless
$C_{\beta_0}$ vanishes at $\rho$ to at least the order of that zero, while the left side
is holomorphic. Hence $C_{\beta_0}(\rho)=0$, against the hypothesis, and $\beta_0$ is
not transparent.

\end{proof}

What the lemma needs from outside is one zero of $\zeta$ that is not a zero of the corresponding
Dirichlet series. Such a zero exists for every fixed exponent in the required interval without
selecting or evaluating any zero numerically.

\begin{proposition}\label{prop:C_noncommon_zero}
For every fixed $\beta\in(\tfrac12,1)$ there is a simple zero
$\rho=\tfrac12+i\gamma$, $\gamma>0$, of $\zeta$ such that $C_\beta(\rho)\neq0$.
\end{proposition}

\begin{proof}
Put $a=1-\beta$ and
$c_\beta(n)=n^{1-\beta}-(n-1)^{1-\beta}$. Then $c_\beta(1)=1$ and
$c_\beta(n)=\mathcal O_\beta(n^{-\beta})$. Hence the series
\[
 C_\beta(s)=\sum_{n\ge1}c_\beta(n)n^{-s}
\]
is holomorphic on $\Re s>a$.
Choose $S>\tfrac12$ so large that
$\sum_{n\ge2}c_\beta(n)n^{-S}\le\tfrac12$. Uniformly in real $t$,
\[
 |C_\beta(S+it)|\ge\tfrac12.
\]
Write $d=S-\tfrac12$ and $g=\beta-\tfrac12$, and choose
\[
 r=d+\frac g3,\qquad R=d+\frac{2g}{3},\qquad R_*=d+\frac{5g}{6}.
\]
Since $S-R_*=a+g/6>a$, every closed disk of centre $S+it$ and radius $R_*$ lies in
$\Re s>a$. On each such disk
\[
 |C_\beta(s)|\le M_\beta:=\sum_{n\ge1}c_\beta(n)n^{-(S-R_*)}<\infty,
\]
independently of $t$. Choose an outer radius $R_t\in(R,R_*)$ whose circle contains no zero.
For clarity, Jensen's calculation is as follows. Put $f_t(w)=C_\beta(S+it+w)$,
and list its zeros $w_\ell$ in $|w|<R_t$, with multiplicity. There are finitely many,
none is zero, and factoring them out leaves a function holomorphic and nonzero on
$|w|\le R_t$. The mean-value identity for its logarithmic modulus, together with
\[
 \frac1{2\pi}\int_0^{2\pi}\log|R_te^{i\theta}-w_\ell|\,d\theta=\log R_t,
\]
gives
\[
 \sum_{|w_\ell|<R_t}\log\frac{R_t}{|w_\ell|}
 =\frac1{2\pi}\int_0^{2\pi}\log|f_t(R_te^{i\theta})|\,d\theta
  -\log|f_t(0)|
 \le\log(2M_\beta).
\]
The displayed mean for a linear factor follows by expanding
$\log(1-(w_\ell/R_t)e^{-i\theta})$ in its uniformly convergent power series and
integrating each nonconstant Fourier term. Every zero with $|w_\ell|\le r$
contributes at least $\log(R/r)$. Thus the disk of radius $r$ contains at most
\[
 K_\beta=\frac{\log(2M_\beta)}{\log(R/r)}
\]
zeros, counted with multiplicity. Its intersection with the critical line has
half-length $h_\beta=(r^2-d^2)^{1/2}>0$. Centres $S+i\ell h_\beta$, for
$0\le\ell\le\lceil T/h_\beta\rceil$, cover the critical segment up to height $T$.
Consequently its number of zeros of $C_\beta$, counted with multiplicity, is at most
$K_\beta(2+T/h_\beta)=\mathcal O_\beta(1+T)$.

The positive proportion of zeros on the critical line goes back to
\nm{Selberg}{A.}~\cite{Selberg1942}, and \nm{Levinson}{N.}~\cite{Levinson1974} carried it past one
third by a method that already yields a positive proportion of simple zeros. What is used here is
the quantitative form of Conrey, which is what makes the count below explicit.

Let $N(T)$ count the nontrivial zeros of $\zeta$ with $0<\Im\rho<T$, with
multiplicity, and let $N_0^*(T)$ count those which are simple and on the critical
line. These are the quantities in \cite[p.~3, (8), (11)--(12)]{Conrey1989}.\index[names]{Conrey, J. B.}
The unconditional Theorem~1 on page~4 of that paper gives
\[
 \liminf_{T\to\infty}\frac{N_0^*(T)}{N(T)}\ge0.401.
\]
The zero-counting formula recorded in \cite[p.~3, (8)]{Conrey1989} is
\[
 N(T)=\frac{T}{2\pi}\log\frac{T}{2\pi e}+\mathcal O(\log T).
\]
Hence $N_0^*(T)\ge0.4N(T)\gg T\log T$ for all sufficiently large $T$.
These zeros are distinct. At most $\mathcal O_\beta(1+T)$ can be zeros of
$C_\beta$, so at least one is not shared. The argument holds for each fixed
$\beta\in(\tfrac12,1)$ and assumes no zero-free region.
\end{proof}

\begin{corollary}\label{cor:phi_half}
$\tau(\Phi)\le\tfrac12$, under no hypothesis. Consequently, if $\Phi$ is a function of good
variation then $\alpha(\Phi)=\tau(\Phi)\le\tfrac12$.
\end{corollary}

\begin{proof}
Let $c>\tfrac12$ and choose $\beta\in(\tfrac12,\min(c,1))$. Proposition~\ref{prop:C_noncommon_zero}
gives a zero $\rho$ with $\Re\rho=\tfrac12>1-\beta$ and $C_\beta(\rho)\neq0$.
Lemma~\ref{lem:phi_sharpness} shows that $\beta$ is not transparent. Thus no
$c>\tfrac12$ belongs to the set whose supremum defines $\tau(\Phi)$.
\end{proof}

The upper bound is qualitative, and it uses no list of zeros and no numerical lower bound for
$C_\beta$.

Both halves are now in place, and they were bought at very different prices. That the threshold
does not exceed one half is Corollary~\ref{cor:phi_half}, and it costs no hypothesis. For every
$c>\tfrac12$ it produces a non-transparent exponent below $c$, out of a simple zero of $\zeta$ on
the critical line that the forcing does not share, and no such zero is ever named. Transparency
below one half costs the hypothesis itself, entering through the transfer law, whose analytic
input at positive exponents is the Littlewood bound on the summatory M\"obius function. The
converse runs on the single exponent $\beta=1$, where the identity of Meissel makes the solution
$\mu(n)/n$ and absorption returns that same bound. Lemma~\ref{lem:littlewood} therefore serves
twice, its direct half as the analytic input and its converse as the last step.

\begin{theorem}\label{thm:tauberian_rh}
The Riemann hypothesis holds if and only if the Ingham function is a function of good
variation of index $\alpha(\Phi)=1/2$.
\end{theorem}

\begin{proof}
Assume the Riemann hypothesis. For the defining equation $r\equiv0$, case (I) of
Theorem~\ref{thm:transfer_law} gives transparency for $\beta<0$. At $\beta=0$ the unique
solution is $a=(1,0,0,\ldots)$, since $\Phi(1/n)=1$, and hence $A(x)=1$ directly.
Case (III) gives transparency for $0<\beta<\tfrac12$. Every exponent below $\tfrac12$ is
therefore transparent, which is $\tau(\Phi)\ge\tfrac12$, and
Corollary~\ref{cor:phi_half} gives the reverse inequality, so the threshold is finite and
equal to $\tfrac12$. Case (IV) gives $A(x)\ll_\eps x^{-1/2+\eps}$ for every
$\beta\ge\tfrac12$, which is the absorption clause of
Definition~\ref{def:reg_index_fgv} at and above the threshold just identified. Both clauses
of that definition therefore hold, so the Ingham function is a function of good
variation\index[terms]{function of good variation} and its regularity index is
$\alpha(\Phi)=\tau(\Phi)=\tfrac12$.

Conversely, suppose that $\Phi$ is a function of good variation with
$\alpha(\Phi)=\tfrac12$ and take the finite forcing exponent $\beta=1$. The defining equation is
\[
 \sum_{k\le n}k a_k\Big\lfloor\frac nk\Big\rfloor=1.
\]
The identity
$\sum_{k\le n}\mu(k)\lfloor n/k\rfloor=1$ and triangular uniqueness give
$ka_k=\mu(k)$. Therefore
\[
 A(x)=\sum_{n\le x}\frac{\mu(n)}n,
 \qquad
 M(x)=xA(x)-\int_1^x A(t)\,dt.
\]
Absorption at $\beta=1$ gives $A(x)\ll_\eps x^{-1/2+\eps}$, hence
$M(x)\ll_\eps x^{1/2+\eps}$. Lemma~\ref{lem:littlewood} concludes the Riemann hypothesis.
\end{proof}

Both implications are paper proofs. Their external inputs are the verified theorems of
Conrey\index[names]{Conrey, J. B.} and Titchmarsh\index[names]{Titchmarsh, E. C.} cited above, and case (II) of Theorem~\ref{thm:transfer_law} is not used.

\begin{remark}
For every continuous function of good variation settled in this volume the arithmetic
index $\alpha(g)$ agrees with the analytic index\index[terms]{analytic index} $\eta(g)$, an agreement that is not
automatic and reflects the fit between the arithmetic grid and the analytic structure of
the kernel. The conjectured exception is $g_{\sqrt2}$, and describing it here anticipates
Chapter~\ref{chap:diophantine}, where it is treated, and the research
dossier~\ref{app:dossier_sqrt2}, which carries its proofs. There $\eta(g_{\sqrt2})=1$ is
proved, transparency fails at the exponent $\tfrac12$, and the transparency frontier\index[terms]{transparency frontier} is pinned
between $1-\log_2(177/100)$ and $\tfrac12$ by
Theorem~\ref{thm:w10-separation}. The Hardy-Littlewood-Ramanujan criterion fails at the
critical forcing, a failure that does not exclude good variation, and the existence of
the index together with the exact value $\tfrac12$ rests on the LOW and ABS estimates,
still open, the separation being produced by the Diophantine relation $(\sqrt2)^2=2$. Under the hypothesis $\Phi$ satisfies $\alpha(\Phi)=\eta(\Phi)=1/2$, so the
two indices coincide there as well.
\end{remark}

The equivalence of Theorem~\ref{thm:tauberian_rh} ties the regularity index of the Ingham
function to the single open problem the present volume does not resolve, and which it
states once here. The floor-function identities and the case analysis in $\beta$ are
developed at length in \cite{CloitreFloor}.

\begin{conjecture}[Riemann hypothesis]
\label{conj:rh}
All nontrivial zeros of the Riemann zeta function lie on the line $\Re s=\tfrac12$.
Equivalently, by Theorem~\ref{thm:tauberian_rh}, $\alpha(\Phi)=\tfrac12$.
\end{conjecture}

\section{The HLR criterion beyond the Ingham function}\label{sec:hlr}
The Hardy-Littlewood-Ramanujan criterion\index[terms]{Hardy--Littlewood--Ramanujan criterion} of Definition~\ref{def:HLR} holds for the Ingham function in its strong
form, by Proposition~\ref{prop:hlr_ingham}, and Chapter~\ref{chap:boundary} used that form as
the Tauberian condition of a boundary theorem. In the inversion of the present chapter it plays a
second role, a bound on the coefficients that validates the truncation of Perron's formula. This
section follows the criterion to the arithmetic generalizations of the Ingham function, where its
two forms separate.

\subsection{Generalized Ingham functions}\label{sec:gen_ingham}

Let $u$ be an arithmetic function. The generalized Ingham function attached to $u$ is
\begin{equation}\label{eq:gen_ingham_def}
  \Phi_u(x):=\sum_{1\le k\le x^{-1}}\frac{u(k)}{k}\,\Phi(kx),\qquad \Phi(y)=y\lfloor1/y\rfloor,
\end{equation}
a superposition of dilated copies of the Ingham function. The Dirichlet unit $\delta$, equal to
$1$ at $1$ and to $0$ beyond, returns $\Phi$ itself, every term of the superposition beyond the
first being killed. The constant function returns a different kernel,
\[
  \Phi_1(x)=x\sum_{m\le 1/x}d(m),
\]
the divisor sum appearing because the pairs $(k,m)$ with $km\le 1/x$ are counted once each. The
construction transports the Mellin transform of $\Phi$ to a Dirichlet series in $u$.

The same computation for a Dirichlet character gives the transform of the twisted kernel. For a
principal character $\Phi_\chi$ grows logarithmically at the origin, as $\Phi_1$ does, and the
integral defining the transform still converges for $\Re z<0$, so the passage to the limit of
Proposition~\ref{prop:mellin_coincidence} covers this case as well.

\begin{proposition}\label{prop:gen_mellin}
For a Dirichlet character\index[terms]{Dirichlet character} $\chi$ the Mellin transform\index[terms]{Mellin transform} of $\Phi_\chi$ is
\begin{equation}\label{eq:gen_mellin}
  \Phi_\chi^{*}(z)=\frac{z}{z-1}\,\zeta(1-z)\,L(\chi,1-z).
\end{equation}
\end{proposition}

\begin{proof}
Since $\Phi(y)=0$ for $y>1$, the sum in \eqref{eq:gen_ingham_def} runs over $k\le1/x$, and for $\Re z<0$,
\[
  \Phi_\chi^{*}(z)=-z\int_0^1\Phi_\chi(t)\,t^{-z-1}\,dt
   =-z\sum_{k\ge1}\frac{\chi(k)}{k}\int_0^{1/k}\Phi(kt)\,t^{-z-1}\,dt.
\]
The substitution $y=kt$ gives $\int_0^{1/k}\Phi(kt)t^{-z-1}\,dt=k^{z}\int_0^1\Phi(y)y^{-z-1}\,dy$, so
\[
  \Phi_\chi^{*}(z)=\Big(\sum_{k\ge1}\chi(k)k^{z-1}\Big)\Big(-z\int_0^1\Phi(y)y^{-z-1}\,dy\Big)
   =L(\chi,1-z)\,\Phi^{*}(z),
\]
and $\Phi^{*}(z)=\frac{z}{z-1}\zeta(1-z)$ gives \eqref{eq:gen_mellin}.
\end{proof}

The transform \eqref{eq:gen_mellin} carries the zeta factor of $\Phi$ multiplied by the $L$-factor of $\chi$. In the open critical strip its zeros, counted with multiplicity, are those of $\zeta(1-z)$ together with those of $L(\chi,1-z)$, a point to which the last subsection returns. On the arithmetic side the generalization is governed by the multiplicativity of $\chi$, and this is where the two forms of the criterion part.

\begin{theorem}\label{thm:gen_hlr}
Let $u$ be an arithmetic function with $u(1)\neq0$, let $u^{-1}$ be its Dirichlet inverse, and let $\Phi_u$ be the generalized Ingham function \eqref{eq:gen_ingham_def} attached to $u$. Then $\Phi_u$ satisfies the HLR criterion, that is
\[
  \sum_{k\le n}a_k\,\Phi_u\!\Big(\frac kn\Big)=n^{-\beta}\ \ (n\ge1)
  \qquad\Longrightarrow\qquad
  na_n=\mathcal{O}_{\beta,\eps}\bigl(n^{\eps}\bigr)\ \ \text{for every }\eps>0
\]
at every $\beta\ge0$, if and only if $u^{-1}(n)=\mathcal{O}_{\eps}(n^{\eps})$ for every $\eps>0$.
\end{theorem}

\begin{proof}
Write $b(n)=na_n^{\Phi}$ for the sequence attached to $\Phi$ by the defining equation with forcing $n^{-\beta}$, so that $b=\mathcal{O}(1)$ by Proposition~\ref{prop:hlr_ingham}, and write $b_u(n)=na_n$ for the sequence attached to $\Phi_u$. Unfolding \eqref{eq:gen_ingham_def} and collecting the double sum along $m=dk$, the identity $\tfrac1m\Phi(m/n)=\tfrac1n\lfloor n/m\rfloor$ turns the defining equation into $\sum_{m\le n}(b_u\star u)(m)\lfloor n/m\rfloor=n^{1-\beta}$, that is $\sum_{t\le n}(\mathbf 1\star u\star b_u)(t)=n^{1-\beta}$. Differencing at the ranks $n$ and $n-1$ leaves $(\mathbf 1\star u\star b_u)(n)=n^{1-\beta}-(n-1)^{1-\beta}$ for $n\ge2$ and the value $1$ at $n=1$, whose inversion by $\mu$ is
\begin{equation}\label{eq:gen_hlr_identity}
u\star b_u=b .
\end{equation}
That identity uses nothing beyond $u(1)\neq0$. Inverting it gives
\[
 na_n=(b\star u^{-1})(n)=\sum_{d\mid n}b(n/d)\,u^{-1}(d).
\]

Suppose $u^{-1}(d)=\mathcal{O}_\eps(d^{\eps})$. Then
\[
  |na_n|\le\sum_{d\mid n}|b(n/d)|\,|u^{-1}(d)|\ll\sum_{d\mid n}d^{\eps}\le \tau(n)\,n^{\eps}\ll n^{2\eps},
\]
since the divisor function\index[terms]{divisor function} satisfies $\tau(n)=n^{o(1)}$, and $\eps>0$ being arbitrary this is the criterion.

Conversely, take $\beta=0$. The difference above is then $n-(n-1)=1$ at every rank $n\ge2$ and $1$ at $n=1$, so $\mathbf 1\star u\star b_u=\mathbf 1$ and $u\star b_u=\mu\star\mathbf 1=\delta$, the unit of Dirichlet convolution. Hence $na_n=u^{-1}(n)$ at every rank, and the criterion at that single exponent already forces the bound on $u^{-1}$.
\end{proof}

The condition bears on the inverse and not on $u$. When $u$ is multiplicative it can be tested
prime by prime, on the coefficients of the reciprocal of the local factor.

\begin{proposition}\label{prop:gen_hlr_local}
Let $u$ be multiplicative with $u(1)=1$, let $v=u^{-1}$, and write
$U_p(X)=\sum_{k\ge0}u(p^{k})X^{k}$, so that $1/U_p(X)=\sum_{k\ge0}v(p^{k})X^{k}$. Suppose that for
every $\eps>0$ there are constants $A\ge0$ and $B\ge1$, independent of the prime and of the
exponent, with
\begin{equation}\label{eq:gen_hlr_local}
  |v(p^{k})|\le B\,(k+1)^{A}\,p^{k\eps}\qquad(p\ \text{prime},\ k\ge1).
\end{equation}
Then $v(n)=\mathcal{O}_\eps(n^{\eps})$ for every $\eps>0$, and $\Phi_u$ satisfies the HLR criterion.
\end{proposition}

\begin{proof}
The inverse of a multiplicative function is multiplicative, so
$|v(n)|=\prod_{p^{k}\|n}|v(p^{k})|$. Applying \eqref{eq:gen_hlr_local} at each prime power and
using $\prod_{p^{k}\|n}(k+1)=\tau(n)$ together with $\prod_{p^{k}\|n}p^{k\eps}=n^{\eps}$,
\[
  |v(n)|\le B^{\omega(n)}\,\tau(n)^{A}\,n^{\eps}.
\]
Now $\omega(n)\ll\log n/\log\log n$ gives $B^{\omega(n)}=n^{o(1)}$, and $\tau(n)=n^{o(1)}$ gives
$\tau(n)^{A}=n^{o(1)}$, so $|v(n)|\le n^{\eps+o(1)}$. As $\eps>0$ is arbitrary the bound follows,
and Theorem~\ref{thm:gen_hlr} applies.
\end{proof}

The hypothesis is a growth condition on the coefficients of the inverse local factors, and it is
satisfied in the three cases the theory needs. For $u$ completely multiplicative the local factor
is $U_p(X)=(1-u(p)X)^{-1}$, its reciprocal is the polynomial $1-u(p)X$, and
$v=\mu u$, so \eqref{eq:gen_hlr_local} holds with $A=0$ and $B=C_\eps$ under the Ramanujan
bound\index[terms]{Ramanujan bound}, the constant being the one that bound carries, and with $B=1$ for a Dirichlet character\index[terms]{Dirichlet character}, the example Proposition~\ref{prop:gen_mellin} uses, whose values have modulus at most one. For an
Euler product\index[terms]{Euler product} of degree at most $r$ independent of the prime,
$U_p(X)=\prod_{i\le r}(1-\alpha_i(p)X)^{-1}$, the reciprocal is a polynomial of degree $r$ whose
coefficients are the elementary symmetric functions of the $\alpha_i(p)$, and the Ramanujan bound
forces $\max_i|\alpha_i(p)|\le p^{\eps}$ through the radius of convergence of $U_p$, so
\eqref{eq:gen_hlr_local} holds with $A=0$ and $B=2^{r}$. The normalized
Ramanujan\index[names]{Ramanujan, S.} tau function\index[terms]{Ramanujan tau function} is the case $r=2$, worked out in
Remark~\ref{rem:hlr_tau}. Every $L$-function whose Euler product has bounded degree and whose
coefficients obey the Ramanujan bound therefore falls under
Proposition~\ref{prop:gen_hlr_local}, which is the form the second volume needs.

The recursion $v(p^{k})=-\sum_{j=1}^{k}u(p^{j})v(p^{k-j})$ can produce geometric growth in $k$
even when the values $u(p^{k})$ are bounded, and \eqref{eq:gen_hlr_local} is exactly what excludes
it. The next remark exhibits the two behaviors side by side.

\begin{remark}\label{rem:hlr_counterexample}
Multiplicativity of $u$ together with the Ramanujan bound does not by itself deliver the condition, and the local series is where the matter is decided. Two multiplicative multipliers make the point, both obeying $u(n)=\mathcal{O}_\eps(n^{\eps})$ and agreeing at every prime.

Take first $u(p)=2$ with $u(p^{j})=0$ for every $j\ge2$, so that $u$ vanishes at every argument divisible by a square and equals $2^{\omega(n)}$ elsewhere. Its local series is the polynomial $1+2x$, whose reciprocal is the infinite series $\sum_{k\ge0}(-2)^{k}x^{k}$, so $u^{-1}(p^{k})=(-2)^{k}$. Along $n=2^{k}$ this gives $|u^{-1}(n)|=n$, and by the converse half of Theorem~\ref{thm:gen_hlr} the sequence at $\beta=0$ is $na_n=u^{-1}(n)$ itself. The criterion fails by a full power of $n$.

Take next $u(p^{j})=2$ for every $j\ge1$, so that $u(n)=2^{\omega(n)}$ at every argument. Its local series is $(1+x)/(1-x)$, whose reciprocal is $(1-x)/(1+x)=1+2\sum_{k\ge1}(-1)^{k}x^{k}$, so $u^{-1}(n)=(-1)^{\Omega(n)}2^{\omega(n)}$, of the same modulus as $u$ itself. Here the criterion holds, and $u$ is not completely multiplicative.

The two multipliers differ only beyond the primes, which is exactly where the reciprocal of the local series is decided. The first local series vanishes at $x=-\tfrac12$, inside the unit disc, and the coefficients of its reciprocal grow like $2^{k}$. The second vanishes at $x=-1$, on the unit circle, and the coefficients of its reciprocal stay bounded. Under Proposition~\ref{prop:gen_hlr_local} with a bounded Euler degree the local series has no zero at all, being the reciprocal of a polynomial, which is the extreme case and the reason an Euler product of bounded degree settles the question.
\end{remark}

The bound of Theorem~\ref{thm:gen_hlr} cannot be sharpened to the strong form. For the character $\chi_4$ modulo $4$ the sequence $na_n$ is unbounded, and the growth is visible along the integers with many prime factors in the residue class $1\bmod4$.

The mechanism is a Dirichlet inversion turning the sequence into a divisor sum over a squarefree
modulus, so the first step is that inversion.

\begin{lemma}\label{lem:gen_chi4_inversion}
Let $\chi=\chi_4$, let $\beta$ be real, and let $(a_n)$ solve the defining equation for
$\Phi_\chi$ at the forcing $n^{-\beta}$. Put $\Delta_\beta(m)=m^{1-\beta}-(m-1)^{1-\beta}$ for
$m\ge2$ and $\Delta_\beta(1)=1$. Then
\begin{equation}\label{eq:gen_chi4_inversion}
n\,a_n=(\Delta_\beta\star h)(n),\qquad h=\mu\star(\mu\chi),
\end{equation}
where $h$ is multiplicative with $h(p)=-1-\chi(p)$, $h(p^{2})=\chi(p)$ and $h(p^{j})=0$ for
$j\ge3$. At every prime $p\equiv1\bmod4$ this gives $h(p)=-2$, so $h(n)=(-2)^{\omega(n)}$
whenever $n$ is squarefree with all its prime factors in that class.
\end{lemma}

\begin{proof}
Since $\Phi(y)=y\lfloor1/y\rfloor$ vanishes for $y>1$, multiplying the defining equation by $n$
and unfolding \eqref{eq:gen_ingham_def} gives
\[
n\sum_{k\le n}a_k\Phi_\chi\Bigl(\frac kn\Bigr)
=\sum_{jk\le n}\chi(j)\,ka_k\Bigl\lfloor\frac n{jk}\Bigr\rfloor
=\sum_{ijk\le n}\chi(j)\,ka_k
=\sum_{t\le n}(\mathbf 1\star\chi\star b)(t),
\]
where $b(k)=ka_k$. The left side is $n^{1-\beta}$ at every rank, so differencing at the ranks
$n$ and $n-1$ leaves $(\mathbf 1\star\chi\star b)(n)=\Delta_\beta(n)$. The character is
completely multiplicative, so its Dirichlet inverse is $\mu\chi$, and inverting $\mathbf 1$ by
$\mu$ gives \eqref{eq:gen_chi4_inversion}. The values of $h$ come from
$h(p^{j})=\sum_{r+s=j}\mu(p^{r})\mu(p^{s})\chi(p^{s})$, where only the terms with
$r\le1$ and $s\le1$ contribute.
\end{proof}

Along a modulus built from that residue class the inversion collapses to one divisor sum, and
the whole question becomes the size of that sum.

\begin{lemma}\label{lem:gen_chi4_split}
Let $p_1<p_2<\cdots$ be the primes congruent to $1\bmod4$, let $n_k=p_1\cdots p_k$, and let
$0<\beta<1$. Then $n_ka_{n_k}=(-2)^{k}Y_k$ with
$Y_k=\sum_{d\mid n_k}\Delta_\beta(d)(-2)^{-\omega(d)}$, and, writing
$\Delta_\beta(m)=(1-\beta)m^{-\beta}+u_\beta(m)$ with $u_\beta(1)=\beta$,
\[
Y_k=(1-\beta)\,\Pi_k+V_k,\qquad
\Pi_k=\prod_{i\le k}\Bigl(1-\tfrac12p_i^{-\beta}\Bigr),\qquad
V_k=\sum_{d\mid n_k}u_\beta(d)\,(-2)^{-\omega(d)} ,
\]
the product decreasing to zero and the sum converging absolutely to a limit $V_\infty$.
There are moreover a rank $k_0$ and constants $0<c_1\le c_2$ depending on $\beta$ alone with
\begin{equation}\label{eq:gen_chi4_tail}
c_1\,\sigma_k\ \le\ V_k-V_\infty\ \le\ c_2\,\sigma_k\qquad(k\ge k_0),
\qquad \sigma_k=\sum_{j>k}p_j^{-1-\beta} .
\end{equation}
\end{lemma}

\begin{proof}
Every divisor of $n_k$ is squarefree with all its prime factors in the class, so
$h(n_k/d)=(-2)^{k-\omega(d)}$ by Lemma~\ref{lem:gen_chi4_inversion} and
\eqref{eq:gen_chi4_inversion} gives the stated form of $Y_k$.

For $m\ge2$ the increment is $\Delta_\beta(m)=(1-\beta)\int_{m-1}^{m}t^{-\beta}\,dt$,
so $u_\beta(m)=(1-\beta)\int_{m-1}^{m}(t^{-\beta}-m^{-\beta})\,dt$, and on the range
of integration $t^{-\beta}-m^{-\beta}=\beta\int_t^{m}s^{-1-\beta}\,ds$ lies between
$\beta m^{-1-\beta}(m-t)$ and $\beta(m-1)^{-1-\beta}(m-t)$. Integrating and using
$\int_{m-1}^{m}(m-t)\,dt=\tfrac12$,
\begin{equation}\label{eq:gen_chi4_ubound}
\tfrac12(1-\beta)\beta\,m^{-1-\beta}\ \le\ u_\beta(m)\ \le\ \tfrac12(1-\beta)\beta\,(m-1)^{-1-\beta},
\end{equation}
so $u_\beta(m)=\tfrac12(1-\beta)\beta\,m^{-1-\beta}(1+\theta_m)$ with $0\le\theta_m\le16/m$,
the last bound coming from $(1-1/m)^{-1-\beta}-1\le(1+\beta)2^{2+\beta}/m$ for $m\ge2$ and
$0<\beta<1$.

The split $\Delta_\beta=(1-\beta)(\cdot)^{-\beta}+u_\beta$ holds at $d=1$ as well, since
$1=(1-\beta)+\beta$, and the power is completely multiplicative, so its divisor sum against
$(-2)^{-\omega}$ factors into $\Pi_k$. By the prime number theorem in the progression
$1\bmod4$, Proposition~\ref{prop:pnt_progressions}, one has $\sum_i p_i^{-1}=\infty$. Since $0<\beta<1$, the sum
$\sum_i p_i^{-\beta}$ dominates it, and the product decreases to zero. Absolute convergence of the
remainder follows from \eqref{eq:gen_chi4_ubound} together with
$(d-1)^{-1-\beta}\le2^{1+\beta}d^{-1-\beta}$, the sum $\sum_dd^{-1-\beta}2^{-\omega(d)}$ over
the squarefree integers built from the class being the convergent product
$\prod_p(1+\tfrac12p^{-1-\beta})$.

For the tail, group the squarefree integers built from the primes $p_i$ that do not divide
$n_k$ by their largest prime factor. Writing such an integer as $p_je$ with $j>k$ and
$e\mid p_1\cdots p_{j-1}$ gives
$(-2)^{-\omega(p_je)}=-\tfrac12(-2)^{-\omega(e)}$, so
$V_k-V_\infty=\tfrac12\sum_{j>k}T_j$ with $T_j=\sum_{e}(-2)^{-\omega(e)}u_\beta(p_je)$. The
expansion above turns $T_j$ into $\tfrac12(1-\beta)\beta\,p_j^{-1-\beta}(Q_j+R_j)$, where
$Q_j=\prod_{i<j}(1-\tfrac12p_i^{-1-\beta})$ lies between the convergent product
$P(\beta)=\prod_i(1-\tfrac12p_i^{-1-\beta})$ and $1$, while $|R_j|\le C_1(\beta)/p_j$ with
$C_1(\beta)=16\prod_p(1+\tfrac12p^{-1-\beta})$. Choosing $k_0$ so large that
$p_j\ge2C_1(\beta)/P(\beta)$ for every $j>k_0$ keeps $Q_j+R_j$ between $\tfrac12P(\beta)$ and
$2$, and summing over $j>k$ gives \eqref{eq:gen_chi4_tail} with
$c_1=\tfrac18(1-\beta)\beta P(\beta)$ and $c_2=(1-\beta)\beta$.
\end{proof}

The three regimes of the forcing exponent now separate, and in each of them the divisor sum
stays away from zero.

\begin{theorem}\label{thm:gen_unbounded}
Let $\chi_4$ be the nontrivial character modulo $4$, let $\beta>0$, and let $n_k=p_1\cdots p_k$
be the product of the first $k$ primes congruent to $1\bmod4$. The sequence attached to
$\Phi_{\chi_4}$ satisfies $|n_ka_{n_k}|\to\infty$. Quantitatively there are a rank $k_1(\beta)$
and a constant $c(\beta)>0$ with
\begin{equation}\label{eq:gen_unbounded_rate}
|n_ka_{n_k}|\ \ge\ c(\beta)\,2^{k}\min\Bigl(1,\ \sum_{j>k}p_j^{-\gamma(\beta)}\Bigr)
\qquad(k\ge k_1),\qquad
\gamma(\beta)=\begin{cases}1+\beta,&0<\beta\le1,\\[2pt] \beta,&\beta>1,\end{cases}
\end{equation}
and the right side tends to infinity. In particular $\Phi_{\chi_4}$ does not satisfy the strong
Hardy-Littlewood-Ramanujan criterion\index[terms]{Hardy--Littlewood--Ramanujan criterion} at any positive exponent.
\end{theorem}

\begin{proof}
Take first $0<\beta<1$ and write $Y_k=(1-\beta)\Pi_k+V_\infty+(V_k-V_\infty)$ as in
Lemma~\ref{lem:gen_chi4_split}. The first term is positive and the third is at least
$c_1\sigma_k>0$ for $k\ge k_0$. If $V_\infty\ge0$ the three terms carry the same sign and
$Y_k\ge c_1\sigma_k$. If $V_\infty<0$ the other two tend to zero, so $Y_k$ tends to $V_\infty$
and $|Y_k|\ge|V_\infty|/2$ beyond some rank. In both cases $|n_ka_{n_k}|=2^{k}|Y_k|$ is
bounded below by a positive multiple of $2^{k}\min(1,\sigma_k)$.

At $\beta=1$ the increment $\Delta_1$ is the unit of Dirichlet convolution, so
\eqref{eq:gen_chi4_inversion} reads $n_ka_{n_k}=h(n_k)=(-2)^{k}$ and the modulus is exactly
$2^{k}$.

For $\beta>1$ put $v(m)=-\Delta_\beta(m)=(\beta-1)\int_{m-1}^{m}t^{-\beta}\,dt$, which
is positive at every $m\ge2$ and satisfies $v(m)=(\beta-1)m^{-\beta}(1+\theta'_m)$ with
$0\le\theta'_m\le C_0(\beta)/m$, so that
$Y_k=1-\sum_{2\le d\mid n_k}v(d)(-2)^{-\omega(d)}$. The series converges absolutely because
$\beta>1$, hence $Y_k$ tends to a limit $Y_\infty$. Grouping the omitted squarefree integers by
their largest prime factor gives
\[
Y_\infty-Y_k
=\frac{\beta-1}{2}\sum_{j>k}p_j^{-\beta}\bigl(Q'_j+R'_j\bigr),\qquad
Q'_j=\prod_{i<j}\Bigl(1-\tfrac12p_i^{-\beta}\Bigr),
\qquad |R'_j|\le \frac{C_2(\beta)}{p_j}.
\]
Since $\beta>1$, the products $Q'_j$ decrease to the positive limit
$P'(\beta)=\prod_i(1-\tfrac12p_i^{-\beta})$. Thus $Q'_j+R'_j\ge P'(\beta)/2$ beyond some
rank, and consequently
$Y_\infty-Y_k\ge c_3(\beta)\sum_{j>k}p_j^{-\beta}$. If
$Y_\infty\le0$ then $|Y_k|\ge c_3\sum_{j>k}p_j^{-\beta}$, and if $Y_\infty>0$ then
$|Y_k|\ge Y_\infty/2$ beyond some rank.

The exponent appearing in the tail sum is in each regime the $\gamma(\beta)$ of
\eqref{eq:gen_unbounded_rate}. Divergence of the right side follows from
$\sum_{j>k}p_j^{-\gamma}\ge p_{k+1}^{-\gamma}$. The prime number theorem in arithmetic
progressions, Proposition~\ref{prop:pnt_progressions} or \cite[\S~II.8]{Tenenbaum2015}, gives
$p_{k+1}\sim2k\log k$ for the primes $1\bmod4$, so the lower bound exceeds a positive multiple of
$2^{k}(k\log k)^{-\gamma}$ and tends to infinity.
\end{proof}

The growth is located on a sparse set of integers, where one part of the transfer law is
explicit and can be compared with the answer just obtained.

\begin{proposition}
\label{prop:gen_unbounded_leading}
With the notation of Theorem~\ref{thm:gen_unbounded}, and for $0<\beta<1$, the leading part of
$na_n$ in
the decomposition of Theorem~\ref{thm:transfer_law} has modulus at least
$(1-\beta)\prod_{i\le k}\bigl(2-p_i^{-\beta}\bigr)$ along $n_k=p_1\cdots p_k$, hence it tends to
infinity.
\end{proposition}

\begin{proof}
From $na_n=(b\star\chi_4^{-1})(n)$ with $\chi_4^{-1}=\mu\chi_4$, and from the decomposition $b=(1-\beta)J_{-\beta}+(\mu\star u)+(\mu\star w)$ of Theorem~\ref{thm:transfer_law}, the leading part of $na_n$ is $(1-\beta)\,(J_{-\beta}\star\chi_4^{-1})(n)$. The function $J_{-\beta}\star\chi_4^{-1}$ is multiplicative, and at a prime $p\equiv1\bmod4$, where $\chi_4(p)=1$,
\[
  (J_{-\beta}\star\chi_4^{-1})(p)=\big(p^{-\beta}-1\big)-\chi_4(p)=p^{-\beta}-2\longrightarrow-2\qquad(p\to\infty).
\]
Each factor has modulus at least $2-p^{-\beta}>1$, the inequality using $\beta>0$, and the
constant $1-\beta$ is positive because $\beta<1$. Along $n_k=p_1\cdots p_k$ the leading part
therefore has modulus at least $(1-\beta)\prod_{i\le k}(2-p_i^{-\beta})\to\infty$. Both ends of
the range are needed. At $\beta=1$ the constant vanishes and the bound is empty, and at $\beta>1$
it changes sign, so the statement is confined to $0<\beta<1$ although
Theorem~\ref{thm:gen_unbounded} covers every $\beta>0$.
\end{proof}

\begin{remark}[Where the growth actually sits]\label{rem:gen_unbounded_leading}
The proposition and the theorem measure two different things, and comparing them locates the
source of the growth. The leading part is $(1-\beta)(-2)^{k}\Pi_k$ along $n_k$, and $\Pi_k$
tends to zero by Lemma~\ref{lem:gen_chi4_split}, so that part is of order $o(2^{k})$ while
$|n_ka_{n_k}|$ itself is of order $2^{k}$ whenever $V_\infty\neq0$. The transfer law splits
$na_n$ into three pieces and it is not the first that carries the growth here, it is the
remainder $V_\infty(-2)^{k}$. The numerical values agree with that reading. At $\beta=0.3$ the
sequence takes the values $-1.55$, $+2.65$, $-4.66$ at $n_1=5$, $n_2=65$, $n_3=1105$, so the
sign alternates as the product of the factors $p^{-\beta}-2$ predicts, while the magnitudes
exceed the leading part alone and follow $2^{k}$ instead.
\end{remark}

In all three regimes used in the proof, the normalized divisor sums
$Y_k=(-2)^{-k}n_ka_{n_k}$ are $\mathcal O_\beta(1)$: this follows from
Lemma~\ref{lem:gen_chi4_split} when $0<\beta<1$, is exact when $\beta=1$, and follows from
absolute convergence when $\beta>1$. Hence
\[
|n_ka_{n_k}|\ll_\beta2^{k}=2^{\omega(n_k)}=n_k^{o(1)}.
\]
Theorem~\ref{thm:gen_unbounded} also shows that these values tend to infinity. They therefore
violate the strong form $\mathcal{O}(1)$ while remaining within the weak form
$\mathcal{O}(n^{\eps})$ of Theorem~\ref{thm:gen_hlr}. The separation of the two forms is thus
an arithmetic fact about $\chi_4$, and the exponent $\eps$ is not a slack in the argument but
the trace of the Ramanujan bound on the multiplier.

\begin{remark}\label{rem:euler}
The passage from the strong form to the weak form is the passage from the trivial multiplier to a multiplier carrying an Euler product, and the condition isolated above is that the Dirichlet inverse inherit the Ramanujan bound. An Euler product of bounded degree delivers it, by Proposition~\ref{prop:gen_hlr_local}, because the reciprocal of such a local factor is a polynomial and the inverse is therefore carried by finitely many prime powers at each prime. Ordinary multiplicativity is not enough, as Remark~\ref{rem:hlr_counterexample} shows with a multiplier whose local series is a polynomial and whose inverse grows. The Ramanujan bound of the Selberg\index[names]{Selberg, A.} class is in this reading the arithmetic residue of unique factorization, and the character twists, whose multipliers are completely multiplicative, are exactly the case where the criterion can be established. This is one heuristic reason to expect the Riemann hypothesis for the $L$-functions that carry both a functional equation\index[terms]{functional equation} and an Euler product, and not for Dirichlet series lacking the second. The oldest example on the second side is the Epstein\index[names]{Epstein, P.} zeta function\index[terms]{Epstein zeta function} of a positive binary quadratic form of class number greater than one, which has a functional equation, no Euler product, and infinitely many zeros in the half plane of absolute convergence \cite{DavenportHeilbronn1936}. The kernels built from such series are the ones Chapter~\ref{chap:diophantine} removes.
\end{remark}

Between the two cases lies a third, and it is the one that carries the most arithmetic.

\begin{remark}[The Ramanujan tau function]\label{rem:hlr_tau}
The normalized Ramanujan\index[names]{Ramanujan, S.} tau function\index[terms]{Ramanujan tau function}
$u(n)=\tau(n)/n^{11/2}$ is multiplicative without being completely multiplicative, and it is the
first case beyond degree one of Proposition~\ref{prop:gen_hlr_local}. What stands in for complete
multiplicativity is the shape of the Euler factor, and the degree is two. \nm{Mordell}{L. J.} proved that $\tau$ is multiplicative and obeys
$\tau(p^{k+1})=\tau(p)\tau(p^{k})-p^{11}\tau(p^{k-1})$ at every prime power
\cite{Mordell1917}, so the Dirichlet series of $u$ has the degree two Euler
product\index[terms]{Euler product} $\prod_p\big(1-u(p)p^{-s}+p^{-2s}\big)^{-1}$, and its Dirichlet inverse is the
multiplicative function determined by
\[
u^{-1}(p)=-u(p),\qquad u^{-1}(p^{2})=1,\qquad u^{-1}(p^{k})=0\ \ (k\ge3).
\]
The inverse is therefore carried by the cube free integers alone, and the bound
$|u(p)|\le2$ of \nm{Deligne}{P.}~\cite{Deligne1974} gives
$|u^{-1}(n)|\le2^{\omega(n)}=n^{o(1)}$. That is condition \eqref{eq:gen_hlr_local} with $A=0$ and
$B=2$, so Proposition~\ref{prop:gen_hlr_local} applies and a degree two Euler product secures the
inheritance as surely as a degree one factor does. What carries the case is the finite degree of
the Euler factor and not multiplicativity alone. The tau multiplier appears in
\cite[\S~6]{Cloitre2016}.
\end{remark}
\begin{numobs}\label{numobs:tau_inverse}
The three displayed values were recomputed for this volume from the exact integer expansion
$q\prod_{n\ge1}(1-q^{n})^{24}=\sum_{n\ge1}\tau(n)q^{n}$ up to $n=400$, with the Dirichlet
inverse formed in sixty digit arithmetic. The inverse vanishes at every non cube free integer
to within $5\cdot10^{-61}$, the ratio $|u^{-1}(n)|\,2^{-\omega(n)}$ stays below $1$ with
maximum $0.9594$ attained at $n=103$, and $\max_{n\le400}|u(n)|=1.9188$, consistent with
Deligne\index[names]{Deligne, P.}'s bound.
\end{numobs}

\subsection{The generalized equivalence}\label{sec:gen_equiv}

The transform \eqref{eq:gen_mellin} for $\chi_4$ reads $\Phi_{\chi_4}^{*}(z)=\frac{z}{z-1}\zeta(1-z)\beta(1-z)$, where $\beta(s)=L(s,\chi_4)$ is the Dirichlet beta function. The factor $z/(z-1)$ has a zero at $z=0$ cancelled by the pole of $\zeta(1-z)$, so $\Phi_{\chi_4}^{*}(0)=\pi/4\neq0$. The zeros in the strip come from the two automorphic factors. Writing $\Theta_\zeta=\sup\{\Re\rho:\zeta(\rho)=0,\,0<\Re\rho<1\}$ and $\Theta_\beta=\sup\{\Re\rho:\beta(\rho)=0,\,0<\Re\rho<1\}$, each nontrivial zero\index[terms]{nontrivial zero} of $\zeta$ or of $\beta$ contributes a zero of $\Phi_{\chi_4}^{*}$ at $z=1-\rho$, so the analytic index is
\begin{equation}\label{eq:gen_analytic_index}
  \inf\{\Re z:\Phi_{\chi_4}^{*}(z)=0\}=1-\max(\Theta_\zeta,\Theta_\beta).
\end{equation}
This equals $1/2$ precisely when $\Theta_\zeta=\Theta_\beta=1/2$, that is, when the Riemann hypothesis holds for $\zeta$ and for $\beta$ at once.

The equivalence is established by the same route as Theorem~\ref{thm:tauberian_rh}, with the
single Dirichlet series $1/\zeta$ replaced by $1/(\zeta L)$. Two preparations are needed. The
first is the arithmetic identity that carries the defining equation to a convolution.

\begin{lemma}\label{lem:chi4_convolution}
Let $(a_n)$ satisfy $\sum_{k\le n}a_k\Phi_{\chi_4}(k/n)=n^{-\beta}$ for $n\ge1$, and put
$\Delta_\beta(1)=1$ and $\Delta_\beta(n)=n^{1-\beta}-(n-1)^{1-\beta}$ for $n\ge2$. Then
\[
  n\,a_n=\bigl(\Delta_\beta\star\mu\star(\mu\chi_4)\bigr)(n),
  \qquad
  \sum_{n\ge1}\frac{n\,a_n}{n^{s}}=\frac{C_\beta(s)}{\zeta(s)\,L(s,\chi_4)},
\]
where $C_\beta(s)=\sum_{n\ge1}\Delta_\beta(n)n^{-s}$.
\end{lemma}

\begin{proof}
Unfolding \eqref{eq:gen_ingham_def} gives
$\Phi_{\chi_4}(k/n)=\tfrac kn\sum_{j\le n/k}\chi_4(j)\lfloor n/(jk)\rfloor$, and collecting the
double sum along $m=jk$ turns the defining equation into
$\sum_{m\le n}b(m)\lfloor n/m\rfloor=n^{1-\beta}$ with $b=(k\mapsto ka_k)\star\chi_4$. The left
side is $\sum_{N\le n}(b\star\mathbf 1)(N)$, so differencing at the ranks $n$ and $n-1$ leaves
$(b\star\mathbf 1)(n)=\Delta_\beta(n)$, the value at $n=1$ being $1$. M\"obius inversion gives
$b=\Delta_\beta\star\mu$. The character being completely multiplicative, its Dirichlet inverse is
$\mu\chi_4$, and $na_n=b\star(\mu\chi_4)$ follows. The series identity is the transcription of the
convolution, the series of $\mu$ being $1/\zeta$ and that of $\mu\chi_4$ being $1/L(\cdot,\chi_4)$.
\end{proof}

The second preparation is the zero that forbids transparency above one half.

\begin{proposition}\label{prop:gen_chi4_noncommon}
For every fixed $\beta\in(\tfrac12,1)$ there is a zero $\rho=\tfrac12+i\gamma$, $\gamma>0$, of
$\zeta(s)L(s,\chi_4)$ with $C_\beta(\rho)\neq0$.
\end{proposition}

\begin{proof}
Take the zero $\rho$ supplied by Proposition~\ref{prop:C_noncommon_zero}.
The partial sums $B(u)=\sum_{n\le u}\chi_4(n)$ satisfy $0\le B(u)\le1$.
Partial summation therefore gives
\[
 L(s,\chi_4)=s\int_1^\infty B(u)u^{-s-1}\,du\qquad(\Re s>0),
\]
where the integral is locally uniformly convergent. Thus $L(\cdot,\chi_4)$
is holomorphic at $\rho$ and
\[
 \zeta(\rho)L(\rho,\chi_4)=0,\qquad C_\beta(\rho)\ne0.
\]
The product need not have a simple zero at $\rho$.
\end{proof}

\begin{theorem}\label{thm:gen_equiv}
The generalized Ingham function $\Phi_{\chi_4}$ is a function of good variation\index[terms]{function of good variation}
with $\alpha(\Phi_{\chi_4})=\tfrac12$ if and only if the Riemann hypothesis holds for $\zeta$ and
for $L(\cdot,\chi_4)$.
\end{theorem}

\begin{proof}
Write $A_1(x)=\sum_{n\le x}na_n$ and $A(x)=\sum_{n\le x}a_n$, and use
Lemma~\ref{lem:chi4_convolution} throughout.

\emph{Sufficiency of the two hypotheses.} Assume both.
For $\beta\ge0$ the coefficients admit a divisor bound. By
Lemma~\ref{lem:chi4_convolution} they factor as
$n a_n=(\Delta_\beta\star\mu)\star(\mu\chi_4)$.
The first factor is bounded by Proposition~\ref{prop:hlr_ingham}
for $\beta>0$, and equals the Dirichlet unit when $\beta=0$.
Since $|\mu\chi_4|\le1$,
\begin{equation}\label{eq:chi4_divisor_bound}
 |n a_n|\le\sum_{d\mid n}|(\Delta_\beta\star\mu)(d)|
             \ll_\beta\tau(n)\qquad(\beta\ge0).
\end{equation}
The restriction on $\beta$ matters. Negative forcing exponents are treated
by absolute convergence, as follows.

For $\beta<0$, let $(b_n)$ be the exact Ingham solution at the same
exponent and put $B(x)=\sum_{n\le x}b_n$.
The convolution identity gives the exact relation
\[
 A(x)=\sum_{d\le x}\frac{\mu(d)\chi_4(d)}{d}B(x/d).
\]
Lemma~\ref{lem:ingham_beta_neg} gives
$t^\beta B(t)\to1/\Phi^*(\beta)$ and bounds this expression
for all $t\ge1$. After multiplication by $x^\beta$, the summand
has modulus at most $C_\beta d^{\beta-1}$, a summable sequence.
Passing to the limit by dominated convergence therefore gives
\[
 x^\beta A(x)\longrightarrow
 \frac1{\Phi^*(\beta)}\sum_{d\ge1}\mu(d)\chi_4(d)d^{\beta-1}
 =\frac1{\Phi^*(\beta)L(1-\beta,\chi_4)}
 =\frac1{\Phi_{\chi_4}^*(\beta)}.
\]
This calculation uses neither zero hypothesis.

At $\beta=0$, one has $\Delta_0=\mathbf1$ and
$na_n=\mu(n)\chi_4(n)$. Corollary~\ref{cor:perron-mobius-halfplane},
applied to $L(\cdot,\chi_4)$ with $\theta=1/2$, proves
$M_{\chi_4}(x)\ll_\eps x^{1/2+\eps}$.
Partial summation makes its Dirichlet series holomorphic in
$\Re s>1/2$, where it agrees with $1/L(s,\chi_4)$ by continuation.
In particular it converges at $s=1$, and, for $0<\eps<1/2$,
\[
 A(x)=\frac1{L(1,\chi_4)}+
       \frac{M_{\chi_4}(x)}x
       -\int_x^\infty\frac{M_{\chi_4}(t)}{t^2}\,dt
 =\frac4\pi+\mathcal O_\eps(x^{-1/2+\eps}).
\]
Here $L(1,\chi_4)=\sum_{m\ge0}(-1)^m/(2m+1)=\pi/4$.
The last equality follows by integrating the finite geometric sum for
$(1+t^2)^{-1}$ on $[0,1]$, whose integrated remainder tends to zero.
The limit is $1/\Phi_{\chi_4}^*(0)$, as required.

For $\beta>0$ write
\[
 \mathcal B_\beta(s)=\frac{C_\beta(s)}{\zeta(s)L(s,\chi_4)},
 \qquad C_\beta(s)=(1-\beta)\zeta(s+\beta)+U(s).
\]
The coefficients of $U$ satisfy $u(n)=\mathcal O_\beta(n^{-\beta-1})$,
including $u(1)=\beta$, so $U$ is bounded and holomorphic on every
closed half plane strictly inside $\Re s>-\beta$.
Fix $0<\eps<1/4$. Choose
\[
 \delta=
 \begin{cases}
 \min\bigl(\eps/4,(1/2-\beta)/2\bigr),&0<\beta<1/2,\\
 \eps/4,&\beta\ge1/2,
 \end{cases}
 \qquad b=\frac12+\delta,\qquad c=\frac54.
\]
Both reciprocals are holomorphic on this strip under their respective
zero hypotheses, with a removable zero at $s=1$ for $1/\zeta$.
Lemma~\ref{lem:perron-zerofree-growth} applies to the three factors
$\zeta(s+\beta)$, $1/\zeta(s)$ and $1/L(s,\chi_4)$.
Applying it with exponent $\nu/3$ to each gives
\[
 |\mathcal B_\beta(\sigma+it)|\ll_{\beta,\nu}(1+|t|)^\nu
 \qquad(b\le\sigma\le c)
\]
at large heights. The bounded $U$ term satisfies this estimate as well.

Put $y=N+1/2$, $T=y^2$ and $\nu=\eps/8$.
The elementary divisor estimate $\tau(n)\ll n^{1/8}$ makes
\eqref{eq:chi4_divisor_bound} a coefficient bound of the form required
by Corollary~\ref{cor:perron-subpower}. To recall its elementary proof,
write $n=\prod p^a$. For $p\ge2^8$ use
$a+1\le2^a\le p^{a/8}$. For each of the finitely many smaller primes,
$\sup_{a\ge0}(a+1)p^{-a/8}<\infty$. Their product gives the constant.
Thus the vertical side, both horizontal sides and the truncation contribute
respectively
\[
 \mathcal O_{\beta,\eps}(y^{1/2+\delta+2\nu}),\qquad
 \mathcal O_{\beta,\eps}\!\left(\frac{y^{-3/4+2\nu}}{\log y}\right),
 \qquad
 \mathcal O_\beta(y^{-3/4}+y^{-7/8}\log y),
\]
with their full integral estimates proved in Appendix~\ref{app:perron}.
All three are $\mathcal O_{\beta,\eps}(y^{1/2+\eps})$.

If $0<\beta<1/2$, the only pole crossed is $s=1-\beta>b$ and
\[
 \Res_{s=1-\beta}\left(\mathcal B_\beta(s)\frac{y^s}{s}\right)
 =\frac{y^{1-\beta}}{\zeta(1-\beta)L(1-\beta,\chi_4)}.
\]
If $\beta\ge1/2$, no pole is crossed. In particular at $\beta=1/2$
the numerator pole is to the left of the chosen line, so there is no
boundary residue to assign. At $\beta=1$, $C_1=1$.
Replacing $y$ by $N$ changes the possible main term by
$\mathcal O_\beta(N^{-\beta})$ and does not change the partial sum.
Consequently
\[
 A_1(N)=\frac{N^{1-\beta}}{\zeta(1-\beta)L(1-\beta,\chi_4)}
       +\mathcal O_{\beta,\eps}(N^{1/2+\eps})
       \qquad(0<\beta<1/2),
\]
\[
 A_1(N)=\mathcal O_{\beta,\eps}(N^{1/2+\eps})
       \qquad(\beta\ge1/2).
\]
Abel summation\index[terms]{Abel summation} turns each display into a main term, a constant and a
remainder,
\[
  A(x)=\frac{A_1(x)}{x}+\int_1^{x}\frac{A_1(t)}{t^{2}}\,dt
  =\frac{1}{\Phi^{*}_{\chi_4}(\beta)}\,x^{-\beta}+c_0+\mathcal{O}_{\beta,\eps}(x^{-1/2+\eps})
  \quad(0<\beta<\tfrac12),
\]
\[
  A(x)=c_0+\mathcal{O}_{\beta,\eps}(x^{-1/2+\eps})\quad(\beta\ge\tfrac12),
\]
the coefficient in the first being $\bigl(1-\tfrac1\beta\bigr)/(\zeta(1-\beta)L(1-\beta,\chi_4))$,
which is $1/\Phi^{*}_{\chi_4}(\beta)$ by \eqref{eq:gen_mellin}. For each fixed $\beta>0$, the constant $c_0=c_0(\beta)$ is the limit of $A(x)$.
Its value must be identified in each regime.

It vanishes, for every $\beta>0$, by the Abelian argument proved in case (III) of
Theorem~\ref{thm:transfer_law}. The series
$\sum_n a_nn^{-w}=\mathcal B_\beta(1+w)$ converges absolutely for $\Re w>0$ by
\eqref{eq:chi4_divisor_bound}. Since $A(x)\to c_0$, partial summation and the split of
$w\int_1^\infty(A(t)-c_0)t^{-w-1}\,dt$ at a fixed large abscissa show that
$\sum_n a_nn^{-w}\to c_0$ as real $w\to0^+$. Therefore
\begin{equation}\label{eq:chi4_constant_vanishes}
  c_0=\lim_{w\to0^{+}}\mathcal B_\beta(1+w)
   =\lim_{w\to0^{+}}\frac{C_\beta(1+w)}{\zeta(1+w)\,L(1+w,\chi_4)}=0 .
\end{equation}
The numerator tends to the finite value $C_\beta(1)$, the series
$\sum_n\Delta_\beta(n)n^{-1-w}$ converging at $w=0$ because
$|\Delta_\beta(n)|\ll_\beta n^{-\beta}$ for $n\ge2$ and $\beta>0$.
This bound includes $\beta=1$, where all these increments vanish. The factor $L(1+w,\chi_4)$ tends to $L(1,\chi_4)=\pi/4$, which is not zero. And
$\zeta(1+w)\to\infty$, $\zeta$ having its pole at $1$. The quotient therefore tends to zero. This
is the mechanism already used in the proof of Theorem~\ref{thm:tauberian_rh}, where the same pole
returns $\sum_{n\ge1}\mu(n)/n=0$.

With $c_0=0$, choose $\eps<1/2-\beta$ in the first display. Every $\beta<\tfrac12$ is transparent, with transparent
coefficient $1/\Phi^{*}_{\chi_4}(\beta)$, and for every $\beta\ge\tfrac12$ the second display
becomes $A(x)\ll_{\beta,\eps}x^{-1/2+\eps}$, which is absorption at the rate one half.
Proposition~\ref{prop:gen_chi4_noncommon} forbids transparency above one half,
by the argument of Corollary~\ref{cor:phi_half} applied to $\zeta L$ in place of $\zeta$. Hence
$\tau(\Phi_{\chi_4})=\tfrac12$ with absorption above it, which is membership with index one half.

\emph{Necessity.} Assume $\alpha(\Phi_{\chi_4})=\tfrac12$ and take $\beta=1$. Then
$\Delta_1(1)=1$ and $\Delta_1(n)=0$ for $n\ge2$, so $\Delta_1$ is the Dirichlet unit and
$na_n=(\mu\star\mu\chi_4)(n)$, of Dirichlet series $1/(\zeta(s)L(s,\chi_4))$. Absorption at
$\beta=1$ gives $A(x)\ll_\eps x^{-1/2+\eps}$, and Abel summation gives
$A_1(x)\ll_\eps x^{1/2+\eps}$. The series $1/(\zeta(s)L(s,\chi_4))$ therefore converges on
$\Re s>\tfrac12$, where it is consequently holomorphic, so $\zeta(s)L(s,\chi_4)$ has no zero in
that half plane. The functional equations reflect the nontrivial zeros of each factor
about the critical line, so every such zero has real part $1/2$. That is the Riemann
hypothesis for $\zeta$ and for $L(\cdot,\chi_4)$ at once.
\end{proof}

\begin{proofstatus}{The theorem is unconditional as a statement of equivalence, and it asserts
neither of the two hypotheses. Appendix~\ref{app:perron} proves the finite inversion,
the contour estimates, the vertical control under each zero hypothesis and the bound for
$M_{\chi_4}$. The two complex analysis inequalities used in that appendix are cited in their
exact forms, with their disk and annulus hypotheses checked. The argument excluding
transparency above one half is Proposition~\ref{prop:gen_chi4_noncommon}.
The equivalence is the terminal result of \S\ref{sec:gen_equiv}, and it is not used as an
assumption in any later proof.}
\end{proofstatus}

Theorem~\ref{thm:gen_equiv} places the two hypotheses on the same footing as the single index of one arithmetic function. It is the generalization toward which the criterion points. The Ingham function reads $\zeta$ alone, through its analytic index $\eta(\Phi)=1-\Theta_\zeta$, and the character-twisted function reads $\zeta$ and $\beta$ together through \eqref{eq:gen_analytic_index}, by the same reading of the zeros of the transform. The equivalences of this chapter identify each arithmetic index with its analytic reading at the critical value one half. The Euler product of $\chi_4$, which makes $\Phi_{\chi_4}$ satisfy the criterion, is in this sense the arithmetic condition that binds the two zero sets to a single regularity index.

\section{Counting problems as an application}\label{sec:counting}

The floor function and the M\"obius\index[terms]{M\"obius function} function together produce counting
formulas, and a counting formula is a right-hand side. That is the whole of what this section
observes, and it turns the equivalence of the chapter into a way of estimating partial sums of
arithmetic functions whose input is a count rather than a Dirichlet series\index[terms]{Dirichlet series}.

Three families of such formulas are standard. The sum $\sum_{k\le n}\mu(k)\lfloor n/k\rfloor^{m}$
counts the $m$-tuples of integers at most $n$ whose greatest common divisor is one, the case $m=1$
being Meissel\index[names]{Meissel, E.}'s identity. The sum $\sum_{k\le n}\mu(k)\lfloor n/k^{p}\rfloor$ counts the $p$-free
integers up to $n$. And for a squarefree $P=p_1\cdots p_m$, the sum
$(-1)^{m}\sum_{k\le n}\mu(Pk)\lfloor n/k\rfloor$ counts the integers up to $n$ all of whose prime
factors divide $P$. Several of these counts are tabulated \cite{OEIS}.

What makes such a count usable here is its rate. When the count is slowly varying\index[terms]{slowly varying} and its
increments are bounded, the defining equation\index[terms]{defining equation} of this theory is being forced by
$n^{-\beta}L_0(n)$ rather than by a pure power, and the transfer of
\S\ref{sec:master_equiv} extends to that forcing. The extension used below is the following one,
and it is quoted with the hypotheses under which it is proved.

\begin{lemma}[Transfer under a slowly varying forcing]\label{lem:sv_forcing}
Let $L_0$ be slowly varying in the sense of Karamata\index[names]{Karamata, J.}, increasing on
$[1,\infty)$, and such that $0\le L_0(n+1)-L_0(n)<C$ for some constant $C$. Let $a_1=1$ and let
$(a_k)$ satisfy, for $n\ge2$,
\[
\sum_{k\le n}a_k\,\Phi\Big(\frac kn\Big)=n^{-\beta}L_0(n).
\]
Then, under the Riemann hypothesis\index[terms]{Riemann hypothesis}, for $\beta\ge\tfrac12$ and
every $\eps>0$,
\[
\sum_{n\le x}a_n\ll_\eps x^{-1/2+\eps}\qquad(x\to\infty).
\]
\end{lemma}
\begin{proofstatus}{Quoted from \cite[Theorem~5.1]{CloitreFloor} and not proved here. The three
hypotheses on $L_0$ are the ones that the proof uses, monotonicity and bounded increments entering
through an Abel summation\index[terms]{Abel summation} on the differenced factor, and slow
variation entering through $L_0(n)=\mathcal{O}(n^{\eps})$. The general statement for an arbitrary
kernel is announced in the same place and is not proved there, so it is not quoted here. The
lemma is used here only in the two examples immediately below, both at $\beta=1$, the sums of
$\mu(6n)$ and of $(-1)^{n-1}\mu(n)$. They are examples rather than a classification, and the broader
RAF encounter with combinatorial number theory remains the open direction stated after them.}
\end{proofstatus}

Take $P=6$. The sum $\sum_{k\le n}\mu(6k)\lfloor n/k\rfloor$ counts the integers up to $n$
whose prime factors are two and three. \nm{Ramanujan}{S.} showed that this count is
asymptotic to $\log(2n)\log(3n)/(2\log2\log3)$, and it is slowly varying with increments
in $\{0,1\}$.
Dividing by $n$ writes the count as a forcing of the Ingham kernel\index[terms]{Ingham kernel},
\[
\sum_{k\le n}\frac{\mu(6k)}{k}\,\Phi\Big(\frac kn\Big)=n^{-1}L_0(n),
\]
and under the Riemann hypothesis\index[terms]{Riemann hypothesis}, the exponent one lying above
$\tfrac12$, Lemma~\ref{lem:sv_forcing} gives $S(x):=\sum_{n\le x}\mu(6n)/n\ll x^{-1/2+\eps}$.

The passage back to $\sum_{n\le x}\mu(6n)$ is an Abel summation\index[terms]{Abel summation}, and
the constant it leaves has to be accounted for. Writing $T(x)=\sum_{n\le x}\mu(6n)$ and summing
$\mu(6n)=n\cdot\mu(6n)/n$ by parts,
\[
T(x)=x\,S(x)-\int_1^{x}S(t)\,dt .
\]
Had $S$ tended to a limit $\ell$, the first term would contribute $\ell x$ and would dominate
every power below one, so the estimate needs $\ell=0$. Here the bound on $S$ supplies it directly,
$S(x)\ll x^{-1/2+\eps}$ forcing $S(x)\to0$, and the same bound makes the integral
$\mathcal{O}(x^{1/2+\eps})$. Hence
\[
\sum_{n\le x}\mu(6n)\ll x^{1/2+\eps} .
\]
The vanishing is visible on the analytic side as well, the Dirichlet series below taking the value
$\prod_{p\ge5}(1-p^{-1})=0$ at $s=1$ by Mertens's theorem\index[names]{Mertens, F.}.

A second case runs the same way. The identity
$\sum_{k\le n}(-1)^{k-1}\mu(k)\lfloor n/k\rfloor=1+2\lfloor\log n/\log2\rfloor$ has a right hand
side that is increasing, slowly varying and of increments in $\{0,2\}$, so
Lemma~\ref{lem:sv_forcing} applies at the exponent one. The Abel step is the same, the partial
sums of $(-1)^{n-1}\mu(n)/n$ tending to zero because they are $\mathcal{O}(x^{-1/2+\eps})$, and it
returns $\sum_{n\le x}(-1)^{n-1}\mu(n)\ll x^{1/2+\eps}$.

Neither estimate is new. The first also follows from the Dirichlet series
\[
\sum_{n\ge1}\frac{\mu(6n)}{n^{s}}=\frac{1}{(1-2^{-s})(1-3^{-s})\zeta(s)}=\prod_{p\ge5}\bigl(1-p^{-s}\bigr),
\]
whose rightmost singularities under the hypothesis are the zeros of $\zeta$ on the critical
line\index[terms]{critical line}. What the route above changes is the input. The analytic route asks
for a Dirichlet series and its continuation, the route here asks for a counting formula and a rate,
and counting formulas of that shape are produced by elementary means in quantity.

Which arithmetic functions admit a counting formula whose count is slowly varying with bounded
increments is not answered here, and the two cases above are the only ones this volume records.
The question is where the theory would meet combinatorial number theory, and it is left open in
that form, without any claim that a subject is being founded.

\section{Intermediate depths and a positive index}\label{sec:intermediate_depths}

\subsection{A zero-free half plane}

The master equivalence $\mathrm{RH} \iff \alpha(\Phi) = 1/2$ admits the following refinement, which connects zero-free half planes of $\zeta$ to lower bounds on the transparency threshold and, under membership, to values of the regularity index. It is the second of the three levels of \S\ref{sec:boundary_to_strip}, between the boundary, which returns a limit, and the critical line.

\begin{corollary}\label{cor:zerofree}
Let $1/2<\sigma<1$. If the Riemann zeta function has no zeros in the half plane $\Re s>\sigma$,
then $\tau(\Phi)\ge1-\sigma$. Conversely, if $\Phi$ is a function of good variation with
$\alpha(\Phi)\ge1-\sigma$, then $\zeta$ has no zeros in that half plane.
\end{corollary}

\begin{proofstatus}{The corollary is proved as stated, but its two directions are not inverse to
one another. The first returns a
transparency threshold\index[terms]{transparency threshold} and leaves membership open, since it
delivers absorption\index[terms]{absorption} at the rate $1-\sigma$ and not at the rate
$\tau(\Phi)$. The second assumes membership. A biconditional would require the first direction to
produce a function of good variation, which it does not. The first direction is invoked in the
discussion following Conjecture~\ref{conj:fgv_membership}, and the converse is used in
Section~\ref{sec:commensurable} to compare $\alpha(\Phi)$ with $\eta(\Phi)$. No theorem assumes
the missing biconditional.}
\end{proofstatus}

\begin{proof}
Assume first that $\zeta$ has no zero in $\Re s>\sigma$.
Corollary~\ref{cor:perron-mobius-halfplane} gives
$M(x)\ll_{\sigma,\eps}x^{\sigma+\eps}$, with a finite Perron proof
for this precise hypothesis. The coefficient decomposition
of Theorem~\ref{thm:transfer_law} is applied to the exact forcing, so $w=0$.

For $0<\beta<1-\sigma$, choose
\[
 0<\delta\le\eps/4,\qquad
 \delta<(1-\sigma-\beta)/2,\qquad b=\sigma+\delta<1-\beta.
\]
The quotient $\zeta(s+\beta)/(s\zeta(s))$ has just the pole
$s=1-\beta$ between $b$ and $5/4$.
Lemma~\ref{lem:perron-zerofree-growth} supplies the subpower bounds
for both factors under the assumed zero-free half plane.
Corollary~\ref{cor:perron-subpower}, with error exponent $\eps/2$,
therefore gives
\[
 (1-\beta)\sum_{n\le x}J_{-\beta}(n)
 =\frac{x^{1-\beta}}{\zeta(1-\beta)}
      +\mathcal O_{\beta,\sigma,\eps}(x^{\sigma+\eps}).
\]
As before the formula first holds at half integers, and the main term
changes by only $\mathcal O_\beta(x^{-\beta})$ on replacing the
abscissa by any real $x$ with the same integer part.
Since $u(k)\ll_\beta k^{-\beta-1}$,
\[
 \sum_{n\le x}(\mu\star u)(n)
 \ll x^{\sigma+\eps}\sum_{k\le x}k^{-\beta-1-\sigma-\eps}
 \ll_{\beta,\sigma,\eps}x^{\sigma+\eps}.
\]
Thus $A_1(x)=x^{1-\beta}/\zeta(1-\beta)
+\mathcal O(x^{\sigma+\eps})$.
Taking $\eps<1-\sigma-\beta$, Abel summation and the Abelian
vanishing of the constant, proved in case (III) of
Theorem~\ref{thm:transfer_law}, give
\[
 A(x)=\frac{x^{-\beta}}{\Phi^*(\beta)}
       +\mathcal O_{\beta,\sigma,\eps}(x^{-(1-\sigma)+\eps}).
\]
Negative $\beta$ are covered by Lemma~\ref{lem:ingham_beta_neg},
and $\beta=0$ has the exact solution $(1,0,0,\ldots)$.
Every $\beta<1-\sigma$ is therefore transparent, which proves
$\tau(\Phi)\ge1-\sigma$.

For completeness, the absorption rate supplied by this hypothesis is
also explicit. Take $0<\eps<1-\sigma$. If $\beta\ge1-\sigma$, then
\[
 \sum_{n\le x}J_{-\beta}(n)
 =\sum_{k\le x}k^{-\beta}M(x/k)
 \ll x^{\sigma+\eps}\sum_{k\le x}k^{-\beta-\sigma-\eps}
 \ll_{\beta,\sigma,\eps}x^{\sigma+\eps},
\]
since $\beta+\sigma+\eps>1$. The $u$ term obeys the same bound.
Abel summation, again with the constant zero, gives
$A(x)\ll x^{-(1-\sigma)+\eps}$.
This is absorption at the rate $1-\sigma$, rather than necessarily
at $\tau(\Phi)$, so it does not by itself prove membership.

Conversely suppose that $\Phi$ is a function of good variation with
$\alpha(\Phi)\ge1-\sigma$. Corollary~\ref{cor:phi_half} gives
$\alpha(\Phi)\le1/2$, so $\beta=1$ is in the absorption range.
Its exact coefficients are $a_n=\mu(n)/n$. Hence
\[
 A(x)\ll_\eps x^{-\alpha(\Phi)+\eps},\qquad
 M(x)=xA(x)-\int_1^x A(t)\,dt
       \ll_\eps x^{1-\alpha(\Phi)+\eps}.
\]
Partial summation makes $\sum\mu(n)n^{-s}$ holomorphic in
$\Re s>1-\alpha(\Phi)$, where it continues $1/\zeta(s)$.
The zeta function has no zero there, and $1-\alpha(\Phi)\le\sigma$
gives the stated zero-free half plane.
\end{proof}

\begin{remark}
Corollary~\ref{cor:zerofree} places the threshold on a scale of fixed zero-free half planes. A
fixed half plane $\Re s>\sigma$ free of zeros, with $\sigma<1$, yields the lower bound
$\tau(\Phi)\ge1-\sigma$, and a better half plane yields a better bound. The scale is one of half
planes and not of regions, and the classical zero-free regions of $\zeta$, which narrow towards
the line $\Re s=1$ without leaving a half plane free, do not enter it.

In the other direction membership is what carries arithmetic information back. If $\Phi$ is a
function of good variation with $\alpha(\Phi)>0$, then $\zeta$ has no zeros in the fixed half
plane $\Re s>1-\alpha(\Phi)$, and the prime number theorem follows, a theorem that
Chapter~\ref{chap:boundary} obtains from the boundary with no hypothesis on the index. Such
membership is not presently known, and the endpoint $\alpha(\Phi)=1/2$ is equivalent to the
Riemann hypothesis.
\end{remark}

\subsection{What a positive index would imply}

Chapter~\ref{chap:boundary} obtained the prime number theorem from the boundary, with no
hypothesis on the index. A positive index would give it again, and much more, since by
Corollary~\ref{cor:zerofree} it gives a half plane free of zeros. The implication below is
recorded to place the prime number theorem among the weakest consequences of the index, and not
as a route to it.

\begin{proposition}\label{prop:pnt_fgv}
Suppose $\Phi$ is a FGV with $\alpha(\Phi) > 0$. Then $M(x)=o(x)$.
\end{proposition}

\begin{proof}
Use the finite forcing exponent $\beta=1$. Multiplying the defining equation by $n$ gives
\[
 \sum_{k\le n}k a_k\Big\lfloor\frac nk\Big\rfloor=1.
\]
The identity
\[
 \sum_{k\le n}\mu(k)\Big\lfloor\frac nk\Big\rfloor
 =\sum_{m\le n}\sum_{k\mid m}\mu(k)=1
\]
and triangular uniqueness give $na_n=\mu(n)$. Hence
$A(x)=\sum_{n\le x}\mu(n)/n$.

Only one case arises. Corollary~\ref{cor:phi_half} gives $\tau(\Phi)\le\tfrac12$ without any
hypothesis, so a positive index satisfies $0<\alpha(\Phi)\le\tfrac12$ and the exponent $\beta=1$
is always in the absorption regime. Absorption there with $\eps=\alpha(\Phi)/2$ gives
$A(x)=\mathcal{O}(x^{-\alpha(\Phi)/2})=o(1)$. No separate transparent case has to be examined.
Abel summation now gives
\[
 M(x)=xA(x)-\int_1^xA(t)\,dt=o(x),
\]
the integral being $o(x)$ by Ces\`aro averaging. No homogeneous equation and no passage to
$\beta=\infty$ is involved.
\end{proof}

\begin{remark}
The statement $M(x)=o(x)$ is one of the elementary equivalents of the prime number theorem, see
\nm{Tenenbaum}{G.}~\cite[\S~II.4]{Tenenbaum2015}. The proposition remains conditional only in the ordinary
sense that the hypothesis $\alpha(\Phi)>0$ is not presently known. The implication from that
hypothesis is proved above using a finite forcing exponent. The hypothesis gives a half plane free
of zeros, by Corollary~\ref{cor:zerofree}, while the conclusion needs only the line $\Re s=1$ and
is proved outright in Proposition~\ref{prop:pnt_boundary}.
\end{remark}
What the chapter settles is a single number. The transparency threshold of the Ingham weight does
not exceed one half, and that much costs no hypothesis. Under the hypothesis it is exactly one
half, with absorption beyond it, so the weight is a function of good variation of that index. The
converse returns the hypothesis from the same value. The two meet at one half, and the meeting is
Theorem~\ref{thm:tauberian_rh}.

What the chapter does not settle is as clear. Membership of $\Phi$ in the class is not proved
here, and it is the exact value of the index, not membership by itself, that the theorem makes
equivalent to the hypothesis. No zero of $\zeta$ is located and none is evaluated, the argument
producing one and never naming it. The same route, run on the character twisted weight, reaches
the two zero sets of $\zeta$ and $L(\cdot,\chi_4)$ at once and asserts neither. The equivalence
moves the question. It does not answer it.

\chapter{Examples of functions of good variation}
\label{chap:examples}

A definition earns its place by what it separates. This chapter puts the definitions of the
preceding chapters to work on the two examples that mark out the possibilities, then surveys the
one variable weights of the gallery.

The first example is affine, continuous and monotone, with a rational transform and a single
zero. Everything about it is exact, the recurrence, its solution, and the three regimes below,
at and above the index. It is the archetype against which the rest of the theory is read, and
it returns in every later chapter as the case where an independent check is available.

The second is the broken harmonic weight\index[terms]{broken harmonic function} at scale two, discontinuous on a lattice of dyadic
points. Its index is again the first zero of its transform, but the individual terms behave
differently, and the Hardy-Littlewood-Ramanujan\index[terms]{Hardy--Littlewood--Ramanujan criterion}\index[names]{Ramanujan, S.} criterion fails at the exponent zero, where the
partial sums against the coordinate are the binary digit sum. Between them the two examples
show that the index can be read off the transform in the continuous case and in the
discontinuous one alike, while the behavior of the individual terms is not determined by the
transform at all.

The gallery that follows is not decoration. The general statements of this volume were found on
its entries before they were stated in general, and each entry is worked constructively rather
than by the shortest route available.

\section{Two canonical examples}
\label{sec:canonical_examples}

Two examples, worked before the gallery, mark out the possibilities. The first is continuous and elementary, the second discontinuous and more delicate, and between them they show smooth transparency, broken harmonics, and the part played by the HLR criterion\index[terms]{Hardy--Littlewood--Ramanujan criterion}.

\subsection*{Example 1: The affine function $g(x) = (1-\lambda)x + \lambda$}

The simplest function of good variation beyond a constant. Its index is $\lambda$, the unique zero of its Mellin transform\index[terms]{Mellin transform}, and it serves as the archetype of the theory, continuous and monotone on $(0,1]$ with $g(1)=1$ and a rational Mellin transform.

\begin{theorem}\label{thm:affine_body}
Let $\lambda \in (0,1)$ and $g(x) = (1-\lambda)x + \lambda$. Then $g$ is a FGV with regularity index $\alpha(g) = \lambda$ and Mellin transform\index[terms]{Mellin transform}
\[
g^*(z) = \frac{z - \lambda}{z - 1}.
\]
The unique zero of $g^*$ is $z = \lambda = \alpha(g)$, so $\alpha(g) = \eta(g)$. More precisely, for any sequence $(a_n)$ satisfying
\[
\sum_{k=1}^n a_k\left((1-\lambda)\frac{k}{n} + \lambda\right) = n^{-\beta},
\]
the partial sums $A(n) = \sum_{k=1}^n a_k$ satisfy:
\begin{itemize}
\item if $\beta < \lambda$: $A(n) \sim \dfrac{1}{g^*(\beta)}\,n^{-\beta}$ \quad (transparency);
\item if $\beta = \lambda$: $A(n) \sim (1-\lambda)\,n^{-\lambda}\log n$ \quad (resonance);
\item if $\beta > \lambda$: $A(n) \sim C\,\dfrac{\Gamma(n+1-\lambda)}{\Gamma(n+1)} \sim C'\,n^{-\lambda}$ \quad (absorption).
\end{itemize}
\end{theorem}

\begin{proof}
Write $S(n)=\sum_{k\le n}k\,a_k$. Multiplying the defining equation\index[terms]{defining equation} by $n$ gives
$(1-\lambda)S(n)+\lambda\,nA(n)=n^{1-\beta}$, and subtracting the same identity at rank $n-1$
removes $S$, since $S(n)-S(n-1)=n\,a_n=n\bigl(A(n)-A(n-1)\bigr)$. What is left is the exact
first order recurrence
\[
A(n)=\Bigl(1-\frac\lambda n\Bigr)A(n-1)+\frac{n^{1-\beta}-(n-1)^{1-\beta}}{n},
\qquad A(1)=1 .
\]
Integrating it by variation of constants\index[terms]{variation of constants} gives a product that the Euler-Gauss formula
evaluates as $\Gamma(n+1-\lambda)/\bigl(\Gamma(2-\lambda)\Gamma(n+1)\bigr)$, of order
$n^{-\lambda}$, against a sum whose generic term is of order $k^{\lambda-\beta-1}$. The three
ranges of $\beta$ are the three behaviors of that sum, divergent, logarithmic and convergent.
The computation is carried out in full in Appendix~\ref{app:A}, where the constants are
identified.
\end{proof}

\begin{remark}
The transparency constant $1/g^*(\beta) = (1-\beta)/(\lambda-\beta)$ is positive for $\beta < \lambda$ (since $g^*(\beta) > 0$) and blows up as $\beta \to \lambda^-$, announcing the resonance. The absorption\index[terms]{absorption} mode $n^{-\lambda}$ is the Gamma-Euler mode attached to the zero of $g^*$, and the spectral information encoded in $g^*$ governs the arithmetic behavior of $(a_n)$ in both cases.
\end{remark}

\subsection*{Example 2: The broken harmonic function $g_2(x) = x \cdot 2^{\lfloor -\log_2 x \rfloor}$}

This is the first discontinuous example. The function $g_2$ is piecewise linear on dyadic intervals, self-similar, and discontinuous at every negative power of 2. Despite these discontinuities it is expected to be a FGV, and the dyadic recurrence argument below, which has no continuous analogue, carries the exact part of that picture.

\begin{definition}
For $x \in (0,1]$, define
\[
g_2(x) := x \cdot 2^{\lfloor -\log_2 x \rfloor}.
\]
Equivalently, $g_2(x) = x \cdot 2^j$ for $x \in (2^{-j-1}, 2^{-j}]$, $j \ge 0$, so on each dyadic interval $g_2$ increases linearly from just above $\tfrac12$ to $1$, attains the value $1$ exactly at the points $x=2^{-j}$, and drops to the right hand limit $\tfrac12$ across each of them, oscillating between $\tfrac12$ and $1$.
\end{definition}

Its Mellin transform\index[terms]{Mellin transform} is
\begin{equation}\label{eq:g2_mellin}
g_2^*(z) = \frac{z}{z-1}\cdot\frac{2^{z-1}-1}{2^z - 1},
\end{equation}
which has zeros at $z = 1 + 2\pi i k/\log 2$ for $k \in \Z \setminus \{0\}$, all on the line $\Re(z) = 1$.

The transform of the dyadic kernel has all its zeros on one vertical line.

\begin{proposition}\label{prop:brokenharm_body}
The transform of $g_2$ is \eqref{eq:g2_mellin}, with zeros $z=1+2\pi ik/\log2$ for
$k\in\Z\setminus\{0\}$, all on the line $\Re(z)=1$, so $\eta(g_2)=1$. The weighted
sums $A_1(n)=\sum_{k\le n}ka_k$ satisfy the exact dyadic recurrence
\[
A_1(n) = A_1\!\left(\left\lfloor\frac{n}{2}\right\rfloor\right) + w(n),
\qquad w(n) = n^{1-\beta} - 2\lfloor n/2\rfloor^{1-\beta}.
\]
At the critical forcing $\beta=1$ the recurrence telescopes exactly to
$A_1(n)=1-\lfloor\log_2 n\rfloor$, and $A(n)=2^{-\lfloor\log_2 n\rfloor}$, so that
$nA(n)$ oscillates with $\limsup_{n\to\infty} nA(n)=2$ and
$\liminf_{n\to\infty} nA(n)=1$.
\end{proposition}

\begin{proof}
The recurrence follows by splitting the defining relation along the self-similar structure
of $g_2$ at $\lfloor n/2\rfloor$. The critical case is the instance $\lambda=2$ of
Appendix~\ref{app:D}, Lemma~\ref{lem:D_critical}.
\end{proof}

The index of that kernel follows, together with its three regimes.

\begin{theorem}\label{thm:g2_index}
The function $g_2$ is a function of good variation with
$\alpha(g_2)=1=\eta(g_2)$. For $\beta<1$ the exponent is transparent, with
$A(n)\sim n^{-\beta}/g_2^{*}(\beta)$, and for $\beta>1$,
$A(n)=\mathcal{O}(n^{-1+\eps})$ for every $\eps>0$.
\end{theorem}

\begin{proof}
The dyadic recurrence and the critical value are the content of
Proposition~\ref{prop:brokenharm_body}. The transparency and absorption\index[terms]{absorption} branches pass
from the weighted sums to $A(n)$ through the finite Abel identity\index[terms]{Abel identity} of
Appendix~\ref{app:D}, whose terminal constant is computed in
Lemma~\ref{lem:D_telescope}, and both then follow from Theorem~\ref{thm:D_index}
at $\lambda=2$, as does $\alpha(g_2)=1=\eta(g_2)$.
\end{proof}

At the exponent zero the weighted sums of the dyadic kernel are a familiar arithmetic function,
and that alone decides the criterion there.

\begin{proposition}
\label{prop:g2_digitsum}
At $\beta=0$ the solution attached to $g_2$ satisfies
\begin{equation}\label{eq:g2_digitsum}
\sum_{k\le n}k\,a_k=s_2(n),
\end{equation}
where $s_2(n)$ is the number of ones in the binary expansion of $n$. Hence
$na_n=s_2(n)-s_2(n-1)$, so that $na_{2^j}=-(j-1)$ along the powers of two and the
sequence $(na_n)$ is unbounded. The strong Hardy-Littlewood-Ramanujan criterion fails
for $g_2$ at $\beta=0$.
\end{proposition}

\begin{proof}
Write $b_k=ka_k$ and $S(m)=\sum_{k\le m}b_k$. Since
$g_2(x)=x\,2^{\lfloor-\log_2x\rfloor}$, the defining equation\index[terms]{defining equation} at $\beta=0$ reads
\[
\sum_{k\le n}b_k\,2^{\lfloor\log_2(n/k)\rfloor}=n .
\]
For $t\ge1$ the geometric sum $\sum_{j\ge0,\,2^j\le t}2^j=2^{\lfloor\log_2t\rfloor+1}-1$
gives $2^{\lfloor\log_2t\rfloor}=\tfrac12\bigl(1+\sum_{2^j\le t}2^j\bigr)$. Applying it
to $t=n/k$ and exchanging the two summations, the equation becomes $T(n)+S(n)=2n$, where
\[
T(n)=\sum_{j\ge0}2^{j}S\bigl(\lfloor n/2^{j}\rfloor\bigr).
\]
Splitting off the term $j=0$ and using
$\lfloor\lfloor n/2\rfloor/2^{i}\rfloor=\lfloor n/2^{i+1}\rfloor$ gives
$T(n)=S(n)+2T(\lfloor n/2\rfloor)$. Eliminating $T$ between the two relations leaves
\[
S(n)=S\bigl(\lfloor n/2\rfloor\bigr)+\bigl(n-2\lfloor n/2\rfloor\bigr),
\qquad S(0)=0,
\]
which is the recursion of the binary digital sum, so $S=s_2$. Differencing gives
$na_n=s_2(n)-s_2(n-1)$, and at $n=2^{j}$ one has $s_2(2^{j})=1$ and $s_2(2^{j}-1)=j$,
whence $na_{2^{j}}=1-j$.
\end{proof}

Identity \eqref{eq:g2_digitsum} is proved above and is self-contained. It is also, as the
research dossier~\ref{app:dossier_sqrt2} records much later, the case $\lambda=2$ of an
endpoint identity holding at every integer scale, and nothing in the present chapter depends
on that reading. On the positive range the criterion holds, and the same dyadic recurrence
delivers it in closed form with an explicit constant.

\begin{theorem}\label{thm:g2_hlr}
Let $\beta$ be complex with $\Re\beta>0$. The solution attached to $g_2$ satisfies, at every
rank,
\begin{equation}\label{eq:g2_hlr_bound}
|n\,a_n|\;\le\;1+|1-\beta|\left(1+\frac{2^{\Re\beta}}{2^{\Re\beta}-1}\right).
\end{equation}
The strong Hardy-Littlewood-Ramanujan criterion\index[terms]{Hardy--Littlewood--Ramanujan criterion} therefore holds for $g_2$ on the whole positive
range.
\end{theorem}

\begin{proof}
Write $A_1(n)=\sum_{k\le n}ka_k$ and $L=\lfloor\log_2n\rfloor$. The dyadic recurrence of
Proposition~\ref{prop:brokenharm_body} reads
$A_1(n)=A_1(\lfloor n/2\rfloor)+n^{1-\beta}-2\lfloor n/2\rfloor^{1-\beta}$. Iterating it along
the ranks $n_i=\lfloor n/2^{i}\rfloor$ down to $n_L=1$, where $A_1(1)=1$, telescopes the two
families against each other and leaves the closed form
\begin{equation}\label{eq:g2_closed}
A_1(n)=n^{1-\beta}-\sum_{i=1}^{L}\Bigl\lfloor\frac{n}{2^{i}}\Bigr\rfloor^{1-\beta} .
\end{equation}
Put $\varphi(m)=m^{1-\beta}-(m-1)^{1-\beta}$ for $m\ge2$ and $\varphi(1)=1$. Along the
positive real axis, the fundamental theorem of calculus gives
\[
\varphi(m)=(1-\beta)\int_{m-1}^{m}t^{-\beta}\,dt,
\qquad
|\varphi(m)|\le |1-\beta|\int_{m-1}^{m}t^{-\Re\beta}\,dt
\le |1-\beta|(m-1)^{-\Re\beta}.
\]
The bound is immediate for $n=1$. Assume henceforth that $n\ge2$. Differencing
\eqref{eq:g2_closed} at the ranks $n$ and $n-1$ uses only
$\lfloor n/2^{i}\rfloor-\lfloor(n-1)/2^{i}\rfloor=\mathbf 1_{\{2^{i}\mid n\}}$, so with
$n=2^{v}m_0$ and $m_0$ odd only the ranks $i\le v$ survive.

When $m_0\ge3$ the two ranks $n$ and $n-1$ carry the same $L$, and
$na_n=\varphi(n)-\sum_{i=1}^{v}\varphi(m_02^{v-i})$. Setting $j=v-i$ gives
$m_02^{j}-1\ge2^{j+1}$, so the surviving terms are bounded by
$|1-\beta|2^{-(j+1)\Re\beta}$ and their sum by
$|1-\beta|/(2^{\Re\beta}-1)$, which together with $|\varphi(n)|\le|1-\beta|$ gives
\eqref{eq:g2_hlr_bound} with room to spare.

When $n=2^{L}$ the rank $n-1$ carries one dyadic level less, the term $i=L$ of
\eqref{eq:g2_closed} is unmatched and contributes $-1$, and
$na_n=\varphi(n)-\sum_{i=1}^{L-1}\varphi(2^{L-i})-1$. With $j=L-i$, the inequality
$2^{j}-1\ge2^{j-1}$ for $j\ge1$ bounds the sum by
$|1-\beta|2^{\Re\beta}/(2^{\Re\beta}-1)$, and \eqref{eq:g2_hlr_bound}
follows. That second bound dominates the first, so it holds at every rank.
\end{proof}

\begin{remark}\label{rem:g2_hlr_endpoints}
The closed form \eqref{eq:g2_closed} carries both endpoints of the criterion. At $\beta=1$ it
reduces to $A_1(n)=1-\lfloor\log_2n\rfloor$, so $a_1=1$ and
$na_n=-\mathbf 1_{\{n\text{ is a power of }2\}}$ for $n\ge2$. At $\beta=0$ it degenerates into the digital sum
\eqref{eq:g2_digitsum}, each floor contributing the parity of its rank, and the bound
\eqref{eq:g2_hlr_bound} degrades like $1/(\beta\log2)$ as $\beta$ decreases to zero. That
degradation is forced, since $na_{2^{L}}$ tends to $1-L$ at the exponent zero by
Proposition~\ref{prop:g2_digitsum}. The success on the positive range and the failure at the
origin therefore come out of one identity.
\end{remark}

\begin{remark}
Two structural differences stand out. First, the index $\alpha(g_2)=1$ is tied to no single zero, being determined by the spectral geometry of the zeros of $g_2^*$, all on $\Re(z)=1$. Second, the critical case $\beta=1$, which is exact, produces not a logarithmic resonance but an oscillation, the limit does not exist and $nA(n)$ oscillates between 1 and 2. This is a signature of broken harmonics that has no continuous counterpart.

Both examples have $\alpha=\eta$ proved, and for $g_2$ the HLR clause is proved as well by Theorem~\ref{thm:g2_hlr}, so that entry carries no conjectural clause. The one worked kernel where the two indices are proved distinct appears in the research dossier~\ref{app:dossier_sqrt2} with $g_{\sqrt{2}}$, where the proved separation places the transparency frontier\index[terms]{transparency frontier} at or below $1/2$, strictly below $\eta=1$, the conjectural value $\alpha=1/2$ resting on the open LOW and ABS estimates, driven by the algebraic resonance $(\sqrt{2})^2 = 2 \in \N$.
\end{remark}

\section{A gallery of one-variable kernels}
\label{sec:zoo_table}

The seven one-variable profiles of Appendices~\ref{app:A} to~\ref{app:G} are collected below,
each with its regularity index $\alpha$, its analytic index $\eta$, and a short remark. Six are
proved to be functions of good variation, while the last is conjectural at its index. Two of them,
\ref{app:B} and~\ref{app:F}, are unbounded at the origin and belong to the recurrence-admissible
extension fixed in \S\ref{sec:fgv_slowly_varying}. Two further profiles, the fractional part\index[terms]{fractional part}
kernel of Appendix~\ref{app:O} and the affine family of Appendix~\ref{app:P}, are long enough
to be read on their own and are not repeated here. A reader new to the theory may look over a few
entries before turning to the analytic theory. The two-variable kernels are gathered in the gallery
of Chapter~\ref{chap:raf}. The broken harmonic kernel at $\sqrt2$ is not a gallery entry, it is the
research dossier of Chapter~\ref{chap:diophantine}.

In the table each kernel appears in canonical form as a profile $g(x)$, and the two indices are
listed with their proof status, open entries remaining open.

\medskip
\renewcommand{\arraystretch}{1.45}
{\small
\begin{tabular}{cp{3.4cm}p{2.4cm}p{2.0cm}p{3.4cm}}
\toprule
{App.} & {Kernel} & $\alpha$ & $\eta$ & {Remark} \\
\midrule
A & $g(x)=(1-\lambda)x+\lambda$ & $\lambda$ & $\lambda$ & Affine archetype. $\alpha=\eta$. Exact Gamma-Euler solution. \\
B & $g(x)=1-\lambda\log x$ & $\lambda$ & $\lambda$ & Logarithmic recurrence. Sharpness through the connection coefficient. Unbounded at $0$ on purpose, the singularity being slowly varying. \\
C & $g(x)=\mu$ on $(0,\lambda]$, $1$ on $(\lambda,1]$ & $\log(1-\mu)/\log\lambda$ & same & Bipolar step. Dilative recurrence with the terminal ranks controlled, two sided bound above the threshold. \\
D & $g(x)=x\lambda^{\lfloor-\log_\lambda x\rfloor}$, $\lambda\in\Z_{\ge2}$ & $1$ & $1$ & Broken harmonic, integer base. A telescoping identity closes the return to $A(n)$. \\
E & $g(x)=|x-\frac12|+\frac12$ & $\eta$ & $\eta$ & Zeros of $2^z-z$ ordered by an exact phase argument. Indicial equation $2^\gamma=\gamma$. The negative band, the expansion $c_\beta n^{-\beta}+2\Re(\kappa n^{-z_0})+\mathcal O(n^{-7/8})$, and sharpness from the meromorphic pole of $\kappa$ are proved without numerical input. \\
F & $g(x)=x-\log x$ & $1/2$ & $1/2$ & Zeros at $\tfrac12\pm\tfrac{\sqrt3}{2}i$. Exact recurrence of order two whose indicial polynomial is the numerator of $g^{*}$, factorized by a discrete Riccati equation. Unbounded at $0$, slowly varying. \\
G & $g(x)=1-x/2$ on $(0,\frac12)$, $x$ on $[\frac12,1]$ & open & $2$ & Every zero of the numerator on the line $\Re z=2$, an infinite family given by an explicit phase equation. Scalar reduction with a parity term, two parity multipliers of mean $-P$. \\
\bottomrule
\end{tabular}
\renewcommand{\arraystretch}{1}
}

\medskip
\noindent The appendices carry a proof of the stated dichotomy except where a
\texttt{proofstatus} box marks an open gap. Entries \ref{app:A} to \ref{app:F} are
proof-complete, and \ref{app:G} remains open at its regimes. The coincidence $\alpha=\eta$ is recorded
only where it is proved. Appendices~\ref{app:B} and~\ref{app:F} use kernels unbounded at $0$
deliberately. Their logarithmic singularity is slowly varying\index[terms]{slowly varying}, and their exact recurrences lie in
the recurrence-admissible class introduced above. They are proved instances of the slowly varying
branch of Conjecture~\ref{conj:fgv_membership}. The power profiles of
Theorem~\ref{thm:fgv_power_singularity} form a further unbounded family beyond slow variation,
with the boundary index $\alpha=\eta=0$.
That closes the first part. Five chapters have taken a weight out of a summation method, made it
the kernel of an equation, attached to it a threshold and a transform, read it at the boundary of
that transform, where it returns the prime number theorem, and shown that for one weight the value of the index is the Riemann hypothesis. The gallery adds seven profiles where
the threshold is computed by hand, six of them completely. What none of it supplies is a
criterion. Nothing so far says which weights have a finite threshold with absorption above it,
and the next part takes
that question up by widening the kernel rather than by narrowing the class.

\part{From functions of good variation to regular arithmetic functions}
\label{part:two}

\rafepigraph{Il y avait loin encore, cependant, de l'arithmétique, où règne le discontinu, à la théorie des fonctions au sens classique.}{There was still a long way, however, from arithmetic, where the discontinuous reigns, to the theory of functions in the classical sense.}{André Weil, \emph{De la métaphysique aux mathématiques} (1960)~\cite{WeilMetaphysique1960}}

Part~\ref{part:one} fixed one number on one weight. What it did not say is what that number
belongs to, and the answer is not the weight but the equation.

Chapter~\ref{chap:raf} widens the kernel from a function of the single ratio $k/n$ to a function of
$n$ and $k$ apart, which is the regular arithmetic function of the title, and builds the class
with an existence theorem that produces an index from the kernel alone.
Chapter~\ref{chap:principles} assembles what makes membership decidable.
Chapters~\ref{chap:volterra} and~\ref{chap:discrete_volterra} invert the operator, the first writing the
Ingham operator as a Volterra equation on the half line and reading the equivalence of
Part~\ref{part:one} a second time through its resolvent, the second building the resolvent of the
discrete equation in the generality the rest of the volume uses.

\chapter{Regular arithmetic functions}
\label{chap:raf}

The kernels of Part~I were functions of a single ratio, $g(k/n)$, and for that form the
arithmetic average of the kernel was an integral, its transform the classical Mellin transform
of $g$. The theory does not end at that form. A general averaging couples the two indices $n$
and $k$ without folding them into one ratio, the integral behind the transform is gone, and the
object that decided the index in Part~I, the location of a zero, may fail to exist. What
survives the passage is not a representation but a property, the turning of the partial sums at
a single exponent.

The chapter opens by setting these equations beside the classical ones they resemble, and by
naming what the theory takes as primary that those do not. It then fixes that property as the
definition of a regular arithmetic function, the class from which the theory takes its name, and
reads the index on a kernel whose transform is constant, where no zero is available and the index
is arithmetic through and through. From there it turns to existence in general, a counterexample
first, then a theorem that grants the index to a whole class of kernels at once, and two families
that theorem delivers, the M\"untz sums and a lacunary profile whose transform stops at a natural
boundary. The last four sections carry machinery and inventory. The scalar reduction turns the
equation into a first order recurrence, an algebra combines indices, a conditional theorem states
the contour argument that would carry the analytic index to the arithmetic one together with every
hypothesis it needs, and a gallery collects the two-variable kernels the appendices prove.

\section{Classical ingredients and a distinct primary object}\label{sec:primary_object}

The equations considered in this book are not new merely as equations. For each fixed exponent
$\beta$, the defining lower triangular system is a discrete Volterra equation\index[terms]{Volterra equation} of the first kind,
and its solution by forward substitution or through a Green matrix\index[terms]{Green's matrix} belongs to classical difference
equation and Volterra theory\index[terms]{Volterra operator}. When $G(n,k)=g(k/n)$, Abel summation\index[terms]{Abel summation} also places the problem next to
Mellin convolution, regular variation\index[terms]{regular variation} and the inversion theorems of Tauberian theory\index[terms]{Tauberian theorem}. In special homogeneous cases a
renewal type formulation may likewise be available. These connections are real and they are used
throughout the book.

Several forcing exponents are compared in what follows, so the subscript is kept on $A_\beta$ and on the sequences attached to it. Elsewhere in the volume the exponent is fixed by its context and the plain $a_n$ and $A(n)$ of Definition~\ref{def:reg_index_fgv} are used.

The primary object extracted from the equations is nevertheless different. The kernel $G$ is fixed
while the complete family of power forcings
\[
\sum_{k\le n}a_{\beta,k}\,G(n,k)=n^{-\beta}
\]
is allowed to vary with $\beta$, and the observable is the ordinary partial sum
$A_\beta(n)=\sum_{k\le n}a_{\beta,k}$. The question is not only whether one fixed solution has a
prescribed asymptotic behavior. It is to determine the maximal interval of exponents on which the
partial sums reproduce the imposed power, and whether all responses beyond that interval are
absorbed\index[terms]{absorption} at a common limiting rate. The resulting number $\alpha(G)$ is a global response index of
the normalized triangular kernel and of its power forcing scale.

Nor is $\alpha(G)$ defined as a zero of a transform. The analytic index\index[terms]{analytic index}
$\eta(G)=\inf\{\Re\rho:G^{*}(\rho)=0\}$ is introduced separately, so that an equality
$\alpha(G)=\eta(G)$, wherever it holds, is a theorem and not a convention. The distinction is
substantive. The rational kernel of Chapter~\ref{chap:raf} has $G^{*}(z)\equiv1$, hence no analytic
index at all, and yet has the arithmetic index $\alpha(G)=2$. In the other direction, for the broken
square root kernel of Chapter~\ref{chap:diophantine} the transparency frontier\index[terms]{transparency frontier} lies at or below
$\tfrac12$ while the analytic index is $1$. General two-variable kernels $G(n,k)$ need not possess
any one-variable Mellin convolution representation at all.

The novelty claimed is therefore limited and precise. Thresholds, indicial roots\index[terms]{indicial root}, resonance,
resolvents\index[terms]{resolvent}, regular variation\index[terms]{regular variation} and Tauberian transfer\index[terms]{Tauberian} are classical. What appears not to have been
isolated in the preceding theories is their organization into a single characteristic of a general
arithmetic triangular kernel, a maximal range of transparent power forcings, supplemented by an
absorption law beyond the transition, and compared with rather than defined by a separate analytic
zero index. No exact antecedent for this combined construction has been identified.

The nearest neighbour is the Tauberian theory of \nm{Bingham}{N. H.} and \nm{Inoue}{A.} for
arithmetic sums and for systems of kernels\index[terms]{Tauberian theorem} \cite{BinghamInoue2000a,BinghamInoue2000b}, and who study the
kernel $\lfloor x\rfloor/x$ itself. Their question is the transfer of regular variation\index[terms]{regular variation} at a fixed
index under a non-vanishing hypothesis on the transform. The question here is the maximal range of
indices over which that transfer holds, together with what happens beyond it, and the two are
compared rather than identified.

Where the two part company is visible in their own account of the kernel. The transform
$\zeta(z+1)/(z+1)$ converges absolutely on $\Re z>0$, which in the parameter of this volume is the
range $\beta<0$, and there the non-vanishing hypothesis costs nothing, the Euler product and the
classical zero-free region keeping $\zeta(1+s)$ away from zero. On that range the two accounts
agree, and \S\ref{sec:beta_negative} settles it here without complex analysis. Beyond it they
record the obstruction themselves, the best known zero-free region for $\zeta$ containing
$\Re s\ge1$ but no half plane $\Re s>1-\eps$, which is what their method would ask for at the
boundary \cite[Remark~4, p.~182]{BinghamInoue2000b}.

Their device against the zeros is the one this volume meets from the other side. A single scale
$\lambda$ puts the zeros of the transform on a vertical line in arithmetic progression,
$\rho+2\pi in/\log\lambda$, which defeats the hypothesis their theorems need, and they escape it
by taking two scales whose logarithms are incommensurable, so that the two progressions meet only
at the index. They attribute that condition to Ingham's\index[names]{Ingham, A. E.} method
\cite[Remark~2, p.~182]{BinghamInoue2000b}. The same progression is what
Chapter~\ref{chap:diophantine} meets with the scale given rather than chosen, where the
arithmetic of $\log\lambda$ is not an obstacle to be avoided but the object under study, and where
the broken kernel at an irrational scale is the case that separates the two indices.

Neither reading attaches a number to the kernel. A Mercerian theorem concludes at one index, so a
hypothesis on the zeros of $\zeta$ would do no more than widen the domain on which such a theorem
can be stated. Here the index is itself the object, the boundary of the range on which the
transfer holds, and it is that boundary that the hypothesis decides.

The comparison with the closure criteria of Nyman\index[names]{Nyman, B.} and Beurling\index[names]{Beurling, A.} is made at the opening of
Chapter~\ref{chap:equivalence}, where the equivalence it stands beside is proved. Two
consequences of that comparison belong here instead.

That writing also places one entry of the gallery beside the Ingham kernel. The profile of
Appendix~\ref{app:O} is $g(x)=1-\{x\}$, the affine function $1-x$ on $(0,1)$ carrying a single
jump at the diagonal, and its transform $1/(1-z)$ has no zero at all. The two are the same
syntax with the fractional part taken of $x$ in one case and of $1/x$ in the other, which turns
one discontinuity into a discontinuity at every $1/m$. Neither index is handed over by the zero
set of the transform. For the first there is no zero to read, and the index $\tfrac14$ comes from
the single jump. For the second the reading is available and its correctness is the Riemann
hypothesis.

The two shapes meet once inside this volume, in Chapter~\ref{chap:ortho}. The orthorecursive
expansion\index[terms]{orthorecursive expansion} of unity is an approximation problem in
$L^{2}([0,1])$, a distance driven down by a greedy projection, and it is settled here by reading
the zeros of the arithmetic Mellin transform\index[terms]{arithmetic Mellin transform} of the
kernel $2/(1+t)$, located and confined completely in that chapter. That is the shape of a closure criterion, an $L^{2}$ distance decided by
a zero set, and it is decidable there because the zero set is available in closed form. For the
Ingham kernel\index[terms]{Ingham kernel} the shape is the same with the zeros of $\zeta$ in place
of those of $2/(1+t)$, which is what puts the two families in one room and what keeps them apart.

In short, RAF theory has the triangular syntax of a difference equation and the resolvent mechanism
of Volterra theory, and where the arithmetic Mellin transform\index[terms]{arithmetic Mellin transform} of the kernel does not vanish it
inherits the conclusions of Tauberian theory\index[terms]{Tauberian}. Its primary object is the phase
transition index of the whole forced family.

\section{The defining equation and the regularity index}

Fix a function $G:\N^{*}\times\N^{*}\to\R$ with $G(n,n)\neq0$ for every $n$, and a function
$r:\N^{*}\to\R$. The defining equation
\begin{equation}\label{eq:defining_relation_gen}
\sum_{k=1}^{n} a_k\,G(n,k)=r(n)\qquad(n\ge1)
\end{equation}
determines $(a_n)_{n\ge1}$ uniquely, by the same isolation of the last term as in the
one-variable case, the relation at $n$ fixing $a_nG(n,n)$ from $a_1,\dots,a_{n-1}$. The forcing
$r(n)=n^{-\beta}$ is again the case carried forward, and $A(n):=\sum_{k=1}^{n} a_k$ is the
partial sum under study.

The transform of Part~I keeps its arithmetic definition here, where no integral is available.

\begin{definition}\label{def:mellin}
For $\Re(z)<0$, when the limit exists, the arithmetic Mellin transform of $G$ is
\[
G^{*}(z):=\lim_{n\to\infty}\Big(-\tfrac{z}{n}\Big)\sum_{k=1}^{n} G(n,k)\Big(\tfrac kn\Big)^{-z-1},
\]
continued where possible beyond $\Re(z)=0$. Its analytic index is
$\eta(G)=\inf\{\Re(\rho):G^{*}(\rho)=0\}$.
\end{definition}

For a kernel $G(n,k)=g(k/n)$ this limit is the Riemann sum of
Proposition~\ref{prop:mellin_coincidence}, and $G^{*}$ is the classical Mellin transform of
$g$, carried by an integral. For a general $G$ no such integral is available, and $G^{*}$ is
fixed by the arithmetic limit alone. It may be constant, with no zero to locate, and then the
analytic index $\eta(G)$ is undefined.

\begin{definition}[Regular arithmetic function and regularity index]\label{def:reg_index}
The exponent $\beta$ is transparent for the kernel $G$ when there is a constant $\Xi_G(\beta)$
with
\begin{equation}\label{eq:transparency_additive}
A(n)=\Xi_G(\beta)\,n^{-\beta}+o(n^{-\beta}),
\end{equation}
a constant that is then unique, by the argument of Definition~\ref{def:transparency_fgv}. The
transparency threshold\index[terms]{transparency threshold} of $G$ is
\[
\tau(G)=\sup\big\{\,a\in\R:\ \text{every }\beta<a\ \text{is transparent for }G\,\big\},
\]
defined for every kernel whose triangular inverse exists. The kernel is a regular arithmetic
function\index[terms]{regular arithmetic function}, or RAF, when $\tau(G)$ is finite and the partial sums are absorbed above it,
$A(n)=\mathcal{O}(n^{-\tau(G)+\eps})$ for every $\beta\ge\tau(G)$ and every $\eps>0$, and its
regularity index is then $\alpha(G)=\tau(G)$. As in Definition~\ref{def:reg_index_fgv}, the notation
$\alpha(G)$ is reserved for that case, so that writing $\alpha(G)$ asserts that $G$ is a regular
arithmetic function.
\end{definition}

The transparent coefficient\index[terms]{transparent coefficient} is computed by the transform. Wherever the identification is
established, and it is established in this chapter for each kernel where it is used, it reads
\[
\Xi_G(z)=\frac{1}{G^{*}(z)},
\]
read as a meromorphic reciprocal. At a pole of $G^{*}$ the value of $\Xi_G$ is zero, and
transparency at such a point asserts $A(n)=o(n^{-\beta})$. At a zero of $G^{*}$ the
reciprocal is not finite and the point is not transparent. Open
Problem~\ref{op:xi_reciprocal} asks which hypotheses can replace the weighted mass condition,
Proposition~\ref{prop:xi_mass_needed} showing that the continuation hypotheses do not suffice on
their own, and the two statements are kept apart for that reason. The identification bears on the transparent range,
the open half line below the threshold, and it says nothing at the threshold itself, where the
definition records absorption. The binomial harmonic\index[terms]{binomial harmonic kernel} kernel of
Appendix~\ref{app:Q} is where that distinction is visible, since it satisfies
$\Xi_G=1/G^{*}$ at every exponent below its index and follows an exact power with the different
constant $1$ at the index. The identification is a theorem
where it is proved and nowhere a convention. Equation~\eqref{eq:transparency_additive} is the form
in which the proofs of this chapter verify transparency.

\begin{remark}\label{rem:index_three_conditions}
The supremum definition unfolds into three conditions, and the proofs below establish
exactly these. Transparency holds at every $\beta<\alpha(G)$. Absorption holds at and
above the index. And the index is sharp, every interval $[\alpha(G),\alpha(G)+\delta)$
with $\delta>0$ contains a nontransparent exponent, since a wholly transparent interval of
that form would raise the supremum. Conversely a number satisfying the three conditions
equals the supremum, so the two descriptions agree.
\end{remark}

At $\beta=\alpha(G)$ the second estimate applies and may carry a logarithmic correction. A
kernel whose partial sums follow the forcing at every real $\beta$ has $\alpha(G)=+\infty$. A function of good variation is a regular arithmetic function of one particular type, the type
in which the kernel is induced by a profile, $G(n,k)=g(k/n)$ with $g$
recurrence-admissible\index[terms]{recurrence-admissible profile} in the sense of
Chapter~\ref{chap:fgv}. Nothing further separates the two classes. Admissibility asks that the
profile be finite at every point of $(0,1]$ and nonzero at $1$, and the two logarithmic profiles
of Appendices~\ref{app:B} and~\ref{app:F}, unbounded at the origin, sit in the class with a
determined index. On that type the
transform is carried by an integral, which is what makes the analytic side available, and on a
general kernel it is not. The two properties need not be settled together, and the passage from
the first to the second is what this chapter carries out.

\begin{remarkx}[Notation, and it follows the type]
Where the kernel is induced by a profile, $G(n,k)=g(k/n)$, this volume writes $g$ and treats it
as a function of the single real variable $x=k/n$ on $(0,1]$, and writes $g^{*}$, $\tau(g)$,
$\alpha(g)$ and $\eta(g)$ for its transform and its three indices. Where the kernel is not of
that form, the volume writes $G(n,k)$, and $G^{*}$, $\tau(G)$, $\alpha(G)$ and $\eta(G)$
accordingly. The quantities are the same in both cases and the letter records which object
carries them. A property of one real variable, continuity for instance, is asserted of a profile
and never of a kernel of two integer variables.
\end{remarkx}

The arithmetic Mellin transform locates the constant of the first estimate and, in favorable
cases, the index itself through its first zero. It enters neither definition above. The existence of $\alpha(G)$ is an arithmetic property of the discrete
equation, whereas the equality $\alpha(G)=\eta(G)$ between the two indices is a theorem to be
proved wherever it holds. The next section exhibits a kernel where the transform is constant, so
that this equality is empty and the index rests on the arithmetic alone.

\section{A regular arithmetic function without an analytic index}
\label{sec:rational_kernel}

The first kernel outside the good variation class is a rational one. Fix two parameters $x,y>0$
with $x\neq y$ and set
\[G(n,k)=\frac{n+k+x}{n+k+y}\cdot\frac{2n+y}{2n+x},\qquad 1\le k\le n.
\]
The second factor depends on $n$ alone and normalizes the diagonal to $G(n,n)=1$. This kernel is
not a function of the ratio $k/n$. The additive shifts $x$ and $y$ break the homogeneity in $n$
and $k$, so $G$ leaves the good variation class while staying smooth and rational, as far from
an arithmetic weight as a kernel of the theory can be.

Its transform is constant. The two factors tend to $1$ uniformly in $k/n$ as $n\to\infty$, so
the limit of Definition~\ref{def:mellin} is the transform of the constant kernel $1$,
\[
G^{*}(z)=1\qquad\text{for all }z .
\]
The transform has no zero, so the analytic index $\eta(G)$ is undefined and the analytic route of
Part~I, which reads the index from the first zero of the transform, has nothing to act on. The
regularity index, if it exists, must come from the equation itself.

It does. The forward difference of the weight in $k$ is
\[
G(n,k+1)-G(n,k)=\frac{2n+y}{2n+x}\cdot\frac{y-x}{(n+k+y)(n+k+y+1)},
\]
and Abel summation\index[terms]{Abel summation} of the defining equation \eqref{eq:defining_relation_gen} with forcing
$n^{-\beta}$ turns it into an exact recurrence for the partial sums,
\begin{equation}\label{eq:rational_recurrence}
A(n)=n^{-\beta}+(y-x)\,\frac{2n+y}{2n+x}\sum_{k=1}^{n-1}
\frac{A(k)}{(n+k+y)(n+k+y+1)}.
\end{equation}
The feedback kernel is of size $n^{-2}$, and reading \eqref{eq:rational_recurrence} against the
forcing separates two regimes. For $\beta<2$ the feedback sum is $o(n^{-\beta})$, the partial
sums follow the forcing and $A(n)\sim n^{-\beta}$, the transparent regime. For $\beta>2$ the
forcing decays faster than the feedback allows, the iterated sums $AA(n)=\sum_{k\le n}A(k)$ tend
to a finite limit $\ell(\beta)$, depending on the pair as well as on the exponent, and
\[
A(n)\sim(y-x)\,\ell(\beta)\,n^{-2},
\]
the absorbed regime, in which the partial sums settle at the intrinsic rate $n^{-2}$ carried by
the kernel. The two regimes meet at
\[
\alpha(G)=2,
\]
the regularity index of the kernel, for every admissible pair. Sharpness asks that $\ell$ be
nonzero at exponents accumulating at two from above, and Theorem~\ref{thm:H_index} obtains it
from a family of shifted impulse solutions whose sums are positive at every high enough rank,
together with the uniqueness of Dirichlet\index[names]{Dirichlet, P. G. L.} coefficients. Appendix~\ref{app:H} carries the
estimates in full.

The intrinsic rate $n^{-2}$ is the rate of the homogeneous solution, the sequence left by the
forcing $r\equiv0$. With a constant transform there is no zero to carry the index, and the index
is read instead from the decay of that homogeneous solution. This is the homogeneity principle\index[terms]{homogeneity principle},
stated for constant transform kernels in \S\ref{sec:homogeneity} of
Chapter~\ref{chap:principles} and established there for the present kernel by the case
analysis of Appendix~\ref{app:H}. What the general principle needs beyond the homogeneous
decay, and why the decay alone cannot carry it, is taken up in \S\ref{sec:existence}.

The rational kernel returns once more. In Chapter~\ref{chap:gauge},
\S\ref{sec:equilibrium_rational}, it is the kernel that stays in equilibrium under a wide class
of gauges\index[terms]{gauge}, its index held at $2$ under deformation for $y>x$, while gauges\index[terms]{gauge} of fast growth lower
the index and mark the boundary of the phenomenon. The same kernel is thus a first regular
arithmetic function, an instance of the homogeneity principle, and an equilibrium kernel, met
here at the start of the theory it runs through.

\section{An existence theorem}
\label{sec:existence}

The gallery that closes this chapter, and the appendices behind it, prove existence one
kernel at a time. This section takes up the question in general, asking what property of a
kernel guarantees that the regularity index exists. The natural guess reads the index off
the homogeneous solution alone, and the first result below is a warning, the guess is
false. What replaces it is a criterion on the inverse of the summation operator, a bound on
all of its columns at once, and the criterion then delivers the index for a whole class of
kernels in a single theorem.

Fix a kernel $G$ with nonvanishing diagonal and put
\[
 d_n=G(n,n),
 \qquad
 c_{n,k}=G(n,k)-G(n,k+1).
\]
The Abelized operator is
\[
 (\mathcal L_Gx)(n)
 =d_nx(n)+\sum_{k=1}^{n-1}c_{n,k}x(k).
\]
Abel summation has turned the system on the coefficients $a_n$ into a triangular system on the partial sums, and the kernel introduced now is what inverts that second system. When every $d_n$ is nonzero, its inverse Green kernel\index[terms]{Green kernel} is defined by
\begin{align}
 H_G(m,m)&=d_m^{-1},
 \label{eq:ex-green-diagonal}\\
 H_G(n,m)&=-d_n^{-1}
 \sum_{k=m}^{n-1}c_{n,k}H_G(k,m)
 \qquad(m<n).
 \label{eq:ex-green-recursion}
\end{align}
The normalized homogeneous trajectory $h_G$ is fixed by
\[
 h_G(1)=1,
 \qquad
 (\mathcal L_Gh_G)(n)=0
 \quad(n\geq2).
\]
For a fixed column $m$, an endpoint limit
\begin{equation}
 Q_m=\lim_{n\to\infty}n^\gamma H_G(n,m)
 \label{eq:ex-column-limit}
\end{equation}
defines the connection coefficient\index[terms]{connection coefficient} of a forcing $f$ by
\[
 \mathcal C_G(f)=\sum_{m\geq1}Q_mf(m)
\]
when the series converges absolutely.

For the power profile $v_\beta(k)=k^{-\beta}$, the exact finite probe is
\[
 \mathcal G_n(\beta)
 =n^\beta(\mathcal L_Gv_\beta)(n)
 =d_n+\sum_{k<n}c_{n,k}
 \left(\frac{k}{n}\right)^{-\beta}.
\]
This probe is an operator consistency test.  Its convergence alone does not
control the inverse columns.

The inverse kernel met here in triangular form returns in its general setting in
Chapters~\ref{chap:volterra} and~\ref{chap:discrete_volterra}.

\begin{lemma}
\label{lem:ex-abel-green}
The defining equation \eqref{eq:defining_relation_gen} is equivalent to
\begin{equation}
 \mathcal L_GA(n)=n^{-\beta}.
 \label{eq:ex-abel-equation}
\end{equation}
For every forcing $f$, the unique solution of
$\mathcal L_GA=f$ is
\begin{equation}
 A(n)=\sum_{m=1}^{n}H_G(n,m)f(m).
 \label{eq:ex-green-convolution}
\end{equation}
The normalized homogeneous trajectory satisfies
\begin{equation}
 h_G(n)=d_1H_G(n,1).
 \label{eq:ex-first-column}
\end{equation}
Thus a one-parameter homogeneous estimate controls only the first column of
the inverse matrix.
\end{lemma}

\begin{proof}
Writing $a_k=A(k)-A(k-1)$ and shifting the second sum gives
\[
 \sum_{k=1}^{n}a_kG(n,k)
 =d_nA(n)+\sum_{k<n}c_{n,k}A(k).
\]
This proves \eqref{eq:ex-abel-equation}.  Forward substitution in the lower
triangular matrix gives \eqref{eq:ex-green-diagonal} and
\eqref{eq:ex-green-recursion}.  A second forward substitution gives
\eqref{eq:ex-green-convolution}.

The first inverse column solves
$\mathcal L_GH_G(\mathord\cdot,1)=e_1$.  Multiplication by $d_1$ makes its
first value equal to one and leaves every later row homogeneous.  Uniqueness
therefore gives \eqref{eq:ex-first-column}.
\end{proof}

The warning comes first.

\begin{theorem}
\label{thm:ex-counterexample}
There is a bounded real triangular kernel with
\[
 G(n,n)=1,
 \qquad
 G^*(z)=1,
 \qquad
 h_G(n)=n^{-2},
\]
where the transform is entire and zero-free, but $G$ is not a RAF with any
finite index.

One such kernel is obtained from
\begin{equation}
 h_n=n^{-2},
 \qquad
 b_1=1,
 \qquad
 b_n=h_n-h_{n-1}
 \quad(n\geq2).
 \label{eq:ex-counter-h}
\end{equation}
Set $G(1,1)=1$.  At ranks two and three set
\[
 G(n,1)=1-h_n,
 \qquad
 G(n,k)=1
 \quad(2\leq k\leq n).
\]
For $n\geq4$, set
\[
 G(n,1)=1-h_n+5b_{n-1},
 \qquad
 G(n,n-1)=-4,
\]
and set all other entries in that row equal to one.

The exact finite probe for $n\geq4$ is
\begin{equation}
 \mathcal G_n(z)=1+\xi_nn^z
 {}+5\left[
 \left(1-\frac2n\right)^{-z}
 -\left(1-\frac1n\right)^{-z}
 \right],
 \qquad
 \xi_n=-h_n+5b_{n-1}.
 \label{eq:ex-counter-probe}
\end{equation}
It converges locally uniformly to one on $\Re z<2$.

For the forcing exponent $\beta=1$, write
$A_1(n)=h_n+q_n$.  Then
\begin{equation}
 q_1=0,
 \qquad
 q_2=\frac12,
 \qquad
 q_3=\frac13,
 \qquad
 q_n=5q_{n-1}-5q_{n-2}+\frac1n
 \quad(n\geq4).
 \label{eq:ex-counter-forced}
\end{equation}
In particular,
\begin{equation}
 |q_n|\geq\frac{263}{60}2^{n-5}
 \qquad(n\geq5).
 \label{eq:ex-counter-growth}
\end{equation}
For every $m\geq3$, the shifted Green column is
\begin{equation}
 H_G(n,m)=
 \frac{r_+^{n-m+1}-r_-^{n-m+1}}{\sqrt5},
 \qquad
 r_\pm=\frac{5\pm\sqrt5}{2}.
 \label{eq:ex-counter-green}
\end{equation}
The shifted columns grow exponentially although the first column equals
$n^{-2}$.
\end{theorem}

\begin{proof}
Every entry of the displayed kernel lies between $-4$ and $1$, apart from
the first entry which is also bounded.  Its diagonal is one.  At ranks at
least four, the Abelized equation is
\begin{equation}
 (\mathcal L_GA)(n)
 =A(n)+\xi_nA(1)+5A(n-2)-5A(n-1).
 \label{eq:ex-counter-abel}
\end{equation}
The coefficient sequence in \eqref{eq:ex-counter-h} telescopes to $h_n$.
Substitution into \eqref{eq:ex-counter-abel} gives
\[
 h_n+(-h_n+5b_{n-1})+5h_{n-2}-5h_{n-1}=0.
\]
The two earlier ranks give the same conclusion directly.  Hence the
normalized homogeneous trajectory is exactly $h_n=n^{-2}$.

Apply the finite-probe definition to the three nonzero Abel coefficients in
\eqref{eq:ex-counter-abel}.  This gives
\eqref{eq:ex-counter-probe}.  Since $\xi_n=\mathcal O(n^{-2})$, the first
error is locally $\mathcal O(n^{\Re z-2})$.  The bracket is locally
$\mathcal O(n^{-1})$ by the mean value theorem on $[1/2,1]$.  The probes
therefore converge locally uniformly to one on $\Re z<2$.  On the defining
half plane $\Re z<0$, the transform is one.  Its continuation is the entire
constant one.

Since $h$ is homogeneous after rank one, substitution of $A_1=h+q$ into
the forced equation gives \eqref{eq:ex-counter-forced}.  Direct calculation
gives
\[
 q_4=-\frac7{12},
 \qquad
 q_5=-\frac{263}{60}.
\]
Suppose $q_{n-1}$ and $q_{n-2}$ are negative, with
$x=|q_{n-1}|\geq2|q_{n-2}|$ and $x\geq263/60$.  Then
\[
 q_n\leq-5x+\frac52x+\frac1n\leq-2x.
\]
Induction from $n=6$ proves \eqref{eq:ex-counter-growth}.

For a shifted impulse at $m\geq3$, the values at ranks $m$ and $m+1$
are one and five.  The later values satisfy
\[
 H_G(n,m)=5H_G(n-1,m)-5H_G(n-2,m).
\]
Solving this recurrence gives \eqref{eq:ex-counter-green}.

If $\alpha(G)>1$, the fixed definition requires transparency at
$\beta=1$, contrary to exponential growth.  If $\alpha(G)\leq1$, the RAF
definition requires absorption at that exponent, with the same contradiction.
No finite RAF index is possible.
\end{proof}

A workable sufficient condition asks two things of the kernel, a bound on the diagonal and a
decay rate on the columns of its inverse. The next statement fixes the two parameters and reads
the index off them.

\begin{theorem}
\label{thm:ex-green-criterion}
Assume $0<\gamma<\infty$ and
\[
 D=\sup_{n\geq1}|d_n^{-1}|<\infty.
\]
For every $0<s<\gamma$, assume
\begin{equation}
 |H_G(n,m)|
 \leq \frac{C_s}{m}
 \left(\frac mn\right)^s
 \qquad(1\leq m<n).
 \label{eq:ex-green-subcritical}
\end{equation}
For every real $\beta<\gamma$, assume that $\Xi_G(\beta)$ is finite and
that there is a comparison sequence $P_\beta$ such that
\begin{align}
 P_\beta(n)&=\Xi_G(\beta)n^{-\beta}+o(n^{-\beta}),
 \label{eq:ex-comparison-profile}\\
 \epsilon_{\beta,n}
 &=n^\beta\bigl(\mathcal L_GP_\beta(n)-n^{-\beta}\bigr)\longrightarrow0.
 \label{eq:ex-profile-consistency}
\end{align}
Then every $\beta<\gamma$ is transparent.

At a point where $0<|G^*(\beta)|<\infty$, the finite-probe condition
\[
 \mathcal G_n(\beta)\longrightarrow G^*(\beta)
\]
supplies these hypotheses with
$P_\beta(n)=\Xi_G(\beta)n^{-\beta}$.  At a pole of $G^*$, a sequence
$P_\beta=o(n^{-\beta})$ satisfying
\eqref{eq:ex-profile-consistency} is a separate hypothesis.  The Green
bound alone does not imply this cancellation.
Moreover,
\begin{equation}
 A_\beta(n)=
 \mathcal O_{\beta,\eps}
 \left(n^{-\gamma+\eps}\right)
 \qquad(\beta\geq\gamma).
 \label{eq:ex-candidate-absorption}
\end{equation}
This is absorption at the candidate exponent.  It does not by itself prove
that $G$ is a RAF.

If the endpoint estimate
\begin{equation}
 |H_G(n,m)|
 \leq\frac{C_\gamma}{m}
 \left(\frac mn\right)^\gamma
 \qquad(m<n)
 \label{eq:ex-green-endpoint}
\end{equation}
also holds, then
\begin{align}
 A_\gamma(n)&=\mathcal O(n^{-\gamma}\log n),
 \label{eq:ex-endpoint-log}\\
 A_\beta(n)&=\mathcal O_\beta(n^{-\gamma})
 \qquad(\beta>\gamma).
 \label{eq:ex-endpoint-above}
\end{align}

Assume in addition that every limit in \eqref{eq:ex-column-limit} exists.
For $\beta>\gamma$, define
\begin{equation}
 \mathcal C_G(\beta)=
 \sum_{m\geq1}Q_mm^{-\beta}.
 \label{eq:ex-power-connection}
\end{equation}
If there are exponents $\beta_j\downarrow\gamma$ for which
$G^*(\beta_j)$ is finite and nonzero and
\begin{equation}
 \mathcal C_G(\beta_j)\neq0,
 \label{eq:ex-green-nonannihilation}
\end{equation}
then sharpness holds, $G$ is a RAF, and
\[
 \alpha(G)=\gamma.
\]
\end{theorem}

\begin{proof}
Fix $\beta<\gamma$ and use the profile in
\eqref{eq:ex-comparison-profile}.  Equation
\eqref{eq:ex-profile-consistency} is the exact identity
\begin{equation}
 \mathcal L_GP_\beta(n)-n^{-\beta}
 =n^{-\beta}\epsilon_{\beta,n}.
 \label{eq:ex-exact-defect}
\end{equation}
Choose $s$ with $\max(0,\beta)<s<\gamma$.  Inversion of
\eqref{eq:ex-exact-defect} and the Green bounds give
\begin{align*}
 n^\beta|A(n)-P_\beta(n)|
 \leq{}&D|\epsilon_{\beta,n}|\\
 &+C_sn^{-(s-\beta)}
 \sum_{m<n}m^{s-\beta-1}|\epsilon_{\beta,m}|.
\end{align*}
The first term tends to zero.  The second tends to zero by weighted Cesaro
summation because $s-\beta>0$.  Together with
\eqref{eq:ex-comparison-profile}, this proves
\eqref{eq:transparency_additive}.  No rate beyond
$\epsilon_{\beta,n}\to0$ is used.  If the finite probe converges to a
finite nonzero value, direct substitution gives
$\epsilon_{\beta,n}=\Xi_G(\beta)\mathcal G_n(\beta)-1$.  When
$\Xi_G(\beta)=0$, that substitution gives the constant defect $-1$,
which explains the separate comparison-profile assumption at a transform
pole.

For $\beta\geq\gamma$, Green inversion applied directly to the forcing
gives
\[
 |A(n)|
 \leq Dn^{-\beta}
 +C_sn^{-s}\sum_{m<n}m^{s-\beta-1}.
\]
Choose $s<\gamma$ close enough to $\gamma$.  The sum is bounded because
$s<\beta$.  This proves \eqref{eq:ex-candidate-absorption}, including all
$\eps>0$ by decreasing $s$ when necessary.

Under \eqref{eq:ex-green-endpoint}, the same sum is harmonic at
$\beta=\gamma$ and bounded above that exponent.  This proves
\eqref{eq:ex-endpoint-log} and \eqref{eq:ex-endpoint-above}.

The endpoint estimate also gives
\[
 |Q_m|\leq C_\gamma m^{\gamma-1}.
\]
Thus \eqref{eq:ex-power-connection} converges absolutely when
$\beta>\gamma$.  Dominated convergence in
\eqref{eq:ex-green-convolution}, with the diagonal term treated separately,
gives
\[
 n^\gamma A(n)\longrightarrow\mathcal C_G(\beta).
\]
At every exponent in \eqref{eq:ex-green-nonannihilation}, this behavior is
incompatible with a nonzero multiple of $n^{-\beta_j}$.  Every number
larger than $\gamma$ contains such a failed exponent.  Thus
the sharpness condition of Remark~\ref{rem:index_three_conditions} holds at $\gamma$.
Subcritical transparency and candidate absorption are the other two conditions of that
remark.  The three conditions therefore prove the RAF assertion and $\alpha(G)=\gamma$.
\end{proof}

A condition constraining a weighted row sum of the Green kernel instead of each of its entries
tolerates individual entries larger than \eqref{eq:ex-green-subcritical} allows, provided they
are compensated along the row. The family of powers is enough to carry the transparency and
absorption clauses.

\begin{theorem}
\label{thm:ex-green-schur}
Keep the setting of Theorem~\ref{thm:ex-green-criterion}, with $D=\sup_n|d_n^{-1}|$ finite, a
parameter $0<\gamma<\infty$, and for every $\beta<\gamma$ a finite $\Xi_G(\beta)$ together
with a comparison sequence satisfying \eqref{eq:ex-comparison-profile} and
\eqref{eq:ex-profile-consistency}. Assume, in place of \eqref{eq:ex-green-subcritical}, the
weighted row bound
\begin{equation}
\label{eq:ex-green-schur}
\sum_{m<n}|H_G(n,m)|\,m^{-s}\ \le\ C_s\,n^{-s}
\qquad(n\ge2)
\end{equation}
for every $s<\gamma$. Then every $\beta<\gamma$ is transparent, and
$A(n)=\mathcal O_{\beta,\eps}(n^{-\gamma+\eps})$ for every
$\beta\ge\gamma$. If \eqref{eq:ex-green-schur} holds at $s=\gamma$ as well, that second
bound improves to $A(n)=\mathcal O_\beta(n^{-\gamma})$. The pointwise estimate
\eqref{eq:ex-green-subcritical} implies \eqref{eq:ex-green-schur}, so the transparency and
absorption clauses of Theorem~\ref{thm:ex-green-criterion} are recovered.
\end{theorem}

\begin{proof}
That \eqref{eq:ex-green-subcritical} implies \eqref{eq:ex-green-schur} follows by choosing
$\max(0,s)<s'<\gamma$, since then
\[
 \sum_{m<n}|H_G(n,m)|\,m^{-s}
 \le C_{s'}\,n^{-s'}\sum_{m<n}m^{s'-s-1}
 \le C'\,n^{-s}.
\]

Fix $\beta<\gamma$. Inversion of the exact defect \eqref{eq:ex-exact-defect} gives
\[
 A(n)-P_\beta(n)=-\sum_{m\le n}H_G(n,m)\,m^{-\beta}\epsilon_{\beta,m}.
\]
Given
$\delta>0$, choose $M$ with $|\epsilon_{\beta,m}|\le\delta$ for $m>M$ and split that sum at
the diagonal, at the ranks up to $M$, and at the rest. After multiplication by $n^{\beta}$ the
diagonal contributes at most $D|\epsilon_{\beta,n}|$, which tends to zero, and the last part
at most $\delta C_\beta$ by \eqref{eq:ex-green-schur} at $s=\beta$. For the middle part fix
$\beta'$ between $\beta$ and $\gamma$. A single term of \eqref{eq:ex-green-schur} at
$s=\beta'$ gives the columnwise bound $|H_G(n,m)|\le C_{\beta'}m^{\beta'}n^{-\beta'}$, so
that part is at most
$C_{\beta'}n^{\beta-\beta'}\sum_{m\le M}m^{\beta'-\beta}|\epsilon_{\beta,m}|$, which tends to zero at fixed
$M$ because $\beta<\beta'$. Hence
$\limsup_{n\to\infty}n^{\beta}|A(n)-P_\beta(n)|\le\delta C_\beta$ for every $\delta$, and
\eqref{eq:ex-comparison-profile} turns that into transparency at $\beta$.

For $\beta\ge\gamma$ and any $s<\gamma$ the same inversion applied to the forcing itself
gives $|A(n)|\le Dn^{-\beta}+\sum_{m<n}|H_G(n,m)|m^{-s}\le Dn^{-\beta}+C_sn^{-s}$,
since $m^{-\beta}\le m^{-s}$ on that range. Taking $s=\gamma-\eps$ gives the
absorption bound, and the same estimate at $s=\gamma$ gives the version without loss.
\end{proof}

\begin{openproblem}[Boundary of the Green criterion]
\label{op:ex-green-boundary}
Three questions remain.
\begin{enumerate}[label=\textup{(\roman*)}]
\item Can the family of power-weight estimates \eqref{eq:ex-green-schur} be replaced by a
single positive weight $w\colon\mathbb N^*\to(0,\infty)$, independent of $\beta$, through a
Schur estimate
\[
 \sum_{m<n}|H_G(n,m)|w(m)\le Cw(n),
\]
while still implying transparency for every $\beta<\gamma$ and absorption for every
$\beta\ge\gamma$?
\item Under the row hypothesis \eqref{eq:ex-green-schur}, the column limits
\eqref{eq:ex-column-limit} need not exist. Can the sharpness clause of
Theorem~\ref{thm:ex-green-criterion} instead be recovered from a rowwise non-annihilation
condition? If the column limits do exist, can the dominating envelope
$\sup_n n^{\gamma}|H_G(n,m)|\le\bar Q_m$, with
$\sum_m\bar Q_m m^{-\beta}<\infty$, be dropped?
\item Which, if any, of the hypotheses in Theorems~\ref{thm:ex-green-criterion}
and~\ref{thm:ex-green-schur} are necessary?
\end{enumerate}
The counterexample above proves only that control of a single column is insufficient.
\end{openproblem}

For a kernel $G(n,k)=g(k/n)$ the Green estimate\index[terms]{Green estimate}\index[names]{Green, G.} can be verified on the
function $g$ itself.

\begin{lemma}
\label{lem:ex-fgv-green}
Let $g$ be bounded on $(0,1]$, let $d=g(1)\neq0$, and put
\[
 \delta_{n,k}=g\left(\frac{k+1}{n}\right)
 -g\left(\frac{k}{n}\right).
\]
Then
\begin{equation}
 (\mathcal L_gA)(n)
 =dA(n)-\sum_{k<n}\delta_{n,k}A(k).
 \label{eq:ex-fgv-operator}
\end{equation}
Fix $0<s<\gamma$.  Suppose that
\begin{align}
 |\delta_{n,k}|&\leq
 \frac{D_s}{k}\left(\frac{k}{n}\right)^s,
 \label{eq:ex-fgv-local}\\
 \sum_{k<n}|\delta_{n,k}|
 \left(\frac{k}{n}\right)^{-s}
 &\leq |d|-\kappa_s
 \qquad(n\geq N_s)
 \label{eq:ex-fgv-gap}
\end{align}
for some $D_s,\kappa_s>0$.  The inverse kernel then satisfies
\[
 |H_g(n,m)|\leq
 \frac{C_s}{m}\left(\frac mn\right)^s
 \qquad(m<n).
\]

Suppose further that $g$ is absolutely continuous on $[0,1]$ and
\begin{equation}
 |g'(t)|\leq Ct^{a-1}
 \label{eq:ex-weighted-derivative}
\end{equation}
almost everywhere for some $a>0$.  For real $\beta<a$,
\[
 g^*(\beta)=g(1)-\int_0^1t^{-\beta}g'(t)\,dt
\]
and the exact probe defect is
\begin{align}
 \mathcal G_n(\beta)-g^*(\beta)
 ={}&\int_0^{1/n}t^{-\beta}g'(t)\,dt
 \notag\\
 &+\sum_{k=1}^{n-1}
 \int_{k/n}^{(k+1)/n}
 \left[
 t^{-\beta}-\left(\frac{k}{n}\right)^{-\beta}
 \right]g'(t)\,dt.
 \label{eq:ex-fgv-probe-defect}
\end{align}
For $\beta\neq0$, this gives
\begin{equation}
 |\mathcal G_n(\beta)-g^*(\beta)|
 \ll_\beta
 \begin{cases}
 n^{-(a-\beta)} & 0<a-\beta<1,\\
 n^{-1}\log n & a-\beta=1,\\
 n^{-1} & a-\beta>1.
 \end{cases}
 \label{eq:ex-fgv-probe-rate}
\end{equation}
At $\beta=0$, the defect is $\mathcal O(n^{-a})$.
\end{lemma}

\begin{proof}
Discrete Abel summation gives \eqref{eq:ex-fgv-operator}.  Its inverse obeys
\[
 H_g(n,m)=d^{-1}\sum_{k=m}^{n-1}\delta_{n,k}H_g(k,m)
 \qquad(m<n).
\]
Separate the term $k=m$, where $H_g(m,m)=d^{-1}$.  Strong induction and
\eqref{eq:ex-fgv-gap} give
\begin{align*}
 |H_g(n,m)|
 \leq{}&\frac{D_s}{|d|^2m}
 \left(\frac mn\right)^s\\
 &+\left(1-\frac{\kappa_s}{|d|}\right)
 \frac{C_s}{m}\left(\frac mn\right)^s.
\end{align*}
A large enough $C_s$ closes the induction.  The finitely many ranks before
$N_s$ are absorbed into that constant.

Integration by parts first in $\Re z<0$ gives
\[
 -z\int_0^1g(t)t^{-z-1}\,dt
 =g(1)-\int_0^1t^{-z}g'(t)\,dt.
\]
The derivative bound continues the right side to $\Re z<a$.  Splitting
the integral into cells and subtracting the left-endpoint sum gives
\eqref{eq:ex-fgv-probe-defect}.  For $\beta\neq0$, the mean value theorem
and \eqref{eq:ex-weighted-derivative} bound its cell sum by
\[
 n^{\beta-a}\sum_{k<n}k^{a-\beta-2}.
\]
The three elementary power-sum estimates give
\eqref{eq:ex-fgv-probe-rate}.  At zero, telescoping gives
$\mathcal G_n(0)=g(1/n)$, whose difference from $g(0)$ is
$\mathcal O(n^{-a})$.
\end{proof}

The criterion applies to a class wide enough to be worth isolating, that of the absolutely
continuous profiles positive at the origin whose derivative obeys a one-sided power bound.

\begin{theorem}
\label{thm:ex-fgv-existence}
Let $g$ be absolutely continuous on $[0,1]$.  Assume that for some
$a>0$,
\[
 g(0)>0,
 \qquad
 0\leq g'(t)\leq Ct^{a-1}
\]
almost everywhere.  Suppose that there is $\gamma$ with
$0<\gamma<a$ such that
\begin{equation}
 g(1)=\int_0^1t^{-\gamma}g'(t)\,dt.
 \label{eq:ex-fgv-crossing}
\end{equation}
Then the kernel $G(n,k)=g(k/n)$ is a RAF and
\[
 \alpha(g)=\eta(g)=\gamma.
\]
Here $\eta(g)$ is the infimum of the real parts of the zeros of the
continuation in $\Re z<a$.
\end{theorem}

Theorem~\ref{thm:ex-fgv-existence} does not require $C^2$ regularity.
An ordinary $C^1$ profile with bounded nonnegative derivative fits the case
$a=1$.  A continuous piecewise $C^1$ profile is also covered when it is
globally absolutely continuous and obeys the same weighted derivative bound.
Finite jumps and general bounded variation\index[terms]{bounded variation} do not supply the local Green
estimate used here.

\begin{proof}
Let
\[
 F(s)=\int_0^1t^{-s}g'(t)\,dt
 \qquad(s<a).
\]
This function is strictly increasing on the real interval relevant to
\eqref{eq:ex-fgv-crossing}.  Also $F(0)=g(1)-g(0)<g(1)$.  Hence
\[
 g^*(s)=g(1)-F(s)>0
 \qquad(s<\gamma).
\]

For $0<s<a$, the derivative bound gives
\[
 0\leq\delta_{n,k}
 \leq C_an^{-a}k^{a-1}
 \leq\frac{C_a}{k}\left(\frac{k}{n}\right)^s.
\]
For $s<\gamma$, the weighted increment sum tends to $F(s)<g(1)$.
It therefore has the gap required by Lemma~\ref{lem:ex-fgv-green}.
The lemma and Theorem~\ref{thm:ex-green-criterion} prove transparency below
$\gamma$ and absorption at that candidate exponent.

For real $s$ with $\gamma<s<a$, one has $g^*(s)<0$.  On the other
hand, the exact recurrence
\[
 A_s(n)=\frac1{g(1)}
 \left[n^{-s}+\sum_{k<n}\delta_{n,k}A_s(k)\right]
\]
shows by induction that every $A_s(n)$ is positive.  Suppose $s$ transparent,
so that $A_s(n)=\Xi'n^{-s}+o(n^{-s})$ for some constant $\Xi'$.  Applying
$\mathcal L_g$ to the two parts identifies that constant.  On the power
profile the finite probe converges for $s<a$ by
Lemma~\ref{lem:ex-fgv-green} and contributes $\Xi'g^*(s)n^{-s}(1+o(1))$, and
on the remainder the positivity of the increments together with the
convergence of $\sum_k\delta_{n,k}(k/n)^{-s}$ gives $o(n^{-s})$ by the
Toeplitz argument already used for transparency.  The defining equation then
forces $\Xi'g^*(s)=1$, so $\Xi'=1/g^*(s)<0$, against the positivity of
$A_s(n)$.  Hence no $s\in(\gamma,a)$ is transparent, such failures occur
arbitrarily close to $\gamma$, and $\alpha(g)=\gamma$.  The candidate absorption is now
the RAF absorption bound.

It remains to locate the first transform zero.  If $g^*(\rho)=0$ and
$\Re\rho<\gamma$, then
\[
 g(1)=\left|\int_0^1t^{-\rho}g'(t)\,dt\right|
 \leq F(\Re\rho)<F(\gamma)=g(1),
\]
which is impossible.  Equation \eqref{eq:ex-fgv-crossing} gives a real zero
at $\gamma$.  Thus $\eta(g)=\gamma$.
\end{proof}

A cell containing a jump has an order-one increment and does not satisfy
\eqref{eq:ex-fgv-local}, so the pointwise estimate is out of reach there. The weighted row bound
of Theorem~\ref{thm:ex-green-schur} is within reach, and finitely many upward jumps are covered
by it.

\begin{theorem}
\label{thm:ex-fgv-jumps}
Let $g_0$ be absolutely continuous on $[0,1]$ with $g_0(0)>0$ and $0\le g_0'(t)\le Ct^{a-1}$
almost everywhere for some $a>0$. Let $0<\lambda_1<\cdots<\lambda_r<1$ and $J_1,\dots,J_r>0$,
and set
\[
 g(t)=g_0(t)+\sum_{j\,:\,\lambda_j<t}J_j
 \qquad(0<t\le1),
\]
so that $g$ is nondecreasing with an upward jump $J_j$ at each $\lambda_j$. Put
\[
 F(s)=\int_0^1t^{-s}g_0'(t)\,dt+\sum_{j\le r}J_j\lambda_j^{-s}
 \qquad(s<a),
\]
and suppose there is $\gamma$ with $0<\gamma<a$ and $F(\gamma)=g(1)$. Then $G(n,k)=g(k/n)$ is
a RAF and $\alpha(g)=\eta(g)=\gamma$. For $r=0$ this is
Theorem~\ref{thm:ex-fgv-existence}, and for $g_0$ constant equal to $\mu$ with $r=1$ and
$J_1=1-\mu$ it is the bipolar step kernel of Appendix~\ref{app:C}, whose index
$\log(1-\mu)/\log\lambda_1$ it returns.
\end{theorem}

\begin{proof}
Write $d=g(1)>0$ and keep the increments $\delta_{n,k}$ of Lemma~\ref{lem:ex-fgv-green}, which
are nonnegative here. The measure $dg$ on $(0,1)$ is $g_0'(t)\,dt$ together with the point masses
$J_j$. For $\Re z<a$, integration by parts gives the holomorphic continuation
\[
 g^{*}(z)=d-\int_0^1t^{-z}g_0'(t)\,dt-\sum_{j\le r}J_j\lambda_j^{-z}.
\]
On the real axis this is $g^{*}(s)=d-F(s)$. The crossing hypothesis forces $dg$ to have
positive mass in $(0,1)$, so $F$ is strictly increasing. If $g^{*}(\rho)=0$ and
$\Re\rho<\gamma$, then
\[
 d=\left|\int_{(0,1)}t^{-\rho}\,dg(t)\right|
 \le \int_{(0,1)}t^{-\Re\rho}\,dg(t)
 =F(\Re\rho)<F(\gamma)=d,
\]
a contradiction. Since $g^{*}(\gamma)=0$, it follows that $\eta(g)=\gamma$.

Fix $s<\gamma$ and put $\kappa=g(1)-F(s)>0$. The weighted increment sums
\[
 W_s(n)=\sum_{k<n}\delta_{n,k}\Bigl(\frac kn\Bigr)^{-s}
\]
converge to $F(s)$, the absolutely continuous part by Lemma~\ref{lem:ex-fgv-green} and the jump
part because the single cell containing $\lambda_j$ has its left endpoint tending to
$\lambda_j$. Choose $N_s$ beyond which
$W_s(n)\le d-\kappa/2$. Inverting \eqref{eq:ex-fgv-operator} gives the Green recursion
\[
 H_g(n,m)=d^{-1}\sum_{k=m}^{n-1}\delta_{n,k}H_g(k,m),
 \qquad H_g(m,m)=d^{-1},
\]
so the weighted row sums $S_s(n)=\sum_{m<n}|H_g(n,m)|m^{-s}$ obey
$S_s(n)\le d^{-1}\sum_{k<n}\delta_{n,k}\bigl[S_s(k)+d^{-1}k^{-s}\bigr]$. Assume
inductively that $S_s(k)\le C_sk^{-s}$ below $n$, which holds at the finitely many ranks under
$N_s$ once $C_s$ is large enough. Then
\[
 S_s(n)\le d^{-1}(C_s+d^{-1})\,n^{-s}W_s(n)
 \le(C_s+d^{-1})\Bigl(1-\frac{\kappa}{2d}\Bigr)n^{-s}\le C_sn^{-s}
\]
as soon as $C_s\ge\frac2\kappa\bigl(1-\frac{\kappa}{2d}\bigr)$, and that is
\eqref{eq:ex-green-schur} at every $s<\gamma$.

For $\beta<\gamma$ the finite probes converge to $g^{*}(\beta)=g(1)-F(\beta)$, a finite
positive number, being Riemann-Stieltjes sums of $t^{-\beta}$ against the integrator $g$ split
at the jumps. The comparison sequence $P_\beta=n^{-\beta}/g^{*}(\beta)$ therefore has defect
tending to zero, and Theorem~\ref{thm:ex-green-schur} gives transparency below $\gamma$
together with $A_\beta(n)=\mathcal O(n^{-\gamma+\eps})$ at and above it.

Sharpness follows as in Theorem~\ref{thm:ex-fgv-existence}. For real $s$ between $\gamma$ and
$a$ one has $g^{*}(s)<0$, while the recurrence
$A_s(n)=d^{-1}\bigl[n^{-s}+\sum_{k<n}\delta_{n,k}A_s(k)\bigr]$ makes every $A_s(n)$ positive
since the increments are nonnegative. Were $s$ transparent with coefficient $\Xi'$, the
Toeplitz argument of that proof applies as follows. For each fixed $k$, all jump points lie
outside $[k/n,(k+1)/n)$ once $n$ is large, and the derivative bound gives
$\delta_{n,k}=\mathcal O_k(n^{-a})$. Hence
\[
 \delta_{n,k}\Bigl(\frac kn\Bigr)^{-s}=\mathcal O_k(n^{s-a})\longrightarrow0.
\]
Together with the convergence of the row sums
$\sum_k\delta_{n,k}(k/n)^{-s}$, this is precisely the pair of Toeplitz conditions used there.
It would force $\Xi'g^{*}(s)=1$ and therefore $\Xi'<0$, against that positivity. So no
exponent of $(\gamma,a)$ is transparent, such exponents accumulate at $\gamma$, and with the
absorption just obtained the three conditions of
Remark~\ref{rem:index_three_conditions} hold, which gives the RAF property and
$\alpha(g)=\gamma$.
\end{proof}

\begin{remark}[The ridge of the Green kernel]\label{rem:ex-fgv-ridge}
When $r\ge1$, the pointwise estimate fails while the row bound does not. For the rank $m$ such
that $\lambda_j\in[m/n,(m+1)/n)$, the $k=m$ term in the Green recursion is at least
$d^{-2}J_j$. Thus the Green kernel carries an entry of order one along
$m\approx\lambda_jn$, which is the discrete shadow of the resolvents with point masses of
Chapter~\ref{chap:volterra}. One entry of order one at $m\asymp\lambda n$ weighs
$\mathcal O(n^{-s})$ in \eqref{eq:ex-green-schur} and is admissible there, while
\eqref{eq:ex-green-subcritical} would ask it to be $\mathcal O(n^{-1})$. This is the
compensation along the row that Open Problem~\ref{op:ex-green-boundary} was stated to reach.
\end{remark}

\begin{openproblem}[FGV classes with jumps]
\label{op:ex-fgv-jumps}
Three cases remain beyond Theorem~\ref{thm:ex-fgv-jumps}. For a profile with finitely many
jumps of arbitrary signs and $g(1)\neq0$, put
\[
 V(s):=\int_0^1t^{-s}|g_0'(t)|\,dt
       +\sum_j|J_j|\lambda_j^{-s},
 \qquad
 \widetilde\gamma_V:=\sup\{s\in\mathbb R:V(s)<|g(1)|\}.
\]
The function $V$ is nondecreasing. The same absolute-value argument therefore gives the row
condition \eqref{eq:ex-green-schur} for every $s<\widetilde\gamma_V$. When the
comparison-profile hypotheses also hold, Theorem~\ref{thm:ex-green-schur} gives transparency
below $\widetilde\gamma_V$, but not sharpness. If the RAF index exists, only
$\widetilde\gamma_V\le\alpha(g)$ follows.
\begin{enumerate}[label=\textup{(\roman*)}]
\item For finitely many jumps of arbitrary signs, which additional sign or non-annihilation
hypotheses determine $\alpha(g)$, and when is the boundary $\widetilde\gamma_V$ sharp?
\item Can the existence theory be extended to countably many jumps when $V(s)$ diverges, by
using cancellation or the arrangement of the jumps instead of absolute variation? This includes
the relevant positive ranges for the Ingham and broken kernels.
\item For a singular continuous component of $dg$, which weighted-variation condition ensures
the Green row bound, and what additional condition gives sharpness?
\end{enumerate}
\end{openproblem}

Finite sums of positive powers with a positive constant term meet the hypotheses of
Theorem~\ref{thm:ex-fgv-existence}, which gives the first family of examples.

\begin{corollary}
\label{cor:ex-muntz}
Let
\[
 g(t)=c_0+\sum_{j=1}^{r}c_jt^{a_j},
 \qquad
 c_0>0,
 \qquad
 c_j>0,
 \qquad
 a_j>0.
\]
There is a unique $\gamma\in(0,\min_ja_j)$ such that
\begin{equation}
 \sum_{j=1}^{r}\frac{c_ja_j}{a_j-\gamma}
 =c_0+\sum_{j=1}^{r}c_j.
 \label{eq:ex-muntz-root}
\end{equation}
The associated FGV is a RAF and $\alpha(g)=\eta(g)=\gamma$.
For $g(t)=c_0+c_1t^a$,
\begin{equation}
 \gamma=\frac{ac_0}{c_0+c_1}.
 \label{eq:ex-one-monomial}
\end{equation}
The affine kernel\index[terms]{affine kernel} is recovered by $a=1$, $c_0=\lambda$, and
$c_1=1-\lambda$.
\end{corollary}

\begin{proof}
The hypotheses of Theorem~\ref{thm:ex-fgv-existence} hold with
$a=\min_ja_j$.  The crossing equation is exactly
\eqref{eq:ex-muntz-root}.  Its left side is continuous and strictly
increasing from $\sum_jc_j$ to infinity, while the right side exceeds its
initial value by $c_0$.  This proves existence and uniqueness.  Solving the
single-monomial equation gives \eqref{eq:ex-one-monomial}.
\end{proof}

\section{A profile with a natural boundary}
\label{sec:ex-lacunary}

The M\"untz sums of Corollary~\ref{cor:ex-muntz} have a rational transform, with poles at the
exponents $a_j$, and the Ingham function continues through the zeta function. Neither the
criterion of \S\ref{sec:existence} nor the index it produces needs anything of the sort. The
family below produces a proved RAF whose transform is holomorphic in a half plane and extends
holomorphically to no point of the edge of that half plane. It grew from a lacunary profile
constructed by Tenenbaum\index[names]{Tenenbaum, G.}, which is its member $c=1/2$, $D=0$ and is
treated in Remark~\ref{rem:ex-lacunary-original} and Open
Problem~\ref{op:ex-lacunary-original}.

Write $u=\log(1/t)$ and set
\begin{equation}
 N(u)=\#\{j\ge0:2^{j}<u\},
 \qquad
 S(u)=\sum_{j\ge0}(u-2^{j})_+=N(u)\,u-2^{N(u)}+1 .
 \label{eq:ex-lac-NS}
\end{equation}
For $c>0$ and $D\ge0$ the profile of the family is
\begin{equation}
 g_{c,D}(t)=1+t^{c}\bigl(D+S(\log(1/t))\bigr)
 \quad(0<t\le1),
 \qquad
 g_{c,D}(0)=1 .
 \label{eq:ex-lac-profile}
\end{equation}
The lacunarity\index[terms]{lacunary series} sits in $S$, whose breakpoints are the doubly
exponential $t=e^{-2^{j}}$.

\begin{lemma}[Transform of the lacunary profile]
\label{lem:ex-lac-transform}
Let $c>0$ and $D\ge0$. The profile \eqref{eq:ex-lac-profile} is absolutely continuous on
$[0,1]$, with
\begin{equation}
 g_{c,D}'(t)=t^{c-1}\bigl[c\bigl(D+S(u)\bigr)-N(u)\bigr],
 \qquad u=\log(1/t),
 \label{eq:ex-lac-derivative}
\end{equation}
outside the breakpoints, and for every $a\in(0,c)$ there is $C_a$ with
$|g_{c,D}'(t)|\le C_at^{a-1}$. Put
\begin{equation}
 \Lambda(w)=\sum_{j\ge0}e^{-w2^{j}}
 \qquad(\Re w>0),
 \label{eq:ex-lac-Lambda}
\end{equation}
holomorphic on that half plane. Then $\gstar_{c,D}$ is holomorphic in $\Re z<c$ and equal there
to
\begin{equation}
 \gstar_{c,D}(z)=1+D-\frac{cD}{c-z}-\frac{z\,\Lambda(c-z)}{(c-z)^{2}} .
 \label{eq:ex-lac-transform}
\end{equation}
In particular $\gstar_{c,D}(0)=1$.
\end{lemma}

\begin{proof}
On $(2^{m-1},2^{m}]$ with $m\ge1$ one has $N(u)=m$ and $S(u)=mu-2^{m}+1$, while $N=S=0$ on
$(0,1]$. So $g_{c,D}$ is piecewise smooth in $u$, and differentiating
$1+e^{-cu}(D+S(u))$ with $dt=-t\,du$ gives \eqref{eq:ex-lac-derivative}. For $u\ge1$ the
definition of $N$ gives $N(u)\le1+\log_{2}u$ and $S(u)\le N(u)u$, so the bracket in
\eqref{eq:ex-lac-derivative} has modulus at most $cD+(cu+1)(1+\log_{2}u)$, and it equals $cD$
for $u<1$. Writing $t^{c-1}=t^{a-1}e^{-(c-a)u}$ with $a<c$, the exponential absorbs that
bracket and leaves the stated bound. It is integrable on $(0,1)$, and $g_{c,D}(t)\to1$ as
$t\to0$, so $g_{c,D}$ is absolutely continuous on $[0,1]$.

Absolute continuity allows the integration by parts used in the proof of
Theorem~\ref{thm:ex-fgv-existence},
\[
 \gstar_{c,D}(z)=-z\int_{0}^{1}g_{c,D}(t)t^{-z-1}\,dt
 =g_{c,D}(1)-\int_{0}^{1}t^{-z}g_{c,D}'(t)\,dt
 \qquad(\Re z<0),
\]
and $g_{c,D}(1)=1+D$ because $S(0)=0$. Substituting $t=e^{-u}$ and writing $w=c-z$,
\[
 \int_{0}^{1}t^{-z}g_{c,D}'(t)\,dt
 =\int_{0}^{\infty}e^{-wu}\bigl[cD+cS(u)-N(u)\bigr]\,du .
\]
The terms of $S$ and $N$ are nonnegative and $\Re w>0$, so their moduli have a finite total
integral and the sums may be integrated term by term. Then
\[
 \int_{0}^{\infty}e^{-wu}\,du=\frac1w,
 \qquad
 \int_{0}^{\infty}e^{-wu}(u-b)_+\,du=\frac{e^{-wb}}{w^{2}},
 \qquad
 \int_{0}^{\infty}e^{-wu}\mathbf 1_{u>b}\,du=\frac{e^{-wb}}{w}
\]
give $cD/w+(c/w^{2}-1/w)\Lambda(w)$. Since $c-w=z$, this is $cD/w+z\Lambda(w)/w^{2}$, which
proves \eqref{eq:ex-lac-transform} on $\Re z<0$. The series \eqref{eq:ex-lac-Lambda} converges
uniformly on $\Re w\ge\sigma$ for every $\sigma>0$, so $\Lambda$ is holomorphic on $\Re w>0$,
and $z\mapsto c-z$ carries $\Re z<c$ onto that half plane. The right side of
\eqref{eq:ex-lac-transform} is therefore holomorphic in $\Re z<c$ and supplies the
continuation. At $z=0$ it equals $1+D-D=1$.
\end{proof}

The monotone members of the family fall under the existence theorem, and the threshold of
monotonicity is explicit.

\begin{theorem}[The monotone members are regular]
\label{thm:ex-lacunary}
Let $c>0$ and put
\begin{equation}
 D_{c}=\max_{m\ge1}\Bigl(\frac mc-1-(m-2)2^{m-1}\Bigr),
 \label{eq:ex-lac-Dc}
\end{equation}
a maximum attained at the least $m$ with $m2^{m-1}\ge1/c$. The profile $g_{c,D}$ is
nondecreasing on $[0,1]$ if and only if $D\ge D_{c}$. For such a $D$ the equation
\begin{equation}
 \frac{cD}{c-\gamma}+\frac{\gamma\,\Lambda(c-\gamma)}{(c-\gamma)^{2}}=1+D
 \label{eq:ex-lac-crossing}
\end{equation}
has a unique root $\gamma_{c,D}$ in $(0,c)$, the kernel $G(n,k)=g_{c,D}(k/n)$ is a RAF, and
\[
 \alpha(g_{c,D})=\eta(g_{c,D})=\gamma_{c,D} .
\]
For $c=1/2$ the threshold is $D_{1/2}=3$ and
\begin{equation}
 \gamma_{1/2,3}=0.07135125878143746087\ldots
 \label{eq:ex-lac-gamma}
\end{equation}
\end{theorem}

\begin{proof}
Put $\varphi(u)=c(D+S(u))-N(u)$, so that $g_{c,D}'(t)=t^{c-1}\varphi(u)$ by
\eqref{eq:ex-lac-derivative} and monotonicity is the nonnegativity of $\varphi$. On $(0,1]$ one
has $\varphi=cD$. On $(2^{m-1},2^{m}]$ with $m\ge1$, the value $S(u)=mu-2^{m}+1$ makes
$\varphi$ affine with positive slope $cm$, so its infimum over that interval is its right limit
at $2^{m-1}$,
\[
 \varphi(2^{m-1}{+})=c\bigl(D+(m-2)2^{m-1}+1\bigr)-m ,
\]
nonnegative exactly when $D\ge h(m):=m/c-1-(m-2)2^{m-1}$. Since $h(1)=1/c>0$, the condition
$cD\ge0$ on $(0,1]$ is implied. Now $h(m+1)-h(m)=1/c-m2^{m-1}$ and $m\mapsto m2^{m-1}$ is
strictly increasing, so $h$ increases up to the least $m$ with $m2^{m-1}\ge1/c$ and does not
increase after it. The maximum \eqref{eq:ex-lac-Dc} is attained there and finite, and
$g_{c,D}$ is nondecreasing if and only if $D\ge D_{c}$. For $c=1/2$ that least $m$ is $2$,
where $h(2)=4-1-0=3$.

Let $D\ge D_{c}$. Then $g_{c,D}(0)=1>0$, and $0\le g_{c,D}'(t)\le C_at^{a-1}$ for every
$a\in(0,c)$ by Lemma~\ref{lem:ex-lac-transform}. The crossing function of
Theorem~\ref{thm:ex-fgv-existence} was computed in the proof of that lemma,
\[
 F(s)=\int_{0}^{1}t^{-s}g_{c,D}'(t)\,dt
 =\frac{cD}{c-s}+\frac{s\,\Lambda(c-s)}{(c-s)^{2}}
 \qquad(s<c).
\]
It is continuous and strictly increasing there, the integrand being nonnegative and not almost
everywhere zero, with $F(0)=g_{c,D}(1)-g_{c,D}(0)=D$ and, since $D\ge D_c>0$, with
$F(s)\to\infty$ as $s\uparrow c$. As $g_{c,D}(1)=1+D$ exceeds $F(0)$, the crossing equation
\eqref{eq:ex-lac-crossing} has a unique root $\gamma_{c,D}$, and it lies in $(0,c)$. Choosing
$a\in(\gamma_{c,D},c)$ puts the profile under Theorem~\ref{thm:ex-fgv-existence}, which gives
the conclusion. The constant \eqref{eq:ex-lac-gamma} is the root for $c=1/2$ and $D=3$. The
terms of \eqref{eq:ex-lac-Lambda} decay doubly exponentially, and eight of them already fix the
printed digits.
\end{proof}

\begin{theorem}[The edge is a natural boundary]
\label{thm:ex-lac-boundary}
Let $c>0$ and $D\ge0$. The function $\gstar_{c,D}$ extends holomorphically to no neighborhood
of any point of the line $\Re z=c$, which is therefore a natural
boundary\index[terms]{natural boundary}.
\end{theorem}

\begin{proof}
The exponents $2^{j}$ are integers, so the substitution $q=e^{-w}$ turns
\eqref{eq:ex-lac-Lambda} into the lacunary series
\[
 \lambda(q)=\sum_{j\ge0}q^{2^{j}}
 \qquad(|q|<1),
 \qquad
 \Lambda(w)=\lambda(e^{-w}) .
\]

First, $\lambda$ is unbounded along the radius through every dyadic root of unity. Let
$q_{0}=e^{2\pi ip/2^{k}}$ with $p$ and $k\ge0$ integers. For $j\ge k$ the exponent $2^{j}$ is a
multiple of $2^{k}$, so $q_{0}^{2^{j}}=1$ and
\[
 \lambda(rq_{0})=\sum_{j<k}(rq_{0})^{2^{j}}+\sum_{j\ge k}r^{2^{j}}
 \qquad(0<r<1),
\]
the first sum having modulus at most $k$. Take $r=1-1/M$ with an integer $M\ge2$. Each $j$ with
$2^{j}\le M$ contributes $r^{2^{j}}\ge(1-1/M)^{M}\ge1/4$, and there are
$\lfloor\log_{2}M\rfloor+1$ such $j$, at most $k$ of them below $k$. Hence
\[
 |\lambda(rq_{0})|\ge\frac{\lfloor\log_{2}M\rfloor+1-k}{4}-k
 \longrightarrow\infty
 \qquad(M\to\infty),
\]
so $\lambda$ is unbounded near $q_{0}$ and extends holomorphically to no neighborhood of it.
The dyadic roots of unity are dense in the unit circle, and the set of boundary points across
which $\lambda$ does extend is open in that circle, so it is empty. This is the gap phenomenon
of Hadamard\index[names]{Hadamard, J.}~\cite{Hadamard1892}.

The same holds for $\Lambda$ on the imaginary axis. If $\Lambda$ extended holomorphically to a
neighborhood of a point $w_{0}$ with $\Re w_{0}=0$, then, $w\mapsto e^{-w}$ being a local
biholomorphism there, composition with its local inverse would extend $\lambda$ to a
neighborhood of $q_{0}=e^{-w_{0}}$, which the previous paragraph excludes.

Let now $\Re z_{0}=c$. If $z_{0}=c$, then for real $s\uparrow c$ the term $j=0$ makes
$\Lambda(c-s)\ge e^{-1}$ once $c-s\le1$, so \eqref{eq:ex-lac-transform} sends
$\gstar_{c,D}(s)$ to $-\infty$ and no holomorphic extension exists at $c$. If $z_{0}\neq c$,
choose $\eps<\min(|z_{0}-c|,c)$. On the disk of radius $\eps$ about $z_{0}$ one has $z\neq c$,
and $|z|\ge\Re z>c-\eps>0$, so \eqref{eq:ex-lac-transform} can be solved for
\[
 \Lambda(c-z)=\frac{(c-z)^{2}}{z}
 \Bigl(1+D-\frac{cD}{c-z}-\gstar_{c,D}(z)\Bigr).
\]
An extension of $\gstar_{c,D}$ across $z_{0}$ would extend $\Lambda$ across $c-z_{0}$, a point
of the imaginary axis.
\end{proof}

\begin{remark}[Regularity without global continuation]
\label{rem:ex-lac-meaning}
For $D\ge D_{c}$ the two theorems bear on the same profile. Its transform is holomorphic in
$\Re z<c$, has its first zero at $\gamma_{c,D}<c$, and extends past no point of $\Re z=c$,
while the arithmetic index exists and equals that zero. Arithmetic regularity therefore does
not require a meromorphic continuation of the transform to the plane. The criterion of
\S\ref{sec:existence} reads the profile and its derivative, and the transform enters only to
name the zero.
\end{remark}

\begin{remark}[The original profile]
\label{rem:ex-lacunary-original}
The member $c=1/2$, $D=0$ is the profile Tenenbaum\index[names]{Tenenbaum, G.} constructed. It
falls outside Theorem~\ref{thm:ex-lacunary}, since $D_{1/2}=3$, and it is bimodal in the
variable $u$: the bracket $\varphi$ of that proof is $S(u)/2-N(u)$, which vanishes
identically on $(0,1]$ and, beyond $u=1$, exactly at $u=7/2$ and at $u=13/3$. So $g_{1/2,0}$
stays at $1$ up to $u=1$, rises to a maximum at $u=7/2$, falls to a local minimum at $u=4$,
rises to a second local maximum at $u=13/3$ and decreases to $1$ beyond. Adding the four monotone increments, with $S(7/2)=4$, $S(4)=5$ and $S(13/3)=6$, gives
the total variation
\begin{equation}
 \int_{0}^{1}|g_{1/2,0}'(t)|\,dt
 =8e^{-7/4}+12e^{-13/6}-10e^{-2}
 =1.41154484\ldots,
 \label{eq:ex-lac-variation}
\end{equation}
larger than $|g_{1/2,0}(1)|=1$. The weighted variation $V$ of Open
Problem~\ref{op:ex-fgv-jumps} already exceeds that value at $s=0$, its threshold being
$\widetilde\gamma_V=-0.07466576906\ldots$, so that route reaches no positive exponent here. Theorem~\ref{thm:ex-lac-boundary} applies unchanged,
and Lemma~\ref{lem:ex-lac-transform} gives
\[
 \gstar_{1/2,0}(z)=1-\frac{z\,\Lambda(\tfrac12-z)}{(\tfrac12-z)^{2}} .
\]
On the real axis $\Lambda(\tfrac12-s)$ and $(\tfrac12-s)^{-2}$ are positive and increasing in
$s<\tfrac12$, so the subtracted term is nonpositive for $s\le0$ and increases strictly from $0$
to $+\infty$ on $[0,\tfrac12)$. Hence $\gstar_{1/2,0}$ is at least $1$ on $s\le0$ and decreases
strictly from $1$ to $-\infty$ on $[0,\tfrac12)$, and it has exactly one real zero in
$\Re z<\tfrac12$, at
\begin{equation}
 \gamma_{1/2,0}=0.10960172249149471740\ldots
 \label{eq:ex-lac-gamma-original}
\end{equation}
\end{remark}

\begin{openproblem}[The original lacunary profile]
\label{op:ex-lacunary-original}
Let $g_{1/2,0}$ be the profile of Remark~\ref{rem:ex-lacunary-original} and $\gamma_{1/2,0}$
its real zero \eqref{eq:ex-lac-gamma-original}.
\begin{enumerate}[label=\textup{(\roman*)}]
\item Is $\gamma_{1/2,0}$ the leftmost zero of $\gstar_{1/2,0}$, so that
$\eta(g_{1/2,0})=\gamma_{1/2,0}$? A numerical computation finds a pair of zeros near
$0.3842\pm1.1016i$, and a numerical search finds none with $\Re z\le\gamma_{1/2,0}$, which is not
a proof that none exists. A proof has to confine any such zero to a compact set and then clear
that set.
\item Is $g_{1/2,0}$ a function of good variation, and if it is, does $\alpha(g_{1/2,0})$ equal
$\gamma_{1/2,0}$?
\end{enumerate}
The natural boundary settles neither question. Since $\gamma_{1/2,0}<\tfrac12$, a strip about
the candidate zero lies inside the half plane of holomorphy, so
Conditional Theorem~\ref{cthm:ex-alpha-eta} is not blocked by the boundary. What is missing for
it is the arithmetic input it assumes, transparency below the candidate and absorption at and
above it, the zero-free half plane of (i), and the three hypotheses that theorem isolates, the
exact Mellin-Perron inversion, the two growth bounds and the nonvanishing of the numerator at
the candidate.
\end{openproblem}

\section{The scalar reduction}
\label{sec:scalar_reduction}

Many kernels of the theory reduce, after one Abel summation, to a first order recurrence
for the partial sums. The following theorem packages that situation once and for all, and
the affine kernel calibrates it.

\begin{proposition}\label{thm:ex-first-order}
Let $\gamma,\delta>0$.  After a fixed initial rank $n_0$, suppose that
the power-forced solutions satisfy
\[
 A_\beta(n)=u_nA_\beta(n-1)+v_{\beta,n},
\]
where $u_n>0$ and
\begin{align*}
 u_n&=1-\frac\gamma n+\mathcal O(n^{-1-\delta}),
 \\
 v_{\beta,n}&=b_\beta n^{-\beta-1}
 +\mathcal O_\beta(n^{-\beta-1-\delta}).
 \end{align*}
Put $\theta=\min(\delta,1)$ and
\[
 P_n=\prod_{j=n_0+1}^{n}u_j,
 \qquad
 q_n=n^\gamma P_n.
\]
Then $q_n$ tends to a number $q_\infty>0$, with
\begin{equation}
 q_n=q_\infty\left(1+\mathcal O(n^{-\theta})\right),
 \label{eq:ex-first-order-q}
\end{equation}
and
\begin{equation}
 A_\beta(n)=n^{-\gamma}q_n
 \left[
 A_\beta(n_0)+
 \sum_{m=n_0+1}^{n}\frac{m^\gamma v_{\beta,m}}{q_m}
 \right].
 \label{eq:ex-first-order-variation}
\end{equation}
The three forcing cases are
\begin{align}
 A_\beta(n)&=\frac{b_\beta}{\gamma-\beta}n^{-\beta}
 +o_\beta(n^{-\beta})
 &&(\beta<\gamma),
 \label{eq:ex-first-order-below}\\
 A_\gamma(n)&=b_\gamma n^{-\gamma}\log n
 +C_\gamma n^{-\gamma}
 +\mathcal O(n^{-\gamma-\theta}\log n),
 \label{eq:ex-first-order-critical}\\
 A_\beta(n)&=C_\beta Q_n n^{-\gamma}
 +\frac{b_\beta}{\gamma-\beta}n^{-\beta}
 +\mathcal O_\beta(n^{-\beta-\theta})
 &&(\beta>\gamma),
 \label{eq:ex-first-order-above}
\end{align}
where $Q_n=q_n/q_\infty$ and
\begin{equation}
 C_\beta=q_\infty
 \left[
 A_\beta(n_0)+
 \sum_{m=n_0+1}^{\infty}
 \frac{m^\gamma v_{\beta,m}}{q_m}
 \right]
 \qquad(\beta>\gamma).
 \label{eq:ex-first-order-connection}
\end{equation}

Suppose this recurrence comes from \eqref{eq:defining_relation_gen}, and suppose
\begin{equation}
 \frac{b_\beta}{\gamma-\beta}=\Xi_G(\beta)
 \qquad(\beta<\gamma).
 \label{eq:ex-first-order-match}
\end{equation}
Then the additive transparency relation \eqref{eq:transparency_additive} holds
below $\gamma$, including any exponent where $\Xi_G(\beta)=0$, while
absorption holds at the candidate exponent.  If there are
$\beta_j\downarrow\gamma$ with $C_{\beta_j}\neq0$ and finite nonzero
$G^*(\beta_j)$, then sharpness holds, the kernel is a RAF, and
$\alpha(G)=\gamma$.
\end{proposition}

\begin{proof}
The exact product iteration gives
\[
 A_\beta(n)=P_n
 \left[A_\beta(n_0)+
 \sum_{m=n_0+1}^{n}\frac{v_{\beta,m}}{P_m}\right],
\]
which is \eqref{eq:ex-first-order-variation}.  Moreover,
\[
 \log\frac{q_n}{q_{n-1}}
 =\gamma\log\frac n{n-1}+\log u_n
 =\mathcal O(n^{-1-\theta}).
\]
The series of logarithms converges.  Positivity of $u_n$ makes its limit
finite and positive, which proves \eqref{eq:ex-first-order-q}.

Set
\[
 r_{\beta,m}=
 \frac{m^\gamma v_{\beta,m}}{q_m}
 -\frac{b_\beta}{q_\infty}m^{\gamma-\beta-1}.
\]
Then
\[
 r_{\beta,m}=\mathcal O_\beta
 \left(m^{\gamma-\beta-1-\theta}\right).
\]
If $\beta<\gamma$, power summation in
\eqref{eq:ex-first-order-variation} proves
\eqref{eq:ex-first-order-below}.  At equality, the main sum is harmonic and
the error sum converges.  This proves \eqref{eq:ex-first-order-critical},
with $C_\gamma$ equal to the resulting finite part.

If $\beta>\gamma$, the series in
\eqref{eq:ex-first-order-connection} converges absolutely.  Subtract its tail
from the finite sum in \eqref{eq:ex-first-order-variation}.  The estimate
\[
 \sum_{m>n}m^{\gamma-\beta-1}
 =\frac{n^{\gamma-\beta}}{\beta-\gamma}
 +\mathcal O(n^{\gamma-\beta-1})
\]
gives \eqref{eq:ex-first-order-above}.  Retaining $Q_n$ is necessary when
the correction to the homogeneous product is larger than the forced term.

The coefficient match \eqref{eq:ex-first-order-match} gives additive
transparency.  If its value is zero, the first line of
\eqref{eq:ex-first-order-below} states exactly that
$A_\beta(n)=o(n^{-\beta})$.
The critical logarithm is bounded by $n^{-\gamma+\eps}$, and the
above-threshold formula is $\mathcal O(n^{-\gamma})$.  A nonzero connection
at exponents decreasing to $\gamma$ prevents transparency at all of them.
This is sharpness at $\gamma$, and the absorption estimate proves the RAF
assertion.
\end{proof}

The affine kernel is the case where every constant of the reduction is explicit, and it
calibrates the statement above.

\begin{corollary}
\label{cor:ex-affine}
For
\[
 g_\lambda(t)=\lambda+(1-\lambda)t,
 \qquad
 0<\lambda<1,
\]
the exact recurrence is
\begin{equation}
 A(n)=
 \left(1-\frac\lambda n\right)A(n-1)
 +\frac{n^{1-\beta}-(n-1)^{1-\beta}}{n}.
 \label{eq:ex-affine-recurrence}
\end{equation}
It satisfies Proposition~\ref{thm:ex-first-order} with
\begin{equation}
 \gamma=\lambda,
 \qquad
 b_\beta=1-\beta,
 \qquad
 g_\lambda^*(z)=\frac{\lambda-z}{1-z}.
 \label{eq:ex-affine-data}
\end{equation}
Consequently,
\[
 \alpha(g_\lambda)=\eta(g_\lambda)=\lambda.
\]
\end{corollary}

\begin{proof}
Discrete Abel summation gives
\[
 A(n)-\frac{1-\lambda}{n}
 \sum_{k<n}A(k)=n^{-\beta}.
\]
Subtracting consecutive ranks after multiplication by the rank gives
\eqref{eq:ex-affine-recurrence}.  Its forcing increment is
\[
 \frac{n^{1-\beta}-(n-1)^{1-\beta}}n
 =(1-\beta)n^{-\beta-1}
 +\mathcal O_\beta(n^{-\beta-2}).
\]
Direct integration gives the transform in \eqref{eq:ex-affine-data}, and
\[
 \frac{1-\beta}{\lambda-\beta}
 =\frac1{g_\lambda^*(\beta)}.
\]
For $\lambda<\beta<1$, the recurrence has positive multiplier, positive
forcing increment, and positive initial value.  Its connection coefficient
is positive.  Such exponents decrease to $\lambda$, so
Proposition~\ref{thm:ex-first-order} proves the arithmetic index.  The transform
has its first zero at $\lambda$, which proves the analytic index separately.
\end{proof}

\section{An algebra of indices}
\label{sec:index_algebra}

The positive class of Theorem~\ref{thm:ex-fgv-existence} is stable under products and
under positive sums, and the index moves in a controlled way. The boundary term at the
origin drives both statements.

\begin{lemma}
\label{lem:ex-zero-boundary-transform}
Let $g$ be absolutely continuous on $[0,1]$, and suppose that
$|g'(t)|\leq Ct^{a-1}$ almost everywhere for some $a>0$.  The
continuation to $\Re z<a$ satisfies
\begin{equation}
 g^*(z)=g(1)-\int_0^1t^{-z}g'(t)\,dt
 =g(0)-z\int_0^1\bigl(g(t)-g(0)\bigr)t^{-z-1}\,dt.
 \label{eq:ex-zero-boundary-transform}
\end{equation}
In particular, $g^*(0)=g(0)$.
\end{lemma}

\begin{proof}
Put $u(t)=g(t)-g(0)$.  Absolute continuity and the derivative bound give
\[
 |u(t)|\leq\frac Ca t^a.
\]
Hence $u(t)t^{-z}\to0$ at zero and the second integral in
\eqref{eq:ex-zero-boundary-transform} converges absolutely when
$\Re z<a$.  Integration by parts on $[\eps,1]$ gives
\[
 -z\int_\eps^1u(t)t^{-z-1}\,dt
 =u(1)-u(\eps)\eps^{-z}
 -\int_\eps^1t^{-z}g'(t)\,dt.
\]
Letting $\eps$ tend to zero and adding $g(0)$ proves both
representations.  They agree with the defining transform on $\Re z<0$, so
they give its continuation.
\end{proof}

The class is stable under multiplication, and the statement below says how the index moves under
a product.

\begin{proposition}\label{thm:ex-positive-product}
For $i=1,2$, let $g_i$ be absolutely continuous on $[0,1]$ and
suppose
\begin{equation}
 g_i(0)>0,
 \qquad
 0\leq g_i'(t)\leq C_it^{a_i-1}
 \label{eq:ex-product-hypotheses}
\end{equation}
almost everywhere.  Assume that $0<\gamma_i<a_i$ and
$g_i^*(\gamma_i)=0$.  Put $h=g_1g_2$.  Then $h$ satisfies the
hypotheses of Theorem~\ref{thm:ex-fgv-existence} with
$a=\min(a_1,a_2)$.  Its transform has a unique real zero $\gamma_h$ in
\[
 0<\gamma_h<\min(\gamma_1,\gamma_2).
\]
Consequently,
\begin{equation}
 \alpha(g_1g_2)=\eta(g_1g_2)=\gamma_h
 <\min\bigl(\alpha(g_1),\alpha(g_2)\bigr).
 \label{eq:ex-product-index}
\end{equation}
\end{proposition}

\begin{proof}
Each $g_i$ is positive and increasing.  Thus
\[
 0\leq h'(t)=g_1'(t)g_2(t)+g_1(t)g_2'(t)
 \leq\bigl(C_1g_2(1)+C_2g_1(1)\bigr)t^{a-1}.
\]
Assume $\gamma_1\leq\gamma_2$.  Then
$\gamma_1<a_1$ and $\gamma_1\leq\gamma_2<a_2$, so
$\gamma_1<a$.  Put $d_i=g_i(0)$ and $u_i=g_i-d_i$.  The transform
lemma and $g_1^*(\gamma_1)=0$ give
\[
 d_1=\gamma_1\int_0^1u_1(t)t^{-\gamma_1-1}\,dt.
\]
Since
$h-d_1d_2=d_1u_2+d_2u_1+u_1u_2$, the $d_2u_1$ term cancels and
\[
 h^*(\gamma_1)=-\gamma_1\int_0^1
 \bigl(d_1u_2(t)+u_1(t)u_2(t)\bigr)t^{-\gamma_1-1}\,dt<0.
\]
The crossing hypothesis makes $g_2$ nonconstant.  Thus $u_2$ is
positive on a set of positive measure, which proves the strict sign.  Also
$h^*(0)=d_1d_2>0$.  A zero lies in $(0,\gamma_1)$.  For real $s<a$,
\[
 (h^*)'(s)=-\int_0^1(-\log t)t^{-s}h'(t)\,dt<0.
\]
The zero is unique.  Theorem~\ref{thm:ex-fgv-existence} proves
\eqref{eq:ex-product-index}.

Equation \eqref{eq:ex-constant-factor} follows by scaling the triangular
equation and the transform.  The finite-product statement follows by
induction.  At each multiplication the new crossing lies below the previous
crossing and the next factor crossing.  It therefore remains in every domain
needed by the next application of the theorem.
\end{proof}

If $C>0$, then
\begin{equation}
 (Cg)^*=Cg^*,
 \qquad
 A_\beta^{Cg}=C^{-1}A_\beta^g,
 \qquad
 \Xi_{Cg}=C^{-1}\Xi_g.
 \label{eq:ex-constant-factor}
\end{equation}
Thus a positive constant factor preserves every defined arithmetic and
analytic index.  This is the exceptional product case.  Finite products with
at least two nonconstant factors, each satisfying all factor hypotheses of
Proposition~\ref{thm:ex-positive-product}, have index strictly below the least
factor index.  Positive constant factors may be removed.  In particular,
\[
 \alpha(g^r)=\eta(g^r)<\alpha(g)
 \qquad(r\geq2)
\]
for every nonconstant $g$ covered by the theorem.

For a sum the index moves the other way, and the next statement locates it between the two.

\begin{proposition}\label{thm:ex-positive-sum}
Let $g_1,g_2$ satisfy \eqref{eq:ex-product-hypotheses}.  Assume that
$0<\gamma_i<a_i$ and $g_i^*(\gamma_i)=0$.  Let
$\lambda_1,\lambda_2>0$, and assume
the common-domain condition
\[
 \max(\gamma_1,\gamma_2)<a:=\min(a_1,a_2).
\]
For $h=\lambda_1g_1+\lambda_2g_2$, there is a unique real crossing
$\gamma_h$, and
\[
 \min(\gamma_1,\gamma_2)
 \leq\gamma_h\leq
 \max(\gamma_1,\gamma_2).
\]
Both inequalities are strict when $\gamma_1\neq\gamma_2$.  When the two
indices agree, $\gamma_h$ equals their common value.  Therefore
\begin{equation}
 \alpha(h)=\eta(h)=\gamma_h.
 \label{eq:ex-sum-index}
\end{equation}
The same result holds for a finite positive sum under
$\max_j\gamma_j<\min_j a_j$.  Its index lies between the least and greatest
summand indices and lies strictly between them unless all indices agree.  No
barycentric formula is asserted.
\end{proposition}

\begin{proof}
The differential hypotheses follow from positivity and linearity.  The same
linearity gives
\begin{equation}
 h^*(s)=\lambda_1g_1^*(s)+\lambda_2g_2^*(s).
 \label{eq:ex-sum-transform}
\end{equation}
Every $g_i^*$ is strictly decreasing on its real continuation interval
because
\[
 (g_i^*)'(s)=-\int_0^1(-\log t)t^{-s}g_i'(t)\,dt<0.
\]
Suppose $\gamma_1<\gamma_2$.  The common-domain condition permits both
endpoint evaluations and gives
\[
 h^*(\gamma_1)=\lambda_2g_2^*(\gamma_1)>0,
 \qquad
 h^*(\gamma_2)=\lambda_1g_1^*(\gamma_2)<0.
\]
Strict monotonicity places the unique zero between the endpoints.  If the
crossings agree, \eqref{eq:ex-sum-transform} vanishes at their common value.
Theorem~\ref{thm:ex-fgv-existence} proves \eqref{eq:ex-sum-index}.

For finitely many summands, evaluate the sum transform at the least and
greatest crossings.  All terms have the required weak sign.  At least one
term has a strict sign when the crossings are not all equal.  This proves the
finite statement.
\end{proof}

Multiplication by the coordinate is the simplest instance of the product rule and is worth
recording on its own.

\begin{corollary}
\label{cor:ex-coordinate-product}
Let $g$ be a bounded Riemann-integrable FGV that is a RAF with finite
$\alpha(g)>0$, and put $h(t)=tg(t)$.  Then $h$ is a RAF and
\[
 \alpha(h)=0.
\]
More precisely, let $B_\delta$ denote the partial sum for $g$ at exponent
$\delta$.  The partial sum for $h$ satisfies the exact identity
\begin{equation}
 A_\beta^h(n)=\frac{B_{\beta-1}(n)}n
 +\sum_{m<n}\frac{B_{\beta-1}(m)}{m(m+1)}.
 \label{eq:ex-coordinate-abel}
\end{equation}
The transform relation is
\begin{equation}
 h^*(z)=\frac z{z-1}g^*(z-1).
 \label{eq:ex-coordinate-transform}
\end{equation}
For $\beta<0$,
\begin{equation}
 A_\beta^h(n)\sim
 \frac{\beta-1}{\beta g^*(\beta-1)}n^{-\beta}
 =\Xi_h(\beta)n^{-\beta}.
 \label{eq:ex-coordinate-below}
\end{equation}
At zero,
\begin{equation}
 A_0^h(n)\sim\frac1{g^*(-1)}\log n,
 \qquad
 g^*(-1)\neq0.
 \label{eq:ex-coordinate-critical}
\end{equation}
For every $\beta>0$, the sequence $A_\beta^h(n)$ has a finite limit and
is $\mathcal O_\beta(1)$.
\end{corollary}

\begin{proof}
Let $(a_k)$ solve the equation for $h$ at exponent $\beta$, and put
$b_k=ka_k$.  Multiplication of that equation by $n$ gives
\[
 \sum_{k\leq n}b_kg(k/n)=n^{1-\beta}.
\]
Triangular uniqueness identifies $(b_k)$ with the coefficient sequence for
$g$ at exponent $\delta=\beta-1$.  Abel summation of
$A_\beta^h(n)=\sum_{k\leq n}b_k/k$ proves
\eqref{eq:ex-coordinate-abel}.  Direct integration in $\Re z<1$ gives
\eqref{eq:ex-coordinate-transform}, and continuation preserves the identity.

Suppose first that $\beta<0$.  Since $\beta-1<\alpha(g)$, transparency
for $g$ gives
\[
 B_{\beta-1}(n)=c_\beta n^{1-\beta}+o(n^{1-\beta}),
 \qquad
 c_\beta=\frac1{g^*(\beta-1)}.
\]
The coefficient is nonzero.  The transform is holomorphic in this left
half plane, and a zero would not be a transparent point.  Abel summation and
the elementary power sum give
\[
 \sum_{m<n}\frac{B_{\beta-1}(m)}{m(m+1)}
 =-\frac{c_\beta}{\beta}n^{-\beta}+o(n^{-\beta}).
\]
Combining this with the boundary term proves
\eqref{eq:ex-coordinate-below}.

At $\beta=0$, transparency of $g$ at $-1$ gives
\[
 B_{-1}(n)=cn+o(n),
 \qquad
 c=\frac1{g^*(-1)}\neq0.
\]
The harmonic form of \eqref{eq:ex-coordinate-abel} proves
\eqref{eq:ex-coordinate-critical}.  Formula
\eqref{eq:ex-coordinate-transform} shows that $h^*$ has a zero at zero,
so zero is not transparent.

Now let $\beta>0$.  If $\beta-1<\alpha(g)$, transparency gives
$B_{\beta-1}(m)=\mathcal O_\beta(m^{1-\beta})$.  The series in
\eqref{eq:ex-coordinate-abel} converges absolutely.  If
$\beta-1\geq\alpha(g)$, choose $0<\eps<\alpha(g)$.  Absorption
gives
\[
 B_{\beta-1}(m)=
 \mathcal O_{\beta,\eps}
 \left(m^{-\alpha(g)+\eps}\right),
\]
and the series again converges absolutely.  The boundary term tends to zero
in both cases.  This proves the bounded absorbed branch.  Subcritical
transparency, absorption for every $\beta\geq0$, and failure at zero are
the three conditions for $\alpha(h)=0$.
\end{proof}

\begin{remark}[The passage to $xg$ does not reverse]
\label{rem:ex-coordinate-converse}
Nothing established here reverses the passage from $g$ to $xg$. Knowing that
$xg$ is a function of good variation of index zero constrains the leading term
of the partial sums of $g$ and leaves its absorbed branch free. This is recorded
as a direction and not as a problem, since any statement of it would have to
quantify the secondary asymptotic required of the absorbed branch, and no
candidate quantification is available.
\end{remark}

\section{From the analytic to the arithmetic index}
\label{sec:alpha_eta_transfer}

The equality $\alpha(G)=\eta(G)$ is a theorem wherever it is proved and nowhere a
definition. The positive class of \S\ref{sec:existence} realizes it unconditionally. Beyond
that class the passage from a first zero of the transform to the index is a contour
argument whose every hypothesis must be stated, and the following conditional theorem
lists them.

\begin{conditionaltheorem}[Mellin-Perron transfer from a first zero]
\label{cthm:ex-alpha-eta}
Assume that Theorem~\ref{thm:ex-green-criterion} has already supplied
transparency for every $\beta<\gamma$ and absorption
$A_\beta(n)=\mathcal O_{\beta,\eps}(n^{-\gamma+\eps})$ for every
$\beta\ge\gamma$ and every $\eps>0$, where $0<\gamma<\infty$.
Assume that $G^*$ has a meromorphic continuation to
$\{\Re z<\gamma+\delta\}$ for some $\delta>0$, has no zero in
$\Re z<\gamma$, and has a simple real zero at $z=\gamma$.

Let $\gamma<\beta_j<\gamma+\delta$ decrease to $\gamma$ and suppose that
$G^*(\beta_j)$ is finite and nonzero. Put $c_j=1/G^*(\beta_j)$.
For every $j$, assume an exact Mellin-Perron representation
\begin{equation}
 A_{\beta_j}(n)-c_jn^{-\beta_j}
 =\lim_{T\to\infty}\frac1{2\pi i}
 \int_{\sigma_--iT}^{\sigma_-+iT}
 \frac{F_j(z)}{G^*(z)}n^{-z}\,dz,
 \label{eq:ex-perron}
\end{equation}
where
\[
 \sigma_-<\gamma<\sigma_{+,j}<\beta_j.
\]
The constant $c_j$ specifies the subtracted term. Its definition makes no
assertion of transparency at $\beta_j$.
Assume that $F_j$ is holomorphic on a neighborhood of the closed strip,
that $F_j(\gamma)\neq0$, and that $F_j/G^*$ has no other pole there.
Suppose there is $M_j\geq0$ such that, away from a fixed disk about $\gamma$,
\begin{align}
 |1/G^*(\sigma+it)|&\leq
 C_j(1+|t|)^{M_j},
 \label{eq:ex-perron-growth}\\
 |F_j(\sigma+it)|&\leq
 C_j(1+|t|)^{-M_j-2}
 \label{eq:ex-perron-decay}
\end{align}
uniformly for $\sigma_-\le\sigma\le\sigma_{+,j}$.
Then the contour can be moved from $\sigma_-$ to $\sigma_{+,j}$, and
\begin{equation}
 A_{\beta_j}(n)-c_jn^{-\beta_j}
 =-\frac{F_j(\gamma)}{(G^*)'(\gamma)}n^{-\gamma}
 +\mathcal O_j(n^{-\sigma_{+,j}}).
 \label{eq:ex-perron-residue}
\end{equation}
Therefore $G$ is a RAF and
\[
 \alpha(G)=\eta(G)=\gamma.
\]
The exact inversion, the two growth bounds and the nonzero numerator at
$\gamma$ are separate hypotheses to verify in an application.
\end{conditionaltheorem}

\begin{proof}
Fix $j$ and write $H_j=F_j/G^*$ and
\[
 I_{\sigma,T}(n)=\frac1{2\pi i}
 \int_{\sigma-iT}^{\sigma+iT}H_j(z)n^{-z}\,dz.
\]
The two growth bounds give $|H_j(\sigma+it)|\ll_j(1+|t|)^{-2}$
for all sufficiently large $|t|$, uniformly across the strip.
On each boundary line the remaining compact interval contains no pole,
so both infinite vertical integrals converge absolutely, with tail bounds
\[
 |I_{\sigma,\infty}(n)-I_{\sigma,T}(n)|
 \ll_j\frac{n^{-\sigma}}{T}
 \qquad(\sigma=\sigma_-,\ \sigma_{+,j}).
\]
For $n\ge2$, each horizontal side at height $\pm T$ has modulus at most
\[
 C_jT^{-2}\int_{\sigma_-}^{\sigma_{+,j}}n^{-\sigma}\,d\sigma
 =C_jT^{-2}\frac{n^{-\sigma_-}-n^{-\sigma_{+,j}}}{\log n}.
\]
Thus, with both vertical segments oriented upward, the residue theorem on
the finite rectangle gives
\[
 I_{\sigma_-,T}(n)
 =I_{\sigma_{+,j},T}(n)
  -\frac{F_j(\gamma)}{(G^*)'(\gamma)}n^{-\gamma}
  +\mathcal O_j\!\left(\frac{n^{-\sigma_-}}{T^2\log n}\right).
\]
The minus sign comes from moving an upward line to the right.
Letting $T\to\infty$ at each fixed $n$, the two tail bounds justify the
passage to the infinite lines. Moreover,
\[
 |I_{\sigma_{+,j},\infty}(n)|
 \le\frac{n^{-\sigma_{+,j}}}{2\pi}
       \int_{\mathbb R}|H_j(\sigma_{+,j}+it)|\,dt
 \ll_j n^{-\sigma_{+,j}}.
\]
This proves \eqref{eq:ex-perron-residue}. The residue coefficient is nonzero,
and both $n^{-\beta_j}$ and $n^{-\sigma_{+,j}}$ are $o(n^{-\gamma})$.
Consequently $n^{\beta_j}A_{\beta_j}(n)$ has no finite limit, so
$\beta_j$ is not transparent for any $j$. Every interval
$(\gamma,\gamma+\varepsilon)$ contains such an exponent, giving
$\tau(G)\le\gamma$. Subcritical transparency gives the reverse inequality,
and absorption at every $\beta\ge\gamma$ proves that $G$ is a RAF with
$\alpha(G)=\gamma$. The zero-free half plane $\Re z<\gamma$ and the zero at
$\gamma$ give $\eta(G)=\gamma$ by Definition~\ref{def:mellin}.
\end{proof}

The Ingham function and the orthorecursive kernel\index[terms]{orthorecursive kernel} stay outside the unconditional reach of
this section. For the Ingham function the identification $\alpha(\Phi)=\tfrac12$ is the
master equivalence\index[terms]{master equivalence} itself, and for the orthorecursive kernel the pointwise transfer keeps
its separate hypotheses, carried in Chapter~\ref{chap:ortho}.

\section{A gallery of two-variable arithmetic kernels}
\label{sec:raf_gallery}

The second gallery keeps nine entries. Appendices~\ref{app:H} to~\ref{app:N}
and~\ref{app:Q} contain
two-variable kernels $G(n,k)$ that do not reduce to a single ratio $k/n$, while
Appendix~\ref{app:O} returns to a ratio kernel with a diagonal jump\index[terms]{diagonal jump}, its whole difficulty
lying in its value at the diagonal. The table records for each the regularity index $\alpha$ where it
is proved, the analytic index $\eta$ where the transform has a zero, and a short remark.
The analytic index is marked $-$ where the transform has no zero, and an entry marked
open is exactly that, an open case of the theory.

\medskip
\renewcommand{\arraystretch}{1.45}
{\small
\begin{tabular}{cp{4.4cm}p{0.8cm}p{0.8cm}p{4.1cm}}
\toprule
{App.} & {Kernel $G(n,k)$} & $\alpha$ & $\eta$ & {Remark} \\
\midrule
\ref{app:H} & $\dfrac{n+k+x}{n+k+y}\cdot\dfrac{2n+y}{2n+x}$ & $2$ & $-$ & Rational kernel. $G^{*}\equiv1$. Exact feedback with telescoping weight. Proved for every $x\neq y$, by shifted impulses and Dirichlet uniqueness. \\
\ref{app:I} & $1$ if $2\mid n-k$, else $k/n$ & $1/2$ & $1/2$ & Parity kernel. Closed by a first order generating equation, the endpoint $-1$ treated by Laguerre and Fej\'er. Two affine channels. \\
\ref{app:J} & $\tfrac12\big(1+\gcd(n,k)/n\big)$ & open & $-$ & Divisor kernel. $G^{*}\equiv\tfrac12$. Exact double Dirichlet form, the section family being the obstruction. \\
\ref{app:K} & $\dfrac{n^{2}+k}{n^{2}+n}$ & $2$ & $-$ & Quadratic rational. $G^{*}\equiv1$. Closed Gamma formula. \\
\ref{app:L} & $\dfrac{\log(n+k)}{\log(2n)}$ & $1$ & $-$ & Logarithmic kernel. $G^{*}\equiv1$. Moving threshold $1-1/\log(2n)$, self-gauged. Positive Abel form, sharpness from positivity, proved. \\
\ref{app:M} & $\tfrac12\big(1+\tfrac{1+\sqrt{k}}{1+\sqrt{n}}\big)$ & $1/4$ & $1/4$ & Square-root kernel. Exact scalar reduction, proved. \\
\ref{app:N} & $\dfrac{\log(n^{2}+k^{2})}{\log(2n^{2})}$ & $2$ & $-$ & Quadratic logarithm. $G^{*}\equiv1$. Moving threshold $2-2/\log(2n^{2})$, self-gauged. The same mechanism read on squares, proved. \\
\ref{app:O} & $1-\{k/n\}$ & $\tfrac14$ & $-$ & Diagonal-jump kernel. Zero-free transform $\tfrac{1}{1-z}$. Laguerre basis and an oscillatory threshold, proved. Twin of the divisor problem. \\
\ref{app:Q} & $\tfrac12\big(1+\binom nk^{-1}\big)$ & $1$ & $-$ & Binomial harmonic kernel. $G^{*}\equiv\tfrac12$. The index comes from the two edge ranks alone, as the crossing of a forced and an edge asymptotic tower, proved. \\
\bottomrule
\end{tabular}
\renewcommand{\arraystretch}{1}
}
\medskip

\noindent Seven of these kernels have a transform without zeros, constant for
\ref{app:H}, \ref{app:J}, \ref{app:K}, \ref{app:L}, \ref{app:N} and \ref{app:Q}, and equal
to $1/(1-z)$ for \ref{app:O}, so that no analytic index is available to them and
their index, where it is known, comes from the equation itself. Every entry is proved by the exact analysis of its
appendix, except \ref{app:J}, whose index is open. For \ref{app:L} and \ref{app:N}
the index is the limit of a threshold that moves with the rank, and the limit is reached by
direct estimates on a positive feedback, sharpness coming from the positivity of the partial
sums alone. The index of \ref{app:J} alone is open, its appendix recording what is exact and
where the gap lies, the divisor kernel carrying an exact double
Dirichlet form whose section family is the obstruction. The kernels gathered here show the reach of the class, from
the rational and the logarithmic to the divisor, the square-root and the binomial, none of them
a function of good variation except the closing diagonal-jump kernel,
of index a quarter.

\chapter{General principles and exactly solvable kernels}
\label{chap:principles}

The definition of the preceding chapter fixes a property. It does not say how to establish it,
and on a general kernel there is no reason for the property to be decidable at all. This chapter
assembles what makes it decidable, the family
where every step is exact, the inversion that supplies the resolvent in closed form, the
classical equation the whole shape descends from, and the two principles that read an index
without a transform or refuse to let two kernels be multiplied.

The chapter is machinery throughout. The affine family is the case where every step is exact
and the two indices coincide. The exact Volterra\index[terms]{discrete Volterra} inversion supplies the resolvent\index[terms]{resolvent} in closed
form. The equidimensional equation of Euler\index[names]{Euler, L.} and Cauchy\index[names]{Cauchy, A.-L.} shows where the shape of the defining
equation\index[terms]{defining equation} comes from, one chapter after that shape was fixed. The homogeneity
principle\index[terms]{homogeneity principle} reads the index off the homogeneous solution when the transform has no zero and so
carries no analytic index\index[terms]{analytic index} at all. Products of kernels then show the multiplicative
obstruction\index[terms]{multiplicative obstruction} that blocks the multiplicative route to the Riemann hypothesis\index[terms]{Riemann hypothesis} and motivates the
gauge\index[terms]{gauge} deformations of Chapter~\ref{chap:gauge}.

\section{The archetype of smoothness, affine functions}
\label{sec:affine_archetype}

Nontrivial FGVs exist constructively, without complex integration, and the
affine family is the proven case of the theory.

\begin{theorem}\label{thm:affine}
Let $g(x)=(1-\lambda)x+\lambda$ with $\lambda\in(0,1)$. Then $g$ is a FGV of
index $\alpha(g)=\lambda$, and its arithmetic Mellin transform is
\[
g^*(z)=\frac{z-\lambda}{z-1},
\]
with a single zero at $z=\lambda$, so that $\alpha(g)=\eta(g)=\lambda$. It is
proved in \cite[Theorem~1.1]{Cloitre2016} and reproduced in full in
Appendix~\ref{app:A}.
\end{theorem}

The mechanism is worth displaying, because the whole problem reduces to a
first order recurrence.

\begin{proposition}\label{prop:affine_recurrence}
For the affine kernel\index[terms]{affine kernel} with forcing $n^{-\beta}$, the partial sums satisfy
exactly
\begin{equation}\label{eq:affine_recurrence}
A(n)=\Big(1-\frac{\lambda}{n}\Big)A(n-1)
+\frac{n^{1-\beta}-(n-1)^{1-\beta}}{n},\qquad n\ge2 .
\end{equation}
\end{proposition}

\begin{proof}
With $S(n):=\sum_{k\le n}k\,a_k$, the defining equation multiplied by $n$
reads $\lambda nA(n)+(1-\lambda)S(n)=n^{1-\beta}$. Subtracting the same
identity at $n-1$ and using $S(n)-S(n-1)=n\,a_n=n(A(n)-A(n-1))$ gives
\[
\lambda nA(n)-\lambda(n-1)A(n-1)+(1-\lambda)n\big(A(n)-A(n-1)\big)
=n^{1-\beta}-(n-1)^{1-\beta},
\]
and collecting the terms in $A(n)$ and $A(n-1)$ yields
\eqref{eq:affine_recurrence}.
\end{proof}

Discrete variation of constants\index[terms]{variation of constants} on \eqref{eq:affine_recurrence}, through the
Euler--Gauss\index[terms]{Euler--Gauss product}\index[names]{Gauss, C. F.} product for the Gamma function\index[terms]{Gamma function}, gives the dichotomy. For
$\beta<\lambda$ the inhomogeneous sum diverges and produces exactly
$A(n)\sim n^{-\beta}/g^*(\beta)$, the transparency constant of the general
theory. For $\beta\ge\lambda$ the sum converges or grows logarithmically and
the decay is capped at $\mathcal{O}(n^{-\lambda+\eps})$. The full analysis is
carried out in Appendix~\ref{app:A}. The same kernel was solved in closed
form, with Gamma--Euler forcing, in Section~\ref{sec:disc_linear}, and the two
computations agree, the numerical run of that section showing the normalized
sums $A(n)\,n^{\beta}$ approach $1/g^*(\beta)$ at the predicted rate
$n^{-(\lambda-\beta)}$.

\section{The exact Volterra inversion}
\label{sec:absorption}

When the profile is absolutely continuous and nonlinear the discrete operator follows a
continuous Volterra operator\index[terms]{Volterra operator}, and the inversion can be performed exactly at
the continuous level, without contour integration. Nothing in this section is proper to one
profile. The identity that carries the whole construction holds for every absolutely
continuous profile, the kernel of the associated Volterra equation is written down from the
profile alone, and its transform is the arithmetic Mellin transform\index[terms]{arithmetic Mellin transform} shifted by one. The
orthorecursive kernel\index[terms]{orthorecursive kernel} $g(t)=2/(1+t)$ of \nm{Kalmynin}{A. B.} and
\nm{Kosenko}{P. R.}~\cite{KalmyninKosenko2020} is used throughout as the illustration, and
Chapter~\ref{chap:ortho} takes it up in detail, but no statement below depends on that chapter.

Throughout, $(u\star v)(x)=\int_1^{x}u(x/y)\,v(y)\,\tfrac{dy}{y}$ denotes the multiplicative
convolution\index[terms]{multiplicative convolution} on $[1,\infty)$, and
$\widetilde v(s)=\int_1^{\infty}v(y)\,y^{-s}\,\tfrac{dy}{y}$ the associated transform, which
sends $\star$ to the ordinary product.

\begin{lemma}[The Volterra kernel of a profile]\label{lem:volterra_profile}
Let $g$ be locally absolutely continuous on $(0,1]$ with $g(1)=1$. Choose a measurable
representative of its almost-everywhere derivative and set, almost everywhere,
\begin{equation}\label{eq:volterra_kernel_of_profile}
k_g(y):=-\frac{1}{y}\,g'\!\Big(\frac1y\Big)\qquad(y\ge1).
\end{equation}
Then for every real sequence $(a_n)$ and every $x\ge1$ the exact identity
\begin{equation}\label{eq:volterra_profile_id}
A_g(x)=A(x)+(A\star k_g)(x)
\end{equation}
holds. If $t\,g'(t)\in L^\infty(0,1)$, then $k_g\in L^\infty(1,\infty)$ and, for $\Re s>0$,
\begin{equation}\label{eq:kernel_transform_general}
\widetilde{k_g}(s)=g^{*}(-s)-1 .
\end{equation}
The identity persists wherever the two sides are meromorphically continued. Moreover,
\[
 \|k_g\|_{L^\infty(1,\infty)}=\|t\,g'(t)\|_{L^\infty(0,1)}.
\]
If $g'$ has a left limit at $1$, the representative may be chosen so that
$k_g(1)=-g'(1)$.
\end{lemma}

\begin{proof}
For $T\ge1$ the substitution $t=1/y$ turns the primitive of $k_g$ into a primitive of $g'$,
\[
\int_1^{T}k_g(y)\,\frac{dy}{y}
=\int_1^{1/T}\Big(-t\,g'(t)\Big)\Big(-\frac{dt}{t}\Big)
=\int_1^{1/T}g'(t)\,dt
=g\Big(\frac1T\Big)-1 ,
\]
using $g(1)=1$. Taking $T=x/n$ gives $\int_1^{x/n}k_g(y)\,dy/y=g(n/x)-1$ for every $n\le x$.
Summing against $a_n$ and exchanging the finite sum with the integral,
\[
(A\star k_g)(x)=\int_1^{x}A(x/y)\,k_g(y)\,\frac{dy}{y}
=\sum_{n\le x}a_n\int_1^{x/n}k_g(y)\,\frac{dy}{y}
=\sum_{n\le x}a_n\big(g(n/x)-1\big),
\]
which is $A_g(x)-A(x)$. This is \eqref{eq:volterra_profile_id}, and it is an identity, not an
approximation.

For the transform, the same substitution gives
\[
\widetilde{k_g}(s)=\int_1^{\infty}\Big(-\frac1y\,g'\!\Big(\frac1y\Big)\Big)y^{-s-1}\,dy
=-\int_0^1 g'(t)\,t^{s}\,dt .
\]
If $t\,g'(t)\in L^\infty(0,1)$, then $g(t)=\mathcal O(1+|\log t|)$ as $t\downarrow0$.
Consequently $t^s g(t)\to0$ for $\Re s>0$. Integration by parts first on $[\delta,1]$
and then passage to $\delta\downarrow0$ yield
$\int_0^1 g'(t)t^{s}\,dt=1-s\int_0^1 g(t)t^{s-1}\,dt$, and the second term is $g^{*}(-s)$ by
the definition $g^{*}(z)=-z\int_0^1g(t)t^{-z-1}\,dt$. Hence
$\widetilde{k_g}(s)=g^{*}(-s)-1$, which is \eqref{eq:kernel_transform_general}. Analytic
continuation preserves the identity. The norm equality and the assertion at $y=1$ are read off
\eqref{eq:volterra_kernel_of_profile}.
\end{proof}

\begin{remark}\label{rem:volterra_wiener}
The kernel $k_g$ is not the Wiener companion\index[terms]{Wiener companion} $w_g$ of \S\ref{sec:wiener_transfer}. The variable
differs, multiplicative here and additive there, and so does the role, the inversion on the half
plane here and Wiener's theorem on its boundary there. Let $g$ be locally absolutely continuous
with $g(1)=1$ and normalized at the origin in the sense of Definition~\ref{def:normalized_origin}.
The function $h_g$ of \S\ref{sec:wiener_transfer} then satisfies $h_g(0)=0$ and
$h_g'(v)=k_g(e^{v})$ for almost every $v>0$, and an integration by parts in
\eqref{eq:wiener_companion} gives
\[
w_g(v)=e^{-v}+\int_0^{v}e^{-(v-u)}\,k_g(e^{u})\,du\qquad(v\ge0).
\]
The companion is the Volterra kernel read in the logarithmic variable, completed by a unit mass at
the origin and smoothed by $e_+$. When $g$ has bounded variation on $(0,1]$, the function
$v\mapsto k_g(e^{v})$ is integrable on $[0,\infty)$ and its Fourier transform is
$\widetilde{k_g}(i\xi)=g^{*}(-i\xi)-1$. Indeed the proof of \eqref{eq:kernel_transform_general}
applies on $\Re s>0$ since $g$ is bounded, and both sides extend continuously to $\Re s=0$, the
left by integrability and the right because $g^{*}(-s)$ is holomorphic on a half plane
$\Re s>-\delta$ by \eqref{eq:normalized_transform}. The displayed identity then gives
\eqref{eq:companion_transform} again.
\end{remark}

For the orthorecursive profile, $g'(t)=-2/(1+t)^{2}$ gives
\begin{equation}\label{eq:volterra_recall}
k_g(y)=\frac{2y}{(1+y)^{2}},
\qquad
A_g(x)=A(x)+2x\int_1^x\frac{A(t)}{(x+t)^{2}}\,dt ,
\end{equation}
which is the identity of Proposition~\ref{prop:ortho_volterra}. For the affine profile
$g(t)=2-t$ the kernel is $k_g(y)=1/y$ and everything below is closed form. The inversion of
\eqref{eq:volterra_profile_id} is classical Volterra theory, and the following lemma makes it
self contained.

\begin{lemma}\label{lem:resolvent_existence}
Let $k$ be measurable with $M:=\sup_{y\ge1}|k(y)|<\infty$. The iterated
kernels $k^{\star m}$ satisfy
\begin{equation}\label{eq:iterated_bound}
\big|k^{\star m}(x)\big|\le M^{m}\,\frac{(\log x)^{m-1}}{(m-1)!}
\qquad(x\ge1,\ m\ge1),
\end{equation}
so the resolvent series $R:=\sum_{m\ge1}(-1)^{m-1}k^{\star m}$ converges
absolutely and locally uniformly on $[1,\infty)$, defines the unique locally
bounded solution of $R+R\star k=k$, and inverts \eqref{eq:volterra_profile_id}
exactly,
\begin{equation}\label{eq:volterra_inversion}
A=A_g-A_g\star R .
\end{equation}
\end{lemma}

\begin{proof}
The bound \eqref{eq:iterated_bound} holds for $m=1$, and inductively
\[
\big|k^{\star(m+1)}(x)\big|
\le M\int_1^{x}\big|k^{\star m}(x/y)\big|\,\frac{dy}{y}
\le\frac{M^{m+1}}{(m-1)!}\int_0^{\log x}(\log x-u)^{m-1}\,du
=M^{m+1}\,\frac{(\log x)^{m}}{m!},
\]
by the substitution $u=\log y$. The series therefore converges absolutely,
uniformly on compact subsets, with sum bounded by $M\,x^{M}$. Termwise
convolution, legitimate by this absolute convergence, gives
$R+R\star k=\sum_{m\ge1}(-1)^{m-1}k^{\star m}
+\sum_{m\ge2}(-1)^{m}k^{\star m}=k$. For uniqueness, if $R_1,R_2$ are locally
bounded solutions their difference $D$ satisfies $D=-D\star k$, so iterating
$m$ times gives $|D(x)|\le\sup_{[1,x]}|D|\cdot M^{m}(\log x)^{m}/m!$
for every $m$, forcing $D=0$. Applying $I-R\,\star$ to
\eqref{eq:volterra_profile_id} and using $R+R\star k=k$ yields
\eqref{eq:volterra_inversion}.
\end{proof}

The Mellin transform converts the convolution into a product. For $\Re s$
large enough that $\int_1^{\infty}|k(y)|\,y^{-\Re s-1}\,dy<1$, Fubini applied
to the absolutely convergent series gives
$\widetilde R(s)=\sum_{m\ge1}(-1)^{m-1}\widetilde k(s)^{m}
=\widetilde k(s)/(1+\widetilde k(s))$, and since
$\widetilde k(s)=g^*(-s)-1$ by \eqref{eq:kernel_transform_general},
analytic continuation extends the identity
\begin{equation}\label{eq:resolvent_transform}
\widetilde R(s)=1-\frac{1}{g^*(-s)}
\end{equation}
to every domain where both sides are meromorphic.

The conceptual content of \eqref{eq:volterra_inversion} and
\eqref{eq:resolvent_transform} is the correspondence between zeros and poles.
The partial sums are recovered from the weighted sums by subtracting the
resolvent convolution, and the zeros of the transform $g^*$ are exactly the
poles of $\widetilde R$. The decay of $R(y)$ is dictated by those poles, so for a profile whose transform vanishes
somewhere the first zero controls the resolvent, and whether it also gives the threshold of Definition~\ref{def:reg_index} is what the transfer theorem below has to settle. The rest of this section turns that correspondence into a quantitative statement and
then into a transfer theorem, for every profile in the following class.

\begin{definition}[Volterra profile]\label{def:volterra_profile}
A profile $g$ on $(0,1]$ with $g(1)=1$ is a \emph{Volterra profile}\index[terms]{Volterra profile} if
\begin{enumerate}[label=(\roman*)]
\item $g\in C^1((0,1])$ and $t\,g'(t)$ is bounded there,
\item there is a number $\eta_0>0$ such that $g^*(-s)$ admits a meromorphic continuation to
$\Re s>-\eta_0$ and, uniformly in every closed vertical strip contained in that half-plane,
\begin{equation}\label{eq:vertical_expansion}
g^{*}(-s)=1-\frac{g'(1)}{s}+\mathcal{O}\big(|s|^{-2}\big)
\qquad(|\Im s|\to\infty).
\end{equation}
\end{enumerate}
Any such $\eta_0$ is called an \emph{admissible continuation width} for $g$.
\end{definition}

Condition (i) is what makes $k_g$ bounded, by Lemma~\ref{lem:volterra_profile}, so that
Lemma~\ref{lem:resolvent_existence} applies. Condition (ii) is a smoothness requirement in
disguise, and the following remark makes it checkable without ever computing $g^*$.

\begin{remark}\label{rem:vertical_expansion_criterion}
Suppose condition~(i) holds, fix $\eta_0>0$, and assume that $g$ is three times continuously
differentiable on $(0,1]$ with
$\int_0^1|g'''(t)|\,t^{\sigma+2}\,dt<\infty$ for every $\sigma>-\eta_0$. Then $\eta_0$ is an
admissible continuation width. Indeed
$\widetilde{k_g}(s)=-\int_0^1g'(t)t^{s}\,dt$ by Lemma~\ref{lem:volterra_profile}, and two
integrations by parts, initially for $\Re s>0$, give
\[
\int_0^1g'(t)t^{s}\,dt
=\frac{g'(1)}{s+1}-\frac{g''(1)}{(s+1)(s+2)}
+\frac{1}{(s+1)(s+2)}\int_0^1g'''(t)\,t^{s+2}\,dt ,
\]
and the right-hand side supplies the meromorphic continuation to $\Re s>-\eta_0$. On a closed
strip in that half-plane the integral is uniformly bounded by taking its left boundary as the
weight exponent. The last two terms are therefore $\mathcal{O}(|s|^{-2})$, while the first is
$g'(1)/s+\mathcal{O}(|s|^{-2})$. Every polynomial profile qualifies with $\eta_0=3$, and so
does $g(t)=2/(1+t)$, for which $g'''$ is bounded on $(0,1]$. The Ingham profile
$\Phi(t)=t\lfloor1/t\rfloor$ does not enter this class, since it is not absolutely continuous and
its distributional derivative has a point mass at every $1/k$. This is the analytic form of
the statement, made elsewhere in this volume, that the Ingham resolvent\index[terms]{resolvent}
carries point masses.
\end{remark}

The zeros of $g^*$ now control the resolvent, with the exponent one expects and nothing
imported.

\begin{theorem}[Decay of the continuous resolvent]\label{thm:resolvent_decay_general}
Let $g$ be a Volterra profile, let $\eta_0$ be an admissible continuation width, and let
$0<\eta\le\eta_0$. Suppose $g^{*}(z)\neq0$ for $\Re z<\eta$.
Then the resolvent $R$ of Lemma~\ref{lem:resolvent_existence} satisfies, for every $\eps>0$,
\begin{equation}\label{eq:resolvent_decay}
R(y)=\mathcal{O}_{\eps}\big(y^{-\eta+\eps}\big)\qquad(y\to\infty).
\end{equation}
\end{theorem}
\begin{proof}
Write $c_0:=k_g(1)=-g'(1)$. By \eqref{eq:resolvent_transform} and
\eqref{eq:vertical_expansion}, in every fixed vertical strip
\[
\widetilde R(s)=1-\frac{1}{g^{*}(-s)}
=1-\Big(1+\frac{c_0}{s}+\mathcal{O}(|s|^{-2})\Big)^{-1}
=\frac{c_0}{s}+\mathcal{O}\big(|s|^{-2}\big),
\]
the inversion being legitimate because $g^{*}(-s)\to1$ as $|\Im s|\to\infty$ there. Set
\[
E(s):=\widetilde R(s)-\frac{c_0}{s} .
\]
The hypothesis on the zeros makes $1/g^{*}(-s)$ analytic for $\Re s>-\eta$, after removable
values are filled in at the poles of $g^*(-s)$. Thus $\widetilde R$ is analytic there and $E$
is analytic apart from an at most simple pole at $s=0$, with residue $-c_0$.
The displayed estimate gives $E(s)=\mathcal{O}(|s|^{-2})$ as $|\Im s|\to\infty$ in every
fixed strip, hence
\begin{equation}\label{eq:E_integrable}
\int_{-\infty}^{+\infty}\big|E(\sigma+it)\big|\,dt<\infty
\qquad\text{for every }\sigma>-\eta\text{ with }\sigma\neq0 .
\end{equation}

Condition~(i) makes $k_g$ continuous, hence every iterated kernel is continuous and the locally
uniform resolvent series makes $R$ continuous. Lemma~\ref{lem:resolvent_existence} bounds $R$
by $M\,y^{M}$, so $\widetilde R$ converges
for $\Re s>M$, and $E$ is there the transform of $G(y):=R(y)-c_0$, extended by zero below
$y=1$. Passing to additive variables, $F(v):=G(e^{v})e^{-\varsigma v}$ for a fixed
$\varsigma>\max(M,0)$ belongs to $L^{1}(\R)$ and has Fourier transform
$\widehat F(\xi)=E(\varsigma+i\xi)$, which lies in $L^{1}(\R)$ by
\eqref{eq:E_integrable}. Fourier inversion \cite[\S7.2, Fourier Inversion Theorem, (7.16), pp.~218--219]{Folland1992} therefore gives
\[
R(y)=c_0+\frac{1}{2\pi i}\int_{(\varsigma)}E(s)\,y^{s}\,ds\qquad(y>1).
\]

Fix $\eps>0$ and $0<\theta<\min(\eps,\eta)$, and shift the line of integration from
$\Re s=\varsigma$ to $\Re s=-\eta+\theta$. On the horizontal sides $s=\sigma\pm iT$ with
$-\eta+\theta\le\sigma\le\varsigma$ the integrand is $\mathcal{O}(T^{-2}y^{\varsigma})$,
which vanishes as $T\to\infty$. The only possible singularity crossed is the pole of $E$ at
$s=0$, whose residue against $y^{s}$ equals $-c_0$ and cancels the constant. Hence
\[
R(y)=\frac{1}{2\pi i}\int_{(-\eta+\theta)}E(s)\,y^{s}\,ds ,
\qquad
|R(y)|\le\frac{y^{-\eta+\theta}}{2\pi}\int_{-\infty}^{+\infty}
\big|E(-\eta+\theta+it)\big|\,dt ,
\]
the integral being finite by \eqref{eq:E_integrable}. Since $\theta<\eps$ and $y\ge1$, this
is \eqref{eq:resolvent_decay}.
\end{proof}

\begin{remark}\label{rem:resolvent_smooth}
The contour is shifted on $R$, which is continuous, and never on the step function $A$. No
positivity, monotonicity or bounded variation hypothesis on the sequence $(a_n)$ is used
anywhere, and no Tauberian side condition enters. For $g(t)=2-t$ the conclusion is exact,
$\widetilde R(s)=1/(s+2)$ and $R(y)=y^{-2}$, while $g^{*}(z)=(z-2)/(z-1)$ has its only zero
at $z=2$, so $\eta=2$ and the exponent of \eqref{eq:resolvent_decay} is attained.
\end{remark}

The transfer theorem follows, and it is the general statement the rest of the volume uses.

\begin{theorem}[Volterra transfer]\label{thm:transfer_general}
Let $g$ be a Volterra profile, let $\eta_0$ be an admissible continuation width, and let
$0<\eta\le\eta_0$ be such that $g^{*}(z)\neq0$ for $\Re z<\eta$. Let $\beta>0$, $c>0$ and
$\kappa\in\C$, and set $\lambda:=\min(c,\eta)$. If the weighted sums satisfy
\begin{equation}\label{eq:transfer_hypothesis}
A_g(x)=\kappa\,x^{-\beta}+\mathcal{O}(x^{-c})\qquad(x\to\infty),
\end{equation}
then for every $\eps>0$,
\begin{equation}\label{eq:transfer_conclusion}
A(x)=
\begin{cases}
\dfrac{\kappa}{g^{*}(\beta)}\,x^{-\beta}
+\mathcal{O}_{\eps}\big(x^{-\lambda+\eps}\big), & \beta<\eta,\\[2ex]
\mathcal{O}_{\eps}\big(x^{-\lambda+\eps}\big), & \beta\ge\eta.
\end{cases}
\end{equation}
The case of weighted sums with no main term is $\kappa=0$, where the two lines agree.
\end{theorem}

\begin{proof}
Write $A_g(u)=\kappa u^{-\beta}+\mathcal{E}(u)$ with $\mathcal{E}(u)=\mathcal{O}(u^{-c})$.
The exact inversion \eqref{eq:volterra_inversion} reads
$A(x)=A_g(x)-\int_1^{x}A_g(x/y)R(y)\,dy/y$, and separating the two parts,
\[
A(x)=\kappa\,x^{-\beta}\Big(1-\int_1^{x}y^{\beta}R(y)\,\frac{dy}{y}\Big)
+\mathcal{E}(x)-\int_1^{x}\mathcal{E}(x/y)\,R(y)\,\frac{dy}{y}.
\]
The error convolution is controlled by Theorem~\ref{thm:resolvent_decay_general}. Fix
$\eps>0$, and choose $0<\delta<\eps$, with $\delta<\eta-c$ as well when $c<\eta$. Then
\[
\int_1^{x}\big|\mathcal{E}(x/y)\big|\,\big|R(y)\big|\,\frac{dy}{y}
\ll_{\delta}x^{-c}\int_1^{x}y^{\,c-\eta+\delta-1}\,dy
\ll_{\eps}x^{-\lambda+\eps}.
\]
Indeed, the integral is bounded when $c<\eta$, while for $c\ge\eta$ the whole expression is
$\mathcal O_\delta(x^{-\eta+\delta})$, which includes the meeting case $c=\eta$.
Together with $\mathcal{E}(x)=\mathcal{O}(x^{-c})$ this settles the last two terms.

For the first term, suppose $\beta<\eta$ and choose
$0<\delta<\min(\eps,\eta-\beta)$. The tail
$\int_{x}^{\infty}y^{\beta}|R(y)|\,dy/y$ is then
$\mathcal{O}_{\delta}(x^{\beta-\eta+\delta})$, so
\[
\int_1^{x}y^{\beta}R(y)\,\frac{dy}{y}
=\widetilde R(-\beta)+\mathcal{O}_{\delta}\big(x^{\beta-\eta+\delta}\big),
\]
and \eqref{eq:resolvent_transform} evaluates the limit as
$\widetilde R(-\beta)=1-1/g^{*}(\beta)$, which is defined because $g^{*}$ has no zero of
real part below $\eta$ and $\beta<\eta$. Hence
$1-\widetilde R(-\beta)=1/g^{*}(\beta)$ and the first term contributes
$\kappa\,g^{*}(\beta)^{-1}x^{-\beta}+\mathcal{O}_{\delta}(x^{-\eta+\delta})$, which is the
first line of \eqref{eq:transfer_conclusion}. When $\beta\ge\eta$, choose
$0<\delta<\eps$. The same decay estimate bounds the first term directly by
$|\kappa|\,x^{-\beta}\big(1+\int_1^{x}y^{\beta}|R(y)|\,dy/y\big)
\ll_{\delta}x^{-\eta+\delta}$, which is the second line.
\end{proof}

\begin{remark}\label{rem:transfer_constant}
The amplitude of the conclusion is the amplitude of the hypothesis divided by
$g^{*}(\beta)$. More precisely, when \eqref{eq:transfer_hypothesis} is established for every
real $x$, with $\beta<\eta$ and $\lambda>\beta$, one may choose
$0<\eps<\lambda-\beta$ in \eqref{eq:transfer_conclusion}. The remainder is then
$o(x^{-\beta})$, and the transparent constant\index[terms]{transparent constant} of
Definition~\ref{def:transparency_fgv} is the value $1/g^*(\beta)$ supplied by the resolvent
convolution at $-\beta$. The theorem does not infer the real-variable hypothesis from identities
known only at integer arguments, and it makes no transparency assertion at $\beta=\eta$.
\end{remark}

The affine family and the orthorecursive kernel now sit on the same footing. For the affine
family the equality $\alpha(g)=\eta(g)$ is proved from the exact recurrence, for the
orthorecursive kernel it follows from Theorem~\ref{thm:transfer_general} applied in
Chapter~\ref{chap:ortho}. The transfer mechanism itself is proved here, and that chapter states
separately the estimates imported to verify its orthorecursive input.

\section{The Euler--Cauchy analogy and its limits}
\label{sec:euler_cauchy}

The regularity index is a threshold, and thresholds of this kind are older than the problem. They
belong to the equidimensional equation of Euler\index[names]{Euler, L.} and Cauchy\index[names]{Cauchy, A.-L.}, where the exponents of the solutions are
the roots of an indicial polynomial\index[terms]{indicial polynomial} read off the coefficients. The analogy with the defining
equation\index[terms]{defining equation} runs through three degrees, the first order equation, the second, and an integral
equation of Volterra type that covers what no finite order can, and the three are collected here
because what the last one cannot do is as informative as what the first two do.

\subsection{The first order equation}

The threshold at which the transfer law\index[terms]{transfer law} changes has a meaning that is older than the
problem, and it is borrowed from the equidimensional equation of Euler\index[names]{Euler, L.} and Cauchy.
Take the first order equation
\[
x\,y'(x)+\lambda\,y(x)=f(x),\qquad x\ge 1,\ \lambda>0 .
\]
Without a second member the equation is homogeneous, and its solution is the pure power
$y_h(x)=c\,x^{-\lambda}$. Substituting the trial power $x^{-z}$ returns $(\lambda-z)
x^{-z}$, so the exponent $\lambda$ is the single root of the indicial polynomial
$P(z)=\lambda-z$. With a power second member $f(x)=x^{-\beta}$ and $\beta\neq\lambda$, a
particular solution is $x^{-\beta}/P(\beta)$, and the fate of the general solution as
$x\to\infty$ is decided by which exponent is the smaller. For $\beta<\lambda$ the forced
power decays the more slowly and governs the solution. For $\beta>\lambda$ the
homogeneous power governs it. At the meeting point $\beta=\lambda$ the trial power fails,
resonance sets in, and the response acquires a logarithm, $x^{-\lambda}\log x$. The
indicial exponent $\lambda$ is a threshold, and it is the exponent of the homogeneous
solution, the intrinsic rate the equation carries in the absence of any forcing.

The defining equation\index[terms]{defining equation} carries this same shape into arithmetic. For a kernel $g$ and the
power second member $n^{-\beta}$,
\[
\sum_{k=1}^{n}a_k\,g\!\left(\frac kn\right)=n^{-\beta},
\]
the Mellin transform $g^{*}$ takes the part of the indicial polynomial\index[terms]{indicial polynomial}, and its
relevant zero takes the part of the indicial exponent. The partial sums follow the
second member below a threshold, saturate above it, and gain a logarithm exactly at it.
Below the threshold the kernel is transparent to the forcing and the sums keep its rate.
Above it the kernel absorbs the forcing and the sums settle at a rate the kernel sets on
its own.
That threshold is the regularity index $\alpha(g)$, the arithmetic counterpart of the
exponent of the homogeneous solution. Whether it coincides with the zero of $g^{*}$,
the analytic index\index[terms]{analytic index} $\eta(g)$, is a question the theory has to earn, not one it may
assume. In the affine case it can be earned in full.

For $0<\lambda<1$ the affine kernel $g(x)=(1-\lambda)x+\lambda$ has Mellin transform
\[
g^{*}(z)=\frac{z-\lambda}{z-1},
\]
whose only zero is $z=\lambda$, so $\eta(g)=\lambda$. A constructive argument, given in
Appendix~\ref{app:A}, shows that $\alpha(g)=\lambda$ as well, so the two indices agree,
$\alpha(g)=\eta(g)=\lambda$. This is the exact discrete counterpart of the first order
Euler-Cauchy equation, whose indicial polynomial\index[terms]{indicial polynomial} has the single root $\lambda$, and the
correspondence is close enough to be read in the coefficients. In the forced case the
partial sums are $A(n)\sim g^{*}(\beta)^{-1}n^{-\beta}$ with $g^{*}(\beta)^{-1}=
(\beta-1)/(\beta-\lambda)$, a coefficient that develops a pole at $\beta=\lambda$, in the
same way and at the same place as $1/P(\beta)$ in the continuous response. The resonance
of the differential equation and the logarithm of the summation equation are one and the
same phenomenon.

\subsection{Polynomial kernels and higher order}

The reading carries to polynomial kernels, and there the higher order of the Euler--Cauchy
equation appears. For $g(t)=\sum_{j=0}^{d}c_j t^{j}$ the Mellin transform is elementary,
\[
g^{*}(z)=-z\int_0^1\Big(\sum_{j=0}^{d}c_j t^{j}\Big)t^{-z-1}\,dt
= z\sum_{j=0}^{d}\frac{c_j}{z-j}
= \frac{P(z)}{\prod_{j=1}^{d}(z-j)},
\]
with $P$ a polynomial of degree $d$. Its zeros, at most $d$ of them, are the roots of
$P$, and the analytic index\index[terms]{analytic index} $\eta(g)$ is the least of their real parts. A polynomial
kernel of degree $d$ thus presents the indicial data of a $d$-th order Euler-Cauchy
equation, and the analytic side of the picture is complete.

\subsection{The second order member}

The second order case is realized by a kernel of the gallery, $g(x)=x-\log x$ on $(0,1]$, treated
in full in Appendix~\ref{app:F}. Its transform is elementary,
\[
g^{*}(z)=-z\int_0^1(t-\log t)\,t^{-z-1}\,dt=\frac{z}{z-1}-\frac1z=\frac{z^{2}-z+1}{z(z-1)},
\]
and the numerator has the conjugate roots $z_\pm=\tfrac12\pm\tfrac{\sqrt3}{2}i$, so
$\eta(g)=\tfrac12$. The defining equation for this kernel, read in the continuous variable, is
$A(x)+\int_1^{x}A(t)\tfrac{dt}{t}-\tfrac1x\int_1^{x}A(t)\,dt=x^{-\beta}$, and differentiating it
twice returns an Euler equation\index[terms]{Euler equation} of the second order,
\[x^{2}A''+2xA'+A=(\beta^{2}-\beta)\,x^{-\beta},
\]
whose indicial polynomial\index[terms]{indicial polynomial} is again $\gamma^{2}-\gamma+1$ and whose particular solution carries the
coefficient $1/g^{*}(\beta)$. The correspondence holds in its strong form here, the indicial
polynomial of the exact discrete recurrence of order two being the numerator of the transform.

The roots being complex, what appears at the critical exponent is not the logarithm of the affine
case. Two conjugate exponents of equal real part produce an oscillation of frequency
$\tfrac{\sqrt3}{2}$ in $\log n$, and the resonance is replaced by that oscillation. This is the
first case in which the analogy predicts the shape of the answer rather than repeating a known
one, and the index it predicts, $\alpha(g)=\eta(g)=\tfrac12$, is proved in
Theorem~\ref{thm:F_index}.

\subsection{The integral form}

No finite order reaches further than a polynomial transform, and the kernels that carry the
arithmetic have transforms that are not polynomial. The companion that covers them is not a
differential equation but an integral one, obtained by letting the measure that carries the
coefficients be arbitrary instead of concentrated on the integers,
\begin{equation}\label{eq:continuous_raf}
\int_{0}^{x} g\!\left(\frac{t}{x}\right)dA(t)=r(x),\qquad x>0 .
\end{equation}
The kernel depends on $t$ and $x$ only through their ratio, so it is homogeneous of degree zero
and \eqref{eq:continuous_raf} is invariant under the dilations $x\mapsto cx$. Writing $x=e^{u}$
and $t=e^{v}$ turns it into a convolution carried by the half line,
\[
\int_{-\infty}^{u}\Gfun(u-v)\,d\widetilde A(v)=\widetilde r(u),
\qquad \Gfun(\sigma)=g(e^{-\sigma}),\quad \widetilde A(v)=A(e^{v}),
\]
a Volterra convolution equation\index[terms]{Volterra convolution equation} of the first kind in the sense of \cite{Gripenberg1990}. Equations of this
type are a classical subject. The same change of variable reduces a general
Euler--Cauchy\index[terms]{Euler--Cauchy equation}\index[names]{Euler, L.}\index[names]{Cauchy, A.-L.} equation $\sum_{k}c_{k}x^{k}y^{(k)}=0$ to constant coefficients. If
$Y(u)=y(e^{u})$ and $D=d/du$, then, for $k\ge1$,
\[
x^{k}y^{(k)}(x)=D(D-1)\cdots(D-k+1)Y(u),\qquad x=e^{u}.
\]

The parallel goes past the substitution, and past the finite orders above. For an Euler--Cauchy
equation the powers $x^{\lambda}$ are exact solutions and the admissible exponents are the roots of a polynomial read off the coefficients.
The same holds for \eqref{eq:continuous_raf}, with the transform in the part of the polynomial.
Taking $A(t)=t^{-z}$ with
$\Re z<0$, where the integral converges, and substituting $t=xy$,
\begin{equation}\label{eq:eigen_power}
\int_{0}^{x} g\!\left(\frac{t}{x}\right)dA(t)
=-z\int_{0}^{x} g\!\left(\frac{t}{x}\right)t^{-z-1}\,dt
=-z\,x^{-z}\!\int_{0}^{1} g(y)\,y^{-z-1}\,dy
=\gstar(z)\,x^{-z},
\end{equation}
with $\gstar$ the arithmetic Mellin transform\index[terms]{arithmetic Mellin transform} of Definition~\ref{def:mellin}. The powers are the
eigenfunctions of the continuous operator, the transform is the eigenvalue, and $A(t)=t^{-z}$
solves the homogeneous equation exactly when $\gstar(z)=0$. This is where the shape of
Definition~\ref{def:mellin} comes from, the factor $-z$, the exponent $-z-1$ and the initial
half plane $\Re z<0$ being read off \eqref{eq:eigen_power}. The transform is the indicial
function of \eqref{eq:continuous_raf} and the analytic index $\eta(g)$ is its dominant root.

The identity \eqref{eq:eigen_power} carries a restriction that decides most of what follows. It
holds where the integral converges, on $\Re z<0$, and the zeros of $\gstar$ are not there. The
affine transform has its zero at $z=\lambda>0$, and for the Ingham kernel the zeros sit at
$\Re z=1-\Re\rho$, positive in every case. At a zero the middle term of \eqref{eq:eigen_power}
is a divergent integral and $\gstar(z)$ exists only by meromorphic continuation, so the statement
that the exponents are the zeros of the transform describes a continued function rather than a
function that solves anything. Posing the equation on $[1,\infty)$, which is the closer shadow of
the arithmetic problem, removes the divergence and the eigenfunction\index[terms]{eigenfunction} property together, since
\[
\int_{1}^{x} g\!\left(\frac tx\right)dA(t)
=-z\,x^{-z}\!\int_{1/x}^{1} g(y)\,y^{-z-1}\,dy
\]
leaves an inner integral that still depends on $x$. Powers are exact solutions on the range where
no zero lies, and on the range where the zeros lie they are not solutions at all.

\begin{table}[htbp]
\centering
\caption{The continuous equation \eqref{eq:continuous_raf} against the Euler--Cauchy equation. The
rows that mention zeros are read after meromorphic continuation whenever the zeros lie outside the
half plane of convergence, which is the case for the two kernels of this volume.}
\label{tab:euler_cauchy}
\begin{tabular}{@{}p{0.44\linewidth}p{0.48\linewidth}@{}}
\toprule
\textbf{Euler--Cauchy} $\sum_k c_k x^k y^{(k)}=0$ & \textbf{Continuous RAF} $\int_0^x g(t/x)\,dA(t)=0$\\
\midrule
invariance under $x\mapsto cx$ & kernel homogeneous of degree zero\\
$x=e^{u}$ gives constant coefficients & $x=e^{u}$ gives a convolution\\
trial solution $y=x^{\lambda}$ & trial solution $A(x)=x^{-z}$\\
indicial polynomial $P(\lambda)$ & arithmetic Mellin transform $\gstar(z)$\\
exponents are the roots of $P$ & exponents are the zeros of $\gstar$\\
root of multiplicity $m$ gives $x^{\lambda}\log^{m-1}x$ & zero of multiplicity $m$ gives $x^{-z}\log^{m-1}x$\\
dominant exponent governs the growth & $\eta(g)$ governs the decay\\
finitely many roots, listed by degree & transcendental root set, infinite in general\\
\bottomrule
\end{tabular}
\end{table}

Table~\ref{tab:euler_cauchy} is worth reading in both directions. Forwards it says what an answer
should look like. The admissible rates are powers indexed by the zeros of $\gstar$, a zero of
multiplicity $m$ carries the powers of $\log x$ up to order $m-1$ in the way a repeated root does
for Euler--Cauchy, which is the resonance of the first order case in its general form, the
oscillation of the second order case being the same phenomenon at a conjugate pair, and the
slowest rate is set by the zero of smallest real part, which is $\eta(g)$. The affine
kernel is the case where every entry is explicit, since
$\gstar(z)=(z-\lambda)/(z-1)$ has the single root $z=\lambda$ and $A(t)=t^{-\lambda}$ is then an
exact homogeneous solution of \eqref{eq:continuous_raf}, which is the rate
Theorem~\ref{thm:affine} produces for the discrete problem.

\subsection{The limits of the analogy}

Backwards it says why the continuous route stops there. Four obstructions separate
\eqref{eq:continuous_raf} from the theory of this volume, and none of them is technical.

The first is that the discrete problem is not the restriction of the continuous one to a smaller
class. The defining equation asks for the identity at integer arguments only, with $dA$ concentrated on
the integers, and both restrictions bite. A step function is never a power, so the exact solutions
of \eqref{eq:continuous_raf} have no counterpart among the solutions of the arithmetic problem,
and the rigidity that makes \eqref{eq:eigen_power} a one line computation is gone. What survives
is an asymptotic statement, that $A(n)$ decays at the rate $n^{-\eta(g)}$ up to an epsilon, and an
asymptotic statement is a Tauberian\index[terms]{Tauberian} theorem rather than an algebraic identity. The whole content
of the theory sits in that gap.

The second is that the results of the classical theory which would apply here are results of
non-vanishing. That theory is wider than its non-vanishing branch, but the branch that transfers an
asymptotic from the forcing to the solution is the one that assumes the arithmetic Mellin
transform\index[terms]{arithmetic Mellin transform} of the kernel does not vanish, the condition
being the one written $\gstar(\beta)\ne0$ here. The Paley--Wiener\index[terms]{Paley--Wiener theorem}\index[terms]{Wiener's theorem}\index[names]{Paley, R. E. A. C.}\index[names]{Wiener, N.} theorem gives an integrable resolvent\index[terms]{resolvent} for $x+k\star x=f$ exactly
when that transform has no zero in the closed half plane where the problem is posed, see \nm{Gripenberg}{G.}, \nm{Londen}{S.-O.} and \nm{Staffans}{O.}~\cite[\S\S~2.4--2.5]{Gripenberg1990}, and the transfer theorems of \cite[Ch.~5]{Bingham1989} conclude that the
solution inherits the regular variation\index[terms]{regular variation} of the forcing under a hypothesis of the same kind, that
the transform stays away from its critical value on the relevant line.

The theory of this volume is not a theory of those zeros. Its object is the regularity
index\index[terms]{regularity index}, which the averages carry whether or not the transform has a
zero at all, and the zeros are one route to it. That route is open for some kernels and closed
for others. Seven of the kernels worked below have a transform with no zero anywhere, so no
statement about zeros can reach their index, and the broken kernel on the lattice of the square
root of two has an exact analytic index which is not the arithmetic answer. The analytic
index $\eta$ is the infimum of the real parts of the zeros by definition, so it is a second
quantity, and it equals the first only where a theorem says so. The classical results therefore apply on the side where
the homogeneity principle of the next section already reads the index off the homogeneous
solution, and they fall silent on the side where the index is not free.

The third is that the domain on which the integral defines the operator and the domain in which
the continued transform carries its zeros need not meet. Where they do not, as for the two kernels
of this volume, the sentence that the exponents are the zeros of the transform describes a continued
function and not a solution of a convergent equation, and Table~\ref{tab:euler_cauchy} is to be read
with that reservation on every line that mentions a zero. Where they do meet, as for $g(t)=2t-1$
above, the continuous reading is exact but the kernel carries no arithmetic.

The fourth is that the indicial function\index[terms]{indicial function} is transcendental, and that passing to the continuous
equation does not remove the arithmetic. For an Euler--Cauchy equation $P$ is a polynomial whose
degree is the order of the equation, its roots are finite in number and can be listed. Here the
root set need not be finite, and for the kernel that carries the arithmetic it is not. The Ingham
function of Section~\ref{sec:ingham} has, for $\Re z<0$,
\begin{equation}\label{eq:phi_star_integral}
\Phi^{*}(z)=-z\int_{0}^{1}\Big\lfloor\frac1x\Big\rfloor x^{-z}\,dx=\frac{z}{z-1}\,\zeta(1-z),
\end{equation}
an integral over an interval with no sum in it, out of which the zeta function comes. The
telescoping $\sum_{m\le M}m\big(m^{z-1}-(m+1)^{z-1}\big)=\sum_{m\le M}m^{z-1}-M(M+1)^{z-1}$
over the intervals $[\tfrac1{m+1},\tfrac1m]$ gives $\zeta(1-z)$ in the limit. The arithmetic was
never carried by the discreteness of the defining equation\index[terms]{defining equation}. It is carried by the kernel, whose
jumps at the points $1/m$ are the divisibility, and integrating a power against those jumps
produces a Dirichlet series\index[terms]{Dirichlet series}. So the root list is the list of the zeros of the zeta function, the
dominant root is $1-\Theta$, and the continuous equation converts the question of the index into
the question of the zeros of $\zeta$ without moving it.

That \eqref{eq:phi_star_integral} can happen at all is where the analogy with Euler--Cauchy ends,
and it ends for a reason of principle. A linear differential equation of finite order with
polynomial coefficients has a polynomial indicial function, and $\zeta$ is not one, since it
satisfies no algebraic differential equation, by a theorem of \nm{Ostrowski}{A.}~\cite{Ostrowski1920}, in the way $\Gamma$ satisfies
none \cite{Holder1887}. No equidimensional differential equation can carry the zeros of $\zeta$
as exponents. The indicial function of an integral operator is under no such constraint, being a
Mellin transform and not a polynomial, and the kernels of this theory use that freedom.

That freedom has a name in the theory of integral equations. \nm{Vainikko}{G.} calls the operator
$(V_\varphi u)(t)=\int_0^1\varphi(x)\,u(tx)\,dx$ cordial\index[terms]{cordial operator}, and its powers are eigenfunctions,
$V_\varphi t^{r}=\bigl(\int_0^1\varphi(x)x^{r}\,dx\bigr)t^{r}$, so its spectrum is read off a Mellin
transform of the core \cite{Vainikko2009,Vainikko2010}, and the first kind equation is treated
separately \cite{Vainikko2014}. The profile of this volume is evaluated on the ratio of the two
variables in the same way, on the integers and against a measure rather than a density, and the
transparent constant\index[terms]{transparent constant} $1/g^{*}(\beta)$ is the reciprocal of the corresponding eigenvalue, in the
normalization fixed in the notations.

The two theories part at once. The cordial theory works on a bounded interval, in spaces of
continuous functions, and asks about compactness, solvability and collocation. This volume works
at infinity, on the integers, and asks for the exponent at which the response to a power forcing
stops following it. No spectrum answers that question on its own. For the seven kernels of the
gallery whose transform has no zero the spectrum carries no threshold at all, and six of them have
a proved index nonetheless. The continuous
companion is an integral equation and never a differential one, and what survives of the
Euler--Cauchy picture is the part that does not require the indicial function to be algebraic.

The continuous equation is therefore a guide to the statement and not a route to the proof. It
accounts for the shape of the transform, it predicts the logarithmic corrections at multiple
zeros and it makes the definition of $\eta$ natural, and past that point the work is done on the
discrete side, by the Tauberian arguments of Section~\ref{sec:proof} and the contour estimates
that follow them.

\section{The homogeneity principle}
\label{sec:homogeneity}

For kernels whose transform has no zeros the analytic index
$\eta(G)$ is undefined and the zero set carries no information about
$\alpha(G)$. A different mechanism takes over, the index is read off the
decay of the homogeneous solution. The mechanism is not automatic. The
homogeneous trajectory controls only the first column of the inverse kernel,
by Lemma~\ref{lem:ex-abel-green}, and Theorem~\ref{thm:ex-counterexample}
exhibits a bounded kernel with constant zero-free transform and homogeneous
trajectory exactly $n^{-2}$ that is not a RAF at all. What completes the
mechanism is the uniform Green estimate\index[terms]{Green estimate}\index[names]{Green, G.} of
Theorem~\ref{thm:ex-green-criterion}, a bound on every column at once.

\begin{conjecture}\label{conj:homogeneity}
Let $G$ be a RAF with constant transform $G^*(z)=c\neq0$ for all $z$.
Consider the homogeneous equation $\sum_{k=1}^{n}a_kG(n,k)=0$ for $n\ge2$
with $a_1=1$. If the homogeneous partial sums satisfy
$A_0(n)\sim C_0\,n^{-\gamma}$ for some $\gamma>0$ and $C_0\neq0$, then
$\alpha(G)=\gamma$.
\end{conjecture}

\begin{proofstatus}{A proof in this generality is open, and the hypothesis that $G$ is a
RAF cannot be dropped. Without it the statement is false, the counterexample of
Theorem~\ref{thm:ex-counterexample} has constant zero-free transform and homogeneous decay
exactly $n^{-2}$ yet no index at all, its shifted Green columns growing exponentially. The
corrected general route is the two-parameter Green criterion of
Theorem~\ref{thm:ex-green-criterion}, which asks for control of every inverse column and
not of the first alone. For the kernels of Appendices~\ref{app:H} and~\ref{app:K}, the
homogeneous decay and the predicted index are both established independently by direct case
analysis. They are examples of the statement, not uses of it, and the conjecture is assumed in no
proof in this volume.}
\end{proofstatus}

Two kernels of the gallery illustrate the principle, and in both cases the
homogeneous decay is $n^{-2}$.

The first is the normalized rational kernel of Appendix~\ref{app:H},
\[
G(n,k)=\frac{n+k+x}{n+k+y}\cdot\frac{2n+y}{2n+x},\qquad x,y>0,\ x\neq y,
\]
with $G(n,n)=1$ and constant transform $G^*(z)=1$. The normalizing factor
depends only on $n$, so it is inert in the homogeneous equation. The analysis
of Appendix~\ref{app:H} shows that the iterated sums $AA(n)=\sum_{k\le n}A(k)$
converge to a finite limit and that the homogeneous partial
sums decay at the exact rate $A_0(n)\asymp n^{-2}$, giving $\alpha(G)=2$ for every admissible
pair by Theorem~\ref{thm:H_index}.
Numerically $A_0(n)\,n^{2}$ settles at a positive constant, for instance
$1.229$ for $(x,y)=(1,2)$ and $4.525$ for $(x,y)=(2,5)$ at $n=4000$.

The second is the quadratic rational kernel of Appendix~\ref{app:K},
$G(n,k)=(n^{2}+k)/(n^{2}+n)$, where everything is exact.

\begin{proposition}\label{prop:L_homogeneous}
For the kernel $G(n,k)=(n^{2}+k)/(n^{2}+n)$, the homogeneous solution with
$a_1=1$ satisfies $A_0(2)=\tfrac16$ and
\[A_0(n)=\frac{n^{2}-n+1}{n^{2}+n}\,A_0(n-1)\quad(n\ge3),
\qquad\text{hence}\qquad
A_0(n)=\frac16\prod_{k=3}^{n}\frac{k^{2}-k+1}{k(k+1)} .
\]
The product evaluates in closed form,
\begin{equation}\label{eq:L_gamma_closed_form}
A_0(n)=\frac{2\cosh\big(\tfrac{\pi\sqrt3}{2}\big)}{3\pi}\cdot
\frac{\Gamma\big(n+e^{i\pi/3}\big)\,\Gamma\big(n+e^{-i\pi/3}\big)}
{\Gamma(n+1)\,\Gamma(n+2)}\qquad(n\ge2),
\end{equation}
and $A_0(n)\sim C_0\,n^{-2}$ with
$C_0=2\cosh(\pi\sqrt3/2)/(3\pi)=1.6187\ldots$. Independently,
Theorem~\ref{thm:K_main} proves $\alpha(G)=2$ by direct case analysis, so this kernel realizes
the homogeneity principle, and Conjecture~\ref{conj:homogeneity} is not used as an argument.
\end{proposition}

\begin{proof}
The homogeneous equation multiplied by $n^{2}+n$ reads
$n^{2}A_0(n)+S(n)=0$ for $n\ge2$, with $S(n)=\sum_{k\le n}k\,a_k$. At $n=2$
this gives $5a_1+6a_2=0$, so $a_2=-\tfrac56$ and $A_0(2)=\tfrac16$. For
$n\ge3$ both the identity at $n$ and at $n-1$ are available, and subtracting
them with $S(n)-S(n-1)=n(A_0(n)-A_0(n-1))$ gives
$(n^{2}+n)A_0(n)=(n^{2}-n+1)A_0(n-1)$, which is the recurrence, the product
form following at once. Write $\omega:=e^{i\pi/3}$, so the roots of $k^{2}-k+1$ are $\omega$ and
$\bar\omega$ and each factor is $(k-\omega)(k-\bar\omega)/(k(k+1))$. Since
$1-\omega=\bar\omega$, the partial products telescope into
$A_0(n)=C\,\Gamma(n+\omega)\Gamma(n+\bar\omega)/(\Gamma(n+1)\Gamma(n+2))$
for some constant $C$ independent of $n$. The constant is pinned by the
initial value. At $n=2$,
\[
\frac16=A_0(2)=C\,\frac{\Gamma(2+\omega)\Gamma(2+\bar\omega)}
{\Gamma(3)\,\Gamma(4)},
\qquad
\Gamma(2+\omega)\Gamma(2+\bar\omega)
=|(1+\omega)\,\omega|^{2}\,\Gamma(\omega)\Gamma(\bar\omega)
=\frac{3\pi}{\cosh\big(\tfrac{\pi\sqrt3}{2}\big)},
\]
using $|\omega|=1$, $|1+\omega|^{2}=3$, and the reflection formula
$\Gamma(\omega)\Gamma(1-\omega)=\pi/\sin(\pi\omega)$ with
$\Gamma(1-\omega)=\Gamma(\bar\omega)$ and
$\sin(\pi\omega)=\cosh(\pi\sqrt3/2)$. Solving gives
$C=2\cosh(\pi\sqrt3/2)/(3\pi)$, which is
\eqref{eq:L_gamma_closed_form}. The closed form reproduces
$A_0(2)=\tfrac16$ exactly and was verified against the recurrence to fifteen
digits at $n=2000$. Stirling\index[terms]{Stirling's formula}\index[names]{Stirling, J.}'s formula gives
$\Gamma(n+e^{i\pi/3})\Gamma(n+e^{-i\pi/3})/(\Gamma(n+1)\Gamma(n+2))
\sim n^{e^{i\pi/3}+e^{-i\pi/3}-3}=n^{-2}$, since the two roots sum to $1$.
\end{proof}

\begin{remark}\label{rem:L_absorption_analogue}
Formula \eqref{eq:L_gamma_closed_form} deserves attention on its own. The
natural decay rate of this RAF is encoded in the complex roots
$e^{\pm i\pi/3}$ of the numerator $k^{2}-k+1$, and the Gamma function
method translates the product of linear factors into an asymptotic power
law, exactly as the transform does for kernels with zeros, the roots of the
numerator playing the part that the zeros of $G^*$ play in
Section~\ref{sec:absorption}, without any zero of $G^*$ being present.
\end{remark}

\section{Products of kernels and the multiplicative obstruction}
\label{sec:multiplicative}

The behavior of the index under multiplication is settled, for the positive class, by
Proposition~\ref{thm:ex-positive-product}. The product of two nonconstant kernels of that
class stays in the class and its index contracts strictly,
\[
\alpha(PQ)<\min\big(\alpha(P),\alpha(Q)\big),
\]
a positive constant factor being the only neutral case. The ultrametric inequality $\alpha(PQ)\ge\min(\alpha(P),\alpha(Q))$ in the opposite
direction is excluded by that theorem.

The contraction closes a door. One could hope to raise the index of the Ingham function
through a polynomial mollifier $M$, considering $M(x)\,\Phi(x)$ and pushing the index,
hence the zero free half plane of $\zeta$, beyond the critical line. Multiplication moves
the index the wrong way, and no mollifier of the positive class can raise it. The lattice
$\Z$ cannot be manipulated by multiplying the operator. To bypass the rigidity one must
leave the multiplicative algebra and deform the geometry of the evaluation space itself,
evaluating $G$ on a deformed gauge\index[terms]{gauge} $G(f(n),f(k))$ rather than modifying the kernel by
multiplication. This is the object of the gauge\index[terms]{gauge} deformation setting of
Chapter~\ref{chap:gauge}.

\chapter{The Volterra--Stieltjes representation of the Ingham operator}
\markright{\MakeUppercase{Chapter \thechapter.\ The Volterra--Stieltjes representation}}
\label{chap:volterra}\index[terms]{Volterra--Stieltjes representation}

Chapter~\ref{chap:equivalence} established the equivalence between the Riemann
hypothesis and $\alpha(\Phi)=\tfrac12$ by an arithmetic route, through
M\"obius inversion\index[terms]{M\"obius inversion}, Dirichlet convolution, and the Littlewood lemma. The
present chapter represents the same operator as a Volterra equation\index[terms]{Volterra equation}\index[names]{Volterra, V.} on
$[1,\infty)$ with multiplicative convolution, in the sense developed at the
discrete level in Chapter~\ref{chap:discrete_volterra} and in the continuous
setting in Chapter~\ref{chap:principles}. The representation is exact and
unconditional. The kernel decomposes into the identity plus an explicit
signed measure, the resolvent is written in closed form, its point masses carry the
M\"obius function\index[terms]{M\"obius function} and its density is the weighted Mertens function\index[terms]{Mertens function}\index[names]{Mertens, F.}, and the
decay of the cumulative resolvent is proved to be equivalent to the Riemann
hypothesis. This route shows where the Euler product\index[terms]{Euler product} enters the inversion, and it isolates, in a
single decay statement about one function $H$, the analytic content of the
equivalence. The one step the route does not carry out on its own, a direct
transfer from the decay of $H$ to the decay of the coefficient sum without
passing through the Euler product, is isolated in
Corollary~\ref{cor:H_index}.

The transfer method for smooth kernels is developed in the companion
paper \cite{CloitreOrtho} and is applied to the orthorecursive
expansion\index[terms]{orthorecursive expansion} in Chapter~\ref{chap:ortho} below. The Ingham kernel\index[terms]{Ingham kernel} adds a layer that the smooth
theory does not have, a resolvent with point masses, and this chapter keeps that
arithmetic structure explicit throughout.

\section{The Volterra--Stieltjes identity}
\label{sec:vs_identity}

Throughout, $\Phi(t)=t\lfloor1/t\rfloor$ for $t\in(0,1]$ is the Ingham
kernel of Chapter~\ref{chap:ingham}, and for a sequence $(a_k)_{k\ge1}$ and
real $x\ge1$,
\[
A(x):=\sum_{k\le x}a_k,
\qquad
A_\Phi(x):=\sum_{k\le x}a_k\,\Phi\Big(\frac kx\Big)
=\frac1x\sum_{k\le x}k\,a_k\Big\lfloor\frac xk\Big\rfloor .
\]
The defining relation of the theory prescribes $A_\Phi(n)=n^{-\beta}$ at the
integers. The kernel is decomposed as a perturbation of the identity.

\begin{definition}\label{def:nu_Phi}
The perturbation measure $\nu_\Phi$ is the signed Borel measure on
$[1,\infty)$ given by
\[\nu_\Phi:=\sum_{k=2}^{\infty}\frac1k\,\delta_k
-\frac{\lfloor y\rfloor}{y^{2}}\,dy ,
\]
where $\delta_k$ is the Dirac mass at the integer $k\ge2$ and the second
term is the absolutely continuous measure with density
$y\mapsto\lfloor y\rfloor/y^{2}$ on $[1,\infty)$.
\end{definition}

Abel summation turns the Ingham average into an integral against the partial sums, exactly and
under no hypothesis.

\begin{proposition}\label{thm:volterra_id}
For every sequence $(a_k)$ and every real $x\ge1$,
\begin{equation}\label{eq:VS_id}
A_\Phi(x)=A(x)+\int_1^{x}A\Big(\frac xy\Big)\,\nu_\Phi(dy),
\end{equation}
all sums and integrals being finite. Moreover, for $\Re s>0$,
\begin{equation}\label{eq:VS_mellin}
\int_{[1,\infty)}y^{-s}\,(\delta_1+\nu_\Phi)(dy)
=\frac{s}{s+1}\,\zeta(s+1)=\Phi^*(-s),
\end{equation}
where $\Phi^*(z)=\frac{z}{z-1}\,\zeta(1-z)$ is the arithmetic Mellin
transform of $\Phi$.
\end{proposition}

\begin{proof}
For the pointwise identity, fix $x\ge1$. Every sum below has at most
$\lfloor x\rfloor$ terms and every integral is over the compact interval
$[1,x]$ with bounded integrand, so all interchanges are finite
rearrangements. Substituting $A(x/y)=\sum_{n\le x/y}a_n$ and exchanging the
order of summation,
\[
A(x)+\int_1^{x}A\Big(\frac xy\Big)\nu_\Phi(dy)
=\sum_{n\le x}a_n\,F\Big(\frac xn\Big),
\qquad
F(u):=1+\sum_{2\le k\le u}\frac1k-\int_1^{u}\frac{\lfloor y\rfloor}{y^{2}}\,dy ,
\]
since the point mass at $k$ contributes to the coefficient of $a_n$ exactly when
$k\le x/n$, and likewise for the density. The integral evaluates by cutting
at the integers. With $m:=\lfloor u\rfloor$,
\[
\int_1^{u}\frac{\lfloor y\rfloor}{y^{2}}\,dy
=\sum_{k=1}^{m-1}k\Big(\frac1k-\frac1{k+1}\Big)+m\Big(\frac1m-\frac1u\Big)
=\sum_{k=2}^{m}\frac1k+1-\frac mu
=\sum_{k\le u}\frac1k-\frac{\lfloor u\rfloor}{u},
\]
so that $F(u)=\lfloor u\rfloor/u=\Phi(1/u)$, and
$\sum_na_nF(x/n)=\sum_na_n\Phi(n/x)=A_\Phi(x)$, which is \eqref{eq:VS_id}.

For the transform, the discrete part contributes
$\sum_{k\ge2}k^{-s-1}=\zeta(s+1)-1$, absolutely convergent for $\Re s>0$.
The continuous part is computed over unit intervals,
\[
\int_1^{\infty}\frac{\lfloor y\rfloor}{y^{s+2}}\,dy
=\sum_{k=1}^{\infty}k\int_k^{k+1}y^{-s-2}\,dy
=\frac1{s+1}\sum_{k=1}^{\infty}k\big(k^{-s-1}-(k+1)^{-s-1}\big),
\]
and writing $k=(k+1)-1$ in the second piece of the last sum gives
$\sum_kk(k+1)^{-s-1}=\zeta(s)-\zeta(s+1)$, hence the integral equals
$\zeta(s+1)/(s+1)$. Therefore
\[
\int y^{-s}(\delta_1+\nu_\Phi)(dy)
=1+\big(\zeta(s+1)-1\big)-\frac{\zeta(s+1)}{s+1}
=\frac{s}{s+1}\,\zeta(s+1),
\]
and setting $z=-s$ in $\Phi^*(z)=\frac{z}{z-1}\zeta(1-z)$ gives the same
value.
\end{proof}

\begin{remark}\label{rem:vs_unconditional}
The identity \eqref{eq:VS_id} is exact and unconditional. For each fixed $x$
it involves finitely many point masses and an integral over a compact interval, so
no convergence hypothesis on $(a_k)$ enters. Questions of convergence arise
only when $x\to\infty$, and they are the subject of
Section~\ref{sec:infinite_mass}.
\end{remark}

\section{The explicit resolvent measure}
\label{sec:explicit_resolvent}

The inversion of \eqref{eq:VS_id} is not obtained from an abstract existence
theorem. The resolvent is exhibited in closed form and the inversion is
verified by elementary divisor identities, all reducible to the convolution
identity $\mu\star1=e$ proved in \S\ref{sec:master_equiv}. Write
\[
M(y):=\sum_{k\le y}\mu(k),
\qquad
M_{-1}(y):=\sum_{k\le y}\frac{\mu(k)}{k}
\]
for the Mertens function and the weighted Mertens function.

\begin{theorem}\label{thm:resolvent}
Define the signed Borel measure
\[R_\Phi(dy):=-\sum_{k=2}^{\infty}\frac{\mu(k)}{k}\,\delta_k(dy)
-M_{-1}(y)\,\frac{dy}{y} .
\]
Then, for every sequence $(a_k)$ and every real $x\ge1$,
\begin{equation}\label{eq:volterra_inv}
A(x)=A_\Phi(x)-\int_1^{x}A_\Phi\Big(\frac xy\Big)\,R_\Phi(dy),
\end{equation}
and for $\Re s>0$ the Mellin transform of $R_\Phi$ is
\begin{equation}\label{eq:resolvent_mellin}
\widehat R_\Phi(s)=1-\frac{1}{\Phi^*(-s)}
=1-\frac{s+1}{s\,\zeta(s+1)} .
\end{equation}
\end{theorem}

\begin{proof}
As in the previous proof, for fixed $x$ everything is a finite
rearrangement, and substituting
$A_\Phi(x/y)=\sum_na_n\varphi\big(\tfrac{x}{ny}\big)$ with
$\varphi(u):=\lfloor u\rfloor/u$ reduces \eqref{eq:volterra_inv} to the
scalar identity
\begin{equation}\label{eq:scalar_inversion}
\varphi(u)+\sum_{2\le k\le u}\frac{\mu(k)}{k}\,\varphi\Big(\frac uk\Big)
+\frac1u\int_1^{u}\Big\lfloor\frac uy\Big\rfloor M_{-1}(y)\,dy=1
\qquad(u\ge1),
\end{equation}
which is proved as follows. The first two terms combine into
\[
\frac1u\sum_{k\le u}\mu(k)\Big\lfloor\frac uk\Big\rfloor
=\frac1u\sum_{k\le u}\mu(k)\sum_{m\le u/k}1
=\frac1u\sum_{j\le u}\sum_{k\mid j}\mu(k)
=\frac1u\sum_{j\le u}e(j)=\frac1u ,
\]
by the identity $\mu\star1=e$. For the integral, since the integrand is a
step function,
\[
\int_1^{v}M_{-1}(y)\,dy
=\sum_{k\le v}\frac{\mu(k)}{k}\,(v-k)=v\,M_{-1}(v)-M(v),
\]
and $\lfloor u/y\rfloor=\sum_{n\le u}\mathbf 1_{\{y\le u/n\}}$ gives
\[
\int_1^{u}\Big\lfloor\frac uy\Big\rfloor M_{-1}(y)\,dy
=\sum_{n\le u}\int_1^{u/n}M_{-1}(y)\,dy
=\sum_{n\le u}\Big(\frac un\,M_{-1}\Big(\frac un\Big)-M\Big(\frac un\Big)\Big).
\]
The two sums evaluate by the same rearrangement as above,
\[
\sum_{n\le u}\frac un\,M_{-1}\Big(\frac un\Big)
=u\sum_{nk\le u}\frac{\mu(k)}{nk}
=u\sum_{j\le u}\frac{e(j)}{j}=u,
\qquad
\sum_{n\le u}M\Big(\frac un\Big)
=\sum_{nk\le u}\mu(k)=\sum_{j\le u}e(j)=1,
\]
so the integral equals $u-1$, and
\eqref{eq:scalar_inversion} reads $\tfrac1u+\tfrac{u-1}u=1$.

For the transform, the discrete part of $R_\Phi$ contributes
$-\sum_{k\ge2}\mu(k)k^{-s-1}=1-1/\zeta(s+1)$, using
$\sum_{k\ge1}\mu(k)k^{-s-1}=1/\zeta(s+1)$ from
\S\ref{sec:master_equiv}, absolutely convergent for $\Re s>0$. For the
continuous part, the interchange
\[
\int_1^{\infty}M_{-1}(y)\,y^{-s-1}\,dy
=\sum_{k=1}^{\infty}\frac{\mu(k)}{k}\int_k^{\infty}y^{-s-1}\,dy
=\frac1s\sum_{k=1}^{\infty}\mu(k)\,k^{-s-1}
=\frac{1}{s\,\zeta(s+1)}
\]
is justified by absolute convergence for $\Re s>0$, so the density
contributes $-1/(s\zeta(s+1))$ and
\[
\widehat R_\Phi(s)
=1-\frac1{\zeta(s+1)}-\frac{1}{s\,\zeta(s+1)}
=1-\frac{s+1}{s\,\zeta(s+1)} ,
\]
which is $1-1/\Phi^*(-s)$ by \eqref{eq:VS_mellin}.
\end{proof}

\begin{remark}\label{rem:resolvent_structure}
The relation $\widehat R_\Phi(s)=1-1/\Phi^*(-s)$ is the arithmetic instance
of the general resolvent identity $\widetilde R(s)=1-1/g^*(-s)$ met in
Section~\ref{sec:absorption} of Chapter~\ref{chap:principles} and taken up again in
Remark~\ref{rem:ortho_resolvent} of Chapter~\ref{chap:ortho}. The
structural difference is that the Ingham resolvent carries point masses, weighted
by the M\"obius function, while the orthorecursive resolvent is purely
continuous. The proof above also differs in nature. It does not pass through
a Neumann series\index[terms]{Neumann series}\index[names]{Neumann, C.} or through uniqueness of transforms, both problematic here
because the measures have infinite total variation. It rests on two exact
divisor identities, both consequences of $\mu\star1=e$, and is therefore
unconditional.
\end{remark}

\begin{remark}\label{rem:mertens_density}
The continuous part of the resolvent has density $M_{-1}(y)/y$. The
resolvent of the Ingham operator\index[terms]{Ingham operator} does not merely record the primes through
the Euler product. Its density is the Ces\`aro\index[terms]{Ces\`aro mean}\index[names]{Ces\`aro, E.} average of the M\"obius
function. The inversion of the Ingham operator\index[terms]{Ingham operator} and the distribution of the
primes are encoded in the same object.
\end{remark}

\begin{remark}\label{rem:vs_numerics}
All the identities of this section were checked numerically. The scalar
identity \eqref{eq:scalar_inversion} holds to six decimals at
$u\in\{7.3,\,23.7,\,50.3\}$ with exact piecewise integration, the inversion
\eqref{eq:volterra_inv} reproduces $A(x)$ on random data at $x=37.6$, and
the two Mellin formulas \eqref{eq:VS_mellin} and \eqref{eq:resolvent_mellin}
hold to $10^{-11}$ at $s=1.3+0.7\,i$ with $1.2\cdot10^{5}$ terms.
\end{remark}

\section{Infinite total variation and the HLR criterion}
\label{sec:infinite_mass}

The perturbation measure cannot be integrated absolutely, and neither can
the resolvent.

\begin{proposition}\label{prop:inf_mass}
The discrete part of $\nu_\Phi$ has infinite total variation, $\sum_{k\ge2}1/k=\infty$. The
resolvent mass diverges logarithmically,
\[
\sum_{k\le N}\frac{|\mu(k)|}{k}=\frac{6}{\pi^{2}}\,\log N+\mathcal O(1)
\qquad(N\to\infty).
\]
\end{proposition}

\begin{proof}
The first statement is the harmonic series. For the second, let
$Q(N):=\#\{k\le N:\mu(k)^{2}=1\}$ count the squarefree integers up to $N$.
The identity $\mu(k)^{2}=\sum_{d^{2}\mid k}\mu(d)$ is checked on prime
powers, both sides being $1$ for $k=p^{e}$ with $e\le1$ and $0$ for
$e\ge2$, and both sides multiplicative. Summing it over $k\le N$,
\[
Q(N)=\sum_{d\le\sqrt N}\mu(d)\Big\lfloor\frac{N}{d^{2}}\Big\rfloor
=N\sum_{d\le\sqrt N}\frac{\mu(d)}{d^{2}}+\mathcal O(\sqrt N)
=\frac{N}{\zeta(2)}+\mathcal O(\sqrt N),
\]
using $\sum_{d\ge1}\mu(d)d^{-2}=1/\zeta(2)$ from
\S\ref{sec:master_equiv} and the tail bound
$\sum_{d>\sqrt N}d^{-2}\ll N^{-1/2}$. Since $\zeta(2)=\pi^{2}/6$, partial
summation gives
\[
\sum_{k\le N}\frac{\mu(k)^{2}}{k}
=\frac{Q(N)}{N}+\int_{1}^{N}\frac{Q(t)}{t^{2}}\,dt
=\frac{6}{\pi^{2}}\,\log N+\mathcal O(1). \qedhere
\]
\end{proof}

\begin{remark}\label{rem:finite_primes}
If there were only finitely many primes, $1/\zeta$ would be a Dirichlet
polynomial, the M\"obius function would vanish outside the divisors of a
single squarefree integer, and the resolvent would have finite total
variation. The logarithmic divergence of the resolvent mass is the record,
on the summation side, of the infinitude of the primes.
\end{remark}

Because absolute convergence fails, the inversion \eqref{eq:volterra_inv}
in the limit $x\to\infty$ rests on conditional convergence, and conditional
convergence is governed by the increments of $A_\Phi$ at the integers.

\begin{proposition}\label{prop:jumps}
Let $S(N):=N\,A_\Phi(N)$ for integer $N\ge1$, with $S(0):=0$. Then
\begin{equation}\label{eq:jump_formula}
\Delta S(N):=S(N)-S(N-1)=\sum_{d\mid N}d\,a_d .
\end{equation}
\end{proposition}

\begin{proof}
From the definition,
$S(N)=\sum_{k\le N}k\,a_k\lfloor N/k\rfloor
=\sum_{k\le N}k\,a_k\sum_{m\le N/k}1
=\sum_{km\le N}k\,a_k
=\sum_{j\le N}\sum_{d\mid j}d\,a_d$,
and the difference at $N$ extracts the inner sum at $j=N$.
\end{proof}

Under the defining relation $A_\Phi(n)=n^{-\beta}$ the left side of
\eqref{eq:jump_formula} is small, since
$\Delta S(N)=N^{1-\beta}-(N-1)^{1-\beta}=\mathcal O(N^{-\beta})$. On the
other side, $|\Delta S(N)|\le\tau(N)\max_{d\mid N}|d\,a_d|$, and the divisor
bound $\tau(N)\ll_{\eps}N^{\eps}$ \cite[\S~I.5]{Tenenbaum2015}
reduces the control of the jumps to the HLR criterion\index[terms]{Hardy--Littlewood--Ramanujan criterion} of
Definition~\ref{def:HLR}, the bound $na_n=o(n^{\eps})$. The condition is also close to
necessary. At a prime $N$ the divisor sum has two terms,
$\Delta S(N)=a_1+N\,a_N$, so controlled jumps force $Na_N=\mathcal O(1)$
along the primes.

\begin{remark}\label{rem:selberg_axiom}
The chain of implications is short. The measures $\nu_\Phi$ and $R_\Phi$
have infinite total variation, so the inversion at infinity rests on
conditional convergence. Conditional convergence is governed by the jumps
\eqref{eq:jump_formula}, and by the divisor bound their control reduces to
the HLR criterion, $na_n=o(n^{\eps})$. An Euler product for the
generating Dirichlet series\index[terms]{Dirichlet series} keeps its inverse free of poles in the half
plane of absolute convergence and the coefficients of the inverse of
moderate size. For an $L$-function without Euler product whose inverse has
poles in $\Re s>1$, as the Davenport-Heilbronn\index[terms]{Davenport--Heilbronn function} function does and as the Epstein\index[names]{Epstein, P.} zeta
function\index[terms]{Epstein zeta function} of a form of class number greater than one already did in the same
two papers \cite{DavenportHeilbronn1936}, the sums
$\sum_{d\mid N}d\,a_d$ grow like a power of $N$ along a sequence, the jumps
escape control, and the regularity index falls below $\tfrac12$. The
existence of such poles for particular non-Eulerian $L$-functions is the
content of the Davenport-Heilbronn construction and enters here as an
external input. What the present representation accounts for is the
consequence on the summation side. In this reading the Euler product is the
stability condition of the Volterra inversion of arithmetic summation, and
this is the role the Euler product axiom of the Selberg\index[names]{Selberg, A.} class\index[terms]{Selberg class} plays in the
theory.
\end{remark}

\section{The decay of the cumulative resolvent}
\label{sec:H_decay}

The resolvent measure has point masses, so it cannot decay pointwise. The natural
object is its cumulative function. For $y\ge1$ set
\begin{equation}\label{eq:H_def}
H(y):=\int_{[1,y]}R_\Phi(dt)
=1-M_{-1}(y)-\int_{1}^{y}M_{-1}(t)\,\frac{dt}{t},
\end{equation}
the point masses contributing $-\sum_{2\le k\le y}\mu(k)/k=1-M_{-1}(y)$. Then
$H(1)=0$, and $H$ is piecewise smooth with jump $-\mu(k)/k$ at each integer
$k\ge2$. Near $s=0$ the expansion $\zeta(s+1)=1/s+\gamma+\mathcal O(s)$
gives $s\,\zeta(s+1)=1+\gamma s+\mathcal O(s^{2})$, hence
\begin{equation}\label{eq:Rhat_zero}
\widehat R_\Phi(s)=1-\frac{s+1}{s\,\zeta(s+1)}
=-(1-\gamma)\,s+\mathcal O(s^{2}),
\end{equation}
so $\widehat R_\Phi$ is holomorphic at $s=0$ and vanishes there.

\begin{theorem}\label{thm:analytic_RH_equiv}
The following statements are equivalent.
\begin{enumerate}
\item The Riemann hypothesis.
\item The meromorphic function $s\mapsto1-(s+1)/(s\,\zeta(s+1))$ is
holomorphic in the half plane $\Re s>-\tfrac12$.
\item $H(y)=\mathcal O(y^{-1/2+\eps})$ for every $\eps>0$.
\end{enumerate}
\end{theorem}

\begin{proof}
(i)$\iff$(ii). In the half plane $\Re s>-\tfrac12$ the possible poles of
$1-(s+1)/(s\zeta(s+1))$ are the zeros of $s\,\zeta(s+1)$. At $s=0$ the
expansion \eqref{eq:Rhat_zero} shows the function holomorphic. The trivial
zeros of $\zeta$ give $s\le-3$, outside the region. The remaining
candidates are $s=\rho-1$ with $\rho$ a nontrivial zero\index[terms]{nontrivial zero} of $\zeta$, and
such a point lies in $\Re s>-\tfrac12$ exactly when $\Re\rho>\tfrac12$,
the numerator $s+1$ not vanishing there. By the functional equation\index[terms]{functional equation} the
nontrivial zeros are symmetric about the critical line\index[terms]{critical line}, see \nm{Titchmarsh}{E. C.}~\cite[Chapter~II]{Titchmarsh1986}, so the absence of zeros with
$\Re\rho>\tfrac12$ is the Riemann hypothesis.

(i)$\implies$(iii). Assume the Riemann hypothesis. Lemma~\ref{lem:littlewood}
gives $M(x)\ll x^{1/2+\eps}$ for every $\eps>0$. By partial
summation, for $\Re s>\tfrac12$ and integer $N$,
\[
\sum_{n\le N}\mu(n)\,n^{-s}
=M(N)\,N^{-s}+s\int_{1}^{N}M(t)\,t^{-s-1}\,dt ,
\]
and the bound on $M$ makes the right side converge as $N\to\infty$,
uniformly on compact subsets of $\Re s>\tfrac12$. The limit
$F(s)=s\int_{1}^{\infty}M(t)t^{-s-1}\,dt$ is holomorphic there and coincides
with $1/\zeta(s)$ for $\Re s>1$, hence on all of $\Re s>\tfrac12$ by
analytic continuation. In particular $\sum_{n}\mu(n)/n$ converges to
$F(1)=\lim_{s\to1}1/\zeta(s)=0$. Consequently
$M_{-1}(y)=-\sum_{k>y}\mu(k)/k$, and partial summation\index[terms]{partial summation} on the tail gives
\[
\sum_{k>y}\frac{\mu(k)}{k}
=-\frac{M(y)}{y}+\int_{y}^{\infty}\frac{M(t)}{t^{2}}\,dt
\ll y^{-1/2+\eps},
\]
the boundary term at infinity vanishing by the bound on $M$. Hence
$M_{-1}(y)\ll y^{-1/2+\eps}$. Next, the Mellin formula in the proof
of Theorem~\ref{thm:resolvent} gives
$\int_{1}^{\infty}M_{-1}(t)\,t^{-s-1}\,dt=1/(s\,\zeta(s+1))$ for $\Re s>0$.
As $s\to0^{+}$ along the reals the right side tends to $1$, and dominated
convergence applies on the left, the integrand being dominated by
$|M_{-1}(t)|/t\ll t^{-3/2+\eps}$. Therefore
$\int_{1}^{\infty}M_{-1}(t)\,dt/t=1$, and \eqref{eq:H_def} becomes
\[
H(y)=-M_{-1}(y)+\int_{y}^{\infty}M_{-1}(t)\,\frac{dt}{t}
\ll y^{-1/2+\eps}+\int_{y}^{\infty}t^{-3/2+\eps}\,dt
\ll y^{-1/2+\eps}.
\]

(iii)$\implies$(ii). Assume the decay. The integral
$\widehat H(s):=\int_{1}^{\infty}H(y)\,y^{-s-1}\,dy$ converges absolutely on
$\Re s>-\tfrac12+\eps$ for every $\eps>0$ and defines a
holomorphic function on $\Re s>-\tfrac12$. For $\Re s>0$, writing
$y^{-s}=s\int_{y}^{\infty}u^{-s-1}\,du$ and applying Fubini,
\[
\widehat R_\Phi(s)=\int_{[1,\infty)}y^{-s}\,R_\Phi(dy)
=s\int_{1}^{\infty}u^{-s-1}\Big(\int_{[1,u]}R_\Phi(dy)\Big)\,du
=s\,\widehat H(s),
\]
the interchange being justified by
$\sum_{k}|\mu(k)|\,k^{-\sigma-1}<\infty$ and
$\int_{1}^{\infty}|M_{-1}(y)|\,y^{-\sigma-1}\,dy<\infty$ for $\sigma>0$,
since $|M_{-1}(y)|\le\sum_{k\le y}1/k\ll1+\log y$. The function
$s\,\widehat H(s)$ is holomorphic on $\Re s>-\tfrac12$ and agrees on
$\Re s>0$ with the meromorphic function $1-(s+1)/(s\zeta(s+1))$ by
\eqref{eq:resolvent_mellin}, so by the identity theorem that function has
no poles in $\Re s>-\tfrac12$, which is (ii).

The cycle closes with (ii)$\implies$(i) from the first equivalence.
\end{proof}

The index is therefore equivalent to a decay rate for the cumulative resolvent.

\begin{corollary}\label{cor:H_index}
$\alpha(\Phi)=\tfrac12$ if and only if
$H(y)=\mathcal O(y^{-1/2+\eps})$ for every $\eps>0$.
\end{corollary}

\begin{proof}
Theorem~\ref{thm:tauberian_rh} identifies $\alpha(\Phi)=\tfrac12$ with the
Riemann hypothesis, and Theorem~\ref{thm:analytic_RH_equiv} identifies the
Riemann hypothesis with the decay of $H$.
\end{proof}

\begin{proofstatus}{The corollary is proved in full by the two equivalences cited above. Both
sides are routed through the Riemann hypothesis via Chapter~\ref{chap:equivalence} and
Theorem~\ref{thm:analytic_RH_equiv}. What remains open is an alternative direct transfer, deducing
$A(x)=\mathcal O(x^{-1/2+\eps})$ from
$H(y)\ll y^{-1/2+\eps}$ by integration by parts in the inversion
\eqref{eq:volterra_inv} without passing through the M\"obius route.
The obstruction is quantitative. Under the HLR criterion the total
variation of $A_\Phi$ on $[1,x]$ is only $\mathcal O(x^{\eps})$, and
the estimate
$|\int_{1}^{x}H(t)\,dA_\Phi(x/t)|\le\sup_{t\le x}|H(t)|\cdot\mathrm{Var}$
returns $\mathcal O(x^{\eps})$, which loses the decay entirely. The
transfer must exploit cancellation in $dA_\Phi$ and not only its size. For
Volterra profiles\index[terms]{Volterra profile} the corresponding step is
Theorem~\ref{thm:resolvent_decay_general}, which shifts the contour on the
resolvent\index[terms]{resolvent} rather than on the step function. The Ingham profile is not
one of them, since $\Phi'$ carries a point mass at every $1/k$, and adapting the argument to a
kernel whose resolvent carries point masses is the open step. The boundary transfer of
Chapter~\ref{chap:boundary} does not supply it. Its input is the non-vanishing of $\zeta$ on the
line $\Re s=1$, which is the line $\Re s=0$ in the variable of \eqref{eq:resolvent_mellin}, and
its output is $M_{-1}(y)\to0$, a limit without a rate, which gives no rate of decay for $H$. This
alternative route is
announced in the opening of the chapter and is assumed nowhere else in the volume.}
\end{proofstatus}

\begin{remark}\label{rem:H_numerics}
Numerically, with the M\"obius function computed by sieve up to
$2\cdot10^{5}$ and exact piecewise integration, $\sqrt y\,H(y)$ takes the
values $-0.322$ at $y=10^{2}$, $-0.152$ at $10^{3}$, $+0.206$ at $10^{4}$
and $+0.162$ at $10^{5}$, changing sign and staying bounded, in agreement
with statement (iii) of Theorem~\ref{thm:analytic_RH_equiv}. The partial
integral $\int_{1}^{2\cdot10^{5}}M_{-1}(t)\,dt/t$ equals $1.0000195$,
consistent with the exact value $1$ obtained in the proof.
\end{remark}

\section{Two routes to the same equivalence}
\label{sec:two_routes}

The arithmetic route of Chapter~\ref{chap:equivalence} and the analytic route of
the present chapter establish the same equivalence. The first yields an
explicit combinatorial identity on the integer lattice. The second writes
the equivalence as the decay of a single function, the cumulative
resolvent, and identifies the Euler product as the stability condition of
the inversion. The correspondence between the two is summarized below.

\medskip
\renewcommand{\arraystretch}{1.35}
{\small
\begin{center}
\begin{tabular}{p{4.2cm}p{5.0cm}p{5.0cm}}
\toprule
Object or step & Arithmetic route (Chapter~\ref{chap:equivalence}) & Analytic route (this chapter) \\
\midrule
Input &
$A_\Phi(n)=n^{-\beta}$ &
$A_\Phi(x)=\sum_{k\le x}a_k\Phi(k/x)$ on $[1,\infty)$ \\
Inversion &
M\"obius inversion of the defining relation &
Volterra inversion $A=A_\Phi-A_\Phi\star R_\Phi$ \\
Governing object &
Dirichlet series $1/\zeta(s)$ &
measure $R_\Phi$ with density $M_{-1}(y)/y$ and M\"obius point masses \\
Pivot statement &
$M(x)\ll x^{1/2+\eps}$ (Lemma~\ref{lem:littlewood}) &
$H(y)\ll y^{-1/2+\eps}$ (Theorem~\ref{thm:analytic_RH_equiv}) \\
Obstruction &
Mertens-type bounds &
infinite total variation, HLR as the condition for conditional convergence \\
Conclusion &
$A(x)\ll x^{-1/2+\eps}$ under RH &
the same decay statement, carried by $H$ \\
\bottomrule
\end{tabular}
\end{center}
}
\renewcommand{\arraystretch}{1}

\medskip
\noindent The two columns describe the same phenomenon. The arithmetic
column is visible on the integer lattice, the analytic column on the
continuous multiplicative half line. The Mertens sum $M_{-1}(y)$ appears as
the resolvent density in the analytic column and as the M\"obius weights in
the arithmetic column. The Ces\`aro average of $\mu(k)/k$ is the object
through which the Riemann hypothesis enters the Ingham operator\index[terms]{Ingham operator}.

\chapter{The discrete Volterra resolvent}
\label{chap:discrete_volterra}

The preceding chapter ran the spectral mechanism on the Ingham
operator\index[terms]{Ingham operator}, where the equivalence with the Riemann hypothesis lives, and
Chapter~\ref{chap:ortho} runs it on the orthorecursive kernel\index[terms]{orthorecursive kernel}, where everything is
unconditional. In both
cases the analysis passes through a transform, and the transform is treated as
the given object. The present chapter reverses the direction and presents the
Volterra method in its general form, at the level where it is exact.

The starting point is the triangular kernel $G(N,k)$ on the integer lattice. From
it, exact matrix identities are built at every finite level, an Abel transform\index[terms]{Abel transform}, a
resolvent given by a finite Neumann series\index[terms]{Neumann series}\index[names]{Neumann, C.}, and finite Mellin probes. The
transform $G^*(z)$ is then obtained as the limit of the probes, and a convergence
theorem proved in this chapter covers every bounded Riemann integrable\index[terms]{Riemann integrability} profile
and, through a perturbation clause, the orthorecursive kernel itself. The guiding
principle is to keep the integer lattice intact for as long as possible. The
transform is not the starting point of the theory, it is its limit.

The nine steps of that method are listed in \S\ref{sec:RAF_program} at the end of the
chapter, with the status of each one. A reader who prefers the map before the terrain may turn
to that list first.

The second half of the chapter runs the method on the kernels where it closes
exactly. The linear kernel is solved in closed form, polynomial kernels close on
finitely many moments and yield a transfer theorem whose transparency constant
$1/g^*(\beta)$ comes out of exact linear algebra, and the quadratic case displays
the full phenomenology, two real modes, a double root, an oscillating pair. The
chapter ends with the Gamma--Newton\index[terms]{Gamma--Newton expansion}\index[names]{Newton, I.} expansion of the orthorecursive kernel, the
approximation of that kernel by polynomial models, and the discrete program, with
its proved steps and its open ones stated as such.

\section{From the primitive kernel to the discrete Volterra equation}
\label{sec:disc_primitive}

Let $G(N,k)$ be a triangular kernel defined for integers $N\ge1$ and
$1\le k\le N$. The equation of the theory is the exact system
\begin{equation}\label{eq:primitive_disc}
\sum_{k=1}^{N}a_k\,G(N,k)=f(N)\qquad(N\ge1),
\end{equation}
where $f$ is given and $(a_n)$ is unknown. Define the partial sums
$A(N):=\sum_{k=1}^{N}a_k$ with $A(0):=0$.

One Abel summation is the tool used throughout this chapter.

\begin{proposition}\label{prop:abel}
For every $N\ge1$,
\[\sum_{k=1}^{N}a_k\,G(N,k)
=A(N)\,G(N,N)-\sum_{k=1}^{N-1}A(k)\big(G(N,k+1)-G(N,k)\big).
\]
Hence, defining the discrete Volterra\index[terms]{discrete Volterra} kernel
\[K(N,k):=G(N,k)-G(N,k+1)\qquad(1\le k\le N-1),
\]
equation \eqref{eq:primitive_disc} is equivalent to the discrete Volterra\index[terms]{discrete Volterra}
equation
\begin{equation}\label{eq:volterra_general}
A(N)\,G(N,N)+\sum_{k=1}^{N-1}K(N,k)\,A(k)=f(N).
\end{equation}
\end{proposition}

\begin{proof}
Writing $a_k=A(k)-A(k-1)$ and shifting the index in the second piece,
\[
\sum_{k=1}^{N}a_kG(N,k)
=\sum_{k=1}^{N}A(k)G(N,k)-\sum_{k=1}^{N-1}A(k)G(N,k+1),
\]
and the terms with $k=N$ collect into $A(N)G(N,N)$ while the remaining ones give
the stated sum.
\end{proof}

\begin{remark}
The kernel $K(N,k)=G(N,k)-G(N,k+1)$ is the discrete partial derivative of $G$ in
its second argument. When $G(N,k)\approx g(k/N)$ one has
$K(N,k)\approx-g'(k/N)/N$, the Riemann sum of the differential $dg$, so the Abel
transform is a discrete integration by parts, with $K$ playing the part of the
kernel of the continuous Volterra identity of Chapter~\ref{chap:volterra}.
\end{remark}

In the orthorecursive case $G(N,k)=(2N-1)/(N+k-1)$ one has $G(N,N)=1$ and
\[K(N,k)=\frac{2N-1}{(N+k-1)(N+k)},
\]
so that \eqref{eq:volterra_general} becomes the exact discrete Volterra\index[terms]{discrete Volterra} equation
of the second kind
\begin{equation}\label{eq:volterra_ortho}
A(N)+\sum_{k=1}^{N-1}K(N,k)\,A(k)=\delta_{N,1}.
\end{equation}

\section{The discrete resolvent by finite Neumann series}
\label{sec:disc_resolvent}

Let $T$ be the strictly triangular operator
$(Tu)(N):=\sum_{k=1}^{N-1}K(N,k)\,u(k)$. In the orthorecursive case,
\eqref{eq:volterra_ortho} reads $(I+T)A=e_1$ with $e_1=(1,0,0,\ldots)$. The
algebraic fact driving everything is that $T$ is strictly triangular, so
$(T^{m}u)(N)$ vanishes as soon as $m\ge N$, and the Neumann series is finite row
by row.

\begin{definition}\label{def:discrete_resolvent}
The discrete resolvent of $K$ is the triangular kernel
\[R(N,k):=\sum_{m\ge1}(-1)^{m-1}\,(T^{m})(N,k),
\]
the sum being finite for each fixed pair $(N,k)$.
\end{definition}

The triangular operator is invertible by a Neumann series that terminates at each rank.

\begin{proposition}\label{prop:exact_inversion}
$(I+T)^{-1}=I-R$ exactly. In the orthorecursive case,
\begin{equation}\label{eq:A_resolvent_disc}
A(N)=\delta_{N,1}-R(N,1).
\end{equation}
\end{proposition}

\begin{proof}
Row by row the operator $T$ is nilpotent, so
$(I+T)\sum_{m\ge0}(-1)^{m}T^{m}=I$ termwise, every sum being finite, and
$\sum_{m\ge0}(-1)^{m}T^{m}=I-R$. Applying $I-R$ to $e_1$ gives
\eqref{eq:A_resolvent_disc}.
\end{proof}

The kernel of the inverse satisfies a recurrence of its own.

\begin{proposition}\label{prop:resolvent_eq_disc}
For $1\le k\le N-1$, the kernel $R(N,k)$ satisfies
\begin{equation}\label{eq:resolvent_eqn_disc}
R(N,k)+\sum_{j=k+1}^{N-1}K(N,j)\,R(j,k)=K(N,k).
\end{equation}
\end{proposition}

\begin{proof}
From the definition, $R=T-T\sum_{m\ge1}(-1)^{m-1}T^{m}=T-TR$, and reading the
$(N,k)$ entry of $TR$, which is $\sum_{j=k+1}^{N-1}K(N,j)R(j,k)$ by
triangularity, gives the display.
\end{proof}

This is the exact discrete form of the classical Volterra resolvent equation.
The continuous identity of Chapter~\ref{chap:volterra} is its limit.

\section{Arithmetic Mellin probes at finite level}
\label{sec:disc_probes}

Spectral information is extracted from the exact structure through finite
probes, the discrete counterparts of the arithmetic Mellin transform.

\begin{definition}\label{def:probes}
For $z\in\C$, $N\ge1$ and $1\le k\le N$, set $u_{N,z}(k):=(k/N)^{-z}$ and
$u_{N,z}(0):=0$. The finite probes attached to $G$, $K$, $R$ are
\begin{align*}
\mathcal{G}_N(z)&:=\sum_{k=1}^{N}\big(u_{N,z}(k)-u_{N,z}(k-1)\big)\,G(N,k),
\\
\mathcal{K}_N(z)&:=\sum_{k=1}^{N-1}u_{N,z}(k)\,K(N,k),\\
\mathcal{R}_N(z)&:=\sum_{k=1}^{N-1}u_{N,z}(k)\,R(N,k).\end{align*}
For fixed $N$ each probe is a finite sum of exponentials in $z$, hence entire.
\end{definition}

The finite probe of the kernel and the finite probe of its inverse are linked by a single
identity.

\begin{proposition}\label{prop:bridge}
For every $N\ge1$ and every $z\in\C$,
\begin{equation}\label{eq:bridge_GK}
\mathcal{G}_N(z)=G(N,N)+\mathcal{K}_N(z).
\end{equation}
In the orthorecursive case $G(N,N)=1$ this gives
$\mathcal{G}_N(z)=1+\mathcal{K}_N(z)$.
\end{proposition}

\begin{proof}
Apply Proposition~\ref{prop:abel} to $a_k=u_{N,z}(k)-u_{N,z}(k-1)$, whose
partial sums are $A(k)=u_{N,z}(k)$, and note $u_{N,z}(N)=1$.
\end{proof}

Read on the resolvent, the same identity gives a recurrence for its probe.

\begin{proposition}\label{prop:probe_resolvent}
For every $N\ge1$ and every $z\in\C$,
\begin{equation}\label{eq:RN_eqn_book}
\mathcal{R}_N(z)+\sum_{j=2}^{N-1}K(N,j)\Big(\frac jN\Big)^{-z}\mathcal{R}_j(z)
=\mathcal{K}_N(z).
\end{equation}
\end{proposition}

\begin{proof}
Multiply \eqref{eq:resolvent_eqn_disc} by $(k/N)^{-z}$ and sum over
$1\le k\le N-1$. The first and last terms give $\mathcal{R}_N(z)$ and
$\mathcal{K}_N(z)$. In the middle term the two finite sums are interchanged and
the profile is factored,
\[
\sum_{k=1}^{N-1}\Big(\frac kN\Big)^{-z}\sum_{j=k+1}^{N-1}K(N,j)R(j,k)
=\sum_{j=2}^{N-1}K(N,j)\Big(\frac jN\Big)^{-z}
\sum_{k=1}^{j-1}\Big(\frac kj\Big)^{-z}R(j,k),
\]
the inner sum being $\mathcal{R}_j(z)$.
\end{proof}

Identity \eqref{eq:RN_eqn_book} expresses the resolvent probe at level $N$
through the resolvent probes at strictly lower levels. It is the exact discrete
form of the continuous relation $R^*=K^*/(1+K^*)$.

\section{The transform as a limit object}
\label{sec:disc_limit}

The probes are exact finite objects, and their significance lies in the limit
$N\to\infty$. The next theorem settles the convergence of the primitive probe
for a class of kernels that covers every example of this monograph.

\begin{theorem}\label{thm:probe_convergence}
Let $g:(0,1]\to\R$ be bounded and Riemann integrable, and let $\Re z<0$.
\begin{enumerate}
\item[(i)] If $G(N,k)=g(k/N)$, then
\[
\mathcal{G}_N(z)\longrightarrow
G^*(z):=-z\int_0^1g(t)\,t^{-z-1}\,dt\qquad(N\to\infty).
\]
\item[(ii)] The same limit holds whenever
$\sup_{1\le k\le N}\big|G(N,k)-g(k/N)\big|=\mathcal{O}(1/N)$.
\end{enumerate}
In both cases the convergence is locally uniform on compact subsets of the half
plane $\Re z<0$, and $\mathcal{K}_N(z)\to G^*(z)-g(1)$ whenever
$G(N,N)\to g(1)$.
\end{theorem}

\begin{proof}
Since $\tfrac{d}{dt}t^{-z}=-z\,t^{-z-1}$ and $t^{-z}\to0$ as $t\to0^{+}$ for
$\Re z<0$, each profile increment is an integral,
\[
u_{N,z}(k)-u_{N,z}(k-1)=-z\int_{(k-1)/N}^{k/N}t^{-z-1}\,dt ,
\]
so that, with the step function $g_N(t):=G(N,\lceil Nt\rceil)$,
\[
\mathcal{G}_N(z)=-z\int_0^1g_N(t)\,t^{-z-1}\,dt .
\]
In case (i), $g_N(t)=g(\lceil Nt\rceil/N)$ converges to $g(t)$ at every
continuity point of $g$, hence almost everywhere by the Lebesgue criterion for
Riemann integrability, while $|g_N|\le\sup|g|$ and $|t^{-z-1}|=t^{-\Re z-1}$ is
integrable on $(0,1)$. Dominated convergence gives the limit. In case (ii) the
extra term is bounded by
\[
\sup_k\big|G(N,k)-g(k/N)\big|\cdot
\sum_{k=1}^{N}\big|u_{N,z}(k)-u_{N,z}(k-1)\big|
\le\frac{C}{N}\cdot|z|\int_0^1t^{-\Re z-1}\,dt
=\frac{C\,|z|}{N\,|\Re z|},
\]
which tends to $0$. For the local uniformity, on a compact subset of
$\Re z<0$ the probes are uniformly bounded, since
$|\mathcal{G}_N(z)|\le\sup_N\sup_k|G(N,k)|\cdot|z|/|\Re z|$, so pointwise
convergence of these entire functions upgrades to locally uniform convergence
by Vitali's theorem. The statement for $\mathcal{K}_N$ follows from the bridge
identity \eqref{eq:bridge_GK}.
\end{proof}

The orthorecursive kernel falls under case (ii) with $g(t)=2/(1+t)$, since
\[
G(N,k)-g\Big(\frac kN\Big)
=\frac{2N-1}{N+k-1}-\frac{2N}{N+k}
=\frac{N-k}{(N+k-1)(N+k)},
\qquad
\Big|G(N,k)-g\Big(\frac kN\Big)\Big|\le\frac1N ,
\]
and the Ingham kernel\index[terms]{Ingham kernel} falls under case (i) with $g=\Phi$, which is bounded and
Riemann integrable. In the orthorecursive case the limit is the transform of
Chapter~\ref{chap:ortho}, and the bridge identity converges to
$G^*(z)=1+K^*(z)$.

The probe equation \eqref{eq:RN_eqn_book} converges, formally, to the
continuous resolvent identity
\begin{equation}\label{eq:Rstar_limit}
R^*(z)=\frac{K^*(z)}{1+K^*(z)},
\end{equation}
which is $\widetilde R(s)=1-1/G^*(-s)$ in the notation of
Chapter~\ref{chap:volterra}. The convergence of the resolvent probes
$\mathcal{R}_N$ is not established in general, and Section~\ref{sec:RAF_program}
records it among the open steps of the program.

\begin{remark}\label{rem:probe_numerics}
The rate in case (ii) is visible numerically. For the orthorecursive kernel at
$z=-0.7+0.3i$ the gap $|\mathcal{G}_N(z)-g^*(z)|$ multiplied by $N$ stabilizes
near $0.18$ for $N$ between $10^{2}$ and $1.6\cdot10^{3}$, the bridge identity
holds to machine precision at every level, and the resolvent probe
$\mathcal{R}_N(z)$ agrees with $K^*/(1+K^*)$ within $2\cdot10^{-4}$ at $N=300$.
The limiting resolvent identity, though unproved, is thus supported at finite
level.
\end{remark}

\begin{remark}\label{rem:conceptual_inversion}
The traditional presentation postulates the transform and derives the
asymptotics from its analytic properties. The discrete method inverts the
logic. The transform is the limit of exact finite objects, the resolvent
identity is the limit of an exact matrix identity, and the spectral content of
the zeros of $G^*$ is already present, discretely, at every finite level.
\end{remark}

\section{The linear kernel in the $\lambda$ formalism}
\label{sec:disc_linear}

The method is now run on the kernels where it closes exactly, beginning with
the model where every step is explicit. The linear case already exhibits the
three asymptotic behaviors of the theory, transparency, absorption, and the
critical logarithm, and it fixes the shape of the transfer theorem to come.

Let $0<\lambda<1$ and consider the linear kernel
\[G(n,k)=\lambda+(1-\lambda)\frac kn,\qquad1\le k\le n,
\]
with diagonal value $G(n,n)=1$. As forcing term\index[terms]{forcing term} take the exact Gamma--Euler
mode
\[b_n^{(\beta)}:=\frac{\Gamma(n-\beta)}{\Gamma(1-\beta)\,\Gamma(n)},
\]
which satisfies $b_{n+1}^{(\beta)}/b_n^{(\beta)}=(n-\beta)/n$, hence
$n(b_{n+1}^{(\beta)}-b_n^{(\beta)})=-\beta\,b_n^{(\beta)}$, and
$b_n^{(\beta)}\sim n^{-\beta}/\Gamma(1-\beta)$.

\subsection{Exact Abel form}

Since $G(n,k)-G(n,k+1)=-(1-\lambda)/n$, Proposition~\ref{prop:abel} gives
\begin{equation}\label{eq:linear_Abel_book}
A_n-\frac{1-\lambda}{n}\sum_{k=1}^{n-1}A_k=b_n^{(\beta)} .
\end{equation}
Introducing $B_n:=\sum_{j=1}^{n}A_j$, this becomes the first order recurrence
\begin{equation}\label{eq:linear_B_book}
B_n=\Big(1+\frac{1-\lambda}{n}\Big)B_{n-1}+b_n^{(\beta)} .
\end{equation}

\subsection{The transform and the spectral threshold}

The transform of $g(t)=\lambda+(1-\lambda)t$ is
\[g^*(z)=\frac{\lambda-z}{1-z},
\]
with a single zero at $z=\lambda$. The threshold between the asymptotic behaviors
is exactly that zero.

\subsection{The homogeneous mode}

For $b_n^{(\beta)}\equiv0$, iterating \eqref{eq:linear_B_book} gives
$B_n=C\,\Gamma(n+2-\lambda)/(\Gamma(2-\lambda)\Gamma(n+1))$ and hence
\[A_n=B_n-B_{n-1}
=C'\,\frac{\Gamma(n+1-\lambda)}{\Gamma(1-\lambda)\,\Gamma(n+1)},
\]
the exact Gamma--Euler mode attached to the spectral exponent $\lambda$.

\subsection{Exact solution away from resonance}

Away from resonance the linear kernel is solved exactly.

\begin{theorem}\label{thm:linear_exact_book}
Assume $\beta\neq\lambda$. Then
\begin{align}
B_n&=\frac{1}{(\beta+1-\lambda)(\beta-\lambda)}
\left[\frac{\Gamma(n+2-\lambda)}{\Gamma(1-\lambda)\,\Gamma(n+1)}
-\frac{\big((\beta+1-\lambda)n+1-\lambda\big)\,\Gamma(n+1-\beta)}
{\Gamma(1-\beta)\,\Gamma(n+1)}\right],\label{eq:linear_B_exact_book}\\[2pt]
A_n&=\frac{1-\lambda}{(\beta+1-\lambda)(\beta-\lambda)}
\frac{\Gamma(n+1-\lambda)}{\Gamma(1-\lambda)\,\Gamma(n+1)}
+\left(\frac{\beta-1}{\beta-\lambda}
+\frac{\beta(1-\lambda)}{n(\beta+1-\lambda)(\beta-\lambda)}\right)
\frac{\Gamma(n-\beta)}{\Gamma(1-\beta)\,\Gamma(n)}.\label{eq:linear_A_exact_book}
\end{align}
\end{theorem}

\begin{proof}
Set $P_n:=\prod_{j=1}^{n}\big(1+\tfrac{1-\lambda}{j}\big)
=\Gamma(n+2-\lambda)/(\Gamma(2-\lambda)\Gamma(n+1))$, so that
\eqref{eq:linear_B_book} reads
$B_n/P_n-B_{n-1}/P_{n-1}=b_n^{(\beta)}/P_n$, with
\[
\frac{b_n^{(\beta)}}{P_n}
=\frac{\Gamma(2-\lambda)}{\Gamma(1-\beta)}\cdot
\frac{n\,\Gamma(n-\beta)}{\Gamma(n+2-\lambda)} .
\]
Write $p:=\beta+1-\lambda$, $q:=1-\lambda$, $D:=p(\beta-\lambda)$ and define
\[
T_n:=\frac{(pn+q)\,\Gamma(n+1-\beta)}{D\,\Gamma(n+2-\lambda)} .
\]
Bringing $T_{n-1}$ and $T_n$ over the common denominator
$D\,\Gamma(n+2-\lambda)$, using $\Gamma(n+1-\lambda)=\Gamma(n+2-\lambda)/(n+1-\lambda)$
and $\Gamma(n+1-\beta)=(n-\beta)\Gamma(n-\beta)$, the numerator of
$T_{n-1}-T_n$ equals
\[
\big(p(n-1)+q\big)(n+q)-\big(pn+q\big)(n-\beta)
=pn(q-1+\beta)+q(q+\beta-p)=Dn,
\]
since $q-1+\beta=\beta-\lambda$ and $q+\beta-p=0$. Hence
$T_{n-1}-T_n=n\,\Gamma(n-\beta)/\Gamma(n+2-\lambda)$, the sum telescopes to
$B_n/P_n=\tfrac{\Gamma(2-\lambda)}{\Gamma(1-\beta)}(T_0-T_n)$, and multiplying
back by $P_n$ gives \eqref{eq:linear_B_exact_book}. Formula
\eqref{eq:linear_A_exact_book} follows from
$A_n=b_n^{(\beta)}+\tfrac{1-\lambda}{n}B_{n-1}$, which is
\eqref{eq:linear_Abel_book}.
\end{proof}

The two regimes are read off the exact formula \eqref{eq:linear_A_exact_book} by the
asymptotic of a Gamma quotient, with no further input.

\begin{corollary}\label{cor:linear_cases_book}
With $\Gamma(n+a)/\Gamma(n+b)\sim n^{a-b}$ applied to \eqref{eq:linear_A_exact_book},
\begin{enumerate}
\item[(i)] if $\beta<\lambda$ then
$A_n\sim\dfrac{1-\beta}{(\lambda-\beta)\,\Gamma(1-\beta)}\,n^{-\beta}$,
which is $A_n\sim b_n^{(\beta)}/g^*(\beta)$ since
$g^*(\beta)=(\lambda-\beta)/(1-\beta)$,
\item[(ii)] if $\beta>\lambda$ then
$A_n\sim\dfrac{1-\lambda}{(\beta+1-\lambda)(\beta-\lambda)\,\Gamma(1-\lambda)}
\,n^{-\lambda}$.
\end{enumerate}
The transparency constant in (i) is exactly $1/g^*(\beta)$, the first
appearance of the general rule.
\end{corollary}

\subsection{The resonant case}

At resonance the solution acquires a logarithm.

\begin{theorem}\label{thm:linear_resonance_book}
Assume $\beta=\lambda$. Then, with $P_n$ as above,
\begin{equation}\label{eq:linear_resonant_exact}
\frac{B_n}{P_n}=(1-\lambda)\Big[\lambda\big(\psi(n+1-\lambda)-\psi(1-\lambda)\big)
+(1-\lambda)\big(\psi(n+2-\lambda)-\psi(2-\lambda)\big)\Big],
\end{equation}
and consequently
\begin{equation}\label{eq:linear_resonant_book}
A_n\sim\frac{1-\lambda}{\Gamma(1-\lambda)}\,n^{-\lambda}\log n .
\end{equation}
\end{theorem}

\begin{proof}
With $\beta=\lambda$ the ratio computed in the previous proof becomes
\[
\frac{b_n^{(\lambda)}}{P_n}
=(1-\lambda)\,\frac{n}{(n-\lambda)(n+1-\lambda)}
=(1-\lambda)\Big(\frac{\lambda}{n-\lambda}+\frac{1-\lambda}{n+1-\lambda}\Big),
\]
the last step by partial fractions. Summing over $1\le m\le n$ and using
$\sum_{m=1}^{n}\tfrac1{m+c}=\psi(n+1+c)-\psi(1+c)$, which is the difference
form of the series for $\psi$ used in Chapter~\ref{chap:ortho}, gives
\eqref{eq:linear_resonant_exact}. Since $\psi(x)=\log x+\mathcal{O}(1)$ and
$P_n\sim n^{1-\lambda}/\Gamma(2-\lambda)$, this yields
$B_n\sim\tfrac{1-\lambda}{\Gamma(2-\lambda)}n^{1-\lambda}\log n$, and
$A_n=b_n^{(\lambda)}+\tfrac{1-\lambda}{n}B_{n-1}$ turns it into
\eqref{eq:linear_resonant_book}.
\end{proof}

\begin{remark}
The linear kernel already carries the full trichotomy. Below the zero of the
transform the forcing passes with the constant $1/g^*(\beta)$, above it the
homogeneous Gamma--Euler mode dominates, and at the zero a logarithm appears.
This is the transfer theorem in miniature, and the closed form
\eqref{eq:linear_resonant_exact} was checked against the recurrence to
$10^{-13}$ at $n=2000$.
\end{remark}

\section{A transfer theorem for polynomial kernels}
\label{sec:disc_poly}

The linear case hides the finite dimensional mechanism behind a one
dimensional closure. Polynomial kernels are the first class where the
algebraic structure becomes visible. After exact Abel summation\index[terms]{Abel summation} the problem
closes on finitely many moments, and the transform appears as the indicial
object of the resulting system.

Let $G(n,k)=P(k/n)$ with $P(t)=\sum_{j=0}^{d}c_jt^{j}$ and $P(1)=1$. Define
$A_n:=\sum_{k=1}^{n}a_k$ and the moments $M_m(n):=\sum_{k=1}^{n}k^{m}A_k$ for
$0\le m\le d-1$.

For a polynomial kernel the Abel transform closes on a finite system.

\begin{proposition}\label{prop:poly_system}
The exact Abel transform of the equation with forcing $b_n^{(\beta)}$ reads
\begin{equation}\label{eq:A_exact_moments_book}
A_n-\sum_{m=0}^{d-1}q_m(n)\,M_m(n-1)=b_n^{(\beta)},\qquad
q_m(n)=\sum_{j=m+1}^{d}c_j\binom jm n^{-j},
\end{equation}
and the normalized moments $U_m(n):=M_m(n)/n^{m}$ satisfy
\begin{equation}\label{eq:poly_U_system}
U(n)=\Big(I+\frac1nM+\mathcal{O}\Big(\frac1{n^{2}}\Big)\Big)U(n-1)
+b_n^{(\beta)}\,\mathbf 1+\mathcal{O}\Big(\frac{b_n^{(\beta)}}{n}\Big),
\end{equation}
where $M:=\mathbf 1\alpha^{\top}-D$ with $\alpha_m:=(m+1)c_{m+1}$ and
$D:=\operatorname{diag}(0,1,\ldots,d-1)$.
\end{proposition}

\begin{proof}
The Volterra kernel\index[terms]{Volterra kernel} is
\[
K(n,k)=P\Big(\frac kn\Big)-P\Big(\frac{k+1}n\Big)
=-\sum_{j=1}^{d}c_jn^{-j}\big((k+1)^{j}-k^{j}\big)
=-\sum_{m=0}^{d-1}\Big(\sum_{j>m}c_j\binom jm n^{-j}\Big)k^{m}
\]
by the binomial theorem, and summing against $A_k$ gives
\eqref{eq:A_exact_moments_book} since $G(n,n)=P(1)=1$. For the system, from
$M_m(n)=M_m(n-1)+n^{m}A_n$, hence
$U_m(n)=(1-\tfrac mn+\mathcal{O}(n^{-2}))U_m(n-1)+A_n$, while
\eqref{eq:A_exact_moments_book} with $M_m(n-1)=(n-1)^{m}U_m(n-1)$ and
$q_m(n)(n-1)^{m}=(m+1)c_{m+1}n^{-1}+\mathcal{O}(n^{-2})$ gives
$A_n=b_n^{(\beta)}+\tfrac1n\sum_m\alpha_mU_m(n-1)+\mathcal{O}(n^{-2})\|U(n-1)\|$.
Substituting produces \eqref{eq:poly_U_system} with the stated $M$.
\end{proof}

The order term in \eqref{eq:poly_U_system} is the truncation of an exact
correction. Theorem~\ref{thm:trace-poly-finite} carries the exact transfer matrix
$T_n=I+M/n+N/n^{2}$, with the block $N/n^{2}$ nonzero and not commuting with $M$
in general, together with the finite resolvent identity it produces.

\begin{proposition}\label{prop:char_identity_book}
With the notation of Proposition~\ref{prop:poly_system},
\[\det\big((1-z)I-M\big)=(1-z)(2-z)\cdots(d-z)\;g^*(z),
\qquad
g^*(z)=\sum_{j=0}^{d}c_j\,\frac{-z}{j-z},
\]
so, once any factor common to numerator and denominator is cancelled, the
effective zeros $\rho$ of the reduced transform correspond to the eigenvalues
$1-\rho$ of $M$. A factor shared with a pole at $z\in\{1,\dots,d\}$ carries no
spectral mode. In the quadratic case $P(1)=-b$ and $P(2)=2c$, so such a
cancellation occurs exactly at $b=0$ or at $c=0$.
\end{proposition}

\begin{proof}
The rank one update formula for determinants gives
$\det((1-z)I+D-\mathbf 1\alpha^{\top})
=\det((1-z)I+D)\big(1-\alpha^{\top}((1-z)I+D)^{-1}\mathbf 1\big)$.
The first factor is $\prod_{m=0}^{d-1}(m+1-z)$. For the second,
\[
\alpha^{\top}\big((1-z)I+D\big)^{-1}\mathbf 1
=\sum_{m=0}^{d-1}\frac{(m+1)c_{m+1}}{m+1-z}
=\sum_{j=1}^{d}\frac{j\,c_j}{j-z}=1-g^*(z),
\]
the last equality since
$g^*(z)=c_0+\sum_{j\ge1}c_j\tfrac{-z}{j-z}$ and $\sum_jc_j=1$ give
$1-g^*(z)=\sum_{j\ge1}c_j\big(1+\tfrac z{j-z}\big)=\sum_{j\ge1}c_j\tfrac j{j-z}$.
\end{proof}

Under separation of the zeros the transfer is complete.

\begin{theorem}\label{thm:poly_transfer_book}
Assume the zeros $\rho_1,\ldots,\rho_d$ of $g^*$ are simple, have distinct real
parts, no two differ by an integer, and $g^*(\beta)\neq0$. Then
\begin{equation}\label{eq:poly_transfer_book}
A_n=\frac{1}{g^*(\beta)}\,b_n^{(\beta)}
+\sum_{j=1}^{d}C_j\,\frac{\Gamma(n+1-\rho_j)}{\Gamma(n+1)}
+o\big(n^{-m_*}\big),
\qquad m_*:=\min\big(\beta,\Re\rho_1,\ldots,\Re\rho_d\big).
\end{equation}
\end{theorem}
\begin{proof}
The system \eqref{eq:poly_U_system} has the form
$U(n)=(I+M/n+V(n))U(n-1)+F(n)$ with $\sum\|V(n)\|<\infty$. The asymptotic
theory of such perturbed linear difference systems is the discrete Levinson\index[names]{Levinson, N.}
theory of Benzaid\index[names]{Benzaid, Z.} and Lutz\index[names]{Lutz, D. A.} \cite{BenzaidLutz1987}. What it supplies is the following.
The principal matrix produces one mode per eigenvalue, and a summable perturbation leaves the
exponents of those modes unchanged, altering only their amplitudes. Their dichotomy condition
holds here because the diagonal entries $1+\mu_j/n$ built from the eigenvalues
$\mu_j=1-\rho_j$ of $M$ separate as powers, the ratio of the products up to $n$
being of exact order $n^{\Re(\mu_i-\mu_j)}$, and the perturbation is summable.
The homogeneous system therefore has a basis of solutions
\[
U^{(j)}(n)=\big(v_j+o(1)\big)\prod_{k=1}^{n}\Big(1+\frac{\mu_j}{k}\Big)
=\big(v_j+o(1)\big)\,
\frac{\Gamma(n+1+\mu_j)}{\Gamma(1+\mu_j)\,\Gamma(n+1)},
\]
where $v_j$ is the eigenvector of $M$ at $\mu_j$, and passing from $U$ back to
$A_n$ through \eqref{eq:A_exact_moments_book} turns these into the
Gamma--Euler modes of \eqref{eq:poly_transfer_book}.

The amplitude of the forcing mode is exact algebra. Seek a particular solution
$U(n)=n\,b_n^{(\beta)}\,\kappa$ plus lower order terms. Since
$(n-1)b_{n-1}^{(\beta)}=b_n^{(\beta)}\big(n+\beta-1+\mathcal{O}(1/n)\big)$,
inserting the trial form into \eqref{eq:poly_U_system} and matching the leading
terms gives $\big((1-\beta)I-M\big)\kappa=\mathbf 1$. The Sherman--Morrison
formula for the rank one update $M=\mathbf 1\alpha^{\top}-D$ then gives, with
$S:=(1-\beta)I+D$,
\[
\kappa=\big(S-\mathbf 1\alpha^{\top}\big)^{-1}\mathbf 1
=\frac{S^{-1}\mathbf 1}{1-\alpha^{\top}S^{-1}\mathbf 1}
=\frac{S^{-1}\mathbf 1}{g^*(\beta)},
\]
the denominator by the computation in the proof of
Proposition~\ref{prop:char_identity_book}. Feeding this back into
$A_n=b_n^{(\beta)}+\tfrac1n\alpha^{\top}U(n-1)+\mathcal{O}(n^{-2})\|U\|$
yields the amplitude
\[
1+\alpha^{\top}\kappa
=1+\frac{\alpha^{\top}S^{-1}\mathbf 1}{g^*(\beta)}
=1+\frac{1-g^*(\beta)}{g^*(\beta)}=\frac{1}{g^*(\beta)},
\]
and the remainder estimates of \nm{Benzaid}{Z.} and \nm{Lutz}{D. A.}~\cite{BenzaidLutz1987} for the variation of
constants give the error term $o(n^{-m_*})$.
\end{proof}

Under the genericity assumed above the discrete Levinson\index[names]{Levinson, N.} dichotomy of
\cite{BenzaidLutz1987} applies to \eqref{eq:poly_U_system} and the principal
products give the stated modes. The excluded configurations, equal real parts,
integer separations, and multiple roots, together with a self-contained
treatment, are the exact finite resolvent identity of
Theorem~\ref{thm:trace-poly-finite} and the quadratic spectral expansion of
Theorems~\ref{thm:trace-quad-simple} and~\ref{thm:trace-quad-exceptional}.

Reading the principal products of the preceding paragraph according to which exponent
dominates gives the two regimes at once.

\begin{corollary}\label{cor:poly_cases}
Under the genericity assumed above, if $\beta<\min_j\Re\rho_j$ then
$A_n\sim b_n^{(\beta)}/g^*(\beta)$, the transparent case. If $\beta>\min_j\Re\rho_j$ the dominant contribution comes
from the zeros of $g^*$ of smallest real part, the absorbed case.
\end{corollary}

When the forcing exponent meets a zero the same computation gives the resonant form.

\begin{proposition}\label{thm:poly_resonance_book}
If $\beta=\rho_\ell$ is a simple zero of $g^*$, the remaining zeros are simple
with distinct real parts and no integer separation from $\rho_\ell$ or from one
another, and the forcing has a nonzero component on the corresponding
eigendirection, then
\begin{equation}\label{eq:poly_resonance_book}
A_n=\widetilde C\,\frac{\Gamma(n+1-\beta)}{\Gamma(n+1)}\,\log n
+\sum_{j=1}^{d}\widetilde C_j\,\frac{\Gamma(n+1-\rho_j)}{\Gamma(n+1)}
+o\big(n^{-\Re\beta}\log n\big),
\qquad\widetilde C\neq0 .
\end{equation}
\end{proposition}

\begin{proof}
Diagonalize the leading matrix. On the non resonant eigendirections the
analysis of Theorem~\ref{thm:poly_transfer_book} applies unchanged. On the
resonant direction the scalar component satisfies
$y_n=\big(1+\tfrac{1-\beta}{n}+\mathcal{O}(n^{-2})\big)y_{n-1}
+c\,b_n^{(\beta)}+\mathcal{O}\big(b_n^{(\beta)}n^{-1}\big)$
with $c\neq0$ by hypothesis, and dividing by the integrating factor
$P_n=\prod_{k\le n}(1+\tfrac{1-\beta}{k})\asymp n^{1-\beta}$ reduces it to the
computation of Theorem~\ref{thm:linear_resonance_book}, the ratio
$b_n^{(\beta)}/P_n$ being of exact order $1/n$, whose sum is the logarithm.
Passing back to $A_n$ gives \eqref{eq:poly_resonance_book}.
\end{proof}

The resonance at a cancelled factor, a double root, or an integer separation of
the roots falls outside the genericity assumed here. Those cases, with the
connection coefficients built by Wronskians\index[terms]{Wronskian} and the simple and double resonances
obtained as finite parts, are Theorem~\ref{thm:trace-quad-exceptional}.

\section{The quadratic kernel and the discrete Euler--Cauchy mechanism}
\label{sec:disc_quadratic}

The quadratic case is the smallest laboratory where the second order geometry
becomes visible, two spectral roots, interaction, multiplicity, oscillation.
It is the bridge between the fully solved linear model and the orthorecursive
kernel.

Let $G(n,k)=a+b(k/n)+c(k/n)^{2}$ with $a+b+c=1$. The transform factors as
\begin{equation}\label{eq:quadratic_transform_book}
g^*(z)=\frac{(z-\rho_1)(z-\rho_2)}{(z-1)(z-2)},
\end{equation}
where $\rho_1,\rho_2$ are the roots of $z^{2}-(3-b-2c)z+2a=0$.

The quadratic case closes on the partial sums alone.

\begin{proposition}\label{prop:quad_indicial_book}
With $B_n:=\sum_{k\le n}A_k$, $C_n:=\sum_{k\le n}kA_k$ and $D_n:=C_n/n$, the
quadratic equation closes on the pair $(B_n,D_n)$ and takes the form
\[\binom{B_n}{D_n}
=\Big(I+\frac1nM+\mathcal{O}\Big(\frac1{n^{2}}\Big)\Big)\binom{B_{n-1}}{D_{n-1}}
+b_n^{(\beta)}\binom11+\mathcal{O}\Big(\frac{b_n^{(\beta)}}{n}\Big),
\qquad
M=\begin{pmatrix}b&2c\\ b&2c-1\end{pmatrix},
\]
and the eigenvalues of $M$ are $\mu_j=1-\rho_j$ with $\rho_1,\rho_2$ the zeros
of $g^*$.
\end{proposition}

\begin{proof}
Here $d=2$, $q_0(n)=\tfrac bn+\tfrac c{n^{2}}$ and $q_1(n)=\tfrac{2c}{n^{2}}$,
so \eqref{eq:A_exact_moments_book} reads
$A_n=b_n^{(\beta)}+\big(\tfrac bn+\tfrac c{n^{2}}\big)B_{n-1}
+\tfrac{2c}{n^{2}}C_{n-1}$. With $C_{n-1}=(n-1)D_{n-1}$ this is
$A_n=b_n^{(\beta)}+\tfrac bnB_{n-1}+\tfrac{2c}nD_{n-1}+\mathcal{O}(n^{-2})
\big(|B_{n-1}|+|D_{n-1}|\big)$. Then $B_n=B_{n-1}+A_n$ gives the first row,
and $D_n=\tfrac{C_{n-1}+nA_n}{n}=\big(1-\tfrac1n\big)D_{n-1}+A_n$ gives the
second, with the stated matrix. Its characteristic polynomial is
$\mu^{2}-(b+2c-1)\mu-b$, and the substitution $\mu=1-z$ together with
$a+b+c=1$ turns it into $z^{2}-(3-b-2c)z+2a=(z-\rho_1)(z-\rho_2)$, the
numerator of \eqref{eq:quadratic_transform_book}.
\end{proof}

The matrix $M$ and its eigenvalues $\mu_j=1-\rho_j$ are exact. The displayed
system is only asymptotic, the exact system with its propagator and summable
correction being Theorem~\ref{thm:trace-poly-finite}.

Three families span the phenomenology. Two real roots separated by an integer,
$\rho_1=\tfrac12$ and $\rho_2=\tfrac32$ with $(a,b,c)=(\tfrac38,\tfrac14,\tfrac38)$,
give a dominant Gamma--Euler mode $M_{1/2}$ and a subordinate mode $M_{3/2}$ that
acquires a logarithm from the integer gap. A double root, $\rho_1=\rho_2=\tfrac12$,
produces the logarithmic companion. A complex pair, $\rho_{1,2}=\tfrac12\pm i$
with $(a,b,c)=(\tfrac58,-\tfrac54,\tfrac{13}8)$, produces oscillation in $\log n$.
Two real roots with non-integer separation give two clean modes with no
logarithm. The complete treatment of the four cases, with the connection
coefficients and the resonances, is Theorems~\ref{thm:trace-quad-simple}
and~\ref{thm:trace-quad-exceptional}.

\begin{remark}\label{rem:quad_numerics}
The three behaviors were verified on the recurrences. For the real pair,
$A_n\Gamma(1-\beta)n^{\beta}$ at $\beta=0.2$ climbs to within seven percent of
the predicted constant $1/g^*(0.2)$ by $n=2500$ at the expected rate
$n^{-0.3}$, and at $\beta=0.9$ the normalized sums $A_nn^{1/2}$ settle at
$0.710$ to three digits. For the complex pair at $\beta=0.9$ the same
normalized sums oscillate between $0.85$ and $1.27$ with the period predicted
by $\Im\rho_1=1$, the signature that no real spectral model can reproduce.
\end{remark}

\section{Gamma--Newton expansion of the orthorecursive kernel}
\label{sec:disc_GN}

The orthorecursive kernel admits an exact discrete expansion in a basis
adapted to the diagonal $k\approx N$, stable under Abel differences, and whose
blocks recover the transform term by term.

\subsection{A Newton--Gamma basis centered at the diagonal}

For integers $N\ge2$, $1\le k\le N$ and $m\ge0$, define
\[\Phi_m(N,k):=\frac{(N-k)^{\underline m}}{(2N-2)^{\underline m}},
\]
where $x^{\underline m}:=x(x-1)\cdots(x-m+1)$ is the falling factorial and
$\Phi_0(N,k)=1$. Since $(N-k)^{\underline m}=0$ for $m>N-k$, the family is
finite at each lattice point.

\begin{proposition}\label{prop:GN_exact}
For every $N\ge2$ and $1\le k\le N$,
\begin{equation}\label{eq:G_exact_series}
G(N,k)=\sum_{m=0}^{N-k}\Phi_m(N,k).
\end{equation}
\end{proposition}

\begin{proof}
Set $x:=N-k$ and $A:=2N-2$, so that $G(N,k)=(A+1)/(A+1-x)$ and the claim is
\begin{equation}\label{eq:finite_newton}
\sum_{m=0}^{x}\frac{x^{\underline m}}{A^{\underline m}}=\frac{A+1}{A+1-x},
\end{equation}
to be proved for all integers $0\le x<A$ by induction on $x$, the identity
being established for every admissible $A$ at once. For $x=0$ both sides equal
$1$. Assuming it at level $x$ for every parameter, split the sum at level
$x+1$ into the term $m=0$ and the rest, where
$(x+1)^{\underline m}=(x+1)\,x^{\underline{m-1}}$ and
$A^{\underline m}=A\,(A-1)^{\underline{m-1}}$, so that
\[
\sum_{m=0}^{x+1}\frac{(x+1)^{\underline m}}{A^{\underline m}}
=1+\frac{x+1}{A}\sum_{m=0}^{x}\frac{x^{\underline m}}{(A-1)^{\underline m}}
=1+\frac{x+1}{A}\cdot\frac{A}{A-x}
=\frac{A+1}{A-x},
\]
which is \eqref{eq:finite_newton} at level $x+1$ with parameter $A$.
\end{proof}

\begin{remark}
The expansion is the Newton expansion of $x\mapsto(2N-1)/(2N-1-x)$ in the
variable $x=N-k$, centered at the diagonal $k=N$.
\end{remark}

\subsection{Exact Abel differences}

The differences of the truncated kernels are computable in closed form.

\begin{proposition}\label{prop:Phi_difference}
For every $m\ge1$,
\[\Phi_m(N,k)-\Phi_m(N,k+1)
=\frac{m}{2N-2}\cdot
\frac{(N-k-1)^{\underline{m-1}}}{(2N-3)^{\underline{m-1}}}.
\]
\end{proposition}

\begin{proof}
With $x=N-k$, the falling factorials satisfy
$x^{\underline m}-(x-1)^{\underline m}
=(x-1)^{\underline{m-1}}\big(x-(x-m)\big)=m\,(x-1)^{\underline{m-1}}$,
since $x^{\underline m}=x\,(x-1)^{\underline{m-1}}$ and
$(x-1)^{\underline m}=(x-1)^{\underline{m-1}}(x-m)$. Dividing by
$(2N-2)^{\underline m}=(2N-2)\,(2N-3)^{\underline{m-1}}$ gives the display.
\end{proof}

The kernel of the inverse expands accordingly.

\begin{corollary}\label{cor:K_expansion}
For $1\le k\le N-1$,
\begin{equation}\label{eq:K_expansion}
K(N,k)=\frac{1}{2N-2}\sum_{m=0}^{N-k-1}(m+1)\,
\frac{(N-k-1)^{\underline m}}{(2N-3)^{\underline m}}.
\end{equation}
\end{corollary}

\begin{proof}
Apply \eqref{eq:G_exact_series} to
$K(N,k)=\sum_{m\ge1}\big(\Phi_m(N,k)-\Phi_m(N,k+1)\big)$, insert
Proposition~\ref{prop:Phi_difference}, and shift the index.
\end{proof}

\subsection{Finite truncations and their continuous blocks}

For $M\ge0$ define the truncated kernel
$G^{(M)}(N,k):=\sum_{m=0}^{\min(M,N-k)}\Phi_m(N,k)$, with $G^{(M)}(N,N)=1$ and
$G^{(M)}(N,k)\uparrow G(N,k)$ as $M\to\infty$, and the corresponding truncated
Volterra kernel $K^{(M)}$ obtained by cutting \eqref{eq:K_expansion} at
$m\le M-1$. These are exact finite models, each fully compatible with the
formalism of Sections~\ref{sec:disc_primitive} and~\ref{sec:disc_resolvent}.

As $N\to\infty$ with $k/N\to t$, each factor of $\Phi_m(N,k)$ satisfies
$(N-k-i)/(2N-2-i)\to(1-t)/2$, so $\Phi_m(N,k)\to2^{-m}(1-t)^{m}$. The
continuous block attached to level $m$ is therefore
\[\Phi_m^*(z):=-z\int_0^1 2^{-m}(1-t)^{m}\,t^{-z-1}\,dt,\qquad\Re z<0 .
\]

Each truncation has a transform expressed through Gamma functions.

\begin{proposition}\label{prop:phi_transform}
For $m\ge0$ and $\Re z<0$,
\[\Phi_m^*(z)=2^{-m}\,
\frac{\Gamma(1-z)\,\Gamma(m+1)}{\Gamma(m+1-z)} .
\]
\end{proposition}

\begin{proof}
The Beta integral \cite[Eq.~5.12.1]{DLMF} gives
$\int_0^1(1-t)^{m}t^{-z-1}\,dt=\Gamma(-z)\Gamma(m+1)/\Gamma(m+1-z)$, and
$-z\,\Gamma(-z)=\Gamma(1-z)$.
\end{proof}

Summing the truncations recovers the transform of the kernel itself.

\begin{corollary}\label{cor:gstar_recovered}
For $\Re z<0$,
\[g^*(z)=-z\int_0^1\frac{2}{1+t}\,t^{-z-1}\,dt
=\sum_{m=0}^{\infty}\Phi_m^*(z)
=\sum_{m=0}^{\infty}2^{-m}\,
\frac{\Gamma(1-z)\,\Gamma(m+1)}{\Gamma(m+1-z)} .
\]
\end{corollary}

\begin{proof}
The expansion $2/(1+t)=\sum_{m\ge0}2^{-m}(1-t)^{m}$ has terms bounded by
$2^{-m}$ on $[0,1]$, so it converges uniformly there, termwise integration
against the integrable weight $t^{-z-1}$ is legitimate, and
Proposition~\ref{prop:phi_transform} identifies each term.
\end{proof}

\begin{remark}
The basis is not a formal device. It reconstructs, block by block, the same
transform that governs Chapter~\ref{chap:ortho}, and it is the natural
candidate for the discrete refinements of the program of
Section~\ref{sec:RAF_program}, in particular its last step.
\end{remark}

\section{The orthorecursive kernel as a limit of polynomial models}
\label{sec:disc_limit_models}

Truncating the geometric expansion of the profile places the orthorecursive
kernel at the end of a canonical sequence of polynomial models. For $M\ge0$
let $P_M(t):=\sum_{m=0}^{M}2^{-m}(1-t)^{m}$, so that $P_M(1)=1$ and
$\sup_{[0,1]}|2/(1+t)-P_M(t)|\le2^{-M}$, and let
\[g_M^*(z)=\sum_{m=0}^{M}2^{-m}\,
\frac{\Gamma(1-z)\,\Gamma(m+1)}{\Gamma(m+1-z)}
\longrightarrow g^*(z)\qquad(M\to\infty),
\]
locally uniformly on $\C\setminus\{1,2,3,\ldots\}$, since on a compact set
avoiding the poles the terms are $\mathcal{O}(2^{-m}m^{\Re z})$ by the ratio
asymptotics of the Gamma function\index[terms]{Gamma function}, and the Weierstrass test applies. By
Hurwitz's theorem every simple zero $\rho$ of $g^*$ is the limit of zeros
$\rho_M$ of $g_M^*$.

\begin{remark}\label{rem:hurwitz_numerics}
The approach is fast. The zero of $g_M^*$ near $\rho_1$ sits at
$1.3631+1.0590\,i$ for $M=8$, at $1.34661+1.05527\,i$ for $M=16$, and at
$1.3465165+1.0551601\,i$ for $M=32$, nine digits of $\rho_1$ recovered from a
degree $32$ polynomial model.
\end{remark}

Theorem~\ref{thm:poly_transfer_book} gives a complete transfer formula for
each truncated solution $A^{(M)}$. Whether the truncated solutions converge to
the orthorecursive one in a norm strong enough to carry the expansion over is
the stability question, and the honest state of that question is the
following.

For $\sigma\in\R$ define the Gamma--Euler weight
$w_\sigma(n):=\Gamma(n+1-\sigma)/\Gamma(n+1)\sim n^{-\sigma}$ and the Banach
space $\mathcal{X}_\sigma$ of sequences with finite norm
$\|u\|_\sigma:=\sup_{n\ge1}|u_n|/w_\sigma(n)$.

\begin{definition}\label{def:uniform_stability_book}
The family $(R_M)_{M\ge0}$ of discrete resolvents of the truncated kernels is
uniformly stable on $\mathcal{X}_\sigma$ if
\begin{equation}\label{eq:uniform_stability_book}
\sup_{M\ge0}\,\sup_{n\ge1}\,\frac{1}{w_\sigma(n)}
\sum_{k=1}^{n-1}\big|R_M(n,k)\big|\,w_\sigma(k)\le C_\sigma<\infty .
\end{equation}
\end{definition}

Truncating the kernel perturbs the solution, and the perturbation obeys an equation of the same
shape.

\begin{proposition}\label{prop:error_abstract_book}
Let $A$ and $A^{(M)}$ solve the exact and the truncated equations with the
same right side, and let $D^{(M)}:=A-A^{(M)}$. Then
\begin{equation}\label{eq:error_resolvent_book}
D_n^{(M)}=E_n^{(M)}-\sum_{k=1}^{n-1}R_M(n,k)\,E_k^{(M)},
\qquad
E_n^{(M)}:=\sum_{k=1}^{n-1}\big(K_M(n,k)-K(n,k)\big)A_k ,
\end{equation}
and if $(R_M)$ is uniformly stable on $\mathcal{X}_\sigma$ then
$\|D^{(M)}\|_\sigma\le(1+C_\sigma)\,\|E^{(M)}\|_\sigma$.
\end{proposition}

\begin{proof}
Write the two equations as $(I+T)A=f$ and $(I+T_M)A^{(M)}=f$. Subtracting,
$(I+T_M)D^{(M)}=(T_M-T)A=E^{(M)}$, and Proposition~\ref{prop:exact_inversion}
applied to $T_M$ gives $D^{(M)}=(I-R_M)E^{(M)}$, which is
\eqref{eq:error_resolvent_book}. The norm bound follows from the triangle
inequality and \eqref{eq:uniform_stability_book}.
\end{proof}

The truncated kernels approach the true one at an explicit rate.

\begin{proposition}\label{prop:consistency_book}
For all $M\ge0$, $N\ge2$ and $1\le k\le N-1$,
\[\big|K_M(N,k)-K(N,k)\big|\le\frac{2(M+2)}{2N-2}\,2^{-M},
\]
and consequently, if $\sum_k|A_k|<\infty$, then
$|E_n^{(M)}|\ll(M+2)\,2^{-M}n^{-1}$ and
$\|E^{(M)}\|_\sigma\ll(M+2)\,2^{-M}$ for every $\sigma\le1$.
\end{proposition}

\begin{proof}
The difference is the tail of \eqref{eq:K_expansion}, and each ratio of
falling factorials obeys
\[
\frac{(N-k-1)^{\underline m}}{(2N-3)^{\underline m}}
=\prod_{i=0}^{m-1}\frac{N-k-1-i}{2N-3-i}
\le\Big(\frac{N-k-1}{2N-3}\Big)^{m}\le2^{-m},
\]
each factor being at most the first and the first being at most $\tfrac12$
for $k\ge1$. Hence
\[
\big|K_M(N,k)-K(N,k)\big|
\le\frac{1}{2N-2}\sum_{m\ge M}(m+1)\,2^{-m}
=\frac{2(M+2)}{2N-2}\,2^{-M},
\]
the sum being evaluated exactly from the geometric series and its derivative.
The bound on $E_n^{(M)}$ follows by summing against $|A_k|$, and for
$\sigma\le1$ the quotient
$|E_n^{(M)}|/w_\sigma(n)\ll(M+2)2^{-M}n^{\sigma-1}$ stays bounded in $n$.
\end{proof}

Under a uniform stability hypothesis the truncation converges to the solution.

\begin{theorem}\label{thm:ortho_conditional_book}
Fix $\sigma\le1$ and assume the family $(R_M)$ is uniformly stable on
$\mathcal{X}_\sigma$. Then $\|A-A^{(M)}\|_\sigma\to0$ as $M\to\infty$, and in
particular $A_n=\mathcal{O}(n^{-\sigma})$ along the purely discrete route.
\end{theorem}

\begin{proof}
Combine Propositions~\ref{prop:error_abstract_book} and~\ref{prop:consistency_book}.
The summability $\sum_k|A_k|<\infty$ follows from Theorem~\ref{thm:ortho_main} of
Chapter~\ref{chap:ortho}: choose $\eps<\alpha_1-1$ in its bound
$|A_k|\ll_\eps k^{-\alpha_1+\eps}$, so the exponent is strictly larger than one.
The weaker coefficient estimate $c_n=\mathcal O(n^{-3/2})$ alone would give only a
$\mathcal O(n^{-1/2})$ tail and is not used here. Since
$A^{(M)}\in\mathcal{X}_\sigma$ for each $M$ by Theorem~\ref{thm:poly_transfer_book}, the limit
$A$ lies in $\mathcal{X}_\sigma$ as well, subject to the stated uniform stability hypothesis.
\end{proof}

\begin{proofstatus}{The implication is proved under the stated uniform-stability hypothesis.
That hypothesis, \eqref{eq:uniform_stability_book}, is open for every $\sigma$ and is the only
unproved input of the theorem. The consistency defect of
Proposition~\ref{prop:consistency_book} decays like $n^{-1}$, so the chain
stability plus consistency cannot reach $\mathcal{X}_\sigma$ beyond
$\sigma=1$. This is not a second gap in the stated implication but the reason for its restriction
to $\sigma\le1$. Reaching the spectral exponent $\alpha_1\approx1.3465$ would require a weighted consistency estimate,
a defect of order $n^{-\alpha_1}$ at least, which the plain tail bound does
not give. The numerical values below are evidence only. For the exact kernel
the weighted sums of \eqref{eq:uniform_stability_book} stay below $1.67$ at
$\sigma=1$ and below $4.1$ at $\sigma=\alpha_1$ up to $n=240$ with no visible
growth, and the truncations $M=1,3$ sit below the exact kernel, so uniform
stability is plausible. Neither the theorem nor either open extension is invoked in a later proof.}
\end{proofstatus}

\section{The RAF program at the discrete level}\label{sec:RAF_program}

The exact identities of this chapter organize into a program.

\begin{enumerate}
\item Start with an exact triangular kernel $G(N,k)$ on the integer lattice.
\item Apply the exact Abel transform to produce the discrete Volterra kernel
$K(N,k)=G(N,k)-G(N,k+1)$.
\item Construct the exact discrete resolvent $R(N,k)$ by the finite Neumann
series.
\item Define the finite arithmetic Mellin probes $\mathcal{G}_N(z)$,
$\mathcal{K}_N(z)$, $\mathcal{R}_N(z)$.
\item Prove convergence of the probes, locally uniformly in $\Re z<0$.
\item Identify the limits $G^*$, $K^*$, $R^*$ and prove the limiting resolvent
identity $R^*=K^*/(1+K^*)$.
\item Continue $G^*$ analytically and study its zeros and poles as spectral
objects.
\item Prove a spectral transfer theorem, the zeros of $G^*$ governing the
asymptotic decay of $A(N)$.
\item Refine toward a trace formula, the discrete resolvent encoding not only
the spectral exponents but their amplitudes.
\end{enumerate}

Steps 1 to 4 are the exact propositions of
Sections~\ref{sec:disc_primitive} to~\ref{sec:disc_probes}. Step 5 is settled
for $\mathcal{G}_N$ and $\mathcal{K}_N$ by
Theorem~\ref{thm:probe_convergence}, for every bounded Riemann integrable
profile and every kernel within $\mathcal{O}(1/N)$ of one, which covers both
the Ingham and the orthorecursive kernels. Within step 6 the limits $G^*$ and
$K^*$ are identified, while the convergence of $\mathcal{R}_N$ and the
limiting resolvent identity remain open, with the numerical support of
Remark~\ref{rem:probe_numerics}. Steps 7 and 8 are achieved for the
orthorecursive kernel in \cite{CloitreOrtho} and
Chapter~\ref{chap:ortho}, for polynomial kernels by
Theorem~\ref{thm:poly_transfer_book}, and for the Ingham kernel they are the
content of Chapter~\ref{chap:equivalence}, conditioned on the Riemann hypothesis
exactly as the theory predicts. The continuous resolvent level of Step 9 is established for the
orthorecursive kernel in Chapter~\ref{chap:trace}, Theorem~\ref{thm:w04d-infinite-trace}, while its
discrete pointwise level stays open, the reason for keeping the exact discrete structure alive rather
than passing at once to the continuous limit.

The class of kernels to which the program should apply can be circumscribed,
tentatively, as follows. The conditions below are those of the resolvent, and the name records
that. They are not the definition of a regular arithmetic function\index[terms]{regular arithmetic function}, which is
Definition~\ref{def:reg_index}.

\begin{definition}\label{def:resolvent_regular}
A triangular kernel $G(N,k)$ is resolvent regular if
\begin{enumerate}[label=(R\arabic*)]
\item there exists a bounded Riemann integrable $g:(0,1]\to\C$ with
$\sup_{1\le k\le N}|G(N,k)-g(k/N)|=\mathcal{O}(1/N)$,
\item the Volterra kernel $K(N,k)$ satisfies uniform summability estimates
controlling $\mathcal{K}_N(z)$ and $\mathcal{R}_N(z)$ in half planes
$\Re z\le-\delta$,
\item the resolvent probes $\mathcal{R}_N(z)$ converge locally uniformly in
$\Re z<0$ and the limiting identity \eqref{eq:Rstar_limit} holds.
\end{enumerate}
\end{definition}

Under (R1), the convergence of $\mathcal{G}_N$ and $\mathcal{K}_N$ is supplied
by Theorem~\ref{thm:probe_convergence}, so the substantive requirements are
(R2), the summability controlling the resolvent probes, and (R3), the
convergence of $\mathcal{R}_N$ together with the limiting identity. The Ingham
kernel and the orthorecursive kernel both satisfy (R1), and
Remark~\ref{rem:probe_numerics} supports (R3) for the orthorecursive kernel at
finite level.

\begin{proofstatus}{The convergence of the primitive and Volterra probes in
$\Re z<0$ is proved in this chapter, Theorem~\ref{thm:probe_convergence}, for
every kernel satisfying (R1). The convergence of the resolvent probes
$\mathcal{R}_N$ and the limiting identity \eqref{eq:Rstar_limit} are open in
general. For the orthorecursive kernel they hold numerically within
$2\cdot10^{-4}$ at $N=300$, and for the Ingham kernel the M\"obius structure
of the resolvent developed in Chapter~\ref{chap:volterra} is the expected
route. The general criterion (R2) isolating the resolvent regular class is part of
the next stage of the program. Thus resolvent regularity is a programmatic class at this stage:
no theorem in this volume assumes Definition~\ref{def:resolvent_regular}, and no kernel is claimed
here to satisfy (R1)--(R3) in full.}
\end{proofstatus}
That is where the second part stops. The operator has been inverted twice, on the half line for
the Ingham kernel and on the integer lattice for every kernel, and in both cases the inversion is
exact before any limit is taken. The transform is no longer given but produced, as the limit of
finite probes, for every bounded Riemann integrable profile. What the part does not deliver is
the matching statement for the resolvent. Their convergence is open in general, and the class
that would make them converge is named here and assumed by no theorem of this volume.

\part{Why the index is arithmetic and regular}
\label{part:three}

\rafepigraph{A l'usage, on se fait des bouts de dictionnaire, qui permettent de passer assez souvent d'une colonne à la colonne voisine.}{With practice one puts together scraps of a dictionary, which allow one to pass often enough from one column to the neighbouring column.}{André Weil, \emph{De la métaphysique aux mathématiques} (1960)~\cite{WeilMetaphysique1960}}

The two preceding parts computed an index and inverted an operator, and in both the index
followed the transform. This part establishes the two things the title names. The index is an
arithmetic quantity before it is anything else, and the word regular is earned rather than
borrowed.
Chapter~\ref{chap:diophantine} proves that the index is arithmetic before it is analytic, by
exhibiting two kernels where the two readings come apart, and that separation is what gives the
equivalence of Theorem~\ref{thm:tauberian_rh} its weight. A reader who has met only kernels where
the two agree may take the regularity index for the first zero of a transform under another name.
The agreement at the Ingham kernel is a fact about that kernel and not a consequence of the
definitions, and it is where the difficulty of the hypothesis sits.
Chapter~\ref{chap:regularity} accounts for the other word, widening the forcing to
$n^{-\beta}L(n)$ with $L$ slowly varying, placing the class beside Karamata regular variation, and
reading the defining equation as the Mercerian problem it is.

\chapter{The arithmetic nature of the index}
\label{chap:diophantine}

The regularity index was defined on the arithmetic side of the theory, from the
decay of the partial sums forced through the defining equation. The analytic index
was defined on the other side, from the first zero of the arithmetic Mellin
transform. Wherever both are available and both are proved, they agree, and that
agreement carries most of the examples of this volume. It is a theorem in each case
and never a law, and this chapter establishes that the two are distinct quantities.

Two kernels show that the two sides are distinct quantities. The rational kernel of
\S\ref{sec:rational_kernel} has a transform without zeros, so it carries no analytic
index at all, and its regularity index rests on the arithmetic alone. The broken
harmonic kernel\index[terms]{broken harmonic function} at $\sqrt2$ is the sharper case. Its transform has infinitely many
zeros, all of them on the vertical line $\Re z=1$, so the analytic index exists, is
unambiguous, and equals $1$. The exponent $\tfrac12$ is nevertheless not transparent,
and the transparency frontier\index[terms]{transparency frontier} is proved to lie at or below $\tfrac12$. The analytic
data is present and it is exact, and it is not the arithmetic answer. The two
readings are therefore not two names for one number.

What decides between the two situations is where the jump set of the profile meets
the evaluation grid. The index is read on the values $g(k/n)$, so it depends on the
arithmetic of that meeting and not on the analytic structure of the kernel alone. For the broken harmonic kernel that meeting is governed by the relation
$(\sqrt2)^{2}=2$, which pairs the layers of the kernel along the powers of two without ever
aligning them with the integer grid.
The Ingham kernel\index[terms]{Ingham kernel} jumps at the points $1/m$, which is the arithmetic of the grid
itself, and the discrepancy\index[terms]{discrepancy} the jumps create is resolved by M\"obius inversion\index[terms]{M\"obius inversion}. The
kernel broken on the lattice of $\sqrt2$ jumps at points the grid never reaches, and
the algebraic relation $(\sqrt2)^2=2$ pairs the layers of the kernel two by two along
the powers of $2$ without ever aligning them with the grid.

The Riemann hypothesis is a statement about the zeros of $\zeta$ and about the
regularity index of the Ingham kernel at once, and the equivalence of
Theorem~\ref{thm:tauberian_rh} carries both readings. What licenses the passage
between them is not analytic. It is the commensurability recorded in
\S\ref{sec:commensurable}, which places the two indices under one analytic object,
the zero set of $\zeta$. The kernel at $\sqrt2$ removes that compatibility and leaves
the analysis untouched, and the zeros of its transform then govern nothing on the
arithmetic side. The mechanism that makes the analytic reading of the hypothesis
available is itself arithmetic.

The three configurations the chapter records follow from this. For the continuous
classes and examples where it is proved, the two indices agree, as the continuous
chapter that follows this one sets out. For kernels whose transform is constant no
zero is available, and the homogeneous mode and the Green control take the place of
the transform. For kernels broken on an irrational lattice the two separate, at the
level of the transparency frontier, the case $\lambda=\sqrt2$ being the prototype and
its complete index remaining open. A distinct exclusion principle, the anti-HLR\index[terms]{anti-HLR statement}
statement of \cite{Cloitre2016}, removes kernels synthesized from $L$-functions
without Euler product\index[terms]{Euler product}, and the Davenport-Heilbronn\index[terms]{Davenport--Heilbronn function} function illustrates both its open
part and its provable part. The Diophantine input used along the way is classical,
the metric theory of continued fractions of Khinchin\index[names]{Khinchin, A. Ya.} \cite{Khinchin1964}
and the approximation theory collected by Schmidt\index[names]{Schmidt, W. M.} \cite{Schmidt1980}.

\section{Commensurable jumps}
\label{sec:commensurable}

The Ingham function $\Phi(x)=x\lfloor1/x\rfloor$ jumps exactly at the
rational points $x=1/m$. The evaluation grid and the jump set then live
in the same arithmetic, and the discrepancy created by the jumps is
resolved exactly by M\"obius inversion and Dirichlet convolution, as
carried out in Chapter~\ref{chap:equivalence} and again, on the continuous
half line, in Chapter~\ref{chap:volterra}.

In this commensurable situation both indices are tied to the same
object. The zeros of $\Phi^*(z)=\frac{z}{z-1}\zeta(1-z)$ are the points
$1-\rho$ with $\rho$ a nontrivial zero\index[terms]{nontrivial zero} of $\zeta$, so
\[
\eta(\Phi)=1-\sup_\rho\Re\rho\le\tfrac12 ,
\]
the inequality because zeros on the critical line\index[terms]{critical line} exist, by Hardy's\index[names]{Hardy, G. H.}
theorem \cite[Chapter~X]{Titchmarsh1986}. On the arithmetic side
$\tau(\Phi)\le\tfrac12$ as well, by Corollary~\ref{cor:phi_half}. Its input is the
existence of a simple critical zero not shared with $C_\beta$, supplied qualitatively by
Conrey\index[names]{Conrey, J. B.}'s positive-proportion theorem, and no zero-free region is used.
In the other direction, if $\Phi$ is a function of good variation\index[terms]{function of good variation} of
index $\alpha(\Phi)\in(0,\tfrac12)$, then Corollary~\ref{cor:zerofree} read with
$\sigma=1-\alpha(\Phi)$ places every zero of $\zeta$ in $\Re s\le1-\alpha(\Phi)$, so that
\[
\alpha(\Phi)\le\eta(\Phi)\le\tfrac12 .
\]
The opposite inequality would follow from a zero-free region\index[terms]{zero-free region}
$\Re s>\sigma$ with $\sigma<1$, and no such region is known. The equality
$\alpha(\Phi)=\eta(\Phi)=\tfrac12$ is equivalent to the
Riemann hypothesis, on both sides at once. The commensurable case does
not decide the equality. It guarantees that a single analytic object,
the zero set of $\zeta$, governs both indices.

\section{Incommensurable jumps: the kernel $g_{\sqrt2}$}
\label{sec:sqrt2}

The separation between the analytic index $\eta(g)$ and the arithmetic
response is proved at the level of the transparency frontier\index[terms]{transparency frontier} when the kernel jumps
on the lattice of $\sqrt2$, and expected in general. This section states what is
established for that kernel and what is not. Every statement below is proved in the
research dossier~\ref{app:dossier_sqrt2}, which also carries the exact block
recurrences, the finite certificates, and the two remaining global obstructions. Each
statement records in its own terms what it leaves open.

The arithmetic side of the comparison is a frontier rather than an index. It reads the
largest range of forcings that pass through the kernel without distortion, and it is
defined without assuming that a regularity index exists.

\begin{definition}[Transparency frontier]
\label{def:w10-frontier}
For a discrete kernel $G$ for which the triangular inverse is defined, let
\begin{equation}
\label{eq:w10-frontier}
\tau_{\mathrm{tr}}(G)=
\sup\bigl\{c\in\R:\text{ every }\beta<c\text{ is transparent}\bigr\}.
\end{equation}
Transparency at $\beta$ means
\[
A(N)=\Xi_G(\beta)N^{-\beta}+o(N^{-\beta}),
\qquad
\Xi_G(\beta)=\frac1{G^*(\beta)},
\]
where the reciprocal is read meromorphically. The frontier is more demanding than the
regularity index of Definition~\ref{def:reg_index}, which asks only that some constant
exist, so that
\[
\tau_{\mathrm{tr}}(G)\le\alpha(G)
\]
for every kernel, with equality wherever the identification
$\Xi_G=1/G^{*}$ is established. The frontier uses the transparency branch alone and
carries no absorption clause, so it does not by itself make $G$ a RAF.
For $G_q$, the set in \eqref{eq:w10-frontier} is nonempty because every
$\beta<0$ is transparent.
\end{definition}

The analytic side is read off the transform, and for the self similar broken kernels it
is the same number for every scale. For $\lambda>1$ let
\begin{equation}\label{eq:g_lambda_def}
g_\lambda(x):=x\,\lambda^{\lfloor-\log_\lambda x\rfloor}
\qquad(0<x\le1),
\end{equation}
the self similar broken harmonic function\index[terms]{broken harmonic function} with jump points
$x=\lambda^{-i}$ and slopes $\lambda^{i-1}$ on
$(\lambda^{-i},\lambda^{-i+1}]$. Its arithmetic Mellin transform is computed directly
from the definition. On the piece $(\lambda^{-i},\lambda^{-i+1}]$ the kernel is
$\lambda^{i-1}x$, so for $\Re z<0$ the defining integral splits over the layers, and the
geometric series $\sum_{i\ge1}\lambda^{iz}=\lambda^{z}/(1-\lambda^{z})$ sums the pieces
into a closed form which then continues meromorphically.

\begin{proposition}\label{thm:w10-transform}
For $\operatorname{Re}z<0$,
\begin{equation}
\label{eq:w10-transform}
\boxed{
g_\lambda^*(z)=
\frac{z}{z-1}\frac{\lambda^{z-1}-1}{\lambda^z-1}.}
\end{equation}
The right side is the meromorphic continuation to $\C$. The points
$z=0$ and $z=1$ are removable, with
\begin{equation}
\label{eq:w10-removable-values}
g_\lambda^*(0)=\frac{1-\lambda^{-1}}{\log\lambda},
\qquad
g_\lambda^*(1)=\frac{\log\lambda}{\lambda-1}.
\end{equation}
The nonremovable poles are the simple points
\[
z=\frac{2\pi i m}{\log\lambda}
\quad(m\in\Z\setminus\{0\}),
\]
and the zeros are the simple points
\[
z=1+\frac{2\pi i k}{\log\lambda}
\quad(k\in\Z\setminus\{0\}).
\]
There are no further cancellations. Under the analytic-index convention of
the monograph,
\begin{equation}
\label{eq:w10-eta}
\eta(g_\lambda)=1
\quad(\lambda>1),
\qquad
\eta(g_{\sqrt2})=1.
\end{equation}
The transform does not record the least positive integer $k$ for which
$\lambda^k$ is integral. Formula \eqref{eq:w10-eta} does not determine an
arithmetic transition.
\end{proposition}

The noncancelled zeros lie on the vertical line $\Re z=1$ and the point $z=1$ is
removable, so the analytic index is $1$ for every scale. The transform records nothing
about the arithmetic of $\lambda$, and the content of this section is that the
discretization does.

For $\lambda=\sqrt2$ the algebraic coincidence $\lambda^{2}=2$ pairs the layers of the
kernel two by two when $g(k/n)$ is evaluated on the integers. The notation that follows
carries that pairing.

\subsection*{Notation for the kernel at $\sqrt2$}

The profile \eqref{eq:g_lambda_def} induces the discrete kernel
\begin{equation}
\label{eq:w10-profile}
G_\lambda(n,k)=g_\lambda(k/n)
\qquad(\lambda>1).
\end{equation}
The value of $g_\lambda$ at zero may be chosen arbitrarily. It has no effect on the
Riemann integral below. For a real forcing exponent $\beta$, let
$(a_n^{(\beta)})_{n\geq1}$ be the unique triangular solution of
\begin{equation}
\label{eq:w10-forcing-equation}
\sum_{k=1}^{n}a_k^{(\beta)}g_\lambda(k/n)=n^{-\beta}.
\end{equation}
Write
\[
A_\beta(N)=\sum_{n\leq N}a_n^{(\beta)},
\qquad
g_\lambda^*(z)=-z\int_0^1g_\lambda(t)t^{-z-1}\,dt.
\]

Several forcing exponents are compared in what follows, so the subscript is kept on $A_\beta$ and on the sequences attached to it. Elsewhere in the volume the exponent is fixed by its context and the plain $a_n$ and $A(n)$ of Definition~\ref{def:reg_index_fgv} are used.

The arithmetic part of this section concerns
\[
q=\sqrt2,
\qquad
g=g_q,
\qquad
p=1-\beta,
\qquad
F_p(n)=n^p\ (n\geq1),
\qquad F_p(0)=0.
\]
Put
\[
\delta_p(n)=F_p(n)-F_p(n-1),
\quad b_n=na_n^{(\beta)},
\quad C(N)=\sum_{n\leq N}b_n,
\quad s(n)=\left\lfloor\frac n q\right\rfloor,
\quad B(m)=\lceil qm\rceil.
\]
For any sequence $W$, write $\Delta W(n)=W(n)-W(n-1)$.
The block variable is
\begin{equation}
\label{eq:w10-Q-definition}
Q_0=0,
\qquad
Q_n=n\sum_{j=0}^{v_2(n)}a_{n/2^j}^{(\beta)}
\quad(n\geq1).
\end{equation}
All occurrences of $Q$, $A$, $b$, and $C$ in the arithmetic results use
the same fixed value of $\beta$.

For $N\geq2$, define
\[
H=\lfloor N/2\rfloor,
\quad x_0=N,
\quad x_{j+1}=s(x_j),
\quad
\eta_n=\lfloor n/2\rfloor-s^2(n),
\quad
\mathcal Z_p(n)=q\eta_nQ_{\lfloor n/2\rfloor}.
\]
Here $\eta_n\in\{0,1\}$ and all sums along $(x_j)$ are finite. The
decisions may be written without nested floors as
\begin{equation}
\label{eq:w10-eta-decisions}
\eta_{2m}=1\quad(m\geq1),
\qquad
\eta_{2m+1}=\mathbf1_{\{\{mq\}<1-1/q\}}
\quad(m\geq1),
\qquad
\eta_1=0.
\end{equation}
Set
\[
R_p(N)=\sum_{j\geq0}\mathcal Z_p(x_j),
\qquad R_p(0)=0.
\]

The centered defect $G_\beta(N)$ used below is the quantity \eqref{eq:w10-G} of
Theorem~\ref{thm:w10-Abel} in the research dossier, the part of $A(N)$ that
carries the orbit trace $R_p$.

\subsection*{The exact reduction}

The pairing of layers turns the forcing equation into a three-term lag identity on the
block variable, and a rescaling by $n^{1/2}$ renders the principal descent branch
marginal. This identity is the mechanism behind every statement that follows.

\begin{lemma}
\label{lem:w10-reduction}
For every real $\beta$ and every $n\geq1$, put
\begin{equation}
\label{eq:w10-chi}
\chi_n=\mathbf1_{\{n=B(s(n))\}}.
\end{equation}
Then
\begin{equation}
\label{eq:w10-recurrence}
\boxed{
Q_n=\delta_p(n)+q\mathbf1_{2\mid n}Q_{n/2}
-(q-1)\chi_nQ_{s(n)}.}
\end{equation}
The decisions are exact, since
\begin{equation}
\label{eq:w10-integer-test}
\chi_n=1
\quad\Longleftrightarrow\quad
(n-1)^2<2s(n)^2<n^2.
\end{equation}
Moreover
\begin{equation}
\label{eq:w10-bQ}
\boxed{b_n=Q_n-2\mathbf1_{2\mid n}Q_{n/2}}
\end{equation}
and, for $N\geq1$,
\begin{equation}
\label{eq:w10-half-block}
\boxed{A(N)=
\sum_{\lfloor N/2\rfloor<n\leq N}\frac{Q_n}{n}.}
\end{equation}
If $S(N)=\sum_{n\leq N}Q_n$ and $S(0)=0$, then
\begin{equation}
\label{eq:w10-S-recurrence}
S(N)=F_p(N)+qS(\lfloor N/2\rfloor)-(q-1)S(s(N)),
\qquad
C(N)=S(N)-2S(\lfloor N/2\rfloor).
\end{equation}
All domains and initial values are included in these formulas.
\end{lemma}

What the reduction yields by analytic means alone is the negative range of forcings and
a smoothed control at every order.

\begin{lemma}\label{thm:w10-acquired-analysis}
Every $\beta<0$ is transparent for $G_q$ and
\begin{equation}
\label{eq:w10-negative-transparency}
A(N)=\frac{N^{-\beta}}{g^*(\beta)}+o(N^{-\beta}).
\end{equation}
For every $\beta\geq0$ and $\eps>0$,
\begin{equation}
\label{eq:w10-hilbert}
\sum_{n\geq1}\frac{|Q_n|^2}{n^{2+2\eps}}<\infty.
\end{equation}
If
\[
\mathcal R_{k,\beta}(x)=
\sum_{n\leq x}a_n^{(\beta)}(1-n/x)^k,
\]
then for every integer $k\geq2$ and every $\eps>0$,
\begin{align}
\mathcal R_{k,\beta}(x)
&=\frac{k!}{\prod_{j=1}^k(j-\beta)}
\frac{x^{-\beta}}{g^*(\beta)}
+O_{\beta,k,\eps}(x^{-1/2+\eps})
&& (0\leq\beta<1/2),
\label{eq:w10-Riesz-low}\\
\mathcal R_{k,\beta}(x)
&=O_{\beta,k,\eps}(x^{-1/2+\eps})
&& (\beta\geq1/2).
\label{eq:w10-Riesz-high}
\end{align}
These are smoothed estimates. They do not give the corresponding raw
estimates.
\end{lemma}

\subsection*{The critical exponent and the separation}

The marginal branch sits at $\beta=1/2$, and the reduction is sharp enough to decide
that exponent.

\begin{theorem}
\label{thm:w10-nontransparency}
At $\beta=1/2$, define
\[
V_r=q^{-r}Q_{2^r}.
\]
Then
\begin{equation}
\label{eq:w10-positive-L}
\boxed{V_r\longrightarrow L,
\qquad L>\frac{1293}{5000}>\frac14.}
\end{equation}
Moreover
\begin{equation}
\label{eq:w10-critical-A-limit}
\sqrt{2^r}A_{1/2}(2^r)
\longrightarrow q^{1/2}-(q-1)L.
\end{equation}
Since $1/g^*(1/2)=q^{1/2}$, the exponent $\beta=1/2$ is not
transparent. In particular,
\[
\tau_{\mathrm{tr}}(G_q)\leq\frac12.
\]
\end{theorem}

The two indices then separate, with explicit numerical bounds on both sides.

\begin{theorem}
\label{thm:w10-separation}
The broken square-root kernel satisfies
\begin{equation}
\label{eq:w10-central-separation}
\boxed{
1-\log_2(177/100)\leq\tau_{\mathrm{tr}}(G_{\sqrt2})
\leq\frac12<1=\eta(g_{\sqrt2}).}
\end{equation}
This compares the transparency frontier with the analytic index. It does
not assume that a complete RAF index exists. The transform has infinitely
many zeros on $\operatorname{Re}z=1$, but it does not determine the
arithmetic frontier.
\end{theorem}

The two-sided bound is the separation announced above. Its upper half is the critical
exponent, its lower half comes from the finite gauge certificate stated further down.
The same critical limit also decides the Hardy-Littlewood-Ramanujan\index[terms]{Hardy--Littlewood--Ramanujan criterion}\index[names]{Ramanujan, S.} criterion.

\begin{theorem}
\label{thm:w10-HLR}
At the critical forcing $\beta=1/2$, the exact identity
\eqref{eq:w10-bQ} and the limit \eqref{eq:w10-positive-L} give
\begin{equation}
\label{eq:w10-HLR-limit}
\boxed{
\frac{b_{2^r}}{2^{r/2}}
=V_r-qV_{r-1}
\longrightarrow(1-q)L\ne0.}
\end{equation}
Consequently
\begin{equation}
\label{eq:w10-HLR-asymptotic}
|2^ra_{2^r}^{(1/2)}|
\sim (q-1)L\,2^{r/2}.
\end{equation}
For every $0<\eps<1/2$,
\begin{equation}
\label{eq:w10-HLR-ratio}
\frac{|2^ra_{2^r}^{(1/2)}|}{(2^r)^\eps}
\longrightarrow+\infty.
\end{equation}
The HLR criterion\index[terms]{Hardy--Littlewood--Ramanujan criterion} requires
$na_n^{(\beta)}=o(n^\eps)$ for every $\eps>0$ and every
$\beta\geq0$. Hence $g_{\sqrt2}$ is not HLR and is not strong HLR.
This failure does not exclude the RAF property. It only prevents a direct
transfer of the HLR inversion used for the Ingham kernel.
\end{theorem}

\subsection*{What a complete index would require}

A frontier is not an index. Turning it into one asks for two raw estimates, one below
the critical exponent and one above it.

\begin{definition}[LOW and ABS]
\label{def:w10-low-abs}
The remaining raw estimates are
\begin{align*}
\mathsf{LOW}&\quad
G_\beta(N)=o(N^{-\beta})
&& (0\leq\beta<1/2),
\\
\mathsf{ABS}&\quad
G_\beta(N)=O_{\beta,\eps}(N^{-1/2+\eps})
&& (\beta\geq1/2,\ \eps>0).
\end{align*}
They are respectively equivalent to
\begin{align*}
R_p(N)&=o(N^p)&& (1/2<p\leq1),
\\
R_p(N)&=O_{\beta,\eps}(N^{1/2+\eps})
&& (p\leq1/2).
\end{align*}
The statements are global in $N$, hence global in the odd cores of $N$.
The letters ABS name the upper raw branch. They do not assert that
$R_p$ is defined by an absolute sum.
\end{definition}

A proper initial part of LOW is proved by a finite certificate. The transfer matrices
$\mathcal M_t$ and the orthant conjugations $\Sigma_t$ are those of the mixed binary
closure, Theorem~\ref{thm:w10-binary} and Lemma~\ref{lem:w10-mixed-certificate} of the
research dossier, and the gauge vectors are tabulated there.

\begin{proposition}\label{thm:w10-partial-LOW}
Let
\[
\mu=\frac{177}{100},
\qquad
\theta=\log_2\mu=0.8237493603\ldots.
\]
For each of the twenty-eight mixed cells there is a positive rational
vector $h_t\in(2000^{-1}\Z)^{14}$. The vectors are listed in
Table~\ref{tab:w10-mixed-gauges}. They satisfy
\[
\frac1{2000}\leq(h_t)_i\leq\frac{2001}{2000},
\qquad
(\Sigma_t\mathcal M_t\Sigma_s)h_s\leq\mu h_t
\]
for every coordinate and every one of the $136$ edges $s\to t$.
The smallest algebraic slack is
\begin{equation}
\label{eq:w10-low-min-slack}
\frac{283023}{200000}-\frac{2001}{2000}q>0,
\qquad
\left(\frac{283023}{200000}\right)^2
-2\left(\frac{2001}{2000}\right)^2
=\frac{21998529}{40000000000}>0.
\end{equation}
Consequently, uniformly for $p\leq1$,
\begin{equation}
\label{eq:w10-theta-Q}
Q_n=O(n^\theta),
\qquad
R_p(N)=O(N^\theta).
\end{equation}
It follows that
\begin{equation}
\label{eq:w10-partial-transparency}
A(N)=\frac{N^{-\beta}}{g^*(\beta)}
+O_\beta(N^{\theta-1})
\quad(0\leq\beta<1-\theta).
\end{equation}
At $\beta=0$, the smaller floor error $O(\log N/N)$ is absorbed. In
particular,
\begin{equation}
\label{eq:w10-frontier-two-sided}
1-\log_2(177/100)
\leq\tau_{\mathrm{tr}}(G_q)\leq\frac12.
\end{equation}
This proves a proper initial part of LOW, not the
full statement.
\end{proposition}

Granting both estimates closes the case, and the implication itself is proved.

\begin{conditionaltheorem}[Exact RAF closure]
\label{cthm:w10-conditional-RAF}
If LOW and ABS both hold, then $G_{\sqrt2}$ is a RAF with
\[
\alpha(G_{\sqrt2})=\frac12.
\]
Indeed, \eqref{eq:w10-negative-transparency} and LOW give transparency for
every $\beta<1/2$. ABS gives the absorbed branch for every
$\beta\geq1/2$. Theorem~\ref{thm:w10-nontransparency} gives sharpness at
the threshold. This implication is proved. Its two hypotheses remain open.
\end{conditionaltheorem}

\begin{openproblem}[Global closure]
\label{op:w10-global}
\sloppy
LOW and ABS remain open. Equivalently, by
Lemma~\ref{lem:w10-trace-equivalence}, write $Q_n^{(p)}$ for the block variable
$Q_n$ attached to the indicated value of $p=1-\beta$. The missing estimates are
\begin{align*}
Q_n^{(p)}&=o(n^p)&& (1/2<p\leq1),\\
Q_n^{(p)}&=O_{p,\eps}(n^{1/2+\eps})&& (p\leq1/2).
\end{align*}
The available Hilbert estimate\index[terms]{Hilbert estimate}\index[names]{Hilbert, D.}, shell-mass bounds, and orbit algebra have not
yielded the first estimate. The finite mixed gauges have multiplier
$177/100>q$, and Proposition~\ref{prop:w10-expanding-cycle} rules out the
full-state critical product estimate. Two concrete obstacles remain, an
observable estimate that removes the expanding direction and critical-scale
control of the initial germ uniformly as the odd core varies.
\end{openproblem}

\begin{conjecture}[Collapse of the index for $g_{\sqrt2}$]
\label{conj:sqrt2_collapse}
$\alpha(g_{\sqrt2})=\tfrac12$.
\end{conjecture}
\begin{proofstatus}{The analytic index $\eta(g_{\sqrt2})=1$ is proved above.
Transparency fails at the critical exponent $\beta=\tfrac12$ by
Theorem~\ref{thm:w10-nontransparency}, the transparency frontier is pinned to
$1-\log_2(177/100)\le\tau_{\mathrm{tr}}(G_{\sqrt2})\le\tfrac12$ by
Theorem~\ref{thm:w10-separation}, and the Hardy--Littlewood--Ramanujan criterion fails at the
critical forcing by Theorem~\ref{thm:w10-HLR}. What remains open is exactly the pair of
estimates LOW and ABS isolated in Open Problem~\ref{op:w10-global}. Together they give the
existence of the regularity index and the value $\tfrac12$ by Conditional
Theorem~\ref{cthm:w10-conditional-RAF}, and without them neither conclusion is proved. The numerical
profile mentioned below is evidence only. This conjecture is a terminal synthesis and is assumed
in no proof in this volume.}
\end{proofstatus}

Since $\eta(g_{\sqrt2})=1$, the statement asserts a strict separation
$\alpha(g_{\sqrt2})=\tfrac12<1=\eta(g_{\sqrt2})$. The continuous
transform carries no trace of the layer pairing, so the separation, now proved at the
level of the transparency frontier, is invisible to $g^*$ and belongs entirely to the
discretization. The full-state all-path Green estimate\index[terms]{Green estimate}\index[names]{Green, G.} once proposed for the family is
false, by the realized expanding cycle of Proposition~\ref{prop:w10-expanding-cycle} in
the research dossier, so the remaining task is an observable estimate that removes its
expanding direction, together with global summation over the odd cores. For $\sqrt2$ the
failure of transparency at $\tfrac12$ replaces the earlier Riesz indication, the
numerical profile points to equality, and general radical claims remain unproved.

\section{The anti-HLR statement and Davenport-Heilbronn}
\label{sec:anti_hlr}

A second exclusion principle operates independently of the jump lattice.
It concerns kernels synthesized from Dirichlet series\index[terms]{Dirichlet series} that satisfy a
functional equation\index[terms]{functional equation} of Riemann\index[names]{Riemann, B.} type but lack an Euler product, and it
explains why such kernels fall outside the HLR class. This is
the anti-HLR statement of \cite{Cloitre2016}, formulated for the
broken harmonic functions\index[terms]{broken harmonic function}.

\begin{definition}[{Broken harmonic function, \cite[\S~2.1]{Cloitre2016}}]
\label{def:BHF}
A bounded function $g$ on $(0,1]$ is a broken harmonic function\index[terms]{broken harmonic function} if there
exist a real sequence $1=u_1>u_2>\cdots$ with $u_n\to0$ and an
increasing sequence $v_n>0$ with $\sup_nu_nv_n<\infty$, such that
$g(x)=v_nx$ for $u_{n+1}<x\le u_n$. The Ingham function is the case
$u_i=1/i$, $v_i=i$, and $g_\lambda$ of \eqref{eq:g_lambda_def} is the
case $u_i=\lambda^{1-i}$, $v_i=\lambda^{i-1}$.
\end{definition}

\begin{remark}[The HLR criterion and its strong form]\label{rem:hlr_forms}
The HLR criterion of Definition~\ref{def:HLR}, following
\cite[\S~2.2]{Cloitre2016}, asks $na_n=o(n^{\eps})$ for every
$\eps>0$. Its strong form asks $na_n=\mathcal O(1)$, which the
Ingham function satisfies by Proposition~\ref{prop:hlr_ingham}. The
strong form implies the HLR criterion, so a violation of the HLR
criterion violates both. The statement below concerns the HLR criterion
in the sense of Definition~\ref{def:HLR}.
\end{remark}

\begin{conjecture}[{Anti-HLR, \cite[\S~3]{Cloitre2016}}]
\label{conj:anti_hlr}
Let $g$ be a broken harmonic function\index[terms]{broken harmonic function} for which the limit
$\lim_{x\to0}g(x)$ exists and is nonzero, and for which the product
$(1-z)g^*(z)$ satisfies a functional equation of Riemann\index[names]{Riemann, B.} type. If $g^*$ has a zero in
the half plane $\Re z<\tfrac12$, then $g$ does not satisfy the HLR
criterion.
\end{conjecture}
\begin{proofstatus}{The conjecture is quoted from \cite[\S~3]{Cloitre2016} and is not proved
here. Proposition~\ref{prop:anti_hlr_right} and Corollary~\ref{cor:DH_not_HLR} establish
unconditionally the failure of HLR for the Davenport--Heilbronn kernel by using zeros with
$\Re s>1$, and they do not use this conjecture. Its open content is the critical-strip range
$\tfrac12<\Re s\le1$, where it is used below only to predict failure of HLR, never as a hypothesis
of a theorem. Moreover Remark~\ref{rem:DH_numerics} records that the Davenport--Heilbronn kernel
satisfies Definition~\ref{def:BHF} only after the monotonicity of the slopes is relaxed. It is
therefore a guiding model for the conjecture, not a literal instance of its stated class.}
\end{proofstatus}

The contrapositive, recorded as Corollary~3.1 in \cite{Cloitre2016},
reads as follows. If such a $g$ satisfies the HLR criterion, then the
nontrivial zeros of $g^*$ lie on the critical line\index[terms]{critical line}, the zeros being
symmetric about that line by the functional equation. The statement
thereby turns the HLR criterion into a zero location principle for the
class of transforms it governs.

\subsection*{The Davenport-Heilbronn kernel}

The illustration proposed in \cite[\S~3.2]{Cloitre2016} starts from the
function of Davenport\index[names]{Davenport, H.} and Heilbronn\index[names]{Heilbronn, H.} \cite{DavenportHeilbronn1936}. The
two papers under that title treat two objects, and the first of them is the Epstein\index[names]{Epstein, P.} zeta
function\index[terms]{Epstein zeta function} of a positive binary quadratic form of class number greater than one, shown there to have
infinitely many zeros in the half plane of absolute convergence. The second is the periodic
series recalled below. Both carry a functional equation of Riemann type and neither carries an
Euler product\index[terms]{Euler product}, and it is that pairing, and not the periodicity, which the
anti-HLR statement addresses. Let
\[
\xi=\frac{-2+\sqrt{10-2\sqrt5}}{\sqrt5-1}=0.284079\ldots
\]
and let $h$ be the $5$-periodic sequence
$h(1),\ldots,h(5)=1,\ \xi,\ -\xi,\ -1,\ 0$. The Dirichlet series
$H(s)=\sum_{n\ge1}h(n)n^{-s}$ continues to an entire function satisfying
a functional equation of Riemann\index[names]{Riemann, B.} type, and it has infinitely many
nontrivial zeros off the critical line\index[terms]{critical line} \cite{DavenportHeilbronn1936}. It
also has infinitely many zeros in the half plane of absolute convergence
$\Re s>1$, whose real parts are studied in
\nm{Bombieri}{E.} and \nm{Ghosh}{A.}~\cite[Theorem~7]{BombieriGhosh2011}. Both features enter below, at two
different levels.

\begin{figure}[p]
\centering
\rafwithfig{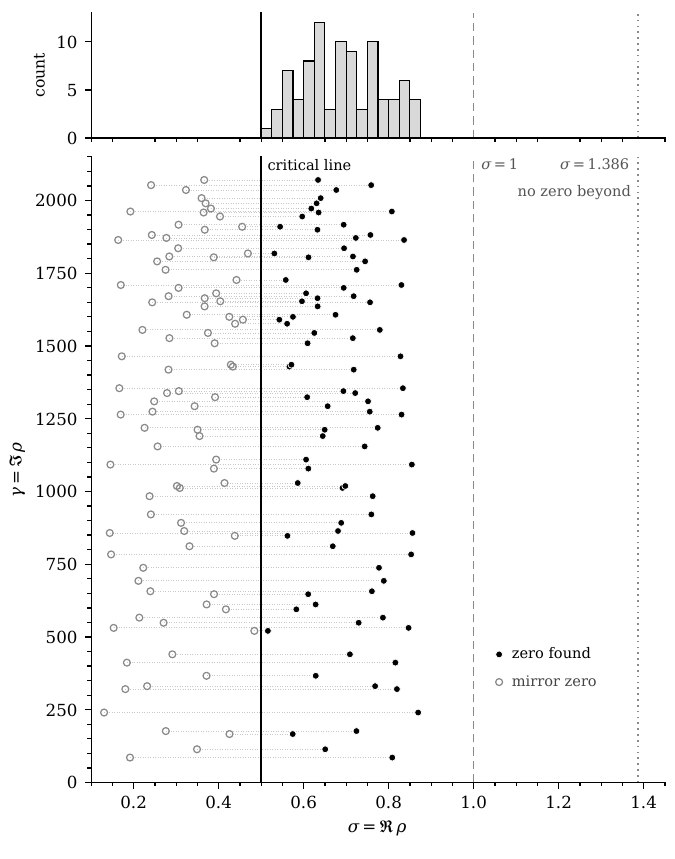}{\includegraphics[width=0.84\textwidth]}%
  {\fbox{\parbox[c][5.11in][c]{0.84\textwidth}{\centering\small
     The figure of the Davenport--Heilbronn zeros is missing. Place
     \texttt{dh\_zeros\_plane.pdf} next to the source file or in a subfolder named
     \texttt{figures}, \texttt{fig} or \texttt{img}. The box has the height of the figure it replaces, so the
     pagination of the volume is unaffected.}}}
\caption{The zeros $\rho=\sigma+i\gamma$ of $H$ with $\sigma>\tfrac12$ and
$0<\gamma<2100$, eighty-eight of them, each drawn as a filled disc. The open circle at
the same height is the mirror zero $1-\bar\rho$, supplied by the functional equation.
Each is accompanied by its mirror, so one hundred and seventy-six zeros lie off the
critical line\index[terms]{critical line} over this range. The strip carries $2147$ zeros in all there, so one zero in
twelve is off the line, a ratio that depends on the window and cannot be constant, the
total growing like $T\log T$ while the number off the line is at most of order $T$. The abscissas run from $0.5159$, at
$\gamma=520.94$, to $0.8695$, at $\gamma=240.40$, and the top panel is their
distribution. None reaches one over this range, although zeros with $\sigma>1$ exist in
infinite number \cite{DavenportHeilbronn1936} and their real parts are described in
\cite[Theorem~7]{BombieriGhosh2011}. No zero lies beyond $\sigma=1.3952$, since
$\sum_{n\ge2}|h(n)|\,n^{-\sigma}<1$ there and the constant term dominates, the sum
reaching one at $\sigma=1.39513\ldots$. The zeros
were located by the argument principle, the count of the strip being compared with the
sign changes of the Hardy\index[names]{Hardy, G. H.} function attached to $H$, and each
one was then refined until $|H(\rho)|$ reached the working precision.}
\label{fig:dh_zeros}
\end{figure}

Figure~\ref{fig:dh_zeros} records what the anomaly looks like over the range computed. The
zeros that leave the critical line\index[terms]{critical line} show no accumulation along a
second vertical there, their abscissas spreading over $\left[0.5159,\,0.8695\right]$ with no
visible preferred value, and they are sparse, one zero in twelve. The two features
just named sit at different scales. The zeros off the line are within reach of a moderate
computation and are the ones drawn. The zeros in the half plane of absolute convergence
are rarer, none occurs below height $2100$, and their existence is taken from
\cite{DavenportHeilbronn1936} rather than from the picture.

The associated kernel is the superposition
\[\Phi_H(x):=x\sum_{1\le k\le1/x}h(k)\Big\lfloor\frac1{kx}\Big\rfloor
=\sum_{k\ge1}\frac{h(k)}{k}\,\Phi(kx)
\qquad(0<x\le1),
\]
the second form because $\Phi(kx)=kx\lfloor1/(kx)\rfloor$ vanishes for
$kx>1$. The change of variables $u=kx$ in the defining integral of the
transform multiplies each term by $k^{-(1-z)}$, so for $\Re z<0$,
\begin{equation}\label{eq:PhiH_transform}
\Phi_H^*(z)=\Phi^*(z)\,H(1-z)
=\frac{z}{z-1}\,\zeta(1-z)\,H(1-z),
\end{equation}
with meromorphic continuation. Every zero of $H$ at $\rho$ with
$\Re\rho>\tfrac12$ produces a zero of $\Phi_H^*$ at $z=1-\rho$ with
$\Re z<\tfrac12$, so the hypothesis of
Conjecture~\ref{conj:anti_hlr} is met, and the statement predicts that
$\Phi_H$ is not HLR.

\begin{remark}[Normalization of the transform]\label{rem:transform_norm}
The transform printed in \cite[\S~3.2]{Cloitre2016} is
$\zeta(1-z)H(1-z)/(1-z)$, which differs from
\eqref{eq:PhiH_transform} by the factor $-z$ attached there to a variant
normalization. The form \eqref{eq:PhiH_transform} is the one consistent
with the normalization $g^*(0)=\lim_{x\to0}g(x)$ used throughout this
monograph, here $\Phi_H^*(0)=H(1)=0.9228\ldots$, and it was checked
against a direct numerical evaluation of the transform to $10^{-7}$ at
two test points. The zero sets relevant to the statement agree in both
normalizations.
\end{remark}

For zeros of $H$ to the right of the convergence line the prediction is
a theorem, by an argument that requires nothing beyond Dirichlet
convolution.

\begin{proposition}
\label{prop:anti_hlr_right}
Let $h$ be a bounded arithmetic function with $h(1)\ne0$, and put
$H(s)=\sum_{n\ge1}h(n)n^{-s}$. Suppose that $H(\rho_0)=0$ for some
$\rho_0$ with $\sigma_0:=\Re\rho_0>1$. Let $(a_n)$ be the solution of the
defining relation $A_{\Phi_H}(n)=1$ for all $n\ge1$, at $\beta=0$, and
set $b_n:=na_n$. Then
\[
\limsup_{n\to\infty}\,|b_n|\,n^{1-\sigma_0+\eps}=+\infty
\qquad\text{for every }\eps>0 .
\]
In particular $(b_n)$ is unbounded and $\Phi_H$ satisfies neither the
HLR criterion nor its strong form.
\end{proposition}

\begin{proof}
The kernel evaluation rearranges into a triple Dirichlet convolution.
Since $\lfloor n/(jk)\rfloor$ counts the integers $i$ with $ijk\le n$,
\[
n\,A_{\Phi_H}(n)
=\sum_{k\le n}k\,a_k\sum_{j\le n/k}h(j)\Big\lfloor\frac n{jk}\Big\rfloor
=\sum_{ijk\le n}b_k\,h(j)
=\sum_{m\le n}(b\star h\star\mathbf 1)(m),
\]
where $\mathbf 1$ is the constant function $1$ and $\star$ is Dirichlet
convolution. The defining relation at $\beta=0$ makes the left side
equal to $n$ for every $n$, so $(b\star h\star\mathbf 1)(m)=1$ for every
$m$, that is $b\star h\star\mathbf 1=\mathbf 1$. Convolving both sides
with $\mu$ and using $\mu\star\mathbf 1=e$ from
\S\ref{sec:master_equiv} gives
\[
b\star h=e ,
\]
the sequence $(b_n)$ is the Dirichlet inverse of $(h_n)$. Now suppose
that $b_n\ll n^{\sigma_0-1-\eps_0}$ for some $\eps_0>0$.
Then $B(s)=\sum_nb_nn^{-s}$ converges absolutely for
$\Re s>\sigma_0-\eps_0$, and $H(s)$ converges absolutely for
$\Re s>1$ since $h$ is bounded. On the half plane
$\Re s>\max(1,\sigma_0-\eps_0)$, which contains $\rho_0$ because
$\sigma_0>1$, the convolution identity gives $B(s)H(s)=1$. Evaluating at
$s=\rho_0$ yields $1=B(\rho_0)\,H(\rho_0)=0$, a contradiction. Hence
$b_n\ll n^{\sigma_0-1-\eps}$ fails for every $\eps>0$,
which is the stated limsup. Since $\sigma_0-1>0$ the sequence is
unbounded, violating the strong HLR criterion, and
$b_n=o(n^{\eps})$ fails for $\eps<\sigma_0-1$, violating
the HLR criterion itself.
\end{proof}

Applied to the Davenport-Heilbronn coefficients, whose Dirichlet series has zeros to the right
of the line of absolute convergence, this settles the case at hand.

\begin{corollary}\label{cor:DH_not_HLR}
The Davenport-Heilbronn kernel\index[terms]{Davenport--Heilbronn function} $\Phi_H$ does not satisfy the HLR
criterion.
\end{corollary}
\begin{proof}
The function $H$ has zeros with $\Re s>1$
 by \nm{Bombieri}{E.} and \nm{Ghosh}{A.}~\cite[Theorem~7]{BombieriGhosh2011}, $h$ is bounded, and $h(1)=1$, so
Proposition~\ref{prop:anti_hlr_right} applies.
\end{proof}

The proposition settles the region right of $1$ unconditionally. The
content of Conjecture~\ref{conj:anti_hlr} is the critical strip\index[terms]{critical strip}. A zero of
$H$ at $\rho$ with $\tfrac12<\Re\rho\le1$ forces, by the same
convolution identity, the abscissa of convergence of $B$ up to
$\Re\rho$, but $\Re\rho-1\le0$, so the pole only excludes decay of
$(b_n)$ faster than $n^{\Re\rho-1}$ and does not by itself contradict
boundedness. The statement asserts that the strip zeros nevertheless
break the criterion, through the finer growth of $(b_n)$ discussed in
\cite{Cloitre2016}. The Ingham operator thus excludes non Eulerian
kernels at two levels, the provable one for zeros beyond the convergence
line and the conjectural one inside the critical strip, and the HLR
criterion is in this sense a quantitative expression of unique
factorization.

\begin{remark}[Numerical behavior of the Davenport-Heilbronn solution]
\label{rem:DH_numerics}
At $\beta=0$ the identity $b\star h=e$ was used to compute $(b_n)$ up to
$n=4000$, and the defining relation $A_{\Phi_H}(n)=1$ was verified
exactly along the way. The maxima of $|b_n|$ over dyadic windows grow
steadily, from $4.29$ on $[50,100)$ to $11.6$ on $[400,800)$ and $34.8$
on $[3200,4000)$, in agreement with
Corollary~\ref{cor:DH_not_HLR}. The window means grow slowly, from
$0.70$ to $1.19$ over the same interval, consistent with a small growth
exponent. The convolution identity in the proof of
Proposition~\ref{prop:anti_hlr_right} was checked to machine precision
on random data. The slopes of $\Phi_H$, the values
$S(M)=\sum_{j\le M}h(j)\lfloor M/j\rfloor$ on
$(1/(M+1),1/M]$, are positive over the computed interval but not monotone,
$S(9)<S(8)$ for instance, so $\Phi_H$ satisfies
Definition~\ref{def:BHF} only with the monotonicity of the slopes
relaxed. The statement of Conjecture~\ref{conj:anti_hlr} follows
\cite{Cloitre2016}, where $\Phi_H$ is the guiding example, and the
definitional slack is recorded here for completeness.
\end{remark}
The separation is the point of the chapter, and it is proved. One kernel, recalled from
\S\ref{sec:rational_kernel}, carries an index while its transform has no zero to carry one. Another has all the zeros of its transform on the line
$\Re z=1$ and loses transparency at one half. In the second case the analytic data is present and
exact, and it is not the arithmetic answer. What that leaves open is the arithmetic answer
itself. The exact index of the kernel broken at $\sqrt2$ is not determined here, Open
Problem~\ref{op:w10-global} carrying the two global estimates it would need. The agreement of the
two indices at the Ingham kernel, which holds under the hypothesis, is therefore a fact about
that kernel and not a law, and the difficulty of the hypothesis sits exactly there.

\chapter{Why the theory is called regular}
\label{chap:regularity}

\og\,La math\'ematique est l'art de donner le m\^eme nom \`a des choses diff\'erentes\,\fg,
wrote Henri Poincar\'e\index[names]{Poincar\'e, H.} \cite{PoincareMethode1908}. Mathematics is the art of giving the same name to
different things, provided that those things, different in their material, are alike in form. The
regularity index\index[terms]{regularity index} is introduced in that sense. It is not a new name for one familiar exponent. It
names a common response structure exhibited by kernels of very different kinds, the transmission of
power forcings below a critical threshold and their absorption\index[terms]{absorption} at a kernel determined rate beyond it.
Affine, arithmetic, discontinuous and fully two-variable kernels differ in their construction,
but the same transition can be asked of each of them.

This chapter concerns a single letter, the R of RAF. It stands for regular, and the choice is
not cosmetic. The theory takes its cue from the pure powers $n^{-\beta}$ on the right of the
defining equation\index[terms]{defining equation}, the prototypical regularly varying sequences, and inherits from them a close
tie to the Tauberian\index[terms]{Tauberian} tradition in which regular variation\index[terms]{regular variation} is the central notion. The pure power
was never the general forcing. The right hand side may carry a slowly varying factor,
$n^{-\beta}L(n)$, and the framework is designed to record that generality, which is what the
name records.

The chapter draws that tie out in full, through the work of Karamata\index[names]{Karamata, J.} and
of Bingham\index[names]{Bingham, N. H.}, Goldie\index[names]{Goldie, C. M.} and Teugels\index[names]{Teugels, J. L.}. It then separates three levels that
must not be conflated, namely an arbitrary forcing, a regularly varying model forcing, and an additive
remainder. The response to that remainder is described by an exact weighted norm in
Section~\ref{sec:perturbed-raf}. The chapter also tests the mechanism on the forcing
$n^{-\beta}\log n$. The answer is a second-order asymptotic of exactly the shape regular
variation predicts, carrying a constant equal to the logarithmic derivative of the Mellin
transform. What is established there is a theorem under hypotheses, stated as such, and not a
consequence of smoothness alone.

\section{Notation}
\label{sec:why_regular}

Throughout the sections that follow, $g:(0,1]\to\R$ is a real profile with $g(1)=1$ and
finite sampled values. The algebraic statements require no more. Whenever a Mellin transform
or an integral estimate is used, the corresponding integrability and continuation hypotheses are
stated explicitly. This includes the bounded Riemann-integrable\index[terms]{Riemann integrability}
branch and, whenever those hypotheses hold, the profiles unbounded but slowly varying at the
origin admitted by Proposition~\ref{prop:fgv_slowly_varying}. For a forcing
$r:\N^{*}\to\R$, $(a_n)$ is the unique solution of
\begin{equation}\label{eq:def}
  \sum_{k=1}^{n}a_k\,g\!\left(\frac{k}{n}\right)=r(n),\qquad n\ge 1,
\end{equation}
and $A(x)=\sum_{n\le x}a_n$. Here $\gstar$ is the arithmetic Mellin transform of Chapter~\ref{chap:fgv}, with logarithmic derivative $\gcirc=(\log\gstar)'=(\gstar)'/\gstar$.

\section{Regular variation and the vocabulary of the theory}\label{sec:rv}

The power law $n^{-\beta}$ is the test the theory was built on and not a limitation of it. Here is
the setting the general forcing comes from. A positive measurable function $f$ on $[X,\infty)$ is \emph{regularly varying of
index} $\rho\in\R$, written $f\in\mathrm{RV}_\rho$, if $f(\lambda x)/f(x)\to\lambda^{\rho}$
as $x\to\infty$ for every $\lambda>0$, and \emph{slowly varying} when $\rho=0$.
Every such $f$ factors as $x^{\rho}\ell(x)$ with $\ell$ slowly varying, the slow
factor carrying whatever is not a clean power, of which $\log x$ and $(\log x)^{\gamma}$ are
the basic instances. The subject originates with Karamata~\cite{Karamata1930}, and its
standard reference is Bingham\index[names]{Bingham, N. H.}, Goldie\index[names]{Goldie, C. M.} and Teugels\index[names]{Teugels, J. L.}~\cite{Bingham1989}, where $\mathrm{RV}$ sits among
a family of neighboring classes, each named for a manner of varying under scaling:
slow and rapid variation\index[terms]{rapid variation}, O-regular (dominated) variation, extended regular
variation, and de~Haan's class\index[terms]{de Haan class}\index[names]{Haan, L. de@de Haan, L.} $\Pi$ of functions with a slowly varying
second-order increment~\cite{deHaan1970}. Each qualifier records how tightly the
function is controlled.

Two kinds of theorem organize the subject. Abelian and Tauberian theorems relate
the regular variation\index[terms]{regular variation} of a function to that of a transform, whether Laplace,
Mellin, or a Dirichlet series, with Karamata's Tauberian theorem\index[terms]{Tauberian theorem} as the model.
Inside Tauberian theory sits a branch that runs the other way within a fixed transform,
recovering the regular variation\index[terms]{regular variation} of $f$ from that of $Tf$ for an averaging or convolution
operator $T$. These are the Mercerian theorems\index[terms]{Mercerian theorem}, named after \nm{Mercer}{J.} and his averaging
theorem of 1907. \nm{Korevaar}{J.}~\cite{Korevaar2004} presents them as a sub-theory of the
Tauberian one, and the treatment of \nm{Bingham}{N. H.}, \nm{Goldie}{C. M.} and
\nm{Teugels}{J. L.}~\cite[Ch.~5]{Bingham1989} places them in the setting of regular variation.
Such a theorem is available when the
Mellin transform of $T$ has no zero on the vertical line through the index in question, the
non-vanishing condition of Wiener's\index[names]{Wiener, N.} Tauberian theory. A Mercerian statement
is in this sense a Tauberian one freed of its side condition, the non-vanishing of the
transform doing the work a Tauberian hypothesis would otherwise have to do.

The operator \eqref{eq:def} is of this kind. It carries $(a_n)$ to
$r(n)=\sum_{k\le n}a_k g(k/n)$, with a kernel depending on $k$ and $n$ only through
$k/n$, so it is a multiplicative convolution whose Mellin transform is $\gstar$.
Solving \eqref{eq:def} for $(a_n)$ inverts the operator, and recovering the size of
$A(x)$ from that of $r$ is a Mercerian problem, which is what transparency records.
For a pure power the inversion preserves regular variation,
\begin{equation}\label{eq:transp}
  r(n)=n^{-\beta}\in\mathrm{RV}_{-\beta}
    \quad\Longrightarrow\quad
  A(x)\sim\frac{x^{-\beta}}{\gstar(\beta)}\in\mathrm{RV}_{-\beta}
  \qquad(\beta<\alpha(g)),
\end{equation}
the multiplier $1/\gstar(\beta)$ being the reciprocal of the Mellin transform at the
index, exactly as a Mercerian theorem gives. This holds while $\gstar(\beta)\ne0$.
Where the two indices are proved to agree, the index $\alpha(g)$ sits at the first zero of $\gstar$, where the Wiener condition
fails and preservation gives way to absorption,
\begin{equation}\label{eq:absorb}
  \beta\ge\alpha(g)\quad\Longrightarrow\quad A(x)=\mathcal{O}\!\left(x^{-\alpha(g)+\eps}\right)
   \quad(\forall\eps>0).
\end{equation}
The vocabulary is fixed by this analogy. The kernel varies well enough that regular
variation passes through it, up to a critical index and with a computable
multiplier, and good variation names that class, beside regular, slow and bounded. Index,
transparency, absorption, and smoothness are the regular-variation data of the
operator. Smoothness, the coincidence of $\alpha(g)$ with the analytic index $\inf\{\Re\rho:\gstar(\rho)=0\}$, is the condition under which \eqref{eq:transp} is a
theorem.

The naming is an analogy and not a definition, which is the licence the epigraph gives, and the
two words sit on the two sides of the defining equation\index[terms]{defining equation}. On the right the forcing is where
regular variation\index[terms]{regular variation} acts, the envelope $n^{-\beta}L(n)$ of Definition~\ref{def:raf_regular_variation}
being regularly varying in the sense of Karamata\index[names]{Karamata, J.}, and that is the word carried by the R of
regular arithmetic functions\index[terms]{regular arithmetic function}. On the left the kernel is where variation is a constraint
rather than a hypothesis on the data. Bounded variation\index[terms]{bounded variation} is the notion at stake there, and the
kernel this volume is built on sits at its edge, the perturbation measure of the Ingham
kernel\index[terms]{Ingham kernel} carrying infinite total variation by Proposition~\ref{prop:inf_mass}, a divergence that
records the infinitude of the primes. An unbounded kernel is admitted in its place when it varies
slowly at the origin, which is Proposition~\ref{prop:fgv_slowly_varying}, and that side is the
word carried by the variation of functions of good variation\index[terms]{function of good variation}. The same name is given to
different things because they are alike in form, and what these two share is a manner of varying
under which a finite threshold exists.

What the vocabulary names can now be stated. The definition of the index carries over from pure
powers to regularly varying envelopes without a change in either clause, and it does not replace
the earlier one. On the pure powers, where $L\equiv1$, the second clause below asks for
$A(n)\ll n^{-\alpha(G)}$ and is therefore stronger than the absorption
$\mathcal{O}(n^{-\alpha(G)+\eps})$ of Definition~\ref{def:reg_index}, which leaves room for a
logarithm at the index.

\begin{definition}\label{def:raf_regular_variation}
A function $K_{-\beta}:\R^{+}\to\R^{+}$ is regularly varying of index
$-\beta$ if $K_{-\beta}(x)=x^{-\beta}L(x)$ where $L$ is slowly varying\index[terms]{slowly varying}, that
is $L(\lambda x)/L(x)\to1$ as $x\to\infty$ for every $\lambda>0$. The kernel
$G$ is a RAF of index $\alpha(G)$ in the regularly varying sense if the
solutions of
\[
\sum_{k=1}^{n}a_k\,G(n,k)=K_{-\beta}(n)
\]
satisfy the two conditions
\begin{itemize}
\item $\beta<\alpha(G)\implies A(n)\sim C(\beta)\,K_{-\beta}(n)$,
\item $\beta\ge\alpha(G)\implies A(n)\ll n^{-\alpha(G)}\,L(n)$.
\end{itemize}
\end{definition}

The index separates two cases. For $\beta<\alpha(G)$ the partial sums are asymptotically
proportional to the forcing term\index[terms]{forcing term}, with an explicit constant, while for
$\beta\ge\alpha(G)$ they satisfy $A(n)\ll n^{-\alpha(G)}L(n)$ whatever the forcing. These two
cases were already met in Definition~\ref{def:reg_index} of Chapter~\ref{chap:raf}, extended here
from pure powers to regularly varying envelopes.

After a pure power comes the forcing $n^{-\beta}\log n$, computed in this chapter.
Its slowly varying factor $\log n$ keeps it in $\mathrm{RV}_{-\beta}$ while removing it from
the pure powers, making it the first properly second-order datum. In Karamata
theory the second-order behavior of averages is governed by de~Haan's class $\Pi$
and an auxiliary function, with the qualitative signature that a derivative in the
index appears. The theorem below is the arithmetic instance, and the answer keeps the
first-order multiplier $1/\gstar(\beta)$ on the $\log x$ term and adds a constant
proportional to $\gcirc(\beta)=(\log\gstar)'(\beta)$.

One difference from the classical theory bears stating. The operator sums over the
integers $k\le n$ rather than integrating, so $\gstar$ is an arithmetic Mellin
transform and can carry arithmetic content. For Ingham's kernel it is built from
$\zeta$, and its zeros are zeros of $\zeta$. The theorems of good variation are thus
the arithmetic counterparts of the Mellin-convolution Mercerian theorems
of~\cite{Bingham1989}, of the same shape but with a Mellin transform whose zero set is
number-theoretic.

A second difference bears stating beside it, and the object of this volume turns on it. The
Mellin convolution setting is treated in full in \cite[Ch.~4 and~5]{Bingham1989} and
in \cite{Korevaar2004}, and there the function is given and is compared with its transform. The
Abelian direction reads $f(x)\sim x^{\rho}\ell(x)$ into
$(k\star f)(x)\sim\check k(\rho)\,x^{\rho}\ell(x)$, a Tauberian converse recovers the first from the
second under side conditions, and a Mercerian statement passes from the convergence of the ratio
to the regular variation\index[terms]{regular variation} of $f$. The prototypes carry names and sit in that same text, the
theorem of \nm{Drasin}{D.} and \nm{Shea}{D. F.} for a kernel of one sign
at \cite[Th.~5.2.1]{Bingham1989} and \cite{DrasinShea1976}, and that of \nm{Jordan}{G. S.} for a
kernel that changes sign at \cite[Th.~5.3.1]{Bingham1989} and \cite{Jordan1974}, the technical
conditions of the second being removed later by \nm{Bingham}{N. H.} and
\nm{Inoue}{A.}~\cite{BinghamInoue2000c}.

Here the function is not given. The transformed object is prescribed, $A_g(x)=x^{-\beta}$, and the
partial sums are the unknown, so the constant $\check k(\rho)$ that the Abelian direction produces is
the reciprocal of the one the solution has to carry. That reciprocal is the transparent
constant\index[terms]{transparent constant} $1/\gstar(\beta)$ of this volume. The exponent then runs over the line while the
kernel stays fixed, and the threshold at which the solution stops carrying the rate is what no
statement made at a single index can produce.

\section{The theorem}\label{sec:statement}

Throughout this chapter a smooth function of good variation is one whose transform
$\gstar$ continues holomorphically and without zeros across
$\Omega=\{z\in\C:\Re z<\alpha(g)\}$, with polynomial growth of $1/\gstar$ on the
vertical lines of every closed substrip. These are hypotheses of the kind listed in
Conditional Theorem~\ref{cthm:ex-alpha-eta}, not consequences of pointwise regularity of
$g$. On $\Omega$ the functions $1/\gstar$ and $\gcirc$ are then holomorphic. Fix
$\beta<\alpha(g)$.

\begin{conditionaltheorem}\label{cthm:main}
Let $g$ be a smooth function of good variation of index $\alpha=\alpha(g)$, assume that
transparency holds with a power saving, uniformly on compact subsets of $\Omega$, in the precise
form of the remainder estimate \eqref{eq:T} displayed below, and let
$(a_n)$ solve $\sum_{k\le n}a_k\,g(k/n)=n^{-\beta}\log n$ with $\beta<\alpha$. Then,
as $x\to\infty$,
\begin{equation}\label{eq:box}
  \boxed{\;A(x)=\frac{1}{\gstar(\beta)}\,x^{-\beta}\log x
      +\frac{\gcirc(\beta)}{\gstar(\beta)}\,x^{-\beta}+o\!\left(x^{-\beta}\right).\;}
\end{equation}
The error is $\mathcal{O}\!\big(x^{-\beta-\delta}\log x\big)$ for some $\delta>0$ depending
on $\beta$.
\end{conditionaltheorem}

Equation \eqref{eq:box} is the second-order refinement of the Mercerian law
\eqref{eq:transp}. The $\log x$ term reproduces the slowly varying factor of the
forcing with the same multiplier $1/\gstar(\beta)$, while the correction
$\gcirc(\beta)/\gstar(\beta)$ is the auxiliary constant of the second-order theory,
$\gcirc$ being the logarithmic derivative of the Mellin transform. The argument uses
one exact identity together with a single estimate on the transparency remainder,
\eqref{eq:T} below, taken as the hypothesis of the chapter, Section~\ref{sec:contour}
establishing that hypothesis for Ingham under RH with the finite contour estimates of
Appendix~\ref{app:perron}. The case
$\beta\ge\alpha$ is treated in Lemma~\ref{lem:absorb}.

\section{Proof of the theorem}\label{sec:proof}

For each fixed truncation the passage from the forcing to $A(x)$ is a finite linear
map whose coefficients depend only on $g$. This rigidly links the two forcings at
issue, since the solution for $n^{-\beta}\log n$ is, term by term, the $\beta$-derivative
of the solution for $n^{-\beta}$.

The solution of the forced equation is a fixed triangular combination of the forcing, the
combination not depending on the exponent.

\begin{proposition}\label{prop:struct}
There are real numbers $b_{n,k}$ $(1\le k\le n)$, depending only on $g$, with
$b_{n,n}=1$, such that the solution of \eqref{eq:def} is $a_n=\sum_{k\le n}b_{n,k}r(k)$.
For every $x\ge1$,
\[
  A(x)=\sum_{k\le x}C_k(x)\,r(k),\qquad C_k(x):=\sum_{k\le n\le x}b_{n,k}\in\R,
\]
the $C_k(x)$ being independent of the forcing. In particular
$z\mapsto A^{(z)}(x):=\sum_{k\le x}C_k(x)k^{-z}$, the solution for $r(n)=n^{-z}$, is
entire for each $x$, and the solution for $r(n)=n^{-\beta}\log n$ has summatory
function
\[
  A_{\log}(x)=\sum_{k\le x}C_k(x)\,k^{-\beta}\log k
    =-\frac{\partial}{\partial z}A^{(z)}(x)\Big|_{z=\beta}.
\]
\end{proposition}

\begin{proof}
Written as $G\mathbf a=\mathbf r$ with $G=(g(k/n))_{n\ge k\ge1}$, the system is
lower triangular with diagonal $g(1)=1$, hence unit lower triangular and invertible
with unit lower-triangular inverse $(b_{n,k})$ depending only on $g$, so
$a_n=\sum_{k\le n}b_{n,k}r(k)$. Summing and reordering the finite sum gives
$A(x)=\sum_{k\le x}C_k(x)r(k)$. Taking $r(k)=k^{-z}$ gives $A^{(z)}(x)$, entire in
$z$. Taking $r(k)=k^{-\beta}\log k$ and comparing with
$-\partial_z\sum_{k\le x}C_k(x)k^{-z}=\sum_{k\le x}C_k(x)k^{-z}\log k$ gives the last
identity.
\end{proof}

Write $\Psi(x,z):=x^{z}A^{(z)}(x)$, entire in $z$, so that transparency
\eqref{eq:transp} reads $\Psi(x,\beta)\to1/\gstar(\beta)$ for real $\beta<\alpha$ and
Proposition~\ref{prop:struct} becomes the exact identity
\begin{equation}\label{eq:Alog}
  A_{\log}(x)=-\frac{\partial}{\partial z}\big[x^{-z}\Psi(x,z)\big]_{z=\beta}
    =x^{-\beta}\log x\,\Psi(x,\beta)-x^{-\beta}\,\partial_z\Psi(x,\beta).
\end{equation}
The only limit below is $x\to\infty$, entering through $\Psi(x,\beta)$ and
$\partial_z\Psi(x,\beta)$.

The identity fixes what is needed of transparency. Its first term carries $\log x$,
so the constant in \eqref{eq:box} is smaller than the leading term by that factor,
and recovering it requires $\Psi(x,\beta)$ to relative accuracy $o(1/\log x)$
together with control of $\partial_z\Psi(x,\beta)$. Transparency by itself gives only
$\Psi(x,\beta)=1/\gstar(\beta)+o(1)$, which cannot separate the two terms, since an
error of exact order $1/\log x$ would shift the constant. What the proof uses instead
is that on every compact $K\subset\Omega$ there are $C_K<\infty$ and $\delta_K>0$
with
\begin{equation*}\label{eq:T}
  \sup_{z\in K}\Big|\,\Psi(x,z)-\frac{1}{\gstar(z)}\,\Big|\le C_K\,x^{-\delta_K}
  \qquad(x\ge2).\tag{T}
\end{equation*}
This estimate is the substantive hypothesis of the theorem. Appendix~\ref{app:perron}
proves the quantitative contour criterion, which Section~\ref{sec:contour} verifies for Ingham under RH.
For a general profile, local boundedness of the entire functions $\Psi(x,\cdot)$
on $\Omega$ is a separate condition to check. Once it holds, the theorem of Vitali
and Porter~\cite[Thm.~2.1]{Schiff1993} turns convergence along a real segment
with an accumulation point in $\Omega$ into locally uniform convergence there.
Cauchy's formula then gives convergence of the derivatives. The power saving
in \eqref{eq:T} requires the additional quantitative estimates.

Analyticity in the parameter enters at one point. The estimate \eqref{eq:T}, though
it bounds only the values of $\Psi$, controls the $\beta$-derivative as well,
uniformly on compact subintervals, which is what allows the asymptotic to be
differentiated. The same Cauchy bound\index[terms]{Cauchy bound} reappears for $\beta\ge\alpha$.

The expansion is uniform on compact ranges of the exponent.

\begin{lemma}\label{lem:unif}
Assume \eqref{eq:T} and let $[\beta_0,\beta_1]\subset(-\infty,\alpha)$ be compact.
Then, as $x\to\infty$,
\[
  \Psi(x,\beta)\to\frac{1}{\gstar(\beta)}\qquad\text{and}\qquad
  \partial_\beta\Psi(x,\beta)\to-\frac{\gcirc(\beta)}{\gstar(\beta)},
\]
both uniformly on $[\beta_0,\beta_1]$. Writing $\eps(x,z)=\Psi(x,z)-1/\gstar(z)$,
\begin{equation}\label{eq:cauchy}
  |\partial_z\eps(x,z_0)|\le\frac{1}{\rho}\sup_{|z-z_0|=\rho}|\eps(x,z)|
  \qquad(\overline D(z_0,\rho)\subset\Omega).
\end{equation}
\end{lemma}

\begin{proof}
Set $I=[\beta_0,\beta_1]$ and pick $\rho<\operatorname{dist}(I,\partial\Omega)$, so
that $K:=\{z:\operatorname{dist}(z,I)\le\rho\}\subset\Omega$ is compact. On $\Omega$
the function $1/\gstar$ is holomorphic and each $\Psi(x,\cdot)$ is entire, so
$\eps(x,\cdot)$ is holomorphic on $K$. By \eqref{eq:T},
$\sup_{z\in K}|\eps(x,z)|\le C_K x^{-\delta_K}=:\omega(x)\to0$, and restriction to
$I$ gives the first limit uniformly. Fixing $\beta\in I$ and applying Cauchy's
formula for the derivative on $\overline D(\beta,\rho)\subset K$ gives
\eqref{eq:cauchy}, hence $|\partial_z\eps(x,\beta)|\le\omega(x)/\rho$ uniformly in
$\beta\in I$. On $I$ the holomorphic derivative $\partial_z$ agrees with
$\partial_\beta$, and $\partial_z(1/\gstar)=-(\gstar)'/(\gstar)^{2}=-\gcirc/\gstar$,
so $\partial_\beta\Psi(x,\beta)=-\gcirc(\beta)/\gstar(\beta)+\mathcal{O}(\omega(x)/\rho)$
uniformly on $I$.
\end{proof}

Write $\eps(x)=\eps(x,\beta)$ and $\eps'(x)=\partial_z\eps(x,\beta)$ for a compact
interval $I$ containing $\beta$ in its interior, so that Lemma~\ref{lem:unif} gives
$|\eps(x)|\le\omega(x)$ and $|\eps'(x)|\le\omega(x)/\rho$ with $\omega(x)=\mathcal{O}(x^{-\delta})$.
Substituting $\Psi(x,\beta)=1/\gstar(\beta)+\eps(x)$ and
$\partial_z\Psi(x,\beta)=-\gcirc(\beta)/\gstar(\beta)+\eps'(x)$ into \eqref{eq:Alog},
\[
  A_{\log}(x)
   =\frac{x^{-\beta}\log x}{\gstar(\beta)}+\frac{\gcirc(\beta)}{\gstar(\beta)}x^{-\beta}
    +x^{-\beta}\big(\log x\,\eps(x)-\eps'(x)\big),
\]
whose first two terms are those of \eqref{eq:box} and whose remainder is at most
$x^{-\beta}\omega(x)(\log x+\rho^{-1})=\mathcal{O}(x^{-\beta-\delta}\log x)=o(x^{-\beta})$.
This proves \eqref{eq:box}. \qed

\medskip

The same identity disposes of the case $\beta\ge\alpha$, where the pure powers are
already small. Recorded as a lemma, it states that the logarithmic twist cannot
break the ceiling $x^{-\alpha}$ set by absorption, and it follows from the Cauchy
bound of Lemma~\ref{lem:unif} applied to a bound rather than an asymptotic. One uses
it whenever the index lies at or below the forcing exponent.

\begin{lemma}\label{lem:absorb}
Suppose \eqref{eq:absorb} holds uniformly on small complex discs about each
$\beta\ge\alpha$: for every $\eps>0$ there are $\rho_0>0$ and $C_{\eps,\beta}$ with
$|A^{(z)}(x)|\le C_{\eps,\beta}x^{-\alpha+\eps}$ for $x\ge2$ and $|z-\beta|\le\rho_0$.
Then for every $\beta\ge\alpha$ the solution of
$\sum_{k\le n}a_k g(k/n)=n^{-\beta}\log n$ satisfies $A(x)=\mathcal{O}(x^{-\alpha+\eps})$ for
every $\eps>0$, and $A(x)=\mathcal{O}(x^{-\alpha}\log x)$ if the bound holds in the sharp
form $\mathcal{O}(x^{-\alpha})$.
\end{lemma}

\begin{proof}
By Proposition~\ref{prop:struct}, $A_{\log}(x)=-\partial_z A^{(z)}(x)|_{z=\beta}$
with $z\mapsto A^{(z)}(x)$ entire, so Cauchy's estimate on $|z-\beta|=\rho$ with
$\rho=\min(\rho_0,\eps)$ gives
$|A_{\log}(x)|\le\rho^{-1}\sup_{|z-\beta|=\rho}|A^{(z)}(x)|\le C_{\eps,\beta}\rho^{-1}x^{-\alpha+\eps}$.
Arbitrariness of $\eps$ yields $A(x)=\mathcal{O}(x^{-\alpha+\eps})$. Under the sharp bound,
$\rho=1/\log x$ gives $|A_{\log}(x)|\le Cx^{-\alpha}\log x$.
\end{proof}

\section{The Mellin transform and the discrete defect}\label{sec:contour}

A finite Perron integral\index[terms]{Perron's formula} separates the coefficient error from
the estimates on the contour. Appendix~\ref{app:perron} proves those estimates. This section records the discrete
sum-minus-integral term, and establishes \eqref{eq:T} for the Ingham profile under the
Riemann hypothesis.

\begin{proofstatus}{For a general profile, Conditional Theorem~\ref{cthm:main} retains
\eqref{eq:T} as a hypothesis. Lemma~\ref{lem:perron-half-integer} and
Proposition~\ref{prop:perron-rectangle} in Appendix~\ref{app:perron} prove the finite inversion
and the contour bounds from coefficient and vertical estimates. Lemma~\ref{lem:perron-discrete-defect} gives an
explicit domain of holomorphy for the discrete defect under its stated assumptions.
Theorem~\ref{thm:ingham-uniform-power} verifies the required estimates for Ingham and
proves \eqref{eq:T} under RH on every compact subset of $\Re z<1/2$.
The remaining general question is to obtain these estimates from the hypotheses on
an arbitrary profile.}
\end{proofstatus}

The finite inversion is Lemma~\ref{lem:perron-half-integer}, and the
three contour errors are bounded in Proposition~\ref{prop:perron-rectangle},
both proved in Appendix~\ref{app:perron}. The appendix also gives the
admissible height window \eqref{eq:perron-height-window}. Here the task
is to identify the exact Dirichlet quotient of the summation equation and
to verify those hypotheses for it.

\subsection{The discrete Mellin defect}

Let $\mathcal R(s)=\sum_{n\ge1}r(n)n^{-s}$. On a right half plane where all
sums converge absolutely, the defining equation gives
\[
 \mathcal R(s)=\sum_{k\ge1}a_k\gamma_k(s),\qquad
 \gamma_k(s)=\sum_{n\ge k}g(k/n)n^{-s}.
\]
The integral corresponding to this inner sum is
\[
 \int_k^\infty g(k/y)y^{-s}\,dy
   =k^{1-s}\int_0^1g(t)t^{s-2}\,dt
   =k^{1-s}\frac{\gstar(1-s)}{s-1}.
\]
Define their difference by $E_k(s)$ and put
$\mathcal H_r(s)=\sum_{k\ge1}a_kE_k(s)$. Thus
\begin{equation}\label{eq:perron-discrete-series}
 \mathcal R(s)=\frac{\gstar(1-s)}{s-1}\mathcal A(s-1)+\mathcal H_r(s).
\end{equation}
The defect depends on the solution as well as on the profile. The next elementary
bound gives a domain in which it has a holomorphic continuation.

\begin{lemma}[A domain for the discrete defect]\label{lem:perron-discrete-defect}
Suppose $g\in C^1[0,1]$ and $|a_k|\le Ck^{-\theta}$ for some real $\theta$.
Then $E_k$ extends holomorphically to $\Re s>0$, and, writing $\sigma=\Re s$,
\begin{equation}\label{eq:perron-quadrature-bound}
 |E_k(s)|\le k^{-\sigma}\left(
 \int_0^1t^\sigma|g'(t)|\,dt
       +|s|\int_0^1t^{\sigma-1}|g(t)|\,dt\right).
\end{equation}
Consequently $\mathcal H_r$ is holomorphic for
$\Re s>\max(0,1-\theta)$ and is $\mathcal O(1+|s|)$ on each fixed closed
vertical substrip in that half plane.
\end{lemma}

\begin{proof}
Put $h_{k,s}(y)=g(k/y)y^{-s}$. Initially for $\sigma>1$,
\[
 E_k(s)=\sum_{n\ge k}\int_n^{n+1}
                 \big(h_{k,s}(n)-h_{k,s}(y)\big)\,dy.
\]
The absolute value is at most $\int_k^\infty|h'_{k,s}(y)|\,dy$. Since
\[
 h'_{k,s}(y)=-k y^{-s-2}g'(k/y)-s y^{-s-1}g(k/y),
\]
the substitution $t=k/y$ gives \eqref{eq:perron-quadrature-bound}.
The same expression for $E_k$ converges locally uniformly for $\sigma>0$,
as follows from the integrable bound for $h'_{k,s}$ on each compact set,
and hence defines a holomorphic continuation there.
On a compact subset of $\sigma>\max(0,1-\theta)$, the summands
$a_kE_k(s)$ are bounded by a constant times $k^{-\theta-\sigma}$ with
$\theta+\sigma>1$ uniformly. Their series therefore converges locally
uniformly. On a fixed closed vertical substrip the same comparison gives
the stated $\mathcal O(1+|s|)$ bound.
\end{proof}

Set $M_r(z)=\mathcal R(1-z)$ and $\mathcal D_r(z)=\mathcal H_r(1-z)$.
Equation \eqref{eq:perron-discrete-series}, with $w=-z$, gives the exact quotient
\[
 Q_r(z):=\frac{\mathcal A(-z)}{-z}
        =\frac{M_r(z)-\mathcal D_r(z)}{\gstar(z)}.
\]
At the initial Perron line $c'=-c$, the finite representation is therefore
\begin{equation}\label{eq:rep}
 A(x)=\frac1{2\pi i}\int_{c'-iT}^{c'+iT}
             Q_r(z)x^{-z}\,dz+\mathcal E_{
             P}(x,T),\qquad
 \mathcal E_{P}(x,T)\ll\frac{x^c+x^\kappa\log x}{T},
\end{equation}
under the coefficient hypothesis of Lemma~\ref{lem:perron-half-integer}.
Continuation of the quotient and its vertical bound must then be checked before
Proposition~\ref{prop:perron-rectangle} can be applied.

For the two forcings of this chapter,
\[
 M_r(z)=\zeta(1+\beta-z)=-\frac1{z-\beta}+\mathcal O(1),\qquad
 M_r(z)=-\zeta'(1+\beta-z)=\frac1{(z-\beta)^2}+\mathcal O(1),
\]
respectively. If $\mathcal D_r$ is holomorphic near $\beta$ and
$\gstar(\beta)\ne0$, the defect changes neither residue there.
Moving the upward $z$-line to the right contributes the negative residue.
Writing $h(z)=x^{-z}/\gstar(z)$ gives, for the pure power, $h(\beta)$, and
for the logarithmic forcing,
\begin{equation}\label{eq:perron-log-residue}
 -h'(\beta)=\frac{x^{-\beta}\log x}{\gstar(\beta)}
       +\frac{\gcirc(\beta)}{\gstar(\beta)}x^{-\beta}.
\end{equation}
In the second case the Laurent expansion of $-\zeta'$ has no simple-pole
term, so no Stieltjes constant is added. The error following these residues
is precisely the error bounded in \eqref{eq:perron-rectangle-bound}, applied
to the full quotient $Q_r(-w)$. Holomorphy of the defect near $\beta$
alone gives no bound for its integral along the shifted line.

\subsection{Uniform transparency for Ingham}

For Ingham the divisor identity supplies the complete quotient directly, including
the part that the smooth quadrature comparison treats as a defect.

\begin{theorem}[Uniform power saving under RH]\label{thm:ingham-uniform-power}
Assume the Riemann hypothesis. For the exact solutions
\[
 \sum_{k\le n}a_k^{(z)}\Phi(k/n)=n^{-z},\qquad
 A^{(z)}(x)=\sum_{n\le x}a_n^{(z)},
\]
estimate \eqref{eq:T} holds on every compact subset
$K\subset\{z\in\C:\Re z<1/2\}$, with $g=\Phi$.
More precisely, put $B=\max_{z\in K}\Re z$ and
\[
 d=\min\left(\frac14,\frac{1/2-B}{2}\right)>0.
\]
Then, uniformly for $z\in K$ and $x\ge2$,
\begin{equation}\label{eq:ingham-uniform-power}
 A^{(z)}(x)=\frac{x^{-z}}{\Phi^*(z)}
                 +\mathcal O_K\!\left(x^{-\Re z-d/2}\right).
\end{equation}
The formula includes $z=0$, with $\Phi^*(0)=1$.
\end{theorem}

\begin{proof}
Write $\Delta_z(1)=1$ and, for $n\ge2$,
\[
 \Delta_z(n)=n^{1-z}-(n-1)^{1-z}
           =(1-z)\int_{n-1}^n t^{-z}\,dt.
\]
Differencing the defining equation multiplied by $n$, and then applying
M\"obius inversion, gives
\begin{equation}\label{eq:ingham-complex-coefficients}
 n a_n^{(z)}=(\mu\star\Delta_z)(n).
\end{equation}
Uniformly on $K$,
\[
 \Delta_z(n)=(1-z)n^{-z}+u_z(n),\qquad
 |u_z(n)|\le C_K n^{-\Re z-1},
\]
where $u_z(1)=z$. For $n\ge2$ the bound follows by writing
$t^{-z}-n^{-z}=z\int_t^n v^{-z-1}\,dv$ inside the preceding integral.
If $v_z=\max(0,-\Re z)$, equation \eqref{eq:ingham-complex-coefficients}
also gives $|a_n^{(z)}|\le C_K\tau(n)n^{v_z-1}$.
For any $h>0$, $\tau(n)\ll_h n^h$. To verify the bound, factor
$n=\prod p^a$. For $p\ge2^{1/h}$ use $a+1\le2^a\le p^{ha}$.
For each of the finitely many smaller primes,
$\sup_{a\ge0}(a+1)p^{-ha}<\infty$, and the product of these finitely
many suprema is the required constant.

Put
\[
 U_z(s)=\sum_{n\ge1}u_z(n)n^{-s},\qquad
 C_z(s)=(1-z)\zeta(s+z)+U_z(s).
\]
The series for $U_z$ converges absolutely when $\Re(s+z)>0$.
On the half plane of absolute convergence, the convolution identity gives
\begin{equation}\label{eq:ingham-complex-quotient}
 F_z(w):=\frac{\mathcal A_z(w)}{w}
 =\frac{(1-z)\zeta(w+1+z)+U_z(w+1)}{w\zeta(w+1)},\qquad
 \mathcal A_z(w)=\sum_{n\ge1}a_n^{(z)}n^{-w}.
\end{equation}
For $z\in K$ choose
\[
 c_z=v_z+\frac14,\qquad \kappa_z=v_z+\frac18,
 \qquad b_z=-\Re z-d.
\]
These choices satisfy $c_z>\max(0,\kappa_z)$, $b_z>-1/2$ and
$b_z+1+\Re z=1-d\ge3/4$.
The series $U_z(w+1)$ is bounded uniformly on this strip, since
\[
 \sum_{n\ge1}|u_z(n)|n^{-\Re w-1}
       \ll_K\sum_{n\ge1}n^{-\Re(w+z)-2}
       \le\sum_{n\ge1}n^{-2+d}<\infty.
\]

Under RH, $w\zeta(w+1)$ has no zero for $\Re w>-1/2$.
At $w=0$ its value is $1$, by removal of the pole of $\zeta$.
Thus the only pole of $F_z$ in the strip is $w=-z$.
Its residue is
\[
 \frac{1-z}{-z\zeta(1-z)}=\frac1{\Phi^*(z)}.
\]
The formula is read by continuation at $z=0$. In that case
$u_0=0$, $F_0(w)=1/w$, and the residue is exactly $1$.

Lemma~\ref{lem:perron-zerofree-growth}, applied to $\zeta$ with $\theta=1/2$,
gives on every fixed closed half plane to the right of the critical line,
for every $\nu>0$ and large $|\Im s|$,
\[
 |\zeta(s)|+|1/\zeta(s)|\ll_\nu(1+|\Im s|)^\nu.
\]
Both arguments $w+1$ and $w+1+z$ in
\eqref{eq:ingham-complex-quotient} remain a fixed distance to the right
of the critical line, uniformly on $K$. Their imaginary parts differ
from $\Im w$ by a bounded quantity. Applying the displayed estimate
with a smaller exponent to each factor gives
\begin{equation}\label{eq:ingham-uniform-vertical}
 |F_z(\sigma+it)|\ll_{K,\nu}(1+|t|)^{-1+\nu}
 \qquad(b_z\le\sigma\le c_z)
\end{equation}
outside fixed neighborhoods of $w=-z$. On the bounded part of the
left boundary, uniformity follows by compactness, since its real
distance from that pole is $d$.

Apply Lemma~\ref{lem:perron-half-integer} and
Proposition~\ref{prop:perron-rectangle} at $x=N+1/2$ with $T=x^2$ and
$\nu=d/4$. The three errors, after multiplication by $x^{\Re z}$, are
bounded respectively by
\begin{align*}
 x^{b_z+\Re z}T^\nu&=x^{-d+2\nu}=x^{-d/2},\\
 \frac{x^{c_z+\Re z}T^{\nu-1}}{\log x}
   &\le\frac{x^{3/4-2+2\nu}}{\log x},\\
 \frac{x^{c_z+\Re z}+x^{\kappa_z+\Re z}\log x}{T}
   &\ll x^{-5/4}(1+\log x).
\end{align*}
Here $c_z+\Re z=\max(\Re z,0)+1/4<3/4$ and
$0<d\le1/4$, so the last two errors are also
$\mathcal O_K(x^{-d/2})$. Heights below the fixed lower bound needed
for the contour concern a bounded range of $x$ and are absorbed into
the constant. The residue gives \eqref{eq:ingham-uniform-power} at
half integers.

For an arbitrary real $x\ge2$, set $y=\lfloor x\rfloor+1/2$.
Then $A^{(z)}(x)=A^{(z)}(y)$, and
$y^{-z}=x^{-z}+\mathcal O_K(x^{-\Re z-1})$, uniformly on $K$.
This last error is smaller than the one in
\eqref{eq:ingham-uniform-power}, which proves the assertion for all $x$.
\end{proof}

\section{Perturbed equations and remainder transfer}\label{sec:perturbed-raf}

A perturbation of a RAF equation\index[terms]{perturbed RAF equation} is an additive error on its right hand side, while the kernel
and therefore its inverse remain fixed. More precisely, choose a model forcing $r_0$ and write
\begin{equation}\label{eq:perturbed-raf}
  r(n)=r_0(n)+e(n),\qquad
  \sum_{k\le n}a_k\,g\!\left(\frac{k}{n}\right)=r_0(n)+e(n).
\end{equation}
Let $(a_n^{[0]})$ solve the equation with right hand side $r_0(n)$, and put
\[
 d_n:=a_n-a_n^{[0]},\qquad
 A_0(x):=\sum_{n\le x}a_n^{[0]},\qquad
 D(x):=A(x)-A_0(x)=\sum_{n\le x}d_n.
\]
Then $d$ solves the exact error equation
\begin{equation}\label{eq:error-equation}
  \sum_{k\le n}d_k\,g\!\left(\frac{k}{n}\right)=e(n).
\end{equation}
Thus the model, the perturbation and the response are three distinct objects. The first may be
$r_0(n)=n^{-\beta}L(n)$, but the results below do not require that choice.

The cumulative inverse coefficients $C_k(x)$ of Proposition~\ref{prop:struct} give an exact
\index[terms]{remainder transfer}
measure of the largest response that a prescribed error envelope can produce.

\begin{proposition}[Exact remainder-transfer norm]\label{prop:remainder-transfer}
Let $w:\N^*\to(0,\infty)$ be a positive weight and define
\begin{equation}\label{eq:transfer-modulus}
  \mathfrak M_{g,w}(x):=\sum_{k\le x}|C_k(x)|\,w(k),
  \qquad
  \|e\|_{w,x}:=\max_{k\le x}\frac{|e(k)|}{w(k)}.
\end{equation}
For the perturbed equation \eqref{eq:perturbed-raf},
\begin{equation}\label{eq:exact-error-response}
  D(x)=\sum_{k\le x}C_k(x)e(k),
  \qquad
  |D(x)|\le \|e\|_{w,x}\,\mathfrak M_{g,w}(x).
\end{equation}
Moreover the modulus is the exact weighted $\ell^\infty$ operator norm at $x$:
\begin{equation}\label{eq:transfer-norm-exact}
 \mathfrak M_{g,w}(x)
   =\sup_{|u(k)|\le w(k)}
      \left|\sum_{k\le x}C_k(x)u(k)\right|.
\end{equation}
Consequently, if $e(n)=\mathcal O(w(n))$ and
$\mathfrak M_{g,w}(x)=\mathcal O(V(x))$, then $D(x)=\mathcal O(V(x))$.
\end{proposition}

\begin{proof}
Linearity of the triangular inverse and Proposition~\ref{prop:struct}, applied to the error
equation \eqref{eq:error-equation}, give the identity in \eqref{eq:exact-error-response}. The
inequality follows from the triangle inequality. Conversely, for a fixed $x$, choose
$u(k)=w(k)\operatorname{sgn}C_k(x)$, with any value in $[-w(k),w(k)]$ when $C_k(x)=0$.
The resulting sum is $\mathfrak M_{g,w}(x)$, which proves \eqref{eq:transfer-norm-exact}.
\end{proof}

This norm is deliberately a worst-case quantity, since it allows the signs of the perturbation to be
chosen adversarially. It is therefore the correct uniform criterion for a whole error class. A
particular arithmetic perturbation may be much smaller because of cancellation, and that sharper
question belongs to the arithmetic analysis of the profile.

\begin{corollary}[Stability of a regularly varying main term]\label{cor:perturbed-main-term}
Let $K_{-\beta}(x)=x^{-\beta}L(x)$ and suppose the solution for the model forcing
$r_0(n)=K_{-\beta}(n)$ satisfies
\[
  A_0(x)=C(\beta)K_{-\beta}(x)+E_0(x).
\]
If $e(n)=\mathcal O(w(n))$ and
$\mathfrak M_{g,w}(x)=\mathcal O(V(x))$, then the solution of
\eqref{eq:perturbed-raf} satisfies
\begin{equation}\label{eq:perturbed-main-term}
  A(x)=C(\beta)K_{-\beta}(x)+E_0(x)+\mathcal O(V(x)).
\end{equation}
In particular, if $E_0(x)=o(K_{-\beta}(x))$ and $V(x)=o(K_{-\beta}(x))$, the transparent
asymptotic is unchanged. If instead $\beta\ge\alpha(g)$ and
$A_0(x)=\mathcal O(x^{-\alpha(g)}L(x))$, then
\[
  A(x)=\mathcal O\!\left(x^{-\alpha(g)}L(x)+V(x)\right).
\]
\end{corollary}

\begin{proof}
Add the estimate for $D(x)=A(x)-A_0(x)$ from
Proposition~\ref{prop:remainder-transfer} to the stated estimate for $A_0(x)$.
\end{proof}

The direct, Abelian direction has an equally concrete form. It starts from control of the
partial sums of a coefficient error and asks how large the resulting right hand error can be.

\begin{proposition}[Abelian remainder bound]\label{prop:abelian-remainder-bound}
Let $(d_n)$ be any real sequence, $D(x)=\sum_{n\le x}d_n$, and
\[
  e(n):=\sum_{k\le n}d_k\,g\!\left(\frac{k}{n}\right).
\]
Then Abel summation gives the exact identity
\begin{equation}\label{eq:abelian-remainder-identity}
 e(n)=D(n)g(1)
   +\sum_{k=1}^{n-1}D(k)
      \left[g\!\left(\frac{k}{n}\right)-g\!\left(\frac{k+1}{n}\right)\right].
\end{equation}
If $g$ has bounded variation on $(0,1]$, then
\begin{equation}\label{eq:abelian-remainder-bv}
 |e(n)|\le
 \bigl(|g(1)|+\operatorname{Var}_{(0,1]}g\bigr)
 \max_{1\le k\le n}|D(k)|.
\end{equation}
\end{proposition}

\begin{proof}
Write $d_k=D(k)-D(k-1)$, with $D(0)=0$, and telescope. This gives
\eqref{eq:abelian-remainder-identity}. Taking absolute values and bounding the variation of
$g$ along the partition $1/n,2/n,\ldots,1$ gives \eqref{eq:abelian-remainder-bv}.
\end{proof}

Proposition~\ref{prop:abelian-remainder-bound} and
Proposition~\ref{prop:remainder-transfer} display the two sides of remainder theory. The Abelian
side sends coefficient control to a right hand error. The Tauberian, or here Mercerian, side
inverts that error and requires control of the cumulative resolvent. The bounded-variation
estimate does not cover profiles of infinite total variation, including the Ingham profile. In
that arithmetic branch cancellation replaces total variation, as developed in
Chapter~\ref{chap:abelian_densities}. This is the discrete RAF counterpart of Tauberian
\index[terms]{Tauberian remainder theory}
remainder theory in \cite[Ch.~VII, pp.~343--419]{Korevaar2004}.

The contour representation also applies to the error equation. Suppose its
coefficients satisfy $|d_n|\le Cn^{\kappa-1}$ and the Mellin identity
\eqref{eq:perron-discrete-series} is available for that equation. Put
$M_e(z)=\sum_{n\ge1}e(n)n^{z-1}$ and let $\mathcal D_e$ denote its discrete
defect. For $x=N+1/2$ and $c>\max(0,\kappa)$,
\begin{equation}\label{eq:perturbation-contour}
 D(x)=\frac1{2\pi i}\int_{-c-iT}^{-c+iT}
       \frac{M_e(z)-\mathcal D_e(z)}{\gstar(z)}x^{-z}\,dz
       +\mathcal O\!\left(\frac{x^c+x^\kappa\log x}{T}\right).
\end{equation}
This is Lemma~\ref{lem:perron-half-integer} applied to $(d_n)$ and the exact
Mellin identity, so the remainder is quantified. A shift uses the continuation
and vertical bounds of Proposition~\ref{prop:perron-rectangle} for this full
quotient.

At a simple zero $\rho$ of $\gstar$ where the numerator is holomorphic, the
term contributed by a shift to the right is
\[
 -\frac{M_e(\rho)-\mathcal D_e(\rho)}{(\gstar)'(\rho)}x^{-\rho}.
\]
For a zero of multiplicity $m$, expand the quotient in a Laurent series and use
\[
 x^{-z}=x^{-\rho}\sum_{j\ge0}\frac{(-\log x)^j(z-\rho)^j}{j!}.
\]
Only $0\le j\le m-1$ can enter the residue, giving $x^{-\rho}$ times a
polynomial in $\log x$ of degree at most $m-1$. Cancellation by the numerator
can lower the degree or remove the term. For real data, conjugate poles
combine into oscillations in $\log x$. The remaining contour and truncation
errors are bounded by \eqref{eq:perron-rectangle-bound}.
The finite polynomial case is developed in Chapter~\ref{chap:trace_poly},
and the infinite-mode orthorecursive case in Chapter~\ref{chap:trace}.

The general remainder mechanism belongs to the present volume, with algebraic inversion here,
arithmetic Abelian estimates in Chapter~\ref{chap:abelian_densities}, and spectral refinement in
Chapters~\ref{chap:trace_poly} and~\ref{chap:trace}. Volume~II specializes this framework
systematically to the Ingham kernel, its related families, and the Riemann hypothesis in its
classical and generalized forms.
The reciprocal-completion example of Appendix~\ref{app:O} makes this division concrete. Its
direct side carries divisor and circle remainders, while its inverse side carries a distinct
regularity index. Definition~\ref{def:O_defect} records the gap and Open
Problem~\ref{op:O_balance} asks when it vanishes.

\section{An application to the Ingham kernel}\label{sec:ingham}

Under the Riemann hypothesis, Theorem~\ref{thm:ingham-uniform-power} proves
\eqref{eq:T} for Ingham on every compact subset of $\Re z<1/2$.
Conditional Theorem~\ref{cthm:main} therefore applies to the logarithmic forcing
for every real $\beta<1/2$, with its uniformity hypothesis now verified.

For $g=\Phi$,
$\Phi^*(z)=\frac{z}{z-1}\zeta(1-z)$, and differentiation gives
\[
 \Phi^\circ(z)=\frac1z+\frac1{1-z}
                    -\frac{\zeta'(1-z)}{\zeta(1-z)}.
\]
Consequently the solution of
$\sum_{k\le n}a_k\Phi(k/n)=n^{-\beta}\log n$ satisfies, under RH,
\begin{align*}
 A(x)={}&\frac{x^{-\beta}\log x}{\Phi^*(\beta)}\\
 &+\frac1{\Phi^*(\beta)}
 \left(\frac1\beta+\frac1{1-\beta}
              -\frac{\zeta'(1-\beta)}{\zeta(1-\beta)}\right)x^{-\beta}
 +\mathcal O_\beta(x^{-\beta-\delta}\log x)
\end{align*}
for some $\delta>0$. At $\beta=0$ the constants are read by continuation.
Indeed $\zeta(1-z)=-1/z+\gamma+\mathcal O(z)$ gives
\[
 \Phi^*(z)=1+(1-\gamma)z+\mathcal O(z^2),\qquad
 \Phi^\circ(0)=1-\gamma,
\]
and hence the formula becomes
\[
 A(x)=\log x+1-\gamma+\mathcal O(x^{-\delta}\log x).
\]
The derivative estimate needed here follows from Lemma~\ref{lem:unif} on a
small closed complex disc about $\beta$ contained in $\Re z<1/2$.

The location of the analytic threshold can also be read without assuming RH.
Let
\[
 \Theta=\sup\{\Re\rho\mid\zeta(\rho)=0,\ 0<\Re\rho<1\}.
\]
Every nontrivial zero $\rho$ gives a zero $z=1-\rho$ of $\Phi^*$, and
taking the infimum of their real parts gives $1-\Theta$.
The trivial zeros give $z=3,5,7,\ldots$, while the apparent zero at $z=0$
is removed by the pole of $\zeta(1-z)$, as computed above. Thus
\[
 \eta(\Phi)=1-\Theta.
\]
Under RH this value is $1/2$, and
Theorem~\ref{thm:ingham-uniform-power} supplies uniform transparency with a
power saving on its left. Absorption and sharpness belong to the master equivalence
$\mathrm{RH}\iff\alpha(\Phi)=1/2$, established in
Theorem~\ref{thm:tauberian_rh} and Corollary~\ref{cor:phi_half}. The additional
conclusion here is the uniform complex-parameter estimate that permits
differentiation and gives both terms of the logarithmic response.

\section{A numerical check}\label{sec:num}

For the affine kernel everything is explicit. Here $\gstar(z)=(\lambda-z)/(1-z)$ gives
$1/\gstar(\beta)=(1-\beta)/(\lambda-\beta)$ and
$\gcirc(\beta)/\gstar(\beta)=(\lambda-1)/(\lambda-\beta)^{2}$, so that
Conditional Theorem~\ref{cthm:main} predicts
$x^{\beta}A(x)-\frac{1-\beta}{\lambda-\beta}\log x\to\frac{\lambda-1}{(\lambda-\beta)^{2}}$.
As $g(k/n)$ is affine in $k$, the recursion \eqref{eq:def} collapses to two running
sums, $a_n=n^{-\beta}\log n-\frac{1-\lambda}{n}\sum_{k<n}k\,a_k-\lambda\sum_{k<n}a_k$,
and $A$ is computed exactly in $\mathcal{O}(N)$ steps. With $\lambda=0.7$ and
$N=2\times10^{6}$ the residual $x^{\beta}A(x)-\frac{1-\beta}{\lambda-\beta}\log x$ is:

\begin{center}
\renewcommand{\arraystretch}{1.15}
\begin{tabular}{r|ccc}
$x$ & $\beta=-0.5$ & $\beta=0.2$ & $\beta=0.4$\\\hline
$10^{3}$        & $-0.2104$ & $-1.1701$ & $-2.9312$\\
$10^{4}$        & $-0.2086$ & $-1.1905$ & $-3.1315$\\
$10^{5}$        & $-0.2084$ & $-1.1970$ & $-3.2321$\\
$10^{6}$        & $-0.2083$ & $-1.1991$ & $-3.2826$\\
$2\times10^{6}$ & $-0.2083$ & $-1.1993$ & $-3.2921$\\\hline
predicted       & $-0.20833$ & $-1.20000$ & $-3.33333$
\end{tabular}
\end{center}

The residual tends to the predicted second constant, not merely to the leading
behavior, confirming both terms of \eqref{eq:box}. Convergence is fast at
$\beta=-0.5$, well inside the strip, and slow at $\beta=0.4$, near $\alpha=0.7$, as
the error $\mathcal{O}(x^{-\beta-\delta}\log x)$ with $\delta$ shrinking toward the index
would predict.

\section*{Bibliographical notes}
\addcontentsline{toc}{section}{Bibliographical notes}

The vocabulary is that of regular variation carried into an arithmetic setting, the
R of RAF being Karamata's regular. Good variation is the class of kernels for which
the multiplicative-convolution operator \eqref{eq:def} preserves regular variation up
to a critical index, a Mercerian property in the sense of~\cite[Ch.~5]{Bingham1989}. The
index is reached, wherever a theorem identifies it with the analytic one, at the first zero of the Mellin transform, where the non-vanishing condition of
Wiener fails and absorption begins. Logarithmic forcing is the first second-order
datum, and the constant it produces, $\gcirc=(\log\gstar)'$, is the arithmetic
analogue of the auxiliary function of de~Haan's class $\Pi$~\cite{deHaan1970}. Karamata's
papers~\cite{Karamata1930} and the treatise~\cite{Bingham1989} are the sources for the continuous
theory, and Korevaar\index[names]{Korevaar, J.}~\cite[especially Ch.~VII, pp.~343--419]{Korevaar2004} gives the Tauberian,
Mercerian and remainder-theoretic background.
What separates these statements from their classical models is a Mellin transform
whose zero set is number-theoretic, as for Ingham's kernel.

\part{Stability and equilibrium of the index}
\label{part:four}

\rafepigraph{Heureusement il s'est trouvé un domaine intermédiaire entre l'arithmétique et la théorie riemannienne, et qui possède, avec chacune de ces deux dernières théories, des ressemblances beaucoup plus étroites qu'elles n'en ont entre elles\,; il s'agit des fonctions algébriques «\,sur un corps fini\,».}{Fortunately there turned out to be an intermediate domain between arithmetic and Riemannian theory, one bearing far closer resemblances to each of those two theories than they bear to one another. It is the theory of algebraic functions over a finite field.}{André Weil, \emph{De la métaphysique aux mathématiques} (1960)~\cite{WeilMetaphysique1960}}

The index has been computed and shown to be arithmetic in its own right, while the associated
triangular systems have been inverted. Throughout, the index has been read in one coordinate,
the ratio $k/n$. This part changes the coordinate.

A gauge deforms a kernel by measuring the two indices through an auxiliary scale, replacing the
ratio $k/n$ by $f(k)/f(n)$. Stability asks whether the regularity index survives the
deformation, and equilibrium asks whether it is left unchanged. For the Ingham operator the
invariance of the index under the gauges $f_q(x)=q^{x}+1$, with $q\ge2$ an integer, would give the
Riemann hypothesis, a conjectural route to it on the side of deformation rather than of decay.

A power gauge changes nothing, as \S\ref{sec:gauge_transform} shows in two lines, so the
whole polynomial family of coordinate changes is invisible to the theory. An exponential gauge
leaves the theory altogether, the deformed kernel depending on the difference of the two
arguments and no longer on their ratio. Everything in this part happens on that second side.

\chapter{Gauge deformations: stability and equilibrium}
\label{chap:gauge}

The arithmetic separation of Chapter~\ref{chap:diophantine} constrains the
Ingham operator against changes of the kernel, and the multiplicative
obstruction of Chapter~\ref{chap:principles} shows that the zeros of
$\zeta$ reappear in every multiplicative combination. A different
deformation leaves the kernel untouched and changes instead the
coordinates in which the operator is evaluated. This chapter studies the
behavior of the regularity index under such a change of variable, and
formulates the invariance of the index for the Ingham operator as a
conjectural route to the Riemann hypothesis.

The order of the chapter is the following. A first section records where the device of
adjoining an auxiliary parameter comes from, since it is not new and its history says what it
is for. A second section carries the whole apparatus, the deformed equation, the probe that
defines the transform in a gauge from the kernel alone, and the observation that separates the
two regimes, namely that a power gauge leaves the transform where it was while an exponential
gauge turns the kernel into a Toeplitz\index[terms]{Toeplitz} kernel and forces the transform to be rebuilt. The
two notions then follow with an example each, stability with a kernel that keeps an index and
loses its value, equilibrium with the same kernel keeping both. The rational kernel marks the
limits of the phenomenon. The next chapter carries the Ingham function, where the zeros of the
transform gauged by $f(x)=q^{x}+1$ sit on the critical line with no hypothesis on $\zeta$, and where the
invariance of the index becomes a formulation of the Riemann hypothesis\index[terms]{Riemann hypothesis}.

\section{A device from Tauberian theory}
\label{sec:gauge_history}

The idea has a precedent in Littlewood's proof of his Tauberian theorem\index[terms]{Tauberian theorem}, the
converse of Abel's\index[names]{Abel, N. H.} theorem under a boundedness condition on the terms of the
series. Tauber had shown that Abel summability of $\sum a_n$ together with
$na_n=o(1)$ forces convergence \cite{Tauber1897}, and Littlewood replaced the
vanishing condition by the boundedness $na_n=\mathcal{O}(1)$ \cite{Littlewood1911}.
The passage from the small $o$ to the large $\mathcal{O}$ is the whole of the
difficulty, since a bounded sequence of terms may oscillate where a vanishing one
cannot, and the proof must recover convergence from a hypothesis that no longer
controls the individual term. Littlewood carried it through by introducing a
parameter that did not belong to the original statement and by passing from the
single problem to the family it generated, and the structure that this family made
visible was what let the proof close. Dyson\index[names]{Dyson, F. J.}, recalling the theorem, described the
method as ``a tour de force done by invoking a completely unmotivated new parameter
which turned out to be the key to the proof'', and used the same device himself to
obtain the phase transition of the one dimensional Ising ferromagnet
\cite{DysonWoS,Dyson1969}. A gauge plays the same part here. It adjoins to a fixed kernel an
auxiliary scale the kernel did not carry, and the response of the index to that
scale exposes structure that a single kernel keeps hidden.

The device is older and wider than that one proof. Three further instances mark out what it
does.

Karamata proved the Hardy-Littlewood Tauberian theorem by replacing the single discontinuous
test function that the statement asks for by the family of polynomials in $e^{-t}$, for which
the hypothesis is available term by term, and then passing to the limit inside the
family, the device of \nm{Karamata}{J.}~\cite{Karamata1930b}. The parameter is the degree, and it belongs to the proof and not
to the statement.

Wiener\index[terms]{Wiener's theorem}\index[names]{Wiener, N.} replaced a single averaging kernel by the closed linear span of its
translates \cite{Wiener1932}. His theorem says that the conclusion drawn for one kernel passes
to every integrable kernel exactly when the Fourier transform\index[terms]{Fourier transform}\index[names]{Fourier, J.} of the first has no zero, so the
family is the object that carries the argument and the transform is what decides whether the
family is large enough. That non-vanishing condition is the one that
Section~\ref{sec:rv} identifies with transparency\index[terms]{transparency} in the present setting, and it enters
here for a second reason, as the criterion attached to a family rather than to a kernel.

Newman\index[names]{Newman, D. J.} proved the prime number theorem\index[terms]{prime number theorem} by inserting into a contour integral\index[terms]{contour integral} the
factor $1+z^{2}/R^{2}$, whose only role is to carry a free radius $R$ that the original
problem does not mention, and by letting that radius grow at the end, the device of \nm{Newman}{D. J.}~\cite{Newman1980}. The
auxiliary quantity is again introduced for the sole purpose of generating a family.

What a gauge shares with these is the adjunction, to a fixed problem, of a scale the problem
did not carry, and the reading of the answer off the family rather than off the single member.
What separates it is that the family is indexed by the coordinate system in which the equation
is evaluated, so the question it raises has a name of its own. Whether the invariant survives
the change of coordinates is stability, and whether its value survives is equilibrium.

\section{The equation in a gauge and its transform}
\label{sec:gauge_transform}

Let $f:\R\to\R^+$ be increasing and unbounded. The gauged defining
relation replaces the evaluation lattice $k/n$ by $f(k)/f(n)$,
\begin{equation}\label{eq:gauged_recurrence}
\sum_{k=1}^{n}a_k\,G\big(f(n),f(k)\big)=f(n)^{-\beta},
\end{equation}
where $G(N,k)$ stands for the kernel evaluated with second argument
scaled by the first, $G(N,k)=g(k/N)$ in the homogeneous case. The
spectral object attached to \eqref{eq:gauged_recurrence} is the
transform of $G$ read in the gauge $f$.

The transform attached to \eqref{eq:gauged_recurrence} is built from the kernel and the gauge
alone. It is the probe of Definition~\ref{def:probes} with the evaluation lattice $k/n$
replaced by $f(k)/f(n)$, and it uses neither the triangular equation nor any asymptotic form
of its solution.

\begin{definition}[The probe in a gauge]\index[terms]{gauged regularity index}
\label{def:gauge_probe}
Write $f_k=f(k)$. For $\Re z<0$ and $n\ge1$ put $u_{n,z}(0)=0$ and
$u_{n,z}(k)=(f_k/f_n)^{-z}$ for $1\le k\le n$, and set
\begin{equation}\label{eq:gauge_probe}
\mathcal G_{f,n}(z)
=\sum_{k=1}^{n}\bigl(u_{n,z}(k)-u_{n,z}(k-1)\bigr)\,G\bigl(f_n,f_k\bigr).
\end{equation}
The arithmetic Mellin transform\index[terms]{arithmetic Mellin transform} of $G$ in the gauge $f$ is
$G^{*}_f(z)=\lim_{n\to\infty}\mathcal G_{f,n}(z)$ where that limit exists, continued to its
maximal domain, and the analytic index in the gauge $f$ is
$\eta_f(G)=\inf\{\Re\rho:G^{*}_f(\rho)=0\}$.
\end{definition}

At $f(x)=x$ the weights $u_{n,z}(k)$ are $(k/n)^{-z}$ and \eqref{eq:gauge_probe} is the
averaging limit $G^{*}(z)=\lim_n n^{z}\sum_{k\le n}[k^{-z}-(k-1)^{-z}]G(n,k)$ of
Chapter~\ref{chap:fgv}, so the definition extends the one already in force and does not
replace it. What the transform in a gauge is for is the asymptotic
\begin{equation}\label{eq:gstar_f}
A(n)\sim\frac{f(n)^{-\beta}}{G^*_f(\beta)}\qquad(n\to\infty),
\quad\text{for }\beta<\alpha_f(G),
\end{equation}
which is a property to be established for each pair $(G,f)$ and not a definition of
$G^{*}_f$.

Two gauges of opposite type show at once why the question has content.

\begin{proposition}
\label{prop:power_gauge}
Let $G(n,k)=g(k/n)$ and $f(x)=x^{p}$ with $p>0$. Then for every $n$ and every $\Re z<0$,
\begin{equation}\label{eq:power_gauge_probe}
\mathcal G_{f,n}(z)=\mathcal H_{n}(pz),
\end{equation}
where $\mathcal H_n$ is the ungauged probe of the composed profile $h(t)=g(t^{p})$.
Consequently $h^{*}(z)=g^{*}(z/p)$ and
\[
G^{*}_f=g^{*},
\qquad
\eta_f(G)=\eta(g).
\]
Moreover the gauged equation \eqref{eq:gauged_recurrence} at the exponent $\beta$ is the
ordinary defining equation\index[terms]{defining equation} for $h$ at the exponent $p\beta$, so $G$ is stable with respect to
$f$ if and only if $h$ is a function of good variation\index[terms]{function of good variation}, and then
$\alpha_f(G)=\alpha(h)/p$. Every power gauge under which the composed profile is a function
of good variation of index $p\,\alpha(g)$ is an equilibrium.
\end{proposition}

\begin{proof}
Since $f_k/f_n=(k/n)^{p}$, the weights are $u_{n,z}(k)=(k/n)^{-pz}$ and the kernel value is
$g((k/n)^{p})=h(k/n)$, which is \eqref{eq:power_gauge_probe}. Substituting $u=t^{p}$ in
$h^{*}(z)=-z\int_0^1g(t^{p})t^{-z-1}\,dt$ gives
$-\tfrac zp\int_0^1g(u)u^{-z/p-1}\,du=g^{*}(z/p)$, so
$G^{*}_f(z)=h^{*}(pz)=g^{*}(z)$ and the zeros of $G^{*}_f$ are those of $g^{*}$. For the
equation, $G(f(n),f(k))=g(k^{p}/n^{p})=h(k/n)$ and $f(n)^{-\beta}=n^{-p\beta}$, so
\eqref{eq:gauged_recurrence} reads $\sum_{k\le n}a_kh(k/n)=n^{-p\beta}$. Transparency for $h$
at $p\beta$ is $A(n)\sim n^{-p\beta}/h^{*}(p\beta)=f(n)^{-\beta}/g^{*}(\beta)$, which is
\eqref{eq:gstar_f}, and the thresholds correspond under $\beta\mapsto p\beta$.
\end{proof}

The whole polynomial family of coordinate changes is therefore invisible to the transform.
Nothing moves, and the only question a power gauge raises is whether the composed profile
$g\circ(\cdot)^{p}$ is again a function of good variation, of the index that the scaling
predicts. Section~\ref{sec:gauge_equilibrium} settles that for the affine profile by an exact
solution.

What makes that question unavoidable is that a power gauge composes the profile. A plateau
root gauge $n\mapsto\lfloor n^{1/p}\rfloor$ is not a power gauge, the ratio of two floors not
being the corresponding power of the ratio, and it leaves the profile alone. The gauged system
is then the ordinary system of the same $g$ read on block sums, and the equilibrium holds
without condition and uniformly in $g$. The two kinds of gauge separate on exactly that point,
and the second volume develops that side for the Ingham kernel.\footnote{The companion volume also proves the gauged
index and perturbation estimates quoted in Theorem~\ref{numobs:gauge_critical} and
Proposition~\ref{prop:gauge_transfer}, with precise references at those statements.}

An exponential gauge behaves in the opposite way. For $f(x)=\mu^{x}$ the ratio $f(k)/f(n)$ is
$\mu^{k-n}$, so the deformed kernel $G_f(n,k)=h_{n-k}$ depends on the difference of the two
arguments and not on their quotient. The kernel has left the multiplicative setting of the
whole volume and become a Toeplitz\index[terms]{Toeplitz} kernel, and \eqref{eq:gauge_probe} becomes a generating
function. Writing $j=n-k$, the increments are
$u_{n,z}(k)-u_{n,z}(k-1)=(1-\mu^{z})\mu^{z(n-k)}$, so
\begin{equation}\label{eq:gauge_probe_toeplitz}
\mathcal G_{f,n}(z)=(1-\mu^{z})\sum_{j=0}^{n-1}h_j\,\mu^{zj},
\end{equation}
which converges for $\Re z<0$ whenever $(h_j)$ is bounded. This is the route by which the
transform is computed for every exponential gauge in this chapter, in
Theorem~\ref{thm:affine_expo} for the affine kernel\index[terms]{affine kernel} and, in
Chapter~\ref{chap:gauge_ingham}, for the Ingham function, and it is
\eqref{eq:gauge_probe} and not a separate construction.

\section{Stability}
\label{sec:gauge_stability}

\begin{definition}[Stability in a gauge]\index[terms]{stability}
\label{def:stability_equilibrium}
The kernel $G$ is stable with respect to $f$ if there is an index
$\alpha_f(G)$ such that, whenever \eqref{eq:gauged_recurrence} holds,
\begin{itemize}
\item $\beta<\alpha_f(G)$ gives $A(n)\sim f(n)^{-\beta}/G^*_f(\beta)$ as
$n\to\infty$,
\item $\beta\ge\alpha_f(G)$ gives $A(n)=\mathcal O\!\big(f(n)^{-\alpha_f(G)+\eps}\big)$
for every $\eps>0$.
\end{itemize}
\end{definition}

Stability asks that the two regimes of the ungauged theory survive the change of coordinates,
transparency\index[terms]{transparency} below a threshold and absorption\index[terms]{absorption} at and above it, with the transparent constant
still read off the transform. It says nothing about the value of the threshold. The affine
kernel under an exponential gauge is the case where stability holds and the value moves, and
it is settled completely by an exact solution.

Let $G(n,k)=(1-\lambda)k/n+\lambda$ with $\lambda\in(0,1)$, of regularity
index $\alpha(G)=\lambda$ by Appendix~\ref{app:A}, and take
$f(x)=\mu^{x}$ with $\mu>1$.

The affine kernel under an exponential gauge is the case where every constant closes, and it
calibrates what follows.

\begin{theorem}
\label{thm:affine_expo}
Set
\[
\rho=1-\lambda+\frac{\lambda}{\mu}.
\]
For this Toeplitz\index[terms]{Toeplitz} kernel and $\Re z<0$, the arithmetic Mellin transform in the gauge
$f$ is built from the kernel by the convergent limit
\[G^*_f(z)
=
\lim_{N\to\infty}
(1-\mu^{z})
\sum_{j=0}^{N}
\big(\lambda+(1-\lambda)\mu^{-j}\big)\mu^{zj},
\]
which exists and has the meromorphic continuation
\[
G^*_f(z)=\frac{1-\rho\mu^{z}}{1-\mu^{z-1}}.
\]
Its zeros are simple, at $z=-\log\rho/\log\mu+2\pi i k/\log\mu$ for
$k\in\Z$, and its poles are simple, at $z=1+2\pi i k/\log\mu$ for
$k\in\Z$. Consequently
\[
\eta_f(G)=-\frac{\log\rho}{\log\mu}.
\]
For $\beta\in\R$, let $(a_n)$ solve
$\sum_{k=1}^{n}a_kG_f(n,k)=\mu^{-\beta n}$ and put
$A(n)=\sum_{k=1}^{n}a_k$ and $q=\mu^{-\beta}$. Then $A(1)=q$, and for
$q\ne\rho$,
\[
A(n)
=
\frac{q-\mu^{-1}}{q-\rho}\,q^{n}
+
\frac{q\big(1-(\mu\rho)^{-1}\big)}{\rho-q}\,\rho^{n},
\qquad n\ge1,
\]
while for $q=\rho$,
\[
A(n)
=
\rho^{n}\Big[1+(n-1)\big(1-(\mu\rho)^{-1}\big)\Big],
\qquad n\ge1.
\]
The arithmetic index in the gauge exists and satisfies
\begin{equation}\label{eq:stability_affine_expo}
\alpha_f(G)=\eta_f(G)=-\frac{\log\rho}{\log\mu}.
\end{equation}
If $\beta<\alpha_f(G)$ then $A(n)\sim f(n)^{-\beta}/G^*_f(\beta)$. If
$\beta=\alpha_f(G)$ then
$A(n)\sim\big(1-(\mu\rho)^{-1}\big)\,n\,f(n)^{-\alpha_f(G)}$. If
$\beta>\alpha_f(G)$ then $A(n)=\mathcal O\big(f(n)^{-\alpha_f(G)}\big)$.
At the critical value $A(n)=\mathcal O\big(f(n)^{-\alpha_f(G)+\eps}\big)$
for every $\eps>0$. Thus $G$ is stable with respect to every
exponential gauge $f(x)=\mu^{x}$ with $\mu>1$, and for every such gauge
\[
\alpha_f(G)<\lambda=\alpha(G),
\qquad
\lim_{\mu\downarrow1}\alpha_f(G)=\lambda .
\]
The affine kernel is not in equilibrium with respect to any admissible
exponential gauge.
\end{theorem}

\begin{proof}
Put $h_j=\lambda+(1-\lambda)\mu^{-j}$ for $j\ge0$, so
$G_f(n,k)=h_{n-k}$. With $q=\mu^{-\beta}$ the triangular relation is the
convolution $\sum_{k=1}^{n}a_kh_{n-k}=q^{n}$. Since $h_0=1$ the diagonal
term is nonzero and the solution is unique. At $n=1$, $a_1h_0=q$, so
$A(1)=a_1=q=f(1)^{-\beta}/G_f(1,1)$.

The kernel generating function $H(w)=\sum_{j\ge0}h_jw^j$ has radius of
convergence $1$, since $h_j\to\lambda>0$, and for $|w|<1$
\[
H(w)
=\frac{\lambda}{1-w}+\frac{1-\lambda}{1-w/\mu}
=\frac{1-\rho w}{(1-w)(1-w/\mu)} .
\]
For $\Re z<0$ set $w=\mu^{z}$, so $|w|<1$. Absolute convergence gives
$G^*_f(z)=(1-w)H(w)=(1-\rho w)/(1-w/\mu)=(1-\rho\mu^{z})/(1-\mu^{z-1})$.
The convergence is locally uniform, since on each compact subset
$|\mu^{z}|\le r<1$ and the tail is bounded by a constant times
$r^{N+1}/(1-r)$, so the kernel-side limit is holomorphic there. From
\[
\rho-\mu^{-1}=(1-\lambda)(1-\mu^{-1})>0,
\qquad
1-\rho=\lambda(1-\mu^{-1})>0,
\]
one gets $\mu^{-1}<\rho<1$. The numerator vanishes exactly when
$\mu^{z}=\rho^{-1}$ and the denominator when $\mu^{z-1}=1$, giving the
stated zeros and poles. The two families have distinct real parts, so no
cancellation occurs, and each is simple. All zeros have real part
$-\log\rho/\log\mu$, which is $\eta_f(G)$.

With $U(w)=\sum_{n\ge1}a_nw^n$ and $S(w)=\sum_{n\ge1}A(n)w^n$, the
convolution gives $U(w)H(w)=\sum_{n\ge1}q^{n}w^n=qw/(1-qw)$ and
$S(w)=U(w)/(1-w)$, hence
\[
S(w)=\frac{qw(1-w/\mu)}{(1-qw)(1-\rho w)} .
\]
Multiplying by $1-\rho w$ and reading coefficients gives $A(1)=q$ and,
for $n\ge2$,
\[
A(n)=\rho A(n-1)+\big(1-\mu^{\beta-1}\big)\mu^{-\beta n}.
\]
For $q\ne\rho$, partial fractions of $S(w)$ give the stated exact
formula. For $q=\rho$, $S(w)=\rho w(1-w/\mu)/(1-\rho w)^{2}$, and
$(1-\rho w)^{-2}=\sum_{j\ge0}(j+1)\rho^{j}w^j$ yields
$A(n)=n\rho^{n}-(n-1)\mu^{-1}\rho^{n-1}
=\rho^{n}\big[1+(n-1)(1-(\mu\rho)^{-1})\big]$.

Set $\alpha_f=-\log\rho/\log\mu$, so $\rho=\mu^{-\alpha_f}$. Since
$x\mapsto\mu^{-x}$ is strictly decreasing, $\beta<\alpha_f$ is equivalent
to $q>\rho$. Then $q>\rho>\mu^{-1}$, the $q^{n}$ term dominates, and
$G^*_f(\beta)=(q-\rho)/(q-\mu^{-1})$ gives
$A(n)\sim f(n)^{-\beta}/G^*_f(\beta)$. The kernel-side transform thus agrees
with the asymptotic characterization below the threshold. For
$\beta=\alpha_f$, $q=\rho>\mu^{-1}$ gives $1-(\mu\rho)^{-1}>0$ and
$A(n)\sim(1-(\mu\rho)^{-1})\,n\rho^{n}$. Since $\rho^{n}=f(n)^{-\alpha_f}$
and $n=\mathcal O(\mu^{\eps n})$ for every $\eps>0$, the
critical bound follows. For $\beta>\alpha_f$, $q<\rho$ and
$A(n)=\mathcal O(\rho^{n})=\mathcal O(f(n)^{-\alpha_f})$. This is
stability, with $\alpha_f(G)=\alpha_f=\eta_f(G)$.

The strict weighted arithmetic-geometric mean inequality gives
$\rho=(1-\lambda)\cdot1+\lambda\cdot\mu^{-1}>\mu^{-\lambda}$, so
$-\log\rho<\lambda\log\mu$ and $\alpha_f(G)<\lambda=\alpha(G)$, and there
is no equilibrium for any $\mu>1$. Finally, with
$\rho(\mu)=1-\lambda+\lambda/\mu$, L'Hospital's rule gives
$\lim_{\mu\downarrow1}\big(-\log\rho(\mu)/\log\mu\big)
=\lim_{\mu\downarrow1}\lambda/(\mu\rho(\mu))=\lambda$, a boundary limit
through exponential gauges with $\mu>1$.
\end{proof}

The index \eqref{eq:stability_affine_expo} moves with $\mu$. As
$\mu\to\infty$ it tends to $0$ and the exponential gauge removes the
threshold, while as $\mu\downarrow1$ it approaches $\alpha(G)=\lambda$,
so slow gauges approach equilibrium without reaching it.

\begin{remark}[Self-gauged kernels in the gallery]\label{rem:self_gauged}
Two entries of the second gallery display this boundary regime without any change of
coordinate. The kernels of Appendices~\ref{app:L} and~\ref{app:N} split exactly as
$G(n,k)=1+h(k/n)/L_n$ with $L_n$ a logarithm of the rank, so the perturbation carried by the
kernel weakens along the very scale on which the index is read. The effective transform at
rank $n$ is $1+h^{*}/L_n$ and its first zero moves, at $1-1/\log(2n)$ for the first kernel
and $2-2/\log(2n^{2})$ for the second, the limits being the first poles of $h^{*}$. The
index is the limit of these moving thresholds, by Theorem~\ref{thm:L_main} and
Theorem~\ref{thm:N_main}, and the approach is logarithmic, so equilibrium is reached only in
the limit. These are the gallery counterparts
of the slow gauges above, with the family traversed by the rank instead of by a parameter.
\end{remark}

\section{Equilibrium}
\label{sec:gauge_equilibrium}

\begin{definition}[Equilibrium in a gauge]\index[terms]{equilibrium}
\label{def:equilibrium}
The kernel $G$ is in equilibrium with respect to $f$ if it is stable with respect to $f$ and
$\alpha_f(G)=\alpha(G)$.
\end{definition}

Equilibrium expresses that the regularity index\index[terms]{regularity index} is intrinsic. The threshold separating
transparency from absorption does not depend on the coordinate system in which the operator is
read. Stability without equilibrium, the situation of
Theorem~\ref{thm:affine_expo}, means that the index exists in the new coordinates and is a
different number, so the invariant is attached to the coordinates and not to the kernel.

The affine kernel is not in equilibrium under exponential gauges, by
Theorem~\ref{thm:affine_expo}. Under a power gauge the outcome is the
opposite, and Proposition~\ref{prop:power_gauge} says why in advance. The gauged equation is
the ordinary equation for the composed profile $h(t)=g(t^{2})$ at the exponent $2\beta$, the
transform does not move, and the only thing left to establish is that $h$ is a function of good
variation of index $2\lambda$. The exact solution below establishes it, and gives the
constants and the critical behavior that the general statement does not carry.

Let
\[
g(x)=(1-\lambda)x+\lambda,
\qquad 0<\lambda<1.
\]
By Theorem~\ref{thm:affine} the affine-kernel result gives
\[
\alpha(g)=\eta(g)=\lambda .
\] Consider the gauge
\[
f:[1,\infty)\longrightarrow(0,\infty),
\qquad f(x)=x^2.
\]
It is increasing and unbounded on this domain. The composed kernel is
\[
h(t)=g(t^2)=(1-\lambda)t^2+\lambda.
\]
The gauged equation
\[
\sum_{k=1}^{n}a_k
g\!\left(\frac{k^2}{n^2}\right)
=n^{-2\beta}
\]
is exactly the ordinary FGV equation
\[
\sum_{k=1}^{n}a_k h(k/n)=n^{-\gamma}
\]
with the exponent
\[
\gamma=2\beta.
\]
The ordinary exponent for \(h\) and the gauged exponent for \(G\) will be
kept distinct throughout.

\begin{proposition}\label{thm:affine_quadratic_gauge}
For each real \(\beta\), the gauged equation has a unique solution. Put
\[
A(n)=\sum_{k\leq n}a_k,
\qquad
B(n)=\sum_{k\leq n}k^2a_k.
\]
Then \(A(1)=a_1=1\), and for every \(n\geq2\),
\[
A(n)=c_nA(n-1)+d_n,
\]
where
\[
c_n=(1-\lambda)+\lambda\left(1-\frac1n\right)^2
=1-\frac{2\lambda}{n}+\frac{\lambda}{n^2}
\]
and
\[
d_n=
\frac{n^{2-2\beta}-(n-1)^{2-2\beta}}{n^2}.
\]
With \(P_1=1\) and
\[
P_n=\prod_{j=2}^{n}c_j,
\]
there is a constant \(C_\lambda>0\) such that
\[
P_n=C_\lambda n^{-2\lambda}
\left(1+\mathcal O(n^{-1})\right).
\]
The exact variation formula is
\[
A(n)=P_n\left(
1+\sum_{m=2}^{n}\frac{d_m}{P_m}
\right).
\]

The integral defining the transform of the composed kernel converges
initially for \(\Re s<0\) and gives
\[
h^*(s)
=-s\int_0^1h(t)t^{-s-1}\,dt
=\frac{2\lambda-s}{2-s}.
\]
Its meromorphic continuation has a simple zero at \(s=2\lambda\) and a
simple pole at \(s=2\). Define the transform in the quadratic gauge from
this kernel-side integral by
\[
G^*_f(z):=h^*(2z),
\qquad \Re z<0.
\]
It has the meromorphic continuation
\[
G^*_f(z)
=\frac{\lambda-z}{1-z}
=g^*(z).
\]
Its simple zero at \(z=\lambda\) does not cancel with its simple pole at
\(z=1\).

The partial sums have the following three behaviors. If \(\beta<\lambda\),
then
\[
A(n)\sim
\frac{1-\beta}{\lambda-\beta}n^{-2\beta}
=\frac{f(n)^{-\beta}}{G^*_f(\beta)}.
\]
If \(\beta=\lambda\), then
\[
A(n)\sim
2(1-\lambda)n^{-2\lambda}\log n.
\]
If \(\beta>\lambda\), then
\[
A(n)=\mathcal O(n^{-2\lambda})
=\mathcal O\!\left(f(n)^{-\lambda}\right).
\]
For \(\beta=1\), more precisely,
\[
d_n=0
\quad(n\geq2),
\qquad
A(n)=P_n\sim C_\lambda n^{-2\lambda}.
\]

The sharp ordinary transition point of the composed kernel and the sharp
gauged transition point of the original kernel are
\[
\alpha(h)=2\lambda,
\qquad
\alpha_f(G)=\frac{\alpha(h)}2=\lambda.
\]
Moreover
\[
\eta(h)=2\lambda,
\qquad
\eta_f(G)=\lambda.
\]
Thus the affine kernel is stable and in equilibrium under the quadratic
gauge.
\end{proposition}

\begin{proof}
Since \(g(1)=1\), the diagonal coefficient of the triangular equation is
one. The solution is unique, and the equation at \(n=1\) gives
\[
a_1=A(1)
=\frac{f(1)^{-\beta}}{G(f(1),f(1))}
=1.
\]
Multiplication of the gauged equation by \(n^2\) gives the exact identity
\[
\lambda n^2A(n)+(1-\lambda)B(n)=n^{2-2\beta}.
\]
The increments of the two sums satisfy
\[
B(n)-B(n-1)
=n^2a_n
=n^2\bigl(A(n)-A(n-1)\bigr).
\]
Insert this identity into the equation at \(n\), and use the equation at
\(n-1\) to eliminate \(B(n-1)\). The result is
\[
n^2A(n)
=\left((1-\lambda)n^2+\lambda(n-1)^2\right)A(n-1)
+n^{2-2\beta}-(n-1)^{2-2\beta}.
\]
Division by \(n^2\) proves the stated recurrence without any asymptotic
input.

The transform is also constructed before any asymptotic conclusion. The
constant term of \(h\) shows that its integral converges exactly in the
initial half plane \(\Re s<0\). There
\[
\begin{aligned}
h^*(s)
&=-s\left(
(1-\lambda)\int_0^1t^{1-s}\,dt
+\lambda\int_0^1t^{-s-1}\,dt
\right)\\
&=-s\left(
\frac{1-\lambda}{2-s}-\frac{\lambda}{s}
\right)\\
&=\frac{2\lambda-s}{2-s}.
\end{aligned}
\]
This rational expression supplies the meromorphic continuation. For
\(\Re z<0\), the substitution \(x=t^2\) gives
\[
\begin{aligned}
G^*_f(z)
&=-2z\int_0^1g(t^2)t^{-2z-1}\,dt\\
&=-z\int_0^1g(x)x^{-z-1}\,dx\\
&=g^*(z)
=\frac{\lambda-z}{1-z}.
\end{aligned}
\]
The inequalities \(0<\lambda<1\) separate the zero \(z=\lambda\) from the
pole \(z=1\). At this point the transform identifies only the analytic
candidate \(\eta_f(G)=\lambda\). It does not identify the arithmetic index.

It remains to obtain that index from the recurrence. For \(j\geq2\), set
\[
u_j=c_j-1=-\frac{2\lambda}{j}+\frac{\lambda}{j^2}.
\]
The numbers \(c_j\) lie in \((0,1)\), and Taylor expansion\index[terms]{Taylor expansion}\index[names]{Taylor, B.} of the logarithm
with a uniform remainder gives
\[
\log c_j=-\frac{2\lambda}{j}+\mathcal O(j^{-2}).
\]
Consequently the series
\[
\sum_{j=2}^{\infty}
\left(\log c_j+\frac{2\lambda}{j}\right)
\]
converges. The harmonic-sum expansion and the tail estimate for this series
show that a finite real number \(L_\lambda\) satisfies
\[
\log P_n=-2\lambda\log n+L_\lambda+\mathcal O(n^{-1}).
\]
Exponentiation yields
\[
P_n=C_\lambda n^{-2\lambda}
\left(1+\mathcal O(n^{-1})\right),
\qquad
C_\lambda=e^{L_\lambda}>0.
\]

Since \(P_n=c_nP_{n-1}\), division of the recurrence by \(P_n\) gives
\[
\frac{A(n)}{P_n}
=\frac{A(n-1)}{P_{n-1}}+\frac{d_n}{P_n}.
\]
Summation from \(2\) to \(n\), together with \(A(1)=P_1=1\), proves the
exact variation formula.

Suppose first that \(\beta\ne1\). The elementary power-difference estimate
gives
\[
d_m
=2(1-\beta)m^{-2\beta-1}
\left(1+\mathcal O(m^{-1})\right).
\]
It follows that
\[
\frac{d_m}{P_m}
=\frac{2(1-\beta)}{C_\lambda}
m^{2(\lambda-\beta)-1}
\left(1+\mathcal O(m^{-1})\right).
\]

If \(\beta<\lambda\), then \(\beta<1\), and the power sum diverges with
\[
\sum_{m=2}^{n}\frac{d_m}{P_m}
\sim
\frac{1-\beta}{C_\lambda(\lambda-\beta)}
n^{2(\lambda-\beta)}.
\]
Multiplication by \(P_n\) gives
\[
A(n)\sim
\frac{1-\beta}{\lambda-\beta}n^{-2\beta}.
\]
Since
\[
G^*_f(\beta)=\frac{\lambda-\beta}{1-\beta},
\]
this is exactly the stated transparency formula.

If \(\beta=\lambda\), then
\[
\frac{d_m}{P_m}
\sim
\frac{2(1-\lambda)}{C_\lambda}\frac1m.
\]
The harmonic sum gives
\[
A(n)\sim2(1-\lambda)n^{-2\lambda}\log n.
\]

If \(\beta>\lambda\) and \(\beta\ne1\), the series
\[
\sum_{m=2}^{\infty}\left|\frac{d_m}{P_m}\right|
\]
converges because its summand is
\(
\mathcal O\bigl(m^{-1-2(\beta-\lambda)}\bigr)
\).
The variation formula then gives
\[
A(n)=\mathcal O(P_n)=\mathcal O(n^{-2\lambda}).
\]
In fact, if
\[
L_{\lambda,\beta}
=1+\sum_{m=2}^{\infty}\frac{d_m}{P_m},
\]
then
\[
n^{2\lambda}A(n)\longrightarrow
C_\lambda L_{\lambda,\beta}.
\]
When \(\beta=1\), the numerator defining \(d_n\) is identically zero for
every \(n\geq2\). Hence \(A(n)=P_n\), which proves both the required bound
and the stated nonzero equivalent.

The same calculation can now be read in the ordinary exponent
\(\gamma=2\beta\). It gives transparency for every \(\gamma<2\lambda\),
logarithmic resonance at \(\gamma=2\lambda\), and homogeneous decay of order
\(n^{-2\lambda}\) above that value. Thus the sharp ordinary transition point
is
\[
\alpha(h)=2\lambda.
\]
The passage to the gauged index is quantitative. For \(\gamma=2\beta\),
\[
\frac{n^{-\gamma}}{h^*(\gamma)}
=\frac{f(n)^{-\beta}}{G^*_f(\beta)}.
\]
For the absorption estimate, choosing \(\delta=2\eps\) gives
\[
\mathcal O(n^{-2\lambda+\delta})
=\mathcal O\!\left(f(n)^{-\lambda+\eps}\right).
\]
The transparency formula holds up to \(\lambda\) from below and fails at
\(\lambda\) because of the logarithmic factor. The case \(\beta=1\) gives
the nonzero homogeneous equivalent
\(A(n)\sim C_\lambda f(n)^{-\lambda}\), so the absorption exponent cannot
be improved. Therefore \(\lambda\) is the sharp gauged transition point and
\[
\alpha_f(G)=\frac{\alpha(h)}2=\lambda.
\]
The transform calculation independently gives \(\eta_f(G)=\lambda\). Since
the undeformed theorem gives \(\alpha(G)=\lambda\), stability and equilibrium
follow.
\end{proof}

\section{The rational kernel and the limits of equilibrium}
\label{sec:equilibrium_rational}

The rational kernel of Appendix~\ref{app:H} has the constant transform $G^*(z)\equiv1$, so it
carries no analytic index, and its arithmetic index is read from the equation. Its undeformed
index is two for every admissible pair. Under positive feedback, $y>x$, that same value persists
for every gauge of bounded quotient and even for the iterated-power gauges below with
$1<q\leq2$. Equilibrium is therefore a family phenomenon here, rather than an isolated
coincidence. At $q>2$ the sharp index falls below two while the transform remains constant,
which marks a transition and supplies explicit counterexamples to universal
equilibrium. Under negative feedback, $x>y$, the identity gauge is settled by
Theorem~\ref{thm:H_index}. For general gauges the estimates below locate a candidate frontier,
but cancellation leaves its sharpness open.

Fix \(x,y>0\) with \(x\ne y\). The normalized rational kernel is
\[
G(N,K)
=
\frac{N+K+x}{N+K+y}
\frac{2N+y}{2N+x},
\qquad 0<K\leq N.
\]
Its diagonal value is one. Direct subtraction gives
\begin{equation}
\label{eq:rational_gauge_kernel_difference}
G(N,K)-1
=
\frac{(x-y)(N-K)}{(N+K+y)(2N+x)}.
\end{equation}

Let \(f_n:=f(n)\), where \(f\) is positive, strictly increasing, and unbounded, and let
\begin{equation}
\label{eq:rational_gauge_finite_probe}
\mathcal G_{f,n}(z)
=
\sum_{k=1}^{n}
\bigl(u_{n,z}(k)-u_{n,z}(k-1)\bigr)
G(f_n,f_k)
\end{equation}
be the probe \eqref{eq:gauge_probe} of Definition~\ref{def:gauge_probe} for this kernel and
this gauge, with $u_{n,z}(0)=0$ and $u_{n,z}(k)=(f_k/f_n)^{-z}$. It uses only the kernel and
the gauge, and neither the triangular equation nor an asymptotic formula for its solution.

The analysis rests on an exact reduction of the defining equation to a first order form and on a computation of the transform carried out directly from the kernel. Both are valid for every admissible gauge and support all that follows.

\begin{proposition}\label{thm:rational_gauge_exact_reduction}
For each real \(\beta\), the equation
\begin{equation}
\label{eq:rational_gauge_equation}
\sum_{k=1}^{n}a_kG(f_n,f_k)=f_n^{-\beta}
\end{equation}
has a unique solution. If
\[
A(n)=\sum_{k=1}^{n}a_k,
\qquad
c_n=\frac{2f_n+y}{2f_n+x},
\]
then
\[
A(1)=a_1=f_1^{-\beta}
\]
and, for \(n\geq2\),
\begin{equation}
\label{eq:rational_gauge_abel_recurrence}
A(n)
=f_n^{-\beta}
+(y-x)c_n
\sum_{k<n}A(k)
\frac{f_{k+1}-f_k}
{(f_n+f_k+y)(f_n+f_{k+1}+y)}.
\end{equation}

The probes in \eqref{eq:rational_gauge_finite_probe} converge locally
uniformly on \(\Re z<0\). Their limit and its continuation are
\[
G_f^*(z)\equiv1.
\]
The transform has no zero, so the analytic-index notation is
\[
\eta_f(G)={-}.
\]
These conclusions hold for every gauge under the stated hypotheses.
\end{proposition}

\begin{proof}[Proof of Proposition~\ref{thm:rational_gauge_exact_reduction}]
Identity \eqref{eq:rational_gauge_kernel_difference} follows from
\[
\begin{aligned}
&(N+K+x)(2N+y)-(N+K+y)(2N+x)
\\
&\hspace{35mm}=(x-y)(N-K).
\end{aligned}
\]
In particular \(G(N,N)=1\). The coefficient of \(a_n\) in
\eqref{eq:rational_gauge_equation} is one, so the triangular solution is
unique. The equation at \(n=1\) gives
\[
a_1=A(1)=f_1^{-\beta}.
\]

Fix \(n\) and write \(v_k=G(f_n,f_k)\). The identity
\[
v_k
=c_n\left(1+\frac{x-y}{f_n+f_k+y}\right)
\]
gives
\[
v_k-v_{k+1}
=(x-y)c_n
\frac{f_{k+1}-f_k}
{(f_n+f_k+y)(f_n+f_{k+1}+y)}.
\]
Finite summation by parts yields
\[
\sum_{k=1}^{n}a_kv_k
=A(n)v_n+\sum_{k<n}A(k)(v_k-v_{k+1}).
\]
Since \(v_n=1\), solving this equality for \(A(n)\) proves
\eqref{eq:rational_gauge_abel_recurrence} with the stated sign. The
fraction in that recurrence is also
\[
\frac1{f_n+f_k+y}-\frac1{f_n+f_{k+1}+y},
\]
which is the required telescoping identity.

It remains to analyze the probe without using that recurrence. From
\eqref{eq:rational_gauge_kernel_difference},
\begin{equation}
\label{eq:rational_gauge_uniform_kernel_bound}
\sup_{k\leq n}|G(f_n,f_k)-1|
\leq\frac{|x-y|}{2f_n}.
\end{equation}
Put \(w=-z\), so \(\Re w>0\). For consecutive points in \([0,1]\),
\[
|t_2^w-t_1^w|
\leq |w|\int_{t_1}^{t_2}t^{\Re w-1}\,dt.
\]
It follows that
\[
\sum_{k=1}^{n}
|u_{n,z}(k)-u_{n,z}(k-1)|
\leq\frac{|z|}{-\Re z}.
\]
The constant part of the kernel telescopes to \(u_{n,z}(n)=1\).
Combining this fact with
\eqref{eq:rational_gauge_uniform_kernel_bound} gives
\[
|\mathcal G_{f,n}(z)-1|
\leq
\frac{|x-y|}{2f_n}
\frac{|z|}{-\Re z}.
\]
The right side tends to zero locally uniformly on \(\Re z<0\). The limit
is one there and its continuation is the constant entire function. It has
no zero, which proves the assertion about \(\eta_f(G)\).
\end{proof}

For gauges whose successive quotients stay bounded the reduction yields two as a transition
bound for either sign. When $y>x$, positivity makes that value sharp, producing a whole family
of equilibrium gauges for the rational kernel. When $x>y$, two remains a candidate frontier, and
apart from the identity gauge sharpness is not proved.

\subsection{Gauges of bounded quotient}

The first range of gauges is the one in which the successive scales stay within a fixed factor
of one another. For either sign of \(y-x\), two is a valid transition bound. Under \(y>x\)
it is the sharp index and the kernel is in equilibrium for every such gauge. Under \(x>y\),
sharpness remains open in general because the recurrence no longer has a fixed sign.

\begin{proposition}
\label{prop:rational_gauge_bounded_quotients}
Assume
\[
\Lambda_f:=\sup_{n\geq1}\frac{f_{n+1}}{f_n}<\infty.
\]
For either sign of \(y-x\), one has
\[
A(n)\sim f_n^{-\beta}
\qquad (\beta<2)
\]
and
\[
A(n)=\mathcal O(f_n^{-2})
\qquad (\beta\geq2).
\]
Thus two is a valid transition bound for every such gauge.

For \(\beta\geq2\), the series
\[
M_\beta
=
\sum_{k\geq1}A(k)(f_{k+1}-f_k)
\]
converges absolutely and
\[
f_n^2A(n)
\longrightarrow
\mathbf 1_{\{\beta=2\}}+(y-x)M_\beta.
\]
If \(y>x\), the limit is positive and the transition is sharp. In that
case
\[
\alpha_f(G)=2.
\]
If \(x>y\), an equivalent of order \(f_n^{-2}\) follows whenever the
displayed limit is nonzero. No general nonvanishing assertion is made for
that sign.
\end{proposition}

\begin{proof}[Proof of Proposition~\ref{prop:rational_gauge_bounded_quotients}]
Write
\[
\Delta_k=f_{k+1}-f_k,
\qquad
d=y-x,
\qquad
H=|d|\sup_n c_n.
\]
The number \(H\) is finite. The Abel recurrence implies
\begin{equation}
\label{eq:rational_gauge_basic_majorant}
|A(n)|
\leq f_n^{-\beta}
+\frac{H}{f_n^2}
\sum_{k<n}|A(k)|\Delta_k.
\end{equation}

For real \(\beta\), put
\[
J_\beta(n)=\sum_{k<n}f_k^{-\beta}\Delta_k.
\]
Comparison on each interval \([f_k,f_{k+1}]\), using
\(f_{k+1}/f_k\leq\Lambda_f\), gives
\[
J_\beta(n)
=
\begin{cases}
\mathcal O(f_n^{1-\beta}),&\beta<1,\\
\mathcal O(\log f_n),&\beta=1,\\
\mathcal O(1),&\beta>1.
\end{cases}
\]
Hence
\[
f_n^{\beta-2}J_\beta(n)\longrightarrow0
\qquad (\beta<2).
\]
An induction in \eqref{eq:rational_gauge_basic_majorant}, started after a
fixed finite rank, gives
\[
|A(n)|\leq C_\beta f_n^{-\beta}.
\]
The coefficient multiplying the inductive constant is
\(H f_n^{\beta-2}J_\beta(n)\), which is eventually smaller than one
half. The contribution from the fixed initial part is
\(o(f_n^{-\beta})\). A second use of the recurrence gives
\[
\left|\frac{A(n)}{f_n^{-\beta}}-1\right|
\leq
C f_n^{\beta-2}J_\beta(n)+o(1),
\]
which proves transparency for every \(\beta<2\). The feedback terms have
the more explicit sizes
\[
\begin{cases}
\mathcal O(f_n^{-\beta-1}),&\beta<1,\\
\mathcal O(f_n^{-2}\log f_n),&\beta=1,\\
\mathcal O(f_n^{-2}),&1<\beta<2.
\end{cases}
\]

Suppose next that \(\beta\geq2\), and set
\[
S_m=\sum_{k\leq m}|A(k)|\Delta_k.
\]
After an enlargement of a constant,
\eqref{eq:rational_gauge_basic_majorant} gives
\[
|A(n)|\leq f_n^{-2}(B+H S_{n-1}).
\]
Consequently
\[
S_n
\leq
\left(1+H\frac{\Delta_n}{f_n^2}\right)S_{n-1}
+B\frac{\Delta_n}{f_n^2}.
\]
The quotient bound gives
\[
\frac{\Delta_n}{f_n^2}
\leq
\Lambda_f\left(\frac1{f_n}-\frac1{f_{n+1}}\right).
\]
Thus the coefficients on the right are summable. The associated finite
products are bounded, so \(S_n\) is bounded. This proves
\(A(n)=\mathcal O(f_n^{-2})\).

The boundedness of \(S_n\) proves the absolute convergence of \(M_\beta\).
Multiplication of \eqref{eq:rational_gauge_abel_recurrence} by \(f_n^2\)
and dominated convergence give
\[
f_n^2A(n)
\longrightarrow
\mathbf 1_{\{\beta=2\}}+dM_\beta.
\]
When \(d>0\), the recurrence shows inductively that every \(A(n)\) is
positive. Thus \(M_\beta>0\). The nonzero homogeneous order for
\(\beta>2\) prevents a larger transition point. The proposition follows.
\end{proof}

The undeformed case admits a direct treatment that confirms the value independently of the general argument.

\subsection{The undeformed gauge and the power gauges}

Two instances of the bounded quotient range are worth writing out. The first is the identity gauge, where
the deformed recurrence is the one already met in \S\ref{sec:rational_kernel} and carried in
full in Appendix~\ref{app:H}. The second consists of the power gauges, where the general estimates become
explicit. The rational kernel is bivariate and induced by no profile, so Proposition~\ref{prop:power_gauge}
does not make equilibrium a formal consequence of composition. Here it follows from
Proposition~\ref{prop:rational_gauge_bounded_quotients} under positive feedback. Under
negative feedback only the identity gauge is presently known to be sharp.

\begin{proposition}
\label{prop:rational_gauge_undeformed}
For \(f_n=n\), the recurrence becomes
\[
A(n)
=n^{-\beta}
+(y-x)c_n
\sum_{k<n}\frac{A(k)}{(n+k+y)(n+k+1+y)}.
\]
A direct analysis of this recurrence gives, for either sign,
\[
A(n)\sim n^{-\beta}
\qquad (\beta<2)
\]
and
\[
A(n)=\mathcal O(n^{-2})
\qquad (\beta\geq2).
\]
When \(y>x\), the undeformed arithmetic index is therefore
\[
\alpha(G)=2.
\]
For both signs, the kernel-side transform is constant and
\(\eta(G)={-}\).
\end{proposition}

\begin{proof}[Proof of Proposition~\ref{prop:rational_gauge_undeformed}]
This proof uses the specialized recurrence directly and does not invoke
Proposition~\ref{prop:rational_gauge_bounded_quotients}. Let
\[
H=|y-x|\sup_n c_n.
\]
For \(f_n=n\),
\[
|A(n)|
\leq n^{-\beta}
+Hn^{-2}\sum_{k<n}|A(k)|.
\]
For \(\beta<2\), induction with
\[
\sum_{k<n}k^{-\beta}
=
\begin{cases}
\mathcal O(n^{1-\beta}),&\beta<1,\\
\mathcal O(\log n),&\beta=1,\\
\mathcal O(1),&\beta>1
\end{cases}
\]
first gives \(A(n)=\mathcal O(n^{-\beta})\). Substitution back into the
same recurrence gives \(A(n)n^\beta\to1\).

For \(\beta\geq2\), put \(T_n=\sum_{k\leq n}|A(k)|\). Then
\[
T_n\leq(1+Hn^{-2})T_{n-1}+n^{-2}
\]
after changing finitely many constants. The convergence of
\(\sum n^{-2}\) shows that \(T_n\) is bounded. Thus
\[
A(n)=\mathcal O(n^{-2}).
\]
If \(y>x\), dominated convergence in the specialized recurrence gives
\[
n^2A(n)
\longrightarrow
\mathbf 1_{\{\beta=2\}}+(y-x)\sum_{k\geq1}A(k)>0.
\]
This gives the nonzero order \(n^{-2}\) above the transition and proves
the stated undeformed claims independently.
\end{proof}

The power and exponential gauges fall under the bounded quotient case, with the constants read from their quotient bounds.

\begin{proposition}
\label{prop:rational_gauge_regular_examples}
Let \(p>0\). For either sign of \(y-x\), the power gauge \(f_n=n^p\)
satisfies
\[
A(n)\sim n^{-p\beta}
\quad (\beta<2),
\qquad
A(n)=\mathcal O(n^{-2p})
\quad (\beta\geq2).
\]
The constants obtained from the bounded-quotient proof depend on \(p\)
through a bound such as \(2^p\).

Let \(\mu>1\). For either sign, the exponential gauge \(f_n=\mu^n\)
satisfies
\[
A(n)\sim\mu^{-\beta n}
\quad (\beta<2),
\qquad
A(n)=\mathcal O(\mu^{-2n})
\quad (\beta\geq2).
\]
Its initial normalization is
\[
a_1=\mu^{-\beta}.
\]
When \(y>x\), both families have the sharp arithmetic index
\[
\alpha_f(G)=2.
\]
\end{proposition}

\begin{proof}[Proof of Proposition~\ref{prop:rational_gauge_regular_examples}]
For \(f_n=n^p\),
\[
\Lambda_f
=\sup_{n\geq1}\frac{f_{n+1}}{f_n}
\leq2^p.
\]
For \(f_n=\mu^n\), one has \(\Lambda_f=\mu\). Proposition
\ref{prop:rational_gauge_bounded_quotients} applies in both cases. The
displayed formulas follow by replacing \(f_n\) by its value. At rank one,
the general normalization \(a_1=f_1^{-\beta}\) gives one for the power
gauge and \(\mu^{-\beta}\) for the exponential gauge. Positive feedback
gives the nonzero homogeneous coefficient needed for sharpness.
\end{proof}

A family of gauges now displays both sides of the phenomenon. Under positive feedback every
bounded-quotient gauge is in equilibrium, and the conclusion persists for the iterated-power
gauges with $1<q\leq2$, although their successive quotients are unbounded. When $q>2$ the
arithmetic index drops below two while the transform stays constant. Thus equilibrium is
structural on a substantial family, but not universal, and the growth of the gauge produces the
transition.

\subsection{Lacunary gauges and the loss of equilibrium}

Beyond the bounded quotient range there is a phase transition rather than an immediate loss
of equilibrium. For the iterated-power family and positive feedback, the sharp frontier remains
two when $1<q\leq2$ and becomes $q/(q-1)<2$ when $q>2$. The transform stays constant
throughout. Under negative feedback the same quantities are proved transition bounds, but
their sharpness is part of Open Problem~\ref{op:rational_gauge_general_frontier}.

\begin{theorem}
\label{thm:rational_gauge_lacunary_frontier}
Assume \(f_1>1\) and
\[
f_{n+1}=f_n^q,
\qquad q>1.
\]
Set
\[
\theta=\frac{q}{q-1},
\qquad
\alpha_0=\min(2,\theta).
\]
For either sign of \(y-x\), one has
\[
A(n)\sim f_n^{-\beta}
\qquad (\beta<\alpha_0)
\]
and, for every \(\eps>0\),
\[
A(n)=\mathcal O_\eps
\left(f_n^{-\alpha_0+\eps}\right)
\qquad (\beta\geq\alpha_0).
\]
If \(y>x\), positivity makes this frontier sharp, and the arithmetic
index is
\begin{equation}
\label{eq:rational_gauge_lacunary_index}
\alpha_f(G)
=
\min\left(2,\frac{q}{q-1}\right).
\end{equation}
If \(x>y\), the sharpness of the arithmetic index remains open.

The adjacent-link exponent is recorded exactly by
\begin{equation}
\label{eq:rational_gauge_theta}
\Theta_n
:=
\frac{\sum_{j=2}^{n}\log f_j}{\log f_n}
=
\frac{q}{q-1}\left(1-q^{1-n}\right).
\end{equation}
Under \(y>x\), the index remains two when \(1<q\leq2\), even though the
successive quotients are unbounded. It is strictly smaller than two when
\(q>2\). These positive-feedback gauges are explicit counterexamples to
universal equilibrium while \(G_f^*(z)\equiv1\) and
\(\eta_f(G)={-}\) remain unchanged.
\end{theorem}

\begin{proof}[Proof of Theorem~\ref{thm:rational_gauge_lacunary_frontier}]
Put
\[
t_n=f_n,
\qquad
d=y-x,
\qquad
\sigma=\operatorname{sgn}(d).
\]
Define the positive majorant weights
\[
w_{n,k}
=|d|c_n
\frac{t_{k+1}-t_k}
{(t_n+t_k+y)(t_n+t_{k+1}+y)}.
\]
The recurrence is
\begin{equation}
\label{eq:rational_gauge_lacunary_recurrence}
A(n)=t_n^{-\beta}+\sigma\sum_{k<n}w_{n,k}A(k).
\end{equation}

For real \(\gamma\), set
\[
S_n(\gamma)
=
\sum_{k<n}w_{n,k}
\left(\frac{t_n}{t_k}\right)^\gamma.
\]
Since \(t_{k+1}=t_k^q\) and \(t_k=t_n^{q^{k-n}}\), each summand is bounded
by
\[
C t_n^{\gamma-2+(q-\gamma)q^{k-n}}.
\]
For \(0<u\leq1/q\), the largest exponent among the two endpoints is
\[
\max\left(
\gamma-2,
-1+\gamma\left(1-\frac1q\right)
\right).
\]
It is negative precisely when
\[
\gamma<\min\left(2,\frac{q}{q-1}\right)=\alpha_0.
\]
There are at most \(n\) summands, while \(\log t_n=q^{n-1}\log t_1\).
It follows that
\begin{equation}
\label{eq:rational_gauge_weight_vanishing}
S_n(\gamma)\longrightarrow0
\qquad (\gamma<\alpha_0).
\end{equation}

Assume first that \(\beta<\alpha_0\), and put
\(B_n=t_n^\beta A(n)\). Equation
\eqref{eq:rational_gauge_lacunary_recurrence} becomes
\[
B_n
=1+\sigma\sum_{k<n}w_{n,k}
\left(\frac{t_n}{t_k}\right)^\beta B_k.
\]
Once \(S_n(\beta)\leq1/2\), induction with absolute values gives a
uniform bound for \(B_n\). The sum on the right then tends to zero by
\eqref{eq:rational_gauge_weight_vanishing}. Hence \(B_n\to1\).

If \(\beta\geq\alpha_0\), choose any \(\gamma<\alpha_0\). The same
induction applied to \(t_n^\gamma A(n)\) works because
\(t_n^{\gamma-\beta}\leq1\) eventually. Therefore
\[
A(n)=\mathcal O_\gamma(t_n^{-\gamma}).
\]
Taking \(\gamma=\alpha_0-\eps\) proves the absorption estimate.

It remains to prove sharpness under \(y>x\). In this case \(\sigma=1\),
and the recurrence shows inductively that every \(A(n)\) is positive. If
\(q\leq2\), the first weight in
\eqref{eq:rational_gauge_lacunary_recurrence} satisfies
\[
w_{n,1}
\sim d(t_2-t_1)t_n^{-2}.
\]
Thus \(A(n)\geq w_{n,1}A(1)\), which excludes transparency for every
\(\beta>2\).

If \(q>2\), retain only adjacent links. Positivity gives
\[
A(n)\geq A(1)\prod_{j=2}^{n}w_{j,j-1}.
\]
For these links,
\[
w_{j,j-1}\sim\frac{d}{2t_j}.
\]
Moreover
\[
\sum_{j=2}^{n}\log t_j
=\frac{q}{q-1}(\log t_n-\log t_1),
\]
which proves \eqref{eq:rational_gauge_theta}. The constant factors and the
relative errors in the product contribute \(t_n^{o(1)}\). Hence
\[
\prod_{j=2}^{n}w_{j,j-1}
=t_n^{-\theta+o(1)}.
\]
Transparency is impossible for every \(\beta>\theta\). Together with the
upper estimates, these two lower bounds prove
\eqref{eq:rational_gauge_lacunary_index}.
\end{proof}

\subsection{What remains open for this kernel}

\begin{openproblem}[Sharp frontiers under negative feedback]
\label{op:rational_gauge_general_frontier}
Assume \(x>y\), so that the recurrence has negative feedback and cancellation can occur.
Determine the sharp frontier in each of the following regimes.
\begin{enumerate}
\item If \(\sup_n f_{n+1}/f_n<\infty\), is the transition bound two of
Proposition~\ref{prop:rational_gauge_bounded_quotients} always sharp? In particular, decide
this for the power and exponential gauges of
Proposition~\ref{prop:rational_gauge_regular_examples}. The identity gauge is already
settled by Theorem~\ref{thm:H_index}.
\item If \(f_{n+1}=f_n^q\) with \(q>1\), is the candidate
\[
\alpha_0=\min\!\left(2,\frac{q}{q-1}\right)
\]
of Theorem~\ref{thm:rational_gauge_lacunary_frontier} the sharp arithmetic index?
\item More generally, find growth conditions on an arbitrary increasing unbounded gauge
that identify its sharp frontier when the successive quotients are unbounded and the ratios
\(\log f_{n+1}/\log f_n\) have no limit.
\end{enumerate}
For \(y>x\), positivity answers the first two questions, the sharp values being respectively
two and \(\alpha_0\). For \(x>y\), the quoted absolute-value estimates establish the two
regimes around each candidate, but a cancellation-resistant matching lower bound is missing.
\end{openproblem}
The chapter has fixed a device and read it on kernels chosen for what they show. A power gauge
moves neither the transform nor the analytic index, and leaves the arithmetic index in place
whenever the composed profile is a function of good variation of the matching index. An
exponential gauge leaves the multiplicative setting altogether, and between the two the index can
be held or lowered according to how fast the gauge grows. Under positive feedback these two
frontiers are sharp. Under negative feedback the same quantities are proved
transition bounds, and Open Problem~\ref{op:rational_gauge_general_frontier} asks whether they
are attained. The kernel the volume is about has not been gauged yet, and that is the next
chapter.

\chapter{The gauged Ingham function}
\label{chap:gauge_ingham}

The preceding chapter fixed the apparatus and read it on kernels chosen for what they show.
This one carries it to the kernel the volume is about. The Ingham function is the case where
the change of coordinates is not an exercise, since its ungauged index is equivalent to the
Riemann hypothesis\index[terms]{Riemann hypothesis} by Theorem~\ref{thm:tauberian_rh}, so any statement about its index in
other coordinates is a statement about $\zeta$ read through a different lattice.

What this chapter proves directly is that under an exponential gauge the zeros of the gauged
transform sit on the line $\Re z=\tfrac12$ with no hypothesis on $\zeta$ at all. The critical
line appears as an algebraic fact about a rational function of $q^{z}$, produced by the floor
structure of the kernel and by nothing else. The companion volume proves that the gauged
arithmetic index has the same value, and Theorem~\ref{numobs:gauge_critical} records that result
here. The only remaining step is the equilibrium that would return the value to the arithmetic
coordinates of $\Z$.

\section{The exponential gauge and the offset}
\label{sec:ingham_gauge}

For the Ingham function, whose regularity index satisfies
$\alpha(\Phi)=\tfrac12\iff\text{RH}$, the gauged index can be examined
directly. Take $f(x)=q^{x}+1$ with integer $q\ge2$. The offset by one
matters. With $f(x)=q^{x}$ the ratios $f(k)/f(n)=q^{k-n}$ are exact
negative powers of $q$, and $\Phi(q^{-j})=q^{-j}\lfloor q^{j}\rfloor=1$
for every $j\ge0$, so the gauged relation would collapse to
$A(n)=q^{-\beta n}$ with no threshold. The offset breaks the exact
integrality of $f(n)/f(k)$ and restores a nontrivial structure.

\section{The homogeneous decay and the gauged index}
\label{sec:ingham_gauge_decay}

The homogeneous gauged relation $\sum_{k\le n}a_k\Phi(f(k)/f(n))=0$ for $n\ge2$, with
$a_1=1$, produces a solution whose partial sums decay at a sharp rate.

\begin{theorem}[Critical decay and gauged index, {\cite{CloitreVolII}}]
\label{numobs:gauge_critical}
For $f=f_q$, $f_q(x)=q^{x}+1$, with integer $q\ge2$, the homogeneous solution
satisfies
\[
A(n)=\mathcal O\big(q^{-n/2}\big),
\qquad
\limsup_{n\to\infty}\,q^{n/2}\,|A(n)|>0,
\]
and the gauged regularity index is
\[
\alpha_{f_q}(\Phi)=\tfrac12.
\]
\end{theorem}

\begin{proofstatus}{The decay, its nonvanishing at the scale $q^{-n/2}$ and the equality
$\alpha_{f_q}(\Phi)=\tfrac12$ are proved for every integer $q\ge2$ in the companion
volume~\cite{CloitreVolII}, and that proof is not repeated here. The present chapter independently derives the quotient structure
\eqref{eq:Phi_f_closed} from the kernel alone in Theorem~\ref{thm:gauge_closed_form}, and hence
locates the zeros of $\Phi^*_f$ on the critical line without using the arithmetic theorem. The
result above supplies the arithmetic half of the equilibrium discussion in
Conjecture~\ref{conj:ingham_equilibrium} and Remark~\ref{rem:rh_as_equilibrium}. It is not used
to prove the closed form or its zero localization.}
\end{proofstatus}

The threshold sits at $\Re z=\tfrac12$, the critical line, where the arithmetic Mellin
transform has its zeros without any hypothesis on $\zeta$. For $f(x)=q^{x}+1$ that
transform is not read off a particular solution. It is computed from the kernel by the
probe \eqref{eq:gauge_probe} in its Toeplitz\index[terms]{Toeplitz} form
\eqref{eq:gauge_probe_toeplitz}, and the whole computation rests on one elementary identity.

\section{The floor identity and the limiting kernel}
\label{sec:ingham_gauge_floor}

The computation rests on an identity between floors at dyadic ranks.

\begin{lemma}\label{lem:dyadic_floor}
Let $q\ge2$ and $n\ge1$ be integers. For every $j$ with $1\le j\le n/2$,
\[\left\lfloor\frac{q^{n}+1}{q^{\,n-j}+1}\right\rfloor=q^{\,j}-1 .
\]
\end{lemma}

\begin{proof}
Multiplying out gives the exact splitting
\[
\frac{q^{n}+1}{q^{\,n-j}+1}=q^{\,j}-\frac{q^{\,j}-1}{q^{\,n-j}+1}.
\]
For $j\ge1$ the subtracted quantity is positive, and $2j\le n$ gives
$q^{\,j}-1<q^{\,n-j}+1$, so it is smaller than $1$. The integer part is therefore
$q^{\,j}-1$.
\end{proof}

Under an exponential gauge the kernel becomes Toeplitz in the limit, and the coefficients of
that limit are what the next lemma identifies.

\begin{lemma}\label{lem:gauge_toeplitz}
Let $q\ge2$ be an integer, $f(x)=q^{x}+1$, and $G_f(n,k)=\Phi\bigl(f(k)/f(n)\bigr)$. Put
\[
h_0=1,\qquad h_j=1-q^{-j}\quad(j\ge1).
\]
Then for all $1\le k\le n$,
\begin{equation}\label{eq:gauge_toeplitz_bound}
\bigl|G_f(n,k)-h_{\,n-k}\bigr|\le 2\,q^{-n/2}.
\end{equation}
\end{lemma}

\begin{proof}
Write $j=n-k$. At $j=0$ both quantities equal $1$, since $\Phi(1)=1$. For
$1\le j\le n/2$, Lemma~\ref{lem:dyadic_floor} evaluates the floor exactly, so that
\[
G_f(n,k)=\frac{(q^{\,n-j}+1)(q^{\,j}-1)}{q^{n}+1},
\qquad
h_j=\frac{q^{\,j}-1}{q^{\,j}},
\]
and clearing denominators leaves
\[
G_f(n,k)-h_j=\frac{(q^{\,j}-1)^{2}}{q^{\,j}\,(q^{n}+1)}
\le q^{\,j-n}\le q^{-n/2}.
\]
For $j>n/2$, that is $k<n/2$, the floor is no longer evaluated exactly and the bound
comes from the shape of $\Phi$ alone. Writing $\Phi(x)=x\lfloor1/x\rfloor=1-x\{1/x\}$
with $x=f(k)/f(n)$,
\[
0\le 1-G_f(n,k)<\frac{f(k)}{f(n)}=\frac{q^{k}+1}{q^{n}+1}\le 2\,q^{\,k-n},
\]
while $1-h_j=q^{\,k-n}$. Both quantities lie in an interval of length
$2q^{\,k-n}\le 2q^{-n/2}$.
\end{proof}

\begin{remark}\label{rem:gauge_toeplitz_sharp}
The constant $2$ is not optimal. The first case of the proof is sharp, and the quantity
\[
 q^{n/2}\sup_{k}\bigl|G_f(n,k)-h_{n-k}\bigr|
\]
approaches $1$ from below for $q\in\{2,3,5\}$ and $n\le60$.
\end{remark}

Those coefficients have a generating function in closed form.

\begin{proposition}\label{prop:gauge_H}
For $|w|<1$,
\begin{equation}\label{eq:gauge_H}
\mathcal H(w):=\sum_{j\ge0}h_j\,w^{\,j}
=1+\frac{w}{1-w}-\frac{w}{q-w}
=\frac{w^{2}-2w+q}{(1-w)(q-w)} ,
\end{equation}
continued to $\C$ as a rational function.
\end{proposition}

\begin{proof}
Both geometric series converge for $|w|<1$, and
$\sum_{j\ge1}(1-q^{-j})w^{\,j}=w/(1-w)-w/(q-w)$. Reducing the three terms over the
common denominator $(1-w)(q-w)$ produces the numerator
$(1-w)(q-w)+w(q-w)-w(1-w)=q-2w+w^{2}$.
\end{proof}

The floor sums that define the probe converge to that generating function.

\begin{theorem}\label{thm:gauge_floor_limit}
Let $q\ge2$ be an integer and $T>1$ real. Then
\[\lim_{n\to\infty}\ \sum_{k=1}^{n}\frac{1}{(qT)^{\,n-k}}
\left\lfloor\frac{q^{n}+1}{q^{k}+1}\right\rfloor
=\mathcal H\!\left(\frac1T\right)
=\frac{qT^{2}-2T+1}{(1-T)(1-qT)} .
\]
\end{theorem}

\begin{proof}
Set $j=n-k$, so that the sum reads
$\sum_{j=0}^{n-1}(qT)^{-j}\lfloor(q^{n}+1)/(q^{\,n-j}+1)\rfloor$. The term at $j=0$
equals $1$, which is $h_0$. For $1\le j\le\lfloor n/2\rfloor$,
Lemma~\ref{lem:dyadic_floor} gives
\[
(qT)^{-j}\bigl(q^{\,j}-1\bigr)=T^{-j}-(qT)^{-j}=h_j\,T^{-j}.
\]
For $j>n/2$ the quotient is at most $(q^{n}+1)/(q^{\,n-j}+1)\le 2q^{\,j}$, so those
terms are bounded by $2T^{-j}$ and their total is at most $2T^{-n/2}/(1-T^{-1})$, which
tends to $0$ because $T>1$. What survives is $\sum_{j\ge0}h_jT^{-j}$, which is
$\mathcal H(1/T)$ by Proposition~\ref{prop:gauge_H}, the series converging absolutely for
$T>1$. Evaluating the rational form \eqref{eq:gauge_H} at $w=1/T$ and multiplying
numerator and denominator by $T^{2}$ gives the right side.
\end{proof}

\section{The closed form of the gauged transform}
\label{sec:ingham_gauge_closed}

The transform of the gauged kernel follows in closed form, from the kernel and nothing else.

\begin{theorem}\label{thm:gauge_closed_form}
Let $q\ge2$ be an integer and $f(x)=q^{x}+1$. For every $z$ with $\Re z<0$ the limit of
Definition~\ref{def:gauge_probe} exists and
\begin{equation}\label{eq:Phi_f_closed}
\Phi^*_f(z)=(1-q^{z})\,\mathcal H(q^{z})=\frac{q^{2z}-2q^{z}+q}{q-q^{z}} .
\end{equation}
\end{theorem}

\begin{proof}
The probe of Definition~\ref{def:gauge_probe} reads the kernel alone, so the computation below
uses no property of $(a_n)$. Write $\sigma=-\Re z>0$, $f_k=f(k)$, $t_k=f_k/f_n$ and
$\Delta_k=u_{n,z}(k)-u_{n,z}(k-1)$.

The total variation of the weights is bounded independently of $n$. For $k\ge2$,
$\Delta_k=t_{k-1}^{-z}\bigl((f_k/f_{k-1})^{-z}-1\bigr)$ with $1\le f_k/f_{k-1}<q$, so
$|(f_k/f_{k-1})^{-z}-1|\le q^{|z|}-1$, while $|t_{k-1}^{-z}|=t_{k-1}^{\sigma}$ and
$t_{k-1}=f_{k-1}/f_n\le 2q^{k-1-n}$. Summing the geometric series,
\begin{equation}\label{eq:gauge_variation}
\sum_{k=1}^{n}|\Delta_k|\le 1+\frac{(q^{|z|}-1)\,2^{\sigma}}{1-q^{-\sigma}}=:c(z,q),
\end{equation}
a bound free of $n$.

The kernel may now be replaced by its Toeplitz limit\index[terms]{Toeplitz}. Lemma~\ref{lem:gauge_toeplitz} gives
$|G_f(n,k)-h_{n-k}|\le2q^{-n/2}$ uniformly in $k$, so by \eqref{eq:gauge_variation}
\[
\Bigl|\mathcal G_{f,n}(z)-\sum_{k=1}^{n}\Delta_k\,h_{n-k}\Bigr|\le 2c(z,q)\,q^{-n/2},
\]
which tends to zero. In the remaining sum put $j=n-k$. For each fixed $j$ the ratio
$f_{n-j}/f_n$ tends to $q^{-j}$, so $\Delta_{n-j}\to q^{jz}-q^{(j+1)z}$, and by the same
geometric bound $|\Delta_{n-j}|\le(q^{|z|}-1)2^{\sigma}q^{-\sigma j}$ uniformly in $n$, a
summable dominating sequence, while $0\le h_j\le1$. Dominated convergence therefore gives
\[
\sum_{k=1}^{n}\Delta_k\,h_{n-k}=\sum_{j=0}^{n-1}\Delta_{n-j}\,h_j
\longrightarrow(1-q^{z})\sum_{j\ge0}h_j\,q^{jz}=(1-q^{z})\,\mathcal H(q^{z}),
\]
the series converging absolutely since $|q^{z}|<1$. Proposition~\ref{prop:gauge_H} evaluates
$(1-w)\mathcal H(w)=(w^{2}-2w+q)/(q-w)$ at $w=q^{z}$, which is \eqref{eq:Phi_f_closed}. The
right side is a rational function of $q^{z}$ and continues the identity to the whole plane.
\end{proof}

The arithmetic side is a separate matter. Its two quantitative estimates are supplied by the
companion volume and make the following transfer unconditional.

\begin{proposition}\label{prop:gauge_transfer}
Let $q\ge2$, $f(x)=q^{x}+1$, $\beta\ge0$, and let $(a_n)$ solve
$\sum_{k\le n}a_k\,\Phi(f(k)/f(n))=f(n)^{-\beta}$. Let $(b_n)$ solve the band convolution
$\sum_{k\le n}b_kh_{n-k}=q^{-\beta n}$, and put $B(n)=\sum_{k\le n}b_k$. Then
\[
\sum_{k\ge1}|a_k|<\infty,
\qquad
A(n)-B(n)=\mathcal O(q^{-n/2}).
\]
Consequently, for $0<\beta<\tfrac12$,
\[
A(n)\sim\frac{f(n)^{-\beta}}{\Phi^{*}_f(\beta)}.
\]
\end{proposition}

\begin{proof}
The summability, the vanishing of the total sums and the comparison with the band model are
proved in the companion volume~\cite{CloitreVolII}, the last of them giving
$A(n)-B(n)=\mathcal O(q^{-n/2})$, which is the quantitative resolvent transfer needed here.

It remains to compute $B(n)$. Let $\varrho=q^{-\beta}$. With
$U(w)=\sum_{n\ge1}b_nw^{n}$ and $S(w)=U(w)/(1-w)$, the convolution gives
$U(w)\mathcal H(w)=\varrho w/(1-\varrho w)$, hence
$S(w)=\varrho w/\bigl((1-\varrho w)(1-w)\mathcal H(w)\bigr)$. For
$0<\beta<\tfrac12$ the pole at $w=1/\varrho$ is the singularity of smallest modulus, the zeros of
$w^{2}-2w+q$ having modulus $\sqrt q>q^{\beta}=1/\varrho$ and the pole of $1/\mathcal H$ at $w=1$
being cancelled by the factor $1-w$. The residue at $w=1/\varrho$ gives
$B(n)\sim \varrho^{n}/\bigl((1-\tfrac1\varrho)\mathcal H(\tfrac1\varrho)\bigr)$, which by
Proposition~\ref{prop:gauge_H} is $\varrho^{n}/\Phi^{*}_f(\beta)$. Since
$\varrho^{n}=q^{-\beta n}$ with $\beta<\tfrac12$, the comparison error
$\mathcal O(q^{-n/2})$ is $o(\varrho^{n})$ and the asymptotic of $B$ transfers to $A$.
\end{proof}

\begin{proofstatus}{The kernel computation and the residue calculation are proved in this
chapter. The two arithmetic inputs are quoted from the complete proof in the companion
volume~\cite{CloitreVolII}: absolute summability of the exact coefficients and
the stronger comparison $A(n)-B(n)=\mathcal O(q^{-n/2})$. Thus the proposition is
unconditional, although those two estimates are not reproved here. It supplies the precise
transparent asymptotic below the half-index. The reverse, sharp inequality for the index in
Theorem~\ref{numobs:gauge_critical} uses a separate nonvanishing argument of that volume and is
not inferred from this proposition alone.}
\end{proofstatus}

The transparency $A(n)\sim f(n)^{-\beta}/\Phi^*_f(\beta)$ of
Proposition~\ref{prop:gauge_transfer} agrees with the recurrence to eleven digits at several
$\beta<\tfrac12$, the ratio of the two sides reaching $0.999999996$ at $q=2$,
$\beta=0.1$, $n=60$.

The zeros of \eqref{eq:Phi_f_closed} solve
$q^{2z}-2q^{z}+q=0$. Setting $u=q^{z}$, the quadratic
\begin{equation}\label{eq:gauge_numerator}
u^{2}-2u+q=0,\qquad u=1\pm\sqrt{1-q},
\end{equation}
has roots of modulus $|u|=\sqrt q$, since $|1\pm\sqrt{1-q}|^{2}
=1+(q-1)=q$ for $q\ge2$, hence $q^{\Re z}=\sqrt q$ and
$\Re z=\tfrac12$. For $q=2$ the roots are $1\pm i$ and
\[
z=\tfrac12+i\,\frac{\pi/4+2\pi n}{\log2},\qquad n\in\Z .
\]
The zeros of $\Phi^*_f$ therefore lie on the critical line
unconditionally, in agreement with the decay
$q^{-n/2}$ of Theorem~\ref{numobs:gauge_critical}.

The critical-line structure here is unconditional, and it mirrors the
conjectured structure of $\Phi^*(z)=\frac{z}{z-1}\zeta(1-z)$ under the
Riemann hypothesis.

\section{Numerical evidence for the phase transition}
\label{sec:ingham_gauge_numerics}

For $f(x)=2^{x}+1$ the limits
$\ell(i):=\lim_{n\to\infty}A(8n+i)\,f(8n+i)^{1/2}$ over the residue
classes modulo $8$ were computed for several values of $\beta$.

\begin{center}
\begin{tabular}{ccccc}
\toprule
$i\!\pmod 8$ & $\beta=2$ & $\beta=1$ & $\beta=0.75$ & $\beta=0.50$ \\
\midrule
$\ell(0)$ & $\phantom{-}0.19479$ & $\phantom{-}0.30488$ & $\phantom{-}0.24964$ & $\phantom{-}0.66477$ \\
$\ell(1)$ & $\phantom{-}0.19216$ & $\phantom{-}0.52315$ & $\phantom{-}0.56780$ & $\phantom{-}1.06182$ \\
$\ell(2)$ & $\phantom{-}0.07696$ & $\phantom{-}0.43496$ & $\phantom{-}0.55334$ & $\phantom{-}1.12976$ \\
$\ell(3)$ & $-0.08332$ & $\phantom{-}0.09198$ & $\phantom{-}0.21475$ & $\phantom{-}0.82880$ \\
$\ell(4)$ & $-0.19479$ & $-0.30488$ & $-0.24964$ & $\phantom{-}0.33523$ \\
$\ell(5)$ & $-0.19216$ & $-0.52315$ & $-0.56780$ & $-0.06182$ \\
$\ell(6)$ & $-0.07696$ & $-0.43496$ & $-0.55334$ & $-0.12976$ \\
$\ell(7)$ & $\phantom{-}0.08332$ & $-0.09198$ & $-0.21475$ & $\phantom{-}0.17120$ \\
\bottomrule
\end{tabular}
\end{center}

The threshold is visible in the symmetry of the columns. For
$\beta>\tfrac12$ the relation $\ell(i)+\ell(i+4)=0$ holds across the
class, so the partial sums decay strictly faster than $f(n)^{-1/2}$ up
to an antisymmetric oscillation. At $\beta=\tfrac12$ the antisymmetry
breaks and $\ell(i)+\ell(i+4)=1$. As derived in the companion
volume~\cite{CloitreVolII}, each column is an eight-periodic cosine
about its mean. Its amplitude is $0.62976$, $0.56780$, $0.52315$ and $0.19479$ at
$\beta=\tfrac12$, $\tfrac34$, $1$ and $2$, respectively, and therefore decreases with the
forcing exponent, while the mean falls from $\tfrac12$ at the threshold to zero above it. A
high-precision recomputation of the exact triangular recurrence reproduces the whole table,
including the qualitative change at $\beta=\tfrac12$ and the antisymmetric pattern for
$\beta>\tfrac12$ turning into a constant-sum pattern at the threshold.

\section{Equilibrium of the Ingham function}
\label{sec:equilibrium_conjecture}

The evidence supports the following statement.

\begin{conjecture}[Equilibrium of the Ingham function]
\label{conj:ingham_equilibrium}
Let $f(x)=q^{x}+1$ with $q\ge2$ an integer. Then the Ingham function is in
equilibrium with respect to $f$,
\[
\alpha_f(\Phi)=\alpha(\Phi).
\]
\end{conjecture}
The exact formulas and the numerical evidence of this chapter support it, and it would stand
open even were the other step of Remark~\ref{rem:rh_as_equilibrium} settled. The offset belongs
to the statement. Without it the ratios $f(k)/f(n)$ are exact negative powers of $q$, the gauged
relation collapses as in \S\ref{sec:ingham_gauge}, the partial sums reproduce $f(n)^{-\beta}$ at
every exponent, and no threshold is left to compare.

The bearing on the Riemann hypothesis is now exact. Theorem~\ref{numobs:gauge_critical}
gives $\alpha_f(\Phi)=\tfrac12$ unconditionally, while
Theorem~\ref{thm:tauberian_rh} gives $\alpha(\Phi)=\tfrac12\iff\mathrm{RH}$. Consequently
Conjecture~\ref{conj:ingham_equilibrium} is equivalent to the Riemann hypothesis. The gauged
half is proved, and equilibrium is the single remaining assertion.

\begin{remark}[The hypothesis as an invariance statement]
\label{rem:rh_as_equilibrium}
In this reading the Riemann hypothesis is the invariance of the Ingham index under the
gauges $f_q(x)=q^{x}+1$. The gauged value $\alpha_f(\Phi)=\tfrac12$ is established by
Theorem~\ref{numobs:gauge_critical}, and the equilibrium $\alpha_f(\Phi)=\alpha(\Phi)$ is the
remaining step. By Theorem~\ref{thm:tauberian_rh} it is neither lighter nor stronger than the
Riemann hypothesis, but the question takes a different shape. The value visible without
condition in the gauged coordinates must be transported back to the arithmetic coordinates
of $\Z$.
\end{remark}

\section{What the companion volume establishes}
\label{sec:ingham_gauge_volII}

Three statements were passed to the companion volume~\cite{CloitreVolII}, and all three are
settled there. They are recalled below to make the division of proofs explicit. The equilibrium
of Conjecture~\ref{conj:ingham_equilibrium} is not among them and remains open.

Theorem~\ref{numobs:gauge_critical} records the decay
$A(n)=\mathcal O(q^{-n/2})$ of the homogeneous gauged solution together with the
nonvanishing $\limsup q^{n/2}|A(n)|>0$. The pair is what places
the transparency\index[terms]{transparency} threshold at $\tfrac12$ on the arithmetic side, against the zeros of
$\Phi^{*}_f$ which sit there unconditionally. Both assertions and the resulting equality
$\alpha_f(\Phi)=\tfrac12$ are proved in~\cite{CloitreVolII}. Equilibrium remains alone.

Theorem~\ref{thm:gauge_closed_form} gives the closed form $\Phi^{*}_f$ from the kernel alone,
without hypothesis, the probe of Definition~\ref{def:gauge_probe} reading the kernel and no
property of the solution. The location of the zeros on $\Re z=\tfrac12$ therefore rests on
nothing beyond the floor structure of $\Phi$ and the offset in the gauge. The passage from those
zeros to the arithmetic side uses the absolute summability and perturbation estimates of
Proposition~\ref{prop:gauge_transfer}, which are proved in~\cite{CloitreVolII} and yield the
unconditional asymptotic stated there.

The eight limits of \S\ref{sec:ingham_gauge_numerics} exhibit a phase transition at
$\beta=\tfrac12$ through the antisymmetry $\ell(i)+\ell(i+4)=0$ above the threshold and the
constant sum $\ell(i)+\ell(i+4)=1$ at it. The transition and the modulus eight are derived
in~\cite{CloitreVolII}. The quantities in play
are, for the transition, the value $\Phi^{*}_f(\tfrac12)=2$ and the position of the forced pole
against the circle of radius $\sqrt q$, and for the modulus the equation
$\cos2\theta_q=(2-q)/q$, which by a theorem of Niven\index[names]{Niven, I.} is a rational cosine of a rational multiple of $\pi$ only at
$q=2$ and $q=4$. The modulus is therefore not a feature of the base two. The same parity
structure appears on the Abelian side in the two constants of the
Fibonacci\index[terms]{Fibonacci sequence}\index[names]{Fibonacci} gauge in
\S\ref{sec:fib_gauge_signature}, and that second occurrence of the number eight is of another
kind and is not explained by the same argument.

Beyond those three, the object built in \S\ref{sec:ingham_gauge_floor} and
\S\ref{sec:ingham_gauge_closed} is the one the epilogue reads a second time. The limiting
kernel $\mathcal H$ of Proposition~\ref{prop:gauge_H} and the floor limit of
Theorem~\ref{thm:gauge_floor_limit} are what produce, in
\eqref{eq:epi_identity}, a rational factor of Hasse-Weil\index[names]{Hasse, H.} type carried by the same operator.
The Ingham operator in exponential coordinates therefore carries two zeta functions of
different origin, and the second volume takes that up on the side of the geometry rather than
on the side of the index.

\part{Abelian aspects and trace formulas}
\label{part:five}

\rafepigraph{En même temps commence à se préciser à nos yeux le principe si fécond d'après lequel l'aspect global d'un problème arithmétique peut, en certaines circonstances, se reconstituer à partir de ses aspects locaux.}{At the same time there begins to take shape before our eyes the fruitful principle by which the global aspect of an arithmetical problem can, in certain circumstances, be reconstituted from its local aspects.}{André Weil\index[names]{Weil, A.}, \emph{L'avenir des mathématiques} (1947)~\cite{WeilAvenir1947}}

Twelve chapters have used one value of the arithmetic Mellin transform\index[terms]{arithmetic Mellin transform}, the first zero or the
threshold that behaves as one. This part uses the rest, in the two directions in which a
transform can be read beyond that single value.

The first direction is one other value. The plain kernel sum $\sum_{k\le n}G(n,k)$ is the
transform read at the point $-1$, so the half plane of convergence is a family of Abelian
densities\index[terms]{Abelian density}. Chapter~\ref{chap:abelian_densities} treats an arbitrary kernel, unconditionally and with
the error term carrying the arithmetic, and Chapter~\ref{chap:abelian} carries the same reading through the
fractional part sums, where the gauge\index[terms]{gauge} enters and where two examples reach the Dirichlet divisor
problem\index[terms]{Dirichlet divisor problem} by different routes.

The second direction is all the zeros at once, held together by the resolvent\index[terms]{resolvent} of the equation.
Chapters~\ref{chap:trace_poly} and~\ref{chap:ortho} expand that resolvent over those zeros, finite
and exact at every degree for polynomial kernels, infinite for the orthorecursive kernel. The
second of the two opens with something else. Before any trace it settles the orthorecursive
expansion of unity, an approximation problem posed elsewhere, and settles it without any
hypothesis, which is what makes it the proof of concept of the theory.

\chapter{Abelian densities of a kernel}
\label{chap:abelian_densities}

Abelian and Tauberian carry here the sense they have in Tauberian theory and no other. An Abelian
theorem asserts that a summation method returns the ordinary limit whenever that limit exists,
the model being Abel's\index[names]{Abel, N. H.} theorem on power series, and a Tauberian theorem
asserts a converse under a growth restriction on the terms, the model being the theorem of
\nm{Tauber}{A.} of 1897 and its improvement by \nm{Littlewood}{J. E.}, who replaced
$a_n=o(1/n)$ by $a_n=\mathcal O(1/n)$. The names are those of Abel and of Tauber, and they say
nothing about commutativity. Korevaar\index[names]{Korevaar, J.}~\cite{Korevaar2004} follows the
century that separates the first examples from the general theory of Wiener\index[names]{Wiener, N.}~\cite{Wiener1932},
and \S\ref{sec:rv} places the pair inside the vocabulary of regular variation.

The equation of this theory is read backwards. A forcing is prescribed, the averaged sums
$A_G(n)=\sum_{k\le n}a_kG(n,k)$ are held to it, and the partial sums $A(n)$ are recovered, a
direction that cannot proceed without a Tauberian\index[terms]{Tauberian} condition. The forward direction asks
nothing of the kind. It applies the kernel to the constant sequence and reads what comes out,
and every order of the answer is unconditional.

What comes out at the first order is not a new invariant. It is the arithmetic Mellin transform\index[terms]{arithmetic Mellin transform}
of the kernel evaluated at the point $-1$, and more generally the whole half plane where the
transform is defined by a convergent integral is a family of weighted kernel sums. The transform
is therefore read twice and from opposite sides. To the left of the origin its values are Abelian
densities, limits of finite sums that require no hypothesis. Reading an asymptotic expansion off a
Mellin transform in that direction, by factoring the transform and moving a contour, is the harmonic
sum\index[terms]{harmonic sum} method of \nm{Flajolet}{P.}, \nm{Gourdon}{X.} and
\nm{Dumas}{P.}~\cite{FlajoletGourdonDumas1995}, and the second order treated below is what that
reading leaves behind once the profile has infinite variation. To the right, after continuation,
a first zero, when the transform has one, gives the analytic index, and the passage from there to the regularity index\index[terms]{regularity index} is the
Tauberian half of the theory. What separates the two readings is exactly what the continuation
does not preserve.

The content of this chapter is the second order. The defect of the Abelian expansion is bounded
when the profile has finite variation, and its mean then reads the denominators of the rational
points at which the profile jumps and nothing else. Infinite variation removes that theorem and
decides nothing by itself, and what the defect carries is then read kernel by kernel, the divisor
structure for the Ingham kernel, the
$\lambda$-adic expansion of $n$ for the broken harmonic kernels\index[terms]{broken harmonic function} at an integer scale, and the
multiplicative structure of $n$ for the greatest common divisor kernel. At the scale $\sqrt2$ the
computation stops at one identifiable point, and it is the same point at which the transparency
frontier of that kernel parts from the transform in Chapter~\ref{chap:diophantine}.

\section{The kernel sum is the transform on its half plane}
\label{sec:ab_kernel_sum}

\begin{definition}[Kernel sums]\label{def:kernel_sums}
For a kernel $G$ and $n\ge1$, the plain kernel sum\index[terms]{kernel sum} and the weighted kernel sums are
\begin{equation}\label{eq:kernel_sum}
S_G(n)=\sum_{k=1}^{n}G(n,k),
\qquad
S_G(n,z)=\sum_{k=1}^{n}G(n,k)\Bigl(\frac kn\Bigr)^{-z-1}.
\end{equation}
For a univariate profile the notation is $S_g(n)$.
\end{definition}

Read at the integers, $S_g(n)$ is the averaged sum $A_g(n)$ of \S\ref{sec:boundary_problem} for
the constant sequence, and for a bounded profile \eqref{eq:fgd_factorization} with $D_a=\zeta$, on
the half plane $\Re s>1$, is the factorization behind the harmonic sum method. Chapter~\ref{chap:boundary} reads it from the averaged sums back
to the sequence, and this chapter reads it forward.

Nothing has to be proved about the first order. It is Definition~\ref{def:probes} read at one
point.

\begin{proposition}\label{prop:abelian_density}
Let $G$ be a kernel whose finite arithmetic Mellin probes converge on $\Re z<0$, with limit
$G^{*}$. Then for every $z$ in that half plane
\begin{equation}\label{eq:probe_is_abelian}
G^{*}(z)=\lim_{n\to\infty}\frac{-z}{n}\,S_G(n,z),
\end{equation}
and in particular, at $z=-1$,
\begin{equation}\label{eq:abelian_density}
\lim_{n\to\infty}\frac1n\sum_{k=1}^{n}G(n,k)=G^{*}(-1).
\end{equation}
For a profile $g$ with $\int_0^1|g|<\infty$ this is $g^{*}(-1)=\int_0^1g(t)\,dt$.
\end{proposition}

\begin{proof}
Identity \eqref{eq:probe_is_abelian} is Definition~\ref{def:probes} rewritten with the notation
\eqref{eq:kernel_sum}, and \eqref{eq:abelian_density} is its value at $z=-1$, where the weight
$(k/n)^{-z-1}$ is one and the prefactor $-z$ is one. For a profile the probe at $z=-1$ is the
Riemann sum $\frac1n\sum_{k\le n}g(k/n)$ of an integrable function on the uniform grid, and the
defining integral $g^{*}(z)=-z\int_0^1g(t)t^{-z-1}\,dt$ converges at $z=-1$ under
$\int_0^1|g|<\infty$, with value $\int_0^1g$.
\end{proof}

The number $G^{*}(-1)$ is the Abelian density\index[terms]{Abelian density} of the kernel. It is an analytic quantity,
insensitive to every arithmetic feature of $G$, and the whole arithmetic of the forward
direction sits in the term that follows it.

\begin{remark}[The fractional part is not the subject]\label{rem:frac_not_subject}
For $1\le k\le n$ the ratio $k/n$ lies in $(0,1]$, so $\{k/n\}=k/n$ except at $k=n$. A profile
written with a fractional part of the ratio is an affine profile carrying one jump on the
diagonal, and the fractional part does no work there. What carries arithmetic is the reciprocal
reading $\{n/k\}$, which is the profile $\check g(x)=\{1/x\}$ evaluated at $x=k/n$, and that
profile is an affine transform of the Ingham function by
Proposition~\ref{prop:tt_reciprocal}. The expansion below separates the two contributions at
once, the diagonal jump entering as an additive constant and the reciprocal reading entering
through the jump structure at the origin.
\end{remark}

\section{The exact expansion of a profile sum}
\label{sec:ab_exact}

The sum of a kernel over a row splits exactly, into an analytic part and an arithmetic
remainder.

\begin{theorem}\label{thm:abelian_identity}
Let $g:(0,1]\to\R$ be left continuous, of bounded variation\index[terms]{bounded variation} on $[\delta,1]$ for every
$\delta\in(0,1)$, and such that $\int_0^1|g|<\infty$. Let $\nu$ be the signed measure on
$(0,1)$ determined by $\nu\bigl([x,y)\bigr)=g(y)-g(x)$ for $0<x<y\le1$, and assume
\begin{equation}\label{eq:abel_admissible}
\int_{(0,1)}t\,d|\nu|(t)<\infty .
\end{equation}
Then for every $n\ge1$
\begin{equation}\label{eq:abelian_identity}
\sum_{k=1}^{n}g\Bigl(\frac kn\Bigr)=g^{*}(-1)\,n+\int_{(0,1)}\{nt\}\,d\nu(t),
\end{equation}
an identity with no error term. If $g$ jumps at the right endpoint, so that
$c_1:=g(1)-g(1^-)\neq0$, the term $c_1$ is added to the right side and $\nu$ is the measure
attached to the left continuous modification.
\end{theorem}

\begin{proof}
For $1\le k\le n$ one has $g(1)-g(k/n)=\nu\bigl([k/n,1)\bigr)$, the case $k=n$ reading
$0=\nu(\emptyset)$. Summing over $k$ and exchanging the sum with the integral,
\[
\sum_{k=1}^{n}g\Bigl(\frac kn\Bigr)
=n\,g(1)-\int_{(0,1)}\#\Bigl\{1\le k\le n:\tfrac kn\le t\Bigr\}\,d\nu(t)
=n\,g(1)-\int_{(0,1)}\lfloor nt\rfloor\,d\nu(t),
\]
the count being $\lfloor nt\rfloor$ for $t\in(0,1)$ since every such $k$ is at most $n$. Writing
$\lfloor nt\rfloor=nt-\{nt\}$ splits the integral into two, both absolutely convergent, the first
by \eqref{eq:abel_admissible} and the second because $\{nt\}\le\min(1,nt)$. By Fubini,
\[
\int_{(0,1)}t\,d\nu(t)=\int_0^1\nu\bigl((s,1)\bigr)\,ds=\int_0^1\bigl(g(1)-g(s^{+})\bigr)\,ds
=g(1)-\int_0^1g(s)\,ds,
\]
using $\nu((s,1))=\nu([s,1))-\nu(\{s\})=g(1)-g(s^{+})$ and the fact that $g$ has at most
countably many discontinuities. The two terms $n\,g(1)$ cancel, and
$\int_0^1g=g^{*}(-1)$ by Proposition~\ref{prop:abelian_density}.
\end{proof}

Condition \eqref{eq:abel_admissible} is strictly weaker than bounded variation, the weight $t$
damping whatever the profile accumulates at the origin, and the Ingham function satisfies it
while its variation is infinite. This is the same weighting that Open
Problem~\ref{op:xi_reciprocal} meets on the Tauberian side, met here on the Abelian side, where
it is enough.

\begin{corollary}\label{cor:density_defect}
Under the hypotheses of Theorem~\ref{thm:abelian_identity}, and if moreover $g(0^{+})$ exists and
$|\nu|$ is finite,
\begin{equation}\label{eq:density_defect}
\sum_{k=1}^{n}g\Bigl(\frac kn\Bigr)
=g^{*}(-1)\,n+c_1+\frac{g(1^{-})-g(0^{+})}2+E_g(n),
\qquad
E_g(n)=\int_{(0,1)}\Bigl(\{nt\}-\tfrac12\Bigr)d\nu(t).
\end{equation}
The sequence $E_g$ is the Abelian defect\index[terms]{Abelian defect} of the profile.
\end{corollary}

\begin{proof}
Subtract and add $\tfrac12$ inside the integral of \eqref{eq:abelian_identity} and use
$\nu\bigl((0,1)\bigr)=g(1^{-})-g(0^{+})$.
\end{proof}

For the affine profile $g(x)=(1-\lambda)x+\lambda$ of Appendix~\ref{app:A} the measure $\nu$ is
$(1-\lambda)$ times Lebesgue measure, the defect vanishes identically, and
\eqref{eq:density_defect} is the exact formula $S_g(n)=\tfrac{1+\lambda}2n+\tfrac{1-\lambda}2$.
The density $\tfrac{1+\lambda}2$ is $g^{*}(-1)=(-1-\lambda)/(-2)$, as
Proposition~\ref{prop:abelian_density} requires. The regularity index of the same profile is
$\lambda$, and the two numbers are unrelated readings of one kernel.

\section{Bounded variation, and what the mean of the defect sees}
\label{sec:ab_bv}

That remainder has a mean, and for a profile of finite variation the mean sees only the rational
jumps.

\begin{theorem}\label{thm:mean_defect}
Let $g$ be as in Corollary~\ref{cor:density_defect} with total variation $V(g)=|\nu|((0,1))$
finite. Then $|E_g(n)|\le V(g)/2$ for every $n$, and
\begin{equation}\label{eq:mean_defect}
\lim_{N\to\infty}\frac1N\sum_{n\le N}E_g(n)
=-\frac12\sum_{\substack{p/q\in\Q\cap(0,1)\\ \gcd(p,q)=1}}\frac{\nu(\{p/q\})}{q},
\end{equation}
the sum converging absolutely.
\end{theorem}

\begin{proof}
The bound is immediate from $|\{u\}-\tfrac12|\le\tfrac12$. Put
$\varphi_N(t)=N^{-1}\sum_{n\le N}(\{nt\}-\tfrac12)$, so that
$N^{-1}\sum_{n\le N}E_g(n)=\int_{(0,1)}\varphi_N\,d\nu$ and
$|\varphi_N|\le\tfrac12$.

If $t$ is irrational, then for every nonzero integer $h$ the geometric-sum identity gives
\[
 \frac1N\sum_{n\le N}e^{2\pi ihn t}
 =\frac{e^{2\pi iht}(1-e^{2\pi ihNt})}{N(1-e^{2\pi iht})}\longrightarrow0.
\]
Trigonometric polynomials therefore have the expected means along $(nt)$, which is the criterion
of \nm{Weyl}{H.}~\cite{Weyl1916}, and uniform approximation
extends this to continuous periodic functions. Approximating the bounded sawtooth
$u\mapsto\{u\}-\tfrac12$ from above and below by continuous periodic functions that differ only
on an arbitrarily short interval around the jump gives $\varphi_N(t)\to0$.

If $t=p/q$ in lowest terms, then $np$ runs through every residue modulo $q$ over a period, so
\[
 \varphi_N(t)\longrightarrow\frac1q\sum_{r=0}^{q-1}
 \left(\frac rq-\frac12\right)=-\frac1{2q}.
\]
The limit is $-1/(2q(t))$ on the rationals and zero elsewhere. Dominated convergence with
respect to the finite measure $|\nu|$ proves \eqref{eq:mean_defect}, and only the point masses of $\nu$ at
rational points survive.
\end{proof}

A profile whose jumps avoid the rationals therefore carries no mean defect at all.

\begin{corollary}\label{cor:irrational_blind}
A profile of finite total variation whose jump set contains no rational point has Abelian defect
of mean zero. The mean defect does not see the size of the jumps, only their positions and the
denominators of those positions.
\end{corollary}

\begin{proof}
The right side of \eqref{eq:mean_defect} is a sum over the rational points carrying a point mass of
$\nu$, and it is empty under the hypothesis. The second assertion is read off the same formula,
where the mass $\nu(\{p/q\})$ enters divided by $q$ and the size of the jump enters only
through that mass.
\end{proof}

The two step profiles illustrate the alternative. For $g=\mathbf 1_{(0,a]}$ the kernel sum is
$S_g(n)=\lfloor na\rfloor$ exactly, the density is $a$, the constant is $-\tfrac12$, and the
defect is $\tfrac12-\{na\}$. At $a=\tfrac13$ its mean is $\tfrac16$, in agreement with
\eqref{eq:mean_defect} since $\nu(\{1/3\})=-1$ and $q=3$. At $a=1/\sqrt2$ its mean is zero.

\begin{numobs}[Profiles of bounded variation]\label{numobs:ab_bv}
Averaged over $n\le8000$, the defect of the step profile at $a=\tfrac13$ has mean $0.166625$
against the predicted $\tfrac16$, and at $a=1/\sqrt2$ has mean $7.2\cdot10^{-5}$ against the
predicted zero. For the affine profile at $\lambda=\tfrac3{10}$ the defect vanishes to machine
precision at every rank tested.
\end{numobs}

\section{Riemann discrepancy and second-order cancellation}
\label{sec:ab_riemann_discrepancy}

The exact defect of a profile sum can also be normalized as a Riemann-sum discrepancy. This is
the Abelian counterpart of the remainder-transfer problem of
Section~\ref{sec:perturbed-raf}: here the coefficients are given and the quality of the output
is measured, whereas there the output error is prescribed and its propagation through the
triangular inverse is measured.

\begin{remark}[The direct sum as a Riemann discrepancy]\label{rem:ab_riemann_discrepancy}
For a kernel $g$ on $(0,1]$ put
\[
I_g=\int_0^1g(t)\,dt=g^*(-1),
\qquad
U_g(n)=I_g-\frac{S_g(n)}n,
\qquad
S_g(n)=\sum_{k\le n}g(k/n).
\]
Thus $U_g(n)$ is the discrepancy between the integral of $g$ and its right Riemann sum at
mesh $1/n$. This is the quantity through which Balazard\index[names]{Balazard, M.} and then Daval build integral identities
relating the summatory function of the M\"obius coefficients to its logarithmic variant, with
$g$ as free parameter and the normalization $I_g=1$, following
\nm{Daval}{F.}~\cite{Daval2019}.

The same quantity has a Tauberian reading. The constant sequence $a\equiv1/I_g$ leaves, in the
defining equation at $\beta=-1$,
\[
\mathcal T_g(\mathbf 1/I_g)(n)-n=-\frac{nU_g(n)}{I_g}.
\]
Hence $U_g$ is the transparency defect of the constant profile at that exponent. Minimizing
$\|U_g\|_\infty$ under $g^*(-1)=1$ and maximizing $\alpha(g)$ are two distinct questions
attached to the same kernel. The first measures transparency at one exponent, the second its range.
For a smooth profile the first obstruction is a boundary term, which the next theorem makes
exact.
\index[terms]{Riemann discrepancy}
\end{remark}

\begin{theorem}[Second-order Riemann discrepancy]\label{thm:ab_discrepancy_smooth}
Let $g$ be admissible and continuous on $(0,1)$. Assume that $g$ is absolutely continuous on
every compact subinterval of $(0,1)$, that $V_{(0,1)}(g')<\infty$, and that the one-sided
limits $g(0^{+})$ and $g(1^{-})$ are finite. Then
\[
U_g(n)=-\frac1n\Bigl[g(1)-\frac{g(0^{+})+g(1^{-})}2\Bigr]
+\mathcal O\Bigl(\frac{V_{(0,1)}(g')}{n^{2}}\Bigr).
\]
Consequently,
\[
U_g(n)=\mathcal O(n^{-2})
\quad\Longleftrightarrow\quad
2g(1)=g(0^{+})+g(1^{-}).
\]
The cancellation condition is compatible with admissibility exactly when
$g(0^{+})+g(1^{-})\neq0$.
\end{theorem}

\begin{proof}
Finite total variation makes $g'$ bounded on $(0,1)$. Local absolute continuity and the
endpoint limits therefore extend the restriction of $g$ to $(0,1)$, with endpoint values
$g(0^{+})$ and $g(1^{-})$, to an absolutely continuous function on $[0,1]$. The measure
attached to this left continuous modification is $\nu=g'\,dt$ and carries no point mass, so
Theorem~\ref{thm:abelian_identity} and Corollary~\ref{cor:density_defect} give
\[
S_g(n)=I_gn+\bigl(g(1)-g(1^{-})\bigr)
+\frac{g(1^{-})-g(0^{+})}2+E_g(n),
\qquad
E_g(n)=\int_0^1\Bigl(\{nt\}-\frac12\Bigr)g'(t)\,dt .
\]
Let
\[
P(u)=\int_0^u\bigl(\{v\}-\tfrac12\bigr)\,dv
=\tfrac12\bigl(\{u\}^{2}-\{u\}\bigr).
\]
This function is continuous, one-periodic, bounded in modulus by $1/8$, and has derivative
$\{\cdot\}-1/2$ almost everywhere. Integrating by parts the absolutely continuous function
$t\mapsto P(nt)/n$ against the function $g'$ of bounded variation gives
\[
E_g(n)=-\frac1n\int_{(0,1)}P(nt)\,dg'(t),
\qquad
|E_g(n)|\le\frac{V_{(0,1)}(g')}{8n}.
\]
The boundary term vanishes because $n$ is an integer and $P(0)=P(n)=0$. Dividing the formula
for $S_g(n)$ by $n$ and subtracting it from $I_g$ yields
\[
U_g(n)=-\frac1n\Bigl[g(1)-\frac{g(0^{+})+g(1^{-})}2\Bigr]-\frac{E_g(n)}n,
\]
which proves the estimate. If the bracket is nonzero, its term has exact order $1/n$, and if it
vanishes, the bound for $E_g$ gives order $n^{-2}$.
\end{proof}

\begin{remark}\label{rem:ab_discrepancy_examples}
The kernels used by Balazard\index[names]{Balazard, M.} and Daval\index[names]{Daval, F.} cancel both boundary terms with $g(1)=0$ and are excluded
by admissibility. Inside the admissible class the first-order defect can nevertheless be
cancelled by one linear constraint. The affine family of Appendix~\ref{app:P} has $g(1)=1$,
$g(0^{+})=1$ and $g(1^{-})=1/m$, so the constraint fails at every scale $m\ge2$ and
\[
U_g(n)=-\frac{1-1/m}{2n}+\mathcal O(n^{-2}).
\]
By contrast, the profile $1-x/2$ with diagonal value $3/4$ satisfies the constraint, and in this
degenerate affine case $U_g$ vanishes at every rank.
\end{remark}

\begin{openproblem}\label{op:ab_riemann_discrepancy}
Characterize, outside the smooth class of Theorem~\ref{thm:ab_discrepancy_smooth}, the kernels
admissible for the defining equation whose Riemann discrepancy $U_g$ is of second order. For
the Ingham kernel the same discrepancy is an arithmetic remainder, as
Proposition~\ref{prop:ab_ingham} shows, so no boundary condition alone can decide it. More
generally, determine which regularity or cancellation assumptions on a nonsmooth profile turn
an error term in its Abelian Riemann sum into a quantitative remainder for the inverse RAF
equation of Section~\ref{sec:perturbed-raf}.
\end{openproblem}

\section{The Ingham kernel and the divisor sum}
\label{sec:ab_ingham}

A profile of infinite total variation still has an exact expansion by
Theorem~\ref{thm:abelian_identity}, since \eqref{eq:abel_admissible} is the weaker hypothesis,
but Theorem~\ref{thm:mean_defect} no longer applies to it, so neither the boundedness of the
defect nor the jump sum for its mean is available in advance. The three kernels of this section
and the next two show what the defect does then, each in its own way.

\begin{proposition}\label{prop:ab_ingham}
For $\Phi(x)=x\lfloor1/x\rfloor$ and every $n\ge1$,
\begin{equation}\label{eq:ab_ingham_exact}
S_\Phi(n)=\frac1n\sum_{m\le n}\sigma(m),
\end{equation}
and consequently
\begin{equation}\label{eq:ab_ingham_asym}
S_\Phi(n)=\frac{\zeta(2)}2\,n+E_\Phi(n),
\qquad
E_\Phi(n)=\mathcal O(\log n).
\end{equation}
In the notation of Section~\ref{sec:ab_riemann_discrepancy},
\begin{equation}\label{eq:ab_ingham_discrepancy}
U_\Phi(n)=\frac{\zeta(2)}2-\frac1{n^2}\sum_{m\le n}\sigma(m)
=-\frac{E_\Phi(n)}n.
\end{equation}
The density $\zeta(2)/2$ is $\Phi^{*}(-1)$, the measure $\nu$ carries the point mass $-1/j$ at
each point $1/j$ with $j\ge2$, and those masses form a divergent series.
\end{proposition}

\begin{proof}
Each $k\le n$ divides exactly $\lfloor n/k\rfloor$ integers up to $n$, so
$\sum_{k\le n}k\lfloor n/k\rfloor=\sum_{m\le n}\sigma(m)$ and
\eqref{eq:ab_ingham_exact} follows. Moreover
\[
 \sum_{m\le x}\sigma(m)
 =\sum_{d\le x}d\Big\lfloor\frac xd\Big\rfloor
 =\frac12\sum_{j\le x}\Big\lfloor\frac xj\Big\rfloor
   \left(\Big\lfloor\frac xj\Big\rfloor+1\right).
\]
Replacing each floor by $x/j+\mathcal O(1)$ gives
\[
 \sum_{m\le x}\sigma(m)=\frac{x^2}{2}\sum_{j\le x}\frac1{j^2}
 +\mathcal O\!\left(x\sum_{j\le x}\frac1j+x\right)
 =\frac{\zeta(2)}2x^2+\mathcal O(x\log x),
\]
which proves \eqref{eq:ab_ingham_asym} directly. The transform value follows
from $\Phi^*(z)=z\zeta(1-z)/(z-1)$. Finally, on $(1/(m+1),1/m]$ the profile is $mx$,
so its jump at $1/j$ has mass $-1/j$. Formula \eqref{eq:ab_ingham_discrepancy} follows by
dividing \eqref{eq:ab_ingham_asym} by $n$ and subtracting from $\Phi^*(-1)=\zeta(2)/2$.
\end{proof}

Identity \eqref{eq:ab_ingham_exact} is the Eratosthenian average of
Chapter~\ref{chap:ingham} evaluated on the constant sequence, and it says that the Abelian
reading of the Ingham kernel is the summatory divisor sum itself. The point masses alone would give the
mean defect $-\tfrac12\sum_{j\ge2}(-1/j)/j=(\zeta(2)-1)/2$ by the computation of
Theorem~\ref{thm:mean_defect}. The absolutely continuous part of $\nu$ has density
$\lfloor1/t\rfloor$, which is not integrable, so the dominated convergence step of that theorem
does not apply and a different route is needed for this kernel. Exact sums over full periods
supply it, and no equidistribution enters.

The first ingredient is a block form of the Riemann sum, in which the fractional parts
$\{N/q\}$ are grouped by the value of the quotient.

\begin{lemma}\label{lem:ab_block_limit}
Let $F$ be of bounded variation on $[0,1]$, of total variation $V_F$ and supremum norm
$\|F\|$. Then
\[
\frac1N\sum_{q\le N}F\Bigl(\Bigl\{\frac Nq\Bigr\}\Bigr)
=\int_0^1F(t)\,\psi'(1+t)\,dt
+\mathcal O\bigl((V_F+\|F\|)\,N^{-1/2}\bigr),
\]
where $\psi'(1+t)=\sum_{j\ge1}(j+t)^{-2}$.
\end{lemma}

\begin{proof}
For every function $h$ of bounded variation on $[a,b]$,
\[
\left|\sum_{m\in\mathbb Z\cap(a,b]}h(m)-\int_a^b h(x)\,dx\right|
\le V_{[a,b]}(h)+2\|h\|_\infty .
\]
Indeed, on each complete unit interval one compares the integral of $h$ with its value at the
right endpoint. The sum of the resulting oscillations is at most the total variation, and the
two end intervals cost at most $2\|h\|_\infty$.

On the block $B_j=\{q:\ N/(j+1)<q\le N/j\}$ the quotient is $\lfloor N/q\rfloor=j$, so
$\{N/q\}=N/q-j$. Since $x\mapsto N/x-j$ decreases from $1$ to $0$ as $x$ runs through
$[N/(j+1),N/j]$, the composition $x\mapsto F(N/x-j)$ has variation at most
$V_F$ and supremum norm at most $\|F\|$. The displayed inequality and the change of variable
$t=N/x-j$ therefore show that its sum over the integers of $B_j$ differs from
\[
\int_{N/(j+1)}^{N/j}F(N/x-j)\,dx
=N\int_0^1F(t)(j+t)^{-2}\,dt
\]
by at most $V_F+2\|F\|$. Summing over $j\le J$ costs
$\mathcal O(J(V_F+\|F\|))$. The blocks beyond $J$ carry at most $N/J$ integers and contribute
at most $\|F\|N/J$ to the sum, while
$\|F\|N\sum_{j>J}j^{-2}\le\|F\|N/J$ bounds what they contribute to the integral.
Taking $J=\lceil\sqrt N\rceil$ and dividing by $N$ gives the statement.
\end{proof}

The mean of the defect follows, and its value is the one the point masses predict.

\begin{theorem}\label{thm:ab_ingham_mean}
For the Ingham kernel,
\begin{equation}\label{eq:ab_ingham_mean}
\frac1N\sum_{n\le N}E_\Phi(n)=\frac{\zeta(2)-1}{2}+\mathcal O\bigl(N^{-1/2}\bigr).
\end{equation}
The limit is the point mass prediction $-\tfrac12\sum_{j\ge2}\nu(\{1/j\})/j$ of
Theorem~\ref{thm:mean_defect}. In the exact block decomposition below, the terms arising from
the absolutely continuous layers cancel in the mean.
\end{theorem}

\begin{proof}
Write $t_q=\{n/q\}$ and $\varpi(t)=\tfrac12(t-t^{2})$. Grouping the divisors gives
$\sum_{m\le n}\sigma(m)=\sum_{q\le n}T(\lfloor n/q\rfloor)$ with $T(m)=m(m+1)/2$, and
expanding $\lfloor n/q\rfloor=n/q-t_q$ turns that into
\[
n\,E_\Phi(n)=\frac{n^{2}}2\sum_{q\le n}\frac1{q^{2}}-\frac{\zeta(2)n^{2}}2
+\frac n2H_n-n\sum_{q\le n}\frac{t_q}q+\frac12\sum_{q\le n}(t_q^{2}-t_q).
\]
Since $\sum_{q>n}q^{-2}=1/n-1/(2n^{2})+\mathcal O(n^{-3})$, the first two terms combine into
$-n/2+1/4+\mathcal O(n^{-1})$, and dividing by $n$ leaves
\begin{equation}\label{eq:ab_ingham_defect_split}
E_\Phi(n)=-\frac12+\sum_{q\le n}\frac1q\Bigl(\frac12-t_q\Bigr)-\frac1n\sum_{q\le n}\varpi(t_q)
+\rho_n,\qquad |\rho_n|\le\frac Cn .
\end{equation}
Fix $q$ and write $N=Lq+r$, where $L=\lfloor N/q\rfloor$ and $0\le r<q$. Over each complete
period the fractional parts $\{n/q\}$ run once through the values $a/q$ with $0\le a<q$, and
the sum of $\tfrac12-\{n/q\}$ is $\tfrac12$. The remaining $r$ ranks contribute
$r/2-r(r+1)/(2q)=r(q-r-1)/(2q)$. Hence
\[
\sum_{n\le N}\left(\frac12-\left\{\frac nq\right\}\right)
=\frac L2+\frac{r(q-r-1)}{2q}.
\]
Moreover the same sum over $1\le n<q$ is zero. It may therefore be inserted when the finite
sums in the middle term of \eqref{eq:ab_ingham_defect_split} are interchanged. Averaging over
$n\le N$ gives
\[
\frac1N\sum_{n\le N}\sum_{q\le n}\frac1q\Bigl(\frac12-t_q\Bigr)
=\frac1N\sum_{q\le N}\Bigl[\frac{\lfloor N/q\rfloor}{2q}
+\frac12\Bigl\{\frac Nq\Bigr\}\Bigl(1-\Bigl\{\frac Nq\Bigr\}-\frac1q\Bigr)\Bigr].
\]
The first part is $\tfrac12\sum_{q\le N}q^{-2}-\tfrac1{2N}\sum_{q\le N}\{N/q\}/q$, that is
$\zeta(2)/2+\mathcal O(N^{-1}\log N)$. The second is
$\tfrac1N\sum_q\varpi(\{N/q\})-\tfrac1{2N}\sum_q\{N/q\}/q$, which
Lemma~\ref{lem:ab_block_limit} applied to $\varpi$ evaluates as $s+\mathcal O(N^{-1/2})$ with
\[
s=\int_0^1\varpi(t)\,\psi'(1+t)\,dt=\int_0^1t\,\psi(1+t)\,dt
=1-\tfrac12\log(2\pi)=0.081061466\ldots,
\]
the middle equality by parts, the boundary terms vanishing because $\varpi(0)=\varpi(1)=0$ and
$\int_0^1\psi(1+t)\,dt=0$. For the last equality, Euler's reflection formula and
$\int_0^1\log(\sin\pi t)\,dt=-\log2$ give
\[
2\int_0^1\log\Gamma(t)\,dt
=\int_0^1\log\frac{\pi}{\sin\pi t}\,dt=\log(2\pi).
\]
Since $\Gamma(1+t)=t\Gamma(t)$ and $\int_0^1\log t\,dt=-1$, this yields
$\int_0^1\log\Gamma(1+t)\,dt=\tfrac12\log(2\pi)-1$.

The same lemma applied at each rank gives
$\tfrac1n\sum_{q\le n}\varpi(t_q)=s+\mathcal O(n^{-1/2})$ for the third term of
\eqref{eq:ab_ingham_defect_split}, whose average over $n\le N$ is
$s+\mathcal O(N^{-1/2})$. Finally,
$N^{-1}\sum_{n\le N}\rho_n=\mathcal O(N^{-1}\log N)=\mathcal O(N^{-1/2})$, and the other
$\mathcal O(N^{-1}\log N)$ terms above obey the same bound. Collecting the three, the two
occurrences of $s$ cancel and $-\tfrac12+\zeta(2)/2$ remains, with the error asserted in
\eqref{eq:ab_ingham_mean}.
\end{proof}

The cancellation of the two block constants is the precise sense in which the mean defect is
blind to the absolutely continuous layers here, and it does not require a separate dominated
convergence argument for their nonintegrable density.

\begin{numobs}[The mean defect]\label{numobs:ab_ingham}
The averages $N^{-1}\sum_{n\le N}E_\Phi(n)$ take the values $0.322069$, $0.322355$,
$0.322469$, $0.322464$ and $0.322465$ at the ranks $N=10^{3}$, $10^{4}$, $10^{5}$,
$10^{6}$ and $2\cdot10^{6}$, against $(\zeta(2)-1)/2=0.322467033$. Over the same range the defect itself
stays between $-1.26$ and $1.89$, well below the bound of \eqref{eq:ab_ingham_asym}.
\end{numobs}

Numerical Observation~\ref{numobs:ab_ingham} agrees with \eqref{eq:ab_ingham_mean} to six
places. The proof uses more than cancellation inside each layer. It also uses the regular
placement and shrinking amplitude of the Ingham layers. Proposition~\ref{prop:ab_layer_counterexample}
shows that zero mass and a summable diameter--variation product alone do not suffice. The
remaining regular-layer question is stated as Open Problem~\ref{op:ab_mean_defect}.

\section{The broken harmonic kernels at an integer scale}
\label{sec:ab_broken}

For $\lambda>1$ let $g_\lambda(x)=x\lambda^{\lfloor-\log_\lambda x\rfloor}$ be the self similar
broken harmonic profile of \eqref{eq:g_lambda_def}, with jump points $\lambda^{-i}$. Its jumps
are of constant size, so its total variation is infinite at every scale, and yet the kernel sum
is exactly computable at the powers of the scale when the scale is an integer.

\begin{proposition}\label{thm:ab_broken}
For every $\lambda>1$ the Abelian density is
\[g_\lambda^{*}(-1)=\int_0^1g_\lambda=\frac{\lambda+1}{2\lambda},
\]
the measure $\nu$ carries the point mass $-(\lambda-1)/\lambda$ at each point $\lambda^{-i}$
with $i\ge1$, and an absolutely continuous part of density $\lambda^{i-1}$ on
$(\lambda^{-i},\lambda^{-i+1})$. If moreover $\lambda$ is an integer then for every $j\ge0$
\begin{equation}\label{eq:ab_broken_exact}
S_{g_\lambda}(\lambda^{j})
=\frac{\lambda+1}{2\lambda}\,\lambda^{j}+(j+1)\,\frac{\lambda-1}{2\lambda}.
\end{equation}
In particular the Abelian defect of $g_\lambda$ is unbounded, and along the powers of the scale
it grows like $\dfrac{\lambda-1}{2\lambda\log\lambda}\,\log n$.
\end{proposition}

\begin{proof}
On $(\lambda^{-i},\lambda^{-i+1}]$ one has $g_\lambda(x)=\lambda^{i-1}x$, which gives the density
and, at $x=\lambda^{-i}$, the value $1$ against the right limit $1/\lambda$, hence the jump.
Integrating layer by layer,
\[
\int_0^1g_\lambda=\sum_{i\ge1}\lambda^{i-1}\frac{\lambda^{-2i+2}-\lambda^{-2i}}2
=\frac{1-\lambda^{-2}}2\sum_{i\ge1}\lambda^{-i+1}
=\frac{1-\lambda^{-2}}2\cdot\frac\lambda{\lambda-1}=\frac{\lambda+1}{2\lambda}.
\]
Condition \eqref{eq:abel_admissible} holds, the layer $i$ contributing
$\mathcal O(\lambda^{-i})$ to $\int t\,d|\nu|$.

Take now $\lambda$ an integer and $n=\lambda^{j}$, and evaluate
$\int_{(0,1)}(\{nt\}-\tfrac12)\,d\nu$ layer by layer, each layer being the interval
$(\lambda^{-i},\lambda^{-i+1})$ together with the point mass at its left endpoint.

For $1\le i\le j$ the map $t\mapsto\lambda^{j}t$ sends the layer onto
$(\lambda^{j-i},\lambda^{j-i+1})$, an interval with integer endpoints and integer length, so the
integral of $\{\lambda^{j}t\}-\tfrac12$ over the layer vanishes, a whole number of periods being
traversed. At the point $\lambda^{-i}$, $\lambda^{j}\lambda^{-i}=\lambda^{j-i}$ is an integer, so
$\{\lambda^{j}\lambda^{-i}\}-\tfrac12=-\tfrac12$, and the layer contributes
$(-\tfrac12)\cdot(-(\lambda-1)/\lambda)=(\lambda-1)/(2\lambda)$. There are $j$ such layers.

For $i>j$ the whole layer satisfies $\lambda^{j}t<1$, so $\{\lambda^{j}t\}=\lambda^{j}t$. The
absolutely continuous part contributes
$\lambda^{j-i-1}(\lambda^{2}-1)/2-(\lambda-1)/(2\lambda)$ and the point mass contributes
$-(\lambda-1)\lambda^{j-i-1}+(\lambda-1)/(2\lambda)$, so the layer contributes
$\lambda^{j-i-1}(\lambda-1)^{2}/2$. Summing over $i>j$,
\[
\frac{(\lambda-1)^{2}}2\sum_{r\ge1}\lambda^{-r-1}
=\frac{(\lambda-1)^{2}}2\cdot\frac1{\lambda(\lambda-1)}=\frac{\lambda-1}{2\lambda}.
\]
Adding the two contributions gives $j(\lambda-1)/(2\lambda)+(\lambda-1)/(2\lambda)$, and
inserting it in \eqref{eq:abelian_identity} gives \eqref{eq:ab_broken_exact}. Reading $j$ as
$\log n/\log\lambda$ gives the growth rate.
\end{proof}

\begin{remark}[Where the computation fails at $\sqrt2$]\label{rem:ab_sqrt2}
The proof used the integrality of $\lambda$ once, at the point where $\lambda^{j-i}$ was declared
an integer, and the whole contribution of the layers $i\le j$ rests on that. For
$\lambda=\sqrt2$ the quantity $\lambda^{j-i}=2^{(j-i)/2}$ is an integer only when $j-i$ is even,
so half of the jumps sit at points whose orbit under multiplication by $n$ is an irrational
rotation, and no exact evaluation of the layer is available. This is the alternative that governs
the transparency frontier\index[terms]{transparency frontier} of $g_{\sqrt2}$ in Chapter~\ref{chap:diophantine}, met here on the
Abelian side, where no Tauberian hypothesis is in play and where the failure is a failure of
exact summation rather than of a Green estimate\index[terms]{Green estimate}\index[names]{Green, G.}.
\end{remark}

\begin{numobs}[Broken harmonic kernels]\label{numobs:ab_broken}
Formula \eqref{eq:ab_broken_exact} is reproduced to machine precision at $\lambda=2$, $3$, $5$
and $0\le j\le7$. At $\lambda=2$ the defect reaches $(j+1)/4$ at $n=2^{j}$ and descends to
$-0.77$, $-1.01$, $-1.26$, $-1.75$ at $n=63$, $127$, $255$, $1023$, so that
$\max_{n\le N}|E|/\log N$ settles on $0.3616$ against the predicted $1/(4\log2)=0.3607$. At
$\lambda=\sqrt2$ the same ratio is near $0.27$ over $N\le1.6\cdot10^{4}$ and the exact formula
fails from $j=2$ on.
\end{numobs}

\section{The greatest common divisor kernel and Pillai's function}
\label{sec:ab_gcd}

The greatest common divisor kernel gives its row sum in closed form, through the function of
Pillai\index[names]{Pillai, S. S.}.

\begin{proposition}\label{prop:ab_gcd}
For the kernel $G(n,k)=\tfrac12(1+\gcd(n,k)/n)$ of Appendix~\ref{app:J} and every $n\ge1$,
\begin{equation}\label{eq:ab_gcd}
S_G(n)=\frac n2+\frac{P(n)}{2n},
\qquad
P(n)=\sum_{d\mid n}d\,\varphi(n/d).
\end{equation}
The Abelian density is $G^{*}(-1)=\tfrac12$, and the defect $P(n)/(2n)$ is unbounded and
depends only on the factorization of $n$.
\end{proposition}

\begin{proof}
The identity $\gcd(n,k)=\sum_{d\mid n,\,d\mid k}\varphi(d)$ gives
$\sum_{k\le n}\gcd(n,k)=\sum_{d\mid n}\varphi(d)n/d=P(n)$ after replacing $d$ by $n/d$.
This proves \eqref{eq:ab_gcd}, and the transform of the kernel is the constant $1/2$.
Also
\[
 \frac{P(n)}n=\sum_{e\mid n}\frac{\varphi(e)}e
 =\prod_{p^{a}\parallel n}\Bigl(1+a\Bigl(1-\frac1p\Bigr)\Bigr),
\]
the product form coming from the multiplicativity of $e\mapsto\varphi(e)/e$ and from
$\sum_{i=1}^{a}\varphi(p^{i})/p^{i}=a(1-1/p)$. For a squarefree product $n=p_1\cdots p_r$ this equals
$\prod_{j=1}^r(2-1/p_j)\ge(3/2)^r$, and is therefore unbounded. This also proves directly
the asserted dependence on the prime factorization. The function $P(n)/n$ is multiplicative,
and the defect $P(n)/(2n)$ is half of it, with the value $\tfrac12$ at $n=1$. The function $P$ is the one studied by
\nm{Pillai}{S. S.}~\cite{Pillai1933}, and its historical name is not needed for the conclusion.
\end{proof}

The kernel of Appendix~\ref{app:J} is the one whose index is out of reach because the section
family $(A_d)_{d\ge2}$ is not controlled by $A$ alone. Its Abelian reading is by contrast exact
and elementary, and the multiplicative structure that blocks the Tauberian direction is the same
one that makes the Abelian direction transparent.

\begin{numobs}[The gcd kernel]\label{numobs:ab_gcd}
Identity \eqref{eq:ab_gcd} is exact at every rank tested, the largest discrepancy\index[terms]{discrepancy} over
$n\in\{12,60,97,360,1024,2310\}$ being $7\cdot10^{-12}$ in double precision. The defect is
$1-1/(2p)$ at every prime, since $P(p)=2p-1$, and it grows with the number of prime factors of
$n$ counted with multiplicity, each factor $p^{a}$ contributing $1+a(1-1/p)\ge1+a/2$ to $P(n)/n$.
\end{numobs}

\section{Two intermediate regimes}
\label{sec:ab_intermediate}

The defects met so far are zero, bounded, of order $\log n$, or half of a multiplicative function. Two kernels of
the gallery fall between the density and those, which shows that the second order of the Abelian
expansion is not confined to one scale.

\begin{proposition}\label{prop:ab_sqrt}
For the square root kernel $G(n,k)=\tfrac12\bigl(1+(1+\sqrt k)/(1+\sqrt n)\bigr)$ of
Appendix~\ref{app:M},
\begin{equation}\label{eq:ab_sqrt}
S_G(n)=\frac56\,n+\frac{\sqrt n}6+\mathcal O(1),
\end{equation}
and $\tfrac56=G^{*}(-1)$ by the closed form $G^{*}(z)=(4z-1)/(2(2z-1))$.
\end{proposition}

\begin{proof}
Summing the kernel gives
$S_G(n)=\tfrac n2+\bigl(n+\sum_{k\le n}\sqrt k\bigr)/\bigl(2(1+\sqrt n)\bigr)$, and
$\sum_{k\le n}\sqrt k=\tfrac23n^{3/2}+\tfrac12n^{1/2}+\mathcal O(1)$ by Euler-Maclaurin.
Expanding the quotient in powers of $n^{-1/2}$ gives $\tfrac n3+\tfrac{\sqrt n}6+\mathcal O(1)$,
and adding $n/2$ gives \eqref{eq:ab_sqrt}. At $z=-1$ the closed form is $(-5)/(2\cdot(-3))=5/6$.
\end{proof}

A kernel damped by a slowly growing gauge produces a defect of intermediate order.

\begin{proposition}\label{prop:ab_gauged}
Let $G(n,k)=1+h(k/n)/L_n$ with $h$ admissible for Theorem~\ref{thm:abelian_identity} and
$L_n\to\infty$. Then
\begin{equation}\label{eq:ab_gauged}
S_G(n)=n+\frac{h^{*}(-1)}{L_n}\,n+\frac{E_h(n)}{L_n},
\end{equation}
an identity with no error term. For the self gauged logarithmic kernel
$G(n,k)=\log(n+k)/\log(2n)$ of Appendix~\ref{app:L}, where $h(t)=\log\frac{1+t}2$ and
$L_n=\log(2n)$, this reads
\[S_G(n)=n-\frac{(1-\log2)\,n}{\log(2n)}+\mathcal O\Bigl(\frac1{\log n}\Bigr).
\]
\end{proposition}

\begin{proof}
Summing $G$ term by term gives $S_G(n)=n+L_n^{-1}S_h(n)$, and
Corollary~\ref{cor:density_defect} applied to $h$ gives \eqref{eq:ab_gauged}. For the logarithmic
kernel, $\int_0^1\log\frac{1+t}2\,dt=\log2-1$ and $h$ has bounded variation, so $E_h$ is
bounded.
\end{proof}

\begin{numobs}[The two intermediate regimes]\label{numobs:ab_middle}
For the square root kernel, $(S_G(n)-\tfrac56n)/\sqrt n$ takes the values $0.169121$,
$0.167481$, $0.166928$, $0.166798$ at $n=10^{3}$, $10^{4}$, $10^{5}$, $4\cdot10^{5}$, against
$\tfrac16$. For the self gauged logarithmic kernel, $(S_G(n)-n)\log(2n)/n$ takes the values
$-0.306506$, $-0.306818$, $-0.306849$ at $n=10^{3}$, $10^{4}$, $10^{5}$, against
$\log2-1=-0.306853$.
\end{numobs}

\section{Abelian reading of the gallery}
\label{sec:ab_gallery}

Every entry of the gallery has an Abelian density, and it is the value of its transform at $-1$.
What separates the entries is the defect.

\begin{center}
\begin{tabular}{lcc>{\raggedright\arraybackslash}p{6.4cm}}
\toprule
Kernel & $G^{*}(-1)$ & kernel sum & defect \\
\midrule
$(1-\lambda)x+\lambda$ & $\frac{1+\lambda}2$ & exact & zero \\
$\mathbf 1_{(0,a]}$ & $a$ & $\lfloor na\rfloor$ & bounded, mean $\frac1{2q}$ or $0$ \\
$\frac{n^{2}+k}{n^{2}+n}$ & $1$ & $n-\frac12+\frac1{n+1}$ & bounded, limit $-\frac12$ \\
$1-\{x\}$ & $\frac12$ & $\frac{n+1}2$ & constant $\frac12$, from the diagonal jump \\
$\frac12\bigl(1+\frac{1+\sqrt k}{1+\sqrt n}\bigr)$ & $\frac56$ & asymptotic & $\frac{\sqrt n}6$ \\
$\frac{\log(n+k)}{\log 2n}$ & $1$ & asymptotic & $-\frac{(1-\log2)n}{\log 2n}$ \\
$x\lfloor1/x\rfloor$ & $\frac{\zeta(2)}2$ & $\frac1n\sum_{m\le n}\sigma(m)$ & $\mathcal O(\log n)$, arithmetic \\
$g_\lambda$, $\lambda\in\N$ & $\frac{\lambda+1}{2\lambda}$ & exact at $n=\lambda^{j}$ & of order $\log n$ along $\lambda^{j}$ \\
$g_{\sqrt2}$ & $\frac{\sqrt2+1}{2\sqrt2}$ & no closed form & boundedness not established, mixed rational and irrational rotations \\
$\frac12\bigl(1+\frac{\gcd(n,k)}n\bigr)$ & $\frac12$ & $\frac n2+\frac{P(n)}{2n}$ & $\tfrac12\prod_{p^{a}\parallel n}\bigl(1+a(1-\tfrac1p)\bigr)$ \\
\bottomrule
\end{tabular}
\end{center}

The density is analytic and the defect is arithmetic. A kernel of finite variation has a bounded
defect whose mean reads the denominators of its rational jumps, by
Theorem~\ref{thm:mean_defect}. Infinite variation removes that theorem and decides nothing by
itself, and in the three kernels worked above the defect carries the arithmetic of the jump set,
the divisor structure for the Ingham kernel, the
$\lambda$-adic expansion of $n$ for the broken harmonics at an integer scale, and the
multiplicative structure of $n$ for the gcd kernel. In no case does the defect determine the
regularity index, and in no case is a Tauberian hypothesis needed to compute it. The next
chapter takes the same forward direction through the fractional part sums, where a gauge is
inserted between the two arguments and the error term reaches the Dirichlet divisor problem\index[terms]{divisor problem}.

\section{What remains open here}
\label{sec:ab_open}

\begin{proposition}[Zero-mass layers alone do not force a mean]
\label{prop:ab_layer_counterexample}
There is a bounded left continuous profile $g$ satisfying the hypotheses of
Theorem~\ref{thm:abelian_identity}, with infinite total variation, whose Stieltjes measure is a
locally finite sum
\[
\nu=\sum_{j\ge1}\nu_j
\]
of measures of mass zero supported on pairwise disjoint intervals. If $\delta_j$ and $v_j$ are
the diameter and total variation of the $j$th layer, then
\[
\delta_j\longrightarrow0,
\qquad
\sum_{j\ge1}\delta_jv_j<\infty.
\]
Every point mass of $\nu$ sits at an irrational point, so the rational jump sum in
\eqref{eq:mean_defect} is zero, but the Ces\`aro mean of
\[
D_g(n):=S_g(n)-g^{*}(-1)n
\]
does not converge.
\end{proposition}

\begin{proof}
Fix irrational numbers $0<a<b<1$. Choose integers $M_j\to\infty$ so rapidly that the
intervals
\[
I_j=\left(\frac{a}{M_j},\frac{b}{M_j}\right]
\]
are pairwise disjoint, and put
\[
g(x)=\sum_{j\ge1}\mathbf 1_{I_j}(x).
\]
The profile is bounded and left continuous, and it has bounded variation on every
$[\delta,1]$. Its Stieltjes layers are
\[
\nu_j=\delta_{a/M_j}-\delta_{b/M_j}.
\]
They have mass zero, diameter $(b-a)/M_j$ and variation two. Taking the $M_j$ at least
geometrically increasing gives
\[
\int_{(0,1)}t\,d|\nu|(t)=(a+b)\sum_{j\ge1}\frac1{M_j}<\infty,
\qquad
\sum_{j\ge1}\delta_jv_j=2(b-a)\sum_{j\ge1}\frac1{M_j}<\infty.
\]
The total variation is nevertheless $\sum_j2=\infty$, and every point mass sits at an irrational point.

Write
\[
D_j(n)=\left\{\frac{na}{M_j}\right\}-
       \left\{\frac{nb}{M_j}\right\}.
\]
Theorem~\ref{thm:abelian_identity}, applied layer by layer, gives
$D_g(n)=\sum_jD_j(n)$. For each fixed $j$, Weyl equidistribution at the two irrational
endpoints gives
\[
\frac1N\sum_{n\le N}D_j(n)\longrightarrow0.
\]
On the other hand, if $n\le M_j$, both arguments of the fractional parts are smaller than one,
and hence
\[
D_j(n)=-\frac{n(b-a)}{M_j},
\qquad
\frac1{M_j}\sum_{n\le M_j}D_j(n)\longrightarrow-\frac{b-a}{2}.
\]

A diagonal choice makes these two behaviors coexist for the full sum. More precisely, choose
integers
\[
M_1<L_1<M_2<L_2<\cdots
\]
so that at $N=M_j$ the mean of $\sum_{i<j}D_i$ has modulus at most $2^{-j}$, while at
$N=L_j$ the mean of $\sum_{i\le j}D_i$ has modulus at most $2^{-j}$. This is possible by the
fixed-layer limits above. The next scales may simultaneously be chosen so large that
\[
\max(M_j,L_j)\sum_{i>j}\frac1{M_i}\longrightarrow0.
\]
For $N\in\{M_j,L_j\}$ and $i>j$, the interval $I_i$ lies below the smallest grid point
$1/N$, so its entire mean contribution has modulus at most
\[
\frac{(b-a)(N+1)}2\sum_{i>j}\frac1{M_i}=o(1).
\]
Consequently the full mean tends to $-(b-a)/2$ along $N=M_j$ and to zero along $N=L_j$.
It therefore does not converge.
\end{proof}

\begin{openproblem}[Mean defects for regular zero-mass layers]\label{op:ab_mean_defect}
Let $g$ satisfy the hypotheses of Theorem~\ref{thm:abelian_identity}, and suppose its Stieltjes
measure has a locally finite decomposition $\nu=\sum_m\nu_m$, where $\nu_m$ has mass zero and
is supported on an interval $I_m$ tending to the origin. Put
\[
\delta_m=\operatorname{diam}(I_m),
\qquad
v_m=|\nu_m|(I_m),
\]
and assume $\sum_m\delta_mv_m<\infty$. The layer defect is unambiguously
\[
D_g(n):=S_g(n)-g^{*}(-1)n-c_1
=\sum_m\int_{I_m}\left(\{nt\}-\tfrac12\right)d\nu_m(t),
\]
because every layer has mass zero and Theorem~\ref{thm:abelian_identity} makes the equivalent
uncentered integral absolutely convergent.

What additional regularity on the positions and variations of the layers guarantees
\[
\lim_{N\to\infty}\frac1N\sum_{n\le N}D_g(n)
=-\frac12\sum_{\substack{p/q\in\Q\cap(0,1)\\(p,q)=1}}
\frac{\nu(\{p/q\})}{q}?
\]
The series on the right is already absolutely convergent under
\eqref{eq:abel_admissible}, since $1/q\le p/q$ for $p\ge1$. In particular, does the conclusion
follow if $g(0^{+})$ exists and $v_m\to0$? Can it be proved for the Ingham-scale class in which,
with $r_m=\sup I_m$,
\[
v_m=\mathcal O(r_m),
\qquad
\delta_m=\mathcal O(r_m^2)?
\]
The Ingham function belongs to this class with $r_m=1/m$,
$\delta_m=1/(m(m+1))$ and $v_m=2/(m+1)$, and its mean formula is
Theorem~\ref{thm:ab_ingham_mean}. Its proof uses exact period sums for the family
$\{n/q\}$, and a general proof would require a short-interval equidistribution estimate strong
enough to exclude the moving microscopic layers of
Proposition~\ref{prop:ab_layer_counterexample}.
\end{openproblem}
Read forward, the equation asks nothing. The kernel applied to the constant sequence returns a
density, and every order of the answer is unconditional. Where the arithmetic enters is the
defect at the second order, and there it is the jump set of the profile that carries it. What is
not settled is how far that reading extends. Open Problem~\ref{op:ab_riemann_discrepancy} asks
which admissible kernels outside the smooth class have a second-order discrepancy, and Open
Problem~\ref{op:ab_mean_defect} asks what regularity gives a mean formula when the Stieltjes
layers have zero mass and infinite variation. The next chapter keeps the direction and inserts a
gauge between the two arguments, and the arithmetic that appears is classical.

\chapter{Fractional-part sums and gauge signatures}
\label{chap:abelian}

The preceding chapter read the forward direction for an arbitrary kernel and stopped at the
second order, where the defect of the Abelian expansion\index[terms]{Abelian defect} was shown to carry the arithmetic of
the jump set. This chapter takes the same direction and inserts a gauge\index[terms]{gauge} between the two
arguments of the kernel, which is where that arithmetic becomes classical.

The object is the fractional part\index[terms]{fractional part} sum $\sum_{k=1}^{n}\{f(n)/f(k)\}$ of a gauge $f$, and the
constants of its main term are signatures of the gauge. For the exponential gauges
$f(n)=m^{n}+1$ the sum reduces exactly to the Dirichlet divisor problem\index[terms]{divisor problem}, so the optimal error
exponent is the divisor exponent itself and no improvement of the one is available without an
improvement of the other. For the Fibonacci\index[terms]{Fibonacci sequence}\index[names]{Fibonacci} gauge the two parity constants are established, the
odd remainder reducing to a dyadic difference of the Gauss circle error term and the even
remainder to a dyadic combination of the divisor remainder.

The second half of the chapter needs no gauge at all. The weight $1-\{x\}$ of
Appendix~\ref{app:O}, read on the ratios $k/n$, is the one whose index is a quarter, and read
on the reciprocal ratios $n/k$ it produces the summatory divisor function\index[terms]{divisor function}. The two objects are
the two triangular halves of one matrix, and the sections that take them apart, in the discrete
operator, in the Mellin transform and in the Bessel skeleton, are where the Abelian side of the
theory is at its most arithmetic. Sums of fractional parts of this kind are studied for their
own sake by Balazard\index[names]{Balazard, M.}, Benferhat\index[names]{Benferhat, L.} and Bouderbala\index[names]{Bouderbala, M.} \cite{BalazardBB2021}, and the classical divisor
and circle estimates used below are treated by Ivi\'c\index[names]{Ivi\'c, A.} \cite[Ch.~13 and~14]{Ivic2003}, where the
divisor problem\index[terms]{divisor problem} comes first and the circle problem\index[terms]{circle problem} is set beside it as its classical analogue.

\section{The exponential gauge and the Dirichlet divisor problem}
\label{sec:divisor_exponent}

For every integer $m\ge2$ the fractional part sum of the exponential gauge reduces
exactly to the summatory divisor function. Write
\[
D(x):=\sum_{\ell\le x}\tau(\ell)=x\log x+(2\gamma-1)x+\Delta(x)
\]
for the summatory divisor function and its remainder, and
$n_{\mathrm{odd}}:=n/2^{v_2(n)}$ for the odd part of $n$, $v_2$ being the $2$-adic
valuation. The reduction shows that the exponent problem for these sums is not an
analogue of the Dirichlet divisor problem\index[terms]{divisor problem} but the divisor problem\index[terms]{divisor problem} itself.

\begin{proposition}\label{prop:frac_divisor_reduction}
Let $m\ge2$ be a fixed integer and
\[
S_m(n):=\sum_{k=1}^{n}\left\{\frac{m^n+1}{m^k+1}\right\}.
\]
Then, for every $\eps>0$,
\[
\boxed{\;
S_m(n)=(\log 2)\,n+\Delta(n)-2\Delta(n/2)-\tau(n_{\mathrm{odd}})
+\mathcal{O}_{m,\eps}(n^{\eps}).\;}
\]
Equivalently, $S_m(n)=D(n)-2D(n/2)-\tau(n_{\mathrm{odd}})+\mathcal{O}_{m,\eps}(n^{\eps})$.
\end{proposition}

\begin{proof}
For $1\le k\le n$ write $n=qk+r$ with $q=\lfloor n/k\rfloor$ and $0\le r<k$. Since
$m^{k}\equiv-1\pmod{m^{k}+1}$, one has
$m^{n}+1\equiv(-1)^{q}m^{r}+1\pmod{m^{k}+1}$, so that
\[
\left\{\frac{m^n+1}{m^k+1}\right\}=
\begin{cases}
\dfrac{m^{r}+1}{m^{k}+1}, & q\ \text{even},\\[7pt]
0, & q\ \text{odd and }r=0,\\[5pt]
1-\dfrac{m^{r}-1}{m^{k}+1}, & q\ \text{odd and }r>0.
\end{cases}
\]
Set $\eta_k:=\mathbf 1_{\{q\ \mathrm{odd},\,r>0\}}$. In every case
$\bigl|\{(m^n+1)/(m^k+1)\}-\eta_k\bigr|\le 2\,m^{-(k-r)}$. Put $h=k-r=(q+1)k-n$. For a
fixed value of $h$ the last identity forces $k\mid n+h$, so the number of indices $k$
with $k-r=h$ is at most $\tau(n+h)$, and
\[
\sum_{k\le n}m^{-(k-r)}\le\sum_{1\le h\le n}m^{-h}\,\tau(n+h)\ll_{m,\eps}n^{\eps}.
\]
Hence $S_m(n)=\sum_{k\le n}\eta_k+\mathcal{O}_{m,\eps}(n^{\eps})$. The indices removed
from the odd quotients by the condition $r>0$ are exactly the divisors $k$ of $n$ with
$n/k$ odd, and there are $\tau(n_{\mathrm{odd}})$ of them, so
\[
\sum_{k\le n}\eta_k=\sum_{k\le n}\mathbf 1_{\{\lfloor n/k\rfloor\ \mathrm{odd}\}}
-\tau(n_{\mathrm{odd}}).
\]
Counting the odd quotients through $\#\{k\le n:\lfloor n/k\rfloor\ge q\}=\lfloor n/q\rfloor$
gives
\[
\sum_{k\le n}\mathbf 1_{\{\lfloor n/k\rfloor\ \mathrm{odd}\}}
=\sum_{j\ge1}\Bigl(\Bigl\lfloor\frac n{2j-1}\Bigr\rfloor-\Bigl\lfloor\frac n{2j}\Bigr\rfloor\Bigr)
=\sum_{j\ge1}(-1)^{j+1}\Bigl\lfloor\frac nj\Bigr\rfloor
=D(n)-2D(n/2),
\]
and $D(n)-2D(n/2)=(\log2)n+\Delta(n)-2\Delta(n/2)$ follows from the definition of
$\Delta$.
\end{proof}

Under the divisor conjecture the exponent above becomes sharp.

\begin{corollary}\label{cor:divisor_quarter}
The Dirichlet divisor conjecture
\[
\Delta(x)=\mathcal{O}_{\eps}\bigl(x^{1/4+\eps}\bigr)
\]
implies
\[
S_m(n)=(\log2)\,n+\mathcal{O}_{m,\eps}\bigl(n^{1/4+\eps}\bigr).
\]
Conversely, this estimate for one fixed integer $m\ge2$ implies the Dirichlet divisor
conjecture. The optimal exponent in the fractional part sum is therefore the Dirichlet
divisor exponent.
\end{corollary}

\begin{proof}
Only the converse requires comment. Set $E(x):=\Delta(x)-2\Delta(x/2)$. By
Proposition~\ref{prop:frac_divisor_reduction} and $\tau(n_{\mathrm{odd}})=n^{o(1)}$, a
bound $S_m(n)-(\log2)n=\mathcal{O}_{m,\eps}(n^{\theta+\eps})$ implies
$E(n)=\mathcal{O}_{\eps}(n^{\theta+\eps})$. For every integer $J\ge1$, iterating the
definition of $E$ gives
\[
\Delta(n)=2^{-J}\Delta(2^{J}n)-\sum_{j=1}^{J}2^{-j}E(2^{j}n),
\]
and the classical bound $\Delta(x)=\mathcal{O}(x^{1/2})$ makes the first term vanish
as $J\to\infty$, so that
\[
\Delta(n)=-\sum_{j\ge1}2^{-j}E(2^{j}n).
\]
For every $\theta<1$ the assumed bound on $E$ then gives
$\Delta(n)=\mathcal{O}_{\eps}(n^{\theta+\eps})$. Taking $\theta=1/4$ proves the
converse.
\end{proof}

\begin{remark}[Best unconditional exponent]\label{rem:huxley_exponent}
The classical bound $\Delta(x)=\mathcal{O}(x^{1/2})$ already proves
Proposition~\ref{prop:exp_gauge_abelian} for every $m\ge2$, and the published estimate
of Huxley\index[names]{Huxley, M. N.} \cite{Huxley2003} yields
\[
S_m(n)=(\log2)\,n+\mathcal{O}_{m,\eps}\bigl(n^{131/416+\eps}\bigr),\qquad
\tfrac{131}{416}=0.314903\ldots
\]
The preprint of Li\index[names]{Li, X.} and Yang\index[names]{Yang, X.} \cite{LiYang2023} announces the improved exponent
$\theta_{*}=0.314483175974\ldots$, the current preprint record.
\end{remark}
\begin{remark}[Optimality]
Hardy's omega theorem for the divisor remainder, carried through the dyadic summation
of the preceding proof, prevents an admissible exponent strictly smaller than $1/4$.
The exponent $1/4$ is the optimal conjectural exponent here, not an exponent suggested
by analogy.
\end{remark}

\subsection{A general quotient profile mechanism}
\label{sec:quotient_profile}

The preceding argument suggests a general mechanism. Let $f:\N\to\N$ be increasing and
\[
S_f(n):=\sum_{k\le n}\left\{\frac{f(n)}{f(k)}\right\}.
\]
Assume that there exists a bounded quotient profile $\omega:\N\to\R$ with
$\omega(0):=0$ such that
\begin{equation}\label{eq:quotient_profile}
S_f(n)=\sum_{k\le n}\omega\!\left(\Bigl\lfloor\frac nk\Bigr\rfloor\right)
+\mathcal{O}_{\eps}(n^{\eps}),
\end{equation}
the error absorbing the exceptional divisibility or resonance terms. Summation over the
quotient blocks gives the exact identity
\[\sum_{k\le n}\omega\!\left(\Bigl\lfloor\frac nk\Bigr\rfloor\right)
=\sum_{q\ge1}\omega(q)\Bigl(\Bigl\lfloor\frac nq\Bigr\rfloor-\Bigl\lfloor\frac n{q+1}\Bigr\rfloor\Bigr)
=\sum_{q\ge1}\bigl(\omega(q)-\omega(q-1)\bigr)\Bigl\lfloor\frac nq\Bigr\rfloor.
\]
Suppose in addition that $\omega(q)-\omega(q-1)=(\mathbf 1*h)(q)$ with $h$ finitely supported
and
\begin{equation}\label{eq:filter_zero}
\sum_{d}\frac{h(d)}{d}=0.
\end{equation}
Then
\[
\begin{aligned}
\sum_{q\ge1}\bigl(\omega(q)-\omega(q-1)\bigr)\Bigl\lfloor\frac nq\Bigr\rfloor
&=\sum_{d}h(d)\,D(n/d)=C_h\,n+\sum_{d}h(d)\,\Delta(n/d),\\
C_h&:=-\sum_{d}\frac{h(d)\log d}{d},
\end{aligned}
\]
and consequently
\begin{equation}\label{eq:divisor_filter}
\boxed{\;
S_f(n)=C_h\,n+\sum_{d}h(d)\,\Delta(n/d)+\mathcal{O}_{\eps}(n^{\eps}).\;}
\end{equation}
In Dirichlet series language, put $H(s):=\sum_d h(d)\,d^{-s}$. The coefficients behind
\eqref{eq:divisor_filter} have Dirichlet series $\zeta(s)^{2}H(s)$, and condition
\eqref{eq:filter_zero} reads $H(1)=0$. It cancels one order of the double pole of
$\zeta(s)^{2}$, leaving a simple pole whose residue is $H'(1)=C_h$, so the linear main
term is produced by a zero of the dilation filter at $s=1$.

For $f(n)=m^n+1$ the quotient profile is $\omega(q)=\mathbf 1_{\{q\ \mathrm{odd}\}}$,
apart from the resonant divisor terms, so
$\omega(q)-\omega(q-1)=(-1)^{q+1}=1-2\cdot\mathbf 1_{\{2\mid q\}}$, corresponding to $h(1)=1$,
$h(2)=-2$, $H(s)=1-2^{1-s}$ and $H'(1)=\log2$.

\subsection{Divisor filters generated by fractional parts}
\label{sec:divisor_filters}

The reduction above is proved under the quotient profile hypothesis
\eqref{eq:quotient_profile}. Without a growth or structural restriction that
hypothesis admits artificial recursive realizations, as Proposition~\ref{prop:quotient_profile_crt}
will show. The arithmetic problem is to find and classify structured increasing
functions $f$ realizing nontrivial profiles $\omega$. For each integer $p\ge2$ a
natural target filter is
\[
h_p(1)=\frac1{p-1},\qquad h_p(p)=-\frac{p}{p-1},
\]
which gives
\[
a_p(q)=\frac{1-p\,\mathbf 1_{\{p\mid q\}}}{p-1},\qquad
\omega_p(q)=\frac{q\bmod p}{p-1},\qquad
C_p=\frac{\log p}{p-1}.
\]
Any natural function $f_p$ satisfying
$\{f_p(n)/f_p(k)\}\approx(\lfloor n/k\rfloor\bmod p)/(p-1)$ with total defect
$n^{o(1)}$ would therefore satisfy
\[
\sum_{k\le n}\left\{\frac{f_p(n)}{f_p(k)}\right\}
=\frac{\log p}{p-1}\,n+\frac{\Delta(n)-p\,\Delta(n/p)}{p-1}+n^{o(1)}.
\]
For $p=2$ this is the family $f(n)=m^n+1$.

\begin{proposition}[Universal realization without structural control]
\label{prop:quotient_profile_crt}
For every function $\omega:\N\to[0,1]$ there is a strictly increasing prime-valued
function $f:\N\to\N$, with $f(n)>2^n$, such that
\[
\left|
\sum_{k\le n}\left\{\frac{f(n)}{f(k)}\right\}
-\sum_{k\le n}\omega\!\left(\left\lfloor\frac nk\right\rfloor\right)
\right|<2
\qquad(n\ge1).
\]
In particular, every profile $\omega_p$ above is realizable in the sense of
\eqref{eq:quotient_profile}, with its error improved to $\mathcal O(1)$.
\end{proposition}

\begin{proof}
Choose $f(1)>2$ prime. Suppose that the distinct primes $f(1),\ldots,f(n-1)$
have been chosen. For each $k<n$, select an integer
$r_{n,k}\in\{1,\ldots,f(k)-1\}$ such that
\[
\left|\frac{r_{n,k}}{f(k)}
-\omega\!\left(\left\lfloor\frac nk\right\rfloor\right)\right|
\le \frac1{f(k)}.
\]
The Chinese remainder theorem\index[terms]{Chinese remainder theorem} gives a
residue class $R_n$ modulo $M_{n-1}:=\prod_{k<n}f(k)$ satisfying
$R_n\equiv r_{n,k}\pmod{f(k)}$ for every $k<n$. Since no $r_{n,k}$ vanishes,
$(R_n,M_{n-1})=1$. Dirichlet's theorem on primes in arithmetic
progressions\index[names]{Dirichlet, P. G. L.}\index[terms]{prime in an arithmetic progression},
Proposition~\ref{prop:pnt_progressions} with $Q=M_{n-1}$, therefore permits the choice of a prime
\[
f(n)\equiv R_n\pmod{M_{n-1}},
\qquad
f(n)>\max\{f(n-1),2^n\}.
\]
For $k<n$ this congruence gives
\[
\left\{\frac{f(n)}{f(k)}\right\}=\frac{r_{n,k}}{f(k)}.
\]
The term $k=n$ contributes zero on the fractional-part side and at most one on
the profile side. Hence
\[
\left|
\sum_{k\le n}\left\{\frac{f(n)}{f(k)}\right\}
-\sum_{k\le n}\omega\!\left(\left\lfloor\frac nk\right\rfloor\right)
\right|
\le 1+\sum_{k<n}\frac1{f(k)}
<1+\sum_{k\ge1}2^{-k}=2.
\]
\end{proof}

Proposition~\ref{prop:quotient_profile_crt} shows that
\eqref{eq:quotient_profile} alone carries no rigidity when unrestricted growth is
allowed. The content of the realization question begins with a quantitative or
algebraic constraint on the gauge.

\begin{openproblem}[Structured realization of the divisor filters]\label{op:divisor_filters}
For $p\ge3$, exhibit a strictly increasing function $f_p:\N\to\N$ whose quotient
profile is $\omega_p$ in the sense of \eqref{eq:quotient_profile} and which has
finite exponential type,
\[
\log f_p(n)=\mathcal O(n).
\]
More restrictively, can such a realization be chosen to satisfy a fixed linear
recurrence with integer coefficients, as $m^n+1$ does for $p=2$, or can either
possibility be ruled out? The recursive prime construction of
Proposition~\ref{prop:quotient_profile_crt} deliberately supplies no such control.
\end{openproblem}

Formula \eqref{eq:divisor_filter} represents filtered divisor remainders through
fractional part observables rather than through the classical Voronoi-Bessel
expansion, which may suggest methods based on quotient blocks, digital recurrences,
dilation dynamics, positivity, periodic profiles, or Ramanujan expansions. For the
original family one may write
\[
S_m(n)+\tau(n_{\mathrm{odd}})-(\log2)\,n=(I-2U_2)\Delta(n)+n^{o(1)},
\qquad (U_2F)(n):=F(n/2),
\]
so the fractional part sum is a dyadic dilation difference of the divisor remainder,
and any direct improvement of its pointwise error transfers to the classical divisor
problem. The warning is that rewriting the sum in terms of $\Delta$ does not by itself
improve the divisor estimate. A new argument must exploit structure visible before the
reduction, a recurrence, positivity, monotonicity, automatic structure, averaging over
parameters, or a family of jointly invertible dilation filters.

\section{The two triangles of the diagonal-jump kernel}
\label{sec:two_triangles}

The preceding sections reach the divisor problem\index[terms]{divisor problem} by deforming the coordinates until the
divisor remainder appears in a fractional part sum. There is a second entrance, and it needs
no deformation at all. The profile $g(x)=1-\{x\}$ of Appendix~\ref{app:O}, read on the ratios
$k/n$, is the function of good variation whose index is one quarter by
Theorem~\ref{thm:O_quarter}. Read on the reciprocal ratios $n/k$, the same profile produces the
summatory divisor function. The two objects are the two triangular halves of one matrix, and
the six sections that follow take them apart in three independent coordinates, the discrete
operator, the Mellin transform, and the Bessel skeleton.

Extend $g$ to the positive half line by periodicity of the fractional part, so that
$g(x)=1-\{x\}$ for every $x>0$, and consider the full ratio matrix
\[\mathbb H_{n,k}=g\!\Big(\frac kn\Big)\qquad(n,k\ge1).
\]
The equation of Appendix~\ref{app:O} uses its lower triangle $k\le n$ row by row. The divisor
sums use the opposite triangle, since $g(n/k)=\mathbb H_{k,n}$ for $k\le n$.

\section{The lower triangle as an additive Volterra operator}

Let $J$ be the summatory matrix $J_{n,k}=\mathbf 1_{k\le n}$, let $I$ be the identity, and let
$D=\operatorname{diag}(1,2,3,\dots)$.

The operator has a matrix form on the lower triangle.

\begin{proposition}\label{prop:tt_matrix}
On the lower triangle,
\begin{equation}\label{eq:tt_volterra}
\mathbb H_-=I+J-D^{-1}JD .
\end{equation}
\end{proposition}

\begin{proof}
The entry of $D^{-1}JD$ at $(n,k)$ is $k/n$ for $k\le n$, so $J-D^{-1}JD$ has entry $1-k/n$
strictly below the diagonal and zero on it. The diagonal value $g(1)=1$ supplies $I$.
\end{proof}

Formula \eqref{eq:tt_volterra} is the additive face of the kernel. The whole difficulty of
Appendix~\ref{app:O} sits in the isolated $I$, the jump that the first two terms do not see.

\section{The reciprocal profile is an affine transform of the Ingham kernel}

The reciprocal profile and the Ingham function are one object read twice.

\begin{proposition}\label{prop:tt_reciprocal}
For every $0<x\le1$, with $\check g(x)=g(1/x)$ and $\Phi(x)=x\lfloor1/x\rfloor$,
\[x\,\check g(x)=x+\Phi(x)-1 .
\]
\end{proposition}

\begin{proof}
Since $\{1/x\}=1/x-\lfloor1/x\rfloor$, one has
$x(1-\{1/x\})=x-1+x\lfloor1/x\rfloor$.
\end{proof}

The single rupture of the lower triangle, carried by the value at $x=1$, unfolds under
reciprocal reflection into every discontinuity of the Ingham kernel. This is the first exact
location of the divisor algebra inside the bridge, and it explains why the two windows on the
same profile carry such different arithmetic.

\section{The terminal jump is the divisor signal}

The role of the value at the integers is isolated by a parameter. For $\lambda\neq0$ set
\[g_\lambda(x)=
\begin{cases}
1-\{x\}, & x\notin\N,\\
\lambda, & x\in\N .
\end{cases}
\]
On the lower triangle only the point $x=1$ is affected. On the reciprocal triangle every
divisibility hit is affected, and the difference is exactly the divisor function.

\begin{proposition}\label{prop:tt_atom}
For every integer $n\ge1$,
\[\sum_{k\le n}g_\lambda(n/k)=\sum_{k\le n}g(n/k)+(\lambda-1)\,\tau(n),
\qquad\text{so}\qquad
\frac{\partial}{\partial\lambda}\sum_{k\le n}g_\lambda(n/k)=\tau(n).
\]
\end{proposition}

\begin{proof}
The two summands differ precisely when $n/k$ is an integer, that is when $k$ divides $n$, and
there are $\tau(n)$ such indices.
\end{proof}

The same parameter deforms the Laguerre\index[terms]{Laguerre} oscillator on the Tauberian side, without moving the
exponent. Let $(a_n)$ solve $\sum_{k\le n}a_k\,g_\lambda(k/n)=n^{-\beta}$ and let $A$ be
its partial sums. The off-diagonal part is affine, so the reduction of
Proposition~\ref{prop:O_structure} becomes $\lambda a_n+\frac1n\sum_{j<n}A(j)=n^{-\beta}$,
and differencing after multiplication by $n$ gives, for $n\ge2$,
\begin{equation}\label{eq:tt_lambda_rec}
n\,A(n)-\bigl(2n-1-\lambda^{-1}\bigr)A(n-1)+(n-1)A(n-2)
=\lambda^{-1}d_\beta(n),
\end{equation}
with $d_\beta(n)=n^{1-\beta}-(n-1)^{1-\beta}$. The homogeneous part is the ordinary Laguerre\index[terms]{Laguerre}
recurrence at argument $\lambda^{-1}$, so its first solution is $u_n=L_n(\lambda^{-1})$, and
Fej\'er's fixed argument expansion gives
\[L_n(\lambda^{-1})=\frac{e^{1/(2\lambda)}\lambda^{1/4}}{\sqrt\pi}\,n^{-1/4}
\Bigl(\cos\Bigl(2\sqrt{n/\lambda}-\frac\pi4\Bigr)+\mathcal O_\lambda(n^{-1/2})\Bigr).
\]
The height of the rupture changes the frequency of the mode and leaves its exponent at one
quarter. At $\lambda=1$, Theorem~\ref{thm:O_quarter} is the complete forced statement.

Proposition~\ref{prop:tt_atom} and the identification $u_n=L_n(\lambda^{-1})$ are the sharpest
elementary form of the correspondence. One and the same terminal jump is the local divisor
signal on the reciprocal side and the Laguerre\index[terms]{Laguerre} oscillator on the triangular side.

\section{The reciprocal row sum is a regularized square of zeta}

For a bounded measurable $f$ on $(0,1]$ define the direct sampling operator on $x\ge1$ by
$(\mathcal Sf)(x)=\sum_{k\le x}f(k/x)$.

The sampling operator has a Mellin transform in closed form.

\begin{lemma}\label{lem:tt_sampling}
For $\Re s>1$ and bounded $f$,
\[\int_1^{\infty}(\mathcal Sf)(x)\,x^{-s-1}\,dx=\zeta(s)\int_0^1 f(u)\,u^{s-1}\,du .
\]
\end{lemma}

\begin{proof}
Interchange the sum and the integral, then substitute $u=k/x$ in each term, which turns
$\int_k^{\infty}f(k/x)x^{-s-1}\,dx$ into $k^{-s}\int_0^1f(u)u^{s-1}\,du$.
\end{proof}

So does the reciprocal profile, and the zeta function appears there.

\begin{proposition}\label{prop:tt_mellin}
For $\Re s>1$,
\begin{equation}\label{eq:tt_mellin_check}
\int_0^1\check g(x)\,x^{s-1}\,dx=\frac1s\Bigl(\zeta(s)-\frac1{s-1}\Bigr),
\end{equation}
the right side being regular at $s=1$ with value $\gamma$, so that
$\int_0^1\bigl(1-\{1/x\}\bigr)\,dx=\gamma$. Consequently, with
$R(x)=\sum_{k\le x}g(x/k)$,
\begin{equation}\label{eq:tt_factorisation}
\int_1^{\infty}R(x)\,x^{-s-1}\,dx=\frac{\zeta(s)}{s}\Bigl(\zeta(s)-\frac1{s-1}\Bigr).
\end{equation}
\end{proposition}

\begin{proof}
The substitution $y=1/x$ turns the first integral into $\int_1^\infty(1-\{y\})y^{-s-1}\,dy$, and
$1-\{y\}=1-y+\lfloor y\rfloor$ with $\int_1^\infty\lfloor y\rfloor y^{-s-1}\,dy=\zeta(s)/s$ gives
$1/s-1/(s-1)+\zeta(s)/s$, which is the stated form. The Laurent expansion\index[terms]{Laurent expansion}\index[names]{Laurent, P. A.}
$\zeta(s)=1/(s-1)+\gamma+\mathcal O(s-1)$ gives the value at $s=1$. Since
$R=\mathcal S\check g$, formula \eqref{eq:tt_factorisation} follows from
Lemma~\ref{lem:tt_sampling}.
\end{proof}

Identity \eqref{eq:tt_factorisation} shows where each zeta factor comes from. The reciprocal
unfolding of Proposition~\ref{prop:tt_reciprocal} produces one, the Abelian sweep over the
integer lattice produces the other, and the subtraction of the polar part cancels the double
pole at $s=1$ and leaves a simple pole of residue $\gamma$, in agreement with the value of the
integral. The degree two object behind the divisor problem is thus assembled from two
elementary operations on one profile.

On the arithmetic side the same sum is elementary. With $N=\lfloor x\rfloor$ and $H_N$ the
harmonic number,
\[R(x)=D(x)+N-xH_N,
\qquad
R(x)-\gamma x=\Delta(x)+(N-x)+x\bigl(\log x+\gamma-H_N\bigr),
\]
the last two terms being bounded. At the integers this is
Proposition~\ref{prop:O_harmonic}. The conjectural estimate
$\Delta(x)=\mathcal O_\eps(x^{1/4+\eps})$ is therefore equivalent to the same discrepancy\index[terms]{discrepancy}
estimate for the Riemann sums of the reciprocal profile,
\[\sum_{k\le x}\check g(k/x)-x\int_0^1\check g(u)\,du=\mathcal O_\eps\bigl(x^{1/4+\eps}\bigr).
\]

\section{One Bessel skeleton on both sides}

The Green basis of Appendix~\ref{app:O} carries a differential operator that the discrete
recurrence hides. Let $u_n,v_n$ be the two solutions of \eqref{eq:O_homog} normalized as in
Lemma~\ref{lem:O_basis}.

\begin{proposition}\label{prop:tt_borel}
For $t>0$,
\begin{equation}\label{eq:tt_borel_u}
\sum_{n\ge0}u_n\frac{t^{n}}{n!}=e^{t}J_0(2\sqrt t),
\end{equation}
and, writing $\sum_{n\ge0}v_nt^{n}/n!=e^{t}F(t)$, the function $F$ solves
\begin{equation}\label{eq:tt_bessel_forced}
t\,F''+F'+F=e^{-t},
\qquad F(0)=0,\ F'(0)=1,
\end{equation}
whose homogeneous basis is $J_0(2\sqrt t)$ and $Y_0(2\sqrt t)$, of Wronskian\index[terms]{Wronskian} $1/(\pi t)$ in the
variable $t$. Variation of constants gives
\[F(t)=\pi\Bigl[Y_0(2\sqrt t)\int_0^{t}e^{-s}J_0(2\sqrt s)\,ds
-J_0(2\sqrt t)\int_0^{t}e^{-s}Y_0(2\sqrt s)\,ds\Bigr].
\]
The two limiting Laplace integrals are
\begin{equation}\label{eq:tt_laplace}
\int_0^{\infty}e^{-s}J_0(2\sqrt s)\,ds=e^{-1},
\qquad
\int_0^{\infty}e^{-s}Y_0(2\sqrt s)\,ds=\frac{e^{-1}}{\pi}\operatorname{Ei}(1).
\end{equation}
\end{proposition}

\begin{proof}
Identity \eqref{eq:tt_borel_u} is the exponential generating function of the Laguerre\index[terms]{Laguerre}
polynomials at argument one \cite[Eq.~18.12.14]{DLMF}. Multiplying \eqref{eq:O_homog} by
$t^{n-1}/(n-1)!$ and summing gives $tX''+(1-2t)X'+tX=1$ for
$X(t)=\sum v_nt^{n}/n!$, and the substitution $X=e^{t}F$ reduces it to
\eqref{eq:tt_bessel_forced}. The equation $tF''+F'+F=0$ is Bessel's equation of order zero in
the variable $2\sqrt t$, and the Wronskian\index[terms]{Wronskian} follows from $W\{J_0,Y_0\}(w)=2/(\pi w)$ together
with $dw/dt=1/\sqrt t$. The Laplace integrals are classical.
\end{proof}

The second Laplace integral of \eqref{eq:tt_laplace} accounts for the constant
$\operatorname{Ei}(1)$ appearing in the Fej\'er expansion of $v_n$ in
Lemma~\ref{lem:O_basis}. That constant is not an artefact of the saddle point\index[terms]{saddle point method} calculation, it
is the Laplace transform\index[terms]{Laplace transform}\index[names]{Laplace, P.-S.} of the second Bessel solution.

The divisor remainder is governed by the same equation. Vorono\"i\index[terms]{Vorono\"i kernel}\index[names]{Vorono\"i, G.}'s formula expresses
$\Delta$ through Bessel functions of order one at the argument $4\pi\sqrt{mx}$
\cite{Jutila2015,EndresSteiner}, whose leading truncated form is
\begin{equation}\label{eq:tt_voronoi}
\Delta(x)=\frac{x^{1/4}}{\pi\sqrt2}\sum_{m\le M}\frac{\tau(m)}{m^{3/4}}
\cos\Bigl(4\pi\sqrt{mx}-\frac\pi4\Bigr)+\mathcal R_M(x).
\end{equation}
The two orders are one differentiation apart, since
\[\frac{d}{dt}Y_0(2\sqrt t)=-\frac{Y_1(2\sqrt t)}{\sqrt t},
\]
so the oscillator that governs the divisor remainder is the derivative of the one that occurs
in the Borel transform\index[terms]{Borel transform}\index[names]{Borel, E.} of the Green function. The parallel is a common Bessel skeleton and not
a numerical coincidence, as the following table records.

\medskip
\renewcommand{\arraystretch}{1.35}
{\small
\begin{tabular}{p{3.1cm}p{5.0cm}p{5.0cm}}
\toprule
& Triangular side, Appendix~\ref{app:O} & Reciprocal side, divisor problem \\
\midrule
Discrete object & inverse of $I+J-D^{-1}JD$ & direct reciprocal row sum \\
Transform & Borel transform of the Laguerre basis & Mellin transform $\zeta(s)\bigl(\zeta(s)-\frac1{s-1}\bigr)/s$ \\
Oscillator & $J_0,Y_0$ at argument $2\sqrt t$ & $Y_1,K_1$ at argument $4\pi\sqrt{mx}$ \\
Envelope & $n^{-1/4}$ & $x^{1/4}m^{-3/4}$ \\
Phase & $2\sqrt n-\pi/4$ & $4\pi\sqrt{mx}-\pi/4$ \\
Status & proved, Theorem~\ref{thm:O_quarter} & conjectural upper bound \\
\bottomrule
\end{tabular}
}
\renewcommand{\arraystretch}{1}
\medskip

\section{What separates the two quarters}

Theorem~\ref{thm:O_quarter} controls the response to one power forcing. Variation of constants
reduces it to two scalar connection series, and the square root phase makes their convergence
abscissa exactly one quarter. In \eqref{eq:tt_voronoi} the critical oscillator is instead
synthesized over every $m$ with coefficient $\tau(m)m^{-3/4}$, a sequence that is not
absolutely summable. Bounding each mode separately destroys the cancellation on which the
conjecture depends, so the quarter theorem does not imply the divisor estimate, and no
implication is claimed. What the structural identities of this section establish is that the
two problems share one profile, one Mellin transform built from two elementary operations, and one
Bessel Green operator. Their difference is global aggregation.

\paragraph{Research direction (from scalar connection to divisor synthesis).}
Proposition~\ref{prop:tt_borel} and the Vorono\"i expansion exhibit a common Bessel skeleton.
The passage from order zero to order one is explicit, and the unresolved difficulty lies in the
global synthesis of the modes weighted by $\tau(m)m^{-3/4}$. A transference theorem
would have to specify its source and target normed spaces, an exact intertwining identity, a
bound uniform in the truncation parameters $x$ and $M$, and a proof of the divisor-weighted
boundedness independent of any estimate equivalent to the Dirichlet divisor conjecture. Until
these objects are identified, this remains a research programme rather than a formal
conjecture.

\begin{numobs}\label{numobs:tt_check}
The Mellin transform \eqref{eq:tt_mellin_check} agrees with piecewise integration between the
jumps of the reciprocal profile to $5\cdot10^{-16}$ at $s=4$, and the value $\gamma$ of the
integral is met to the truncation error. The Laplace integrals \eqref{eq:tt_laplace} agree with
quadrature at twenty five digits, the second to a measured ratio of one at fifteen decimals.
The identification $u_n=L_n(\lambda^{-1})$ holds with residual zero in exact rational
arithmetic for $\lambda\in\{1,\tfrac12,\tfrac32,2\}$, and the forced recurrence
\eqref{eq:tt_lambda_rec} reproduces forward substitution in the defining equation to
$10^{-14}$.
\end{numobs}

\section{The Fibonacci gauge}
\label{sec:fibonacci_gauge}
\label{sec:fib_gauge_signature}

The Fibonacci sequence gives the gauge $f(n)=F_{n}$, of irrational growth ratio $\varphi = (1+\sqrt{5})/2$. The Abelian analog here is fully proved and reveals a sharp dichotomy between odd and even indices.

\begin{theorem}[Harcos, {\cite{HarcosMO2024}}]\label{thm:fibonacci_frac}
Let $F_{n}$ denote the $n$-th Fibonacci number. Then:
\begin{equation}\label{eq:fibonacci_frac}
\sum_{k=1}^n \left\{\frac{F_{n}}{F_{k}}\right\} =
\begin{cases}
\dfrac{\pi}{8}\,n + \mathcal{O}(\sqrt{n}), & n \text{ odd},\\[6pt]
\dfrac{3\log 2}{4}\,n + \mathcal{O}(\sqrt{n}), & n \text{ even}.
\end{cases}
\end{equation}
\end{theorem}

The sum was the subject of a question on MathOverflow, and the theorem with the proof below
is the answer given there by Gergely Harcos\index[names]{Harcos, G.}~\cite{HarcosMO2024}.

The bound $\mathcal{O}(\sqrt n)$ is the one supplied by the balancing in the proof below. It is far
from optimal, and the next remark identifies what the two remainders really are.

\begin{proof}
The proof rests on the following Lucas identity. For Fibonacci numbers $F$ and Lucas numbers $L$,
\[
F_{u}L_{v} = \begin{cases} F_{u+v} + (-1)^v F_{u-v}, & u \ge v,\\ F_{u+v} - (-1)^u F_{v-u}, & u \le v. \end{cases}
\]
Plugging $(mk, mk \pm r)$ for $(u,v)$ yields the congruences modulo $F_{k}$ (using $F_{k} \mid F_{mk}$):
\begin{align*}
F_{2mk+r} &\equiv F_{r}(-1)^{mk} \pmod{F_{k}},\\
F_{2mk-r} &\equiv F_{r}(-1)^{mk-r-1} \pmod{F_{k}}.
\end{align*}
These determine $\{F_{n}/F_{k}\}$ exactly in each case. Setting $I(n,t) := \N \cap (n/(t+1),\, n/t)$, the cases $n = 2mk+r$ correspond to $k \in I(n,2m)$ and $n = 2mk-r$ to $k \in I(n,2m-1)$.

The contributions of each interval, sorted by parity of $m$, are:
\begin{itemize}
\item odd $m$: both $I(n,2m)$ and $I(n,2m-1)$ contribute $\tfrac{1}{2}|I(n,\cdot)| + \mathcal{O}(1)$;
\item even $m$: $I(n,2m)$ contributes $\mathcal{O}(1)$. For $I(n,2m-1)$, the contribution is $\mathcal{O}(1)$ if $n$ is odd and $|I(n,2m-1)| + \mathcal{O}(1)$ if $n$ is even.
\end{itemize}
Fixing $M \ge 1$ and summing over $m \in \{1,\ldots,M\}$, the partial sum over $k \ge \lceil n/(2M+1)\rceil$ equals, up to $\mathcal{O}(M)$:
\begin{itemize}
\item If $n$ is odd: $\tfrac{1}{2}\sum_{\substack{1\le m\le M\\ m\text{ odd}}} (n/(2m-1) - n/(2m+1)) = \tfrac{\pi}{8}n + \mathcal{O}(n/M)$.
\item If $n$ is even: $\tfrac{1}{2}\sum_{\substack{m\text{ odd}}} (\cdots) + \sum_{\substack{m\text{ even}}} (n/(2m-1) - n/(2m)) = \tfrac{3\log 2}{4}n + \mathcal{O}(n/M)$.
\end{itemize}
The contribution of $k < n/(2M+1)$ is $\mathcal{O}(n/M)$. Choosing $M = \lfloor\sqrt{n}\rfloor$ balances the errors to $\mathcal{O}(M + n/M) = \mathcal{O}(\sqrt{n})$.
\end{proof}

The remainder is not the end of the story. Identified exactly, it is the error term of a
classical problem, and which of the two depends on the parity.

\begin{theorem}[Exact remainders, {\cite[Theorem~1.1]{CloitreFibonacci}}]\label{thm:fibonacci_exact}
Let $\Delta$ be the Dirichlet divisor error of \S\ref{sec:divisor_exponent}, so that
$\sum_{\ell\le x}\tau(\ell)=x\log x+(2\gamma-1)x+\Delta(x)$, and let
$\Delta_C$ be the Gauss circle error, normalized as in Appendix~\ref{app:O} by
$N_C(x)=\pi x+\Delta_C(x)$, where $N_C(x)$ counts the lattice points of the disc of radius
$\sqrt x$. Then for every $\eps>0$,
\[
\sum_{k=1}^{n}\left\{\frac{F_n}{F_k}\right\}
=\frac{\pi}{8}\,n+\frac14\Bigl(\Delta_C(n)-\Delta_C\bigl(\tfrac n2\bigr)\Bigr)+\mathcal{O}(n^{\eps})
\qquad(n\ \text{odd},\ n\ge3),
\]
\[
\sum_{k=1}^{n}\left\{\frac{F_n}{F_k}\right\}
=\frac{3\log 2}{4}\,n+2\Delta\bigl(\tfrac n2\bigr)-5\Delta\bigl(\tfrac n4\bigr)
+2\Delta\bigl(\tfrac n8\bigr)+\mathcal{O}(n^{\eps})
\qquad(n\ \text{even},\ n\ge4).
\]
\end{theorem}

\begin{proofstatus}{Imported result. The two displayed exact remainder formulas are quoted
from \cite[Theorem~1.1]{CloitreFibonacci} and are not proved here. Their proof reduces the odd
remainder through Jacobi's two-square identity and the even remainder through an alternating
divisor identity. Within this volume they sharpen Theorem~\ref{thm:fibonacci_frac}, supply the
forward half of Corollary~\ref{cor:fibonacci_equivalence}, and provide the exact input interpreted
in Remarks~\ref{rem:fibonacci_profile_principle} and the two preceding remarks. No theorem
outside \S\ref{sec:fibonacci_gauge} depends on them.}
\end{proofstatus}

The identification turns a bound into an equivalence.

\begin{corollary}[The two parities against the two classical problems, {\cite{CloitreFibonacci}}]
\label{cor:fibonacci_equivalence}
Let $0\le\theta_C,\theta<1$. Then $\Delta_C(x)=\mathcal{O}(x^{\theta_C+\eps})$ and
$\Delta(x)=\mathcal{O}(x^{\theta+\eps})$ for every $\eps>0$ if and only if
\[
\sum_{k=1}^{n}\left\{\frac{F_n}{F_k}\right\}=\frac{\pi}{8}\,n+\mathcal{O}(n^{\theta_C+\eps})
\quad(n\ \text{odd}),
\qquad
\sum_{k=1}^{n}\left\{\frac{F_n}{F_k}\right\}=\frac{3\log 2}{4}\,n+\mathcal{O}(n^{\theta+\eps})
\quad(n\ \text{even}).
\]
The equivalence carries no hypothesis. In particular the conjecture that $\tfrac14$ is the
optimal exponent in the Gauss circle problem and in the Dirichlet divisor problem is the
conjecture that $\tfrac14$ is the optimal remainder exponent for the Fibonacci sum, separately
along each parity.
\end{corollary}

\begin{proofstatus}{Imported result. The forward implications follow from
Theorem~\ref{thm:fibonacci_exact} and the triangle inequality. The converse implications are
proved in \cite[Proposition~6.1]{CloitreFibonacci}: the circle error is recovered by dyadic
iteration, while the divisor error is recovered from the two-scale recurrence for
$B(m)=\Delta(2m)-2\Delta(m)$. The corollary is used here to identify the two optimal remainder
exponents and to interpret the remarks that follow, and no later proof assumes it.}
\end{proofstatus}

The two error terms have already been met in this volume, on a different object. The
fractional part kernel of Appendix~\ref{app:O}, with its diagonal jump, produces the same pair, the
divisor error on one member of its family and the Gauss circle error on another, and
Corollary~\ref{cor:O_indices} reads the index off them there. A gauge on the Fibonacci sequence
and a jump on the diagonal of a kernel arrive at the same two classical remainders.

\begin{remark}[The two constants and the gauge]
The constants $\pi/8 \approx 0.3927$ and $3\log 2/4 \approx 0.5199$ are not arithmetic accidents. They arise from summing the series $\sum_{m \text{ odd}} (1/(2m-1) - 1/(2m+1)) = \pi/4$ and the related alternating logarithmic series, after accounting for the parity splitting. The error term is not of order $\sqrt n$, and it is not a signature of the index. The odd remainder reduces to a dyadic difference of the Gauss circle error term and the even remainder to a dyadic combination of the Dirichlet divisor error term, so the optimal exponent is the circle exponent along the odd indices and the divisor exponent along the even ones, which is the content of the companion paper~\cite{CloitreFibonacci}. Both are $\mathcal{O}(n^{131/416+\eps})$ by Huxley\index[names]{Huxley, M. N.}, and the preprint of Li\index[names]{Li, X.} and Yang\index[names]{Yang, X.} announces $0.3144831759741\ldots$, exactly as recorded in Remark~\ref{rem:huxley_exponent} for the exponential gauge. The obstruction is cancellation in exponential sums, not the zero-free region\index[terms]{zero-free region} of an $L$-function. Hardy's omega theorems give $\Omega(n^{1/4})$ for the two classical remainders themselves \cite{Hardy1916}, and no formal transfer of that lower bound to their dyadic differences is available, so the value $1/4$ along each parity is expected rather than forced. What is certain is that the exponent is not $1/2$. Definition~\ref{def:O_delta} and Remark~\ref{rem:O_fibonacci} place the two parities in the same frame as the two classical problems.
\end{remark}

The pairing between a parity and a classical problem is not a property of the Fibonacci
sequence. It belongs to the family of second order recurrences, and inside that family it can
be exchanged.

\begin{remark}[The Lucas companion exchanges the two parities, {\cite[Theorem~7.1]{CloitreFibonacci}}]
\label{rem:lucas_parity_exchange}
Fix an integer $c\ge1$ and let $F^{(c)}$ and $L^{(c)}$ obey the same recurrence
$X_{j+2}=cX_{j+1}+X_j$ with the initial values $F^{(c)}_0=0$, $F^{(c)}_1=1$ and
$L^{(c)}_0=2$, $L^{(c)}_1=c$, so that $c=1$ gives the Fibonacci and Lucas
sequences\index[terms]{Lucas sequence}\index[names]{Lucas, E.}. The sums
$\sum_{k\le n}\{F^{(c)}_n/F^{(c)}_k\}$ obey the two formulas of
Theorem~\ref{thm:fibonacci_exact} for every $c$, so the parameter is invisible to the
asymptotics. The companion sums $\mathcal{L}_c(n)=\sum_{k\le n}\{L^{(c)}_n/L^{(c)}_k\}$
reverse the assignment. For every $\eps>0$,
\[
\mathcal{L}_c(n)=\frac{3\log2}{4}\,n+\Delta(n)-3\Delta\bigl(\tfrac n2\bigr)
+3\Delta\bigl(\tfrac n4\bigr)-2\Delta\bigl(\tfrac n8\bigr)+\mathcal{O}(n^{\eps})
\qquad(n\ \text{odd}),
\]
\[
\mathcal{L}_c(n)=\frac{\pi}{8}\,n+\frac14\Delta_C\bigl(\tfrac n2\bigr)
+\mathcal{O}(n^{\eps})
\qquad(n\ \text{even}).
\]
The odd indices now carry the divisor error and the even indices the circle error, the
opposite of Theorem~\ref{thm:fibonacci_exact}. Nothing analytic produces the exchange. The
Lucas congruence modulo $L^{(c)}_k$ carries the sign $(-1)^{s(k+1)}$ where the Fibonacci one
carries $(-1)^{sk}$, which shifts by one the residue class of $\lfloor n/k\rfloor$ that puts a
fractional part near one. The periodic profile changes, and with it the divisor sum attached
to its first difference.
\end{remark}

\begin{remark}[What decides, and where it stops]\label{rem:fibonacci_profile_principle}
The two preceding statements share a mechanism rather than a computation. A sum
$\sum_{k\le n}\{A_n/A_k\}$ reduces, up to $\mathcal{O}(n^{\eps})$, to
$\sum_{k\le n}\omega(\lfloor n/k\rfloor)$ for a periodic $\omega$ as soon as the sequence has a
residue law modulo $A_k$, and the telescoping identity
\[
\sum_{k\le x}\omega\Bigl(\Bigl\lfloor\frac xk\Bigr\rfloor\Bigr)
=\sum_{q\le x}\bigl(\omega(q)-\omega(q-1)\bigr)\Bigl\lfloor\frac xq\Bigr\rfloor
=\sum_{j\le x}\ \sum_{q\mid j}\bigl(\omega(q)-\omega(q-1)\bigr),
\qquad\omega(0)=0,
\]
then turns the first difference of $\omega$ into a divisor
sum. The divisor function, its alternating analogue and the character $\chi_4$ arise this way,
which is why the divisor and the circle problems appear and why the choice between them is a
congruence and not an estimate. The same mechanism is what
Appendix~\ref{app:O} runs on the diagonal of a kernel, and
Corollary~\ref{cor:O_indices} reads an index off the resulting error terms. What is needed is
congruence information about the sequence, not growth. Order three is already beyond the
argument, since the companion matrix of the Tribonacci
sequence\index[terms]{Tribonacci sequence} is three by three and its reduction modulo $T_k$ is
not a scalar sign. Computation of the exact remainders for $n\le4000$ puts
$\sum_{k\le n}\{T_n/T_k\}$ near $n/2$ with no visible separation of the three classes modulo
three, and no asymptotic formula is claimed there \cite{CloitreFibonacci}.
\end{remark}

\begin{remark}[Odd-even dichotomy as a possible Abelian signature]
The split between $\pi/8$ (odd $n$) and $3\log 2/4$ (even $n$) in \eqref{eq:fibonacci_frac} is a direct Abelian manifestation of the parity structure visible in the phase-transition table of Chapter~\ref{chap:gauge_ingham}. There, for $\beta > 1/2$, one has $\ell(i) + \ell(i+4) = 0$ (antisymmetry). At $\beta = 1/2$ the symmetry breaks to $\ell(i) + \ell(i+4) = 1$. The two fractional part constants appear to encode the same phase transition on the Abelian side.
\end{remark}

\section{Further exponential-gauge sums}

For $f(x) = 2^x+1$, define $W_m(n) := \sum_{k=1}^n k^m \Phi(f(k)/f(n))$.

\begin{proposition}[Polynomial correction and parity law]\label{prop:Wm_polynomial}
For every fixed integer $m\geq 0$, define
\[
 c_r:=\sum_{j=1}^{\infty}\frac{j^r}{2^j},
 \qquad
 P_m(X):=\sum_{j=1}^{\infty}\frac{(X-j)^m}{2^j}
       =\sum_{r=0}^{m}(-1)^r\binom{m}{r}c_rX^{m-r}.
\]
Then $P_m\in\Z[X]$ has degree $m$.  The integers $c_r$ form the entry A000629 of the On-Line
Encyclopedia of Integer Sequences\index[terms]{On-Line Encyclopedia of Integer Sequences}~\cite{OEIS}. In particular,
\[
 P_0(X)=1,\qquad P_1(X)=X-2,\qquad
 P_2(X)=X^2-4X+6,
\]
and $P_3(X)=X^3-6X^2+18X-26$.  Moreover,
\begin{align*}
 \lim_{n\to\infty}\frac{2^n}{n^m}
 \left(W_m(2n)-\sum_{k=1}^{2n}k^m+P_m(2n)\right)&=3,\\
 \lim_{n\to\infty}\frac{2^n}{n^m}
 \left(W_m(2n+1)-\sum_{k=1}^{2n+1}k^m+P_m(2n+1)\right)&=2.
\end{align*}
Consequently,
\[
 W_m(N)=\sum_{k=1}^{N}k^m-P_m(N)
        +\mathcal{O}_m\!\left(N^m2^{-N/2}\right).
\]
\end{proposition}

The sums $W_m$ were the subject of a question on MathOverflow, and the reduction to residues on
which the proof rests is the answer given there by Hung-Hsun Yu\index[names]{Yu, Hung-Hsun}~\cite{YuMO2024}.
The proof below includes the details and makes the final parity passage explicit.

\begin{proof}
Let $r_{N,k}$ be the least nonnegative residue of $2^N+1$ modulo $2^k+1$, so that
\[
 \Phi\!\left(\frac{2^k+1}{2^N+1}\right)
 =1-\frac{r_{N,k}}{2^N+1},
\]
and therefore
\[
 D_m(N):=\sum_{k=1}^{N}k^m-W_m(N)
 =\sum_{k=1}^{N-1}\frac{r_{N,k}}{2^N+1}\,k^m.
\]
For $N/2\leq k<N$, reduction of $2^N+1$ modulo $2^k+1$ gives
\[
 r_{N,k}=2^k+2-2^{N-k},
\]
whereas for $N/3<k<N/2$ it gives
\[
 r_{N,k}=2^{N-2k}+1.
\]
For $k\leq N/3$ the elementary bound $0\leq r_{N,k}<2^k+1$ suffices.
Replacing $(2^N+1)^{-1}$ by $2^{-N}$, retaining the first range and
absorbing the other terms therefore yields
\[
 D_m(N)=
 \sum_{1\leq j\leq \lfloor N/2\rfloor}
 \left(2^{-j}-2^{j-N}\right)(N-j)^m
 +\mathcal{O}_m\!\left(N^{m+1}2^{-2N/3}\right).
\]
By the definition of $P_m$ this becomes
\begin{equation}\label{eq:Wm_central_remainder}
 W_m(N)-\sum_{k=1}^{N}k^m+P_m(N)
 =\sum_{j=1}^{\infty}2^{-\max(j,N-j)}(N-j)^m
 +\mathcal{O}_m\!\left(N^{m+1}2^{-2N/3}\right).
\end{equation}

If $N=2n$, put $j=n+r$ in the main term of
\eqref{eq:Wm_central_remainder}.  Dominated convergence gives
\[
 \frac{2^n}{n^m}
 \sum_{j=1}^{\infty}2^{-\max(j,2n-j)}(2n-j)^m
 \longrightarrow \sum_{r\in\Z}2^{-|r|}=3.
\]
If $N=2n+1$, the same substitution gives
\[
 \frac{2^n}{n^m}
 \sum_{j=1}^{\infty}2^{-\max(j,2n+1-j)}(2n+1-j)^m
 \longrightarrow
 \sum_{r\in\Z}2^{-\max(r,1-r)}=2.
\]
The error in \eqref{eq:Wm_central_remainder} is $o(n^m2^{-n})$
on both subsequences, proving the two limits and the stated uniform order.

Finally, expansion of $(X-j)^m$ gives the displayed coefficient formula.
Here $c_0=1$, and for $r\geq1$ a shift of the geometric series gives
\[
 c_r=1+\sum_{s=0}^{r-1}\binom{r}{s}c_s.
\]
Thus every $c_r$ is an integer, and the leading coefficient of $P_m$ is
$c_0=1$.
\end{proof}

The fractional part sums $\sum_{k=1}^n \{(2^n+1)/(2^k+1)\}$ exhibit the same phenomenon. One can establish:
\[\sum_{n/2 < k \le n} \left\{\frac{2^n+1}{2^k+1}\right\} = \frac{n}{2} + \mathcal{O}(1),
\]
together with the exact expressions for the fractional part in two dyadic intervals:
\begin{align*}
\left\lfloor\frac{n+1}{2}\right\rfloor \le k \le n-1 &\implies \left\{\frac{2^n+1}{2^k+1}\right\} = \frac{2^k+2-2^{n-k}}{2^k+1}, \\
\left\lfloor\frac{n+3}{3}\right\rfloor \le k \le \left\lfloor\frac{n}{2}\right\rfloor &\implies \left\{\frac{2^n+1}{2^k+1}\right\} = \frac{2^{n-2k}+1}{2^k+1}. \end{align*}
These exact formulas provide the microscopic structure of the sum, but extending the dyadic decomposition to the full sum requires controlling sub-dyadic intervals.

\begin{proposition}\label{prop:exp_gauge_abelian}
\[
\sum_{k=1}^n \left\{\frac{2^n+1}{2^k+1}\right\} = (\log 2)\,n + \mathcal{O}(n^{1/2}).
\]
\end{proposition}

\begin{proof}
This is the case $m=2$ of Proposition~\ref{prop:frac_divisor_reduction} in
\S\ref{sec:divisor_exponent}, proved there, which identifies the error term exactly as
$\Delta(n)-2\Delta(n/2)-\tau(n_{\mathrm{odd}})+\mathcal{O}_{\eps}(n^{\eps})$ and
sharpens the bound to $\mathcal{O}_{\eps}(n^{131/416+\eps})$ unconditionally.
\end{proof}

The classical model is the fractional part sum $\sum_{k=1}^n \{n/k\} = (1-\gamma)n + \mathcal{O}(n^{1/2})$ of the standard divisor problem. The gauge changes the constant, $1-\gamma$ for the identity gauge, $\log 2$ for the exponential family, and \S\ref{sec:divisor_exponent} identifies the error term of the exponential family exactly as a dilation difference of the divisor remainder.
Two reductions carry the chapter. The exponential gauge sends the fractional part sum exactly to
the Dirichlet divisor problem, so its optimal error exponent is the divisor exponent and no
improvement of the one is available without an improvement of the other. The Fibonacci gauge
splits by parity, the odd remainder into a dyadic difference of the circle error term and the
even remainder into a dyadic combination of the divisor remainder. The weight of
Appendix~\ref{app:O}, read on the reciprocal ratios, then reaches the same problem from the other
side. The gauge moves the constant and not the difficulty. Open
Problem~\ref{op:divisor_filters} asks for a structured realization of the quotient profile at the
remaining filters.

\chapter{Finite trace formulas for polynomial kernels}
\label{chap:trace_poly}

A transform can also be read at all its zeros at once, and the object that holds them together is
the resolvent\index[terms]{resolvent} of the equation, the kernel $R$ inverting the operator, already met in
Chapter~\ref{chap:discrete_volterra}. A trace formula, in this part and throughout, means an
identity that writes $R$ as a sum over the zeros of the transform,
\[
R(y)=2\,\Re\sum_{n}\frac{y^{-\rho_n}}{g^{*\prime}(\rho_n)}+\text{a controlled remainder},
\]
each zero contributing one term and the derivative of the transform at that zero fixing its
amplitude. The name records the shape, a quantity attached to an operator expressed as a sum
over its spectrum.

Appendix~\ref{app:O} supplies a nonpolynomial counterpoint, in which direct arithmetic remainders and
the inverse Laguerre mode coexist in reciprocal readings of a completed profile, but
Proposition~\ref{prop:O_undetermined} shows that one reading does not determine the other.

Two clarifications are owed at once, because the term is loaded. Nothing here is a trace
formula in the sense of Selberg\index[names]{Selberg, A.}. There is no group, no quotient, and no geodesic, and no
such structure is claimed. The comparison to draw is the explicit formula\index[terms]{explicit formula} of Riemann\index[names]{Riemann, B.} and von
Mangoldt\index[names]{Mangoldt, H. von@von Mangoldt, H.}, which writes a counting function as a sum over the zeros of $\zeta$. Here the
transform $g^{*}$ takes the place of $\zeta$ and its zeros the place of the nontrivial zeros,
and the identity is exact with an explicit error rather than asymptotic.

What the reader needs is in this volume. The transform of Chapter~\ref{chap:fgv}, the
resolvent of Chapter~\ref{chap:discrete_volterra}, and nothing further. The two auxiliary
notions the chapters lean on, the indicial polynomial\index[terms]{indicial polynomial} whose roots are the exponents an equation
can carry and the Wronskian\index[terms]{Wronskian} that measures the independence of the modes, are defined where
they first appear.

One picture governs the whole chapter. The transform of a polynomial kernel is a rational
function, and its denominator is an indicial polynomial\index[terms]{indicial polynomial} whose roots are the exponents the
equation can carry. Each root produces a mode, a solution of the homogeneous equation behaving
like a power of the rank. The amplitude attached to each mode is fixed by a Wronskian\index[terms]{Wronskian}, the
discrete determinant that measures how independent the modes are, and the passage from these
modes back to the coefficients of the sequence is a transfer of singularities, read off the
generating function near its singular point. Roots, modes, Wronskians\index[terms]{Wronskian}, transfer. Every formula
below sits at one of those four stations.

For polynomial kernels the finite resolvent identity is exact at every degree. The affine kernel\index[terms]{affine kernel}
is the degree-one laboratory, the general polynomial kernel carries that exact identity, and the
quadratic kernel carries a complete spectral expansion, separated real modes, a complex pair, an
integer gap, a double mode, and the simple and double resonances, under the hypotheses stated for
each case. A complete spectral expansion at arbitrary fixed degree, organized by resonance
classes and carried through every prescribed coefficient order, stays open and closes the chapter.
The passage to a kernel with infinitely many spectral modes is the subject of
the next chapter.

\section{The affine finite trace formula}
\label{sec:trace_affine}

The affine kernel $g(x)=(1-\lambda)x+\lambda$ with $0<\lambda<1$ has arithmetic Mellin transform
$g^*(z)=(z-\lambda)/(z-1)$, a single zero at $z=\lambda$ and a pole at $z=1$. Its regularity index is
$\alpha(g)=\eta(g)=\lambda$, established in Appendix~\ref{app:A} by the Euler--Gauss\index[terms]{Euler--Gauss product}\index[names]{Gauss, C. F.} product, and its
closed-form partial sums are Theorem~\ref{thm:linear_exact_book}. The degree-one case carries the
finite resolvent identity in closed form for an arbitrary forcing, and the two forcings used later in
the chapter are its specializations. Write $M_\lambda(n)=\Gamma(n+1-\lambda)/\Gamma(n+1)\sim n^{-\lambda}$
for the discrete monomial.

\begin{proposition}\label{thm:trace-affine}
Let $R\colon\{1,2,\dots\}\to\mathbb C$ be arbitrary and let $(a_n)$ solve
$\sum_{k\le n}a_k\,g(k/n)=R(n)$ for $n\ge1$. Then the partial sums $A(n)=\sum_{k\le n}a_k$ satisfy,
for every $n\ge1$,
\begin{equation}\label{eq:trace-affine-general}
A(n)=R(n)+(1-\lambda)\,M_\lambda(n)\sum_{m=1}^{n-1}\frac{\Gamma(m+1)}{\Gamma(m+2-\lambda)}\,R(m),
\end{equation}
the sum being empty at $n=1$.
\end{proposition}

\begin{proof}
Writing $g(k/n)=\lambda+(1-\lambda)k/n$ and separating the two sums, the equation reads
$\lambda A(n)+(1-\lambda)n^{-1}\sum_{k\le n}k\,a_k=R(n)$. Eliminating the weighted sum at
ranks $n$ and $n-1$ gives
\[
 nA(n)=(n-\lambda)A(n-1)+nR(n)-(n-1)R(n-1),\qquad A(0)=0.
\]
Its homogeneous solutions are scalar multiples of $M_\lambda(n)$. In the stated domain
$0<\lambda<1$ the exact normalization is
\[
 M_\lambda(n)=\frac{\Gamma(n+1-\lambda)}{\Gamma(n+1)}
 =\Gamma(1-\lambda)\prod_{k=1}^{n}\left(1-\frac\lambda k\right).
\]
Variation of constants against this Gamma-normalized mode, with kernel
$\Gamma(m+1)/\Gamma(m+2-\lambda)=1/[(m+1-\lambda)M_\lambda(m)]$, telescopes to
\eqref{eq:trace-affine-general}. Positive integral values of $\lambda$ lie outside the stated
domain, and there the Gamma factor and the vanishing product require a separate limiting
normalization and are not asserted by this formula.
\end{proof}

Two specializations follow. With $R(n)=n^{-\beta}$ the identity becomes
\begin{equation}\label{eq:trace-affine-finite}
 \boxed{\,A(n)=n^{-\beta}+(1-\lambda)\,M_\lambda(n)\sum_{m=1}^{n-1}
 \frac{\Gamma(m+1)}{\Gamma(m+2-\lambda)}\,m^{-\beta}\,},
\end{equation}
and with $R(n)=M_\beta(n)$ it produces the discrete-monomial forcing used by the polynomial transfer
theorem below, so that Theorem~\ref{thm:trace-poly-finite} at $d=1$ and \eqref{eq:trace-affine-general}
agree on the affine kernel. The homogeneous mode is $M_\lambda(n)$ and the resolvent probe of
Chapter~\ref{chap:discrete_volterra} converges to $\widetilde R(s)=1-1/g^*(-s)=(\lambda-1)/(s+\lambda)$,
whose leftmost singularity sits at $s=-\lambda$, in the sign convention $z=-s$ of
Proposition~\ref{prop:ortho_volterra}.

For the forcing $n^{-\beta}$ the three regularity behaviors follow from \eqref{eq:trace-affine-finite}
by the ratio of the two rank scales. For $\beta<\lambda$ the forcing dominates and
$n^{\beta}A(n)\to 1/g^*(\beta)$, the transparent behavior. For $\beta>\lambda$ the homogeneous mode
dominates and $n^{\lambda}A(n)$ tends to a constant, the absorbed behavior. At the threshold
$\beta=\lambda$ the two scales coincide and $A(n)\sim(1-\lambda)\,n^{-\lambda}\log n$, the resonance.
The value of the index is read from these three behaviors and not from the transform, in keeping
with the standing convention of the volume. The degenerate exponents are exact, $A(n)=P_n$ at
$\beta=1$ and $A(n)=\lambda^{-1}+(1-\lambda^{-1})P_n$ at $\beta=0$, with $P_n=M_\lambda(n)/\Gamma(2-\lambda)$.

The finite resolvent form is the accessible face of the general identity. The polynomial theory that
follows carries it to every degree, and the quadratic theory exhibits the modal interactions that the
affine model cannot show.

\section{The polynomial kernel and its transform}
\label{sec:trace_poly_setup}

Let
\[
 g(t)=\sum_{j=0}^{d}c_jt^j,
 \qquad
 d\geq1,
 \qquad
 \sum_{j=0}^{d}c_j=1,
\]
and consider
\begin{equation}
 \sum_{k\leq n}a_k g\!\left(\frac{k}{n}\right)=M_\beta(n),
 \qquad
 A(n):=\sum_{k\leq n}a_k,
 \qquad n\geq1,
 \label{eq:trace-poly-triangular}
\end{equation}
where
\[
 M_\rho(n):=\frac{\Gamma(n+1-\rho)}{\Gamma(n+1)}
 \sim n^{-\rho}.
\]
The forcing is used only when every value required by
\eqref{eq:trace-poly-triangular} is finite.  In particular, positive
integers at least two require a different rank-one convention and are not
part of the present forcing statement.

The discrete mode is the monomial $M_\lambda(n)$ of
Section~\ref{sec:trace_affine}, and the normalized homogeneous solution is
\[
 P_n=\frac{M_\lambda(n)}{\Gamma(2-\lambda)}.
\]
No other notation for a discrete mode is used.

For \(\Re z<0\), define the arithmetic Mellin transform by
\[
 g^*(z):=-z\int_0^1g(t)t^{-z-1}\,dt.
\]
The rational function obtained from this integral will always be reduced
before its zeros are read.

The quadratic case is
\begin{equation}
 g(t)=a+bt+ct^2,
 \qquad
 a+b+c=1,
 \qquad
 c\neq0.
 \label{eq:trace-quad-kernel}
\end{equation}

\section{An exact resolvent at finite rank}
\label{sec:trace_poly_results}

For a polynomial profile the transform is a finite sum of simple fractions, and the resolvent
that inverts the triangular system closes after finitely many terms. The statement below records
both, exactly and at every rank.

\begin{theorem}
\label{thm:trace-poly-finite}
The arithmetic Mellin transform is
\begin{equation}
 g^*(z)=z\sum_{j=0}^{d}\frac{c_j}{z-j}.
 \label{eq:trace-poly-transform}
\end{equation}
Common factors in its numerator and denominator are to be cancelled.

For \(0\leq r\leq d-1\), set
\[
 S_r(n):=\sum_{k\leq n}k^rA(k),
 \qquad
 U_r(n):=\frac{S_r(n)}{n^r},
\]
and write \(U(n)=(U_0(n),\ldots,U_{d-1}(n))^{\mathsf T}\).
Define
\[
 q_r(n):=\sum_{j=r+1}^{d}c_j\binom jr n^{-j},
 \qquad
 \ell_{n,r}:=q_r(n)(n-1)^r,
 \qquad
 \ell_n:=(\ell_{n,0},\ldots,\ell_{n,d-1}).
\]
Then Abel summation gives the exact identity
\begin{equation}
 A(n)=M_\beta(n)+\sum_{r=0}^{d-1}q_r(n)S_r(n-1).
 \label{eq:trace-poly-abel}
\end{equation}
The exact moment system is
\begin{equation}
 U(n)=T_nU(n-1)+M_\beta(n)\mathbf 1,
 \qquad
 U(1)=M_\beta(1)\mathbf 1,
 \label{eq:trace-poly-system}
\end{equation}
where
\begin{equation}
 (T_n)_{r,s}
 =\delta_{r,s}\left(\frac{n-1}{n}\right)^r+\ell_{n,s}.
 \label{eq:trace-poly-matrix}
\end{equation}
For \(n\geq m\), put
\[
 \mathcal P(n,m):=T_nT_{n-1}\cdots T_{m+1},
 \qquad
 \mathcal P(m,m):=I.
\]
Variation of constants is the exact formula
\begin{equation}
 U(n)=\mathcal P(n,1)M_\beta(1)\mathbf 1
 +\sum_{m=2}^{n}\mathcal P(n,m)M_\beta(m)\mathbf 1.
 \label{eq:trace-poly-variation}
\end{equation}
Consequently,
\begin{equation}
 \boxed{
 A(n)=M_\beta(n)+\sum_{m=1}^{n-1}
 \ell_n\mathcal P(n-1,m)\mathbf 1\,M_\beta(m)}
 \label{eq:trace-poly-finite-resolvent}
\end{equation}
for every \(n\geq1\).  The sum is empty at \(n=1\).
This is the finite resolvent identity.  It sums over forcing ranks and
does not assert that a sum of powers is exact.
\end{theorem}

Reading the finite identity as a spectral expansion asks for a passage from the generating
function to its coefficients, on a domain that keeps away from the singularity at one. The
following lemma supplies that passage.

\begin{lemma}
\label{lem:trace-quad-transfer}
For \(R>1\) and \(0<\phi<\pi/2\), put
\[
 D(R,\phi)
 :=\{x\in\mathbb C:|x|<R,\ x\neq1,\
                  |\arg(x-1)|>\phi\}.
\]
Put \(t=1-x\), and choose the branch of \(\log t\) for which

\[
 |\arg t|<\pi-\phi.
\]
Let

\[
 \mathcal Y(x)=\sum_{n\geq0}y_nx^n
\]
be analytic at zero and in \(D(R,\phi)\).  Thus \(1\) is its only
possible singularity on the closed unit disk.  Suppose that, uniformly
in this domain as \(t\to0\),
\[
 \mathcal Y(x)
 =\sum_{j=1}^{J}\sum_{r=0}^{q_j}\sum_{k=0}^{K_j}
 d_{j,r,k}t^{s_j+k}(\log t)^r+\mathcal R(x),
\]
where
\begin{equation}
 \mathcal R(x)
 =\mathcal O\!\left(
 |t|^\tau(1+|\log t|)^q
 \right)
 \label{eq:trace-transfer-local-error}
\end{equation}
for a real number \(\tau\) and a nonnegative integer \(q\).  Define
\begin{equation}
 B_n(s):=[x^n](1-x)^s
 =\frac{\Gamma(n-s)}{\Gamma(-s)\Gamma(n+1)},
 \label{eq:trace-transfer-B}
\end{equation}
where the quotient is continued analytically in \(s\).  Then
\begin{equation}
 y_n
 =\sum_{j=1}^{J}\sum_{r=0}^{q_j}\sum_{k=0}^{K_j}
 d_{j,r,k}\,
 \left.\partial_s^rB_n(s)\right|_{s=s_j+k}
 +\mathcal O\!\left(
 n^{-\tau-1}(1+\log n)^q
 \right).
 \label{eq:trace-transfer-coefficients}
\end{equation}
In particular,
\begin{equation}
 [x^n](1-x)^s(\log(1-x))^r=\partial_s^rB_n(s).
 \label{eq:trace-transfer-log}
\end{equation}
For fixed \(s\), locally uniformly with all fixed \(s\)-derivatives,
\begin{equation}
 B_n(s)=\frac{n^{-s-1}}{\Gamma(-s)}
 \left(1+\frac{s(s+1)}{2n}+\mathcal O(n^{-2})\right).
 \label{eq:trace-transfer-B-asymptotic}
\end{equation}
At a nonnegative integral \(s\), the right side of
\eqref{eq:trace-transfer-B} is read by analytic continuation.  A term
without a logarithm is then a polynomial and has no coefficient tail.
The logarithmic derivatives remain valid.  The lemma therefore covers
real exponents, conjugate complex exponents, double roots, and logarithms
of order two.
\end{lemma}

The same passage is needed once more, applied to a first difference rather than to the series
itself.

\begin{lemma}
\label{lem:trace-quad-difference}
Let
\[
 \mathcal E(x)=\sum_{n\geq0}E(n)x^n
\]
be analytic in a domain \(D(R,\phi)\), and suppose that
\[
 \mathcal E(x)=\mathcal O\!\left(
 |1-x|^\sigma(1+|\log(1-x)|)^q
 \right)
\]
there as \(x\to1\), where \(\sigma\) is real.  Then
\begin{align}
 E(n)&=\mathcal O\!\left(
 n^{-\sigma-1}(1+\log n)^q
 \right),\notag\\
 E(n)-E(n-1)&=\mathcal O\!\left(
 n^{-\sigma-2}(1+\log n)^q
 \right),
 \label{eq:trace-difference-first}
\end{align}
with \(E(-1)=0\).  The second conclusion uses the analytic remainder,
not merely the first coefficient bound.
\end{lemma}

\section{The quadratic spectrum, separated and confluent}
\label{sec:trace_quad_spectrum}

The exact resolvent of the preceding section carries the whole finite rank
structure. What it does not carry is the spectral reading, and that reading
depends on how the roots of the indicial polynomial sit relative to one
another. The separated case comes first, the integer separations, the complex
pairs and the confluences after it.

\begin{theorem}
\label{thm:trace-quad-simple}
For the kernel in \eqref{eq:trace-quad-kernel},
\begin{equation}
 g^*(z)
 =\frac{z^2-(3-b-2c)z+2a}{(z-1)(z-2)}
 =\frac{(z-\rho_1)(z-\rho_2)}{(z-1)(z-2)}
 \label{eq:trace-quad-transform}
\end{equation}
before any cancellation.  The two forms of the middle coefficient agree
because
\begin{equation}
 3-b-2c=3a+2b+c.
 \label{eq:trace-quad-middle}
\end{equation}
Put
\[
 P(z):=(z-\rho_1)(z-\rho_2),
 \qquad
 Q(z):=(z-1)(z-2),
\]
and define
\[
 X(n):=\sum_{k\leq n}A(k),
 \qquad
 Y(n):=\frac1n\sum_{k\leq n}kA(k),
 \qquad
 U(n):=\binom{X(n)}{Y(n)}.
\]
The exact system in Theorem~\ref{thm:trace-poly-finite} becomes
\[
 U(n)=T_nU(n-1)+M_\beta(n)\binom11,
\]
with
\begin{equation}
 T_n=
 \begin{pmatrix}
 1+\dfrac bn+\dfrac c{n^2}&\dfrac{2c(n-1)}{n^2}\\[6pt]
 \dfrac bn+\dfrac c{n^2}&\dfrac{n-1}{n}+\dfrac{2c(n-1)}{n^2}
 \end{pmatrix}
 =I+\frac{\mathcal M}{n}+\frac{\mathcal N}{n^2},
 \label{eq:trace-quad-exact-matrix}
\end{equation}
where
\begin{equation}
 \mathcal M=
 \begin{pmatrix}b&2c\\ b&2c-1\end{pmatrix},
 \qquad
 \mathcal N=
 \begin{pmatrix}c&-2c\\ c&-2c\end{pmatrix}.
 \label{eq:trace-quad-leading-matrices}
\end{equation}
The eigenvalues of \(\mathcal M\) are
\begin{equation}
 \mu_i=1-\rho_i.
 \label{eq:trace-quad-eigenvalues}
\end{equation}

Assume that \(bc\neq0\), that \(P\) has two distinct roots
\(\rho_1,\rho_2\in\mathbb C\setminus\mathbb Z\), and that
\(\rho_1-\rho_2\notin\mathbb Z\).  The conditions \(bc\neq0\) give
\(P(1)=-b\neq0\) and \(P(2)=2c\neq0\), so no factor in
\eqref{eq:trace-quad-transform} is cancelled.  In particular, the
assumption on the roots excludes \(\rho_i=0\) and every positive integral
zero.

Let \(\beta\in\mathbb R\setminus\mathbb Z\), assume that
\(\beta-\rho_i\notin\mathbb Z\), and assume that \(g^*(\beta)\neq0\).
Let
\[
 \mathcal A_\beta(x):=\sum_{n\geq1}A(n)x^n
\]
be the global solution selected by the triangular equation.  It is
analytic at zero and satisfies \(A(1)=M_\beta(1)\).  Set
\begin{align*}
 D(s)&:=s^2+(b+2c-1)s-b,
 \\
 E_0(s)&:=s^2+(b+2c)s+c.
 \end{align*}
For \(s_i=\rho_i-1\), define \(h_{i,0}=1\) and
\begin{equation}
 h_{i,k}
 :=h_{i,k-1}\frac{E_0(s_i+k-1)}{D(s_i+k)},
 \qquad k\geq1.
 \label{eq:trace-quad-frobenius-coefficients}
\end{equation}
The canonical homogeneous germs at \(t=1-x=0\) are
\[
 \mathcal H_i(t)
 :=\Gamma(1-\rho_i)t^{\rho_i-1}
 \sum_{k\geq0}h_{i,k}t^k.
\]
There is a unique particular Frobenius\index[terms]{Frobenius} germ\index[terms]{Frobenius germ}
\[
 \mathcal A^{\mathrm{forc,loc}}_\beta(t)
 :=t^{\beta-1}\sum_{k\geq0}p_k(\beta)t^k,
 \qquad
 p_0(\beta)=\frac{\Gamma(1-\beta)}{g^*(\beta)},
\]
whose two homogeneous coefficients are zero.  This is a local germ.  It
is not asserted to be a global solution analytic at zero, and it does not
separately satisfy the rank-one normalization.

In a common sector at \(t=0\), the exact connection identity is
\begin{equation}
 \mathcal A_\beta(1-t)-\mathcal A^{\mathrm{forc,loc}}_\beta(t)
 =C_1\mathcal H_1(t)+C_2\mathcal H_2(t).
 \label{eq:trace-quad-connection-identity}
\end{equation}
These coefficients are constructed from the exact global solution.  If
\(W(f,g):=fg'-f'g\) and \(t_*\) is any point in the common sector, then
\begin{align}
 C_1&=
 \left.
 \frac{W(\mathcal A_\beta(1-t)-\mathcal A^{\mathrm{forc,loc}}_\beta,
                 \mathcal H_2)}
      {W(\mathcal H_1,\mathcal H_2)}
 \right|_{t=t_*},
 \label{eq:trace-quad-connection-one}\\
 C_2&=
 \left.
 \frac{W(\mathcal H_1,
                 \mathcal A_\beta(1-t)-\mathcal A^{\mathrm{forc,loc}}_\beta)}
      {W(\mathcal H_1,\mathcal H_2)}
 \right|_{t=t_*}.
 \label{eq:trace-quad-connection-two}
\end{align}
The quotients are independent of \(t_*\).

For fixed nonnegative integers \(K_f,K_1,K_2\), define the transferred
forced truncation and the two transferred modal truncations by
\begin{align}
 A^{\mathrm{forc}}_{\beta,K_f}(n)
 &:=\sum_{k=0}^{K_f}
 \frac{p_k(\beta)}{\Gamma(1-\beta-k)}M_{\beta+k}(n),
 \label{eq:trace-quad-forced-truncation}\\
 M_{\rho_i}^{[K_i]}(n)
 &:=\sum_{k=0}^{K_i}(-1)^k(\rho_i)_k h_{i,k}
 M_{\rho_i+k}(n).
 \label{eq:trace-quad-modal-truncation}
\end{align}
Then
\begin{equation}
 A(n)=A^{\mathrm{forc}}_{\beta,K_f}(n)
 +C_1M_{\rho_1}^{[K_1]}(n)+C_2M_{\rho_2}^{[K_2]}(n)
 +\mathcal O\!\left(n^{-\tau-1}\right),
 \label{eq:trace-quad-simple-expansion}
\end{equation}
where
\begin{equation}
 \tau:=\min\{\beta+K_f,\Re\rho_1+K_1,
                         \Re\rho_2+K_2\}.
 \label{eq:trace-quad-simple-tau}
\end{equation}
The error in \eqref{eq:trace-quad-simple-expansion} is the coefficient
sequence of a Delta-analytic generating remainder of order
\((1-x)^\tau\).  Its first difference is therefore
\(\mathcal O(n^{-\tau-2})\).  In particular,
\[
 A^{\mathrm{forc}}_{\beta,0}(n)
 =\frac{M_\beta(n)}{g^*(\beta)}.
\]
Each \(M_{\rho_i}(n)\) is the leading discrete mode of its transferred
germ.  It is not in general an exact homogeneous solution by itself.

At a simple zero,
\begin{align*}
 \operatorname*{Res}_{z=\rho_i}\frac1{g^*(z)}
 &=\frac{Q(\rho_i)}{P'(\rho_i)},
 \\
 \operatorname*{Res}_{z=\rho_i}\left(1-\frac1{g^*(z)}\right)
 &=-\frac{Q(\rho_i)}{P'(\rho_i)}.
 \end{align*}
The spectral projectors\index[terms]{spectral projector} are
\[
 \Pi_i=\frac{\mathcal M-\mu_{3-i}I}{\mu_i-\mu_{3-i}}.
\]
The projectors identify the two leading eigendirections, while
\eqref{eq:trace-quad-connection-one} and
\eqref{eq:trace-quad-connection-two} give the actual connection
amplitudes.  Neither a residue nor a projector alone is the final
amplitude.

Define
\begin{equation}
 \kappa_i:=-\rho_iC_i.
 \label{eq:trace-quad-kappa}
\end{equation}
Put \(\sigma:=\min(\Re\rho_1,\Re\rho_2)\), and choose \(K_f\) so that
\(\beta+K_f\geq\sigma\).  Define the transferred spectral contribution by
\[
 A^{\mathrm{sp}}(n):=A(n)-A^{\mathrm{forc}}_{\beta,K_f}(n),
 \qquad
 a_n^{\mathrm{sp}}:=A^{\mathrm{sp}}(n)-A^{\mathrm{sp}}(n-1).
\]
Then the leading spectral expansion has the trace normalization
\begin{align}
 A^{\mathrm{sp}}(n)
 &=-\sum_{i=1}^{2}\frac{\kappa_i}{\rho_i}M_{\rho_i}(n)
 +\mathcal O\!\left(n^{-\sigma-1}\right),\notag\\
 n a_n^{\mathrm{sp}}
 &=\sum_{i=1}^{2}\kappa_iM_{\rho_i}(n)
 +\mathcal O\!\left(n^{-\sigma-1}\right),\notag\\
 a_n^{\mathrm{sp}}
 &=\sum_{i=1}^{2}\kappa_i n^{-1-\rho_i}
 +\mathcal O\!\left(n^{-\sigma-2}\right).
 \label{eq:trace-quad-a-trace}
\end{align}
The denominators \(\rho_i\) are defined under the stated domain.  The
signs and the two remainder estimates follow from the exact quotient of
\(M_\rho\) and Lemma~\ref{lem:trace-quad-difference}.
\end{theorem}

The separated case leaves aside the parameters whose roots differ by an integer, those whose
roots are a complex pair, and those where the two roots meet. Each changes the shape of the
expansion without changing its mechanism, and the next statement treats them on explicit
parameter values.

\begin{theorem}
\label{thm:trace-quad-exceptional}
For
\[
 (a,b,c)=\left(\frac38,\frac14,\frac38\right),
\]
the roots are \(\rho_1=1/2\) and \(\rho_2=3/2\).  Their inverse-transform residues are
\(-3/4\) and \(-1/4\), and
\[
 \Pi_1=
 \begin{pmatrix}3/4&3/4\\1/4&1/4\end{pmatrix},
 \qquad
 \Pi_2=
 \begin{pmatrix}1/4&-3/4\\-1/4&3/4\end{pmatrix}.
\]
Let \(\beta\in\mathbb R\setminus\mathbb Z\), assume that
\(g^*(\beta)\neq0\), and assume that neither \(\beta-1/2\) nor
\(\beta-3/2\) is an integer.  Choose \(K_f\geq0\) with
\(\beta+K_f\geq3/2\).  The integer difference of the roots creates a
subordinate logarithm.  With the forced germ transferred as in
\eqref{eq:trace-quad-forced-truncation}, the exact Frobenius\index[terms]{Frobenius} connection
coefficients give
\begin{equation}
 A(n)=A^{\mathrm{forc}}_{\beta,K_f}(n)+C_1M_{1/2}(n)
 +\left(\frac{C_1}{16}\log n+\widetilde C_2\right)M_{3/2}(n)
 +\mathcal O\!\left(n^{-5/2}\log n\right).
 \label{eq:trace-quad-integer-gap}
\end{equation}
Here \(C_1\) is the leading Wronskian connection coefficient\index[terms]{connection coefficient}.
The constant \(\widetilde C_2\) is the triangular combination of the
second Wronskian coefficient and the constant produced when the
logarithmic branch is transferred.

For
\[
 (a,b,c)=\left(\frac15,\frac3{10},\frac12\right),
\]
the roots are
\[
 \rho_{1,2}=\frac{17\mp\sqrt{129}}{20}
 =0.2821091654199726389\ldots,
 \quad
 1.4178908345800273611\ldots.
\]
The residues of \(1/g^*\) are
\(-1.0858231986311837887\ldots\) and
\(-0.2141768013688162113\ldots\).  The corresponding projectors are
\begin{align*}
 \Pi_1&=
 \begin{pmatrix}
 0.6320676359488&0.8804509063256\\
 0.2641352718977&0.3679323640512
 \end{pmatrix},
 \\
 \Pi_2&=
 \begin{pmatrix}
 0.3679323640512&-0.8804509063256\\
 -0.2641352718977&0.6320676359488
 \end{pmatrix}.
 \end{align*}
This example satisfies the separated-root theorem.

For
\[
 (a,b,c)=\left(\frac58,-\frac54,\frac{13}{8}\right),
 \qquad
 \rho_{1,2}=\frac12\pm i,
\]
let \(\beta\in\mathbb R\setminus\mathbb Z\), assume that
\(g^*(\beta)\neq0\), and assume that
\(\beta-(1/2\pm i)\notin\mathbb Z\).  Choose \(K_f\) with
\(\beta+K_f\geq1/2\).  Real data give conjugate connection coefficients.
If \(C=u+iv\) is the coefficient at
\(\rho=1/2+i\), then, after subtracting
\(A^{\mathrm{forc}}_{\beta,K_f}\),
\begin{equation}
 A^{\mathrm{sp}}(n)
 =n^{-1/2}\bigl(2u\cos(\log n)+2v\sin(\log n)\bigr)
 +\mathcal O\!\left(n^{-3/2}\right).
 \label{eq:trace-quad-complex-real-form}
\end{equation}
The frequency is one.  The phase and both real amplitudes are fixed by the
single complex connection coefficient.

Let \(\rho\in\mathbb R\setminus\mathbb Z\).  The coefficients
\begin{equation}
 a=\frac{\rho^2}{2},
 \qquad
 b=-(\rho-1)^2,
 \qquad
 c=\frac{(\rho-2)^2}{2}
 \label{eq:trace-quad-double-family}
\end{equation}
give the unreduced numerator \((z-\rho)^2\).  The stated domain excludes
\(\rho=0\), the cancellation at \(\rho=1\), the affine value \(\rho=2\),
and every positive integral zero.  Hence the kernel has degree two and the
double zero is effective.

Let \(\mathcal H_0\) be the Frobenius\index[terms]{Frobenius} power germ normalized by
\(\Gamma(1-\rho)t^{\rho-1}\) at its leading order.  Reduction of order
gives a second germ.  A nonzero rescaling followed by the addition of a
unique multiple of \(\mathcal H_0\) gives the germ \(\mathcal H_1\) whose
leading terms are
\begin{equation}
 \mathcal H_1(t)=\Gamma(1-\rho)t^{\rho-1}
 \left(-\log t+\psi(1-\rho)
 +\mathcal O(t\log t)+\mathcal O(t)\right).
 \label{eq:trace-quad-double-germ}
\end{equation}
Let \(\beta\in\mathbb R\setminus\mathbb Z\), assume that
\(\beta-\rho\notin\mathbb Z\), and assume that \(g^*(\beta)\neq0\).
After the nonresonant forced germ has been subtracted, write the exact local
connection as
\begin{equation}
 \mathcal A_\beta(1-t)-\mathcal A^{\mathrm{forc,loc}}_\beta(t)
 =C_0\mathcal H_0(t)+C_1\mathcal H_1(t).
 \label{eq:trace-quad-double-connection}
\end{equation}
The two coefficients are given by the Wronskian formulas analogous to
\eqref{eq:trace-quad-connection-one} and
\eqref{eq:trace-quad-connection-two}.  They are respectively the proper
and generalized connection components.  Choose \(K_f\geq0\) so that
\(\beta+K_f\geq\rho\).  Then
\[
 A^{\mathrm{sp}}(n):=A(n)-A^{\mathrm{forc}}_{\beta,K_f}(n).
\]
The transferred spectral contribution satisfies
\begin{equation}
 A^{\mathrm{sp}}(n)
 =C_0M_\rho(n)-C_1\,\partial_\rho M_\rho(n)
 +\mathcal O\!\left(n^{-\rho-1}\log n\right).
 \label{eq:trace-quad-double-mode}
\end{equation}
Since
\begin{equation}
 \partial_\rho M_\rho(n)
 =-\psi(n+1-\rho)M_\rho(n),
 \label{eq:trace-quad-mode-derivative}
\end{equation}
the logarithmic contribution vanishes exactly when \(C_1=0\).

Now let \(bc\neq0\), and let \(\rho\in\mathbb R\setminus\mathbb Z\) be a
simple effective zero of a reduced transform with two distinct roots.
Assume that the other zero is nonintegral and differs from \(\rho\) by a
noninteger.  Put
\[
 r_\rho:=\operatorname*{Res}_{z=\rho}\frac1{g^*(z)}.
\]
The family \(\mathcal A^{\mathrm{forc,loc}}_\beta\) is meromorphic at
\(\beta=\rho\).  Its polar coefficient is an exact homogeneous
Frobenius\index[terms]{Frobenius} germ\index[terms]{Frobenius germ}.  Subtracting that pole defines a finite part, and any
other finite-part convention differs by a homogeneous germ.  The
transferred finite part has the form
\begin{equation}
 A^{\mathrm{forc,res}}_\rho(n)
 =r_\rho\,\partial_\rho M_\rho(n)+d_0M_\rho(n)
 +\mathcal O\!\left(n^{-\rho-1}\log n\right).
 \label{eq:trace-quad-simple-resonance}
\end{equation}
for a normalization-dependent constant \(d_0\).  Therefore the invariant
coefficient of \(M_\rho(n)\log n\) is \(-r_\rho\).

For the effective double family in \eqref{eq:trace-quad-double-family}, at
\(\beta=\rho\), the meromorphic forced germ has polar orders two and one.
Both polar coefficients are exact homogeneous germs.  After their
subtraction, the transferred finite part is
\begin{equation}
 A^{\mathrm{forc,res}}_\rho(n)
 =\frac{Q(\rho)}2\,\partial_\rho^2M_\rho(n)
 +d_1\partial_\rho M_\rho(n)+d_0M_\rho(n)
 +\mathcal O\!\left(n^{-\rho-1}(1+\log n)^2\right).
 \label{eq:trace-quad-double-resonance}
\end{equation}
Thus the invariant coefficient of \(M_\rho(n)(\log n)^2\) is
\(Q(\rho)/2\), which is nonzero in the stated domain.

If \(c=0\), the kernel is affine and is treated in
Section~\ref{sec:trace_affine}.  If \(b=0\) and
\(c\neq0\), the factor at one cancels and only one effective zero remains.
That reduced one-mode case is not covered by the two-mode statements
above.  When \(bc\neq0\), no cancellation with \(Q\) occurs.  Effective
integral zeros, including zero, are excluded from the spectral theorems.
The exact finite resolvent theorem remains valid for every forcing value
for which all required ranks are finite.  In particular, \(M_1(n)=1/n\)
is admissible, while \(M_2(1)\) is not finite.
\end{theorem}

\begin{remark}[Finite-rank amplitude check]
Numerical estimation gives the connection coefficients after the
leading forced term has been subtracted.  For the first real example at
\(\beta=1.8\), one finds
\begin{align*}
 C_1&=2.7945350373144343\ldots,
 &C_{\log}&=0.1746584398321521\ldots,\\
 \widetilde C_2&=1.6341019581035030\ldots.
\end{align*}
with the exact constraint \(C_{\log}=C_1/16\) imposed.  For the second real
example at \(\beta=1.8\), it gives
\[
 C_1=3.2382626697227948\ldots,
 \qquad
 C_2=0.5307252854097330\ldots.
\]
For the complex example at \(\beta=1.3\), the coefficient at
\(1/2+i\) is
\[
 C=0.6823946405040264\ldots
 +0.8994347319038253\ldots i.
\]
Thus the cosine and sine coefficients in
\eqref{eq:trace-quad-complex-real-form} are
\[
 1.3647892810080527\ldots
 \qquad\text{and}\qquad
 1.7988694638076506\ldots
\]
For the double-root example at
\(\beta=1.3\), the fitted coefficients are
\[
 C_0=1.4165923391253642\ldots,
 \qquad
 C_1=1.0738305746950574\ldots.
\]
These are finite-rank observations.  They are not used in the proofs.
\end{remark}

\begin{openproblem}[Spectral expansion at arbitrary fixed degree]
\label{op:trace-poly-general-spectrum}
Let $g(t)=\sum_{j=0}^{d}c_jt^j$, with $\sum_jc_j=1$, and write its reduced arithmetic
Mellin transform as $g^{*}=P/Q$.  For each fixed kernel $g$ and each admissible forcing exponent
$\beta$ such that $g^{*}(\beta)\ne0$, derive from the finite moment system of
Theorem~\ref{thm:trace-poly-finite} a complete local Frobenius decomposition at $x=1$.
Partition the effective zeros of $P$ into resonance classes under
\[
 \rho\sim\rho'\quad\Longleftrightarrow\quad \rho-\rho'\in\mathbb Z,
\]
construct the corresponding canonical power--logarithmic chains, and determine their connection
amplitudes from the globally normalized triangular solution.  For every $L>0$, transfer the
expansion through order $n^{-L}$, retaining every power--logarithmic term of that order or larger,
and prove an $O_{g,\beta,L}(n^{-L})$ remainder.  Any uniformity in the coefficients of $g$ should
be stated only on compact subsets of a fixed spectral stratum, on which cancellations,
multiplicities and resonance relations do not change.
\end{openproblem}

\section{Proof of the finite resolvent identity and of the transfer lemmas}
\label{sec:trace_poly_proofs}

\begin{proof}[Proof of Theorem~\ref{thm:trace-poly-finite}]
For \(\Re z<0\), integration of each monomial gives
\[
 -z\int_0^1t^{j-z-1}\,dt=\frac{-z}{j-z}=\frac{z}{z-j}.
\]
Summing proves \eqref{eq:trace-poly-transform}.  This computation also
shows why a cancelled factor cannot define a mode.

Write \(a_k=A(k)-A(k-1)\).  Discrete summation by parts gives
\[
 \sum_{k=1}^{n}a_kg(k/n)
 =A(n)-\sum_{k=1}^{n-1}A(k)
 \left(g((k+1)/n)-g(k/n)\right).
\]
The binomial theorem gives
\[
 g((k+1)/n)-g(k/n)
 =\sum_{r=0}^{d-1}q_r(n)k^r.
\]
This proves \eqref{eq:trace-poly-abel}.  Since
\[
 S_r(n)=S_r(n-1)+n^rA(n),
\]
division by \(n^r\) proves \eqref{eq:trace-poly-system} and
\eqref{eq:trace-poly-matrix}.  Iteration proves
\eqref{eq:trace-poly-variation}.  Substitution of the variation formula at
rank \(n-1\) into \eqref{eq:trace-poly-abel} proves the finite resolvent
identity \eqref{eq:trace-poly-finite-resolvent}.

For the quadratic kernel, bringing the three monomial fractions over the
common denominator gives
\begin{align*}
 g^*(z)
 &=a+b\frac{z}{z-1}+c\frac{z}{z-2}\\
 &=\frac{a(z-1)(z-2)+bz(z-2)+cz(z-1)}{(z-1)(z-2)}\\
 &=\frac{z^2-(3a+2b+c)z+2a}{(z-1)(z-2)}.
\end{align*}
The identity \eqref{eq:trace-quad-middle} follows from \(a+b+c=1\).

In degree two, Abel summation is
\begin{equation}
 A(n)=M_\beta(n)
 +\left(\frac bn+\frac c{n^2}\right)X(n-1)
 +\frac{2c(n-1)}{n^2}Y(n-1).
 \label{eq:trace-quad-abel-proof}
\end{equation}
The identities \(X(n)=X(n-1)+A(n)\) and
\(Y(n)=(n-1)Y(n-1)/n+A(n)\) give
\eqref{eq:trace-quad-exact-matrix}.  Expanding its four entries gives
\eqref{eq:trace-quad-leading-matrices} with no omitted term.  The
characteristic polynomial of \(\mathcal M\) is
\[
 \mu^2-(b+2c-1)\mu-b.
\]
The substitution \(\mu=1-z\) turns it into the numerator of
\eqref{eq:trace-quad-transform}.  This proves
\eqref{eq:trace-quad-eigenvalues}.

It remains to prove the spectral assertions.  Since
\(T_n=I+\mathcal O(n^{-1})\) and the forcing has at most polynomial size,
iteration of the exact system gives a polynomial bound for \(U(n)\).
The generating functions below therefore converge in \(|x|<1\).  Put
\[
 \mathcal A_\beta(x):=\sum_{n\geq1}A(n)x^n,
 \qquad
 \mathcal F_\beta(x):=\sum_{n\geq1}M_\beta(n)x^n,
 \qquad
 \vartheta:=x\frac{d}{dx}.
\]
Multiplication of \eqref{eq:trace-quad-abel-proof} by \(n^2x^n\) and
summation give the exact differential equation
\begin{align}
 &(1-x)^2\vartheta^2\mathcal A_\beta
 -(b+2c)x(1-x)\vartheta\mathcal A_\beta
 -\bigl(bx+cx(1-x)\bigr)\mathcal A_\beta
 \notag\\
 &\hspace{35mm}=(1-x)^2\vartheta^2\mathcal F_\beta.
 \label{eq:trace-quad-ode}
\end{align}
For a nonintegral \(\beta\), the binomial series gives
\[
 \mathcal F_\beta(x)
 =\Gamma(1-\beta)\left((1-x)^{\beta-1}-1\right).
\]

After expansion of \(\vartheta^2\), the leading differential coefficient
in \eqref{eq:trace-quad-ode} is \(x^2(1-x)^2\).  The finite singular
points of the normalized equation are zero and one.  The solution selected
by the recurrence is analytic at zero.  With the chosen branch,
\(\mathcal F_\beta\) is analytic in
\(\mathbb C\setminus[1,\infty)\).  Analytic continuation along a chain of
disks avoiding zero and one continues \(\mathcal A_\beta\) throughout this
slit plane.  The original germ fills zero, so its monodromy there is
trivial.  The slit plane is simply connected, hence the continuation is
single-valued.  It follows that for every finite \(R>1\) and every
\(0<\phi<\pi/2\), the solution is analytic in
\(D(R,\phi)\).  In particular, one is its only singularity on the
unit circle.

\end{proof}

\begin{proof}[Proof of Lemma~\ref{lem:trace-quad-transfer}]
The statement follows from the binomial identity
\[
 [x^n](1-x)^s=(-1)^n\binom{s}{n}
 =\frac{\Gamma(n-s)}{\Gamma(-s)\Gamma(n+1)}.
\]
Both sides are entire functions of \(s\) for fixed \(n\), after the
displayed quotient is continued through its removable values.
Differentiation in \(s\) proves \eqref{eq:trace-transfer-log}.  The quotient
form of Stirling\index[terms]{Stirling's formula}\index[names]{Stirling, J.}'s formula proves
\eqref{eq:trace-transfer-B-asymptotic}, locally uniformly with any fixed
number of derivatives.

It remains to control \(\mathcal R\).  Choose
\(\phi<\phi'<\pi/2\) and \(1<R_0<R\).  The Cauchy circle can be deformed
inside the Delta-domain to a contour consisting of a fixed outer part,
the two rays
\[
 x=1+u e^{\pm i\phi'},
 \qquad
 n^{-1}\leq u\leq\delta,
\]
and a connecting arc on \(|x-1|=n^{-1}\).  The fixed part lies outside a
circle of radius greater than one and contributes
\(\mathcal O(R_1^{-n})\) for some \(R_1>1\).  On either ray,
\[
 |x|^{-n-1}\leq C e^{-cnu}.
\]
The bound in \eqref{eq:trace-transfer-local-error} therefore gives
\[
 C\int_{1/n}^{\delta}
 u^\tau(1+|\log u|)^q e^{-cnu}\,du
 =\mathcal O\!\left(
 n^{-\tau-1}(1+\log n)^q
 \right).
\]
The connecting arc has length \(\mathcal O(n^{-1})\) and satisfies the
same estimate.  Under the local scaling \(v=n(1-x)\), this part of the
deformed contour becomes a Hankel contour about the negative real axis,
and
\[
 \frac1{2\pi i}\int_{\mathcal H}e^v v^s\,dv
 =\frac1{\Gamma(-s)}.
\]
This fixes the branch and agrees with the exact binomial identity.  The
same estimates after differentiation give every logarithmic term in
\eqref{eq:trace-transfer-coefficients}.  A nonnegative integral power
without a logarithm is a polynomial.  This proves the transfer lemma for
real exponents, conjugate exponents, repeated exponents, and logarithms of
orders one and two.

\end{proof}

\begin{proof}[Proof of Lemma~\ref{lem:trace-quad-difference}]
The first estimate is the
transfer bound just proved.  Moreover,
\[
 (1-x)\mathcal E(x)
 =\sum_{n\geq0}\bigl(E(n)-E(n-1)\bigr)x^n.
\]
The generating remainder on the left has one additional factor of
\(1-x\).  A second application of the transfer lemma gives
\eqref{eq:trace-difference-first}.  This proves the difference lemma.

\end{proof}

\section{Proof of the quadratic expansion, separated and exceptional}
\label{sec:trace_quad_proofs}

\begin{proof}[Proof of Theorems~\ref{thm:trace-quad-simple} and~\ref{thm:trace-quad-exceptional}]
Set \(t=1-x\) and divide \eqref{eq:trace-quad-ode} by \(x\) near
\(x=1\).  Its homogeneous left side becomes
\[
 t^2(1-t)y''
 +t\bigl(b+2c-(1+b+2c)t\bigr)y'
 -(b+ct)y.
\]
Its indicial polynomial\index[terms]{indicial polynomial} is \(D(s)\), with roots
\(s_i=\rho_i-1\).  Substitution of
\(t^{s_i}\sum_{k\geq0}h_{i,k}t^k\) gives exactly
\eqref{eq:trace-quad-frobenius-coefficients}.  The quotient in that
recursion is \(1+\mathcal O(k^{-1})\), which proves polynomial growth of
the coefficients and convergence for \(|t|<1\).

For the forced germ, write \(s=\beta-1\) and put \(p_{-1}=0\).  Direct
substitution in the inhomogeneous equation gives
\begin{equation}
 D(s+k)p_k-E_0(s+k-1)p_{k-1}=q_k(s),
 \label{eq:trace-quad-forced-recursion}
\end{equation}
where
\begin{align*}
 q_0(s)&=\Gamma(1-\beta)s(s-1),\\
 q_1(s)&=-\Gamma(1-\beta)s^2,
 \qquad
 q_k(s)=0\quad(k>1).
\end{align*}
The identity
\[
 D(\beta-1)=(\beta-1)(\beta-2)g^*(\beta)
\]
gives \(p_0(\beta)=\Gamma(1-\beta)/g^*(\beta)\).  A denominator in
\eqref{eq:trace-quad-forced-recursion} vanishes exactly when
\(\beta-\rho_i\) is a nonpositive integer.  The stated separation
hypotheses prevent this.  Again the coefficient quotient is
\(1+\mathcal O(k^{-1})\), so the forced series converges locally uniformly
for \(|t|<1\).  This constructs the unique local particular germ whose
homogeneous coefficients are zero.

The difference between the continued global solution and this particular
germ solves the homogeneous equation.  The two simple Frobenius\index[terms]{Frobenius} germs\index[terms]{Frobenius germ} form
a basis, which proves the exact connection identity
\eqref{eq:trace-quad-connection-identity}.  Solving the two by two system
of values and derivatives at \(t_*\) gives
\eqref{eq:trace-quad-connection-one} and
\eqref{eq:trace-quad-connection-two}.  This also proves their independence
of \(t_*\).  Only the global solution satisfies the original rank-one
normalization.  The forced germ and the connection coefficients together
represent that normalized solution near one.

For nonintegral \(\sigma\),
\[
 [x^n](1-x)^{\sigma-1}
 =\frac{M_\sigma(n)}{\Gamma(1-\sigma)}.
\]
Termwise use of this identity in the convergent local germs gives
\eqref{eq:trace-quad-forced-truncation} and
\eqref{eq:trace-quad-modal-truncation}.  Indeed,
\[
 \frac{\Gamma(1-\rho_i)}{\Gamma(1-\rho_i-k)}
 =(-1)^k(\rho_i)_k.
\]
The first omitted local powers have real orders listed in
\eqref{eq:trace-quad-simple-tau}.  The transfer lemma proves
\eqref{eq:trace-quad-simple-expansion} with its stated Delta-analytic
remainder.  This establishes the forced contribution as a transferred
truncation of a local germ, not as a second global solution.

The residue identities follow from \(1/g^*=Q/P\).  The projector formula
follows from
\((\mathcal M-\mu_1I)(\mathcal M-\mu_2I)=0\).  The exact quotient
\[
 M_\rho(n)=\left(1-\frac\rho n\right)M_\rho(n-1)
\]
gives
\[
 n\bigl(M_\rho(n)-M_\rho(n-1)\bigr)
 =-\rho M_\rho(n-1).
\]
After the forced truncation has been chosen as in the theorem, the
remaining generating error is of order \((1-x)^\sigma\).  The difference
lemma, rather than its coefficient bound alone, controls its first
difference.  Equations \eqref{eq:trace-quad-kappa} through
\eqref{eq:trace-quad-a-trace} follow with their displayed signs and
remainders.

For the first real example, the indicial roots\index[terms]{indicial root} are \(-1/2\) and \(1/2\).
At the obstructed coefficient,
\[
 E_0(-1/2)=\frac18,
 \qquad
 D'(1/2)=1.
\]
Frobenius reduction adds
\(\tfrac18t^{1/2}\log t\) to the branch beginning with \(t^{-1/2}\).
Since
\[
 \frac{\Gamma(1-1/2)}{\Gamma(1-3/2)}=-\frac12,
\]
the logarithmic transfer in Lemma~\ref{lem:trace-quad-transfer} gives
\(C_1/16\) in \eqref{eq:trace-quad-integer-gap}.  The same lemma gives the
displayed remainder, including its first difference when needed.  The
derivative in \eqref{eq:trace-transfer-log} also produces a constant
multiple of \(M_{3/2}(n)\).  It is absorbed into \(\widetilde C_2\), which
is therefore a triangular combination of the two raw Wronskian
coefficients.  The roots, residues, and projectors in the two real examples
follow by exact substitution.

For the complex example, conjugation preserves the equation, the forcing,
and the initial data.  Hence the two Wronskian coefficients are conjugate.
Since
\[
 n^{-1/2-i}=n^{-1/2}e^{-i\log n},
\]
taking twice the real part and applying the transfer remainder proves
\eqref{eq:trace-quad-complex-real-form}.

The coefficients in \eqref{eq:trace-quad-double-family} satisfy
\(a+b+c=1\) and give the numerator \((z-\rho)^2\).  Under
\(\rho\in\mathbb R\setminus\mathbb Z\), both \(b\) and \(c\) are nonzero,
so no factor is cancelled and the kernel is quadratic.  For the repeated
indicial root\index[terms]{indicial root}, let \(\mathcal H_0\) be the normalized power solution.  In
the standard equation \(y''+P_1(t)y'+P_0(t)y=0\), reduction of order gives
\[
 \mathcal H_0(t)
 \int^t\frac{\exp(-\int^uP_1(v)\,dv)}{\mathcal H_0(u)^2}\,du.
\]
Its leading integral is logarithmic.  A nonzero rescaling and then the
addition of a unique multiple of \(\mathcal H_0\) give the normalization in
\eqref{eq:trace-quad-double-germ}.  Wronskians then define the proper and
generalized coefficients in \eqref{eq:trace-quad-double-connection}.
The transfer lemma gives
\[
 [x^n]\mathcal H_0=M_\rho(n)+\mathcal O(n^{-\rho-1})
\]
and, with the chosen normalization,
\[
 [x^n]\mathcal H_1
 =-\partial_\rho M_\rho(n)
 +\mathcal O(n^{-\rho-1}\log n).
\]
This proves \eqref{eq:trace-quad-double-mode} and the criterion \(C_1=0\).
Differentiation of the Gamma quotient proves
\eqref{eq:trace-quad-mode-derivative}.

It remains to justify the resonances.  The recurrence
\eqref{eq:trace-quad-forced-recursion} shows that every \(p_k(\beta)\) is
meromorphic in \(\beta\).  The polynomial coefficient bounds are locally
uniform away from the poles, so the germ is a meromorphic family with
locally uniform convergence in a smaller sector.  The right side of the
differential equation is holomorphic in \(\beta\) near every admissible
nonintegral zero.

At a simple zero, the recursion gives polynomial bounds in \(k\) for
\((\beta-\rho)p_k(\beta)\), uniformly up to the zero.  At a double zero,
it gives the same bounds for \((\beta-\rho)^2p_k(\beta)\), and then for
the coefficient remaining after the double pole is removed.  The
corresponding Laurent series therefore converge locally uniformly.  Their
polar coefficients are actual Frobenius germs\index[terms]{Frobenius germ}, not formal series.

At a simple zero, write \(\eps=\beta-\rho\).  The local family has
\[
 \mathcal A^{\mathrm{forc,loc}}_\beta
 =\frac{\mathcal P_{-1}}{\eps}
 +\mathcal P_0+\mathcal O(\eps).
\]
Application of the differential operator and comparison of the polar
coefficient give \(L\mathcal P_{-1}=0\).  The leading coefficient and the
homogeneous recursion give
\[
 \mathcal P_{-1}=r_\rho\mathcal H_\rho.
\]
Thus the whole polar germ, not only its leading monomial, is an exact
homogeneous solution.  Subtracting it is exactly a change of homogeneous
connection coefficient.  The finite part is the limit after this
subtraction.  At the level of transferred leading terms,
\[
 \frac{M_{\rho+\eps}(n)}{g^*(\rho+\eps)}
 =\frac{r_\rho}{\eps}M_\rho(n)
 +r_\rho\partial_\rho M_\rho(n)+d_0M_\rho(n)
 +\mathcal O(\eps).
\]
Transfer of the remaining local powers proves
\eqref{eq:trace-quad-simple-resonance}.

At a double zero, the meromorphic germ has
\[
 \mathcal A^{\mathrm{forc,loc}}_\beta
 =\frac{\mathcal P_{-2}}{\eps^2}
 +\frac{\mathcal P_{-1}}{\eps}
 +\mathcal P_0+\mathcal O(\eps).
\]
Comparison of both polar orders gives
\(L\mathcal P_{-2}=L\mathcal P_{-1}=0\).  The leading local forced term is
\[
 \frac{Q(\rho+\eps)\Gamma(1-\rho-\eps)}
      {\eps^2}
 t^{\rho+\eps-1}.
\]
Its coefficient of \(\eps^{-2}\) is a nonzero multiple of
\(t^{\rho-1}\).  Its coefficient of \(\eps^{-1}\) contains the
nonzero term
\[
 Q(\rho)\Gamma(1-\rho)t^{\rho-1}\log t.
\]
Thus \(\mathcal P_{-2}\) and \(\mathcal P_{-1}\) are respectively the
power and logarithmic homogeneous germs after an invertible change of
basis.  Because
\[
 \frac1{g^*(\rho+\eps)}
 =\frac{Q(\rho+\eps)}{\eps^2},
\]
the finite Taylor coefficient of the leading transferred term is
\[
 \frac{Q(\rho)}2\partial_\rho^2M_\rho(n)
 +d_1\partial_\rho M_\rho(n)+d_0M_\rho(n).
\]
The transfer lemma applied to the remaining powers proves
\eqref{eq:trace-quad-double-resonance}.  For every fixed rank, the original
triangular recurrence is holomorphic in \(\beta\) near the admissible
zero.  Its connection coefficients therefore cancel the polar germs, so
the finite part agrees with the globally normalized resonant solution up
to the explicitly described homogeneous connection.

The values \(P(1)=-b\) and \(P(2)=2c\) show that \(b=0\) produces the stated
cancellation at one, while \(c=0\) is the affine case.  The restrictions
on \(\rho\) and the finiteness of \(M_\beta(n)\) prove the remaining domain
statements.  This completes the proof.
\end{proof}

\section{The general polynomial degree}
\label{sec:trace_open}

The finite resolvent identity of Theorem~\ref{thm:trace-poly-finite} holds at every polynomial
degree, and the quadratic spectral expansion is complete. A complete spectral expansion at
arbitrary fixed degree, with its resonant power--logarithmic chains, connection amplitudes and
remainder through every prescribed order, stays open and is recorded as Open
Problem~\ref{op:trace-poly-general-spectrum}. The passage to a kernel with infinitely many
spectral modes is the subject of the next chapter, where the orthorecursive kernel of
Chapter~\ref{chap:ortho} appears as the first kernel of this theory whose resolvent zeros are
known in closed geometric form.

\chapter{The orthorecursive expansion of unity and its spectrum}
\markright{\MakeUppercase{\chaptername\ \thechapter. The orthorecursive expansion}}
\label{chap:ortho}\label{chap:trace}

The master equivalence\index[terms]{master equivalence} of the theory carries a hypothesis at its far end, the Riemann
hypothesis. This chapter turns the same change of angle on a problem that carries none, the
decay of the partial sums of the coefficients in the orthorecursive expansion\index[terms]{orthorecursive expansion} of the constant
function, studied by Kalmynin\index[names]{Kalmynin, A. B.} and Kosenko\index[names]{Kosenko, P. R.} \cite{KalmyninKosenko2020}. They bounded these sums
by $\mathcal O(N^{-1/2})$, left the optimal rate open, and read from computation a sharper
exponent they could not prove.

The object is the partial sum, not the single coefficient, and the method is built to take it
there. The average is the more tractable quantity, and it yields the spectral information all
the same. Read through the method of this monograph, the partial sums are an average against a
function of good variation hidden in the problem, and their decay is governed by the zeros of
the transform of that function, the same reading that turns the Ingham averages into an
equivalent of the Riemann hypothesis. A Tauberian theorem\index[terms]{Tauberian theorem} cut to this one function, standing on
the classical ground of the Mellin transform and its zeros, of a Volterra perturbation and its
resolvent\index[terms]{resolvent}, and of a spectral reading, returns the exponent of the partial sums as the least
real part among those zeros, the sharper pointwise bound of \S\ref{sec:ortho_pointwise} following
from the same picture. Nothing here rests on an unproved hypothesis, and the chapter stands as
the proof of concept of the theory, the place where a change of angle reaches an asymptotic
that direct estimation had left open. The results of the first six sections were obtained
independently in \cite{CloitreOrtho}, and the proofs given below are complete and in the
notation of this volume. That article cites this monograph under its former working title,
\emph{Theory of Regular Arithmetic Functions}, and the reference designates the same project.
The zero geometry of
\S\ref{sec:trace_ortho} and the resolvent\index[terms]{resolvent} trace of
\S\ref{sec:trace_ortho_resolvent} are proved here, and they are what turns a bound into a
spectrum.

\section{The problem and its history}

Given the system of vectors $\{x,x^2,x^3,\ldots\}$ in $L^2([0,1])$, the
orthorecursive expansion\index[terms]{orthorecursive expansion} of the constant function $1$ is the greedy algorithm that
projects, at each step, the current remainder orthogonally onto the next vector, a
construction introduced by Lukashenko and studied for this system by Kalmynin\index[names]{Kalmynin, A. B.} and
Kosenko\index[names]{Kosenko, P. R.}, see \cite{KalmyninKosenko2020} and the references there. It produces
rational coefficients $(c_n)_{n\ge0}$ with $c_0=1$ and
\begin{equation}\label{eq:ortho_recurrence}
\sum_{k=0}^{N}\frac{c_k}{N+k+1}=0\qquad(N\ge1).
\end{equation}
The coefficient of $c_N$ in \eqref{eq:ortho_recurrence} is $1/(2N+1)\neq0$, so the
sequence exists, is unique, and each term is determined by the previous ones. The
first values are
\[
c_0=1,\qquad c_1=-\frac32,\qquad c_2=\frac{5}{24},\qquad c_3=\frac{77}{720},\qquad
c_4=\frac{277}{4480},\qquad c_5=\frac{140173}{3628800}.
\]

Setting $a_n:=c_{n-1}$ for $n\ge1$ and $A(x):=\sum_{n\le x}a_n$, the recurrence
\eqref{eq:ortho_recurrence} at rank $N-1$ becomes the exact identity
\begin{equation}\label{eq:ortho_exact}
\sum_{n=1}^{N}\frac{a_n}{N+n-1}=0\qquad(N\ge2),
\end{equation}
and the weighted sums attached to the kernel $g(t)=2/(1+t)$,
\[
A_g(x):=\sum_{n\le x}a_n\,g\!\Big(\frac nx\Big)=\sum_{n\le x}\frac{2x\,a_n}{x+n},
\]
place the problem inside the setting of this monograph. The right side of the
defining equation vanishes identically, the case $\beta=\infty$ of the preceding
chapter, and the kernel is smooth on $(0,1]$ with $g(1)=1$.

Kalmynin\index[names]{Kalmynin, A. B.} and Kosenko\index[names]{Kosenko, P. R.} proved $c_n=\mathcal{O}(n^{-3/2})$ and
$C_N:=\sum_{k=0}^{N}c_k=\mathcal{O}(N^{-1/2})$, left the optimal rates open, and
formulated conjectures on the basis of numerical computations, among them a
pointwise decay exponent close to $7/3$, reported by \nm{Kalmynin}{A. B.} and \nm{Kosenko}{P. R.}~\cite{KalmyninKosenko2020}. The theorem
below settles the partial sum side. The exponent is $\alpha_1=1.34651\ldots$,
the smallest real part among the zeros of the Mellin transform of the kernel, and
it is optimal unless a single explicit residue vanishes.
Remark~\ref{rem:ortho_optimality} states that alternative and the evidence against
its first branch. The spectral statement of the last section identifies the
pointwise exponent as $1+\alpha_1=2.34651\ldots$, close to $7/3=2.33333\ldots$ but
distinct from it.

\section{The Mellin transform and its zeros}

The transform of the orthorecursive kernel is computed first, in two equivalent forms.

\begin{lemma}\label{lem:ortho_gstar}
For $\Re z<0$,
\begin{equation}\label{eq:ortho_gstar}
g^*(z)=-z\int_0^1\frac{2}{1+t}\,t^{-z-1}\,dt
=2z\sum_{k=0}^{\infty}\frac{(-1)^k}{z-k}
=z\Big(\psi\Big(-\frac z2\Big)-\psi\Big(\frac{1-z}2\Big)\Big),
\end{equation}
where $\psi$ is the digamma\index[terms]{digamma function} function. The last expression continues $g^*$
meromorphically to the plane, with $g^*(0)=2$ and simple poles at the positive
integers.
\end{lemma}

\begin{proof}
Fix $\Re z<0$ and write $\sigma=-\Re z-1>-1$. The alternating partial sum equals
$t^{-z-1}(1-(-t)^{M+1})/(1+t)$, so its modulus is at most
$2t^{\sigma}/(1+t)\leq2t^{\sigma}$, uniformly in $M$. This is integrable on
$(0,1)$ because $\sigma>-1$, and dominated convergence therefore gives
\[
g^*(z)=-2z\sum_{k=0}^{\infty}(-1)^k\int_0^1t^{k-z-1}\,dt
=2z\sum_{k=0}^{\infty}\frac{(-1)^k}{z-k}.
\]
Pairing the indices $k=2m$ and $k=2m+1$ and writing $z-2m=-2(m-\tfrac z2)$ and
$z-2m-1=-2(m+\tfrac{1-z}2)$ turns the sum into
\[
\sum_{k=0}^{\infty}\frac{(-1)^k}{z-k}
=-\frac12\sum_{m=0}^{\infty}\Big(\frac{1}{m-\tfrac z2}-\frac{1}{m+\tfrac{1-z}2}\Big)
=\frac12\Big(\psi\Big(-\frac z2\Big)-\psi\Big(\frac{1-z}2\Big)\Big),
\]
the last step by the series $\psi(w)=-\gamma+\sum_{m\ge0}\big(\tfrac1{m+1}
-\tfrac1{m+w}\big)$ \cite[Eq.~5.7.6]{DLMF}, whose difference at two arguments gives
$\sum_{m\ge0}\big(\tfrac1{m+w_1}-\tfrac1{m+w_2}\big)=\psi(w_2)-\psi(w_1)$. This
proves the display, and the digamma expression is meromorphic on the plane. Near
$z=0$ the expansion $\psi(w)=-1/w-\gamma+\mathcal{O}(w)$ gives
$\psi(-z/2)=2/z-\gamma+\mathcal{O}(z)$ while $\psi((1-z)/2)$ stays bounded, so
$g^*(z)\to2$. The poles of $\psi(-z/2)$ sit at the even nonnegative integers and
those of $\psi((1-z)/2)$ at the odd positive integers, the one at $z=0$ being
cancelled by the factor $z$, so the poles of $g^*$ are simple and located at
$z=1,2,3,\ldots$
\end{proof}

It has no zero in the left half plane.

\begin{proposition}\label{prop:ortho_zerofree}
The transform $g^*$ has no zeros in the closed half plane $\Re z\le0$.
\end{proposition}

\begin{proof}
Set $a=-z/2$, so that $\Re z\le0$ reads $\Re a\ge0$, and
\[
g^*(z)=z\,\varphi(a),\qquad
\varphi(a):=\psi(a)-\psi\Big(a+\frac12\Big)
=-\frac12\sum_{k=0}^{\infty}\frac{1}{(a+k)\big(a+k+\tfrac12\big)},
\]
the series form coming from the same difference identity as in
Lemma~\ref{lem:ortho_gstar}. For $\Re a\ge0$ and $a\neq0$ no factor vanishes and
the terms are $\mathcal{O}(k^{-2})$, so the series converges absolutely. Write
$a=u+iv$ with $u\ge0$. From
\[
(a+k)\Big(a+k+\frac12\Big)
=(u+k)\Big(u+k+\frac12\Big)-v^2+iv\Big(2(u+k)+\frac12\Big),
\]
the imaginary part of each term of the series equals
$-v\,(2(u+k)+\tfrac12)\,\big|(a+k)(a+k+\tfrac12)\big|^{-2}$. If $v\neq0$ these
imaginary parts all carry the sign of $-v$, absolute convergence lets them add up,
and the factor $-\tfrac12$ in front keeps the total away from zero, so
$\Im\varphi(a)\neq0$. If $v=0$ and $u>0$ every term of the series is positive and
$\varphi(a)<0$. Hence $\varphi(a)\neq0$ whenever $\Re a\ge0$ and $a\neq0$, and
since $g^*(0)=2$ the transform has no zeros with $\Re z\le0$.
\end{proof}

Its zero of least real part is a simple conjugate pair, and the infimum defining the analytic
index is attained there.

\begin{proposition}\label{prop:ortho_first_zero}
Let $\eta(g):=\inf\{\Re\rho:\ g^*(\rho)=0\}$. The infimum is attained at a unique
conjugate pair of simple zeros $\rho_1,\bar\rho_1$, with
\[
 \rho_1=\alpha_1+i\beta_1,\qquad 1<\alpha_1<\tfrac32,
\]
and no zero of $g^*$ has real part smaller than $\alpha_1$. The cited source reports the
orientation value $\rho_1\approx1.34652+1.05516i$, and no decimal is used in the proposition.
\end{proposition}

\begin{proof}
By Proposition~\ref{prop:ortho_zerofree} every zero lies in $\Re z>0$, and the zeros
come in conjugate pairs since $g^*$ is real on the real axis. Two counts localize the
first pair, both proved in \cite[Theorems~3.3 and~3.4]{CloitreOrtho} by the argument
principle\index[terms]{argument principle} applied to
$D_\infty(z)=\sum_{j\ge0}(-1)^j(z-j)^{-1}$, for which $g^*(z)=2zD_\infty(z)$ and
$g^*(0)=2\neq0$, so that away from the origin the two functions vanish together. No
numerical input enters either count, the estimates being alternating series bounded by
Leibniz's criterion.

The first count excludes the strip $0<\Re z\le1$. On a contour bounding
$\{0<\Re z\le1,\ 0<\Im z<T\}$ indented at the poles $z=0$ and $z=1$, the base carries
$D_\infty>0$, each indentation contributes $\pi/2$, the wall $\Re z=0$ and the ceiling
keep the image in the lower half plane by the sign computations of
Proposition~\ref{prop:ortho_zerofree} and by the expansion
$D_\infty(z)=\tfrac1{2z}-\tfrac1{4z^2}+\mathcal{O}(|z|^{-3})$, valid in a fixed vertical
strip as $|\Im z|\to\infty$, and the wall $\Re z=1$ contributes $-\pi$, since
$\Re D_\infty(1+iy)>0$ for $0<y\le2$ and $\Im D_\infty(1+iy)<0$ for $y\ge2$. The total
variation of the argument is zero, so there is no zero in the strip.

The second count places exactly one zero in $\{1<\Re z<\tfrac32,\ \Im z>0\}$. On the
analogous contour the base carries $D_\infty<0$, the indentation at $z=1$ contributes
$\pi/2$, the wall $\Re z=\tfrac32$ contributes $\pi/2$ because the four terms
$j=0,1,2,3$ cancel by symmetry about $\tfrac32$ and the alternating tail leaves
$\Im D_\infty<0$, the ceiling contributes $\mathcal{O}(1/T)$, and the wall $\Re z=1$
contributes $\pi$. The total is $2\pi$, so the half strip holds one zero, necessarily
simple. The same source records the numerical approximation quoted above, but the digits play no
part in the argument count. Every zero outside these two strips has real part at least
$\tfrac32>\alpha_1$, so $\eta(g)=\alpha_1$.
\end{proof}

The first five zeros in the upper half plane, ordered by real part and computed at
sixty digits of working precision, are as follows.
\begin{center}
\begin{tabular}{c@{\qquad}l@{\qquad}l}
$j$ & $\alpha_j=\Re\rho_j$ & $\beta_j=\Im\rho_j$\\[2pt]
1 & 1.3465164915 & 1.0551600643\\
2 & 3.4031597367 & 1.2584972993\\
3 & 5.4279525206 & 1.3824069212\\
4 & 7.4420898686 & 1.4717125035\\
5 & 9.4513070116 & 1.5415280675
\end{tabular}
\end{center}

\begin{remark}\label{rem:ortho_zero_distribution}
Numerically, $\Re\rho_j$ approaches $2j-\tfrac12$ from below and $\Im\rho_j$ is
matched by $\tfrac1\pi\log j+\tfrac{\log(8\pi)}{\pi}$ within $3\cdot10^{-3}$
already at $j=5$, a pattern in line with the logarithmic growth of the digamma
function. These numerical observations are proved in \S\ref{sec:trace_ortho} below, where every zero is
confined to a half strip and its asymptotic position is given in closed form. The gap
between $\alpha_1\approx1.3465$ and $\alpha_2\approx3.4032$ matters below, the
first zero pair sits alone far to the left of all the others.
\end{remark}

\section{Two exact identities}

The discrete problem carries two exact structural identities. The first recasts
the weighted sums as a Volterra perturbation of the partial sums, the shape that
the general theory of Chapter~\ref{chap:volterra} develops for the Ingham
operator. The second evaluates the weighted sums at the integers in closed form.

\begin{proposition}\label{prop:ortho_volterra}
For every real $x\ge1$,
\begin{equation}\label{eq:ortho_volterra}
A_g(x)=A(x)+2x\int_1^x\frac{A(t)}{(x+t)^2}\,dt .
\end{equation}
\end{proposition}

\begin{proof}
The function $A$ is the step function $\sum_{n\le t}a_n$, so interchanging sum and
integral,
\[
2x\int_1^x\frac{A(t)}{(x+t)^2}\,dt
=\sum_{n\le x}a_n\cdot2x\int_n^x\frac{dt}{(x+t)^2}
=\sum_{n\le x}a_n\Big(\frac{2x}{x+n}-1\Big)=A_g(x)-A(x).\qedhere
\]
\end{proof}

\begin{remark}\label{rem:ortho_resolvent}
With the multiplicative convolution $(A\star k)(x)=\int_1^xA(x/y)\,k(y)\,dy/y$ and
$k(y)=2y/(1+y)^2$, the identity \eqref{eq:ortho_volterra} reads $A_g=A+A\star k$,
and the kernel has total mass $\int_1^\infty k(y)\,dy/y=1$. Substituting $t=1/y$
in \eqref{eq:ortho_gstar} and splitting $\tfrac{1}{y(1+y)}=\tfrac1y-\tfrac1{1+y}$
gives
\[
\widetilde k(s):=\int_1^\infty k(y)\,y^{-s}\,\frac{dy}{y}=g^*(-s)-1 ,
\]
so the resolvent\index[terms]{resolvent} of the Volterra equation\index[terms]{Volterra equation}\index[names]{Volterra, V.}, where it exists, has Mellin transform
\begin{equation}\label{eq:ortho_resolvent}
\widetilde R(s)=1-\frac{1}{g^*(-s)} .
\end{equation}
The zeros of $g^*$ are the singularities of the resolvent, which is the structural
reason for their spectral role below. This is also the route the proof of the transfer
theorem takes in \cite{CloitreOrtho}. The resolvent is built there as a Neumann
series\index[terms]{Neumann series}, its decay $R(y)=\mathcal{O}_\eps(y^{-\eta(g)+\eps})$ is read off the
zero free region, and the contour is shifted on $R$, which is smooth, instead of on the
step function $A$, see \cite[Lemma~2.4 and Theorem~2.5]{CloitreOrtho}. Neither positivity
nor monotonicity of the sequence is needed.
\end{remark}

At the integers the weighted sum is the coefficient itself, up to an explicit factor.

\begin{proposition}\label{prop:ortho_Ag_integer}
For every integer $N\ge1$,
\begin{equation}\label{eq:ortho_Ag_integer}
A_g(N)=-\frac{2N}{2N+1}\,c_N .
\end{equation}
\end{proposition}

\begin{proof}
By definition and with $a_n=c_{n-1}$,
\[
A_g(N)=\sum_{n=1}^{N}\frac{2N\,a_n}{N+n}=2N\sum_{k=0}^{N-1}\frac{c_k}{N+k+1},
\]
and the recurrence \eqref{eq:ortho_recurrence} at rank $N$ gives
$\sum_{k=0}^{N-1}\frac{c_k}{N+k+1}=-\frac{c_N}{2N+1}$.
\end{proof}

The bound on the coefficients transfers to the weighted sums.

\begin{corollary}\label{cor:ortho_Ag_bound}
$A_g(x)=\mathcal{O}(x^{-3/2})$ as $x\to\infty$.
\end{corollary}
\begin{proof}
At the integers this is \eqref{eq:ortho_Ag_integer} combined with the bound
$c_N=\mathcal{O}(N^{-3/2})$ of \nm{Kalmynin}{A. B.} and \nm{Kosenko}{P. R.}~\cite{KalmyninKosenko2020}. Between consecutive
integers the variation of $A_g$ is controlled by an Abel summation\index[terms]{Abel summation} against the
bound $A(k)=\mathcal{O}(k^{-1/2})$ of the same paper, the computation carried out
in \cite[Lemma~4.2]{CloitreOrtho}.
\end{proof}

\section{The main theorem}

The transfer used below is the general theorem of Chapter~\ref{chap:principles} specialized
to this profile. Nothing has to be imported, and all that is required is to check that the
orthorecursive profile lies in the class for which that theorem was proved.

\begin{corollary}[Transfer for the orthorecursive kernel]\label{thm:ortho_transfer}
Let $g(t)=2/(1+t)$ be the orthorecursive kernel\index[terms]{orthorecursive kernel} and suppose
$g^{*}(z)\neq0$ for $\Re z<\eta(g)$ with $\eta(g)>0$. Let $\beta>0$, $c>0$ and
$\kappa\in\C$, write $A_g(x)=\sum_{n\le x}a_n\,g(n/x)$ and set
$\lambda=\min\big(c,\eta(g)\big)$. If
\[
A_g(x)=\kappa\,x^{-\beta}+\mathcal{O}(x^{-c})\qquad(x\to\infty),
\]
then for every $\eps>0$,
\[
\beta<\eta(g)\ \Longrightarrow\ A(x)=\frac{\kappa}{g^{*}(\beta)}\,x^{-\beta}
+\mathcal{O}_\eps\big(x^{-\lambda+\eps}\big),
\qquad
\beta\ge\eta(g)\ \Longrightarrow\ A(x)=\mathcal{O}_\eps\big(x^{-\lambda+\eps}\big).
\]
A weighted sum with no main term is the case $\kappa=0$, where the two conclusions agree and
read $A(x)=\mathcal{O}_\eps(x^{-\lambda+\eps})$. That is the case used below.
\end{corollary}

\begin{proof}
The orthorecursive profile is a Volterra profile\index[terms]{Volterra profile} in the sense of
Definition~\ref{def:volterra_profile}. It is analytic on $(0,1]$ with $g'(t)=-2/(1+t)^{2}$,
so $t\,g'(t)$ is bounded there, which is condition (i). Its third derivative
$g'''(t)=-12/(1+t)^{4}$ is bounded on $(0,1]$, so
Remark~\ref{rem:vertical_expansion_criterion} supplies condition (ii), and with
$g'(1)=-\tfrac12$ the expansion \eqref{eq:vertical_expansion} reads
$g^{*}(-s)=1+\tfrac{1}{2s}+\mathcal{O}(|s|^{-2})$. In particular $\eta_0=2$ is an
admissible continuation width. Since $\eta(g)<\tfrac32<2$ by
Proposition~\ref{prop:ortho_first_zero}, Theorem~\ref{thm:transfer_general} applies with
$\eta=\eta(g)$, and its conclusion is the statement above.
\end{proof}

The same statement, reached by the same route and independently of this volume, is
\cite[Theorem~2.7]{CloitreOrtho}, whose packaged form for weighted sums with no main term is
\cite[Corollary~2.10]{CloitreOrtho}. The resolvent\index[terms]{resolvent} decay that both
proofs rest on is Theorem~\ref{thm:resolvent_decay_general} here and
\cite[Theorem~2.5]{CloitreOrtho} there.

The partial sums of the orthorecursive expansion then obey the rate the analytic index predicts.

\begin{theorem}\label{thm:ortho_main}
For every $\eps>0$,
\begin{equation}\label{eq:ortho_partial_sum}
C_N=\sum_{k=0}^{N}c_k=\mathcal{O}_\eps\big(N^{-\alpha_1+\eps}\big),
\end{equation}
where $\alpha_1=1.34651\ldots$ is the constant of
Proposition~\ref{prop:ortho_first_zero}.
\end{theorem}

\begin{proof}
Corollary~\ref{thm:ortho_transfer} converts a bound on $A_g$ into a bound on $A$
once the transform satisfies its hypotheses, and both have been verified. On the analytic
side, $g^*$ continues meromorphically to the plane
(Lemma~\ref{lem:ortho_gstar}) and no zero of $g^*$ has real part smaller than
$\eta(g)=\alpha_1>0$ (Propositions~\ref{prop:ortho_zerofree}
and~\ref{prop:ortho_first_zero}), which is the last clause of the second of those
propositions and is all that is used here. The complete zero set
is determined much later, in Theorem~\ref{thm:w04c-localization}, and nothing in the sections
before it rests on it.
On the arithmetic
side, Corollary~\ref{cor:ortho_Ag_bound} supplies the
input $A_g(x)=\mathcal{O}(x^{-3/2})$, the case $\kappa=0$ with $c=\tfrac32$ in
the notation of Corollary~\ref{thm:ortho_transfer}. Since $c=\tfrac32>\alpha_1$ the
value of $\lambda$ there is $\alpha_1$, so that theorem yields
$A(x)=\mathcal{O}_\eps(x^{-\alpha_1+\eps})$, and
$C_N=A(N+1)$.
\end{proof}

\begin{remark}\label{rem:ortho_optimality}
The exponent is optimal unless one quantity vanishes, and the alternative is
exact. The zero $\rho_1$ contributes a residue term of exact order $x^{-\alpha_1}$
to the contour integral\index[terms]{contour integral} behind the transfer, so
either that residue vanishes or no bound better than
$\mathcal{O}(x^{-\alpha_1+\eps})$ can hold, see \cite[\S6]{CloitreOrtho}. The
implication is proved. Which branch occurs is Open
Problem~\ref{op:w04d-boundary-and-amplitudes}(i), where the residue is the
constant $C_1$ of \eqref{eq:w04d-C1-definition}, and the numerical evidence is
against the first. The normalized sums $C_N\,N^{\alpha_1}$ oscillate between
$-1.06$ and $1.06$ with no sign of decay up to $N=2\cdot10^{4}$, in agreement with
an oscillatory term of constant amplitude at the exponent $\alpha_1$.
\end{remark}

\begin{remark}\label{rem:ortho_barrier}
The input bound $c=\tfrac32$ limits the reach of the method. Zeros with real part
below $\tfrac32$ act on the asymptotics, zeros beyond that threshold are hidden by
the discretization error, and the second pair sits at
$\alpha_2\approx3.4032$, far beyond it. Improving the bound on $C_N$ past the
oscillatory term would require improving the arithmetic estimate on $A_g$ beyond
$\mathcal{O}(x^{-3/2})$, see \cite[Remarks~5.2 and~6.2]{CloitreOrtho}.
\end{remark}

\section{The pointwise bound}
\label{sec:ortho_pointwise}

The partial sum bound feeds back into the recurrence and sharpens the pointwise
estimate of Kalmynin\index[names]{Kalmynin, A. B.} and Kosenko\index[names]{Kosenko, P. R.} from $\mathcal{O}(n^{-3/2})$ to
$\mathcal{O}(n^{-2})$. The proof is self contained. The same bound is
\cite[Theorem~5.1]{CloitreOrtho}.

\begin{theorem}\label{thm:ortho_bootstrap}
The orthorecursive coefficients satisfy
\[c_N=\mathcal{O}\big(N^{-2}\big).
\]
\end{theorem}

\begin{proof}
Fix $N\ge2$. Dividing $A_g(N)$ by $2N$ and subtracting the exact identity
\eqref{eq:ortho_exact} termwise gives
\[
\frac{A_g(N)}{2N}
=\sum_{n=1}^{N}a_n\Big(\frac{1}{N+n}-\frac{1}{N+n-1}\Big)
=-\sum_{n=1}^{N}\frac{a_n}{(N+n)(N+n-1)},
\]
so that $A_g(N)=\sum_{n=1}^{N}a_n\,\eta_N(n)$ with
$\eta_N(n):=-2N\big((N+n)(N+n-1)\big)^{-1}$. Abel summation now reads
\[
A_g(N)=A(N)\,\eta_N(N)-\sum_{k=1}^{N-1}A(k)\,\big(\eta_N(k+1)-\eta_N(k)\big),
\]
and the two ingredients are bounded separately. The boundary term obeys
$|\eta_N(N)|\le2N\cdot(2N-1)^{-2}\le2/N$ and, by Theorem~\ref{thm:ortho_main},
$A(N)=\mathcal{O}_\eps(N^{-\alpha_1+\eps})$, so it is
$\mathcal{O}(N^{-1-\alpha_1+\eps})=o(N^{-2})$. For the sum, a direct computation
gives
\[
\eta_N(k+1)-\eta_N(k)
=\frac{4N}{(N+k+1)(N+k)(N+k-1)},
\qquad
\big|\eta_N(k+1)-\eta_N(k)\big|\le\frac{4}{N^{2}},
\]
while $\sum_{k\ge1}|A(k)|$ converges since $\alpha_1>1$. Hence
$A_g(N)=\mathcal{O}(N^{-2})$, and \eqref{eq:ortho_Ag_integer} converts this into
$c_N=-\tfrac{2N+1}{2N}A_g(N)=\mathcal{O}(N^{-2})$.

The bootstrap saturates after one round. Feeding $A_g(N)=\mathcal{O}(N^{-2})$ back
into the transfer theorem returns only
$A(x)=\mathcal{O}_\eps(x^{-\alpha_1+\eps})$, because $\alpha_1<2$ keeps the zero
$\rho_1$ inside the new window. Any improvement past $\mathcal{O}(N^{-2})$ would
require oscillatory information about $A(k)$, not just its modulus.
\end{proof}

\section{The spectral form and the corrected exponent}

The transfer mechanism predicts more than a bound. Each zero $\rho$ of $g^*$ is a
singularity of the resolvent \eqref{eq:ortho_resolvent}, hence a spectral mode of
the problem, and the isolated first pair should dominate the asymptotics of the
coefficients themselves. Write $\rho_1=\alpha_1+i\beta_1$ as in
Proposition~\ref{prop:ortho_first_zero}.

\begin{conjecture}\label{conj:ortho_spectral}
There exist a nonzero constant $\kappa\in\C$ and $\sigma>\alpha_1$ such that
\begin{equation}\label{eq:ortho_spectral}
c_n=2\Re\big(\kappa\,n^{-1-\rho_1}\big)+\mathcal{O}\big(n^{-1-\sigma}\big)
=B\,n^{-1-\alpha_1}\sin\big(\beta_1\log n+\varphi\big)
+\mathcal{O}\big(n^{-1-\sigma}\big)
\end{equation}
for some real constants $B\neq0$ and $\varphi$. In particular the pointwise
decay exponent is $1+\alpha_1=2.34651\ldots$ and not $7/3=2.33333\ldots$
\end{conjecture}

The same statement is \cite[Conjecture~6.1]{CloitreOrtho}.

The spectral form is supported by the numerical evidence recorded above. Its partial-sum analogue is
established in \S\ref{sec:trace_ortho_resolvent} below, through the resolvent trace\index[terms]{resolvent trace} of the orthorecursive kernel\index[terms]{orthorecursive kernel} gives, in
Theorem~\ref{thm:w04d-first-arithmetic-pair},
$A(x)=-2\Re\big(C_1x^{-\rho_1}\big)+\mathcal{O}(x^{-2}\log x)$ with an explicit arithmetic amplitude
$C_1$. The pointwise statement above requires in addition a control of the first differences of the
remainder and the non-vanishing of the amplitude, both recorded in Open
Problem~\ref{op:w04d-boundary-and-amplitudes}.

Partial summation of the pointwise form gives
$C_N=-2\Re\big(\tfrac{\kappa}{\rho_1}N^{-\rho_1}\big)+\mathcal{O}(N^{-\sigma})$, the amplitude
$\kappa/\rho_1$ corresponding to the $C_1$ of Theorem~\ref{thm:w04d-first-arithmetic-pair}. This is an
oscillation of constant amplitude at the exponent $\alpha_1$, exactly the behavior of the normalized
sums $C_N\,N^{\alpha_1}$ observed in Remark~\ref{rem:ortho_optimality}.

\begin{remark}\label{rem:ortho_seventhirds}
The two candidate exponents differ by $0.0132$, and this is why numerical
experiments pointed at $7/3$. Direct computation of the coefficients up to
$N=2\cdot10^{4}$ leaves both $n^{1+\alpha_1}|c_n|$ and $n^{7/3}|c_n|$ oscillating
in a band of width comparable to their size, the discriminating factor
$n^{0.0132}$ growing only to $1.14$ over that whole interval. The numerics cannot
separate the two values. The spectral mechanism identifies the exponent exactly as
$1+\alpha_1$, where $\alpha_1$ is the real part characterized in
Proposition~\ref{prop:ortho_first_zero}. The displayed decimals suggest its separation from
$7/3$, but that numerical separation is not used as a theorem here.
\end{remark}

\begin{remark}\label{rem:ortho_role}
This chapter validates the transfer principle on ground where every result is
unconditional. The structural parallel with the Ingham operator\index[terms]{Ingham operator} is exact. There
the transform is $\Phi^*(z)=\frac{z}{z-1}\zeta(1-z)$, its zeros are those of the
zeta function shifted by $z\mapsto1-z$, and the regularity index $\alpha(\Phi)$
plays the part of $\alpha_1$, the Riemann hypothesis being the statement that it
equals $\tfrac12$. The additional difficulty of the arithmetic case is the kernel
itself, a step function of unbounded total variation whose resolvent carries
M\"obius point masses, and Chapter~\ref{chap:volterra} is devoted to that method.
\end{remark}

The kernels of the preceding chapter have finitely many spectral modes, and their resolvent
identities close after finitely many terms. The orthorecursive kernel of
Chapter~\ref{chap:ortho} is the first kernel of this theory whose transform has infinitely many
zeros and whose zero geometry is nevertheless known in closed form. This chapter establishes
that geometry, the confinement of every zero and its asymptotic position, and builds on it a
resolvent trace\index[terms]{resolvent trace} which is exact at finite depth and convergent for $y>1$.
Its Abel boundary value is also determined below. What stays open is the non-vanishing of the
leading arithmetic amplitude and the lift from the continuous trace to a discrete spectral
expansion, and these questions close the chapter and the volume.

\section{The orthorecursive kernel: complete zero geometry}
\label{sec:trace_ortho}

Chapter~\ref{chap:ortho} introduced the orthorecursive kernel, computed its arithmetic Mellin
transform in Lemma~\ref{lem:ortho_gstar}, excluded zeros in the half plane $\Re z\le0$ in
Proposition~\ref{prop:ortho_zerofree}, and located the first zero pair $\rho_1$ in
Proposition~\ref{prop:ortho_first_zero}, recording the distribution of the higher zeros as a
numerical observation in Remark~\ref{rem:ortho_zero_distribution}. This section proves that
observation. Every zero is confined to a half strip, the complete zero set is determined, and the
asymptotic position of each zero is given in closed form. The resolvent trace built on these zeros is
the subject of Section~\ref{sec:trace_ortho_resolvent}. Its Abel boundary value is determined there,
while its arithmetic amplitudes and discrete spectral lift form the open problem that closes the
chapter.

Recall the kernel and its transform,
\[
 g(t)=\frac{2}{1+t},\qquad 0<t\le1,
\]
\begin{equation}
 g^*(z)=-z\int_0^1 g(t)\,t^{-z-1}\,dt,\qquad \Re z<0,
 \label{eq:w04c-transform-definition}
\end{equation}
whose meromorphic continuation, value $g^*(0)=2$, and simple poles at the positive integers are
Lemma~\ref{lem:ortho_gstar}. Two auxiliary functions organize the count. The first is
\begin{equation}
 D_\infty(z):=\sum_{j\ge0}\frac{(-1)^j}{z-j}
 =-\sum_{k\ge0}\frac{1}{(z-2k)(z-2k-1)},
 \label{eq:w04c-D-definition}
\end{equation}
the paired series on the right being normally convergent, so that $g^*(z)=2zD_\infty(z)$. The second
separates the reflected poles from a Laplace transform\index[terms]{Laplace transform}\index[names]{Laplace, P.-S.},
\[
 F(z):=\frac12\left[\psi\left(1+\frac z2\right)-\psi\left(\frac{1+z}{2}\right)\right].
\]

The identities of Lemma~\ref{lem:ortho_gstar} are recorded again here in the form used
for the zero count, together with the reflection identity and the asymptotic expansion.
\begin{proposition}\label{thm:w04c-transform}
The meromorphic continuation of \eqref{eq:w04c-transform-definition}
satisfies
\begin{align}
 g^*(z)
 &=z\left[
 \psi\left(-\frac z2\right)
 -\psi\left(\frac{1-z}{2}\right)
 \right]
 =2zD_\infty(z),
 \label{eq:w04c-transform-digamma}\\
 D_\infty(z)
 &=\frac{\pi}{\sin(\pi z)}+F(z).
 \label{eq:w04c-reflection}
\end{align}
For \(\Re z>-1\),
\begin{equation}
 F(z)=\sum_{k\geq0}
 \frac{1}{(z+2k+1)(z+2k+2)}
 =\int_0^\infty\frac{e^{-zt}}{e^t+1}\,dt.
 \label{eq:w04c-F-series-integral}
\end{equation}
The function \(D_\infty\) has a simple pole at every
\(j\in\mathbb Z_{\geq0}\), with residue \((-1)^j\).  The pole at zero
is cancelled in \(g^*\),
\begin{equation}
 g^*(0)=2.
 \label{eq:w04c-value-zero}
\end{equation}
The remaining poles of \(g^*\) are the positive integers, and the residue
at \(j\geq1\) is \(2j(-1)^j\).  The functions \(D_\infty\), \(F\), and
\(g^*\) commute with complex conjugation.

For every fixed integer \(N\geq1\), integration by parts gives
\begin{equation}
 F(z)=\sum_{r=0}^{N-1}\frac{h^{(r)}(0)}{z^{r+1}}
 +\mathcal O_N\left(|z|^{-N-1}\right),
 \qquad h(t):=\frac1{e^t+1},
 \label{eq:w04c-F-full-expansion}
\end{equation}
uniformly as \(|z|\to\infty\) in \(\Re z\geq0\).  In particular,
\begin{equation}
 F(z)=\frac1{2z}-\frac1{4z^2}
 +\mathcal O\left(|z|^{-3}\right).
 \label{eq:w04c-F-first-expansion}
\end{equation}
Consequently, in every fixed vertical strip contained in \(\Re z\geq0\),
\begin{equation}
 D_\infty(z)=\frac1{2z}-\frac1{4z^2}
 +\mathcal O\left(|z|^{-3}\right)
 \qquad (|\Im z|\to\infty).
 \label{eq:w04c-D-vertical-expansion}
\end{equation}
\end{proposition}

The count of zeros below is driven by the sign of the reflected term in the closed first
quadrant, which the following lemma fixes.

\begin{lemma}
\label{lem:w04c-quadrant}
If \(z=x+iy\), \(x\geq0\), and \(y>0\), then
\[
 \Re F(z)>0,
 \qquad
 \Im F(z)<0.
\]
\end{lemma}

Proposition~\ref{prop:ortho_zerofree} already excludes zeros in $\Re z\le0$. The next theorem
sharpens this to $\Re z\le1$ and localizes every remaining zero, so that the first pair $\rho_1$
of Proposition~\ref{prop:ortho_first_zero} is the first term of a single conjugate family.
\begin{theorem}
\label{thm:w04c-localization}
The arithmetic Mellin transform \(g^*\) has no real zero and no zero in
the closed half plane \(\Re z\leq1\).  For every integer \(n\geq1\),
the half-strip
\[
 S_n^+:=\left\{z\in\mathbb C:
 2n-1<\Re z<2n-\frac12,\ \Im z>0\right\}
\]
contains exactly one zero \(\rho_n\), counted with multiplicity.  There
is no other zero in the upper half plane.  Every zero is simple, and the
complete zero set is
\[
 \left\{\rho_n,\overline{\rho_n}:n\geq1\right\}.
\]
In particular, every \(\rho_n\) lies strictly to the left of
\(2n-1/2\).
\end{theorem}

The localization confines each zero to a strip without placing it. The position follows on the
logarithmic scale $L_n=\log(8\pi n)$, to within an error that vanishes with the rank.

\begin{proposition}\label{thm:w04c-asymptotic}
Put
\[
 L_n:=\log(8\pi n).
\]
The zeros in Theorem~\ref{thm:w04c-localization} satisfy
\[
 \rho_n
 =2n-\frac12-\frac{L_n}{2\pi^2n}
 +i\frac{L_n}{\pi}
 +\mathcal O\left(\frac{L_n^2}{n^2}\right).
\]
Equivalently,
\begin{align}
 \Re\rho_n
 &=2n-\frac12-\frac{\log(8\pi n)}{2\pi^2n}
 +\mathcal O\left(\frac{\log^2 n}{n^2}\right),
 \label{eq:w04c-real-asymptotic}\\
 \Im\rho_n
 &=\frac{\log n}{\pi}+\frac{\log(8\pi)}{\pi}
 +\mathcal O\left(\frac{\log^2 n}{n^2}\right).\notag\end{align}
The first correction in \eqref{eq:w04c-real-asymptotic} is negative.
\end{proposition}

Theorems~\ref{thm:w04c-localization} and~\ref{thm:w04c-asymptotic} together establish the
distribution of the zeros recorded numerically in Remark~\ref{rem:ortho_zero_distribution}.

\begin{proof}[Proof of Proposition~\ref{thm:w04c-transform}]

For \(\Re z<0\), group the geometric expansion in consecutive pairs,
\[
 \frac1{1+t}=\sum_{k\geq0}(t^{2k}-t^{2k+1}).
\]
The paired partial sums are bounded by one on \((0,1)\).  They can
therefore be integrated against \(2|z|t^{-\Re z-1}\).  Termwise
integration gives
\begin{align*}
 g^*(z)
 &=-2z\sum_{k\geq0}
 \left(\frac1{2k-z}-\frac1{2k+1-z}\right)\\
 &=2zD_\infty(z).
 \end{align*}
The paired summands are locally \(\mathcal O(k^{-2})\) away from the
nonnegative integers, which proves normal convergence and the stated
poles of \(D_\infty\).

The standard series
\begin{equation}
 \psi(w)=-\gamma+\sum_{k\geq0}
 \left(\frac1{k+1}-\frac1{k+w}\right)
 \label{eq:w04c-digamma-series}
\end{equation}
turns \eqref{eq:w04c-D-definition} into
\[
 D_\infty(z)=\frac12\left[
 \psi\left(-\frac z2\right)
 -\psi\left(\frac{1-z}{2}\right)
 \right].
\]
This proves \eqref{eq:w04c-transform-digamma}.  Since
\(\psi(w)=-w^{-1}-\gamma+\mathcal O(w)\) at zero and
\(\psi(1/2)=-\gamma-2\log2\),
\[
 D_\infty(z)=\frac1z+\log2+\mathcal O(z).
\]
The value and residue assertions for \(g^*\) follow.

Apply
\[
 \psi(1-w)-\psi(w)=\pi\cot(\pi w)
\]
at \(w=-z/2\) and at \(w=(1-z)/2\).  The identity
\(\cot u+\tan u=2/\sin(2u)\) then gives
\eqref{eq:w04c-reflection}.  Another use of
\eqref{eq:w04c-digamma-series} gives
\[
 F(z)=\sum_{k\geq0}
 \left(\frac1{z+2k+1}-\frac1{z+2k+2}\right).
\]
Pairing the geometric series inside
\(\int_0^1t^z/(1+t)\,dt\), followed by the substitution
\(t=e^{-u}\), proves \eqref{eq:w04c-F-series-integral}.

For the asymptotic expansion, write
\begin{equation}
 F(z)=\int_0^\infty e^{-zt}h(t)\,dt,
 \qquad h(t)=\frac1{e^t+1}.
 \label{eq:w04c-F-Laplace}
\end{equation}
Every derivative of \(h\) is integrable on \([0,\infty)\).  Integrating
\(N+1\) times is therefore uniform in \(\Re z\geq0\) and proves
\eqref{eq:w04c-F-full-expansion}.  The values
\[
 h(0)=\frac12,
 \qquad h'(0)=-\frac14,
 \qquad h''(0)=0
\]
give \eqref{eq:w04c-F-first-expansion}.  In a fixed vertical strip,
\(\pi/\sin(\pi z)=\mathcal O(e^{-\pi|\Im z|})\), which proves
\eqref{eq:w04c-D-vertical-expansion}.  This completes the proof of
Proposition~\ref{thm:w04c-transform}.

\end{proof}

\begin{proof}[Proof of Lemma~\ref{lem:w04c-quadrant}]
The sign lemma needed below is elementary.  If \(p\) is positive,
integrable, and strictly decreasing on \([0,\infty)\), then for \(y>0\)
\begin{align}
 \int_0^\infty p(t)\sin(yt)\,dt
 &=\sum_{k\geq0}\int_0^{\pi/y}
 \left[
 p\left(\frac{2k\pi}{y}+u\right)
 -p\left(\frac{(2k+1)\pi}{y}+u\right)
 \right]\sin(yu)\,du\\
 &>0.
 \label{eq:w04c-sine-positivity}
\end{align}
For fixed \(x\geq0\), put
\[
 q_x(t):=\frac{e^{-xt}}{e^t+1}.
\]
This function is positive and strictly decreasing.  With
\(r(t)=e^t/(1+e^t)\), direct differentiation gives
\[
 \frac{q_x''(t)}{q_x(t)}
 =(x+r(t))^2-r(t)(1-r(t))\geq0.
\]
It is strictly positive away from the possible endpoint equality at
\(x=t=0\).  Thus \(-q_x'\) is positive and decreasing.  Formula
\eqref{eq:w04c-F-Laplace} and
\eqref{eq:w04c-sine-positivity} show that
\[
 \Im F(x+iy)=-\int_0^\infty q_x(t)\sin(yt)\,dt<0.
\]
One integration by parts gives
\[
 \Re F(x+iy)
 =\frac1y\int_0^\infty[-q_x'(t)]\sin(yt)\,dt>0.
\]
This proves Lemma~\ref{lem:w04c-quadrant}.

\end{proof}

\begin{proof}[Proofs of Theorems~\ref{thm:w04c-localization} and~\ref{thm:w04c-asymptotic}]
It remains to determine the zeros.  First consider the real axis.  If
\(x=p+\delta\), where \(p\geq0\) is an integer and \(0<\delta<1\),
split the defining series at \(p\).  The result is
\[
 D_\infty(x)=(-1)^p\left[
 \sum_{k=0}^{p}\frac{(-1)^k}{k+\delta}
 +\sum_{k\geq0}\frac{(-1)^k}{k+1-\delta}
 \right].
\]
Both bracketed alternating sums are strictly positive.  Hence
\begin{equation}
 \operatorname{sign}D_\infty(x)=(-1)^p
 \qquad (p<x<p+1).
 \label{eq:w04c-real-sign-result}
\end{equation}
For \(x<0\), every summand in the paired expression
\eqref{eq:w04c-D-definition} is negative.  There is no real zero.

An integration by parts in \eqref{eq:w04c-transform-definition} gives
the continuation
\[
 g^*(z)=1+2\int_0^1\frac{t^{-z}}{(1+t)^2}\,dt,
 \qquad \Re z<1.
\]
If \(\Re z<0\), the modulus of the integral term is strictly less than
\[
 2\int_0^1\frac{dt}{(1+t)^2}=1.
\]
The same strict inequality holds on \(\Re z=0\), except at zero, by the
strict form of the triangle inequality.  At zero,
\eqref{eq:w04c-value-zero} applies.  Thus \(g^*\) has no zero in
\(\Re z\leq0\).

For \(z=x+iy\), \(y>0\), set
\(T(z)=\pi/\sin(\pi z)\).  A direct calculation gives
\begin{align*}
 \Re T(z)
 &=\frac{\pi\sin(\pi x)\cosh(\pi y)}
 {\sin^2(\pi x)+\sinh^2(\pi y)},
 \\
 \Im T(z)
 &=-\frac{\pi\cos(\pi x)\sinh(\pi y)}
 {\sin^2(\pi x)+\sinh^2(\pi y)}.
 \end{align*}
When \(0<x<1\), the real parts of both \(T(z)\) and \(F(z)\) are
positive.  When \(x=1\), the real part of \(T(z)\) vanishes and that of
\(F(z)\) remains positive.  Conjugation treats \(y<0\).  This proves the
zero-free assertion for \(\Re z\leq1\).

The same signs give the full confinement.  If
\(\sin(\pi x)\geq0\), then \(\Re D_\infty(x+iy)>0\).  If
\(\cos(\pi x)\geq0\), then
\(\Im D_\infty(x+iy)<0\).  The complement in \((1,\infty)\) of the
union of the intervals
\[
 (2n-1,2n-1/2),
 \qquad n\geq1,
\]
is covered by these two sign conditions.  Both inequalities are strict
on the boundary.  Every zero in the upper half plane must therefore
belong to one of the half-strips \(S_n^+\).

Fix \(n\geq1\), and put
\[
 a=2n-1,
 \qquad b=2n-\frac12.
\]
For small \(\eps>0\) and large \(Y\), take the positively
oriented boundary of
\[
 \{a<\Re z<b,\ 0<\Im z<Y\}
\]
with a quarter-disc of radius \(\eps\) removed at the pole
\(a\).  Traverse the bottom from \(a+\eps\) to \(b\), the right
side upward, the top to the left, the left side downward, and the
quarter-circle clockwise back to \(a+\eps\).

On the bottom, \eqref{eq:w04c-real-sign-result} gives
\(D_\infty<0\).  On the right side,
\(\cos(\pi b)=0\), so \(\Im D_\infty<0\).  On the left side,
\(\sin(\pi a)=0\), so \(\Re D_\infty>0\).  Uniformly on the top,
\[
 D_\infty(z)=\frac1{2z}+\mathcal O(Y^{-2}),
\]
and the argument variation tends to zero.  Finally,
\[
 D_\infty(z)=-\frac1{z-a}+\mathcal O(1)
\]
near the odd pole \(a\).

Using continuous arguments on the indicated half planes, the limiting
variations on the five pieces are
\[
 0,\qquad \frac\pi2,\qquad 0,\qquad \pi,
 \qquad \frac\pi2.
\]
Their sum is \(2\pi\).  The winding number\index[terms]{winding number} is integral and is unchanged
once \(\eps\) is small and \(Y\) is large.  The argument
principle gives exactly one zero counted with multiplicity in
\(S_n^+\).  The confinement shows that there are no others.  A zero
with total multiplicity one is simple.  The coefficients of
\(D_\infty\) are real, so conjugation supplies the lower zeros and proves
Theorem~\ref{thm:w04c-localization}.

Only the asymptotic location remains.  Write
\[
 \rho_n=x_n+iy_n,
 \qquad
 w_n=\rho_n+\frac12.
\]
The half-strip localization gives \(x_n=2n+\mathcal O(1)\).  At a zero,
\begin{equation}
 |\sin(\pi\rho_n)|=\frac\pi{|F(\rho_n)|}
 =2\pi|\rho_n|\left(1+\mathcal O(n^{-1})\right).
 \label{eq:w04c-sine-modulus}
\end{equation}
Since
\begin{equation}
 |\sin(\pi(x+iy))|^2
 =\sin^2(\pi x)+\sinh^2(\pi y),
 \label{eq:w04c-sine-identity}
\end{equation}
equations \eqref{eq:w04c-sine-modulus} and
\eqref{eq:w04c-sine-identity} imply
\[
 \sinh(\pi y_n)=2\pi|\rho_n|
 \left(1+\mathcal O(n^{-1})\right).
\]
Since \(|\rho_n|\geq x_n\asymp n\), this gives
\(y_n=\pi^{-1}\log n+\mathcal O(1)\), and hence
\begin{equation}
 y_n=\mathcal O(\log n),
 \qquad
 e^{-2\pi y_n}=\mathcal O(n^{-2}).
 \label{eq:w04c-height-preliminary}
\end{equation}

The two leading terms of \eqref{eq:w04c-F-first-expansion} give the
cancellation
\[
 \frac{2\pi}{F(\rho_n)}
 =4\pi w_n\left(1+\mathcal O(n^{-2})\right).
\]
The zero equation and
\(\sin(\pi\rho_n)=-\cos(\pi w_n)\) now yield
\begin{equation}
 e^{-i\pi w_n}
 =\frac{2\pi}{F(\rho_n)(1+e^{2i\pi w_n})}
 =4\pi w_n\left(1+\mathcal O(n^{-2})\right).
 \label{eq:w04c-exponential-equation}
\end{equation}
Here the second equality uses \eqref{eq:w04c-height-preliminary}.

The principal logarithm is available because \(w_n\) lies in the first
quadrant.  Its real part belongs to \((2n-1/2,2n)\), which fixes the
integer in the logarithm of \eqref{eq:w04c-exponential-equation}.  Thus
\begin{equation}
 w_n=2n+\frac{i}{\pi}
 \left[\operatorname{Log}(4\pi w_n)+\mathcal O(n^{-2})\right].
 \label{eq:w04c-log-equation}
\end{equation}
Set \(u_n=w_n-2n\).  The preliminary height bound gives
\(u_n=\mathcal O(L_n)\), and
\[
 \operatorname{Log}(4\pi w_n)
 =L_n+\frac{u_n}{2n}
 +\mathcal O\left(\frac{L_n^2}{n^2}\right).
\]
Substitution into \eqref{eq:w04c-log-equation} gives
\begin{align*}
 u_n
 &=\frac{iL_n}{\pi}+\frac{iu_n}{2\pi n}
 +\mathcal O\left(\frac{L_n^2}{n^2}\right)\\
 &=\frac{iL_n}{\pi}-\frac{L_n}{2\pi^2n}
 +\mathcal O\left(\frac{L_n^2}{n^2}\right).
 \end{align*}
Since \(\rho_n=w_n-1/2\), this proves
Proposition~\ref{thm:w04c-asymptotic}.
\end{proof}

\section{The orthorecursive resolvent trace}
\label{sec:trace_ortho_resolvent}

The first half of this chapter recast the orthorecursive weighted sums as a multiplicative Volterra
perturbation of the partial sums (Proposition~\ref{prop:ortho_volterra}) and evaluated them at the
integers (Proposition~\ref{prop:ortho_Ag_integer}). Section~\ref{sec:trace_ortho} determined the
complete zero set of the transform. This section reads the resolvent attached to the kernel as a
trace over those zeros, exact at finite depth and, for $y>1$, as a convergent infinite sum, and
extracts from it the first spectral pair of the orthorecursive partial sums.

Consider the orthorecursive kernel and its Volterra kernel\index[terms]{Volterra kernel}
\[
 g(t)=\frac{2}{1+t},\qquad k(y)=\frac{2y}{(1+y)^2},\qquad 0<t\le1\le y.
\]
Multiplicative convolution is written
\begin{equation}
 (f\star h)(y):=\int_1^y f(y/v)h(v)\,\frac{dv}{v}.
 \label{eq:w04d-convolution}
\end{equation}
The resolvent is the locally convergent Volterra series
\begin{equation}
 R:=\sum_{m\ge1}(-1)^{m-1}k^{\star m},
 \label{eq:w04d-resolvent-series}
\end{equation}
the unique locally bounded solution of
\begin{equation}
 R+R\star k=k.
 \label{eq:w04d-resolvent-equation}
\end{equation}
For a function on the multiplicative half line, put
\[
 \widetilde f(s):=\int_1^\infty f(y)y^{-s}\,\frac{dy}{y}.
\]
The arithmetic Mellin transform of the kernel is that of Lemma~\ref{lem:ortho_gstar},
\[
 g^*(z)=-z\int_0^1 g(t)t^{-z-1}\,dt=2zD_\infty(z),
\]
and Section~\ref{sec:trace_ortho} gives its reflection identity
\begin{align}
 D_\infty(z)&=\frac{\pi}{\sin(\pi z)}+F(z),\notag\\
 F(z)&=\int_0^\infty\frac{e^{-zu}}{e^u+1}\,du,\qquad \Re z>-1,
 \label{eq:w04d-F-integral}
\end{align}
together with the simple zero pairs $\rho_n,\overline{\rho_n}$ of
Theorem~\ref{thm:w04c-localization},
\begin{equation}
 2n-1<\Re\rho_n<2n-\frac12,\qquad \Im\rho_n>0,
 \label{eq:w04d-zero-strips}
\end{equation}
and their asymptotics of Proposition~\ref{thm:w04c-asymptotic}, with $L_n=\log(8\pi n)$,
\begin{equation}
 \rho_n=2n-\frac12-\frac{L_n}{2\pi^2n}+i\frac{L_n}{\pi}+\mathcal O\!\left(\frac{L_n^2}{n^2}\right).
 \label{eq:w04d-zero-asymptotic}
\end{equation}
There are no other zeros. The transformed resolvent and its regularized part are
\begin{equation}
 \widetilde R(s)=1-\frac1{g^*(-s)},\qquad E(s):=\widetilde R(s)-\frac1{2s},
 \label{eq:w04d-E-definition}
\end{equation}
the identity $\widetilde R(s)=1-1/g^*(-s)$ being the one of
Chapter~\ref{chap:discrete_volterra}. For a vertical line the arrows record its orientation,
\begin{align*}
 \int_{(a)}^{\uparrow}H(s)\,ds&:=\lim_{T\to\infty}\int_{a-iT}^{a+iT}H(s)\,ds,
 \\
 \int_{(a)}^{\downarrow}H(s)\,ds&:=\lim_{T\to\infty}\int_{a+iT}^{a-iT}H(s)\,ds.
 \end{align*}

The resolvent identity below is the one met for the discrete Volterra\index[terms]{discrete Volterra} operator in
Chapter~\ref{chap:discrete_volterra}, here carried by the orthorecursive kernel.
\begin{proposition}\label{thm:w04d-resolvent-inversion}
The series \eqref{eq:w04d-resolvent-series} converges locally uniformly
and solves \eqref{eq:w04d-resolvent-equation}.  For every \(c>0\), it
also converges in the weighted space associated with
\(y^{-c}\,dy/y\), and
\begin{equation}
 \widetilde R(s)=1-\frac1{g^*(-s)}
 \qquad (\Re s>0).
 \label{eq:w04d-resolvent-transform}
\end{equation}
On every line \(\Re s=c>0\), one has
\begin{equation}
 E(c+it)=\mathcal O_c((1+|t|)^{-2}).
 \label{eq:w04d-E-right-decay}
\end{equation}
For \(y>1\), Fourier inversion in the logarithmic variable gives
\begin{equation}
 R(y)=\frac12+\frac1{2\pi i}
 \int_{(c)}^{\uparrow}E(s)y^s\,ds.
 \label{eq:w04d-regularized-inversion}
\end{equation}
\end{proposition}

Moving the line of \eqref{eq:w04d-regularized-inversion} to the left crosses the poles one by
one, and the passage needs a bound on the displaced line that does not grow faster than the
number of poles crossed. The following lemma supplies it, with explicit constants.

\begin{lemma}\label{thm:w04d-vertical-control}
For every integer \(N\geq1\),
\begin{align}
 \int_{\mathbb R}|E(-2N+it)|\,dt
 \leq{}&\frac{31}{6}
 +\frac{40e^{-\pi}}{3(1-e^{-2\pi})}
 \left(2N+1+\frac1\pi\right)\notag\\
 &+\frac{2\pi}{3N}
 \left(\frac14+\frac5{6\sqrt3}\right)
 <9N.
 \label{eq:w04d-explicit-vertical-bound}
\end{align}
In particular, the dependence on the displaced line is at most linear
in \(N\).
\end{lemma}

The two ingredients now combine. Shifting the line past the first $N$ conjugate pairs and
discarding the displaced integral by the bound just proved gives the trace at finite depth.

\begin{theorem}
\label{thm:w04d-finite-trace}
For every integer \(N\geq1\) and every \(y>1\),
\begin{align*}
 R(y)
 &=2\Re\sum_{n=1}^{N}
 \frac{y^{-\rho_n}}{g^{*\prime}(\rho_n)}
 +\frac1{2\pi i}\int_{(-2N)}^{\uparrow}E(s)y^s\,ds
 \\
 &=2\Re\sum_{n=1}^{N}
 \frac{y^{-\rho_n}}{g^{*\prime}(\rho_n)}
 -\frac1{2\pi i}\int_{(-2N)}^{\downarrow}E(s)y^s\,ds.
 \end{align*}
The second display is the return-side convention for the finite contour.
The remainder satisfies
\begin{equation}
 \left|R(y)-2\Re\sum_{n=1}^{N}
 \frac{y^{-\rho_n}}{g^{*\prime}(\rho_n)}\right|
 \leq\frac{9N}{2\pi}y^{-2N}.
 \label{eq:w04d-finite-remainder}
\end{equation}
At every upper zero,
\begin{equation}
 \mathop{\rm Res}_{s=-\rho_n}\widetilde R(s)
 =\frac1{g^{*\prime}(\rho_n)}.
 \label{eq:w04d-residue-sign}
\end{equation}
The artificial pole of \(E\) at zero has residue \(-1/2\).  Its
contribution cancels the constant \(1/2\) in
\eqref{eq:w04d-regularized-inversion}.
\end{theorem}

Letting the depth grow asks for the size of the residues, which the reflection identity supplies,
and the trace then closes as a convergent sum for $y>1$.

\begin{theorem}
\label{thm:w04d-infinite-trace}
The residues satisfy
\begin{align}
 g^{*\prime}(\rho_n)
 &=-i\pi+\frac{i\pi/2-1}{\rho_n}
 +\mathcal O(n^{-2}),
 \label{eq:w04d-derivative-asymptotic}\\
 \frac1{g^{*\prime}(\rho_n)}
 &=\frac{i}{\pi}
 +\frac{i\pi/2-1}{\pi^2\rho_n}
 +\mathcal O(n^{-2}).
 \label{eq:w04d-residue-asymptotic}
\end{align}
In particular,
\begin{equation}
 \frac1{g^{*\prime}(\rho_n)}
 =\frac{i}{\pi}+\mathcal O(n^{-1}).
 \label{eq:w04d-residue-target}
\end{equation}
For every \(y>1\), the exact infinite trace is
\begin{equation}
 R(y)=2\Re\sum_{n\geq1}
 \frac{y^{-\rho_n}}{g^{*\prime}(\rho_n)}.
 \label{eq:w04d-infinite-resolvent-trace}
\end{equation}
The series converges absolutely.  Its convergence is locally uniform on
every half line \(y\geq1+\delta\), where \(\delta>0\).
\end{theorem}

The arithmetic estimates used below are established in Chapter~\ref{chap:ortho}: the
partial-sum bound behind \eqref{eq:ortho_partial_sum}, the pointwise bound
$c_N=\mathcal O(N^{-2})$ of Theorem~\ref{thm:ortho_bootstrap}, and the integer evaluation
$A_g(N)=-\tfrac{2N}{2N+1}c_N$ of Proposition~\ref{prop:ortho_Ag_integer}.
\begin{theorem}
\label{thm:w04d-first-arithmetic-pair}
Let \(a_n=c_{n-1}\) be the orthorecursive sequence and put
\begin{equation}
 A(x):=\sum_{n\leq x}a_n,
 \qquad
 A_g(x):=\sum_{n\leq x}\frac{2x\,a_n}{x+n}.
 \label{eq:w04d-arithmetic-sums}
\end{equation}
The established arithmetic estimates imply the continuous bound
\begin{equation}
 A_g(x)=\mathcal O(x^{-2}).
 \label{eq:w04d-Ag-continuous-bound}
\end{equation}
The integral
\begin{equation}
 C_1:=\frac1{g^{*\prime}(\rho_1)}
 \int_1^\infty A_g(u)u^{\rho_1}\,\frac{du}{u}
 \label{eq:w04d-C1-definition}
\end{equation}
converges absolutely, and
\begin{equation}
 A(x)=-2\Re\left(C_1x^{-\rho_1}\right)
 +\mathcal O(x^{-2}\log x).
 \label{eq:w04d-first-pair-asymptotic}
\end{equation}
This assertion does not require \(C_1\neq0\).  The residue of the
resolvent is only one factor in \(C_1\).  The other factor is the global
Mellin datum carried by \(A_g\).
\end{theorem}

\begin{proof}[Proof of Proposition~\ref{thm:w04d-resolvent-inversion}]

Pass to the logarithmic coordinate \(u=\log y\), and write
\[
 K(u):=k(e^u)=\frac{2e^u}{(1+e^u)^2}
 =\frac1{2\cosh^2(u/2)}.
\]
Additive convolution on \([0,\infty)\) is the logarithmic form of
\eqref{eq:w04d-convolution}.  Since \(0<K(u)\leq1/2\), induction on
the convolution order gives
\begin{equation}
 |K^{*m}(u)|\leq
 \left(\frac12\right)^m\frac{u^{m-1}}{(m-1)!}.
 \label{eq:w04d-neumann-local-bound}
\end{equation}
The Neumann series\index[terms]{Neumann series}\index[names]{Neumann, C.} is therefore locally uniform.  The usual telescoping
calculation gives \(R+R*K=K\).  If two locally bounded functions solve
this equation, iteration of their difference and
\eqref{eq:w04d-neumann-local-bound} prove uniqueness.

For \(c>0\), put
\[
 q(c):=\int_0^\infty K(u)e^{-cu}\,du.
\]
The total mass of \(K\) is one, and the inequality
\(e^{-cu}<1\) is strict for \(u>0\).  Hence \(q(c)<1\).  The Neumann
series converges in the corresponding weighted \(L^1\) space and gives
\begin{equation}
 \widetilde R(s)
 =\frac{\widetilde k(s)}{1+\widetilde k(s)}.
 \label{eq:w04d-neumann-transform}
\end{equation}
The identity
\[
 \widetilde k(s)=g^*(-s)-1
\]
now proves \eqref{eq:w04d-resolvent-transform}.

The derivatives of \(K\) are integrable after multiplication by
\(e^{-cu}\).  Two integrations by parts, together with
\(K(0)=1/2\) and \(K'(0)=0\), give
\begin{equation}
 \widetilde k(s)=\frac1{2s}+\mathcal O_c(|s|^{-2})
 \qquad (\Re s=c,\ |\Im s|\to\infty).
 \label{eq:w04d-k-transform-decay}
\end{equation}
Substitution in \eqref{eq:w04d-neumann-transform} proves
\eqref{eq:w04d-E-right-decay}.  The function
\(e^{-cu}(R(e^u)-1/2)\) belongs to \(L^1(0,\infty)\).  Its Fourier
transform is \(E(c+it)\), which is also integrable by
\eqref{eq:w04d-E-right-decay}. The locally uniform Neumann series makes
$R$ continuous, so the weighted function is continuous at every $u>0$.
Fourier inversion \cite[\S7.2, (7.16), pp.~218--219]{Folland1992} therefore
proves \eqref{eq:w04d-regularized-inversion} and completes the proof of
Proposition~\ref{thm:w04d-resolvent-inversion}.

\end{proof}

\begin{proof}[Proof of Lemma~\ref{thm:w04d-vertical-control}]
The dependence on the displaced lines requires a separate argument.  Put
\[
 z=2N-it,
 \qquad
 T(z):=\frac{\pi}{\sin(\pi z)},
 \qquad
 p(u):=\frac1{4\cosh^2(u/2)},
\]
and define
\[
 Q(z):=\int_0^\infty e^{-zu}p(u)\,du,
 \qquad H(z):=2zT(z).
\]
Since \(p=-h'\), where \(h(u)=(e^u+1)^{-1}\), one integration by
parts in \eqref{eq:w04d-F-integral} gives
\begin{equation}
 2zF(z)=1-2Q(z),
 \qquad
 g^*(z)=1+H(z)-2Q(z).
 \label{eq:w04d-gstar-QH}
\end{equation}

On the line \(\Re z=2N\), the reflection term \(T(z)\) is purely
imaginary.  Its imaginary part and that of \(F(z)\) have the same sign
when \(t\neq0\).  Indeed,
\begin{equation}
 T(2N-it)=\frac{i\pi}{\sinh(\pi t)},
 \label{eq:w04d-T-even-line}
\end{equation}
and the strict quadrant lemma for \(F\), followed by conjugation, gives
the same sign for its imaginary part together with \(\Re F(z)>0\).
Thus \(|T(z)+F(z)|\geq|F(z)|\).  Also
\[
 |Q(z)|\leq\int_0^\infty e^{-2Nu}p(u)\,du
 \leq\frac1{8N}.
\]
It follows that
\begin{equation}
 |g^*(z)|=2|z|\,|T(z)+F(z)|
 \geq2|z|\,|F(z)|
 =|1-2Q(z)|\geq\frac34.
 \label{eq:w04d-gstar-lower-bound}
\end{equation}
At \(t=0\), the reciprocal of \(g^*\) is understood by continuity at
the pole of \(g^*\).

A second pair of integrations by parts gives
\begin{equation}
 Q(z)=\frac1{4z}+\frac{J(z)}{z^2},
 \qquad
 J(z):=\int_0^\infty e^{-zu}p''(u)\,du.
 \label{eq:w04d-Q-second-expansion}
\end{equation}
The substitution \(v=\tanh(u/2)\) yields
\begin{align*}
 \|p''\|_1
 &=\frac14\int_0^1|1-3v^2|\,dv
 =\frac1{3\sqrt3}.
 \end{align*}
Thus \(|J(z)|\leq1/(3\sqrt3)\).

Write \(\Delta(z)=H(z)-2Q(z)\).  From
\eqref{eq:w04d-E-definition} and \eqref{eq:w04d-gstar-QH},
\begin{equation}
 E(-z)=\frac{\Delta(z)}{1+\Delta(z)}+\frac1{2z}.
 \label{eq:w04d-E-delta}
\end{equation}
For \(|t|\leq1\), equations \eqref{eq:w04d-gstar-lower-bound} and
\eqref{eq:w04d-E-delta} give
\begin{equation}
 |E(-z)|\leq1+\frac43+\frac14=\frac{31}{12}.
 \label{eq:w04d-E-central-bound}
\end{equation}
For \(|t|\geq1\), substitute
\eqref{eq:w04d-Q-second-expansion} into
\eqref{eq:w04d-E-delta}.  The term \(-1/(2z)\) in \(\Delta\)
cancels the last term in \eqref{eq:w04d-E-delta}.  Using
\(|z|\geq2\) and \eqref{eq:w04d-gstar-lower-bound} gives
\begin{equation}
 |E(-z)|
 \leq\frac53|H(z)|
 +\frac4{3|z|^2}
 \left(\frac14+\frac5{6\sqrt3}\right).
 \label{eq:w04d-E-tail-bound}
\end{equation}
Formula \eqref{eq:w04d-T-even-line} gives
\[
 |H(z)|=\frac{2\pi|z|}{|\sinh(\pi t)|}.
\]
For \(t\geq1\), use
\[
 |z|\leq2N+t,
 \qquad
 \frac1{\sinh(\pi t)}
 \leq\frac{2e^{-\pi t}}{1-e^{-2\pi}}.
\]
Integrating \eqref{eq:w04d-E-central-bound} and
\eqref{eq:w04d-E-tail-bound} on the two symmetric tails proves the
first inequality in \eqref{eq:w04d-explicit-vertical-bound}.  After
division by \(N\), every nonconstant contribution is largest at
\(N=1\).  Substitution in the displayed elementary constant gives a
value less than nine.  This proves the second inequality and
Lemma~\ref{thm:w04d-vertical-control}.

\end{proof}

\begin{proof}[Proof of Theorem~\ref{thm:w04d-finite-trace}]
It remains to move the inversion line.  Fix \(c>0\).  For
\(0\leq\Re s\leq c\), the integrations by parts leading to
\eqref{eq:w04d-k-transform-decay} are uniform and give
\(E(s)=\mathcal O_{N,c}(|s|^{-2})\).  For
\(-2N\leq\Re s\leq0\), put \(z=-s\).  Repeated integration by parts
in \eqref{eq:w04d-F-integral}, together with the exponentially small
reflection term on horizontal segments, gives uniformly
\[
 F(z)=\frac1{2z}-\frac1{4z^2}
 +\mathcal O_{N,c}(|z|^{-3}),
 \qquad
 \frac{\pi}{\sin(\pi z)}=\mathcal O_{N,c}(e^{-\pi|\Im z|}).
\]
Hence \(g^*(-s)=1+1/(2s)+\mathcal O_{N,c}(|s|^{-2})\), which again
gives
\[
 E(\sigma+iT)=\mathcal O_{N,c}(T^{-2})
\]
uniformly in \(-2N\leq\sigma\leq c\).  Hence the horizontal sides of
the contour vanish as \(T\to\infty\).

The poles crossed by the contour are zero and the points
\(-\rho_n,-\overline{\rho_n}\) with \(1\leq n\leq N\).  The poles of
\(g^*\) at the positive integers are zeros of its reciprocal and create
no pole of \(E\).  If \(s=-\rho+h\), then
\[
 g^*(-s)=-h\,g^{*\prime}(\rho)+\mathcal O(h^2).
\]
Consequently,
\[
 -\frac1{g^*(-s)}
 =\frac1{h\,g^{*\prime}(\rho)}+\mathcal O(1),
\]
which proves the positive sign in \eqref{eq:w04d-residue-sign}.
The function \(\widetilde R\) is regular at zero because \(g^*(0)=2\).
Therefore \(E\) has residue \(-1/2\) there.

The residue theorem\index[terms]{residue theorem} applied to the upward lines now gives
\[
 \frac1{2\pi i}\int_{(c)}^{\uparrow}E(s)y^s\,ds
 =-\frac12
 +2\Re\sum_{n=1}^N
 \frac{y^{-\rho_n}}{g^{*\prime}(\rho_n)}
 +\frac1{2\pi i}\int_{(-2N)}^{\uparrow}E(s)y^s\,ds.
\]
The first term on the right cancels the first term in
\eqref{eq:w04d-regularized-inversion}.  Reversing the orientation of
the last integral proves both forms of the finite trace.  Lemma
\ref{thm:w04d-vertical-control} gives
\eqref{eq:w04d-finite-remainder}.  This proves
Theorem~\ref{thm:w04d-finite-trace}.

\end{proof}

\begin{proof}[Proof of Theorem~\ref{thm:w04d-infinite-trace}]
The residue asymptotic comes from the same reflection identity.  At a
zero \(\rho=\rho_n\), one has \(T(\rho)=-F(\rho)\) and
\begin{equation}
 g^{*\prime}(\rho)
 =2\rho\left[\pi\cot(\pi\rho)F(\rho)+F'(\rho)\right].
 \label{eq:w04d-gprime-reflection}
\end{equation}
Repeated integration by parts in \eqref{eq:w04d-F-integral} gives,
uniformly at these zeros,
\begin{align*}
 F(\rho)
 &=\frac1{2\rho}-\frac1{4\rho^2}
 +\mathcal O(|\rho|^{-3}),
 \\
 F'(\rho)
 &=-\frac1{2\rho^2}+\frac1{2\rho^3}
 +\mathcal O(|\rho|^{-4}).
 \end{align*}
The zero asymptotic \eqref{eq:w04d-zero-asymptotic} implies
\[
 e^{-2\pi\Im\rho_n}=\mathcal O(n^{-2}),
 \qquad
 \cot(\pi\rho_n)=-i+\mathcal O(n^{-2}).
\]
Substitution into \eqref{eq:w04d-gprime-reflection} proves
\eqref{eq:w04d-derivative-asymptotic}.  Inverting that expansion proves
\eqref{eq:w04d-residue-asymptotic} and
\eqref{eq:w04d-residue-target}.

The residues are bounded.  Equations \eqref{eq:w04d-zero-strips} and
\eqref{eq:w04d-residue-target} show that, for \(y>1\),
\[
 \sum_{n\geq1}\left|
 \frac{y^{-\rho_n}}{g^{*\prime}(\rho_n)}\right|
 \ll\sum_{n\geq1}y^{-(2n-1)}<\infty.
\]
The same majorant is uniform for \(y\geq1+\delta\).  Finally,
\eqref{eq:w04d-finite-remainder} tends to zero as \(N\to\infty\).
This proves the exact trace and completes the proof of
Theorem~\ref{thm:w04d-infinite-trace}.

\end{proof}

\begin{remarkx}[The boundary of the trace]
The trace already has a canonical Abel boundary value. Indeed, the locally
uniform Neumann series of Proposition~\ref{thm:w04d-resolvent-inversion}
shows that $R$ is continuous at $y=1$. Every convolution power of order
at least two vanishes there, and therefore
\[
 R(1)=k(1)=\frac12.
\]
Together with \eqref{eq:w04d-infinite-resolvent-trace}, this gives
\begin{equation}
 \lim_{y\downarrow1}
 2\Re\sum_{n\ge1}\frac{y^{-\rho_n}}{g^{*\prime}(\rho_n)}=\frac12.
 \label{eq:w04d-Abel-boundary}
\end{equation}
The order of these operations matters. Equations
\eqref{eq:w04d-zero-asymptotic} and \eqref{eq:w04d-residue-asymptotic} give
\[
 \frac1{\rho_n}=\frac1{2n}+\mathcal O\!\left(\frac{\log n}{n^2}\right),
 \qquad
 2\Re\frac1{g^{*\prime}(\rho_n)}
 =-\frac1{\pi^2n}+\mathcal O\!\left(\frac{\log n}{n^2}\right).
\]
Consequently the conjugate-pair partial sums obtained by putting $y=1$
termwise satisfy
\[
 2\Re\sum_{n\le N}\frac1{g^{*\prime}(\rho_n)}
 =-\frac1{\pi^2}\log N+\mathcal O(1).
\]
Thus the ordinary paired boundary series diverges logarithmically, whereas
the Abel procedure in \eqref{eq:w04d-Abel-boundary} has the value $1/2$.
\end{remarkx}

\begin{proof}[Proof of Theorem~\ref{thm:w04d-first-arithmetic-pair}]
Only the arithmetic transfer remains.  The established bounds for the
orthorecursive sequence give, for every sufficiently small
\(\eps>0\),
\[
 A(N)=\mathcal O_\eps
 \left(N^{-\Re\rho_1+\eps}\right),
 \qquad
 c_N=\mathcal O(N^{-2}).
\]
Choose \(\eps<\Re\rho_1-1\).  Then
\begin{equation}
 \sum_{N\geq1}|A(N)|<\infty.
 \label{eq:w04d-A-summable}
\end{equation}
For \(N<x<N+1\), differentiation of
\eqref{eq:w04d-arithmetic-sums} gives
\[
 A_g'(x)=\sum_{n\leq N}a_n b_x(n),
 \qquad
 b_x(t):=\frac{2t}{(x+t)^2}.
\]
Abel summation gives
\begin{equation}
 A_g'(x)=A(N)b_x(N)
 -\sum_{n=1}^{N-1}A(n)
 \left(b_x(n+1)-b_x(n)\right).
 \label{eq:w04d-Ag-Abel}
\end{equation}
For \(1\leq t\leq N\leq x\),
\[
 0\leq b_x'(t)=\frac{2(x-t)}{(x+t)^3}
 \leq\frac2{x^2},
 \qquad b_x(N)=\mathcal O(x^{-1}).
\]
Equations \eqref{eq:w04d-A-summable} and
\eqref{eq:w04d-Ag-Abel} imply \(A_g'(x)=\mathcal O(x^{-2})\).
At the integers, the exact identity of Proposition~\ref{prop:ortho_Ag_integer} gives
\[
 A_g(N)=-\frac{2N}{2N+1}c_N=\mathcal O(N^{-2}).
\]
Integration over each unit interval proves
\eqref{eq:w04d-Ag-continuous-bound}.

Since \(\Re\rho_1<3/2<2\), equation
\eqref{eq:w04d-Ag-continuous-bound} proves the absolute convergence in
\eqref{eq:w04d-C1-definition}.  The finite trace with \(N=1\) can be
written
\[
 R(y)=2\Re\left(
 \frac{y^{-\rho_1}}{g^{*\prime}(\rho_1)}\right)+Q_1(y),
 \qquad Q_1(y)=\mathcal O(y^{-2}).
\]
The exact Volterra inversion, the inverse of Proposition~\ref{prop:ortho_volterra}, is
\begin{equation}
 A=A_g-A_g\star R.
 \label{eq:w04d-exact-arithmetic-inversion}
\end{equation}
After the change of variable \(u=x/y\), the first mode in the
convolution equals
\begin{equation}
 2\Re\left[
 \frac{x^{-\rho_1}}{g^{*\prime}(\rho_1)}
 \int_1^xA_g(u)u^{\rho_1}\,\frac{du}{u}
 \right].
 \label{eq:w04d-first-mode-convolution}
\end{equation}
The tail from \(x\) to infinity is \(\mathcal O(x^{-2})\) after
multiplication by \(x^{-\rho_1}\).  The remaining convolution satisfies
\begin{equation}
 |(A_g\star Q_1)(x)|
 \ll x^{-2}\int_1^x\frac{du}{u}
 \ll x^{-2}\log x.
 \label{eq:w04d-Q1-convolution-bound}
\end{equation}
Combining \eqref{eq:w04d-exact-arithmetic-inversion},
\eqref{eq:w04d-first-mode-convolution}, and
\eqref{eq:w04d-Q1-convolution-bound} proves
\eqref{eq:w04d-first-pair-asymptotic}.  The minus sign comes from the
subtraction in \eqref{eq:w04d-exact-arithmetic-inversion}.  This proves
Theorem~\ref{thm:w04d-first-arithmetic-pair}.
\end{proof}

\begin{openproblem}[Arithmetic amplitudes and the discrete spectral lift]
\label{op:w04d-boundary-and-amplitudes}
Complete the arithmetic spectral picture in the following three steps.
\begin{enumerate}[label=\textup{(\roman*)}]
\item Prove that \(C_1\) in \eqref{eq:w04d-C1-definition} is nonzero. In
view of \eqref{eq:w04d-first-pair-asymptotic}, this would show that the
exponent \(\alpha_1\) is attained by the orthorecursive partial sums.
\item Put
\[
 E(x):=A(x)+2\Re\bigl(C_1x^{-\rho_1}\bigr).
\]
Prove, for some \(\sigma>\alpha_1\), the first-difference estimate
\[
 E(n+1)-E(n)=\mathcal O\bigl(n^{-1-\sigma}\bigr).
\]
Since \(c_n=A(n+1)-A(n)\), this would prove
Conjecture~\ref{conj:ortho_spectral}, with \(\kappa=\rho_1C_1\).
\item Determine whether, for every fixed spectral rank \(M\), the continuous
trace admits a discrete finite-rank expansion built from the conjugate
Gamma modes
\[
 \frac{\Gamma(n+1-\rho_j)}{\Gamma(n+1)},\qquad 1\le j\le M,
\]
with globally normalized arithmetic amplitudes and a quantified remainder.
No infinite discrete spectral sum is asserted here.
\end{enumerate}
\end{openproblem}

\begin{proofstatus}{Theorems~\ref{thm:w04d-finite-trace}--\ref{thm:w04d-first-arithmetic-pair}
prove the finite trace, the residue asymptotics and infinite trace for $y>1$, and the first
partial-sum pair. Equation~\eqref{eq:w04d-Abel-boundary} and the preceding boundary discussion
prove Abel summability at $y=1$ and the logarithmic divergence of the ordinary paired series.
These results supply Step~9 before
Definition~\ref{def:resolvent_regular}, the discussion after
Conjecture~\ref{conj:ortho_spectral}, and the epilogue. Open
Problem~\ref{op:w04d-boundary-and-amplitudes} remains because $C_1$ may cancel, a coefficient
trace needs first differences, and no finite-rank discrete Gamma expansion is proved. No theorem
assumes these three clauses.}
\end{proofstatus}
Here the body of the volume ends, and it ends on the kernel that carries no hypothesis. The decay
left open by Kalmynin\index[names]{Kalmynin, A. B.} and Kosenko\index[names]{Kosenko, P. R.} is
settled from above, the exponent bounded by the leftmost zero of the transform of a function of
good variation hidden in the problem, and the resolvent is written as a sum over all of its
zeros. No statement
of the chapter rests on a hypothesis about $\zeta$. What is not settled is which branch of the
alternative holds. The leading amplitude may cancel, and Open
Problem~\ref{op:w04d-boundary-and-amplitudes} carries that question together with the
first-difference control a coefficient trace would need. The method is
the one that turns the Ingham averages into the Riemann hypothesis. Here it runs to the end.

\backmatter

\chapter*{Epilogue: from Volume I to Volume II}

\rafepigraph{Le savant doit ordonner\,; on fait la science avec des faits comme une maison avec des pierres\,; mais une accumulation de faits n'est pas plus une science qu'un tas de pierres n'est une maison.}{The scientist must set things in order. Science is built with facts as a house is built with stones, but an accumulation of facts is no more a science than a heap of stones is a house.}{Henri Poincaré\index[names]{Poincar\'e, H.}, \emph{La Science et l'Hypothèse} (1902)~\cite{Poincare1902}}

\label{chap:epilogue}
\addcontentsline{toc}{chapter}{Epilogue: from Volume I to Volume II}
\markboth{FROM VOLUME I TO VOLUME II}{FROM VOLUME I TO VOLUME II}

Volume I built one equivalence and studied one invariant. The
equivalence is $\alpha(\Phi)=\tfrac12\iff\text{RH}$, established on the
arithmetic side through M\"obius inversion\index[terms]{M\"obius inversion} and the Littlewood bound and,
on the analytic side, through the Volterra resolvent and the decay of
its cumulative function. The boundary of the transform came first, and from the non-vanishing of
$\zeta$ on the line $\Re s=1$ alone it returned Ingham's theorem and the prime number theorem, a
limit without a rate. The invariant is the regularity index
$\alpha(G)$, read across a gallery of kernels, continuous, arithmetically
structured, and diophantinely rough. Seven of those kernels have a transform without a single
zero and six of them carry a proved index, so the invariant is at its clearest where no
analytic route exists at all. The last chapters added a further
element. Under a change of the evaluation coordinates, a gauge $f$, the
analytic index of the gauged Ingham transform\index[terms]{Ingham transform} is $\tfrac12$ for a family of
exponential gauges, its zeros lying on the critical line without any hypothesis
on $\zeta$. This closing
chapter records where that observation leads, and it names the threads
that the next volume carries forward.

\section*{A bridge to finite fields}
\label{sec:epilogue_bridge}

The arithmetic Mellin transform of the Ingham function under the gauge
$f(x)=q^{x}+1$, computed in
Chapter~\ref{chap:gauge_ingham}, is
\[\Phi^*_f(z)=\frac{q^{2z}-2q^{z}+q}{q-q^{z}},\qquad q\ge2\ \text{integer}.
\]
A local zeta function\index[terms]{local zeta function} counts. For a smooth projective curve $C$ over the finite field\index[terms]{finite field} $\mathbb F_q$,
let $N_m$ be the number of points of $C$ with coordinates in the extension $\mathbb F_{q^m}$. The
local zeta function of $C$ is the generating series
\begin{equation}\label{eq:epi_zeta_def}
Z(C/\mathbb F_q,T)=\exp\Bigl(\sum_{m\ge1}\frac{N_m}{m}\,T^{m}\Bigr),
\end{equation}
and the theorem of Weil\index[names]{Weil, A.} says that it is a rational function of $T$,
\begin{equation}\label{eq:epi_weil_form}
Z(C/\mathbb F_q,T)=\frac{P(T)}{(1-T)(1-qT)},
\end{equation}
where $P$ is a polynomial with integer coefficients, of degree twice the genus of $C$, with
$P(0)=1$, whose reciprocal roots all have modulus $\sqrt q$. That last statement is the Riemann
hypothesis for the curve. The polynomial $P$ is the characteristic polynomial of the Frobenius
endomorphism, the map raising coordinates to the power $q$, written in the reciprocal variable.
For an elliptic curve, where the genus is one, $P(T)=1-aT+qT^{2}$ with $a=q+1-N_1$, the
Frobenius trace\index[terms]{Frobenius trace}, which measures the deviation of the point count from its expected
value $q+1$.

Written in the variable $T=q^{-z}$, the transform $\Phi^*_f$ factors through a function of
exactly the shape \eqref{eq:epi_weil_form}. Take $P(T)=1-2T+qT^{2}$, the Frobenius polynomial
of trace $a=2$ in the form given by Hindry\index[names]{Hindry, M.} \cite{Hindry2012}, and write
\begin{equation}\label{eq:epi_hasse_zeta}
Z(C/\mathbb F_q,T)=\frac{qT^{2}-2T+1}{(1-T)(1-qT)}
\end{equation}
for the corresponding right side of \eqref{eq:epi_weil_form}. A direct substitution gives the
identity
\begin{equation}\label{eq:epi_identity}
\Phi^*_f(z)=\Big(1-\frac1T\Big)\,Z(C/\mathbb F_q,T),\qquad T=q^{-z}.
\end{equation}
When $q$ is a power of a prime and an elliptic curve $C$ of trace $2$ exists over
$\mathbb F_q$, the notation is literal and \eqref{eq:epi_hasse_zeta} is the zeta function of
that curve. For a general integer $q\ge2$ the same expression is a formal factor of
Hasse\index[names]{Hasse, H.}-Weil type, with no curve behind it, since a finite field $\mathbb F_q$ exists only
when $q$ is a prime power, and the letter $C$ in the notation is then kept for the form alone.
The Ingham operator, an object of Tauberian\index[names]{Tauber, A.} analysis over $\Z$, and this
rational factor meet in one identity, valid for every integer $q\ge2$.

The meeting is not only formal. The finite-field zeta function is
recovered exactly from the floor structure of the Ingham kernel, by
Theorem~\ref{thm:gauge_floor_limit}. For integer $q\ge2$ and real $T>1$,
\begin{equation}\label{eq:epi_floor_limit}
\lim_{n\to\infty}\sum_{k=1}^{n}\frac{1}{(qT)^{\,n-k}}
\Big\lfloor\frac{q^{n}+1}{q^{k}+1}\Big\rfloor
=\frac{qT^{2}-2T+1}{(1-T)(1-qT)}
=Z(C/\mathbb F_q,T).
\end{equation}
The left side is a weighted sum of the same integer quotients
$\lfloor(q^{n}+1)/(q^{k}+1)\rfloor$ that govern the gauged Ingham
recurrence, and its limit is the local zeta function\index[terms]{local zeta function} on the right.

The location of the zeros is algebraic. The zeros of $\Phi^*_f$ are the zeros of the
numerator $qT^{2}-2T+1$, a quadratic whose two roots have product $1/q$ and are
complex conjugate for $q\ge2$, so each has modulus $q^{-1/2}$. In the variable $z$
this is $\Re z=\tfrac12$. This is a fact about the polynomial, holding for every
integer $q\ge2$, independently of any geometric realization.

When $q$ is a power of a prime, the modulus $\sqrt q$ of the reciprocal roots is the
Riemann hypothesis for a curve over $\mathbb F_q$, proved by Hasse\index[names]{Hasse, H.} for elliptic
curves and by Weil\index[names]{Weil, A.} for curves of any genus \cite{Hindry2012}, and placed by Deligne\index[names]{Deligne, P.}
in the setting of the Weil conjectures \cite{Deligne1974}. Whether a curve of trace
$a=2$ exists over $\mathbb F_q$ is not settled by the Hasse\index[names]{Hasse, H.} bound alone. It is
governed by the classification of Waterhouse\index[names]{Waterhouse, W. C.} \cite{Waterhouse1969}. For an odd prime
power $q$ the trace $a=2$ is admissible, since $\gcd(2,q)=1$ places it in the
ordinary case. For even prime powers only the cases verified in that classification
occur, and $q=8$, for instance, admits no curve of trace $2$. The location of the
zeros of $\Phi^*_f$ on the critical line is therefore an algebraic fact for every
integer $q\ge2$, while its reading as the Riemann hypothesis of a curve is available
only over the admissible prime powers. The identification of the regularity index
$\alpha_f(\Phi)$ with one half is proved in the companion volume and recorded here as
Theorem~\ref{numobs:gauge_critical}. The equilibrium of
Conjecture~\ref{conj:ingham_equilibrium} is therefore the sole remaining step and is equivalent
to the classical Riemann hypothesis.

\begin{remarkx}[Why the two objects meet]\label{rem:epi_numerics}
None of the three statements above is numerical. The floor limit
\eqref{eq:epi_floor_limit} is Theorem~\ref{thm:gauge_floor_limit}, proved from the
dyadic identity of Lemma~\ref{lem:dyadic_floor}, and its value is the generating
function $\mathcal H$ of the limiting gauged kernel read at $w=1/T$, so that
$\mathcal H(1/T)=Z(C/\mathbb F_q,T)$. The identity \eqref{eq:epi_identity} is then
$\Phi^*_f(z)=(1-w)\mathcal H(w)$ rewritten with $w=q^{z}=1/T$, and the modulus of the
numerator roots is the computation of \eqref{eq:gauge_numerator}. The local zeta
function of a curve over $\mathbb F_q$ with Frobenius trace $2$ is the generating
function of the limiting kernel of the Ingham operator under the gauge
$f(x)=q^{x}+1$, and the Frobenius polynomial $qT^{2}-2T+1$ and the numerator
$w^{2}-2w+q$ of $\mathcal H$ are one object read in two reciprocal variables. The
numerical agreement to forty digits at test points records a check and not a
definition.
\end{remarkx}

\begin{remarkx}[On the geometric object]\label{rem:epi_geometry}
The right side of \eqref{eq:epi_hasse_zeta} is written as the zeta
function of an elliptic curve over $\mathbb F_q$ with trace $a=2$. Which
curve realizes this trace for each $q$, and how the family varies with
$q$, is a question of the arithmetic of elliptic curves over finite
fields rather than of the present analytic theory. Pinning the geometric
object precisely, and reading the gauge parameter as a
geometric datum, is left open here. Reading a Riemann hypothesis through a geometry
over a base smaller than $\Z$ is the program recalled by Connes\index[names]{Connes, A.} \cite{Connes2026}.
\end{remarkx}

\section*{Toward Volume II}
\label{sec:epilogue_volII}

Volume I connected two languages. On one side stands the analysis of the
Ingham operator, the Volterra resolvent, the Mellin transform, the
Tauberian passage from averaged sums to partial sums. On the other
stands the arithmetic of $\zeta$, of the M\"obius function\index[terms]{M\"obius function}, of the
distribution of the primes. The regularity index is the quantity through
which the two are read as one. The identity \eqref{eq:epi_identity} adds
a third language. The same operator, in exponential coordinates, carries
the zeta function of a curve over a finite field, an object of algebraic
geometry whose Riemann hypothesis is a theorem.

The identity of the epilogue shows that one operator carries two zeta
functions of different origin, the transcendental $\zeta$ and the finite-field
$Z(C/\mathbb F_q,T)$, and that the change of arithmetic coordinates is what brings
them into the same frame. Building that into a bridge, from the side where the Riemann
hypothesis is a theorem to the side where it is not, is the programme of the second volume.

One comparison is already in view. The resolvent trace of
Chapter~\ref{chap:trace} realizes the resolvent of the orthorecursive kernel as a sum over the
zeros of its transform, exact for $y>1$ and Abel-summable at $y=1$. The non-vanishing of its
leading arithmetic amplitude and its discrete spectral lift are left open in Open
Problem~\ref{op:w04d-boundary-and-amplitudes}. On the geometric side the
identity \eqref{eq:epi_identity} writes the gauged Ingham transform\index[terms]{Ingham transform} through a rational
factor of Hasse\index[names]{Hasse, H.}-Weil type, whose Frobenius trace\index[terms]{Frobenius trace} is the arithmetic counterpart of that
analytic sum over zeros. Setting the two side by side, in the manner of the trace formulas of
Selberg\index[names]{Selberg, A.} and of Connes\index[names]{Connes, A.}, is a comparison this
volume leaves open.

What the volume isolates is not an equation but an object. The triangular system, the resolvent and
the Mellin transform were all available before, and the book uses them. What was not available was the
single characteristic that reads a kernel through the whole family of power forcings at once, a
maximal transparent range together with an absorption law beyond it, compared with rather than
defined by the zeros of a transform. The equations were already there, the common object was not.

The reserved position of Volume I stands on its own. The regularity index is
defined by the discrete equation, its value at the Ingham kernel is equivalent to the
Riemann hypothesis, and the gauge deformations show that the threshold moves under a
change of coordinates while the question of an intrinsic reading remains open. That
question is where the next volume begins.

\vspace{1.2\baselineskip}
\rafepigraph{Combien de temps faudra-t-il encore pour que notre pierre de Rosette, à nous autres arithméticiens, rencontre son Champollion\,?}{How much longer will it take before our Rosetta stone, we arithmeticians, meets its Champollion?}{André Weil, \emph{De la métaphysique aux mathématiques} (1960)~\cite{WeilMetaphysique1960}}

\appendix
\makeatletter\@mainmattertrue\makeatother
\addcontentsline{toc}{part}{Appendices: a gallery of arithmetic kernels}

\chapter*{The gallery and the shape of its proofs}
\addcontentsline{toc}{chapter}{The gallery and the shape of its proofs}
\markboth{THE GALLERY AND THE SHAPE OF ITS PROOFS}{THE GALLERY AND THE SHAPE OF ITS PROOFS}

The seventeen gallery entries and the research dossier that follow are the material of the
theory. Each studies a fixed kernel through its defining equation\index[terms]{defining equation}
for power forcings, determining its regularity index\index[terms]{regularity index} or isolating
the estimates still needed to do so. The proofs exhibit the mechanism of the operator, and
three reasons make this constructive approach necessary.

The first is that for seven of these kernels there is nothing to invoke. The shifted rational kernel
of Appendix~\ref{app:H}, the quadratic rational kernel of Appendix~\ref{app:K}, the two logarithmic
kernels of Appendices~\ref{app:L} and~\ref{app:N}, the fractional part\index[terms]{fractional part} kernel of Appendix~\ref{app:O},
the binomial harmonic kernel of Appendix~\ref{app:Q} and the greatest common divisor kernel of
Appendix~\ref{app:J} have transforms without zeros, so no
analytic index\index[terms]{analytic index} exists for them. Any route that locates the index among the zeros of a transform, a
contour argument\index[terms]{contour integral}, a Tauberian transfer\index[terms]{Tauberian}, an inversion of Paley--Wiener\index[terms]{Paley--Wiener theorem} type, has nothing to
locate. For these kernels the index is arithmetic, and the equation is the only access to it. Six of the
seven carry a proved index, the greatest common divisor kernel of Appendix~\ref{app:J} being the
one left open, so six kernels of this gallery have a proved regularity index and no analytic
index in existence. That is the sharpest form the thesis of this volume takes.

The second is that the general routes to an index that this volume carries are not unconditional. The transparency
expansion of Conditional Theorem~\ref{cthm:main} holds under the estimate~$T$, the transfer of
Conditional Theorem~\ref{cthm:ex-alpha-eta} holds under a Mellin--Perron\index[terms]{Perron's formula} package, and the homogeneity
principle of \S\ref{sec:homogeneity} is a conjecture. Proving a gallery entry by one of these would
replace an established result by a conditional one. For the indices established in the gallery,
the constructive proofs supply unconditional arguments. The remaining conjectures and
conditional conclusions are marked at the entries concerned. Applications of the homogeneity
principle in the body are supported by gallery theorems that do not assume it.

This does not make the general routes idle, and they are used throughout the body. They say why a
threshold sits where it does rather than only that it does, they apply to families where a
construction applies to one kernel, and they predict the shape of an answer before any computation
is made. The gallery is also where they are tested. Every entry whose constructive answer agrees
with what a conditional route would return is evidence for the hypothesis that route carries, and
the two kinds of proof therefore support one another instead of competing. What the gallery does
not allow is the reverse order, an unproved hypothesis carrying a result that is then counted as
established.

The third is that a single technique would hide what the index is. The entries are solved by exact
first order recurrences and telescoping, by a dilative recurrence on a geometric scale, by Gamma
products, by discrete resolvents\index[terms]{resolvent} and Volterra\index[terms]{Volterra operator} inversion, by contour integration, by a second order
recurrence with its Riccati\index[terms]{Riccati equation} factorization and its Casoratian\index[terms]{Casoratian}, by generating functions of Mahler
type, and by direct estimation of the probe. Comparing these methods shows how the same
definition of an index applies to very different operators. A reader who follows one technique
learns a computation, and a reader who follows several sees an invariant.

Constructive here describes proofs that display how the operator acts on the forcing. Such proofs
may use contour integration, special functions and asymptotic expansions. In analytic number
theory, elementary means avoiding complex analysis, a condition several of these proofs do not meet.

\galleryentry{A}{The affine archetype}
 {$g(x)=(1-\lambda)x+\lambda$ on $(0,1]$, with $\lambda\in(0,1)$}
 {function of good variation}
 {$g^{*}(z)=\frac{z-\lambda}{z-1}$, single zero at $z=\lambda$}
 {$\alpha(g)=\lambda=\eta(g)$}
 {proved, Theorem~\ref{thm:affine_body}, stated in the body and proved here}

\rafgalleryfig{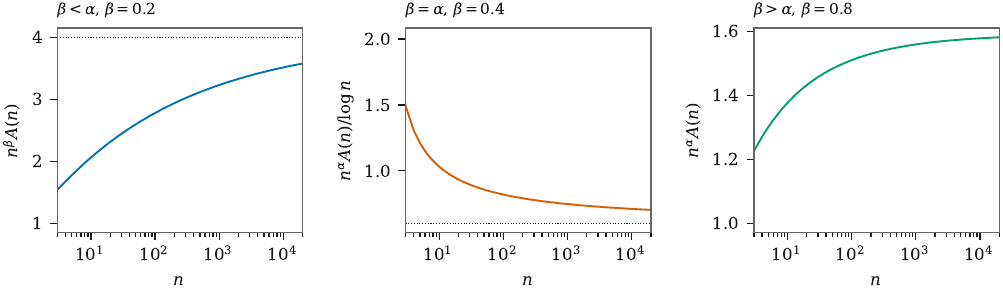}{Above, the affine profile at three values of the parameter, each running from $\lambda$ at the origin to one at the diagonal. Below, the three regimes at $\lambda=\tfrac25$, on a logarithmic scale of ranks. Below the index the partial sums follow the imposed rate and the constant approaches $1/g^{*}(\beta)$, marked by the dotted line. At the index the response carries an extra factor $\log n$, and the quotient approaches $1-\lambda$. Above the index the imposed rate is absorbed and the partial sums settle on $n^{-\alpha}$. The three lower panels are the three clauses of Definition~\ref{def:reg_index_fgv} read on the one kernel where every step is exact.}{fig:app_A}

For $\lambda \in (0,1)$, consider the affine function $g:(0,1] \to \mathbb{R}$ defined by
\[
    g(x) = (1-\lambda)x + \lambda.
\]
Its Mellin transform is given by:
\[
    g^{*}(z) = -z \int_{0}^{1} \left( (1-\lambda)t + \lambda \right) t^{-z-1} dt = (1-\lambda)\frac{-z}{1-z} + \lambda\frac{-z}{-z} = \frac{z-\lambda}{z-1}.
\]
The transform vanishes only at $z=\lambda$. That this value is the regularity index of $g$ is
Theorem~\ref{thm:affine_body} of the body, and the proof given here is constructive. It solves
the defining equation\index[terms]{defining equation} exactly and reads the three forcing ranges off that solution. The
solution rests on the elementary integration of a linear first-order recurrence, recorded first.

\begin{lemma}\label{lem:A_recurrence_solution}
Let $(x_n)_{n\geq 1}$, $(u_n)_{n\geq 2}$, and $(v_n)_{n\geq 2}$ be sequences satisfying
the recurrence $x_n = u_n x_{n-1} + v_n$ for $n \geq 2$, with $u_n\ne0$ for every $n\ge2$.
The solution is given by:
\[
    x_n = \left(\prod_{k=2}^{n} u_k\right) \left(x_1 + \sum_{k=2}^{n} \frac{v_k}{\prod_{i=2}^{k} u_i}\right) \quad \text{for } n \geq 2.
\]
\end{lemma}

\begin{proof}
Induction on $n$. The formula holds at $n=2$, where it reads $x_2=u_2x_1+v_2$. Assuming it at
rank $n-1$ and multiplying by $u_n$ adds the factor $u_n$ to the product and leaves the bracket
unchanged, so that $u_nx_{n-1}$ is the product up to $n$ times the bracket up to $n-1$, and
adding $v_n$ adds the term $k=n$ to the bracket, since
$v_n=\bigl(\prod_{k=2}^{n}u_k\bigr)\,v_n/\prod_{i=2}^{n}u_i$.
\end{proof}

\begin{proof}[Proof of Theorem~\ref{thm:affine_body}]
The aim is to show that $g$ is a function of good variation with regularity index
$\alpha(g)=\lambda$, through the asymptotic behavior of the partial sums
$A(n)=\sum_{k\le n}a_k$ for all real $\beta$.

\medskip
The recurrence for $A(n)$ comes first. The defining relation \cref{eq:defining_relation} is $\sum_{k=1}^{n} a_k g(k/n) = n^{-\beta}$. Substituting $g(x)$ yields:
\[
    \sum_{k=1}^{n} a_k \left( (1-\lambda)\frac{k}{n} + \lambda \right) = n^{-\beta}.
\]
Writing $\sum_{k=1}^n k a_k = nA(n) - \sum_{j=1}^{n-1} A(j)$ (Abel summation) and applying the same at rank $n-1$, after subtraction one obtains the exact first-order recurrence:
\begin{equation}\label{eq:A_recurrence}
    A(n) = \left(1-\frac{\lambda}{n}\right)A(n-1) + \frac{n^{1-\beta} - (n-1)^{1-\beta}}{n}.
\end{equation}

\medskip
Apply Lemma~\ref{lem:A_recurrence_solution} with $u_n = 1 - \lambda/n$ and
$v_n = (n^{1-\beta}-(n-1)^{1-\beta})/n$. Its nonvanishing hypothesis holds because
$0<\lambda<1$ and $n\ge2$ give $u_n>0$. The product is:
\[
    \prod_{k=2}^n u_k = \prod_{k=2}^n \left(1-\frac{\lambda}{k}\right) \sim \frac{C}{n^\lambda} \quad \text{as } n \to \infty,
\]
with $C=1/\Gamma(2-\lambda)$, since
$\prod_{k=2}^{n}(1-\lambda/k)=\Gamma(n+1-\lambda)/\bigl(\Gamma(2-\lambda)\Gamma(n+1)\bigr)$ by the
Euler-Gauss product formula and $\Gamma(n+1-\lambda)/\Gamma(n+1)\sim n^{-\lambda}$.
At $\beta=1$, every $v_k$ vanishes, so $A(n)=\prod_{k=2}^{n}u_k=\mathcal O(n^{-\lambda})$.
For $\beta\ne1$, the generic summand is
\[
\frac{v_k}{\prod_{i=2}^{k}u_i}
=\frac{k^{1-\beta}-(k-1)^{1-\beta}}{k}\cdot\frac{\Gamma(2-\lambda)\Gamma(k+1)}{\Gamma(k+1-\lambda)}
\sim(1-\beta)\,\Gamma(2-\lambda)\,k^{\lambda-\beta-1},
\]
using $k^{1-\beta}-(k-1)^{1-\beta}\sim(1-\beta)k^{-\beta}$.

\medskip
A case analysis on the convergence of the sum $\sum_{k=2}^n k^{-\beta+\lambda-1}$ concludes the argument.

\begin{itemize}
    \item For $\beta<\lambda$ the exponent $\lambda-\beta-1>-1$, so the sum diverges and is
    asymptotic to $(1-\beta)\Gamma(2-\lambda)n^{\lambda-\beta}/(\lambda-\beta)$. Multiplying by
    the product $n^{-\lambda}/\Gamma(2-\lambda)$ leaves
    \[
    A(n)\sim\frac{1-\beta}{\lambda-\beta}\,n^{-\beta}=\frac{n^{-\beta}}{g^{*}(\beta)},
    \]
    the two constants agreeing because $g^{*}(\beta)=(\beta-\lambda)/(\beta-1)$. The
    transparent constant is therefore the reciprocal of the transform, with no free factor left
    to determine.
    \item For $\beta=\lambda$ the sum is
    $\sim(1-\lambda)\Gamma(2-\lambda)\log n$, giving
    $A(n)\sim(1-\lambda)n^{-\lambda}\log n=\mathcal O(n^{-\lambda+\eps})$.
    \item For $\beta>\lambda$, $\beta\ne1$, the sum converges absolutely to a finite limit,
    so $A(n)=\mathcal O(n^{-\lambda})$. The case $\beta=1$ was treated above.
\end{itemize}

\medskip
The logarithm at $\beta=\lambda$ prevents transparency there, so no larger transparent
threshold is possible. Collecting the three ranges gives $\alpha(g)=\lambda=\eta(g)$.
\end{proof}

\begin{numobs}[Numerical control]\label{numobs:A_check}
The recurrence \eqref{eq:A_recurrence} reproduces forward substitution in the defining equation
to $2.9\cdot10^{-15}$ over $n\le400$ at $\beta=\tfrac1{10}$, $\tfrac32$ and $3$. At
$\lambda=\tfrac12$ the transparent constant is met, $n^{\beta}A(n)$ taking the values
$1.498152$, $1.499772$, $1.499973$, $1.499986$ at $n=10^{3}$, $10^{4}$, $10^{5}$,
$2\cdot10^{5}$ for $\beta=-\tfrac12$, against $1/g^{*}(-\tfrac12)=\tfrac32$, and
$1.964322$, $1.988716$, $1.996432$, $1.997477$ at $\beta=0$ against $1/g^{*}(0)=2$, the
approach being of order $n^{\beta-\lambda}$ as the case analysis predicts. At $\lambda=\tfrac13$
the same four ranks give $1.498611$, $1.499855$, $1.499985$, $1.499993$ at $\beta=-1$ against
$\tfrac32$, and $2.113427$, $2.135106$, $2.140827$, $2.141501$ at $\beta=-\tfrac14$ against
$\tfrac{15}7=2.142857$. Above the index the absorbed rate is the one stated, $n^{\lambda}A(n)$
settling on $1.128378$ at $\beta=1$ and on $0.527405$ at $\beta=2$ for $\lambda=\tfrac12$, and
on $1.107732$ at $\beta=1$ for $\lambda=\tfrac13$.
\end{numobs}

\galleryentry{B}{The slowly varying logarithmic profile}
 {$g(x)=1-\lambda\log x$ on $(0,1]$, with $\lambda\in(0,1)$}
 {function of good variation}
 {$g^{*}(z)=\frac{z-\lambda}{z}$, single zero at $z=\lambda$, pole at $z=0$}
 {$\alpha(g)=\lambda=\eta(g)$}
 {proved, Theorem~\ref{thm:B_index}}

\rafgalleryfig{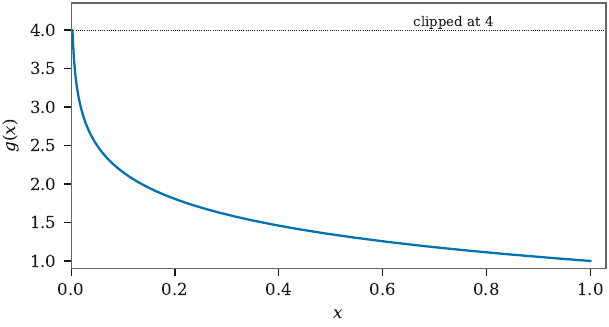}{The profile $g(x)=1-\tfrac12\log x$, drawn with the values above four removed. It is an unbounded, slowly varying profile and a proved instance of the second branch of Conjecture~\ref{conj:fgv_membership}. The power-singular family of Theorem~\ref{thm:fgv_power_singularity} lies beyond this scale.}{fig:app_B}

For $\lambda \in (0,1)$, consider the function $g:(0,1] \to \mathbb{R}$ defined by
\[
    g(x) = 1 - \lambda \log x .
\]
The Mellin transform, using the convention $g^{*}(z) = -z \int_0^1 g(t)\, t^{-z-1} dt$, is
\[
    g^{*}(z) = -z \int_0^1 \left(1 - \lambda \log t\right) t^{-z-1} dt
    = 1 - \frac{\lambda}{z} = \frac{z - \lambda}{z},
\]
the integral converging absolutely on $\Re z<0$ and the value on the right continuing it to
the whole plane. The transform has a single zero at $z=\lambda$ and a pole at $z=0$, so
$\eta(g)=\lambda$.

The threshold of this profile is read directly off the forced equation. Several forcing
exponents are compared in the statement below, so the subscript is kept on $A_\beta$ there,
while the proof, at one fixed exponent, writes $A(n)$ as elsewhere in the volume, in the
notation of Definition~\ref{def:reg_index_fgv}.

\begin{theorem}\label{thm:B_index}
The defining equation $\sum_{k\le n}a_k\,g(k/n)=n^{-\beta}$ has transparency threshold
$\lambda$. For every $\beta<\lambda$,
\[
A_\beta(n)=\frac{1}{g^{*}(\beta)}\,n^{-\beta}+o(n^{-\beta}),
\qquad \frac{1}{g^{*}(\beta)}=\frac{\beta}{\beta-\lambda},
\]
for $\beta=\lambda$ one has $A_\lambda(n)=\mathcal{O}(n^{-\lambda}\log n)$, and for every
$\beta>\lambda$ one has $A_\beta(n)=\mathcal{O}(n^{-\lambda})$. Transparency fails at a set
of exponents accumulating at $\lambda$ from above, so the threshold is attained and
$\alpha(g)=\lambda=\eta(g)$.
\end{theorem}

\begin{proof}
The defining relation reads $A(n)+\lambda\log(n)A(n)-\lambda S(n)=n^{-\beta}$ with
$S(n)=\sum_{k\le n}a_k\log k$. Subtracting the relation at rank $n-1$ and using
$S(n)-S(n-1)=a_n\log n$ together with $a_n=A(n)-A(n-1)$, the terms carrying $\log n$ cancel
and the exact first order recurrence
\begin{equation}\label{eq:B_recurrence}
A(n)=u_n\,A(n-1)+v_{\beta,n},
\qquad
u_n=1-\lambda\log\frac{n}{n-1},
\qquad
v_{\beta,n}=n^{-\beta}-(n-1)^{-\beta},
\end{equation}
holds for $n\ge2$, with $A(1)=1$ since $g(1)=1$. The recurrence of
Appendix~\ref{app:A} is not the one governing this kernel, the forcing increment being here
the difference of consecutive powers rather than its Abel transform\index[terms]{Abel transform}.

Expanding $\log\frac{n}{n-1}=\frac1n+\frac1{2n^2}+\mathcal{O}(n^{-3})$ and
$(n-1)^{-\beta}=n^{-\beta}(1+\beta/n+\mathcal{O}_\beta(n^{-2}))$ gives
\[
u_n=1-\frac{\lambda}{n}+\mathcal{O}(n^{-2}),
\qquad
v_{\beta,n}=-\beta\,n^{-\beta-1}+\mathcal{O}_\beta(n^{-\beta-2}),
\]
so \eqref{eq:B_recurrence} satisfies the hypotheses of Proposition~\ref{thm:ex-first-order} with
$\gamma=\lambda$, $b_\beta=-\beta$ and $\delta=1$. Positivity of $u_n$ holds for every
$n\ge2$, since $\lambda\log\frac{n}{n-1}\le\lambda\log2<1$. The three regimes
\eqref{eq:ex-first-order-below}, \eqref{eq:ex-first-order-critical} and
\eqref{eq:ex-first-order-above} follow, and the coefficient match
\eqref{eq:ex-first-order-match} is the identity
\[
\frac{b_\beta}{\gamma-\beta}=\frac{-\beta}{\lambda-\beta}=\frac{\beta}{\beta-\lambda}
=\frac{1}{g^{*}(\beta)} ,
\]
which gives transparency below $\lambda$, the value $\beta=0$ being the pole of $g^{*}$ where
the transparency constant is read as $0$ by the convention of
Definition~\ref{def:reg_index_fgv}.

Sharpness is the remaining point, and it is settled by the connection coefficient\index[terms]{connection coefficient} of
\eqref{eq:ex-first-order-connection}. Write
\[
P_m=\prod_{j=2}^{m}u_j>0,
\qquad
D(\beta)=1+\sum_{m\ge2}\frac{v_{\beta,m}}{P_m},
\]
so that $C_\beta=q_\infty D(\beta)$ in the notation of Proposition~\ref{thm:ex-first-order} taken
with $n_0=1$. Since $q_m=m^{\lambda}P_m$ converges to a finite positive limit, there are
constants $0<c_1\le c_2$ with $c_1m^{-\lambda}\le P_m\le c_2m^{-\lambda}$ for every $m\ge2$.
The mean value theorem gives $|v_{\beta,m}|\le|\beta|(m-1)^{-\Re\beta-1}$, so on a compact
subset of $\{\Re\beta>\lambda\}$, where $\Re\beta\ge\lambda+\eps$ for some
$\eps>0$, the general term of the series is dominated by a constant multiple of
$m^{-1-\eps}$. The series therefore converges uniformly on compact subsets and
$D$ is holomorphic on the half plane $\{\Re\beta>\lambda\}$, which is connected.

The behavior of $D$ along the real axis at infinity locates one point where it does not
vanish. The term of index $m=2$ equals $(2^{-\beta}-1)/(1-\lambda\log2)$ and tends to
$-1/(1-\lambda\log2)$ as $\beta\to+\infty$. For the remaining terms, the bound above gives
\[
\Big|\sum_{m\ge3}\frac{v_{\beta,m}}{P_m}\Big|
\le\frac{\beta}{c_1}\sum_{m\ge3}m^{\lambda}(m-1)^{-\beta-1}
\le\frac{2^{\lambda}\beta}{c_1}\sum_{j\ge2}j^{\lambda-\beta-1}
\le\frac{2^{\lambda}\beta}{c_1}\Big(2^{\lambda-\beta-1}+\frac{2^{\lambda-\beta}}{\beta-\lambda}\Big),
\]
which tends to $0$ as $\beta\to+\infty$. Hence
\[
\lim_{\beta\to+\infty}D(\beta)=1-\frac{1}{1-\lambda\log2}
=\frac{-\lambda\log2}{1-\lambda\log2}<0 ,
\]
the inequality holding because $0<1-\lambda\log2<1$ for $\lambda\in(0,1)$. So $D$ is a
holomorphic function on a connected open set, not identically zero, and its zero set is
discrete. There are therefore real exponents $\beta_j\downarrow\lambda$ with
$D(\beta_j)\neq0$, hence $C_{\beta_j}\neq0$, while $g^{*}(\beta_j)=(\beta_j-\lambda)/\beta_j$
is finite and nonzero. The sharpness clause of Proposition~\ref{thm:ex-first-order} applies and
gives $\alpha(g)=\lambda$.
\end{proof}

\begin{remark}\label{rem:B_class}
This entry sits in the gallery by design. It is unbounded at the origin, so it does not
belong to the bounded class under which Definition~\ref{def:reg_index_fgv} is stated, and it
nevertheless has a determined index, which is the point it is there to make. Boundedness is
a convenience of the definition and not a mechanism of the theory, and the proof above uses
only the exact recurrence \eqref{eq:B_recurrence}. The scale of the singularity is what
matters. Since $g(1/u)=1+\lambda\log u$ is slowly varying at infinity, the kernel falls under
Proposition~\ref{prop:fgv_slowly_varying}, the transform is given by an absolutely convergent
integral on the whole half plane $\Re z<0$, and the arithmetic Riemann sums converge to it
there. The analytic index is therefore available in the usual sense, and the coincidence
$\alpha=\eta=\lambda$ is read on the same footing as for the bounded entries. Appendix~\ref{app:F}
carries the same logarithmic scale. These two examples support the slowly varying branch of
Conjecture~\ref{conj:fgv_membership}. Theorem~\ref{thm:fgv_power_singularity} supplies a distinct
power-singular family beyond it, without enlarging the hypotheses of that conjecture.
\end{remark}

\begin{numobs}\label{numobs:B_check}
At $\lambda=2/5$ the recurrence \eqref{eq:B_recurrence} reproduces forward substitution in
the defining equation to $10^{-39}$ over $n\le90$ in $40$ digit arithmetic. At
$\beta=-1/2$ the ratio $n^{\beta}A(n)$ equals $0.555559$ at $n=10^{6}$ against
$1/g^{*}(\beta)=5/9$. At $\beta=7/20$ the same ratio follows
$-7+7.86\,n^{\beta-\lambda}$ over three decades. Above the threshold $n^{\lambda}A(n)$
converges to $q_\infty D(\beta)$, the two sides agreeing to six decimals at $\beta=3/2$. The
truncated coefficient $D(\beta)$ stays negative for $\beta$ down to $\lambda+10^{-3}$ and
reaches $-0.3836212740$ at $\beta=60$ against the predicted limit
$-\lambda\log2/(1-\lambda\log2)=-0.3836212740$.
\end{numobs}

\galleryentry{C}{The bipolar step profile}
 {$g(x)=\mu$ on $(0,\lambda]$ and $g(x)=1$ on $(\lambda,1]$, with $0<\lambda<1$ and $0<\mu<1$}
 {function of good variation}
 {$g^{*}(z)=1-q\,\lambda^{-z}$ with $q=1-\mu$, unique real zero at $z_0=\log(1-\mu)/\log\lambda$}
 {$\alpha(g)=z_0=\eta(g)$}
 {proved, Theorem~\ref{thm:C_index}}

Fix $0<\lambda<1$ and $0<\mu<1$. Define the step function $g:(0,1]\to\mathbb{R}$ by
\[
    g(x)=
    \begin{cases}
        \mu, & 0<x\le \lambda,\\[1mm]
        1,   & \lambda < x \le 1 .
    \end{cases}
\]
The Mellin transform is $g^{*}(z) = 1 - q\,\lambda^{-z}$ with $q=1-\mu$, its unique zero on
the real line sits at
\[
z_0=\frac{\log(1-\mu)}{\log\lambda}>0 ,
\]
and $\eta(g)=z_0$.

For this profile the two indices agree, and Theorem~\ref{thm:C_index} is what proves it.
Several forcing exponents are compared in what follows, so the subscript is kept on $A_\beta$
and on the sequences attached to it. Elsewhere in the volume the exponent is fixed by its
context and the plain $a_n$ and $A(n)$ of Definition~\ref{def:reg_index_fgv} are used.

\begin{theorem}\label{thm:C_index}
The kernel $g$ is a function of good variation with $\alpha(g)=z_0=\eta(g)$. For every
$\beta<z_0$,
\[
A_\beta(n)=\frac{1}{g^{*}(\beta)}\,n^{-\beta}
+\mathcal{O}_\beta\!\left(n^{-z_0}+n^{-\beta-1}\log n\right),
\]
and for every $\beta>0$ the partial sums obey the two sided estimate
\begin{equation}\label{eq:C_two_sided}
\lambda^{\beta}\,n^{-z_0}\ \le\ A_\beta(n)\ \le\ C_\beta\,n^{-z_0}
\qquad(\beta\ge z_0,\ n\ge1),
\end{equation}
with a logarithmic factor at $\beta=z_0$. Above the threshold the partial sums are positive
while $1/g^{*}(\beta)$ is negative, so transparency fails at every such exponent.
\end{theorem}

\begin{proof}
Splitting the defining sum at $m=\lfloor\lambda n\rfloor$, the condition $k/n\le\lambda$ is
$k\le\lfloor\lambda n\rfloor$, so the sum equals $\mu A(m)+(A(n)-A(m))$ and the defining
relation becomes the exact dilative recurrence
\begin{equation}\label{eq:C_recurrence}
    A_\beta(n) = q\,A_\beta(\lfloor \lambda n \rfloor) + n^{-\beta},
    \qquad q = 1-\mu .
\end{equation}
Setting $n_0=n$, $n_{j+1}=\lfloor\lambda n_j\rfloor$ and $J=J(n)=\max\{j:n_j\ge1\}$, with the
convention $A_\beta(0)=0$, iteration of \eqref{eq:C_recurrence} is exact and gives
\begin{equation}\label{eq:C_iterated}
    A_\beta(n)=\sum_{j=0}^{J}q^{j}\,n_j^{-\beta} .
\end{equation}
Every term of \eqref{eq:C_iterated} is positive, a fact used twice below.

The trajectory obeys $\lambda n_j-1<n_{j+1}\le\lambda n_j$, whence by induction
\begin{equation}\label{eq:C_traj}
\lambda^{j}n-\frac{1}{1-\lambda}<n_j\le\lambda^{j}n\qquad(0\le j\le J).
\end{equation}
The terminal ranks are the ranks at which the comparison $n_j\approx\lambda^{j}n$ ceases to
be accurate, and they are handled by the geometric decay of the weight rather than by that
comparison. Fix $N_0>1+\frac{2}{1-\lambda}$ and let $J_0=\min\{j:n_j<N_0\}$, so that
$J_0\le J$. For $j\ge J_0$ one has $1\le n_j<N_0$, hence
$n_j^{-\beta}\le\max(1,N_0^{|\beta|})$, and
\[
\Big|\sum_{j\ge J_0}q^{j}n_j^{-\beta}\Big|
\le\max(1,N_0^{|\beta|})\sum_{j\ge J_0}q^{j}
=\mathcal{O}_{\beta,\lambda,\mu}\!\left(q^{J_0}\right).
\]
By \eqref{eq:C_traj} at $j=J_0$, the inequality $\lambda^{J_0}n-\frac1{1-\lambda}<N_0$ gives
$\lambda^{J_0}<N_1/n$ with $N_1=N_0+\frac1{1-\lambda}$, and taking logarithms against the
negative number $\log\lambda$ gives $J_0>\log(n/N_1)/|\log\lambda|$. Since $\log q<0$ and
$\log q/\log\lambda=z_0$,
\begin{equation}\label{eq:C_terminal}
q^{J_0}<\exp\!\Big(\log q\cdot\frac{\log(n/N_1)}{|\log\lambda|}\Big)
=\Big(\frac{n}{N_1}\Big)^{-z_0}=N_1^{z_0}\,n^{-z_0},
\end{equation}
so the terminal ranks contribute $\mathcal{O}_{\beta,\lambda,\mu}(n^{-z_0})$ in total.

For $j<J_0$ one has $n_j\ge N_0$, so \eqref{eq:C_traj} reads $n_j=\lambda^{j}n(1+\theta_j)$
with $-\frac1{(1-\lambda)N_0}\le\theta_j\le0$, and $|\theta_j|\le\frac12$ by the choice of
$N_0$. Hence $n_j^{-\beta}=(\lambda^{j}n)^{-\beta}\bigl(1+\mathcal{O}_\beta(|\theta_j|)\bigr)$
with $|\theta_j|\le\frac{1}{(1-\lambda)\lambda^{j}n}$, and
\[
\sum_{j<J_0}q^{j}n_j^{-\beta}
=n^{-\beta}\sum_{j<J_0}r^{j}
+\mathcal{O}_\beta\Big(\frac{n^{-\beta-1}}{1-\lambda}\sum_{j<J_0}\bigl(r\lambda^{-1}\bigr)^{j}\Big),
\qquad r=q\lambda^{-\beta}.
\]
The error sum is $\mathcal{O}(n^{-\beta-1})$ when $r\lambda^{-1}<1$ and
$\mathcal{O}(n^{-\beta-1}\log n)$ when $r\lambda^{-1}=1$. When $r\lambda^{-1}>1$ it is
$\mathcal{O}\bigl(n^{-\beta-1}(r\lambda^{-1})^{J_0}\bigr)$, and that case forces
$\lambda^{-\beta-1}>1/q>1$, hence $\beta>-1$. The minimality of $J_0$ gives
$n_{J_0-1}\ge N_0$, and $n_{J_0-1}\le\lambda^{J_0-1}n$ then gives
$\lambda^{J_0}\ge\lambda N_0/n$, so that
$\lambda^{-(\beta+1)J_0}\le\bigl(n/(\lambda N_0)\bigr)^{\beta+1}$. Together with
\eqref{eq:C_terminal} this yields
\[
(r\lambda^{-1})^{J_0}=q^{J_0}\lambda^{-(\beta+1)J_0}
=\mathcal{O}_{\beta,\lambda,\mu}\!\left(n^{\beta+1-z_0}\right),
\]
so the error is again $\mathcal{O}_\beta(n^{-z_0})$. In every case the error is
$\mathcal{O}_\beta(n^{-\beta-1}\log n+n^{-z_0})$.

The value of $r$ decides the main sum. Since $\log\lambda<0$, the inequality $r<1$ is
$\beta<z_0$. For $\beta<z_0$ the sum $\sum_{j<J_0}r^{j}$ equals
$\frac{1}{1-r}-\frac{r^{J_0}}{1-r}$, and $n^{-\beta}r^{J_0}=q^{J_0}\lambda^{-\beta J_0}
n^{-\beta}=\mathcal{O}_\beta(n^{-z_0})$ by the same two bounds, so
\[
A_\beta(n)=\frac{n^{-\beta}}{1-q\lambda^{-\beta}}
+\mathcal{O}_\beta\!\left(n^{-z_0}+n^{-\beta-1}\log n\right)
=\frac{1}{g^{*}(\beta)}n^{-\beta}+o(n^{-\beta}),
\]
which is transparency, the error being $o(n^{-\beta})$ precisely because $\beta<z_0$. At
$\beta=z_0$ the sum has $J_0=\mathcal{O}(\log n)$ terms equal to $1$, giving
$A_{z_0}(n)=\mathcal{O}(n^{-z_0}\log n)$. For $\beta>z_0$ one has $r>1$ and
$\sum_{j<J_0}r^{j}\le r^{J_0}/(r-1)$, so $A_\beta(n)=\mathcal{O}_\beta(n^{-z_0})$, which is
the upper bound in \eqref{eq:C_two_sided}.

The lower bound comes from positivity. Keeping only the last term of \eqref{eq:C_iterated}
and using $\lambda n_J<1\le n_J$, valid because $n_{J+1}=\lfloor\lambda n_J\rfloor=0$, gives
$n_J^{-\beta}>\lambda^{\beta}$ for $\beta>0$. From \eqref{eq:C_traj}, $1\le n_J\le\lambda^{J}n$
gives $\lambda^{J}\ge1/n$, hence $J\le\log n/|\log\lambda|$, and since $\log q<0$ this bound
on $J$ gives
\[
q^{J}\ \ge\ \exp\!\Big(\log q\cdot\frac{\log n}{|\log\lambda|}\Big)
\ =\ n^{-z_0}.
\]
Therefore
\[
A_\beta(n)\ \ge\ q^{J}n_J^{-\beta}\ \ge\ \lambda^{\beta}\,n^{-z_0}\qquad(\beta>0,\ n\ge1),
\]
which is the lower bound in \eqref{eq:C_two_sided}.

Sharpness follows at once. For $\beta>z_0$ the transform value
$g^{*}(\beta)=1-q\lambda^{-\beta}$ is negative, so the transparency profile
$n^{-\beta}/g^{*}(\beta)$ is negative, whereas $A_\beta(n)$ is positive by
\eqref{eq:C_iterated}. Independently of sign, $n^{\beta}A_\beta(n)\ge\lambda^{\beta}
n^{\beta-z_0}$ tends to infinity. Transparency fails at every exponent above $z_0$, so the
supremum of Definition~\ref{def:reg_index_fgv} is exactly $z_0$, the absorption estimate
$\mathcal{O}(n^{-z_0+\eps})$ holds at and above it, and $g$ is a function of good
variation with $\alpha(g)=z_0=\eta(g)$.
\end{proof}

\begin{numobs}\label{numobs:C_check}
For the pairs $(\lambda,\mu)$ equal to $(1/2,1/2)$, $(3/10,7/10)$, $(0.618,1/4)$ and
$(9/10,1/10)$, the recurrence \eqref{eq:C_recurrence} and the iterated form
\eqref{eq:C_iterated} reproduce forward substitution in the defining equation to machine
accuracy over $n\le400$ and over the six exponents $\beta\in\{-2/5,0,z_0/2,z_0,1.3z_0,3z_0\}$.
At $(\lambda,\mu)=(1/2,1/2)$, where $z_0=1$, the ratio $n^{\beta}A_\beta(n)$ at $n=10^{6}$
equals $1.610122$ against $1/g^{*}(-2/5)=1.610122$, and $2.601160$ against
$1/g^{*}(3/10)=2.601269$. Above the threshold $n^{z_0}A_\beta(n)$ stays positive and bounded
while $1/g^{*}(\beta)$ is negative. Over $n\le3\cdot10^{5}$ the quantity $n^{z_0}A_\beta(n)$
has minimum $1$ for each pair tested, against the proved lower bound $\lambda^{\beta}$.
\end{numobs}

\galleryentry{D}{The integer scale broken kernel}
 {$g(x)=x\,\lambda^{\lfloor-\log_\lambda x\rfloor}$ on $(0,1]$, with $\lambda\ge2$ an integer}
 {function of good variation}
 {$g^{*}(z)=\frac{z}{z-1}\cdot\frac{1-\lambda^{z-1}}{1-\lambda^{z}}$, first zero on the real line at $z=1$}
 {$\alpha(g)=1=\eta(g)$}
 {proved, Theorem~\ref{thm:D_index}}

\rafgalleryfig{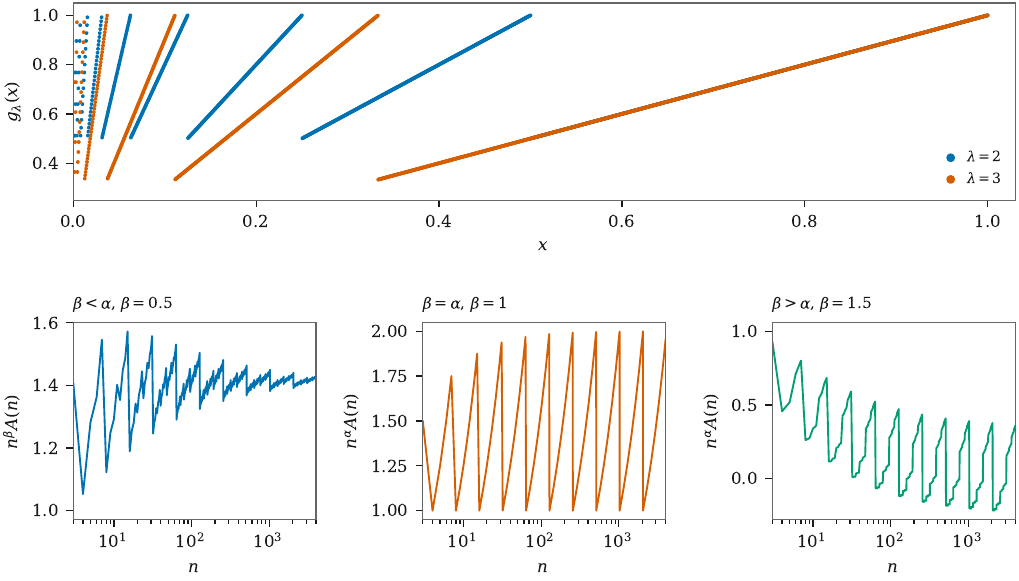}{Above, the profile $g_\lambda$ at the scales $\lambda=2$ and $\lambda=3$, a fan of segments each reaching the value one at $x=\lambda^{-i}$ and dropping to $1/\lambda$ just to the right of that point. Below, the three regimes at $\lambda=2$, where $\alpha=1$. The partial sums follow the imposed rate below the index, oscillate boundedly at it, and settle on the rate $n^{-\alpha}$ above it, the oscillation having period $\log\lambda$ on the logarithmic scale in all three.}{fig:app_D}

Let $\lambda\ge2$ be an integer and $g(x)=x\,\lambda^{\lfloor-\log_\lambda(x)\rfloor}$ on
$(0,1]$. On each layer $(\lambda^{-(j+1)},\lambda^{-j}]$ one has $g(x)=\lambda^{j}x$, so
$\lambda^{-1}<g\le1$ and the kernel is bounded with $g(1)=1$. The Mellin transform is
\[
g^{*}(z) = \frac{z}{z-1}\cdot\frac{1-\lambda^{z-1}}{1-\lambda^{z}},
\]
whose first zero on the real line sits at $z=1$, so $\eta(g)=1$.

The index of this kernel is one, and the three regimes are explicit.

\begin{theorem}\label{thm:D_index}
The kernel $g$ is a function of good variation with $\alpha(g)=1=\eta(g)$. For every
$\beta<1$,
\[
A_\beta(n)=\frac{1}{g^{*}(\beta)}\,n^{-\beta}+o(n^{-\beta}),
\]
the transform being regular at the origin with
$g^{*}(0)=(\lambda-1)/(\lambda\log\lambda)$, at $\beta=1$ one has $A_1(n)=\lambda^{-\lfloor\log_\lambda n\rfloor}$ exactly, and for
$\beta>1$ one has $A_\beta(n)=\mathcal{O}(n^{-1})$ with transparency failing at a set of
exponents accumulating at $1$ from above.
\end{theorem}

The proof runs through the weighted sums $B(n)=\sum_{k\le n}ka_k$, which satisfy an exact
dilative recurrence, and then returns to $A(n)$ through an Abel identity\index[terms]{Abel identity} whose terminal
constant is computed by a telescoping argument.

\begin{lemma}\label{lem:D_Brec}
For $n\ge\lambda$ and $m=\lfloor n/\lambda\rfloor$ one has
$\lfloor-\log_\lambda(k/n)\rfloor=\lfloor-\log_\lambda(k/m)\rfloor+1$ for every $k\le m$, and
consequently
\begin{equation}\label{eq:D_Bn-recurrence}
  B(n) = B(m) + n^{1-\beta} - \lambda\, m^{1-\beta},
  \qquad m=\lfloor n/\lambda\rfloor,
\end{equation}
with $B(n)=n^{1-\beta}$ for $1\le n<\lambda$ and $B(0)=0$.
\end{lemma}

\begin{proof}
The claimed identity is $\lfloor\log_\lambda(n/k)\rfloor=\lfloor\log_\lambda(\lambda m/k)\rfloor$
for $k\le m$. It fails only if some power $\lambda^{i}$ satisfies
$\lambda m<k\lambda^{i}\le n$. From $\lambda m\le n<\lambda m+\lambda$ such a power would give
$m<k\lambda^{i-1}<m+1$. For $i\ge1$ the quantity $k\lambda^{i-1}$ is an integer, which
excludes it. For $i=0$ the requirement $k>\lambda m$ contradicts $k\le m$, and for $i\le-1$
one has $k\lambda^{i}\le m/\lambda<\lambda m$. Multiplying the defining relation by $n$
turns it into $\sum_{k\le n}ka_k\lambda^{j(k,n)}=n^{1-\beta}$ with
$j(k,n)=\lfloor-\log_\lambda(k/n)\rfloor$. The ranks $m<k\le n$ have $j=0$ and contribute
$B(n)-B(m)$, and the ranks $k\le m$ have $j(k,n)=j(k,m)+1$ and contribute $\lambda m^{1-\beta}$
by the same relation at $m$. The initial segment comes from the top layer alone, where
$g(k/n)=k/n$ for every $k\le n$.
\end{proof}

The weighted sums fall into three regimes according to where the forcing exponent sits.

\begin{lemma}\label{lem:D_Bregimes}
Let $r=\lambda^{\beta-1}$ and $c_\beta=\dfrac{1-\lambda^{\beta}}{1-\lambda^{\beta-1}}$. Then
\[
B(n)=c_\beta\,n^{1-\beta}+\mathcal{O}_\beta\!\left(n^{-\beta}+1\right)\quad(\beta<1),
\qquad
B(n)=1-\lfloor\log_\lambda n\rfloor(\lambda-1)\quad(\beta=1),
\]
and $B(n)=\mathcal{O}_\beta(1)$ for $\beta>1$. At $\beta=0$ the solution is the digital sum,
$B(n)=s_\lambda(n)$.
\end{lemma}

\begin{proof}
Write $n_0=n$, $n_{j+1}=\lfloor n_j/\lambda\rfloor$ and let $J$ be the first index with
$n_J<\lambda$. Iterating \eqref{eq:D_Bn-recurrence} gives
$B(n)=n_J^{1-\beta}+\sum_{j<J}w(n_j)$ with
$w(t)=t^{1-\beta}-\lambda\lfloor t/\lambda\rfloor^{1-\beta}$. Since
$\lfloor t/\lambda\rfloor=(t/\lambda)(1+\mathcal{O}(1/t))$, the mean value theorem gives
$w(t)=(1-\lambda^{\beta})t^{1-\beta}+\mathcal{O}_\beta(t^{-\beta})$. The trajectory obeys
$\lambda^{-j}n-\frac{\lambda}{\lambda-1}<n_j\le\lambda^{-j}n$, and the ranks with
$n_j\ge N_0$ satisfy $n_j^{1-\beta}=(\lambda^{-j}n)^{1-\beta}(1+\mathcal{O}_\beta(1/n_j))$,
exactly as in Appendix~\ref{app:C}. For $\beta<1$ the series $\sum_j\lambda^{-j(1-\beta)}$
converges to $1/(1-\lambda^{\beta-1})$, the finitely many terminal ranks contribute
$\mathcal{O}_\beta(1)$, and the stated estimate follows. At $\beta=1$ the forcing is
$n^{0}-\lambda m^{0}=1-\lambda$ at every rank with $m\ge1$, so the iteration is exact and
gives $B(n)=1-\lfloor\log_\lambda n\rfloor(\lambda-1)$. For $\beta>1$ the same series
converges after the substitution and $B(n)$ stays bounded. At $\beta=0$ the forcing is the
remainder $n-\lambda\lfloor n/\lambda\rfloor$, the least significant digit of $n$ in base
$\lambda$, and repeated Euclidean division gives the digital sum $s_\lambda(n)$, whence
$0\le B(n)\le(\lambda-1)(1+\lfloor\log_\lambda n\rfloor)$. No global equivalent
$B(n)\sim C\log n$ holds at $\beta=0$, the quotient $B(n)/\log n$ tending to $0$ along
$n=\lambda^{m}$ and to $(\lambda-1)/\log\lambda$ along $n=\lambda^{m}-1$.
\end{proof}

The passage from $B$ to $A$ rests on the finite Abel identity
\begin{equation}\label{eq:D_abel}
A(n)=\frac{B(n)}{n}+\sum_{k<n}\frac{B(k)}{k(k+1)},
\end{equation}
which is exact for every $n\ge1$. For $\beta>0$ the series on the right converges absolutely
by Lemma~\ref{lem:D_Bregimes}, and its value is what the following lemma computes.

\begin{lemma}\label{lem:D_telescope}
For every $\beta>0$,
\[
T:=\sum_{k\ge1}\frac{B(k)}{k(k+1)}=0,
\qquad\text{hence}\qquad
A(n)=\frac{B(n)}{n}-\sum_{k\ge n}\frac{B(k)}{k(k+1)} ,
\]
while at $\beta=0$ the same series converges to
$T_0=\lambda\log\lambda/(\lambda-1)=1/g^{*}(0)$.
\end{lemma}

\begin{proof}
Absolute convergence holds since $B(k)=\mathcal{O}_\beta(k^{1-\beta}+\log k)$ by
Lemma~\ref{lem:D_Bregimes}, which makes the general term
$\mathcal{O}_\beta(k^{-1-\beta}+k^{-2}\log k)$. The ranks $k$ with
$\lfloor k/\lambda\rfloor=i$ are exactly $\lambda i\le k\le\lambda i+\lambda-1$, and on that
block
\[
\sum_{k=\lambda i}^{\lambda i+\lambda-1}\frac{1}{k(k+1)}
=\frac{1}{\lambda i}-\frac{1}{\lambda(i+1)}
=\frac{1}{\lambda}\cdot\frac{1}{i(i+1)} ,
\]
the sum telescoping. Splitting $T$ at $k=\lambda$ and inserting
\eqref{eq:D_Bn-recurrence} gives
\[
T=\sum_{k<\lambda}\frac{k^{1-\beta}}{k(k+1)}
+\sum_{k\ge\lambda}\frac{B(\lfloor k/\lambda\rfloor)}{k(k+1)}
+\sum_{k\ge\lambda}\frac{k^{1-\beta}-\lambda\lfloor k/\lambda\rfloor^{1-\beta}}{k(k+1)} .
\]
The middle sum equals $\frac1\lambda\sum_{i\ge1}B(i)/(i(i+1))=T/\lambda$ by the block
identity, and the same identity applied to the subtracted part of the third sum gives
\[
\sum_{k\ge\lambda}\frac{\lambda\lfloor k/\lambda\rfloor^{1-\beta}}{k(k+1)}
=\lambda\sum_{i\ge1}i^{1-\beta}\cdot\frac{1}{\lambda\,i(i+1)}
=\sum_{i\ge1}\frac{i^{-\beta}}{i+1} .
\]
The two remaining pieces are the first sum and $\sum_{k\ge\lambda}k^{-\beta}/(k+1)$, whose
total is $\sum_{k\ge1}k^{-\beta}/(k+1)$, the very series just produced. They cancel, and
$T=T/\lambda$. Since $\lambda\ge2$, this forces $T=0$. All rearrangements are legitimate,
each series involved being absolutely convergent for $\beta>0$.

At $\beta=0$ the series $\sum_{k}k^{-\beta}/(k+1)$ diverges and the cancellation has to be
carried out at finite level. There $B(k)=s_\lambda(k)$ and the same splitting gives
$T_0(1-1/\lambda)=W_\lambda$ with
\[
W_\lambda=\sum_{k\ge1}\frac{k-\lambda\lfloor k/\lambda\rfloor}{k(k+1)}
=\lim_{N\to\infty}\Big[\sum_{k=2}^{\lambda N}\frac1k-\sum_{i=2}^{N}\frac1i\Big]
=\lim_{N\to\infty}\bigl(H_{\lambda N}-H_{N}\bigr)=\log\lambda ,
\]
the middle equality using the block identity once more. Hence
$T_0=\lambda\log\lambda/(\lambda-1)$, which is $1/g^{*}(0)$ since
$g^{*}(z)=\frac{z}{z-1}\frac{1-\lambda^{z-1}}{1-\lambda^{z}}$ has the removable value
$(\lambda-1)/(\lambda\log\lambda)$ at the origin.
\end{proof}

The partial sums obey an exact recurrence under dilation by the scale.

\begin{lemma}\label{lem:D_Arec}
For every $\beta>0$ and every $m\ge1$,
\begin{equation}\label{eq:D_Arec}
A(\lambda m)=\frac{A(m)}{\lambda}+m^{-\beta}\bigl(\lambda^{-\beta}-1\bigr)+\rho(m),
\qquad
\rho(m)=\sum_{i=m}^{\lambda m-1}\frac{i^{-\beta}}{i+1}>0 .
\end{equation}
\end{lemma}

\begin{proof}
Write $R(n)=\sum_{k\ge n}B(k)/(k(k+1))$, so that $A(n)=B(n)/n-R(n)$ by
Lemma~\ref{lem:D_telescope}. Repeating the block computation of that lemma on the tail
starting at $\lambda m$ gives $R(\lambda m)=R(m)/\lambda-\rho(m)$. Combining this with
\eqref{eq:D_Bn-recurrence} at $n=\lambda m$, where $\lfloor n/\lambda\rfloor=m$ exactly,
\[
A(\lambda m)=\frac{B(m)+(\lambda m)^{1-\beta}-\lambda m^{1-\beta}}{\lambda m}
-\frac{R(m)}{\lambda}+\rho(m)
=\frac{A(m)}{\lambda}+m^{-\beta}\bigl(\lambda^{-\beta}-1\bigr)+\rho(m). \qedhere
\]
\end{proof}

The critical forcing is settled inside this appendix, by the same weighted sums.

\begin{lemma}\label{lem:D_critical}
At $\beta=1$ the solution is carried by the powers of $\lambda$ alone, with $a_1=1$,
$a_{\lambda^{e}}=-(\lambda-1)\lambda^{-e}$ for $e\ge1$, and $a_n=0$ at every other rank.
Consequently
\begin{equation}\label{eq:D_critical_closed}
A(n)=\lambda^{-\lfloor\log_\lambda n\rfloor}
\qquad(n\ge1),
\end{equation}
and $nA(n)$ runs over $[1,\lambda)$ on each block $\lambda^{e}\le n<\lambda^{e+1}$, with
$\liminf_{n\to\infty}nA(n)=1$ and $\limsup_{n\to\infty}nA(n)=\lambda$.
\end{lemma}

\begin{proof}
Lemma~\ref{lem:D_Bregimes} gives $B(n)=1-(\lambda-1)\lfloor\log_\lambda n\rfloor$ at
$\beta=1$. Since $na_n=B(n)-B(n-1)$ with $B(0)=0$, the coefficient vanishes at every rank
where $\lfloor\log_\lambda n\rfloor$ does not jump, and $\lambda\ge2$ being an integer that
happens at every rank except $n=1$ and the powers $n=\lambda^{e}$ with $e\ge1$, where the
jump is $-(\lambda-1)$. This gives the stated coefficients. Summing them,
\[
A(n)=1-(\lambda-1)\!\!\sum_{1\le e\le\lfloor\log_\lambda n\rfloor}\!\!\lambda^{-e}
=1-\bigl(1-\lambda^{-\lfloor\log_\lambda n\rfloor}\bigr),
\]
which is \eqref{eq:D_critical_closed}. On the block $\lambda^{e}\le n<\lambda^{e+1}$ the
product $nA(n)=n\lambda^{-e}$ increases from $1$, attained at $n=\lambda^{e}$, to
$\lambda-\lambda^{-e}$, attained at $n=\lambda^{e+1}-1$, which gives the two extreme limits.
\end{proof}

\begin{remark}\label{rem:D_no_E}
The comparison with integral scales of the research dossier~\ref{app:dossier_sqrt2},
Proposition~\ref{prop:w10-integer-comparison}, contains the integer case as an instance. The
direction of the dependence is the one recorded here, and neither this appendix nor the body
statements that rest on it need the irrational dossier.
\end{remark}

\begin{proof}[Proof of Theorem~\ref{thm:D_index}]
Take $\beta\in(0,1)$. By Lemma~\ref{lem:D_telescope} and Lemma~\ref{lem:D_Bregimes},
\[
A(n)=\frac{B(n)}{n}-\sum_{k\ge n}\frac{B(k)}{k(k+1)}
=c_\beta n^{-\beta}-c_\beta\sum_{k\ge n}k^{-1-\beta}+o(n^{-\beta})
=c_\beta n^{-\beta}\Big(1-\frac1\beta\Big)+o(n^{-\beta}),
\]
and $c_\beta(\beta-1)/\beta=1/g^{*}(\beta)$ by the closed form of the transform. This is
transparency on $(0,1)$. For $\beta<0$ the series in \eqref{eq:D_abel} diverges and the same
computation is carried on the partial sums, $\sum_{k<n}B(k)/(k(k+1))=c_\beta n^{-\beta}/(-\beta)
+o(n^{-\beta})$, which combined with $B(n)/n=c_\beta n^{-\beta}+o(n^{-\beta})$ gives the same
constant. At $\beta=0$ the two terms of \eqref{eq:D_abel} behave differently, $B(n)/n$ being
$\mathcal{O}(\log n/n)$ by the digital sum bound while the series converges, so
$A(n)\to T_0=1/g^{*}(0)$ by Lemma~\ref{lem:D_telescope}, which is transparency at the
exponent $0$. At $\beta=1$ the exact
value $A_1(n)=\lambda^{-\lfloor\log_\lambda n\rfloor}$ is Lemma~\ref{lem:D_critical}, so
$nA_1(n)$ oscillates between $1$ and $\lambda$. For $\beta>1$ both $B(n)$ and the tail are
$\mathcal{O}_\beta(1/n)$ after division by $n$, so $A(n)=\mathcal{O}_\beta(n^{-1})$.

Sharpness above $1$ is read off \eqref{eq:D_Arec}, which gives for the sequence $mA(m)$
\[
\lambda m\,A(\lambda m)=mA(m)+\lambda m^{1-\beta}\bigl(\lambda^{-\beta}-1\bigr)+\lambda m\rho(m),
\qquad A(1)=1 .
\]
For $\beta>1$ one has $\rho(m)=\mathcal{O}_\beta(m^{-\beta})$, so both increments are
$\mathcal{O}_\beta(m^{1-\beta})$ and the series over $m=\lambda^{i}$ converges. Hence
\[
L(\beta):=\lim_{s\to\infty}\lambda^{s}A(\lambda^{s})
=1+\sum_{i\ge0}\Big[\lambda\,\lambda^{i(1-\beta)}\bigl(\lambda^{-\beta}-1\bigr)
+\lambda^{i+1}\rho(\lambda^{i})\Big]
\]
exists, and the same estimates show that the series converges uniformly on compact subsets of
$\{\Re\beta>1\}$, so $L$ is holomorphic there. As $\beta\to+\infty$ along the real axis, the
index $i=0$ contributes $\lambda(\lambda^{-\beta}-1)+\lambda\rho(1)$ with
$\rho(1)=\sum_{j<\lambda}j^{-\beta}/(j+1)\to\tfrac12$, and every index $i\ge1$ contributes a
quantity bounded by a constant multiple of $\lambda^{2i}\lambda^{-i\beta}$, which tends to
$0$. Hence $L(\beta)\to1-\lambda/2$. For $\lambda\ge3$ that limit is at most $-\tfrac12$, so
$L$ is not identically zero. For $\lambda=2$ the limit vanishes and the next order is needed,
the coefficient of $2^{-\beta}$ being $2$ from the index $i=0$ and $-4+\tfrac43$ from the
index $i=1$, so that $L(\beta)=-\tfrac23\,2^{-\beta}+\mathcal{O}(4^{-\beta})$, again not
identically zero. In both cases $L$ is holomorphic and not identically zero on a connected
open set, so its zeros are isolated and there are real exponents $\beta_j\downarrow1$ with
$L(\beta_j)\neq0$. At such an exponent $\lambda^{s}A(\lambda^{s})$ tends to a nonzero limit,
so $A(n)$ is not $o(n^{-\beta_j})$ and transparency fails. The supremum of
Definition~\ref{def:reg_index_fgv} is therefore $1$, and $\alpha(g)=1=\eta(g)$.
\end{proof}

\begin{numobs}\label{numobs:D_check}
For $\lambda\in\{2,3,5\}$ the recurrence \eqref{eq:D_Bn-recurrence} reproduces forward
substitution in the defining equation to $10^{-12}$ over $n\le300$ and five exponents. The
truncation of $T$ at $4\cdot10^{6}$ equals $+3.382\cdot10^{-1}$ at $\beta=1/10$ and
$+5.367\cdot10^{-3}$ at $\beta=2/5$ for $\lambda=2$, against the tail values
$-c_\beta N^{-\beta}/\beta$ predicted by Lemma~\ref{lem:D_Bregimes}, namely $0.339$ and
$5.376\cdot10^{-3}$, so the limit is zero to three digits in each case. The dilative
recurrence \eqref{eq:D_Arec} holds to $10^{-16}$ over $m\le3000$ for the three values of
$\lambda$. At $\beta=2$ and $\lambda=2$ the limit $L(\beta)=-0.15650828$ agrees with
$nA(n)$ at $n=2^{20}$ to seven decimals, and the computed values of $L$ at
$\beta\in\{6/5,8/5,2,3,6\}$ are nonzero for each $\lambda$ tested. At $\beta=0$ the partial
sums converge to $1.38629366$, $1.64791804$ and $2.01179316$ at $n=3\cdot10^{6}$ for
$\lambda=2,3,5$, against $\lambda\log\lambda/(\lambda-1)$ equal to $1.38629436$,
$1.64791843$ and $2.01179739$.
\end{numobs}

\galleryentry{E}{The folded affine profile}
 {$g(x)=\bigl|x-\tfrac12\bigr|+\tfrac12$ on $(0,1]$}
 {function of good variation}
 {$g^{*}(z)=\frac{2^{z}-z}{1-z}$, zeros the solutions of $2^{z}=z$, first pair
 $z_0=\eta+i\omega$}
 {$\tau(g)=\alpha(g)=\eta(g)=\eta$}
 {transparency for every $\beta<\eta$, absorption for $\beta\ge\eta$, and failures of
 transparency arbitrarily close above $\eta$, all proved in Theorem~\ref{thm:E_index}, and the
 sharpness follows from the nonremovable pole of the meromorphic connection coefficient}

\rafgalleryfig{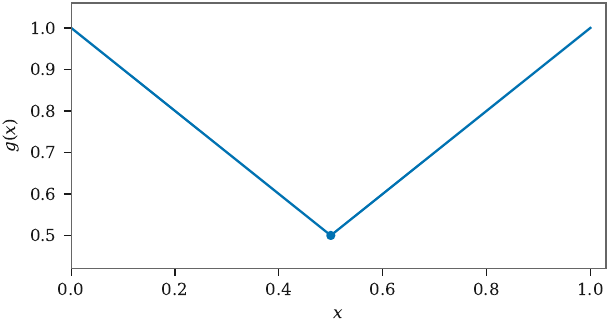}{The folded profile $g(x)=\bigl|x-\tfrac12\bigr|+\tfrac12$,
decreasing to one half at the fold and increasing back to one. The fold produces the first
complex characteristic pair, and Theorem~\ref{thm:E_index} proves that its real part is the
regularity index.}{fig:app_E}

Consider the piecewise linear function $g:(0,1]\to\R$ given by
\[
g(x) = \Bigl|x-\tfrac12\Bigr|+\tfrac12=
\begin{cases}
1-x, & 0<x<\tfrac12,\\[1mm]
x, & \tfrac12\le x\le1,
\end{cases}
\]
bounded with $g(1)=1$, of arithmetic Mellin transform
\[g^{*}(z)=\frac{2^{z}-z}{1-z}.
\]

The defining equation collapses to a recurrence in consecutive ranks carrying one dilated term.

\begin{proposition}\label{prop:E_exact}
Put $m_n=\lfloor(n-1)/2\rfloor$ and $d_\beta(n)=n^{1-\beta}-(n-1)^{1-\beta}$ for $n\ge2$,
with $d_\beta(1)=1$. Then for every $n\ge2$,
\begin{equation}\label{eq:E_recurrence}
A(n)=A(n-1)+\frac{d_\beta(n)-A(m_n)}{n},
\end{equation}
with $A(1)=1$, and with the convention $A(0)=0$ the same relation holds at $n=1$.
\end{proposition}

\begin{proof}
The condition $k/n<\tfrac12$ is $2k<n$, that is $k\le m_n$, so with $B(j)=\sum_{k\le j}ka_k$
the defining relation reads
\[
A(m_n)-\frac{2}{n}B(m_n)+\frac{1}{n}B(n)=n^{-\beta},
\]
the first two terms coming from the branch $1-k/n$ and the third from the branch $k/n$.
Multiplying by $n$ and differencing between ranks $n$ and $n-1$ splits by parity. For even
$n=2p$ one has $m_n=m_{n-1}=p-1$, the terms in $B(m)$ cancel and what remains is
$A(m_n)+na_n=d_\beta(n)$. For odd $n=2p+1$ one has $m_n=p$ and $m_{n-1}=p-1$, so
$nA(p)-(n-1)A(p-1)=A(p-1)+na_p$ and $B(p)-B(p-1)=pa_p$, and since $n-2p=1$ the terms in $a_p$
combine to $a_p$, leaving $A(p-1)+a_p+na_n=d_\beta(n)$, which is again
$A(m_n)+na_n=d_\beta(n)$. Both parities give \eqref{eq:E_recurrence}.
\end{proof}

The recurrence is of first order in consecutive ranks and carries one dilated term. It
determines the sequence from $A(1)$ alone, so the space of homogeneous solutions has dimension
one, where the continuous delay equation obtained from it by summation has a larger solution space.

\begin{proposition}\label{prop:E_zeros}
The zeros of $g^{*}$ are the solutions of $2^{z}=z$. Every solution $z=\sigma+it$ has
$\sigma>0$, and those with $t>0$ are exactly the pairs satisfying
\begin{equation}\label{eq:E_param}
t^{2}=4^{\sigma}-\sigma^{2},
\qquad
t\log2=\arctan\frac{t}{\sigma}+2k\pi,
\qquad k\in\mathbb Z_{\ge0},
\end{equation}
together with their conjugates. There is exactly one solution for each $k$, and their real parts
increase strictly with $k$. If $z_0=\eta+i\omega$ is the solution for $k=0$, then
$\eta<7/8$ and every other upper-half-plane solution has real part greater than one. Consequently
the only zeros of $\Delta(s)=s-1+2^{1-s}$ in $\Re s>0$ are
$s_0=1-z_0$ and $\bar s_0$.
\end{proposition}

\begin{proof}
Taking moduli and arguments in $2^z=z$ gives \eqref{eq:E_param}, and a solution with
$\sigma\le0$ would have $|z|\le1$ and an argument of modulus below $\pi/2$, a contradiction.
For $\sigma>0$ put
\[
t(\sigma)=\sqrt{4^\sigma-\sigma^2},\qquad
f(\sigma)=t(\sigma)\log2-\arctan\frac{t(\sigma)}\sigma .
\]
The inequality $2^\sigma>\sigma$ makes $t$ real. Direct differentiation gives
\[
f'(\sigma)=\frac{4^\sigma}{t(\sigma)}
 \left(\log2-\frac{\sigma}{4^\sigma}\right)^2
 +\frac{t(\sigma)}{4^\sigma}>0.
\]
Since $f(0+)=\log2-\pi/2<0$ and $f(\sigma)\to\infty$, each level $2k\pi$ is met exactly once.
Also $f(1)<2<2\pi$, so the branches $k\ge1$ have real part greater than one.

It remains to locate the first branch without decimals. At $\sigma=7/8$ let
$t^2=2^{7/4}-49/64$. The exact inequalities
$(5/3)^4<8<(17/10)^4$ give $8/5<t<13/8$, while
$\log2>56/81>69/100$. Hence $t\log2>138/125$. With
$y=(7/8)/t>7/13$, the alternating series for $\arctan y$ and $\pi<22/7$ give
\[
\arctan\frac{t}{7/8}=\frac\pi2-\arctan y
<\frac{11}{7}-\left(\frac7{13}-\frac1{3}\left(\frac7{13}\right)^3\right)
<\frac{138}{125}.
\]
Thus $f(7/8)>0$, and strict monotonicity gives $\eta<7/8$. Finally
$s=1-z$ carries precisely the first conjugate pair into $\Re s>0$.
\end{proof}

Powers are eigenfunctions of that recurrence up to an error, which gives the indicial equation
and the continuous reading.

\begin{proposition}\label{prop:E_indicial}
Let $L[A](n)=n\bigl(A(n)-A(n-1)\bigr)+A(m_n)$, the operator of \eqref{eq:E_recurrence}. For
every complex $\gamma$,
\begin{equation}\label{eq:E_indicial}
L\bigl[n^{-\gamma}\bigr]=\bigl(2^{\gamma}-\gamma\bigr)n^{-\gamma}
+\mathcal O_\gamma\bigl(n^{-\gamma-1}\bigr),
\end{equation}
so the indicial equation\index[terms]{indicial equation} of the recurrence is $2^{\gamma}=\gamma$, the zero set of the
numerator of $g^{*}$. In the continuous variable the same equation is $xA'(x)=-A(x/2)$, which
under $x=e^{u}$ becomes the constant delay equation
\begin{equation}\label{eq:E_delay}
y'(u)=-y(u-\log2),
\end{equation}
of characteristic equation $2^{\gamma}=\gamma$ under $y=e^{-\gamma u}$.
\end{proposition}

\begin{proof}
Expanding, $n(n^{-\gamma}-(n-1)^{-\gamma})=-\gamma n^{-\gamma}+\mathcal O_\gamma(n^{-\gamma-1})$,
and $m_n=\tfrac n2\bigl(1-(1+2\theta_n)/n\bigr)$ with $\theta_n\in\{0,\tfrac12\}$ gives
$m_n^{-\gamma}=2^{\gamma}n^{-\gamma}+\mathcal O_\gamma(n^{-\gamma-1})$, whence
\eqref{eq:E_indicial}. The continuous form follows from summing \eqref{eq:E_recurrence} and
replacing the sum over the ranks by an integral, the dilated argument becoming a constant
shift in the logarithmic variable.
\end{proof}

The forcing matches the same way. Since $d_\beta(n)=(1-\beta)n^{-\beta}+\mathcal
O_\beta(n^{-\beta-1})$, the profile $c_\beta n^{-\beta}$ with
\begin{equation}\label{eq:E_constant}
c_\beta=\frac{1-\beta}{2^{\beta}-\beta}=\frac{1}{g^{*}(\beta)}
\end{equation}
satisfies $L[c_\beta n^{-\beta}]=d_\beta(n)+\mathcal O_\beta(n^{-\beta-1})$, so the
transparency constant is fixed by the same computation that fixes the indicial equation.

\begin{proposition}\label{prop:E_below}
For every $\beta<-1$ one has $A(n)=c_\beta n^{-\beta}+o(n^{-\beta})$.
\end{proposition}

\begin{proofstatus}{The proof below gives an error $o(n^{-\beta})$ for every
$\beta<-1$. Its intermediate bound is $\mathcal O(n+n^{-\beta-1})$ when
$\beta\ne-2$, and $\mathcal O(n\log n)$ at $\beta=-2$.
The endpoint and the interval $[-1,-\gamma_0)$ are recovered by
Proposition~\ref{prop:E_negative_band}. The remaining band
$[-\gamma_0,0)$ is settled later by Proposition~\ref{prop:E_negative_band_complete},
using a Dirichlet-series argument for the residual. Thus Proposition~\ref{prop:E_below}
supplies precisely the range $\beta<-1$ in the negative-exponent clause of
Corollary~\ref{cor:E_regimes}, and no later proof depends on its rough error bound.}
\end{proofstatus}

\begin{proof}
Set $E(n)=A(n)-c_\beta n^{-\beta}$, so that
$L[E](n)=\mathcal O_\beta(n^{-\beta-1})$ by the preceding paragraph.
With $M(n)=\max_{k\le n}|E(k)|$, the recurrence gives
\[
 M(n)\le\left(1+\frac1n\right)M(n-1)+C_\beta n^{-\beta-2}.
\]
Dividing by $n+1$ and summing yields
\[
 M(n)\ll_\beta n\left(1+\sum_{j=2}^n j^{-\beta-3}\right)
 \ll_\beta
 \begin{cases}
 n,&-2<\beta<-1,\\
 n\log n,&\beta=-2,\\
 n^{-\beta-1},&\beta<-2.
 \end{cases}
\]
Every branch is $o(n^{-\beta})$ in its stated range, proving transparency
for all $\beta<-1$.
\end{proof}

\subsection*{A priori growth and the arithmetic Dirichlet series}

The passage to the regime $\beta\ge0$ is the content of the rest of the appendix. The route
is the one announced by the indicial equation itself. The partial sums are encoded in a
Dirichlet series, the recurrence becomes a functional equation\index[terms]{functional equation} whose entire denominator is
$\Delta(s)=s-1+2^{1-s}$, the zeros of $\Delta$ are the points $1-z$ with $2^{z}=z$, and the
asymptotics of $A$ are read off the poles by a smoothed Perron summation. Throughout
this appendix put
\[
\gamma_0=\frac{W_0(\log2)}{\log2}=0.641185744\ldots,
\qquad\text{the unique real solution of }\gamma_0=2^{-\gamma_0},
\]
and fix the notation $s_0=1-z_0$, so that $\Re s_0=1-\eta(g)=0.175321\ldots$ and
$\Delta(s_0)=\Delta(\bar s_0)=0$.

Before any asymptotic, the solution needs an a priori bound, valid for every forcing in the
right half plane.

\begin{proposition}\label{prop:E_growth}
For every $\eps>0$ and every complex $\beta$ with $\Re\beta\ge0$,
\[
|A_\beta(n)|\le C_{\eps,\beta}\,n^{\gamma_0+\eps},
\]
and the constant is bounded on compact sets of $\{\Re\beta\ge0\}$.
\end{proposition}

\begin{proof}
Write $\gamma=\gamma_0+\eps$ and note $|d_\beta(n)|\le C(1+|\beta|)n^{-\Re\beta}\le
C(1+|\beta|)$ for $n\ge1$. From \eqref{eq:E_recurrence},
$|A(n)|\le|A(n-1)|+\bigl(C(1+|\beta|)+|A(m_n)|\bigr)/n$. Assume inductively $|A(k)|\le Kk^{\gamma}$
for $k<n$. Since $m_n\le n/2$,
\[
|A(n)|\le K(n-1)^{\gamma}+\frac{C(1+|\beta|)}{n}+\frac{K}{n}\Bigl(\frac n2\Bigr)^{\gamma}
=Kn^{\gamma}\Bigl[\Bigl(1-\tfrac1n\Bigr)^{\gamma}+\frac{2^{-\gamma}}{n}\Bigr]
+\frac{C(1+|\beta|)}{n}.
\]
The bracket is $1-(\gamma-2^{-\gamma})/n+\mathcal O(n^{-2})$, and $\gamma-2^{-\gamma}>0$
because $\gamma>\gamma_0$ and $u\mapsto u-2^{-u}$ is increasing with root $\gamma_0$. Hence
for $n\ge n_0(\eps)$ the bracket is at most $1-\delta_\eps/n$ with $\delta_\eps>0$, and the
induction closes as soon as $K\delta_\eps n^{\gamma-1}\ge C(1+|\beta|)/n$, which holds for
every $n\ge n_0$ since $\gamma>0$. Choosing $K$ large enough to cover the ranks below $n_0$
gives the bound, uniformly for $\beta$ in a compact set.
\end{proof}

For the homogeneous solution the same induction closes with an explicit exponent and no
constant at all, which is worth recording once.

\begin{lemma}\label{lem:E_explicit}
Let $H$ be the solution of the homogeneous recurrence $H(n)=H(n-1)-H(m_n)/n$ with $H(1)=1$.
Then $|H(n)|\le n^{13/20}$ for every $n\ge1$.
\end{lemma}

\begin{proof}
Write $\nu=\tfrac{13}{20}$ and note $\nu-2^{-\nu}>\tfrac1{100}$. The claim holds with equality
at $n=1$, and $|H(2)|=1\le2^{\nu}$. For $n\ge3$, since $m_n\le n/2$ and
$(1-1/n)^{\nu}\le1-\nu/n$ by concavity,
\[
|H(n)|\le|H(n-1)|+\frac{|H(m_n)|}{n}
\le n^{\nu}\Bigl[\Bigl(1-\frac1n\Bigr)^{\nu}+\frac{2^{-\nu}}{n}\Bigr]
\le n^{\nu}\Bigl[1-\frac{\nu-2^{-\nu}}{n}\Bigr]<n^{\nu}. \qedhere
\]
\end{proof}

The same induction, run against an arbitrary bounded second member rather than against the
forcing itself, carries the transparent range of Proposition~\ref{prop:E_below} past its
endpoint.

\begin{lemma}\label{lem:E_growth_forced}
With the convention $X(0)=0$, let $(X(n))_{n\ge1}$ satisfy
$L[X](n)=\varphi(n)$ for $n\ge2$, with
$\Phi_0=\sup_{n\ge2}|\varphi(n)|$ finite. Then for every $\eps>0$,
\[
|X(n)|\le C_\eps\bigl(|X(1)|+\Phi_0\bigr)\,n^{\gamma_0+\eps}\qquad(n\ge1).
\]
\end{lemma}

\begin{proof}
Write $\gamma=\gamma_0+\eps$ and
$\delta_\eps=\tfrac12(\gamma-2^{-\gamma})>0$. Since
\[
\Bigl(1-\frac1n\Bigr)^\gamma+\frac{2^{-\gamma}}n
=1-\frac{\gamma-2^{-\gamma}}n+\mathcal O_\gamma(n^{-2}),
\]
there is an $n_0=n_0(\eps)$ such that the left-hand side is at most
$1-\delta_\eps/n$ for $n\ge n_0$. Put
$M_j=\max_{1\le k\le j}|X(k)|$. For $2\le j<n_0$, the recurrence and $m_j<j$ give
\[
M_j\le\Bigl(1+\frac1j\Bigr)M_{j-1}+\frac{\Phi_0}{j},
\]
and iteration over these finitely many ranks yields
$M_{n_0-1}\le C_\eps\bigl(|X(1)|+\Phi_0\bigr)$.

Choose $K$ at least
\[
\max\left\{\frac{\Phi_0}{\delta_\eps},
\max_{1\le k<n_0}\frac{|X(k)|}{k^\gamma}\right\}.
\]
Assume $|X(k)|\le Kk^\gamma$ for every $k<n$, where $n\ge n_0$. Since $m_n\le n/2$,
\begin{align*}
|X(n)|
&\le |X(n-1)|+\frac{\Phi_0+|X(m_n)|}{n}\\
&\le Kn^\gamma\left[\Bigl(1-\frac1n\Bigr)^\gamma+\frac{2^{-\gamma}}n\right]
   +\frac{\Phi_0}{n}\\
&\le Kn^\gamma\Bigl(1-\frac{\delta_\eps}{n}\Bigr)+\frac{\Phi_0}{n}
\le Kn^\gamma,
\end{align*}
because $K\delta_\eps\ge\Phi_0$ and $n^\gamma\ge1$. Thus the induction closes, and the
preceding finite-rank estimate permits $K\le C_\eps(|X(1)|+\Phi_0)$.
\end{proof}

\begin{proposition}\label{prop:E_negative_band}
For every $\beta\in[-1,-\gamma_0)$ and every $\eps>0$,
\[
A(n)=c_\beta\,n^{-\beta}+\mathcal O_{\beta,\eps}\bigl(n^{\gamma_0+\eps}\bigr).
\]
In particular, on choosing $0<\eps<-\beta-\gamma_0$,
\[
A(n)=\frac{n^{-\beta}}{g^{*}(\beta)}+o(n^{-\beta}),
\]
so every such $\beta$ is transparent for the folded affine profile.
Together with Proposition~\ref{prop:E_below} every exponent below $-\gamma_0$ is therefore
transparent.
\end{proposition}

\begin{proof}
Put $P(k)=c_\beta k^{-\beta}$ with $P(0)=0$ and $E=A-P$, so that $L[E]=r$ with
$r(n)=d_\beta(n)-L[P](n)$. The three expansions
$d_\beta(n)=(1-\beta)n^{-\beta}+\mathcal O_\beta(n^{-\beta-1})$,
$n\bigl(n^{-\beta}-(n-1)^{-\beta}\bigr)=-\beta n^{-\beta}+\mathcal O_\beta(n^{-\beta-1})$ and
$m_n^{-\beta}=2^{\beta}n^{-\beta}\bigl(1+\mathcal O(1/n)\bigr)$, all valid from the rank $3$
on, combine through the identity $(1-\beta)-c_\beta(2^{\beta}-\beta)=0$ of
\eqref{eq:E_constant} and leave $|r(n)|\le C_\beta n^{-\beta-1}$, a bounded quantity exactly
when $\beta\ge-1$. The rank $2$ contributes a constant. Lemma~\ref{lem:E_growth_forced} applied
to $E$, whose value at the first rank is $1-c_\beta$, gives
$E(n)=\mathcal O_{\beta,\eps}(n^{\gamma_0+\eps})$ for every $\eps>0$. Since
$-\beta>\gamma_0$, choosing $0<\eps<-\beta-\gamma_0$ makes this remainder $o(n^{-\beta})$.
\end{proof}

The Dirichlet series of the solution satisfies a functional equation, which continues it past
its abscissa of convergence and locates its poles.

\begin{proposition}\label{prop:E_dirichlet}
Let $\beta$ be complex with $\Re\beta\ge0$ and let
\[
D(s)=\sum_{n\ge1}A(n)\,n^{-s},
\]
absolutely convergent on $\Re s>1+\gamma_0$. Put $\Delta(s)=s-1+2^{1-s}$ and
\begin{align*}
\rho_s(n)&=\int_0^1\bigl[(n+u)^{-s}-n^{-s}+su\,n^{-s-1}\bigr]\,du,\notag\\
\tau_s(m)&=(2m+1)^{-s}+(2m+2)^{-s}-2\,(2m)^{-s}+3s\,(2m)^{-s-1},\notag\\
R_2(s)&=\sum_{n\ge1}A(n)\rho_s(n),
\qquad
T_2(s)=\sum_{m\ge1}A(m)\tau_s(m),\notag\\
\Phi(s)&=\sum_{n\ge1}d_\beta(n)\,n^{-s}
=(1-\beta)\,\zeta(s+\beta)+\Psi_\beta(s),
\qquad
\Psi_\beta(s)=\sum_{n\ge1}\bigl[d_\beta(n)-(1-\beta)n^{-\beta}\bigr]n^{-s}.\notag
\end{align*}
Then $|\rho_s(n)|\le\tfrac16|s(s+1)|\,n^{-\Re s-2}$ and
$|\tau_s(m)|\le\tfrac52|s(s+1)|\,(2m)^{-\Re s-2}$
for $\Re s\ge-2$, the series $R_2$ and $T_2$ converge absolutely on
$\{\Re s>\gamma_0-1\}$, the series $\Psi_\beta$ converges absolutely on
$\{\Re s>-\Re\beta\}$, and on $\Re s>1+\gamma_0$
\begin{equation}\label{eq:E_funceq}
\Delta(s)\,D(s)
=\Phi(s)+s\Bigl[\frac{s-1}{2}+3\cdot2^{-s-1}\Bigr]D(s+1)-(s-1)R_2(s)-T_2(s)
=:N(s).
\end{equation}
Consequently $D$ extends meromorphically to $\{\Re s>\gamma_0-1\}$. Away from the
collisions $\beta=z_0,\bar z_0$, its only possible poles in $\{\Re s>0\}$ are simple,
at $s=1-\beta$ with residue $c_\beta$ and at $s_0,\bar s_0$ with residues
\begin{equation}\label{eq:E_kappa}
\kappa(\beta)=\frac{N_\beta(s_0)}{\Delta'(s_0)},
\qquad
\overline{\kappa(\bar\beta)}\ \text{at }\bar s_0,
\qquad
\Delta'(s_0)=1-z_0\log2\neq0 .
\end{equation}
At $\beta=z_0$ the forcing pole and the zero $s_0$ collide and $D$ has a double pole at
$s_0$, and the conjugate statement holds at $\beta=\bar z_0$.
\end{proposition}

\begin{proof}
The recurrence in the form $n\bigl(A(n)-A(n-1)\bigr)+A(m_n)=d_\beta(n)$, valid for every
$n\ge1$ with $A(0)=0$, is multiplied by $n^{-s}$ and summed. Abel summation and the vanishing
of the boundary term give, for $\Re s>1+\gamma_0$,
\[
\sum_{n\ge1}\bigl(A(n)-A(n-1)\bigr)n^{1-s}
=\sum_{n\ge1}A(n)\bigl[n^{1-s}-(n+1)^{1-s}\bigr]
=(s-1)\sum_{n\ge1}A(n)\int_0^1(n+u)^{-s}\,du,
\]
while the dilated term regroups along $m_n=m$ for $n\in\{2m+1,2m+2\}$ into
$\sum_{m\ge1}A(m)\bigl[(2m+1)^{-s}+(2m+2)^{-s}\bigr]$. Writing
$\int_0^1(n+u)^{-s}\,du=n^{-s}-\tfrac s2n^{-s-1}+\rho_s(n)$ and
$(2m+1)^{-s}+(2m+2)^{-s}=2^{1-s}m^{-s}-3s(2m)^{-s-1}+\tau_s(m)$ and collecting the terms
proportional to $D(s)$ and to $D(s+1)$ yields \eqref{eq:E_funceq}. The two pointwise bounds
follow from Taylor expansion\index[terms]{Taylor expansion}\index[names]{Taylor, B.} with integral remainder. For $\rho_s$ the integrand is
$s(s+1)\int_0^u(u-v)(n+v)^{-s-2}\,dv$ and $(n+v)^{-\Re s-2}\le n^{-\Re s-2}$ for $\Re s\ge-2$,
so that integrating in $u$ over $[0,1]$ produces the factor $\int_0^1\tfrac{u^2}{2}\,du=\tfrac16$.
For $\tau_s$ the same expansion at base point $2m$ with increments $1$ and $2$ leaves the
remainder $s(s+1)\sum_{j=1}^{2}\int_0^j(j-v)(2m+v)^{-s-2}\,dv$, of modulus at most
$\tfrac52|s(s+1)|(2m)^{-\Re s-2}$. With Proposition~\ref{prop:E_growth} the series $R_2$ and
$T_2$ converge absolutely and locally uniformly on $\Re s>\gamma_0-1$. For $\Psi_\beta$,
$d_\beta(n)=(1-\beta)\int_0^1(n-u)^{-\beta}\,du$ for $n\ge2$ gives
$|d_\beta(n)-(1-\beta)n^{-\beta}|\le|\beta(1-\beta)|(n-1)^{-\Re\beta-1}$, which is summable
against $n^{-\Re s}$ on $\Re s>-\Re\beta$.

Continuation is by descent through two strips of width one. On
$\{\Re s>\gamma_0\}$ every term on the right of \eqref{eq:E_funceq} is defined and
meromorphic, $D(s+1)$ living in the initial half plane, so $D=N/\Delta$ extends there. One
more application, now with $D(s+1)$ given by the first extension, reaches
$\{\Re s>\gamma_0-1\}$. The possible poles of $N$ on $\{\Re s>0\}$ are the forcing pole
at $s=1-\beta$ and those inherited from $D(s+1)$, the latter having real part at most zero.
The zeros of $\Delta$ with positive real part are exactly $s_0$ and $\bar s_0$ by
Proposition~\ref{prop:E_zeros}, and they are simple. If the forcing pole is distinct from these
zeros, its residue is
$(1-\beta)/\Delta(1-\beta)=c_\beta$, while the residue at $s_0$ is
$N_\beta(s_0)/\Delta'(s_0)$. If $\beta=z_0$, however, the term
$(1-\beta)\zeta(s+\beta)$ has a simple pole at $s=s_0$ with nonzero residue $1-z_0$,
division by the simple zero of $\Delta$ produces a double pole. Conjugation gives the statements
at $\bar s_0$ and for $\bar\beta$.
\end{proof}

Shifting the contour asks for growth on vertical lines.

\begin{lemma}\label{lem:E_lines}
Fix $\theta\in(0,\tfrac1{20}]$ and $\beta$ with $\Re\beta\ge0$. On the region
$\{\theta\le\Re s\le3\}$ deprived of discs of radius $\tfrac1{10}$ around $1-\beta$, $s_0$
and $\bar s_0$,
\[
|D(s)|\le C_{\theta,\beta}\,(1+|t|)^{5},\qquad s=\sigma+it .
\]
\end{lemma}

\begin{proof}
Partial summation at length $|t|+2$ gives the elementary bound
$|\zeta(\sigma'+it)|\le C_\theta(1+|t|)$ for $\sigma'\ge\theta$ and $|\sigma'+it-1|\ge\tfrac1{10}$,
hence $|\Phi(s)|\le C(1+|t|)$ on the region, the series $\Psi_\beta$ being bounded there. The
bounds of Proposition~\ref{prop:E_dirichlet} give $|R_2(s)|+|T_2(s)|\le C(1+|t|)^{2}$ on
$\Re s\ge\gamma_0-1+\tfrac1{100}$. On $\Re s\ge1+\gamma_0+\tfrac1{100}$ the series for $D$
converges absolutely, so $|D|\le C$ there. Descending, on
$\gamma_0+\tfrac1{100}\le\Re s\le1+\gamma_0+\tfrac1{100}$ the identity \eqref{eq:E_funceq}
gives
\[
|D(s)|\le\frac{C(1+|t|)+C(1+|t|)^{2}\cdot C+C(1+|t|)\,(1+|t|)^{2}}{|\Delta(s)|},
\]
and $|\Delta(s)|\ge|t|-2-2^{1-\theta}\ge\tfrac12|t|$ for $|t|\ge T_0$, while for $|t|\le T_0$
the function $1/\Delta$ is bounded on the region, the zeros of $\Delta$ being excluded by the
discs. Hence $|D|\le C(1+|t|)^{3}$ on that strip. One more descent, with $D(s+1)$ now bounded
by $C(1+|t|)^{3}$, gives $|D|\le C(1+|t|)^{5}$ on $\theta\le\Re s\le\gamma_0+\tfrac1{100}$,
and the three ranges assemble into the claim.
\end{proof}

\subsection*{The asymptotic expansion}

The expansion of the solution follows, with the two conjugate modes and a remainder.

\begin{theorem}\label{thm:E_asymptotic}
Let $\beta\ge0$ be real. With $c_\beta$ as in \eqref{eq:E_constant} and $\kappa(\beta)$ as in
\eqref{eq:E_kappa},
\begin{equation}\label{eq:E_master}
A(n)=c_\beta\,n^{-\beta}+2\,\Re\bigl(\kappa(\beta)\,n^{-z_0}\bigr)
+\mathcal O_\beta\bigl(n^{-7/8}\bigr).
\end{equation}
\end{theorem}

\begin{proof}
Write $P(n)=c_\beta n^{-\beta}$, $M(n)=\kappa n^{-z_0}+\bar\kappa n^{-\bar z_0}$ with
$\kappa=\kappa(\beta)$, and let $r(n)=A(n)-P(n)-M(n)$ be the residual, whose Dirichlet
series is
\[
D_r(s)=D(s)-c_\beta\,\zeta(s+\beta)-\kappa\,\zeta(s+z_0)-\bar\kappa\,\zeta(s+\bar z_0).
\]
By Proposition~\ref{prop:E_dirichlet} the poles of $D$ on $\{\Re s>0\}$ are cancelled by the
subtracted zeta terms, residue by residue, so $D_r$ is holomorphic on $\{\Re s>\theta\}$ for
any fixed $\theta\in(0,\tfrac1{25}]$ chosen with $|1-\beta-\theta|\ge\tfrac1{100}$, and by
Lemma~\ref{lem:E_lines} together with the zeta bound it satisfies
$|D_r(\sigma+it)|\le C(1+|t|)^{5}$ on $\theta\le\sigma\le3$.

Fix $r=15$ and, for $x\ge2$, the Riesz mean
\[
S(x)=\frac{1}{r!}\sum_{n\le x}r(n)\,(x-n)^{r}
=\frac{1}{2\pi i}\int_{(c)}D_r(s)\,x^{s+r}\,\frac{\Gamma(s)}{\Gamma(s+r+1)}\,ds,
\qquad c=2+\gamma_0,
\]
by Lemma~\ref{lem:perron-riesz}, since
$\Gamma(s)/\Gamma(s+r+1)=1/[s(s+1)\cdots(s+r)]$.
The a priori bound for the coefficients, with a fixed $\eps<1$, gives absolute
convergence on $c=2+\gamma_0$. The growth exponent is $M=5<r=15$.
The quotient has no pole between $\theta$ and $c$, since $D_r$ is holomorphic
there and $\theta>0$. Thus the shift formula of that lemma gives
\begin{equation}\label{eq:E_riesz}
S(x)=\mathcal O\bigl(x^{\theta+r}\bigr).
\end{equation}

Let $h=x^{1-\delta}$ with $\delta=\tfrac1{200}$, and let $\Delta_h^{r+1}$ denote the
$(r+1)$-fold forward difference with step $h$ in the variable $x$. For a function
$f\in C^{r+1}$ one has $\Delta_h^{r+1}f(x)=\int_{[0,h]^{r+1}}f^{(r+1)}(x+u_1+\cdots+u_{r+1})\,du$,
and applying this to the truncated power $u\mapsto(u-n)^{r}_{+}/r!$, whose $(r+1)$-st
distributional derivative is the Dirac mass at $n$, gives
\[
V(n):=\frac{1}{h^{r+1}}\,\Delta_h^{r+1}\Bigl[\frac{(x-n)_+^{r}}{r!}\Bigr]
=\frac{1}{h^{r+1}}\,\varphi_{r+1}(n-x)\ \ge0,
\]
where $\varphi_{r+1}$ is the $(r+1)$-fold convolution of the indicator of $[0,h]$ with
itself, of total mass $h^{r+1}$ and supported on $x<n<x+(r+1)h$. Applying $\Delta_h^{r+1}$ to the identity defining $S$ and
dividing by $h^{r+1}$,
\begin{equation}\label{eq:E_window}
\sum_{n}V(n)\,r(n)=\frac{\Delta_h^{r+1}S(x)}{h^{r+1}}
=\mathcal O\bigl(x^{\theta+r}h^{-r-1}\bigr)
=\mathcal O\bigl(x^{\theta-1+(r+1)\delta}\bigr),
\end{equation}
each of the $2^{r+1}$ terms in the difference obeying \eqref{eq:E_riesz}. The total mass
$\mu=\sum_nV(n)$ is evaluated by the same computation applied to the constant sequence.
Euler and Maclaurin give $\sum_{n\le x}(x-n)^{r}/r!=x^{r+1}/(r+1)!-x^{r}/(2\,r!)+\mathcal
O(x^{r-1})$, the difference operator annihilates polynomials of degree at most $r$ and maps
$x^{r+1}/(r+1)!$ to $h^{r+1}$ exactly, so $\mu=1+\mathcal O(x^{r-1}h^{-r-1})=1+\mathcal O(x^{-1})$.

The oscillation of $r$ across the window is controlled by the recurrence. The exact relation
$n\bigl(A(n)-A(n-1)\bigr)+A(m_n)=d_\beta(n)$, the matching
$L[P](n)=d_\beta(n)+\mathcal O_\beta(n^{-\beta-1})$ of \eqref{eq:E_constant}, and the indicial
identity \eqref{eq:E_indicial} at $\gamma=z_0$ and $\bar z_0$, where $2^{z_0}-z_0=0$, give
after subtraction
\[n\bigl(r(n)-r(n-1)\bigr)+r(m_n)=\mathcal O_\beta\bigl(n^{-\beta-1}+n^{-\eta-1}\bigr),
\]
where $\eta=\eta(g)$. The main and oscillating profiles thus cancel structurally inside the
increments, which is the point of the whole construction. Consequently, for every $k$ in the
window,
\begin{align*}
|r(k)-r(\lceil x\rceil)|
&\le\sum_{j}\Bigl[\frac{|r(m_j)|}{j}
       +C_\beta\bigl(j^{-\beta-2}+j^{-\eta-2}\bigr)\Bigr]\\
&\le C\Bigl[x^{-\delta}\sup_{y\in[\frac x2-1,\,\frac{11x}{20}]}|r(y)|
+x^{-\beta-1-\delta}+x^{-\eta-1-\delta}\Bigr],
\end{align*}
the sum over $j$ having at most $(r+1)h$ terms. Combining with \eqref{eq:E_window} and
$\mu\ge\tfrac12$,
\begin{equation}\label{eq:E_bootstrap}
|r(n)|\le C_\beta\Bigl[x^{\theta-1+(r+1)\delta}
+x^{-\delta}\sup_{y\in[\frac x2-1,\,\frac{11x}{20}]}|r(y)|
+x^{-\beta-1-\delta}+x^{-\eta-1-\delta}\Bigr],
\qquad n=\lceil x\rceil .
\end{equation}
Choose $\theta=\theta(\beta)\in(0,1/25]$ with
$|1-\beta-\theta|\ge1/100$, as required above, and keep $r=15$, $\delta=1/200$.
Then $\theta-1+16\delta\le-0.88$.
For all sufficiently large $n$, take $x=n$ in \eqref{eq:E_bootstrap}.
The delayed indices then lie between $n/3$ and $3n/5$. If
$r(m)=\mathcal O(m^a)$ for some real $a$, their supremum is
$\mathcal O(n^a)$, so that inequality improves the bound to
\[
 r(n)=\mathcal O_\beta\bigl(n^{\max(-0.88,a-\delta)}\bigr).
\]
The two remaining forcing terms decay faster than $n^{-0.88}$ because
$\beta\ge0$ and $\eta>0$. Proposition~\ref{prop:E_growth}, together with
the boundedness of $P$ and $M$, supplies the initial exponent
$a=\gamma_0+\eps$. A fixed finite number of repetitions therefore gives
$r(n)=\mathcal O_\beta(n^{-0.88})$. Since $0.88>7/8$, this implies
\eqref{eq:E_master}.
\end{proof}

The missing negative band follows by applying the same smoothed argument to the residual rather
than to the original forcing.

\begin{proposition}\label{prop:E_negative_band_complete}
Every real $\beta$ with $-\gamma_0\le\beta<0$ is transparent:
\[
A(n)=c_\beta n^{-\beta}+o(n^{-\beta}).
\]
\end{proposition}

\begin{proof}
Put $q=-\beta$, so $0<q\le\gamma_0<1$, and set
$P(n)=c_\beta n^q$, $E(n)=A(n)-P(n)$. The calculation of
Proposition~\ref{prop:E_negative_band} gives
\[
L[E](n)=r(n)=\mathcal O_\beta(n^{q-1}).
\]
The forcing is bounded, and Lemma~\ref{lem:E_growth_forced} gives
$E(n)=\mathcal O(n^{\gamma_0+\varepsilon_0})$.

Let $F(s)=\sum E(n)n^{-s}$ and $R(s)=\sum r(n)n^{-s}$. Then $R$ converges normally for
$\Re s>q$. The Abel calculation in Proposition~\ref{prop:E_dirichlet}, applied to the residual,
gives
\[
\Delta(s)F(s)=R(s)+M(s)F(s+1)-(s-1)R_{2,E}(s)-T_{2,E}(s),
\]
where $M(s)=s[(s-1)/2+3\cdot2^{-s-1}]$ and the two weighted series converge normally for
$\Re s>\gamma_0-1+\varepsilon_0$. Choose
\[
\max(q,\Re s_0)<\theta<1,
\quad \gamma_0-1+\varepsilon_0<\theta,
\quad 0<16\delta<\min(1,1+q-\theta).
\]
The identity is used twice, each use gaining one unit to the left. Since
$E(n)=\mathcal O(n^{\gamma_0+\varepsilon_0})$, the series $F$ converges absolutely on
$\Re s>1+\gamma_0+\varepsilon_0$. Read as
\begin{equation}\label{eq:E_descent}
F(s)=\frac{R(s)+M(s)F(s+1)-(s-1)R_{2,E}(s)-T_{2,E}(s)}{\Delta(s)},
\end{equation}
the identity expresses $F$ at $s$ through $F$ at $s+1$. The first descent
reaches $\Re s>\max(q,\gamma_0+\varepsilon_0)$ and the second reaches
$\Re s>\max(q,\gamma_0-1+\varepsilon_0)$, meromorphically at possible zeros
of $\Delta$. The restriction $\Re s>q$ comes from $R$, and the other
restriction from $R_{2,E}$ and $T_{2,E}$. The choice of $\theta$ lies strictly
inside their common domain.

The division in \eqref{eq:E_descent} is legitimate on the whole region and costs nothing. By
Proposition~\ref{prop:E_zeros} the only zeros of $\Delta$ in $\Re s>0$ are $s_0$ and
$\bar s_0$, and $\theta>\Re s_0$, so $\Delta$ does not vanish on $\Re s\ge\theta$. Writing
$s=\sigma+it$ with $\sigma\ge\theta$, the bounds $|2^{1-s}|=2^{1-\sigma}\le2^{1-\theta}$ and
$|s-1|\ge|t|$ give
\[
|\Delta(s)|\ \ge\ |t|-2^{1-\theta}\ \ge\ \tfrac12|t|\ \ge\ \tfrac12T_0
\qquad\bigl(|t|\ge T_0:=2^{2-\theta}\bigr),
\]
while on the compact set $\{\theta\le\sigma\le3,\ |t|\le T_0\}$, which carries no zero, the
continuous function $1/\Delta$ is bounded. Hence $1/\Delta=\mathcal O(1)$ on
$\theta\le\Re s\le3$.

The polynomial bound follows by the same two steps. Taking the initial growth
bound with $\varepsilon_0/2$ makes the series $F$ bounded on
$\Re s\ge1+\gamma_0+\varepsilon_0$. On the region $\Re s\ge\theta>q$,
$|R(s)|=\mathcal O(1)$,
$|M(s)|=\mathcal O((1+|t|)^{2})$ and, by the weight estimates of
Proposition~\ref{prop:E_dirichlet},
$|R_{2,E}(s)|+|T_{2,E}(s)|=\mathcal O((1+|t|)^{2})$. The numerator of \eqref{eq:E_descent} is
then $\mathcal O((1+|t|)^{3})$, the dominant term being $(s-1)R_{2,E}(s)$, so the first descent
gives $F(s)=\mathcal O((1+|t|)^{3})$ on
$\gamma_0+\varepsilon_0\le\Re s\le1+\gamma_0+\varepsilon_0$. In the second descent
$M(s)F(s+1)=\mathcal O((1+|t|)^{5})$ dominates, and
$F(s)=\mathcal O\bigl((1+|\Im s|)^{5}\bigr)$ on $\theta\le\Re s\le3$.

For the Riesz mean of order $15$,
\[
S(x)=\frac1{15!}\sum_{n\le x}E(n)(x-n)^{15},
\]
Lemma~\ref{lem:perron-riesz}, with $M=5<15$ and initial line
$c=2+\gamma_0$ (choose $0<\varepsilon_0<1$), gives
$S(x)=\mathcal O(x^{\theta+15})$. No pole is crossed, because $F$ is
holomorphic on $\Re s\ge\theta>0$. With $h=x^{1-\delta}$, apply the forward difference of order $16$ and divide by
$h^{16}$. The resulting weights are the ones computed before \eqref{eq:E_window} at $r=15$,
namely $V(n)=h^{-16}\varphi_{16}(n-x)$ with $\varphi_{16}$ the sixteenfold convolution of the
indicator of $[0,h]$ with itself. They are therefore nonnegative, supported on $x<n<x+16h$,
and of mass $1+\mathcal O(x^{-1})$, the mass being evaluated there by Euler and Maclaurin on
the constant sequence. Thus
\[
\sum_nV(n)E(n)=\mathcal O(x^{\theta-1+16\delta}).
\]
The exact recurrence of the residual gives, throughout that window,
\[
|E(k)-E(\lceil x\rceil)|\le C\left(
x^{-\delta}\sup_{x/2-1\le y\le11x/20}|E(y)|+x^{q-1-\delta}\right).
\]
The same pointwise iteration as in Theorem~\ref{thm:E_asymptotic} now applies.
Put
\[
 p=\max(\theta-1+16\delta,q-1-\delta)<q.
\]
Nonnegative weights of mass tending to one convert the last two displays
into
\[
 |E(n)|\le C\left(n^p+n^{-\delta}
                 \sup_{n/3\le m\le3n/5}|E(m)|\right)
\]
for all sufficiently large $n$. A bound $E(m)=\mathcal O(m^a)$ is therefore
improved to
\[
 E(n)=\mathcal O(n^{\max(p,a-\delta)}).
\]
Starting with $a=\gamma_0+\varepsilon_0$ and repeating a fixed number
$k>(\gamma_0+\varepsilon_0-p)/\delta$ of times yields
$E(n)=\mathcal O(n^p)=o(n^q)$. This includes $q=\gamma_0$ and proves the claim.
\end{proof}

\begin{corollary}\label{cor:E_regimes}
For every real $\beta<\eta(g)$,
\[
A(n)=\frac{n^{-\beta}}{g^*(\beta)}+o(n^{-\beta}),
\]
and for every $\beta\ge\eta(g)$,
$A(n)=\mathcal O_\beta(n^{-\eta(g)})$.
\end{corollary}

\begin{proof}
Propositions~\ref{prop:E_below}, \ref{prop:E_negative_band}, and
\ref{prop:E_negative_band_complete} cover $\beta<0$. For
$0\le\beta<\eta$, Theorem~\ref{thm:E_asymptotic} and $\eta<7/8$ make both residual modes
$o(n^{-\beta})$. For $\beta\ge\eta$, every term of \eqref{eq:E_master} is
$\mathcal O(n^{-\eta})$.
\end{proof}

\subsection*{Sharpness and the connection coefficient}

The complex resonance that prevents global holomorphy supplies the missing nontriviality.

\begin{proposition}\label{prop:E_kappa}
The coefficient $\kappa(\beta)$ of \eqref{eq:E_kappa} is meromorphic on
$\{\Re\beta>0\}$ and holomorphic on the connected domain
$\Omega=\{\Re\beta>0\}\setminus\{z_0\}$. At $\beta=z_0$,
\[
\kappa(\beta)=\frac{1-z_0}{\Delta'(s_0)}\frac1{\beta-z_0}+\mathcal O(1).
\]
Consequently, for every $\varepsilon>0$ there is a real
$\beta\in(\eta,\eta+\varepsilon)$ with $\kappa(\beta)\neq0$, and transparency fails at every
such exponent.
\end{proposition}

\begin{proof}
In $N_\beta(s_0)$ the term $(1-\beta)\zeta(s_0+\beta)$ has, as $\beta\to z_0$, the
principal part $(1-z_0)/(\beta-z_0)$. All remaining terms are holomorphic there.
Indeed $\Psi_\beta(s_0)$, $R_{2,\beta}(s_0)$, and $T_{2,\beta}(s_0)$ converge normally on
compact subsets, while $D_\beta(s_0+1)$ is obtained from \eqref{eq:E_funceq} at $s_0+1$,
its denominator is $\Delta(s_0+1)=1-z_0/2\neq0$ and every series on the right converges
normally near $z_0$. Division by $\Delta'(s_0)\neq0$ gives the displayed principal part,
which is nonremovable.

Thus $\kappa$ is not identically zero on $\Omega$. If it vanished throughout a real interval
$(\eta,\eta+\varepsilon)$, the identity theorem on the connected domain $\Omega$ would make it
identically zero, contradicting the pole. Hence nonzero real points occur in every such interval.

At one of these points, write $z_0=\eta+i\omega$ and
$\kappa=|\kappa|e^{i\phi}$. The expansion \eqref{eq:E_master} gives
\[
n^\eta\left(A(n)-c_\beta n^{-\beta}\right)
=2|\kappa|\cos(\phi-\omega\log n)+o(1).
\]
Taking $n$ to be the nearest integer to
$\exp((\phi+2k\pi)/\omega)$ makes the cosine tend to one. The limsup is therefore
$2|\kappa|>0$, whereas transparency at $\beta>\eta$ would force $n^\eta A(n)\to0$.
\end{proof}

\subsection*{An independent mode functional}

At $\beta=1$ the differenced forcing vanishes, $d_1(n)=0$ for $n\ge2$, so
$H=A_{1}$ is the fundamental solution of the homogeneous recurrence with $H(1)=1$, and
$\kappa(1)$ is its mode coefficient. The following functional computes it as a limit of
finite windows, in exact analogy with the conserved spectral projections of the delay
equation \eqref{eq:E_delay}.

Several forcing exponents are compared in what follows, so the subscript is kept on $A_\beta$ and on the sequences attached to it. Elsewhere in the volume the exponent is fixed by its context and the plain $a_n$ and $A(n)$ of Definition~\ref{def:reg_index_fgv} are used.

\begin{proposition}\label{prop:E_functional}
With $\gamma=z_0$ define
\[
I(n)=n^{\gamma}H(n)-\gamma\sum_{k=m_n+1}^{n}k^{\gamma-1}H(k).
\]
Then $I(2p+2)-I(2p)=\mathcal O(p^{-2})$, the sequence $I(n)$ converges, and
\begin{equation}\label{eq:E_mode_limit}
\lim_{n\to\infty}I(n)=\Delta'(s_0)\,\kappa(1)=N(s_0)\big|_{\beta=1},
\qquad
\Delta'(s_0)=1-z_0\log2 .
\end{equation}
In particular $\kappa(1)\neq0$ if and only if $I(n)$ has a nonzero limit.
\end{proposition}

\begin{proof}
Substituting $H(n)=H(n-1)-H(m_n)/n$ in $I(n)-I(n-1)$ and tracking the window
$(m_n,n]$, which gains the rank $n$ and loses the rank $m_n$ exactly when $n$ is odd, gives
the exact drift
\[
I(n)-I(n-1)=P(n)\,H(n-1)-Q(n)\,H(m_n),
\]
with $P(n)=n^{\gamma}-\gamma n^{\gamma-1}-(n-1)^{\gamma}=\mathcal O(n^{\gamma-2})$ and
$Q(n)=n^{\gamma-1}-\gamma n^{\gamma-2}-\gamma\,\mathbf 1_{\{2\nmid n\}}\,m_n^{\gamma-1}$. For
odd $n$ the identity $\gamma\,2^{1-\gamma}=2$, a rewriting of $2^{\gamma}=\gamma$, turns
$\gamma m_n^{\gamma-1}$ into $2(n-1)^{\gamma-1}(1+\mathcal O(1/n))$, so $Q$ alternates at
size $n^{\gamma-1}$ between parities. Pairing the odd rank $2p+1$ with the even rank $2p+2$,
which share $m_n=p$, cancels the alternating part to leading order and leaves
\[
I(2p+2)-I(2p)
=\mathcal O\bigl(p^{\eta-2}\bigr)\bigl(|H(2p)|+|H(2p+1)|+|H(p)|\bigr).
\]
Corollary~\ref{cor:E_regimes} at $\beta=1$ gives $|H(n)|\le Cn^{-\eta}$, so the paired drift
is $\mathcal O(p^{-2})$, summable, and $I$ converges along even ranks, the odd gaps being
$\mathcal O(n^{-1})$ singly. For the value of the limit insert
\eqref{eq:E_master} at $\beta=1$, where $c_1=0$, into $I(n)$. The $\kappa n^{-z_0}$ part
contributes $\kappa\bigl(1-\gamma\sum_{m_n<k\le n}k^{-1}\bigr)\to\kappa(1-z_0\log2)$. The
conjugate part contributes
$\bar\kappa\,n^{\gamma-\bar\gamma}\bigl[1-\gamma(1-2^{\bar\gamma-\gamma})/(\gamma-\bar\gamma)\bigr]
+\mathcal O(n^{-1})$, and the bracket vanishes identically because
$2^{\bar\gamma}=\bar\gamma$, so the conjugate mode is annihilated in the limit. The remainder
$\mathcal O(n^{-7/8})$ contributes $\mathcal O(n^{\eta-7/8})\to0$. Hence
$\lim I=\kappa(1)(1-z_0\log2)=\kappa(1)\Delta'(s_0)$, which equals $N(s_0)$ at $\beta=1$ by
\eqref{eq:E_kappa}.
\end{proof}

The functional proves convergence and identifies its limit, but it does not by itself decide the
particular value $\kappa(1)$. No such value is needed for sharpness after
Proposition~\ref{prop:E_kappa}.
The regimes now assemble into the statement of the entry.

\begin{theorem}\label{thm:E_index}
The folded affine profile is a function of good variation. Transparency holds for every real
$\beta<\eta(g)$, while for every $\beta\ge\eta(g)$,
$A_\beta(n)=\mathcal O_\beta(n^{-\eta(g)})$. There are nontransparent exponents arbitrarily
close to $\eta(g)$ from above. Consequently
\[
\tau(g)=\alpha(g)=\eta(g)=\Re z_0.
\]
\end{theorem}

\begin{proof}
Corollary~\ref{cor:E_regimes} gives transparency below $\eta$ and absorption at and above it.
Proposition~\ref{prop:E_kappa} gives, in every interval $(\eta,\eta+\varepsilon)$, an exponent
where transparency fails. Hence the supremum defining $\tau$ is exactly $\eta$, and the absorption
clause gives membership and $\alpha=\tau$.
\end{proof}

\begin{remark}
The meromorphic argument proves that $\kappa$ is not the zero function, and that is all the
sharpness clause consumes. It does not decide the particular value $\kappa(1)$, which is a
question of its own and enters no proof here.
\end{remark}

\begin{numobs}\label{numobs:E_check}
Exploratory computations agree with the formulas above over the ranges tested, and they are
recorded because they show the reader what the entry looks like, not because anything rests on
them. The recurrence \eqref{eq:E_recurrence} reproduces forward substitution in the defining
equation to $2\cdot10^{-14}$ over $n\le3000$ at four exponents including $\beta=\eta(g)$. The
functional equation \eqref{eq:E_funceq} is verified at interior points to $5\cdot10^{-13}$. The
expansion \eqref{eq:E_master} is verified at
$\beta\in\{0,\tfrac25,1,\tfrac32,\tfrac52\}$, the residual times $n^{7/8}$ tending to zero
through values below $10^{-4}$ at $n=4\cdot10^{6}$, and the connection coefficient at those
exponents comes out as $-0.775+0.625i$, $-0.552+1.098i$, $0.033+1.406i$, $0.429+1.345i$ and
$0.745+1.061i$, nonzero in every case tested. Above the threshold $n^{\eta}|A(n)|$ is bounded
without decay, its supremum matching $2|\kappa(\beta)|$, for instance $2.812$ at $\beta=1$, and
the number of sign changes of $A$ up to $2\cdot10^{6}$ is seven against $7.24$ predicted by the
phase $\Im(z_0)\log n/\pi$. The period of the oscillation in $\log n$ is
$2\pi/\Im z_0=4.0086$, a factor $55.07$ in $n$, which is why short numerical windows do not
display it. None of this enters a proof, and the nontriviality that
Theorem~\ref{thm:E_index} consumes comes instead from the meromorphic pole established in
Proposition~\ref{prop:E_kappa}.
\end{numobs}
\galleryentry{F}{The linear logarithmic profile}
 {$g(x)=x-\log x$ on $(0,1]$}
 {function of good variation}
 {$g^{*}(z)=\frac{z^{2}-z+1}{z(z-1)}$, zeros $z_\pm=\tfrac12\pm\tfrac{\sqrt3}{2}i$}
 {$\alpha(g)=\tfrac12=\eta(g)$}
 {proved, Theorem~\ref{thm:F_index}}

For $x\in(0,1]$ put $g(x)=x-\log x$. Then $g(1)=1$ and, the integral converging absolutely on
$\Re z<0$,
\[g^{*}(z)=-z\int_0^1(t-\log t)\,t^{-z-1}\,dt=\frac{z}{z-1}-\frac1z=\frac{z^2 - z + 1}{z(z-1)}.
\]
The zeros of the numerator are $z_\pm=\tfrac12\pm\tfrac{\sqrt3}{2}i$, so $\eta(g)=\tfrac12$.
The kernel is unbounded at the origin, where it grows like $-\log x$, a slowly varying scale,
so \S\ref{sec:fgv_slowly_varying} applies and both the transform and the arithmetic averages
are available on the whole half plane $\Re z<0$. This entry and Appendix~\ref{app:B} are the
two unbounded members of the gallery.

\begin{theorem}\label{thm:F_index}
The kernel $g(x)=x-\log x$ is a function of good variation with
$\alpha(g)=\tfrac12=\eta(g)$. For every $\beta<\tfrac12$,
\[
A(n)=\frac{1}{g^{*}(\beta)}\,n^{-\beta}+\mathcal{O}_\beta\!\left(n^{-1/2}+n^{-\beta-1}\right)
=\frac{1}{g^{*}(\beta)}\,n^{-\beta}+o(n^{-\beta}),
\]
and for every $\beta\ge\tfrac12$ one has $A(n)=\mathcal{O}_\beta(n^{-1/2})$, with no
logarithmic loss at the critical exponent. Transparency fails at a set of exponents
accumulating at $\tfrac12$ from above, so the threshold is attained.
\end{theorem}

The proof occupies the rest of this appendix. Two Abel summations turn the defining equation
into a closed relation between the partial sums alone, differencing twice makes it a recurrence
of order two, a discrete Riccati\index[terms]{Riccati equation}\index[names]{Riccati, J.} equation factorizes that recurrence into two of order one, and
the three regimes then follow by variation of constants\index[terms]{variation of constants}.

\subsection*{The exact equation and its order}

The defining equation is equivalent to a recurrence of the second order.

\begin{proposition}\label{prop:F_exact}
For every $n\ge1$ the defining equation $\sum_{k\le n}a_k\,g(k/n)=n^{-\beta}$ is equivalent to
\begin{equation}\label{eq:F_exact}
A(n)+\sum_{r<n}A(r)\Big[\log\Big(1+\frac1r\Big)-\frac1n\Big]=n^{-\beta}.
\end{equation}
\end{proposition}

\begin{proof}
Put
\[
 B(n)=\sum_{k\le n}ka_k,
 \qquad
 S(n)=\sum_{k\le n}a_k\log k.
\]
The defining equation reads $A(n)\log n+B(n)/n-S(n)=n^{-\beta}$, and Abel summation
gives
\[
 B(n)=nA(n)-\sum_{r<n}A(r),
 \qquad
 S(n)=A(n)\log n-\sum_{r<n}A(r)\log(1+1/r).
\]
Substituting both, the two terms carrying $\log n$ cancel and \eqref{eq:F_exact}
follows.
\end{proof}

The bracket in \eqref{eq:F_exact} is the whole content of the kernel. It is not a multiplier
acting on $A(n-1)$, and the equation does not reduce to first order. Differencing twice removes
the accumulated sum.

\begin{proposition}\label{prop:F_order2}
Put
\[
c_n=n(n-1)\log\frac{n}{n-1}-(n-1),
\qquad
\Delta_n=n^{-\beta}-(n-1)^{-\beta},
\qquad
\varphi(n)=n(n-1)\Delta_n .
\]
Then $c_n=\tfrac12-\tfrac1{6n}+\mathcal{O}(n^{-2})$, and for every $n\ge2$
\begin{equation}\label{eq:F_order2}
n(n+1)A(n+1)-\big[2n^{2}-c_{n+1}\big]A(n)+\big[n(n-1)-c_n+1\big]A(n-1)
=\varphi(n+1)-\varphi(n).
\end{equation}
Written in normalized form $A(n+1)=T_nA(n)-D_nA(n-1)+f_n$, the coefficients satisfy
\begin{equation}\label{eq:F_TD}
T_n=2-\frac2n+\frac{3}{2n^{2}}+\mathcal{O}(n^{-3}),
\qquad
D_n=1-\frac2n+\frac{5}{2n^{2}}+\mathcal{O}(n^{-3}),
\qquad
f_n=(\beta^{2}-\beta)\,n^{-\beta-2}+\mathcal{O}_\beta(n^{-\beta-3}).
\end{equation}
\end{proposition}

\begin{proof}
Let $U(n)=\sum_{r\le n}A(r)$. Taking the difference of \eqref{eq:F_exact} at ranks $n$ and
$n-1$, and using $U(n-1)=U(n-2)+A(n-1)$ to rewrite $-U(n-1)/n+U(n-2)/(n-1)$ as
$U(n-2)/(n(n-1))-A(n-1)/n$, one obtains after multiplication by $n(n-1)$ the exact relation
\[
n(n-1)\,a_n+U(n-2)+c_n\,A(n-1)=\varphi(n)\qquad(n\ge2).
\]
Its difference between ranks $n+1$ and $n$ replaces $U(n-1)-U(n-2)$ by $A(n-1)$ and removes the
accumulated sum. Substituting $a_{n+1}=A(n+1)-A(n)$ and $a_n=A(n)-A(n-1)$, and using
$(n+1)n+n(n-1)=2n^{2}$, gives \eqref{eq:F_order2}. The expansion
$\log\frac{n}{n-1}=\frac1n+\frac1{2n^{2}}+\frac1{3n^{3}}+\cdots$ gives the stated form of
$c_n$, and dividing the two brackets of \eqref{eq:F_order2} by $n(n+1)$ gives the expansions of
$T_n$ and $D_n$. For the forcing, expanding $n(n-1)^{1-\beta}$ in powers of $1/n$ gives
$\varphi(n)=-\beta n^{1-\beta}+\tfrac12\beta(1-\beta)n^{-\beta}+\mathcal{O}_\beta(n^{-\beta-1})$,
whence $\varphi(n+1)-\varphi(n)=-\beta(1-\beta)n^{-\beta}+\mathcal{O}_\beta(n^{-\beta-1})$ and
$f_n=(\beta^{2}-\beta)n^{-\beta-2}+\mathcal{O}_\beta(n^{-\beta-3})$.
\end{proof}

The second-order recurrence has an indicial polynomial, whose roots are the exponents of the
homogeneous solutions.

\begin{proposition}\label{prop:F_indicial}
Let $L$ denote the left side of \eqref{eq:F_order2} divided by $n(n+1)$. For every complex
$\gamma$,
\begin{equation}\label{eq:F_indicial}
L\big[n^{-\gamma}\big]=\big(\gamma^{2}-\gamma+1\big)\,n^{-\gamma-2}
+\mathcal{O}_\gamma\!\left(n^{-\gamma-3}\right).
\end{equation}
The indicial polynomial\index[terms]{indicial polynomial} of the recurrence is therefore the numerator of $g^{*}$, and its roots
are $\gamma_\pm=\tfrac12\pm\tfrac{\sqrt3}{2}i$.
\end{proposition}

\begin{proof}
With $(1\pm 1/n)^{-\gamma}=1\mp\gamma/n+\gamma(\gamma+1)/(2n^{2})+\mathcal{O}(n^{-3})$ and the
expansions \eqref{eq:F_TD}, the three terms of $n(n+1)L[n^{-\gamma}]$ contribute
\[
n^{2}+n-\gamma n-\gamma+\tfrac{\gamma(\gamma+1)}2,
\qquad
-2n^{2}+\tfrac12,
\qquad
n^{2}-n+\tfrac12+\gamma n-\gamma+\tfrac{\gamma(\gamma+1)}2,
\]
each up to $\mathcal{O}_\gamma(1/n)$, the correction $c_n-\tfrac12=\mathcal{O}(1/n)$ entering
only at that order. The coefficients of $n^{2}$ and of $n$ cancel identically, and what remains
is $\gamma(\gamma+1)-2\gamma+1=\gamma^{2}-\gamma+1$. Dividing by $n(n+1)$ gives
\eqref{eq:F_indicial}.
\end{proof}

The cancellation of the term of order $n$ for every $\gamma$ is the degeneracy of this kernel.
The frozen characteristic equation of \eqref{eq:F_order2} has the double root $1$, so the
theorem of Poincar\'e\index[names]{Poincar\'e, H.} and Perron\index[names]{Perron, O.} gives nothing, and the two exponents separate only at order
$n^{-2}$. The factorization below is what replaces it.

\subsection*{Factorization through a discrete Riccati equation}

Each root is realized by a solution of an associated Riccati recursion, from a finite rank on.

\begin{lemma}\label{lem:F_riccati}
Let $\gamma$ be either root of $\gamma^{2}-\gamma+1=0$. There are $N$ and a complex sequence
$(r_n)_{n\ge N}$ with
\begin{equation}\label{eq:F_riccati}
r_{n+1}=T_n-\frac{D_n}{r_n},
\qquad
r_n=1-\frac{\gamma}{n}+\mathcal{O}(n^{-2}).
\end{equation}
\end{lemma}

\begin{proof}
Write $R_n=1-\gamma/n$ and look for $r_n=R_n+\eps_n$. The equation \eqref{eq:F_riccati} read
backwards is $r_n=D_n/(T_n-r_{n+1})$, and the denominator is
\[
T_n-R_{n+1}=1+\frac{\gamma-2}{n}+\frac{\tfrac32-\gamma}{n^{2}}+\mathcal{O}(n^{-3}),
\]
bounded away from $0$ for $n$ large. Expanding the quotient and using
$\gamma^{2}-\gamma=-1$ together with $\tfrac52-\tfrac32=1$,
\[
\frac{D_n}{T_n-R_{n+1}}
=1-\frac{\gamma}{n}+\frac{\tfrac52-\tfrac32+\gamma^{2}-\gamma}{n^{2}}+\mathcal{O}(n^{-3})
=R_n+\mathcal{O}(n^{-3}),
\]
the coefficient of $n^{-2}$ vanishing precisely because $\gamma$ is an indicial root\index[terms]{indicial root}. Hence
$\eps$ must satisfy
\begin{equation}\label{eq:F_eps}
\eps_n=\kappa_n\,\eps_{n+1}+\sigma_n+\theta_n(\eps_{n+1}),
\qquad
\kappa_n=\frac{D_n}{(T_n-R_{n+1})^{2}}=1+\frac{2-2\gamma}{n}+\mathcal{O}(n^{-2}),
\end{equation}
with $\sigma_n=\mathcal{O}(n^{-3})$ and $|\theta_n(u)|\le C|u|^{2}$ for $|u|\le\tfrac12$, the
quadratic remainder of the same expansion. Since $\Re(2-2\gamma)=1$, the backward products
satisfy $\big|\prod_{j=n}^{m-1}\kappa_j\big|\le C\,m/n$ for $N\le n\le m$.

Consider $X_K=\{\eps=(\eps_n)_{n\ge N}:\ \|\eps\|:=\sup_n n^{2}|\eps_n|\le K\}$ with the
norm $\|\cdot\|$, a complete metric space, and on it the map
\[
(\Psi\eps)_n=\sum_{m\ge n}\Big(\prod_{j=n}^{m-1}\kappa_j\Big)
\big[\sigma_m+\theta_m(\eps_{m+1})\big].
\]
A fixed point of $\Psi$ solves \eqref{eq:F_eps}. The series converges absolutely, since
$(m/n)\,m^{-3}$ and $(m/n)\,m^{-4}$ are summable in $m$, and
\[
n^{2}\big|(\Psi\eps)_n\big|
\le n^{2}\sum_{m\ge n}C\frac mn\Big(\frac{C_\sigma}{m^{3}}+\frac{CK^{2}}{m^{4}}\Big)
\le C_1+\frac{C_2K^{2}}{N},
\]
so $\Psi$ maps $X_K$ into itself for $K=2C_1$ and $N$ large. For two elements of $X_K$ the
quadratic terms differ by at most $2CK\,m^{-2}\|\eps-\eps'\|/m^{2}$, whence
$\|\Psi\eps-\Psi\eps'\|\le C_3KN^{-1}\|\eps-\eps'\|<\tfrac12\|\eps-\eps'\|$ for $N$ large.
The contraction principle gives the fixed point, and $r_n=R_n+\eps_n$ is the required solution.
\end{proof}

The two solutions so obtained form a basis of the homogeneous space.

\begin{proposition}\label{prop:F_basis}
The homogeneous recurrence associated with \eqref{eq:F_order2} has, for $n\ge N$, two solutions
\begin{equation}\label{eq:F_basis}
h_\pm(n)=n^{-\gamma_\pm}\big(1+o(1)\big),
\end{equation}
and they are independent. In particular every homogeneous solution is
$\mathcal{O}(n^{-1/2})$, and the two basis solutions have modulus $\asymp n^{-1/2}$ with phase
$\mp\tfrac{\sqrt3}{2}\log n$.
\end{proposition}

\begin{proof}
Take $\gamma=\gamma_+$ in Lemma~\ref{lem:F_riccati} and set $h_+(n)=\prod_{j=N+1}^{n}r_j$.
Then $h_+(n)=r_nh_+(n-1)$, and with $s_n=D_n/r_n$ the identity $r_{n+1}+s_n=T_n$ and
$s_nr_n=D_n$ give
\[
h_+(n+1)-T_nh_+(n)+D_nh_+(n-1)
=h_+(n-1)\big[r_{n+1}r_n-T_nr_n+D_n\big]=0,
\]
so $h_+$ is a homogeneous solution. Since $r_j=1-\gamma_+/j+\mathcal{O}(j^{-2})$, the series
$\sum_j\log(r_j)+\gamma_+/j$ converges absolutely, so
$h_+(n)=C_+n^{-\gamma_+}(1+o(1))$ with $C_+\neq0$, and after normalization this is
\eqref{eq:F_basis}. Taking $\gamma=\gamma_-$ gives $h_-$ in the same way. The two are
independent because $\gamma_+\neq\gamma_-$ makes the ratio $h_+/h_-$ equal to
$n^{-i\sqrt3}(1+o(1))$, which has no limit. Finally
$|n^{-\gamma_\pm}|=n^{-1/2}$ and $\arg n^{-\gamma_\pm}=\mp\tfrac{\sqrt3}{2}\log n$.
\end{proof}

Their Casoratian is computable in closed form, which is what makes variation of constants
explicit here.

\begin{lemma}\label{lem:F_casoratian}
The Casoratian\index[terms]{Casoratian} $W(n)=h_+(n)h_-(n+1)-h_-(n)h_+(n+1)$ satisfies $W(n+1)=D_nW(n)$ and
$|W(n)|\asymp n^{-2}$.
\end{lemma}

\begin{proof}
The relation $W(n+1)=D_nW(n)$ is immediate from the recurrence. Since
$D_n=1-2/n+\mathcal{O}(n^{-2})$, the series $\sum_j\log D_j+2/j$ converges and
$W(n)=W(N)\,C\,n^{-2}(1+o(1))$. It is nonzero because the two solutions are independent.
\end{proof}

\subsection*{The three regimes}

\begin{proof}[Proof of Theorem~\ref{thm:F_index}]
Fix $\beta$ and put $c_\beta=1/g^{*}(\beta)=(\beta^{2}-\beta)/(\beta^{2}-\beta+1)$, which is
finite for every real $\beta$ since $\beta^{2}-\beta+1>0$. By
Proposition~\ref{prop:F_indicial} with $\gamma=\beta$,
\[
L\big[c_\beta n^{-\beta}\big]=c_\beta(\beta^{2}-\beta+1)n^{-\beta-2}
+\mathcal{O}_\beta(n^{-\beta-3})
=(\beta^{2}-\beta)n^{-\beta-2}+\mathcal{O}_\beta(n^{-\beta-3}),
\]
which is $f_n+\mathcal{O}_\beta(n^{-\beta-3})$ by \eqref{eq:F_TD}. The error
$E(n)=A(n)-c_\beta n^{-\beta}$ therefore satisfies the same recurrence with forcing
$\tilde f_n=f_n-L[c_\beta n^{-\beta}]=\mathcal{O}_\beta(n^{-\beta-3})$.

Discrete variation of constants against the basis of Proposition~\ref{prop:F_basis} gives, for
$n>N$,
\begin{equation}\label{eq:F_green}
E(n)=\lambda_+h_+(n)+\lambda_-h_-(n)
+\sum_{m=N}^{n-1}\frac{h_+(m)h_-(n)-h_-(m)h_+(n)}{W(m)}\,\tilde f_m,
\end{equation}
with $\lambda_\pm$ fixed by the two initial values. The Green kernel\index[terms]{Green kernel} obeys
\[
\Big|\frac{h_+(m)h_-(n)-h_-(m)h_+(n)}{W(m)}\Big|
\le C\,\frac{m^{-1/2}n^{-1/2}}{m^{-2}}=C\,m^{3/2}n^{-1/2},
\]
by Proposition~\ref{prop:F_basis} and Lemma~\ref{lem:F_casoratian}, so the sum in
\eqref{eq:F_green} is bounded by $C\,n^{-1/2}\sum_{m<n}m^{-\beta-3/2}$. That sum converges for
$\beta>-\tfrac12$ and is $\mathcal{O}(n^{-\beta-1/2})$ otherwise, so the whole sum is
$\mathcal{O}_\beta(n^{-1/2}+n^{-\beta-1})$. The homogeneous part of \eqref{eq:F_green} is
$\mathcal{O}(n^{-1/2})$ by Proposition~\ref{prop:F_basis}. Hence
\begin{equation}\label{eq:F_master}
A(n)=\frac{n^{-\beta}}{g^{*}(\beta)}+\mathcal{O}_\beta\!\left(n^{-1/2}+n^{-\beta-1}\right).
\end{equation}
For $\beta<\tfrac12$ both error terms are $o(n^{-\beta})$, which is transparency with the stated
constant. For $\beta\ge\tfrac12$ the first term of \eqref{eq:F_master} is itself
$\mathcal{O}(n^{-1/2})$, so $A(n)=\mathcal{O}_\beta(n^{-1/2})$, with no logarithm at the
critical exponent, the degeneracy of the indicial roots\index[terms]{indicial root} being a pair of simple complex roots
rather than a double real one.

For sharpness, write $\Lambda(\beta)=(\lambda_+(\beta),\lambda_-(\beta))$ for the pair of
connection coefficients in \eqref{eq:F_green}. If $\Lambda(\beta)\neq(0,0)$ then the
homogeneous part of $A$ is $C\,n^{-1/2}$ in modulus along a sequence of ranks, by the
independence of $h_\pm$ and the fact that $h_+/h_-$ has no limit, so
$n^{\beta}A(n)$ does not converge to $1/g^{*}(\beta)$ and transparency fails. The two
coefficients are obtained from the initial data by inverting a fixed nonsingular matrix and
subtracting the tail of the sum in \eqref{eq:F_green}, so each is holomorphic in $\beta$ on
$\{\Re\beta>\tfrac12\}$, the series converging locally uniformly there. As $\beta\to+\infty$
along the real axis the forcing tends to $0$ geometrically and the initial data tend to
$A(1)=1$ and $A(2)=\tfrac12-\log2$, a nonzero vector, so $\Lambda(\beta)$ tends to the nonzero
pair attached to that vector. Hence $\Lambda$ is not identically zero, its zero set is discrete,
and there are real exponents $\beta_j\downarrow\tfrac12$ with $\Lambda(\beta_j)\neq(0,0)$. At
each of them transparency fails, so the supremum of Definition~\ref{def:reg_index_fgv} equals
$\tfrac12$, and $\alpha(g)=\tfrac12=\eta(g)$.
\end{proof}

\begin{remark}[The continuous reading]\label{rem:F_continuous}
The same equation in the continuous variable is
$A(x)+\int_1^{x}A(t)\frac{dt}{t}-\frac1x\int_1^{x}A(t)\,dt=x^{-\beta}$, and differentiating it
twice turns it into the Euler equation\index[terms]{Euler equation}
$x^{2}A''+2xA'+A=(\beta^{2}-\beta)x^{-\beta}$, whose particular solution carries the coefficient
$1/g^{*}(\beta)$ of \eqref{eq:F_master}. This kernel is the second order member of the
Euler--Cauchy analogy, and Section~\ref{sec:euler_cauchy} places it beside the first order case
and the integral form, with the complex pair of exponents accounting for the oscillation of
frequency $\tfrac{\sqrt3}{2}$ in $\log n$ that replaces the resonance of the affine kernel.
\end{remark}

\begin{numobs}\label{numobs:F_check}
Equation \eqref{eq:F_exact} and recurrence \eqref{eq:F_order2} reproduce forward substitution
in the defining equation to $3\cdot10^{-13}$ over $n\le400$ at the exponents
$\beta\in\{-2/5,1/5,1/2,9/10\}$, the recurrence holding from $n=2$ on. The Riccati sequence
built by backward iteration satisfies $n^{2}(r_n-1+\gamma_+/n)\to-0.25$, and the factorization
identity $r_{n+1}+D_n/r_n=T_n$ holds to machine accuracy. The two products obey
$|{\textstyle\prod}r_j|\,n^{1/2}\to1.297918$ and $|{\textstyle\prod}s_j|\,n^{3/2}\to2.623294$,
and $\arg\prod r_j+\tfrac{\sqrt3}{2}\log n$ is constant to five decimals over three decades.
The Casoratian\index[terms]{Casoratian} satisfies $|W(n)|\,n^{2}\to3.404820$. At $\beta=-2/5$ the ratio
$n^{\beta}A(n)$ equals $0.358975$ at $n=10^{6}$ against $1/g^{*}(\beta)=0.358974$. Above the
threshold $n^{1/2}A(n)$ stays between $-1.34$ and $+1.75$ over six decades without tending to
zero, and its sign changes up to $2\cdot10^{6}$ number three against the four predicted by the
phase.
\end{numobs}

\galleryentry{G}{The two branch discontinuous profile}
 {$g(x)=1-x/2$ on $(0,1/2)$ and $g(x)=x$ on $[1/2,1]$}
 {conjectured function of good variation}
 {$g^{*}(z)=\frac{P(z)}{z-1}$ with $P(z)=z-\bigl(1-\frac z4\bigr)2^{z}$, every zero on the line $\Re z=2$}
 {$\eta(g)=2$ and $\tau(g)\le2$ proved, membership and $\alpha(g)=2$ open}
 {open, transparency below $2$ and absorption at and above $2$ are unproved, Conjecture~\ref{conj:G_index_revised}, reduced to Open Problem~\ref{op:G_green}}

\rafgalleryfig{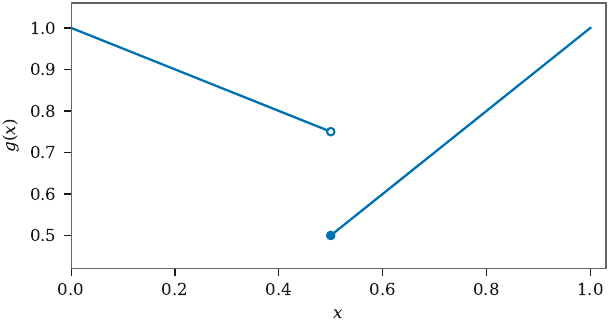}{The two branches, with the value at $x=\tfrac12$ belonging to the lower one. The open circle marks the limit from the left and the filled circle the value taken. The jump is what separates this entry from the affine archetype, and the transform it produces has a pole at $z=1$ and a first zero at $z=2$.}{fig:app_G}

Define $g:(0,1]\to\R$ by
\[
    g(x)=\begin{cases}1-x/2, & 0<x<1/2,\\ x, & 1/2\le x\le 1.\end{cases}
\]
The breakpoint belongs to the second branch, with $g(1/2)=1/2$, and the arithmetic Mellin
transform is
\[g^{*}(z)=\frac{P(z)}{z-1},
\qquad
P(z)=z-\Bigl(1-\frac z4\Bigr)2^{z},
\]
with $P(2)=0$.

One Abel summation turns the defining equation into an exact identity between the two weighted
sums.

\begin{proposition}\label{prop:G_identity}
Put $A(r)=\sum_{k\le r}a_k$, $B(r)=\sum_{k\le r}k\,a_k$ and $m=\lfloor n/2\rfloor$. Then,
exactly at every rank,
\begin{equation}\label{eq:G_identity}
A(m)+\frac{B(n)}{n}-\frac{3B(m)}{2n}-\frac14\,\mathbf{1}_{\{2\mid n\}}\,a_m=n^{-\beta}.
\end{equation}
\end{proposition}

\begin{proof}
When $n$ is odd every index $k\le m$ falls in the first branch and every index $k>m$ in the
second, and the two partial summations assemble into the first three terms. When $n=2m$ is
even the index $k=m$ has $k/n=1/2$ and belongs to the second branch, where $g(1/2)=1/2$, and
carrying it there produces the correction $-a_m/4$. The elementary control at $n=2$ confirms
the parity term. The direct sum is $a_1g(1/2)+a_2g(1)=a_1/2+a_2$, the first three terms of
\eqref{eq:G_identity} alone give $3a_1/4+a_2$, and the corrected identity restores
$a_1/2+a_2$.
\end{proof}

\subsection*{All the zeros lie on one line}

The zeros of the indicial polynomial sit on a single vertical line, and they can be listed.

\begin{theorem}\label{thm:G_zeros}
Every zero of $P$ lies on the vertical line $\Re z=2$. The zeros are exactly
\[
z=2,\qquad z=2\pm it,
\]
where $t>0$ runs through the solutions of the phase equation
\begin{equation}\label{eq:G_phase}
t\log2=2\arctan\frac t2+2k\pi,\qquad k\in\Z_{\ge0}.
\end{equation}
The branch $k=0$ has the single positive root $t_0'=2.685182329\ldots$, each branch $k\ge1$
has exactly one root $t_k$, with $t_1=13.161965\ldots$, $t_2=22.404914\ldots$, and
$t_k=(2k+1)\pi/\log2+\mathcal O(1/k)$. Consequently
\[
\eta(g)=2 .
\]
\end{theorem}

\begin{proof}
Since $z=4$ is not a zero, $P(z)=0$ is equivalent to
\begin{equation}\label{eq:G_equiv}
2^{\,z-2}=\frac{z}{4-z}.
\end{equation}
Write $z=\sigma+it$ and let
\[
u(\sigma,t)=(\sigma-2)\log2-\log\Bigl|\frac{z}{4-z}\Bigr|
=(\sigma-2)\log2-\tfrac12\log\frac{\sigma^{2}+t^{2}}{(4-\sigma)^{2}+t^{2}},
\]
the difference of the logarithmic moduli of the two sides of \eqref{eq:G_equiv}, so that a
zero must satisfy $u=0$ together with the equality of arguments. The identity
$u(2,t)=0$ holds for every $t$, because $|2+it|=|2-it|$, the point $2$ being equidistant
from $0$ and $4$.

Suppose first $|t|\ge\tfrac32$. Then
\[
\frac{\partial u}{\partial\sigma}
=\log2-\frac{\sigma}{\sigma^{2}+t^{2}}-\frac{4-\sigma}{(4-\sigma)^{2}+t^{2}}
\ \ge\ \log2-\frac{1}{2|t|}-\frac{1}{2|t|}
\ \ge\ \log2-\frac23\ >\ 0,
\]
the elementary bound $x/(x^{2}+t^{2})\le1/(2|t|)$ being applied twice. So
$\sigma\mapsto u(\sigma,t)$ is strictly increasing and vanishes only at $\sigma=2$. Every
zero with $|t|\ge\tfrac32$ therefore has $\sigma=2$.

Suppose next $0<|t|\le\tfrac32$, and by conjugation take $t>0$. The argument of the right
side of \eqref{eq:G_equiv} is
\[
\psi(z)=\arg z-\arg(4-z),
\]
and the argument of the left side is $t\log2\in(0,\,\tfrac32\log2]\subset(0,1.04]$. For
$0<\sigma<4$ one has $\psi=\arctan(t/\sigma)+\arctan(t/(4-\sigma))$, with
\[
\frac{\partial\psi}{\partial\sigma}=\frac{-t}{\sigma^{2}+t^{2}}+\frac{t}{(4-\sigma)^{2}+t^{2}},
\]
negative for $\sigma<2$ and positive for $\sigma>2$, so on $(0,4)$ the function $\psi$
attains its minimum in $\sigma$ at $\sigma=2$, where $\psi=2\arctan(t/2)$. The function
$q(t)=2\arctan(t/2)-t\log2$ satisfies $q(0)=0$, is increasing up to
$t^{*}=2\sqrt{1/\log2-1}=1.3307\ldots$ and decreasing after, with
$q(\tfrac32)=0.2472\ldots>0$, so $q>0$ on $(0,\tfrac32]$. Hence for $0<\sigma<4$,
$\sigma\neq2$,
\[
\psi(z)>2\arctan(t/2)>t\log2 ,
\]
strictly. For $\sigma\le0$ one has $\arg z\ge\pi/2$ and $\arg(4-z)<0$, so
$\psi>\pi/2>1.04$, and for $\sigma\ge4$ one has $\arg(4-z)\le-\pi/2$, so again
$\psi\ge\pi/2$. In every case $0<t\log2<\psi(z)$, and since $\psi\le\tfrac{3\pi}2$ in every
region while $t\log2+2\pi>2\pi>\tfrac{3\pi}2$ and $t\log2-2\pi<0$, no determination shifted
by a multiple of $2\pi$ can restore the equality of arguments. There is therefore no zero
with $0<|t|\le\tfrac32$ and $\sigma\neq2$.

On the real axis put $f(\sigma)=2^{\sigma-2}-\sigma/(4-\sigma)$. For $\sigma\le0$ and for
$\sigma>4$ the second term is nonpositive and $f>0$, and $\sigma=4$ is not a zero. On
$(0,4)$ one has $f'(\sigma)=2^{\sigma-2}\log2-4/(4-\sigma)^{2}<0$, because
$(4-\sigma)^{2}2^{\sigma-2}\log2$ is maximal at $\sigma=4-2/\log2$ where it equals
$3.13\ldots<4$. So $f$ is strictly decreasing on $(0,4)$ and vanishes only at $\sigma=2$,
which is the only real zero.

It remains to describe the zeros on the line. At $z=2+it$ equation \eqref{eq:G_equiv} reads
$2^{it}=(2+it)/(2-it)$, an equality of unimodular numbers which holds exactly when the phase
equation \eqref{eq:G_phase} holds. The function $\varphi(t)=t\log2-2\arctan(t/2)$ has
$\varphi(0)=0$, decreases to its minimum $\varphi(t^{*})=-0.2518\ldots$ and then increases to
infinity with slope tending to $\log2$. So $\varphi=0$ has the roots $t=0$ and one further
root $t_0'$, and $\varphi=2k\pi$ has exactly one root $t_k$ for each $k\ge1$. Solving
numerically gives the stated values, and inverting $\varphi(t)=2k\pi$ with
$2\arctan(t/2)=\pi-4/t+\mathcal O(t^{-3})$ gives the asymptotic form of $t_k$.
\end{proof}

The transform of this kernel thus carries infinitely many zeros, all of them on a single vertical line, in contrast with the isolated conjugate pairs met elsewhere in the gallery. The analytic index is their common real part, and whether the regularity index takes that same value is the open question of this entry.

\subsection*{The complete scalar reduction}

The parity of the rank splits the equation into two branches.

\begin{proposition}\label{prop:G_reduction}
Let $d_\beta(n)=n^{1-\beta}-(n-1)^{1-\beta}$. The defining equation is equivalent to
$A(1)=1$ together with the pair of recurrences, for $p\ge1$,
\begin{align}
A(2p)&=A(2p-1)+\frac{d_\beta(2p)-A(p-1)}{2p},\label{eq:G_even}\\[1mm]
A(2p+1)&=A(2p)+\frac{d_\beta(2p+1)-A(p)-\tfrac p2\,\bigl(A(p)-A(p-1)\bigr)}{2p+1}.\label{eq:G_odd}
\end{align}
The pair determines the sequence from $A(1)$ alone, so the space of homogeneous solutions
has dimension one.
\end{proposition}

\begin{proof}
Multiply \eqref{eq:G_identity} by $n$ and difference between consecutive ranks, splitting by
parity. For $n=2p$ one has $m_n=p$ and $m_{n-1}=p-1$, and with $B(p)-B(p-1)=pa_p$ the
difference reads
\[
2pA(p)-(2p-1)A(p-1)+2p\,a_{2p}-\tfrac32\,p\,a_p-\tfrac{2p}{4}\,a_p=d_\beta(2p).
\]
The left side collapses, $2pA(p)-(2p-1)A(p-1)=A(p-1)+2pa_p$ and
$\tfrac32p+\tfrac p2=2p$, so the four terms in $a_p$ cancel and
$A(p-1)+2p\,a_{2p}=d_\beta(2p)$, which is \eqref{eq:G_even}. For $n=2p+1$ one has
$m_n=m_{n-1}=p$, the terms in $B(p)$ cancel, the parity term enters with the opposite sign
from rank $2p$, and
\[
A(p)+(2p+1)a_{2p+1}+\tfrac p2\,a_p=d_\beta(2p+1),
\]
which is \eqref{eq:G_odd}. Conversely the pair reconstructs \eqref{eq:G_identity} by
summation, the value at $n=1$ anchoring the telescoping.
\end{proof}

Each branch carries its own multiplier, and the two combine into the indicial data.

\begin{proposition}\label{prop:G_indicial}
Write \eqref{eq:G_even} and \eqref{eq:G_odd} as $n\bigl(A(n)-A(n-1)\bigr)+F_n[A]=d_\beta(n)$,
where $F_n$ is the dilated feedback of the relevant parity. For every complex $\gamma$,
\[
n\bigl(\Delta n^{-\gamma}\bigr)+F_n[k\mapsto k^{-\gamma}]
=\chi_{\pm}(\gamma)\,n^{-\gamma}+\mathcal O_\gamma(n^{-\gamma-1}),
\qquad
\begin{cases}
\chi_{+}(\gamma)=2^{\gamma}-\gamma & (n\ \text{even}),\\[1mm]
\chi_{-}(\gamma)=\bigl(1-\tfrac\gamma2\bigr)2^{\gamma}-\gamma & (n\ \text{odd}),
\end{cases}
\]
and the half sum of the two parity multipliers is
\[
\tfrac12\bigl(\chi_{+}(\gamma)+\chi_{-}(\gamma)\bigr)
=\Bigl(1-\frac\gamma4\Bigr)2^{\gamma}-\gamma=-P(\gamma),
\]
the negative of the numerator of $g^{*}$.
\end{proposition}

\begin{proof}
On even ranks $n=2p$ the feedback is $A(p-1)$ and
$(p-1)^{-\gamma}=2^{\gamma}n^{-\gamma}(1+\mathcal O(1/n))$. On odd ranks $n=2p+1$ the
feedback is $A(p)+\tfrac p2(A(p)-A(p-1))$, and for the power profile
$A(p)=p^{-\gamma}$ the increment is $A(p)-A(p-1)=-\gamma p^{-\gamma-1}(1+\mathcal O(1/p))$,
so the feedback is $p^{-\gamma}(1-\tfrac\gamma2)+\mathcal O(p^{-\gamma-1})
=(1-\tfrac\gamma2)2^{\gamma}n^{-\gamma}+\mathcal O(n^{-\gamma-1})$. In both cases
$n\Delta n^{-\gamma}=-\gamma n^{-\gamma}+\mathcal O(n^{-\gamma-1})$, and the two multipliers
follow. The half sum is a direct computation.
\end{proof}

The dilation structure is therefore that of Appendix~\ref{app:E} with one new feature. The
feedback alternates between two multipliers along the parity of the rank, and only their mean
reproduces the numerator of the transform. A power profile alone can not be indicial for
both parities at once, the natural unknown being the pair of subsequences of even and odd
rank, and the two exact recurrences above are the corresponding system.

A maximum norm applied to the two recurrences loses the cancellation between the two parity
multipliers. The following generating equation records the same obstruction without
suppressing the odd rank increment.

\begin{proposition}\label{prop:G_mahler}
Put $A(0)=0$ and
\[
F_\beta(z)=\sum_{n\ge1}a_nz^n.
\]
For $|z|<1$ the defining system is equivalent to
\begin{equation}\label{eq:G_mahler}
 zF_\beta'(z)
 +\frac{z}{1-z}F_\beta(z^2)
 +\frac{z^3}{2}F_\beta'(z^2)
 =(1-z)\operatorname{Li}_{\beta-1}(z).
\end{equation}
\end{proposition}

\begin{proof}
The coefficient of $z^n$ in the first term is $na_n$. The coefficient in
$zF_\beta(z^2)/(1-z)$ is $A(p-1)$ when $n=2p$ and $A(p)$ when
$n=2p+1$. The last term contributes $pa_p/2$ at rank $2p+1$ and zero at
even ranks. The coefficient on the right is
$d_\beta(n)=n^{1-\beta}-(n-1)^{1-\beta}$ for $n\ge2$. The coefficient at
rank one gives $a_1=1$. These are exactly \eqref{eq:G_even} and
\eqref{eq:G_odd}.
\end{proof}

At the endpoint no transparent coefficient exists, which is what leaves this entry open.

\begin{proposition}\label{prop:G_endpoint}
There is no real number $\Xi$ such that the solution at $\beta=2$ satisfies
\[
A(n)=\Xi n^{-2}+o(n^{-2}).
\]
Consequently $\tau(g)\le2$. If transparency holds for every $\beta<2$ and the
partial sums are absorbed above that value, then $g$ is a function of good
variation\index[terms]{function of good variation} with $\alpha(g)=2$, and the index is sharp.
\end{proposition}

\begin{proof}
Suppose that $A(r)=\Xi r^{-2}+o(r^{-2})$. The series
\[
S=\sum_{j\ge1}A(j)
\]
is absolutely convergent. Summation by parts gives
\[
B(r)=rA(r)-\sum_{j<r}A(j)
=-S+\frac{2\Xi}{r}+o(r^{-1}).
\]
Put $m=\lfloor n/2\rfloor$. The assumed asymptotic gives
\[
n^2A(m)=4\Xi+o(1),
\qquad
nB(n)=-Sn+2\Xi+o(1),
\]
and
\[
-\frac{3n}{2}B(m)=\frac{3Sn}{2}-6\Xi+o(1).
\]
It also gives
\[
n^2\mathbf 1_{\{2\mid n\}}a_m=o(1).
\]
Multiplication of \eqref{eq:G_identity} by $n^2$ yields
\[
1=\frac{Sn}{2}+(4+2-6)\Xi+o(1)
=\frac{Sn}{2}+o(1).
\]
Division by $n$ gives $S=0$. Substitution in the preceding equality then
gives $1=o(1)$, which is impossible.

Thus $\beta=2$ is not transparent. No number $c>2$ can have every
$\beta<c$ transparent, so $\tau(g)\le2$. If every $\beta<2$ is
transparent and the partial sums are absorbed above that value, then
$\alpha(g)=2$. Every interval $[2,2+\delta)$ contains
the nontransparent point $2$, which gives sharpness.
\end{proof}

The comparison kernel from Appendix~\ref{app:P} at $m=2$ is
\[
q(x)=1-\frac{x}{2}\quad(0<x<1),
\qquad q(1)=1.
\]
Its transform is $(z-2)/(2(z-1))$. It has one real zero at $2$. The kernel
of this appendix has the real zero and the pairs $2\pm it_k$. The numerical
data below separate these two situations. For $q$, the quantity $n^2A(n)$
settles when $\beta>2$. For $g$, it retains the frequencies $t_k$ on the
logarithmic scale.

\begin{conjecture}\label{conj:G_index_revised}
The kernel $g$ has $\alpha(g)=\eta(g)=2$. For every $\beta<2$,
\[
A(n)=\frac{1}{g^*(\beta)}n^{-\beta}+o(n^{-\beta}).
\]
At the endpoint,
\[
A(n)=-\frac{1}{2-2\log2}\,n^{-2}\log n+\mathcal O(n^{-2}).
\]
For every $\beta>2$, one has $A(n)=\mathcal O_\beta(n^{-2})$.
\end{conjecture}

\begin{openproblem}[Green kernel estimates]\label{op:G_green}
Let $\mathcal L$ act on sequences $U=(U_n)_{n\ge0}$ with $U_0=0$ by
\[
(\mathcal LU)_1=U_1,
\]
\[
(\mathcal LU)_{2p}
=2p(U_{2p}-U_{2p-1})+U_{p-1},
\]
and
\[
(\mathcal LU)_{2p+1}
=(2p+1)(U_{2p+1}-U_{2p})
+U_p+\frac p2(U_p-U_{p-1}).
\]
For $1\le j\le n$, define $R(n,j)$ by
\[
R(n,j)=0\quad(n<j),
\qquad
\mathcal L\bigl(R(\,\cdot\,,j)\bigr)(n)=\mathbf1_{\{n=j\}}.
\]
Put
\[
P(z)=z-\left(1-\frac z4\right)2^z.
\]
Prove the following two assertions.

For every $\eps>0$ there is $C_\eps>0$ such that
\begin{equation*}\label{eq:G_green_bound}
|R(n,j)|\le C_\eps jn^{-2+\eps},
\qquad 1\le j\le n.
\tag{G1}
\end{equation*}
For every real $\sigma<2$ one has
\begin{equation*}\label{eq:G_green_limit}
\lim_{n\to\infty}
\sum_{j=2}^nR(n,j)\left(\frac jn\right)^{-\sigma}
=-\frac1{P(\sigma)}.
\tag{G2}
\end{equation*}
\end{openproblem}

\begin{conditionaltheorem}\label{cthm:G_green_closure}
If \eqref{eq:G_green_bound} and \eqref{eq:G_green_limit} hold, then $g$ is
a RAF with $\alpha(g)=2$. Transparency holds below $2$, absorption holds
at and above $2$, and the index is sharp.
\end{conditionaltheorem}

\begin{proof}
The inverse relation is
\[
A(n)=\sum_{j\le n}R(n,j)d_\beta(j),
\]
where $d_\beta(1)=1$ and
\[
d_\beta(j)=(1-\beta)j^{-\beta}
+\mathcal O_\beta(j^{-\beta-1})
\qquad(j\ge2).
\]
For $\beta<2$, \eqref{eq:G_green_limit} gives
\[
A(n)
=\frac{\beta-1}{P(\beta)}n^{-\beta}+o(n^{-\beta})
=\frac{1}{g^*(\beta)}n^{-\beta}+o(n^{-\beta}).
\]
The error is $o(n^{-\beta})$ by \eqref{eq:G_green_bound}. For $\beta=2$,
the same bound produces at most the factor
$\sum_{j\le n}j^{-1}=\mathcal O(\log n)$. For $\beta>2$, the corresponding sum
converges. Hence $A(n)=\mathcal O_\eps(n^{-2+\eps})$ for every
$\beta\ge2$. Proposition~\ref{prop:G_endpoint} gives the upper bound on
the index and its sharpness.
\end{proof}

\begin{proofstatus}{The exact identity, zero localization, parity recurrences and
multipliers, and the Mahler differential equation are proved.
Proposition~\ref{prop:G_endpoint} proves that $\beta=2$ is not transparent
and gives $\tau(g)\le2$. In Conditional
Theorem~\ref{cthm:G_green_closure}, \eqref{eq:G_green_limit} supplies the
main term below $2$, while \eqref{eq:G_green_bound} controls its error and
yields RAF absorption $\mathcal O_\eps(n^{-2+\eps})$ at and above $2$.
These estimates do not prove the sharper endpoint asymptotic or the bound
$\mathcal O_\beta(n^{-2})$ above $2$ stated in
Conjecture~\ref{conj:G_index_revised}. A maximum norm cannot establish
\eqref{eq:G_green_bound}: it removes the cancellation between the two
parity multipliers. Numerical Observation~\ref{numobs:G_check_revised} is
evidence only. Nothing outside this appendix assumes the conjecture or
either Green estimate.}
\end{proofstatus}

\begin{numobs}\label{numobs:G_check_revised}
Exact rational forward substitution at $\beta=3$ gives zero residual for
\eqref{eq:G_identity}, \eqref{eq:G_even}, \eqref{eq:G_odd}, and
\eqref{eq:G_mahler} through rank $120$. At $80$ decimal digits and through
rank $500$, the maximum residual is below $1.5\cdot10^{-79}$ for the
identity and the two recurrences at
$\beta\in\{\tfrac12,\tfrac32,3\}$.

The first positive phase and its period in $\log n$ are
\[
 t_0=2.68518232902163\ldots,
 \qquad
 2\pi/t_0=2.33994736196\ldots
\]
The number $2\pi/\log2$ is the asymptotic spacing of the $t_k$, not a period in
$\log n$.

At $\beta=3$, a fit of $n^2A(n)$ on the first ten frequencies $t_k$ over
$20000\le n\le4\cdot10^6$ detects the infinite family through its first
terms. The coefficient of determination is $0.99379$, and the dominant
fitted amplitude is $21.4247$ at $t_0$. Over
$250000\le n\le4\cdot10^6$, the measured range is
\[
-31.86323\le n^2A(n)\le25.34876.
\]
For the comparison kernel, $n^2A(n)=-3.289871$ at
$n=4\cdot10^6$, with no persistent oscillation. At $\beta=2$, a fit
including the same frequencies gives the coefficient $-1.62410$ for
$\log n$, against $-1/(2-2\log2)=-1.62945\ldots$. These data support
Conjecture~\ref{conj:G_index_revised} and do not prove the Green kernel
estimates.
\end{numobs}
\galleryentry{H}{The shifted rational kernel}
 {$G(n,k)=\frac{n+k+x}{n+k+y}\cdot\frac{2n+y}{2n+x}$, with $x,y>0$ and $x\neq y$}
 {regular arithmetic function}
 {$G^{*}(z)=1$, constant, no zero, so $\eta$ is not defined}
 {$\alpha(G)=2$, $\eta$ not defined}
 {proved, Theorem~\ref{thm:H_index}}

Fix $x,y>0$ with $x\neq y$ and put $c=y-x$. The kernel
$G(n,k)=\frac{n+k+x}{n+k+y}\cdot\frac{2n+y}{2n+x}$ satisfies $G(n,n)=1$ and has constant
arithmetic Mellin transform $G^{*}(z)=1$, so no analytic index is available and the index must
come from the equation.

\begin{proposition}\label{prop:H_exact}
With $w_n=\dfrac{2n+y}{2n+x}$ the defining equation is equivalent to
\begin{equation}\label{eq:H_recurrence-A}
A(n)=n^{-\beta}+c\,w_n\sum_{k<n}\frac{A(k)}{(n+k+y)(n+k+1+y)},
\end{equation}
and for every $k\ge1$
\begin{equation}\label{eq:H_telescope}
\sum_{n>k}\frac{1}{(n+k+y)(n+k+1+y)}=\frac{1}{2k+1+y}.
\end{equation}
\end{proposition}

\begin{proof}
Abel summation of $\sum_{k\le n}a_kG(n,k)$ against $v_k=G(n,k)$ uses
$v_{k+1}-v_k=w_n(y-x)/((n+k+y)(n+k+1+y))$, which follows from expanding the numerator of the
difference, and $v_n=1$. For \eqref{eq:H_telescope}, the summand is
$\frac{1}{n+k+y}-\frac{1}{n+k+1+y}$ and the sum telescopes from $n=k+1$.
\end{proof}

Throughout the appendix put
\[q(n,k)=\frac{1}{(n+k+y)(n+k+1+y)},
\qquad
W=\max\Bigl(1,\frac{4+y}{4+x}\Bigr)=\sup_{n\ge2}w_n,
\]
so that $q(n,k)\le q(n,1)\le n^{-2}$ for $1\le k<n$ and \eqref{eq:H_telescope} reads
$\sum_{n>k}q(n,k)=1/(2k+1+y)$.

Below the index the feedback term is negligible and the response follows the forcing.

\begin{proposition}\label{prop:H_regimes}
For every $\beta<2$ the feedback term of \eqref{eq:H_recurrence-A} is $o(n^{-\beta})$ and
$A(n)\sim n^{-\beta}=n^{-\beta}/G^{*}(\beta)$. For every $\beta>2$ the series
\begin{equation}\label{eq:H_ell}
\ell(\beta)=\sum_{k\ge1}A_\beta(k)
\end{equation}
converges absolutely and
\begin{equation}\label{eq:H_asymptotic}
A(n)=c\,\ell(\beta)\,n^{-2}+\mathcal{O}_\beta(n^{-3}\log n+n^{-\beta}),
\end{equation}
for both signs of $c$. In particular $A(n)=\mathcal{O}(n^{-2})$ at and above the
exponent two.
\end{proposition}

\begin{proof}
Let $T(n)=\sum_{k\le n}|A(k)|$. Bounding the feedback sum of \eqref{eq:H_recurrence-A} by its
largest weight gives $|A(n)|\le n^{-\beta}+|c|\,W\,T(n-1)\,q(n,1)$ with $q(n,1)\le n^{-2}$,
hence $T(n)\le T(n-1)\bigl(1+|c|W n^{-2}\bigr)+n^{-\beta}$, and the product
$\prod_{m\ge2}(1+|c|Wm^{-2})$ converges. Therefore $T(n)\le C_{x,y}\sum_{k\le n}k^{-\beta}$,
so $T$ is bounded for $\beta>1$ and grows at most like $n^{1-\beta}$ for $\beta<1$. Feeding
the bound back into \eqref{eq:H_recurrence-A}, the feedback term is $\mathcal
O(T(n)\,n^{-2})$, which gives $A(n)=\mathcal O(n^{-2})$ for $\beta>1$ and makes the feedback
$o(n^{-\beta})$ for every $\beta<2$, whence the transparency statement. Absolute convergence
of \eqref{eq:H_ell} for $\beta>2$ follows from $A(n)=\mathcal O(n^{-2})$. For the asymptotic,
compare the weight to its diagonal value. For $1\le k<n$ one has
$|q(n,k)-n^{-2}|\le C_y\,(k+1)\,n^{-3}$, since the two factors of $n^{2}q(n,k)^{-1}$ are
$1+\mathcal O((k+1)/n)$ with bounded implied quantities on that range, so
\[
\sum_{k<n}A(k)\,q(n,k)=\frac{1}{n^{2}}\sum_{k<n}A(k)
+\mathcal O\Bigl(n^{-3}\sum_{k<n}(k+1)\,|A(k)|\Bigr)
=\frac{\ell(\beta)}{n^{2}}+\mathcal O\bigl(n^{-3}\log n\bigr),
\]
the discarded tail $\sum_{k\ge n}|A(k)|$ being $\mathcal O(n^{-1})$ and the weighted sum
$\sum_{k<n}(k+1)|A(k)|$ being $\mathcal O(\log n)$, both by $|A(k)|=\mathcal O(k^{-2})$.
With $w_n=1+\mathcal O(1/n)$ this is \eqref{eq:H_asymptotic}.
\end{proof}

By \eqref{eq:H_asymptotic} the sharpness of the exponent two reduces to one question, whether
the function $\ell$, holomorphic on $\{\Re\beta>2\}$ by locally uniform convergence, vanishes
identically. The answer is that it never does, for any admissible pair, and the proof goes
through a family of shifted impulse solutions which expands $\ell$ in a Dirichlet series.

\begin{lemma}\label{lem:H_family}
For $j\ge1$ let $E^{(j)}$ be the solution of
\begin{equation}\label{eq:H_impulse}
E^{(j)}(n)=\mathbf 1_{\{n=j\}}+c\,w_n\sum_{k<n}E^{(j)}(k)\,q(n,k)\qquad(n\ge1),
\end{equation}
so that $E^{(j)}(n)=0$ for $n<j$ and $E^{(j)}(j)=1$. Then the following hold.
\begin{enumerate}[label=(\roman*)]
\item For every $n\ge1$ and every $\beta$,
$A_\beta(n)=\sum_{j\le n}j^{-\beta}E^{(j)}(n)$.
\item There is a constant $C_{x,y}$, independent of $j$, with
$\sum_{n\ge j}|E^{(j)}(n)|\le C_{x,y}$, and moreover
$|E^{(j)}(n)|\le |c|\,W\,C_{x,y}\,n^{-2}$ for $n>j$.
\item The sums $\lambda_j=\sum_{n\ge j}E^{(j)}(n)$ satisfy, for every real $\beta>2$,
\begin{equation}\label{eq:H_dirichlet}
\ell(\beta)=\sum_{j\ge1}\lambda_j\,j^{-\beta},
\end{equation}
the series converging absolutely.
\end{enumerate}
\end{lemma}

\begin{proof}
The recurrence \eqref{eq:H_recurrence-A} determines its solution uniquely from the forcing,
the diagonal term being isolable, and it is linear in the forcing. The right side of (i)
solves the recurrence with forcing $\sum_{j\ge1}j^{-\beta}\mathbf 1_{\{n=j\}}=n^{-\beta}$,
each rank involving finitely many terms since $E^{(j)}(n)=0$ for $j>n$, so (i) follows by
uniqueness. For (ii) put $T_j(n)=\sum_{j\le k\le n}|E^{(j)}(k)|$. For $n>j$ the largest
weight in \eqref{eq:H_impulse} is $q(n,j)\le n^{-2}$, so
$|E^{(j)}(n)|\le|c|\,W\,T_j(n-1)\,n^{-2}$ and
$T_j(n)\le T_j(n-1)(1+|c|Wn^{-2})$ with $T_j(j)=1$. The product
$\prod_{m\ge2}(1+|c|Wm^{-2})=:C^{0}_{x,y}$ converges, so $T_j(n)\le C^{0}_{x,y}$ for all $n$,
which gives both bounds of (ii) with $C_{x,y}=C^{0}_{x,y}\bigl(1+|c|W\sum_{m\ge2}m^{-2}\bigr)$.
For (iii), the double series $\sum_j j^{-\beta}\sum_n|E^{(j)}(n)|\le C_{x,y}\,\zeta(\beta)$
converges for $\beta>1$, so summing (i) over $n$ and exchanging the two sums is legitimate
for $\beta>2$, where $\sum_n A_\beta(n)$ converges absolutely by
Proposition~\ref{prop:H_regimes}.
\end{proof}

The sharpness argument needs positivity of the coefficients from a computable rank on.

\begin{lemma}\label{lem:H_positivity}
Put $\theta_j=\dfrac{|c|\,W}{2j+1+y}$. For every $j$ with $\theta_j<\tfrac12$,
\begin{equation}\label{eq:H_lambda_bound}
\lambda_j\ \ge\ \frac{1-2\theta_j}{1-\theta_j}\ >\ 0 .
\end{equation}
In particular $\lambda_j>0$ for every $j>\tfrac12\bigl(2|c|W-1-y\bigr)$.
\end{lemma}

\begin{proof}
Write $S_j=\sum_{n>j}|E^{(j)}(n)|$, which is finite by Lemma~\ref{lem:H_family}. Taking
absolute values in \eqref{eq:H_impulse}, summing over $n>j$ and exchanging the two sums,
the telescoped weight \eqref{eq:H_telescope} gives
\[
S_j\le|c|\,W\sum_{k\ge j}|E^{(j)}(k)|\sum_{n>k}q(n,k)
=|c|\,W\sum_{k\ge j}\frac{|E^{(j)}(k)|}{2k+1+y}
\le\theta_j\bigl(1+S_j\bigr),
\]
since every rank $k$ carrying a nonzero term satisfies $k\ge j$. When $\theta_j<1$ this
yields $S_j\le\theta_j/(1-\theta_j)$, and then
$\lambda_j\ge E^{(j)}(j)-S_j=1-S_j\ge(1-2\theta_j)/(1-\theta_j)$, which is positive exactly
when $\theta_j<\tfrac12$.
\end{proof}

The index of this kernel is two, whatever the two parameters.

\begin{theorem}\label{thm:H_index}
For every $x,y>0$ with $x\neq y$ the kernel $G$ satisfies $\alpha(G)=2$.
\end{theorem}

\begin{proof}
Transparency below two and absorption at and above two are
Proposition~\ref{prop:H_regimes}, so it remains to show that transparency fails at exponents
accumulating at two from above. Suppose $\ell$ vanished identically on $\{\Re\beta>2\}$. By
\eqref{eq:H_dirichlet} the absolutely convergent Dirichlet series $\sum_j\lambda_jj^{-\beta}$
would vanish for every real $\beta>2$. An induction on $r\ge1$ then forces every
$\lambda_r=0$, since if $\lambda_1=\dots=\lambda_{r-1}=0$ then
$r^{\beta}\sum_{j\ge r}\lambda_jj^{-\beta}=\lambda_r+\sum_{j>r}\lambda_j(j/r)^{-\beta}$ and
the tail tends to $0$ as $\beta\to+\infty$ by dominated convergence, the terms being bounded
by $C_{x,y}(j/r)^{-3}$. This contradicts Lemma~\ref{lem:H_positivity}, which produces ranks
$j$ with $\lambda_j>0$. Hence $\ell$ is not identically zero on the connected half plane
$\{\Re\beta>2\}$. Being holomorphic there, it cannot vanish identically on any real interval
$(2,2+\delta)$, so there are real exponents $\beta_m\downarrow2$ with $\ell(\beta_m)\neq0$.
At such an exponent \eqref{eq:H_asymptotic} gives
$A_{\beta_m}(n)\sim c\,\ell(\beta_m)\,n^{-2}$ with $c\neq0$, which is not $o(n^{-\beta_m})$,
so transparency fails there. The supremum of Definition~\ref{def:reg_index_fgv} is therefore
two.
\end{proof}

\begin{remark}[The constant $\ell_\infty$ and why the shifted family is needed]
\label{rem:H_ellinf}
The limit of $\ell(\beta)$ as $\beta\to+\infty$ along the real axis is
$\ell_\infty(x,y)=\lambda_1$, the sum of the solution carrying the initial impulse alone, by
dominated convergence in \eqref{eq:H_dirichlet}. Reducing sharpness for $x>y$ to the
nonvanishing of this single constant does not work. The contraction condition
$|y-x|W/(3+y)<\tfrac12$ verifies it only on a neighborhood of the diagonal and does not extend
to the whole quadrant. The constant $\lambda_1$ changes
sign along a curve in the $(x,y)$ plane, numerically it vanishes near $(15.67,\,1)$ on the
line $y=1$, and it is negative beyond that curve. Lemma~\ref{lem:H_positivity} explains both
phenomena at once. The contraction ratio attached to the impulse at rank $j$ is
$|c|W/(2j+1+y)$, so no fixed rank can serve every pair, while for each pair all ranks beyond
$|c|W$ are positive, and uniqueness of Dirichlet\index[names]{Dirichlet, P. G. L.} coefficients converts any single positive
rank into the sharpness of the index.
\end{remark}

\begin{numobs}\label{numobs:H_check}
The recurrence \eqref{eq:H_recurrence-A} reproduces forward substitution in the defining
equation to $10^{-15}$ over $n\le400$ at the pairs $(3,1)$ and $(5,\tfrac12)$. The quantity
$n^{2}A(n)$ at $n=1.2\cdot10^{4}$ agrees with $c\,\ell(\beta)$ to three or four digits,
for instance $-1.791522$ against $-1.798107$ at $(3,1)$ and $\beta=\tfrac52$, the residual
gap decreasing like $1/n$. The constant $\ell(\beta)$ equals $0.899$ and $0.693$ at $(3,1)$
for $\beta=\tfrac52$ and $4$, and $0.524$ and $0.368$ at $(5,\tfrac12)$. The Dirichlet
representation \eqref{eq:H_dirichlet} is verified at $(9,1)$ to five digits, for instance
$0.256322$ against $0.256326$ at $\beta=3$. The coefficients $\lambda_j$ at $(41,1)$, where
$|c|=40$, equal $-0.1385$, $-0.0990$, $-0.0020$, $+0.1633$, $+0.4450$, $+0.6440$ at
$j=1,2,4,8,20,40$, in agreement with the bound \eqref{eq:H_lambda_bound} wherever
$\theta_j<\tfrac12$, for instance $\lambda_{40}=0.644\ge0.048$. The sign change of
$\lambda_1=\ell_\infty$ on the line $y=1$ is bracketed by $\lambda_1=+0.0307$ at $x=14$ and
$\lambda_1=-0.0054$ at $x=16$.
\end{numobs}
\galleryentry{I}{The parity switched kernel}
 {$G(n,k)=1$ if $n-k$ is even and $G(n,k)=k/n$ otherwise}
 {regular arithmetic function}
 {$G^{*}(z)=\frac{2z-1}{2(z-1)}$, single zero at $z=\tfrac12$}
 {$\alpha(G)=\eta(G)=\tfrac12$}
 {proved, Theorem~\ref{thm:I_index}}

\rafgalleryfig{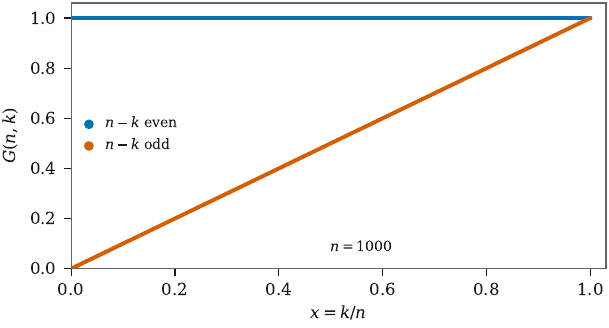}{The kernel at the single rank $n=1000$. It is not of the form $g(k/n)$, and its row splits into two branches according to the parity of $n-k$, the value being one on the even branch and $k/n$ on the odd one. The ratio $k/n$ alone does not determine the value, which is what places this entry outside the univariate class.}{fig:app_I}

The kernel is
\[
G(n,k)=\begin{cases}1 & n-k \text{ even},\\ k/n & n-k \text{ odd},\end{cases}
\]
with $G(n,n)=1$.

The finite probes converge, and their limit is the transform of this two-variable kernel.

\begin{proposition}\label{prop:I_transform}
The finite probes of $G$ converge, locally uniformly on $\Re z<0$, to
\begin{equation}\label{eq:I_transform}
G^{*}(z)=\frac12\Bigl(1+\frac{z}{z-1}\Bigr)=\frac{2z-1}{2(z-1)},
\end{equation}
whose only zero is $z=\tfrac12$, so $\eta(G)=\tfrac12$.
\end{proposition}

\begin{proof}
In the probe $\tfrac{-z}{n}\sum_{k\le n}G(n,k)(k/n)^{-z-1}$ the ranks with $n-k$ even and
those with $n-k$ odd form two interleaved progressions of step two, each contributing one
half of a Riemann sum. The first carries the constant profile $1$ and the second the profile
$t$, so the probe converges to $\tfrac12\bigl(-z\int_0^1t^{-z-1}\,dt\bigr)
+\tfrac12\bigl(-z\int_0^1t\cdot t^{-z-1}\,dt\bigr)=\tfrac12+\tfrac{z}{2(z-1)}$ on $\Re z<0$,
locally uniformly, the two half sums having the same limits as full Riemann sums of step
two.
\end{proof}

The kernel splits into an affine part and a parity part, which is what places it outside the
one-variable class.

\begin{proposition}\label{prop:I_channels}
Identically in $n$ and $k\le n$,
\begin{equation}\label{eq:I_split}
G(n,k)=\frac{1+k/n}{2}+(-1)^{n-k}\,\frac{1-k/n}{2},
\end{equation}
so with $B(n)=\sum_{k\le n}ka_k$, $\widetilde A(n)=\sum_{k\le n}(-1)^{k}a_k$ and
$\widetilde B(n)=\sum_{k\le n}(-1)^{k}ka_k$ the defining equation reads
\begin{equation}\label{eq:I_channels}
S(n)+(-1)^{n}\,\widetilde T(n)=n^{-\beta},
\qquad
S(n)=\frac{A(n)+B(n)/n}{2},
\qquad
\widetilde T(n)=\frac{\widetilde A(n)-\widetilde B(n)/n}{2}.
\end{equation}
The balanced channel $S$ is the affine functional with profile $(1+t)/2$, of transform
\eqref{eq:I_transform} with zero $\tfrac12$, applied to $(a_k)$. The twisted channel
$\widetilde T$ is the affine functional with profile $(1-t)/2$, whose transform
$-1/(2(z-1))$ has no zero and whose diagonal weight vanishes, applied to the twisted
sequence $((-1)^{k}a_k)$.
\end{proposition}

\begin{proof}
For $n-k$ even the right side of \eqref{eq:I_split} is $\tfrac{1+k/n}2+\tfrac{1-k/n}2=1$ and
for $n-k$ odd it is $\tfrac{1+k/n}2-\tfrac{1-k/n}2=k/n$. Summing against $a_k$ and using
$(-1)^{n-k}=(-1)^{n}(-1)^{k}$ gives \eqref{eq:I_channels}.
\end{proof}

\subsection*{The exact local system and its elimination}

The two parities give a local system in two unknowns.

\begin{proposition}\label{prop:I_system}
Put $u_j=a_{2j}$, $v_j=a_{2j-1}$, $U_0=V_0=0$, $U_m=\sum_{j\le m}u_j$,
$V_m=\sum_{j\le m}v_j$, and, for $m\ge2$,
\[
d_e(m)=(2m)^{1-\beta}-(2m-2)^{1-\beta},
\qquad
d_o(m)=(2m-1)^{1-\beta}-(2m-3)^{1-\beta}.
\]
The defining equation is equivalent to $V_1=1$, $U_1=2^{-\beta}-\tfrac12$ together with, for
$m\ge2$,
\begin{align}
2m\,U_m-(2m-2)\,U_{m-1}+(2m-1)\bigl(V_m-V_{m-1}\bigr)&=d_e(m),\label{eq:I_even}\\
(2m-1)\,V_m-(2m-3)\,V_{m-1}+2(m-1)\bigl(U_{m-1}-U_{m-2}\bigr)&=d_o(m),\label{eq:I_odd}
\end{align}
and the partial sums are recovered through $A(2m)=U_m+V_m$ and $A(2m+1)=U_m+V_{m+1}$.
\end{proposition}

\begin{proof}
At rank $n=2m$ the even indices $k=2j$ carry weight $1$ and the odd indices $k=2j-1$ carry
weight $k/n$, so the defining equation reads
$U_m+\tfrac{1}{2m}\sum_{j\le m}(2j-1)v_j=(2m)^{-\beta}$. At rank $n=2m-1$ the parities
exchange roles and $V_m+\tfrac{1}{2m-1}\sum_{j\le m-1}2j\,u_j=(2m-1)^{-\beta}$. Multiplying
by $2m$, respectively $2m-1$, and differencing consecutive values of $m$ removes the
weighted sums, since
$\sum_{j\le m}(2j-1)v_j-\sum_{j\le m-1}(2j-1)v_j=(2m-1)(V_m-V_{m-1})$ and similarly for the
even weights, and produces \eqref{eq:I_even} and \eqref{eq:I_odd}. The initial values are the
ranks $n=1,2$ of the defining equation. Conversely the pair reconstructs the two families of
relations by summation.
\end{proof}

Earlier printed versions of this system mixed two sets of notations and are superseded by
Proposition~\ref{prop:I_system}, which is validated against forward substitution in the
defining equation at both parities.

\begin{proposition}\label{prop:I_elimination}
Eliminating $V$ from \eqref{eq:I_even} and \eqref{eq:I_odd} gives, for $m\ge3$, the exact
third order recurrence $L_J[U](m)=g_\beta(m)$ with
\begin{align}
L_J[U](m)&=(2m-1)\Bigl[(m+1)U_{m+1}-3m\,U_m+(3m-2)\,U_{m-1}-(m-1)\,U_{m-2}\Bigr]\notag\\
&\qquad+2m\,U_m-(2m-2)\,U_{m-1},\label{eq:I_third}\\[1mm]
g_\beta(m)&=d_e(m)+(2m-1)\,\frac{\Delta d(m)-\Delta d(m+1)}{2},
\qquad \Delta d=d_o-d_e.\notag
\end{align}
The frozen characteristic polynomial of \eqref{eq:I_third} is $(r-1)^{3}$, and for every
complex $\gamma$
\begin{equation}\label{eq:I_indicial}
L_J\bigl[m^{-\gamma}\bigr]=-2\,(2\gamma-1)\,m^{-\gamma}+\mathcal O_\gamma(m^{-\gamma-1}),
\end{equation}
so the indicial polynomial is proportional to $2\gamma-1$, the numerator of the transform
\eqref{eq:I_transform}.
\end{proposition}

\begin{proof}
Equation \eqref{eq:I_even} expresses $(2m-1)(V_m-V_{m-1})$ through $U$, and combining it
with \eqref{eq:I_odd} in the form
$(2m-1)V_m-(2m-3)V_{m-1}=(2m-1)(V_m-V_{m-1})+2V_{m-1}$ isolates
\[
V_{m-1}=m\,U_m-(2m-2)\,U_{m-1}+(m-1)\,U_{m-2}+\tfrac12\Delta d(m).
\]
Shifting this identity by one and inserting both expressions back into \eqref{eq:I_even}
eliminates $V$ entirely and yields \eqref{eq:I_third} after collecting terms. Dividing the
bracket by its leading weights, the frozen recurrence is
$r^{3}-3r^{2}+3r-1=(r-1)^{3}$. For the indicial computation write the bracket as
$m\bigl[U_{m+1}-3U_m+3U_{m-1}-U_{m-2}\bigr]+\bigl[U_{m+1}-2U_{m-1}+U_{m-2}\bigr]$. On
$U=m^{-\gamma}$ the first group is a third difference, of size $\mathcal O(m^{-\gamma-2})$
after multiplication by $m$, while the second group expands to
$-\gamma m^{-\gamma-1}+\mathcal O(m^{-\gamma-2})$, so the bracket times $2m-1$ contributes
$-2\gamma m^{-\gamma}$ at leading order. The remaining part contributes
$2m\,U_m-(2m-2)U_{m-1}=2m\bigl(U_m-U_{m-1}\bigr)+2U_{m-1}
=(2-2\gamma)m^{-\gamma}+\mathcal O(m^{-\gamma-1})$. Adding the two gives
\eqref{eq:I_indicial}.
\end{proof}

\subsection*{The generating equation and the index}

The decomposition \eqref{eq:I_channels} shows that a second family lives at the scale of the
indicial root. The twisted channel is forced by $(-1)^{n}\bigl(n^{-\beta}-S(n)\bigr)$, so the
coefficients acquire a parity alternating component. Any reduction to a single scalar
majorant of $|U|+|V|$ discards the cancellation between the two families and cannot see the
threshold, which is why the reduction formerly attempted here is not used. The scalar
generating equation below keeps both families and closes the appendix.

Several forcing exponents are compared in what follows, so the subscript is kept on $A_\beta$ and on the sequences attached to it. Elsewhere in the volume the exponent is fixed by its context and the plain $a_n$ and $A(n)$ of Definition~\ref{def:reg_index_fgv} are used.

\begin{proposition}\label{prop:I_ogf}
Put $A(0)=0$ and
\[
F_\beta(z)=\sum_{n\ge1}A(n)z^n.
\]
For $|z|<1$ one has
\begin{equation}\label{eq:I_ogf}
F_\beta'(z)
-\left(\frac{1}{2(1-z)}
+\frac{1-z}{2(1+z)^2}\right)F_\beta(z)
=\frac{\operatorname{Li}_{\beta-1}(z)}{z}.
\end{equation}
Define
\[
\mu(z)=\sqrt{1-z^2}
\exp\left(\frac{1}{1+z}-1\right),
\]
with the branches positive on $(-1,1)$. Then
\begin{equation}\label{eq:I_integral}
F_\beta(z)=\mu(z)^{-1}
\int_0^z\mu(w)\frac{\operatorname{Li}_{\beta-1}(w)}{w}\,dw.
\end{equation}
\end{proposition}

\begin{proof}
Summation by parts in the two channels gives
\[
2nS(n)=2nA(n)-\sum_{k<n}A(k)
\]
and
\[
2n\widetilde T(n)
=\sum_{k<n}(-1)^k(2n-2k-1)A(k).
\]
Consequently
\[
\sum_{n\ge1}2nS(n)z^n
=2zF_\beta'(z)-\frac{z}{1-z}F_\beta(z)
\]
and
\[
\sum_{n\ge1}2n(-1)^n\widetilde T(n)z^n
=-\frac{z(1-z)}{(1+z)^2}F_\beta(z).
\]
The second equality uses
\[
\sum_{r\ge1}(2r-1)(-z)^r
=-\frac{z(1-z)}{(1+z)^2}.
\]
Since $S(n)+(-1)^n\widetilde T(n)=n^{-\beta}$, division by $2z$
gives \eqref{eq:I_ogf}. Direct differentiation gives
\[
\frac{\mu'(z)}{\mu(z)}
=-\frac{1}{2(1-z)}-\frac{1-z}{2(1+z)^2}.
\]
Multiplication by $\mu$ and integration from zero give
\eqref{eq:I_integral}.
\end{proof}

The generating function has an irregular endpoint, and its behavior there governs the
transition.

\begin{lemma}\label{lem:I_irregular}
Put
\[
q(z)=\frac{1}{2(1-z)}+\frac{1-z}{2(1+z)^2},
\qquad
H(z)=\frac{\exp(1-1/(1+z))}{\sqrt{1-z^2}}.
\]
Let $f$ be holomorphic in a neighborhood of $-1$, and let $Y$ solve
$Y'-qY=f$ in the unit disc. On each closed sector with vertex $-1$ that
avoids the two tangent rays, and for every integer $N\ge2$, there are a
constant $C_N$ and a polynomial $P_N(s)$ in $s=1+z$, beginning with
$s^2$, such that
\[
Y(-1+s)=C_NH(-1+s)+P_N(s)+\mathcal O(s^{N+1}).
\]
Let $\Gamma_-$ be the part inside $|z+1|<r$ of a Darboux contour that
runs on the two sides of these sectors, and define
\[
c_-(n)=\frac{1}{2\pi i}\int_{\Gamma_-}Y(z)z^{-n-1}\,dz.
\]
Then $c_-(n)=\mathcal O(n^{-1/2})$.
\end{lemma}

\begin{proof}
In the variable $s=1+z$ one has
\[
q(z)=s^{-2}-\frac{1}{2s}+\frac{1}{2(2-s)}.
\]
Multiplication of the equation by $s^2$ determines the coefficients of
$P_N$ recursively. The first two coefficients vanish and the coefficient
of $s^2$ is $-f(-1)$. After subtraction of this polynomial, variation of
constants gives
\[
\frac{Y-P_N}{H}
=C_N+\int H^{-1}\mathcal O(s^{N-1})\,ds.
\]
Repeated integration by parts, using $(H^{-1})'=-qH^{-1}$, gives the
stated remainder on every such sector.

Put $w=-z$. The homogeneous term satisfies the exact identity
\[
H(-w)=(1+w)^{-1/2}(1-w)^{-1/2}
\exp\left(-\frac{w}{1-w}\right).
\]
The generating formula
\[
\sum_{n\ge0}L_n^{(\nu)}(x)w^n
=(1-w)^{-\nu-1}
\exp\left(-\frac{xw}{1-w}\right)
\]
and the fixed parameter Fej\'er estimate
\[
L_n^{(\nu)}(x)
=\frac{e^{x/2}}{\sqrt\pi\,x^{\nu/2+1/4}}
n^{\nu/2-1/4}
\left(
\cos\left(2\sqrt{nx}-\frac{\nu\pi}{2}-\frac\pi4\right)
+\mathcal O(n^{-1/2})
\right)
\]
hold for fixed real $\nu$ and fixed $x>0$. With
$\nu=-1/2$ and $x=1$ they give
\[
L_n^{(-1/2)}(1)
=\frac{e^{1/2}}{\sqrt\pi}n^{-1/2}
\left(\cos(2\sqrt n)+\mathcal O(n^{-1/2})\right).
\]
The factor $(1+w)^{-1/2}$ is holomorphic at $w=1$. Its other endpoint is
$z=1$, which is treated separately. The polynomial jets continue through
$-1$. Since $N$ is arbitrary, their remainders give no larger coefficient.
This proves the bound for the endpoint $-1$.
\end{proof}

The three regimes of the parity kernel follow.

\begin{theorem}\label{thm:I_index}
For the kernel of this appendix the following assertions hold.

If $\beta<1/2$, then
\[
A(n)=\frac{2(\beta-1)}{2\beta-1}n^{-\beta}
+o(n^{-\beta}).
\]
At $\beta=1/2$, one has
\[
A(n)=\frac12n^{-1/2}\log n+\mathcal O(n^{-1/2}).
\]
If $\beta>1/2$, then
\[
A(n)=\mathcal O_\beta(n^{-1/2}).
\]
Consequently the kernel is a RAF with
\[
\alpha(G)=\eta(G)=\frac12.
\]
Transparency holds below $1/2$, absorption holds at and above $1/2$,
and the index is sharp.
\end{theorem}

\begin{proof}
Put $t=1-z$. As $z$ tends to one in a dented neighborhood,
\[
\mu(z)=\sqrt2e^{-1/2}t^{1/2}(1+\mathcal O(t)),
\qquad
\mu(z)^{-1}=\frac{e^{1/2}}{\sqrt2}t^{-1/2}(1+\mathcal O(t)).
\]
For every real $\beta$, the function
$\operatorname{Li}_{\beta-1}(z)/z$ is holomorphic at $z=-1$. Equation
\eqref{eq:I_ogf} continues $F_\beta$ through the unit circle away from
$1$ and $-1$. Thus the hypotheses of Lemma~\ref{lem:I_irregular} hold.
For $\beta<1/2$, the polylogarithm has the expansion
\[
\operatorname{Li}_{\beta-1}(z)
=\Gamma(2-\beta)t^{\beta-2}(1+\mathcal O(t))+\mathcal O(1).
\]
Equation \eqref{eq:I_integral} gives
\[
F_\beta(z)
=\frac{\Gamma(2-\beta)}{1/2-\beta}
t^{\beta-1}(1+o(1))+\mathcal O(t^{-1/2}).
\]
The coefficient transfer
\[
[z^n](1-z)^{\beta-1}
\sim\frac{n^{-\beta}}{\Gamma(1-\beta)}
\]
and Lemma~\ref{lem:I_irregular} give
\[
A(n)=\frac{\Gamma(2-\beta)}{(1/2-\beta)\Gamma(1-\beta)}n^{-\beta}
+o(n^{-\beta}).
\]
The coefficient is $2(\beta-1)/(2\beta-1)$, which is $1/G^{*}(\beta)$ by
\eqref{eq:I_transform}.

At $\beta=1/2$, equation \eqref{eq:I_integral} gives
\[
F_{1/2}(z)
=\frac{\sqrt\pi}{2}t^{-1/2}\log(1/t)+\mathcal O(t^{-1/2}).
\]
Coefficient transfer and Lemma~\ref{lem:I_irregular} give the endpoint
formula in the statement.

For $\beta>1/2$, the integral in \eqref{eq:I_integral} converges at one.
The endpoint $1$ therefore contributes $\mathcal O(n^{-1/2})$. The only other
endpoint on the circle of convergence is $-1$, and
Lemma~\ref{lem:I_irregular} gives the same bound there.

The first formula proves transparency for every $\beta<1/2$, with
$\Xi(\beta)=1/G^*(\beta)$. The logarithm proves that $\beta=1/2$ is not
transparent. Thus $\alpha(G)=1/2$ and the index is sharp. The last two formulas
give $A(n)=\mathcal O_\eps(n^{-1/2+\eps})$ for every $\beta\ge1/2$, which
is absorption.
\end{proof}

\begin{remark}\label{rem:I_twisted_phase}
The endpoint $-1$ contributes an alternating term of order $n^{-1/2}$
whose phase is $2\sqrt n+\mathcal O(1)$. The term is visible in $A(n)$ itself. It
is not confined to the two parity subsequences. In the notation of
Proposition~\ref{prop:I_system}, the corresponding phase of $U_m-V_m$ is
$2\sqrt{2m}+\mathcal O(1)$. There is no logarithmic rotation, and no factor
$(-1)^m$ should be placed before $U_m-V_m$.
\end{remark}

\begin{numobs}\label{numobs:I_check_revised}
Exact rational forward substitution at $\beta=3$ gives zero residual for
the channel identity, the local system, the eliminated identity for
$V_{m-1}$, the third order recurrence, and \eqref{eq:I_ogf} through rank
$120$. At $80$ decimal digits and through rank $500$, the largest scaled
residual is below $1.5\cdot10^{-75}$ at
$\beta\in\{\tfrac15,\tfrac12,\tfrac32\}$.

At $\beta=3$ and $500\le n\le4000$, the minimum of
$|a_n|\sqrt n$ is $0.00106138605925$ at $n=3168$, and the maximum is
$1.47860899333$ at $n=658$. A fit at large rank gives the phase
$2\sqrt n$ with coefficient of determination greater than
$0.99999999$. The fitted amplitude of $U_m-V_m$ is
$1.47830\ldots$. The corresponding alternating amplitude in
$\sqrt nA(n)$ is $0.73915\ldots$, so the second mode is visible in the
partial sums.

At $\beta=1/2$, fitting $\sqrt nA(n)$ to a constant, a multiple of
$\log n$, and the alternating phase gives a logarithmic coefficient
$0.50000\ldots$, in agreement with Theorem~\ref{thm:I_index}.
\end{numobs}
\galleryentry{J}{The greatest common divisor kernel}
 {$G(n,k)=\frac{1}{2}\!\left(1+\frac{\gcd(n,k)}{n}\right)$}
 {conjectured regular arithmetic function}
 {$G^{*}(z)=\tfrac12$, constant, no zero, so $\eta$ is not defined}
 {$\alpha(G)=1$ conjectured, $\eta$ not defined}
 {open, stability of $(I-T)^{-1}$ at order $n^{-1+\eps}$ is unproved, Conjecture~\ref{conj:J_main}}

\rafgalleryfig{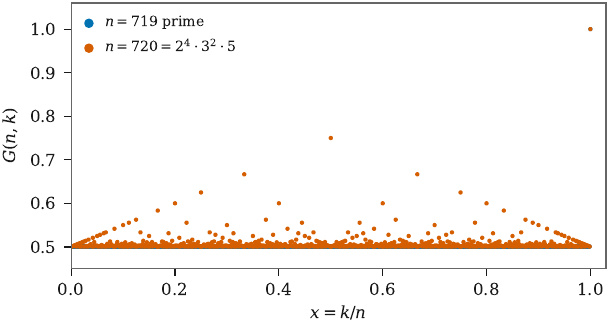}{The kernel at two consecutive ranks. At the prime rank the greatest common divisor is one below the diagonal, so the row is flat at one half with a single point at the diagonal. At the next rank, which has thirty divisors, the row carries a point above one half at every $k$ sharing a factor with $n$. Two adjacent ranks give unrelated rows, and the section family that Conjecture~\ref{conj:J_main} leaves open is that arithmetic.}{fig:app_J}

The arithmetic kernel $G(n,k)=\frac{1}{2}(1+\gcd(n,k)/n)$ has constant Mellin transform
$G^{*}(z)=\tfrac12$, without zero, so no analytic index is available.

\begin{proposition}\label{prop:J_divisor}
Put $A_d(m)=\sum_{\ell\le m}a_{d\ell}$, so that $A_1=A$. The identity
\[
 \gcd(n,k)=\sum_{d\mid n,\ d\mid k}\varphi(d)
\]
turns the defining equation into the exact relation
\begin{equation}\label{eq:J_core}
A(n)+\frac{1}{n}\sum_{d\mid n}\varphi(d)\,A_d(n/d)=2n^{-\beta}.
\end{equation}
\end{proposition}

\begin{proof}
Summing $a_k\gcd(n,k)$ and exchanging the order, each divisor $d$ of $n$ collects the ranks
$k$ that are multiples of $d$, whose partial sum up to $n$ is $A_d(n/d)$.
\end{proof}

The divisor structure of the kernel gives a double Dirichlet series in closed form.

\begin{proposition}\label{prop:J_dirichlet}
Let $\mathcal A_d(w)=\sum_{m\ge1}A_d(m)\,m^{-w}$. For $\Re s$ large enough that all series
converge absolutely,
\begin{equation}\label{eq:J_dirichlet}
\sum_{n\ge1}A(n)\,n^{-s}
+\sum_{d\ge1}\frac{\varphi(d)}{d^{\,s+1}}\,\mathcal A_d(s+1)
=2\,\zeta(s+\beta).
\end{equation}
\end{proposition}

\begin{proof}
Multiply \eqref{eq:J_core} by $n^{-s}$ and sum. In the double sum write $n=dm$, so that
$n^{-s-1}\varphi(d)A_d(n/d)$ regroups as $\varphi(d)d^{-s-1}\cdot A_d(m)m^{-s-1}$, and the
forcing gives $2\zeta(s+\beta)$. The interchanges are legitimate by absolute convergence,
the coefficients being of at most polynomial growth.
\end{proof}

The coupling weight $\varphi(d)d^{-s-1}$ is multiplicative, with
$\sum_d\varphi(d)d^{-s-1}=\zeta(s)/\zeta(s+1)$, so \eqref{eq:J_dirichlet} is the equation of
an infinite multiplicatively graded system in the section series $\mathcal A_d$, of which the
visible part $A=A_1$ is one component.

\begin{conjecture}\label{conj:J_main}
$\alpha(G)=1$. For $0<\beta<1$, $A(n)\sim 2n^{-\beta}$, and for $\beta\ge 1$,
$A(n)=\mathcal{O}_\eps(n^{-1+\eps})$.
\end{conjecture}

\begin{proofstatus}{Propositions~\ref{prop:J_divisor} and
\ref{prop:J_dirichlet} are proved. Summation by parts also gives the exact
scalar reformulation
\[
(Tx)(n)=\frac1{2n}\sum_{k<n}\bigl(\gcd(n,k+1)-\gcd(n,k)\bigr)x(k),
\]
under which the defining equation is $A=p_\beta+TA$, with
$p_\beta(n)=n^{-\beta}$. No estimate for the section family is assumed here.
For every fixed $\gamma>0$, a second summation by parts and
$\gcd(n,j)=\sum_{d\mid n,\ d\mid j}\varphi(d)$ give
\[
(Tp_\gamma)(n)=\frac12(n-1)^{-\gamma}
 +\mathcal O_{\gamma,\eps}(n^{-1+\eps}).
\]
Thus the leading multiplier is one half when $0<\gamma<1$, whereas
$Tp_\gamma(n)=\mathcal O_{\gamma,\eps}(n^{-1+\eps})$ when $\gamma\ge1$.
This proves that the threshold signature is present in $T$, but it does not
control its iterates. The missing step is the stability of the resolvent
$(I-T)^{-1}$ on sequences of order $n^{-1+\eps}$.

Two sufficient routes are known, but no equivalence is asserted. It is enough
to prove $A_d(m)\ll_{\beta,\eps}d^{-1+\eps}m^\eps$ for every $d\ge2$.
The coefficient estimate
$a_n=\mathcal O_{\beta,\eps}(n^{-1+\eps})$, the kernel analogue of the weak
Hardy--Littlewood--Ramanujan\index[terms]{Hardy--Littlewood--Ramanujan criterion}
bound in Definition~\ref{def:HLR}, is a stronger sufficient condition, since summing
it along multiples implies the displayed section estimate. Either route gives
$A(n)=2n^{-\beta}+\mathcal O_{\beta,\eps}(n^{-1+\eps})$ and hence both regimes
of the conjecture.

No estimate controlling the sections by $A$ alone is established. Whether
they tend to nonzero limits is open, and the numerical values still drift at
$m=5\cdot10^7$. Numerical Observation~\ref{numobs:J_check} is evidence only.
Nothing outside this appendix assumes Conjecture~\ref{conj:J_main} or either
sufficient bound.}
\end{proofstatus}

\begin{numobs}\label{numobs:J_check}
The relation \eqref{eq:J_core} reproduces forward substitution to $5\cdot10^{-16}$ over
$n\le800$, and the resummed identity \eqref{eq:J_dirichlet} holds to $5\cdot10^{-15}$ at
$s=3$ under truncation-consistent evaluation. Below the threshold the transparency constant
is met, $n^{\beta}A(n)/2=1.0043$ at $n=4\cdot10^{3}$ for $\beta=\tfrac12$. Above it
$nA(n)$ stays bounded without settling, its values at $n=10^{3}$, $2\cdot10^{3}$,
$4\cdot10^{3}$ being $0.43$, $0.25$, $0.29$ at $\beta=\tfrac32$, while the sections
vary slowly over that window, $A_2$ near $-0.28$ and $A_3$ near $-0.13$ for
$500\le m\le10^{3}$, values that continue to move at larger cutoffs and are not limits.
\end{numobs}
\galleryentry{K}{The quadratic rational kernel}
 {$G(n,k)=\frac{n^{2}+k}{n^{2}+n}$}
 {regular arithmetic function}
 {not computed here, the gallery table of \S\ref{sec:raf_gallery} records $G^{*}\equiv1$, no zero, so $\eta$ is not defined}
 {$\alpha(G)=2$, $\eta$ not defined}
 {proved, Theorem~\ref{thm:K_main}}

The index of this kernel is two as well, and the regimes are explicit on both sides.

\begin{theorem}\label{thm:K_main}
$\alpha(G)=2$. For $0<\beta<2$, $A(n)=n^{-\beta}+o(n^{-\beta})$. For $\beta\ge 2$, $A(n)=\mathcal{O}(n^{-2})$.
\end{theorem}

\begin{proof}
After expanding and applying Abel summation ($\sum k a_k = nA(n)-AA(n-1)$), the defining relation simplifies exactly to:
\begin{equation}\label{eq:K_core}
A(n) = n^{-\beta} + \frac{AA(n-1)}{n(n+1)}, \qquad AA(n):=\sum_{j=1}^n A(j).
\end{equation}

All partial sums are positive. Indeed $A(1)=1$, and if $A(j)>0$ for every $j<n$ then
$AA(n-1)>0$, so \eqref{eq:K_core} gives $A(n)>n^{-\beta}>0$. The feedback term is therefore
positive as well, and the two ranges are read off it.

For $\beta\ge2$ the auxiliary sequence $S(n)=\sum_{j\le n}A(j)$ satisfies
$S(n)\le S(n-1)\bigl(1+1/(n(n+1))\bigr)+n^{-\beta}$. Since $\prod_k\bigl(1+1/(k(k+1))\bigr)$
converges and $\sum_kk^{-\beta}$ converges, $S(n)$ is bounded, and \eqref{eq:K_core} gives
$A(n)\le n^{-\beta}+S(n-1)/(n(n+1))=\mathcal{O}(n^{-2})$.

For $\beta<2$ set $M(n)=\max_{j\le n}j^{\beta}A(j)$, so that $M(1)=1$. By \eqref{eq:K_core},
\[
n^{\beta}A(n)\le1+\frac{M(n-1)}{n(n+1)}\,n^{\beta}\sum_{j<n}j^{-\beta},
\]
and the last factor tends to zero in each of the three subranges. For $\beta<1$ it is at most
$n^{\beta}\cdot n^{1-\beta}/(1-\beta)$ divided by $n(n+1)$, hence at most
$1/\bigl((1-\beta)(n+1)\bigr)$. For $\beta=1$ it is at most $(1+\log n)/(n+1)$. For
$1<\beta<2$ it is at most $\zeta(\beta)n^{\beta}/(n(n+1))$, which tends to zero because
$\beta<2$. In each case there is a rank $n_1$ beyond which the factor is at most $\tfrac12$,
so that $M(n)\le\max\bigl(M(n_1),2\bigr)$ for every $n$ by induction, and $M$ is bounded. The
same three estimates then give $AA(n-1)/(n(n+1))=o(n^{-\beta})$, whence
$A(n)=n^{-\beta}+o(n^{-\beta})$.

The two ranges together give $\alpha(G)=2$.
\end{proof}

The homogeneous equation is solvable in closed form, and its decay exhibits the index directly.

\begin{proposition}\label{prop:K_homogeneous}
The homogeneous equation $\sum_{k\le n}a_kG(n,k)=0$ reduces to the one-step recurrence
$A_0(n)=\frac{n^{2}-n+1}{n^{2}+n}A_0(n-1)$ with $A_0(1)=1$, and
\begin{equation}\label{eq:K_homogeneous}
A_0(n)=\prod_{k=2}^{n}\frac{k^{2}-k+1}{k(k+1)}
=\frac{2\cosh\!\left(\tfrac{\pi\sqrt3}{2}\right)}{\pi}\cdot
\frac{\Gamma\!\left(n+e^{i\pi/3}\right)\Gamma\!\left(n+e^{-i\pi/3}\right)}{\Gamma(n+1)\,\Gamma(n+2)},
\end{equation}
an identity at every rank. By Stirling\index[terms]{Stirling's formula}\index[names]{Stirling, J.}, $A_0(n)\sim C_0n^{-2}$, where
$C_0=2\cosh(\pi\sqrt3/2)/\pi=4.856379584\ldots$
\end{proposition}

\begin{proof}
Factor $k^{2}-k+1=(k-\omega)(k-\bar\omega)$ with $\omega=e^{i\pi/3}$, and note
$1-\omega=e^{-i\pi/3}$. Telescoping the three products gives
$\prod_{k=2}^{n}(k-\omega)=\Gamma(n+1-\omega)/\Gamma(2-\omega)$,
$\prod_{k=2}^{n}k=\Gamma(n+1)$ and $\prod_{k=2}^{n}(k+1)=\Gamma(n+2)/2$, so the product equals
$2\,\Gamma(n+1-\omega)\Gamma(n+1-\bar\omega)/\bigl(\Gamma(2-\omega)\Gamma(2-\bar\omega)
\Gamma(n+1)\Gamma(n+2)\bigr)$, and $n+1-\omega=n+e^{-i\pi/3}$. The prefactor is real, since
$2-\omega=\tfrac32-i\tfrac{\sqrt3}2$, so that
$\Gamma(2-\omega)\Gamma(2-\bar\omega)=\bigl|\Gamma(\tfrac32+i\tfrac{\sqrt3}2)\bigr|^{2}$, and
$\Gamma(\tfrac32+it)\Gamma(\tfrac32-it)=(\tfrac14+t^{2})\,\pi/\cosh(\pi t)$ by the reflection
formula. At $t=\sqrt3/2$ the factor $\tfrac14+t^{2}$ is one, leaving
$\pi/\cosh(\pi\sqrt3/2)$, whose reciprocal times two is $C_0$. Since
$e^{-i\pi/3}+e^{i\pi/3}=1$, the Gamma quotient is asymptotic to $n^{1-3}=n^{-2}$.
\end{proof}

The homogeneity principle\index[terms]{homogeneity principle} of Chapter~\ref{chap:principles} would read $\alpha(G)=2$ off this
decay directly. That principle is open in general, so the case analysis above, and not
\eqref{eq:K_homogeneous}, is what proves the index here.

\begin{numobs}[Numerical control]\label{numobs:K_check}
The recurrence \eqref{eq:K_core} reproduces forward substitution in the defining equation to
$1.7\cdot10^{-15}$ over $n\le400$ at $\beta=\tfrac1{10}$, $\tfrac32$ and $3$. Below the index
$n^{\beta}A(n)$ tends to one, its values at $n=10^{3}$, $10^{4}$, $10^{5}$, $2\cdot10^{5}$ being
$1.002013$, $1.000201$, $1.000020$, $1.000010$ at $\beta=\tfrac12$ and $1.108919$, $1.034938$,
$1.011094$, $1.007849$ at $\beta=\tfrac32$, the approach slowing as $\beta$ rises toward two, as
the error term $\mathcal O(n^{\beta-2})$ requires. Above the index $n^{2}A(n)$ is bounded and
settles, on $2.064332$ at $\beta=\tfrac52$ and on $1.886397$ at $\beta=3$. The closed form
\eqref{eq:K_homogeneous} agrees with the product to $6\cdot10^{-45}$ at
$n=1,2,3,5,10,50,200$ in forty-digit arithmetic, and $n^{2}A_0(n)$ reaches $4.856331$ at
$n=2\cdot10^{5}$ against $C_0=4.856380$.
\end{numobs}

\galleryentry{L}{The self gauged logarithmic kernel}
 {$G(n,k)=\frac{\log(n+k)}{\log(2n)}$}
 {regular arithmetic function}
 {$G^{*}(z)=1$, constant, no zero, so $\eta$ is not defined}
 {$\alpha(G)=1$, $\eta$ not defined}
 {proved, Theorem~\ref{thm:L_main}}

\rafgalleryfig{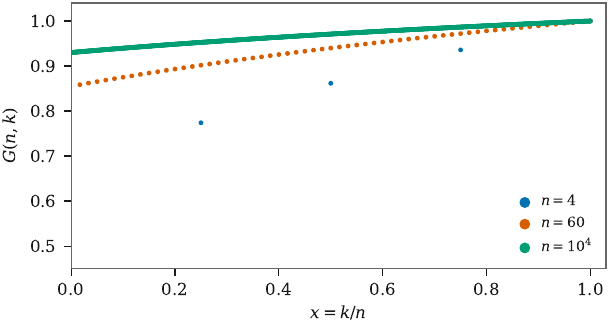}{The kernel at three ranks. Unlike a univariate profile it moves with $n$, rising toward the constant one as the gauge $\log(2n)$ grows, and the first zero of the effective transform $1+h^{*}/L_n$ moves with it. The limiting transform is constant and has no zero, so the index comes from the equation, and Theorem~\ref{thm:L_main} identifies it with the limit of those moving thresholds. The approach is logarithmic, which is why equilibrium is reached only in the limit.}{fig:app_L}

The kernel is increasing in its second variable with unit diagonal, and its transform is
constant. The index cannot come from the transform, and the mechanism that produces it is a
threshold that moves with the rank.

\begin{proposition}\label{prop:L_transform}
The finite probes of $G$ converge, locally uniformly on $\Re z<0$, to $G^{*}(z)=1$. The
analytic index $\eta(G)$ is undefined.
\end{proposition}

\begin{proof}
The kernel increases from $G(n,1)$ to $G(n,n)=1$ along each row, so the total row variation is
$\bigl(\log(2n)-\log(n+1)\bigr)/\log(2n)=\mathcal O(1/\log n)$. On a compact subset of
$\Re z<0$ the powers $(k/n)^{-z}$ in the finite probe have modulus at most one, so the probe
differs from $G(n,n)=1$ by $\mathcal O_z(1/\log n)$ and tends to one.
\end{proof}

The kernel is an exact perturbation of the constant one, by a logarithmic profile divided by a
growing gauge.

\begin{proposition}\label{prop:L_gauge}
Put $h(t)=\log\frac{1+t}{2}$ and $L_n=\log(2n)$. Then, identically in $n$ and $k\le n$,
\begin{equation}\label{eq:L_split}
G(n,k)=1+\frac{h(k/n)}{L_n},
\end{equation}
so the defining equation reads $A(n)+L_n^{-1}\sum_{k\le n}a_k\,h(k/n)=n^{-\beta}$. The
transform of the perturbing profile is
\begin{equation}\label{eq:L_hstar}
h^{*}(z)=-z\int_0^1 h(t)\,t^{-z-1}\,dt
=-\log2-z\sum_{m\ge1}\frac{(-1)^{m+1}}{m(m-z)}
=-\int_0^1\frac{t^{-z}}{1+t}\,dt,
\end{equation}
the first two forms on $\Re z<0$ and by continuation, the third valid on $\Re z<1$,
meromorphic with simple poles at the positive integers, the residue at $z=1$ being $1$.
Consequently the finite-rank effective transform
\[G^{*}_n(z)=1+\frac{h^{*}(z)}{L_n}
\]
has its first positive zero at
\begin{equation}\label{eq:L_threshold}
z_n=1-\frac{1}{L_n}+\mathcal O(L_n^{-2}),
\end{equation}
and $z_n\to1$ as $n\to\infty$.
\end{proposition}

\begin{proof}
Identity \eqref{eq:L_split} is $\log(n+k)=\log(2n)+\log\frac{1+k/n}{2}$ divided by $\log(2n)$.
For \eqref{eq:L_hstar}, split $h=\log(1+t)-\log2$, use $\int_0^1t^{-z-1}\,dt=-1/z$ on $\Re z<0$
and expand $\log(1+t)=\sum_{m\ge1}(-1)^{m+1}t^{m}/m$, each term contributing
$-z(-1)^{m+1}/(m(m-z))$. The series converges locally uniformly away from the positive
integers and exhibits the poles. The third form follows from integration by parts, the
boundary terms vanishing since $h(1)=0$, and the integral converges and is analytic on
$\Re z<1$, where it agrees with the series by continuation. Near $z=1$ the term $m=1$ gives
$z/(z-1)$, of residue $1$, so $h^{*}(z)\sim1/(z-1)$ there, and $1+h^{*}(z)/L_n=0$ reads
$1/(z-1)\simeq-L_n$, which is \eqref{eq:L_threshold}.
\end{proof}

The Abel form of that perturbation has a sign, and the sign is what drives the estimate.

\begin{lemma}\label{lem:L_positive}
Put $\delta_n(k)=h\bigl(\tfrac{k+1}{n}\bigr)-h\bigl(\tfrac kn\bigr)
=\log\bigl(1+\tfrac{1}{n+k}\bigr)$ for $1\le k<n$. The defining equation is equivalent to
\begin{equation}\label{eq:L_abel}
A(n)=n^{-\beta}+\frac{1}{L_n}\sum_{k<n}A(k)\,\delta_n(k),
\qquad
\frac{1}{2(n+k)}\le\delta_n(k)\le\frac{1}{n+k},
\end{equation}
and consequently $A(n)>0$ for every $n$ and every real $\beta$, with
$A(n)\ge\delta_n(1)/L_n\ge1/(4nL_n)$ for $n\ge2$.
\end{lemma}

\begin{proof}
Abel summation of $\sum_{k\le n}a_kh(k/n)$ against the partial sums gives
\[
 A(n)h(1)-\sum_{k<n}A(k)\delta_n(k),
\]
and $h(1)=0$. The bounds on $\delta_n$ are $\log(1+x)\in[\tfrac x2,x]$ for
$0\le x\le1$. Positivity follows by induction from
$A(1)=1$, every term on the right of \eqref{eq:L_abel} being positive, and the term
$k=1$ alone gives the lower bound.
\end{proof}

The index of the kernel follows, together with the three regimes.

\begin{theorem}\label{thm:L_main}
$\alpha(G)=1$. More precisely the following hold.
\begin{enumerate}[label=(\roman*)]
\item For every $\beta<1$,
\[
n^{\beta}A(n)=1-\frac{h^{*}(\beta)}{L_n}+\mathcal O_\beta\bigl(L_n^{-2}\bigr),
\]
so transparency holds with constant $1/G^{*}(\beta)=1$ and the finite-rank constant is
$1/G^{*}_n(\beta)$ at the order $1/L_n$.
\item For every $\beta>1$ there is $B_\beta$ with $A(n)\le n^{-\beta}+B_\beta\,n^{-1}$, and at
$\beta=1$ one has $A(n)=\mathcal O(n^{-1}\log n)$, so absorption holds with no logarithmic
loss above the threshold.
\item For every $\beta>1$, $n^{\beta}A(n)\ge n^{\beta-1}/(4L_n)\to\infty$, so transparency
fails at every exponent above one.
\end{enumerate}
\end{theorem}

\begin{proof}
Everything rests on \eqref{eq:L_abel} and on the exact harmonic mass
\begin{equation}\label{eq:L_harmonic}
\sum_{k=1}^{n-1}\frac{1}{k(n+k)}=\frac{H_{n-1}+H_n-H_{2n-1}}{n},
\qquad
H_{n-1}+H_n-H_{2n-1}\le\log n+\gamma-\log2+\frac{3}{2n},
\end{equation}
the identity by partial fractions and telescoping, the bound by the standard two-sided
estimates for $H_m$ together with $n(n-1)/(2n-1)\le n/2$.

For (ii) fix $\beta>1$ and prove $A(n)\le n^{-\beta}+Bn^{-1}$ by induction, with $B$ chosen
below. Inserting the hypothesis into \eqref{eq:L_abel} and using $\delta_n(k)\le1/(n+k)$,
\[
A(n)\le n^{-\beta}+\frac{1}{L_n}\Bigl[\sum_{k<n}\frac{k^{-\beta}}{n+k}
+B\sum_{k<n}\frac{1}{k(n+k)}\Bigr]
\le n^{-\beta}+\frac{\zeta(\beta)+B\bigl(\log n+\gamma-\log2+\tfrac3{2n}\bigr)}{nL_n}.
\]
Since $L_n=\log n+\log2$, the requirement that the last expression be at most
$n^{-\beta}+B/n$ reads $\zeta(\beta)\le B\bigl(2\log2-\gamma-\tfrac3{2n}\bigr)$, and
$2\log2-\gamma=0.809\ldots$, so $B=2\zeta(\beta)$ works for $n\ge8$, and enlarging $B$ to
cover the ranks below $8$ closes the induction. At $\beta=1$ the hypothesis
$A(k)\le ML_k/k$ gives, with $L_k\le L_n$ and \eqref{eq:L_harmonic},
$A(n)\le n^{-1}+M(\log n-0.11)/n$, which is at most $ML_n/n$ for $M\ge2$ by the same
margin.

For (i) fix $\beta<1$. The bootstrap $|A(k)|\le Mk^{-\beta}$ closes first. Inserting it into
\eqref{eq:L_abel},
\[
A(n)\le n^{-\beta}+\frac{M}{L_n}\sum_{k<n}\frac{k^{-\beta}}{n+k},
\qquad
\sum_{k<n}\frac{k^{-\beta}}{n+k}\le C_\beta\,n^{-\beta},
\]
the last bound by splitting at $n/2$, the near part being at most
$n^{-1}\sum_{k\le n/2}k^{-\beta}$ and the far part at most
$2^{\beta}n^{-\beta}\log2$ when $\beta\ge0$, with the same conclusion for $\beta<0$ since
then $k^{-\beta}\le n^{-\beta}$ termwise. So $A(n)\le n^{-\beta}(1+C_\beta M/L_n)$, and for
$n$ beyond the rank where $L_n\ge2C_\beta$ the choice $M=2$, enlarged to cover the initial
ranks, closes the induction as in the absorption case. With the crude bound established,
write $A(k)=k^{-\beta}+E(k)$, $|E(k)|\le C k^{-\beta}/L_k$, which is the bound just proved
reinserted once into \eqref{eq:L_abel}. Then
\[
\sum_{k<n}A(k)\delta_n(k)
=\sum_{k<n}\frac{k^{-\beta}}{n+k}+\mathcal O\Bigl(\sum_{k<n}\frac{k^{-\beta}}{L_k(n+k)}\Bigr)
+\mathcal O\Bigl(\sum_{k<n}\frac{k^{-\beta}}{(n+k)^{2}}\Bigr),
\]
using $\delta_n(k)=\tfrac{1}{n+k}+\mathcal O((n+k)^{-2})$. The main sum is a Riemann sum,
\[
\sum_{k<n}\frac{k^{-\beta}}{n+k}
=n^{-\beta}\int_0^1\frac{t^{-\beta}}{1+t}\,dt+\mathcal O_\beta\bigl(n^{-1}+n^{-\beta-1}\bigr)
=-h^{*}(\beta)\,n^{-\beta}+\mathcal O_\beta(n^{\max(-1,-\beta-1)}),
\]
by \eqref{eq:L_hstar}. Monotonicity controls the discretization, and the edge ranks
$k\le\sqrt n$ contribute $\mathcal O(n^{-(1+\beta)/2})$. The two error sums are
$\mathcal O(n^{-\beta}/L_n)$, the first because the ranks $k\le\sqrt n$ contribute
$\mathcal O(n^{(\beta-1)/2})$ relatively and the ranks beyond carry $1/L_k\le2/L_n$, the
second being $\mathcal O(n^{-\beta-1}\log n)$. Hence
\[
 \sum_{k<n}A(k)\delta_n(k)=-h^{*}(\beta)\,n^{-\beta}\bigl(1+\mathcal O(1/L_n)\bigr),
\]
and \eqref{eq:L_abel} gives (i).

For (iii), Lemma~\ref{lem:L_positive} gives $A(n)\ge1/(4nL_n)$ outright, so for $\beta>1$
the ratio $A(n)/n^{-\beta}$ tends to infinity and $A(n)=n^{-\beta}/G^{*}(\beta)+o(n^{-\beta})$
is impossible. Combining, transparency holds for every $\beta<1$ by (i) and fails for every
$\beta>1$ by (iii), so the supremum of Definition~\ref{def:reg_index_fgv} is exactly one,
and (ii) is the absorption statement at the index.
\end{proof}

\begin{remark}[The self-gauged reading]\label{rem:L_selfgauged}
The theorem realizes the limit regime of Definition~\ref{def:stability_equilibrium}. The
effective kernel at rank $n$ is the member $\lambda=1/L_n$ of the fixed family
$1+\lambda h$, whose index moves like $1-\lambda$, and Theorem~\ref{thm:L_main} identifies the index of the self-gauged kernel with the limit of the moving thresholds \eqref{eq:L_threshold}. No diagonal argument is
needed, the signed perturbation $h\le0$ making the feedback positive and the three regimes
accessible to direct estimates. The finite-rank constant of (i) is the one to test
numerically, the convergence to the limit constant being logarithmic.
\end{remark}

\begin{numobs}\label{numobs:L_check}
Positivity holds at every rank computed. The transparency constant at rank $n$ matches (i),
at $\beta=\tfrac25$, where $h^{*}(\beta)=-1.268758$ by both the series and the integral of
\eqref{eq:L_hstar}, the measured $n^{\beta}A(n)$ at $n=2\cdot10^{3}$, $4\cdot10^{3}$,
$8\cdot10^{3}$ is $1.189363$, $1.172398$, $1.157854$ against the resummed prediction
$1.180599$, $1.164380$, $1.150835$, the gap being of order $L_n^{-2}$. Above the threshold,
at $\beta=\tfrac85$, the quantity $nA(n)$ equals $0.868$, $0.866$, $0.864$ at the same
ranks, bounded as in (ii), and at $\beta=\tfrac65$ the quantity $n^{\beta}A(n)$ grows
through $6.89$, $8.20$, $10.67$ at $n=2\cdot10^{3}$, $5\cdot10^{3}$, $2\cdot10^{4}$, as
(iii) requires. The thresholds themselves are $z_n=0.868$, $0.931$, $0.965$, $0.986$ at
$n=10^{3}$, $10^{6}$, $10^{12}$, $10^{30}$, so no computation of practical size displays
the limit index directly, which is why the finite-rank constant is the meaningful test.
\end{numobs}
\galleryentry{M}{The square root kernel}
 {$G(n,k)=\frac12\left(1+\frac{1+\sqrt{k}}{1+\sqrt n}\right)$}
 {regular arithmetic function}
 {$G^{*}(z)=\frac{4z-1}{2(2z-1)}$, only zero at $z=\tfrac14$, the point $z=\tfrac12$ being a pole}
 {$\alpha(G)=\eta(G)=\tfrac14$}
 {proved, Theorem~\ref{thm:ex-appendix-n}}

The square-root kernel is settled in one statement, its index and its regimes together.

\begin{theorem}
\label{thm:ex-appendix-n}
Let
\[
 G_N(n,k)=\frac12\left(1+\frac{1+\sqrt{k}}{1+\sqrt n}\right).
\]
Then $A_\beta(1)=1$, and for $n\geq2$,
\begin{equation}
 A_\beta(n)=u_nA_\beta(n-1)+v_{\beta,n},
 \label{eq:ex-n-recurrence}
\end{equation}
where
\begin{align}
 u_n&=\frac{\sqrt n+2+\sqrt{n-1}}{2(1+\sqrt n)},\notag\\
 v_{\beta,n}&=n^{-\beta}
 -\frac{1+\sqrt{n-1}}{1+\sqrt n}(n-1)^{-\beta}.
 \label{eq:ex-n-v}
\end{align}
These coefficients satisfy
\begin{align}
 u_n&=1-\frac1{4n}+\mathcal O(n^{-3/2}),
 \label{eq:ex-n-u-expansion}\\
 v_{\beta,n}&=\left(\frac12-\beta\right)n^{-\beta-1}
 +\mathcal O_\beta(n^{-\beta-3/2}).
 \label{eq:ex-n-v-expansion}
\end{align}
The transform and its reciprocal are
\begin{equation}
 G_N^*(z)=\frac12+\frac z{2z-1}
 =\frac{4z-1}{2(2z-1)},
 \qquad
 \Xi_{G_N}(z)=\frac{2(2z-1)}{4z-1}.
 \label{eq:ex-n-transform}
\end{equation}
The three forcing branches are
\begin{align}
 A_\beta(n)&=\Xi_{G_N}(\beta)n^{-\beta}+o(n^{-\beta})
 &&(\beta<1/4),
 \label{eq:ex-n-below}\\
 A_{1/4}(n)&=\frac14n^{-1/4}\log n+C_{1/4}n^{-1/4}
 +\mathcal O(n^{-3/4}\log n),
 \label{eq:ex-n-critical}\\
 A_\beta(n)&=C_\beta\mathcal Q_n n^{-1/4}
 +\Xi_{G_N}(\beta)n^{-\beta}
 +\mathcal O_\beta(n^{-\beta-1/2})
 &&(\beta>1/4),
 \label{eq:ex-n-above}
\end{align}
where $\mathcal Q_n=1+\mathcal O(n^{-1/2})$ and $C_\beta>0$ for every
$\beta>1/4$.  Therefore
\begin{equation}
 \alpha(G_N)=\eta(G_N)=\frac14.
 \label{eq:ex-n-indices}
\end{equation}
\end{theorem}

\begin{proof}
Put
\[
 D_j=\sqrt{j+1}-\sqrt j,
 \qquad
 c_n=\frac1{2(1+\sqrt n)},
 \qquad
 T_n=\sum_{j\leq n}D_jA_\beta(j).
\]
Abel summation in the defining equation gives the exact Volterra relation
\begin{equation}
 A_\beta(n)=n^{-\beta}+c_nT_{n-1},
 \qquad
 T_0=0.
 \label{eq:ex-n-volterra}
\end{equation}
Eliminating $T_{n-2}$ between consecutive ranks gives
\eqref{eq:ex-n-recurrence}--\eqref{eq:ex-n-v}.  Also
\[
 u_n=1-\frac{\sqrt n-\sqrt{n-1}}{2(1+\sqrt n)}>0.
\]
The expansions
\[
 \frac{\sqrt n-\sqrt{n-1}}{2(1+\sqrt n)}
 =\frac1{4n}+\mathcal O(n^{-3/2})
\]
and
\[
 \frac{1+\sqrt{n-1}}{1+\sqrt n}
 =1-\frac1{2n}+\mathcal O(n^{-3/2})
\]
together with the binomial expansion of $(n-1)^{-\beta}$ prove
\eqref{eq:ex-n-u-expansion} and \eqref{eq:ex-n-v-expansion}.

The row profile is the exact linear combination
\[
 G_N(n,k)=\left(\frac12+\frac1{2(1+\sqrt n)}\right)
 +\frac{\sqrt n}{2(1+\sqrt n)}\left(\frac kn\right)^{1/2}.
\]
The finite probes therefore converge in $\Re z<1/2$ to the transform of
$\tfrac12(1+t^{1/2})$.  Direct integration gives
\eqref{eq:ex-n-transform}.  Its only zero is $1/4$, and its pole at
$1/2$ is not a zero.

Proposition~\ref{thm:ex-first-order} applies with
\[
 \gamma=\frac14,
 \qquad
 \delta=\frac12,
 \qquad
 b_\beta=\frac12-\beta.
\]
The coefficient identity is
\[
 \frac{b_\beta}{\gamma-\beta}
 =\frac{1/2-\beta}{1/4-\beta}
 =\Xi_{G_N}(\beta).
\]
This proves the three formulas
\eqref{eq:ex-n-below}--\eqref{eq:ex-n-above}, apart from the asserted sign
of the connection coefficient.

For that sign, \eqref{eq:ex-n-volterra} gives
\begin{equation}
 T_n=(1+c_nD_n)T_{n-1}+D_nn^{-\beta}.
 \label{eq:ex-n-positive-state}
\end{equation}
Let
\[
 R_n=\prod_{j=1}^n(1+c_jD_j).
\]
Since $c_nD_n=1/(4n)+\mathcal O(n^{-3/2})$, there is an
$r_\infty>0$ such that
\begin{equation}
 R_n=r_\infty n^{1/4}\bigl(1+\mathcal O(n^{-1/2})\bigr).
 \label{eq:ex-n-growing-product}
\end{equation}
Variation of constants in \eqref{eq:ex-n-positive-state} gives
\[
 T_n=R_n\sum_{m=1}^n\frac{D_mm^{-\beta}}{R_m}.
\]
For $\beta>1/4$, this positive series converges to a number
$S_\beta>0$, since its summand is
\[
 \frac1{2r_\infty}m^{-\beta-3/4}
 \bigl(1+\mathcal O(m^{-1/2})\bigr).
\]
Equations \eqref{eq:ex-n-volterra} and
\eqref{eq:ex-n-growing-product} yield
\[
 A_\beta(n)\sim\frac{r_\infty S_\beta}{2}n^{-1/4}.
\]
Thus $C_\beta=r_\infty S_\beta/2>0$.

For $1/4<\beta<1/2$, the forcing $v_{\beta,n}$ is eventually positive
but need not be positive at the first ranks.  More exactly,
\[
 v_{\beta,n}>0
 \quad\Longleftrightarrow\quad
 \beta<
 \frac{\log\bigl((1+\sqrt{n-1})/(1+\sqrt n)\bigr)}
 {\log((n-1)/n)},
\]
and the quotient tends to $1/2$.  Initial ranks require another argument.
Relation \eqref{eq:ex-n-volterra} gives that argument because its forcing,
its weight, and $A_\beta(1)$ are positive.  Hence
$A_\beta(n)>0$ at every rank and for every real $\beta$.  On
$(1/4,1/2)$, one has $G_N^*(\beta)<0$.  If such a $\beta$ were transparent
with a constant $\Xi'$, the argument of Theorem~\ref{thm:ex-fgv-existence}
identifies $\Xi'$ with $1/G_N^*(\beta)<0$, so the partial sums would be
eventually negative, which is impossible.  Such failures occur arbitrarily
near $1/4$.  This proves sharpness.  The scalar bounds give absorption, and the transform of this
profile has its first zero at $1/4$ as well.  Therefore
\eqref{eq:ex-n-indices} follows.  At the pole $\beta=1/2$, one has
$\Xi_{G_N}(1/2)=0$ but $A_{1/2}(n)\sim C_{1/2}n^{-1/4}$.  This point is
not transparent.
\end{proof}

\begin{remark}
The transform has its pole at $1/2$ and its zero at $1/4$. A reading that takes the pole
for the zero produces the value $1/2$ for this kernel, which is wrong. At the pole,
$\Xi_{G_N}(1/2)=0$ while $A_{1/2}(n)\sim C_{1/2}\,n^{-1/4}$, so that point is not
transparent, and the two indices agree here by \eqref{eq:ex-n-indices}.
\end{remark}

\begin{numobs}[Numerical control]\label{numobs:M_check}
The recurrence \eqref{eq:ex-n-recurrence} reproduces forward substitution in the defining
equation to $3.1\cdot10^{-15}$ over $n\le400$ at $\beta=\tfrac1{10}$, $\tfrac32$ and $3$. Below
the index the transparent constant is met, $n^{\beta}A_\beta(n)$ taking the values $1.451536$,
$1.481653$, $1.493270$, $1.495046$ at $n=10^{3}$, $10^{4}$, $10^{5}$, $2\cdot10^{5}$ for
$\beta=-\tfrac14$, against $\Xi_{G_N}(-\tfrac14)=\tfrac32$, and $1.749434$, $1.857592$,
$1.919645$, $1.932399$ at $\beta=0$ against $\Xi_{G_N}(0)=2$, the secondary branch
$n^{-1/4}$ of \eqref{eq:ex-n-above} accounting for the approach of order $n^{\beta-1/4}$.
At the index the four ranks give $n^{1/4}A_{1/4}(n)/\log n$ equal to $0.366978$, $0.338628$,
$0.321303$, $0.317321$, matching \eqref{eq:ex-n-critical} with $C_{1/4}$ near $0.82$. Above the
index $n^{1/4}A_\beta(n)$ settles on $0.910522$ at $\beta=\tfrac12$ and on $0.316748$ at
$\beta=1$. The value at $\beta=\tfrac12$ is the one that separates the pole from the zero,
since $\Xi_{G_N}(\tfrac12)=0$ while the measured limit is positive. Positivity holds at every
rank computed, the smallest value of $A_{1/2}$ over $n\le2\cdot10^{5}$ being
$4.3\cdot10^{-2}$.
\end{numobs}

\galleryentry{N}{The self gauged quadratic logarithmic kernel}
 {$G(n,k)=\frac{\log(n^{2}+k^{2})}{\log(2n^{2})}$}
 {regular arithmetic function}
 {$G^{*}(z)=1$, constant, no zero, so $\eta$ is not defined}
 {$\alpha(G)=2$, $\eta$ not defined}
 {proved, Theorem~\ref{thm:N_main}}

The same mechanism operates with the profile read on squares, and the moving threshold settles at two instead of one. The limiting transform is again constant, so the index comes from the equation.

\begin{proposition}\label{prop:N_transform}
The finite probes of $G$ converge, locally uniformly on $\Re z<0$, to $G^{*}(z)=1$, by the row
variation argument of Appendix~\ref{app:L}, so $\eta(G)$ is undefined. With
$h_2(t)=\log\frac{1+t^{2}}{2}$ and $L_n=\log(2n^{2})$ one has the identity
$G(n,k)=1+h_2(k/n)/L_n$, the transform
\begin{equation}\label{eq:N_hstar}
h_2^{*}(z)=-\log2-z\sum_{m\ge1}\frac{(-1)^{m+1}}{m(2m-z)}
=-2\int_0^1\frac{t^{1-z}}{1+t^{2}}\,dt,
\end{equation}
the integral form valid on $\Re z<2$, is meromorphic with simple poles at the even integers,
the residue at $z=2$ being $2$, and the effective transform $G^{*}_n=1+h_2^{*}/L_n$ has its
first positive zero at $z_n=2-2/L_n+\mathcal O(L_n^{-2})$.
\end{proposition}

\begin{proof}
As in Proposition~\ref{prop:L_gauge}, with $\log(1+t^{2})=\sum_{m\ge1}(-1)^{m+1}t^{2m}/m$, so
that the $m$-th term contributes $-z(-1)^{m+1}/(m(2m-z))$ and the poles sit at the even
integers, the residue at $z=2$ being $2$. The integral form is integration by parts with
$h_2'(t)=2t/(1+t^{2})$ and $h_2(1)=0$.
\end{proof}

The quadratic self-gauged kernel behaves in the same way, with the index moved to two.

\begin{theorem}\label{thm:N_main}
$\alpha(G)=2$. With $\delta_n(k)=\log\bigl(1+\tfrac{2k+1}{n^{2}+k^{2}}\bigr)$ the defining
equation is equivalent to
$A(n)=n^{-\beta}+L_n^{-1}\sum_{k<n}A(k)\delta_n(k)$, all partial sums are positive, and the
following hold.
\begin{enumerate}[label=(\roman*)]
\item For every $\beta<2$,
$n^{\beta}A(n)=1-h_2^{*}(\beta)/L_n+\mathcal O_\beta(L_n^{-2})$.
\item For every $\beta>2$ there is $B_\beta$ with $A(n)\le n^{-\beta}+B_\beta\,n^{-2}$, and
$A(n)=\mathcal O(n^{-2}\log n)$ at $\beta=2$.
\item For every $\beta>2$, $n^{\beta}A(n)\ge c\,n^{\beta-2}/L_n\to\infty$, with $c=1$
admissible, so transparency fails at every exponent above two.
\end{enumerate}
\end{theorem}

\begin{proof}
The Abel form and positivity are as in Lemma~\ref{lem:L_positive}, with
$\delta_n(k)=h_2(\tfrac{k+1}n)-h_2(\tfrac kn)$ and
$\tfrac{2k+1}{2(n^{2}+k^{2})}\le\delta_n(k)\le\tfrac{2k+1}{n^{2}+k^{2}}$, and (iii) follows
from the term $k=1$ alone, $A(n)\ge\delta_n(1)/L_n\ge1/(n^{2}L_n)$.

For (ii) fix $\beta>2$, choose $K=10$, and prove
$A(n)\le n^{-\beta}+B\,n^{-2}\mathbf 1_{\{n>K\}}+A(n)\mathbf 1_{\{n\le K\}}$ by induction on
$n>K$. The three parts of the feedback are estimated separately. The profile part is
$\sum_{k<n}k^{-\beta}(2k+1)(n^{2}+k^{2})^{-1}\le(2\zeta(\beta-1)+\zeta(\beta))n^{-2}$. The
initial part, over $k\le K$, is at most $C_K n^{-2}$ with $C_K=\sum_{k\le K}A(k)(2k+1)$. The
main part is
\[
B\sum_{K<k<n}\frac{2k+1}{k^{2}(n^{2}+k^{2})}
\le\frac{B}{n^{2}}\Bigl[2\bigl(H_{n-1}-\!\!\sum_{k<n}\tfrac{k}{n^{2}+k^{2}}\bigr)+\frac1K\Bigr]
\le\frac{B}{n^{2}}\Bigl[2\log n+2\gamma-\log2+\frac1K+\frac cn\Bigr],
\]
using the partial fraction $\tfrac{1}{k(n^{2}+k^{2})}=\tfrac1{n^{2}}\bigl(\tfrac1k-\tfrac
k{n^{2}+k^{2}}\bigr)$ and $\sum_{k<n}k/(n^{2}+k^{2})\ge\tfrac12\log2-\tfrac cn$. Since
$L_n=2\log n+\log2$, the induction closes when
$B\bigl(2\log2-2\gamma-\tfrac1K-\tfrac cn\bigr)\ge2\zeta(\beta-1)+\zeta(\beta)+C_K$, and
$2\log2-2\gamma-\tfrac1{10}=0.131\ldots>0$, so a large $B$ works for $n$ beyond an explicit
rank, the remaining ranks being absorbed into $B$. At $\beta=2$ the hypothesis
$A(k)\le ML_k/k^{2}$ closes the same way with the margin unchanged.

For (i) the bootstrap $|A(k)|\le Mk^{-\beta}$ closes for $\beta<2$ from
$\sum_{k<n}k^{-\beta}(2k+1)(n^{2}+k^{2})^{-1}\le C_\beta n^{-\beta}$, by splitting at $n/2$,
and the refined constant follows as in Theorem~\ref{thm:L_main}(i) from
\[
\sum_{k<n}k^{-\beta}\,\frac{2k+1}{n^{2}+k^{2}}
=n^{-\beta}\int_0^1\frac{2t^{1-\beta}}{1+t^{2}}\,dt+\text{lower order}
=-h_2^{*}(\beta)\,n^{-\beta}+\text{lower order},
\]
with \eqref{eq:N_hstar}, the error analysis being identical. The three parts assemble into
$\alpha(G)=2$ exactly as in Theorem~\ref{thm:L_main}.
\end{proof}

\begin{numobs}\label{numobs:N_check}
At $\beta=\tfrac25$, where $h_2^{*}(\beta)=-0.903229$ by both forms of \eqref{eq:N_hstar},
the measured constants $1.060741$, $1.055507$, $1.051115$ at $n=2\cdot10^{3}$,
$4\cdot10^{3}$, $8\cdot10^{3}$ match the finite-rank prediction of (i) to three to five
hundredths of one per cent. Above the threshold $n^{2}A(n)$ is bounded, $1.9671$, $1.9640$,
$1.9618$ at $\beta=\tfrac52$ and $n=5\cdot10^{3}$, $10^{4}$, $2\cdot10^{4}$, while at
$\beta=\tfrac{11}5$ the quantity $n^{\beta}A(n)$ grows through $12.19$, $14.56$, $19.07$ at
$n=2\cdot10^{3}$, $5\cdot10^{3}$, $2\cdot10^{4}$, as (iii) requires. Positivity holds at
every rank computed, the smallest value of $A$ over $n\le2\cdot10^{4}$ being
$4.9\cdot10^{-9}$ at $\beta=\tfrac52$. The thresholds are $z_n=1.862$, $1.929$, $1.964$,
$1.986$ at $n=10^{3}$, $10^{6}$, $10^{12}$, $10^{30}$.
\end{numobs}
\galleryentry{O}{The fractional part kernel with a diagonal jump}
 {$g(x)=1-\{x\}$ on $(0,1]$, so $G(n,k)=1-k/n$ for $k<n$ and $G(n,n)=1$}
 {function of good variation}
 {$g^{*}(z)=\frac{1}{1-z}$, single pole at $z=1$ and no zero, so $\eta$ is not defined}
 {$\alpha(g)=\tfrac14$, $\eta$ not defined}
 {proved, Theorem~\ref{thm:O_quarter}}

\rafgalleryfig{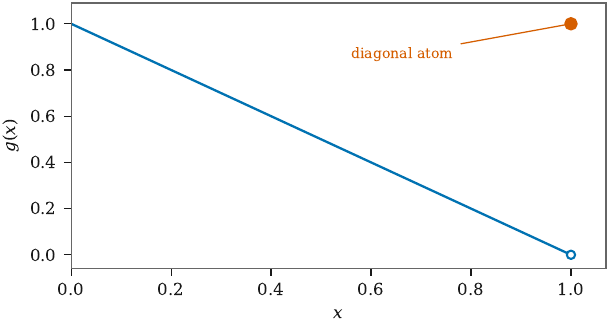}{The profile is the affine function $1-x$ off the diagonal and takes the value one at $x=1$. The open circle marks the value the affine branch would take there and the filled circle the value the kernel takes. The whole difficulty of this entry sits in that single point, which the first two terms of the additive decomposition \eqref{eq:tt_volterra} do not see.}{fig:app_O}

On $(0,1)$ this kernel is the affine function $1-x$, and at the diagonal it jumps to
$g(1)=1-\{1\}=1$. Two evaluations of the same profile must be kept apart. The defining
equation evaluates $g$ at the ratios $x=k/n\le1$, so its triangular kernel is
$G(n,k)=g(k/n)=1-\{k/n\}$, equal to $1-k/n$ for $k<n$ and to $1$ at $k=n$. The reciprocal
row evaluation $g(n/k)=1-\{n/k\}$, defined through the periodic extension of the fractional
part to $(0,\infty)$, is not the kernel of the equation, it is the observable whose sums
produce the divisor identity of Proposition~\ref{prop:O_harmonic} below. The affine
function $1-x$ alone is not admissible, since its value at $1$ vanishes and the defining
equation then loses its diagonal. The fractional part restores the diagonal, the equation
becomes well posed, and the price of that stability is a spectrum. The one profile thus
carries an identified mode on the column side and a classical open problem on the row
side.

\begin{proposition}\label{prop:O_structure}
The defining equation $\sum_{k\le n}a_k\,g(k/n)=n^{-\beta}$ is equivalent to
\begin{equation}\label{eq:O_exact}
\frac{1}{n}\sum_{j<n}A(j)+a_n=n^{-\beta},\qquad A(j)=\sum_{k\le j}a_k .
\end{equation}
The arithmetic Mellin transform of the kernel is
\begin{equation}\label{eq:O_transform}
g^{*}(z)=\frac{1}{1-z},
\end{equation}
with a single pole at $z=1$ and no zero. For $\beta<0$ the continuous equation with the
same kernel, $\int_0^x a(t)\,(1-t/x)\,dt=x^{-\beta}$, has the exact locally integrable
solution $A(x)=(1-\beta)\,x^{-\beta}$, where $A(x)=\int_0^x a(t)\,dt$. The same expression
extends formally in $\beta$, but for $\beta\ge0$ it is no longer an ordinary primitive
based at the origin, and at $\beta=1$ it vanishes identically and fails the integrated
equation outright.
\end{proposition}

\begin{proof}
For $k<n$ the kernel value is $1-k/n$ and at $k=n$ it is $1$, so the defining sum equals
$\sum_{k<n}a_k(1-k/n)+a_n$. Counting the ranks $j$ with $k\le j<n$ gives
$\sum_{j<n}A(j)=\sum_{k<n}a_k(n-k)$, and division by $n$ yields \eqref{eq:O_exact}. For the
transform, $-z\int_0^1(1-t)\,t^{-z-1}\,dt=1+\frac{z}{1-z}=\frac{1}{1-z}$ on $\Re z<0$, and
the value of $g$ at the single point $t=1$ does not affect the integral. For the continuous
equation, integration by parts turns it into $\frac1x\int_0^xA(t)\,dt=x^{-\beta}$, so
$\int_0^xA(t)\,dt=x^{1-\beta}$, and differentiation gives $A(x)=(1-\beta)\,x^{-\beta}$,
which is $x^{-\beta}/g^{*}(\beta)$. A primitive based at the origin must vanish there,
which holds exactly when $\beta<0$.
\end{proof}

The transform is zero free, so no analytic index is available and the formal continuous
response carries no finite spectral threshold, which the exact solution confirms on its
ordinary range $\beta<0$. The discrete equation behaves otherwise. The diagonal jump\index[terms]{diagonal jump} turns
a relation of the first kind into one of the second kind, and the conversion carries a
spectrum.

\begin{proposition}\label{prop:O_laguerre}
Let $(a_n)$ be the homogeneous solution, $a_1=1$ and $\frac1n\sum_{j<n}A(j)+a_n=0$ for
$n\ge2$. Then for $n\ge4$
\begin{equation}\label{eq:O_recurrence}
n\,a_n-(2n-3)\,a_{n-1}+(n-2)\,a_{n-2}=0,
\end{equation}
and
\[a_n=\frac{(-1)^{n+1}}{n!}\,\mathrm{La}_{n-2}(1;1,1)\qquad(n\ge2),
\]
where the monic first associated Laguerre polynomials with parameter $\alpha=1$ are
defined by $\mathrm{La}_{-1}(1;1,x)=0$, $\mathrm{La}_{0}(1;1,x)=1$, and
\[
\mathrm{La}_{m}(1;1,x)
=(x-2m-2)\,\mathrm{La}_{m-1}(1;1,x)
-m(m+1)\,\mathrm{La}_{m-2}(1;1,x)\qquad(m\ge1).
\]
For the corresponding non-monic normalization, see
\nm{Ismail}{M. E. H.}~\cite[p.~160, Thm.~5.6.1, Eq.~(5.6.11)]{Ismail2005}.
The coefficient triangle and the values at $x=1$ are the OEIS\index[terms]{On-Line Encyclopedia of Integer Sequences} entries A201201 and
A201202~\cite{OEIS}.
\end{proposition}

\begin{proof}
Multiplying the homogeneous relation by $n$ and subtracting the relation at rank $n-1$
gives $A(n-1)+n\,a_n-(n-1)\,a_{n-1}=0$. Subtracting the same combination at rank $n-1$ and
using $A(n-1)-A(n-2)=a_{n-1}$ yields \eqref{eq:O_recurrence}, valid as soon as the three
ranks involved are homogeneous, that is for $n\ge4$. Substituting
$a_n=(-1)^{n+1}L_{n-2}/n!$ into \eqref{eq:O_recurrence} and multiplying by
$(-1)^{n+1}(n-1)!$ turns it into
\[
L_{n-2}+(2n-3)\,L_{n-3}+(n-1)(n-2)\,L_{n-4}=0,
\]
which with $m=n-2$ is the three-term recurrence
$L_m+(2m+1)L_{m-1}+m(m+1)L_{m-2}=0$ of the first associated monic Laguerre values at
$x=1$. The initial values agree, $a_2=-\tfrac12$ against $\mathrm{La}_0=1$ and
$a_3=-\tfrac12$ against $\mathrm{La}_1(1;1,1)=-3$, so the two sequences coincide.
\end{proof}

The index concerns the partial sums rather than the coefficients, and they satisfy an
exact recurrence of their own, together with a closed generating function.

\begin{proposition}\label{prop:O_partial_sums}
Let $A_\beta(n)$ denote the partial sums of the forced solution, with $A_\beta(0)=0$ and
$A_\beta(1)=1$. For $n\ge2$,
\begin{equation}\label{eq:O_Arec}
n\,A_\beta(n)-2(n-1)\,A_\beta(n-1)+(n-1)\,A_\beta(n-2)=n^{1-\beta}-(n-1)^{1-\beta},
\end{equation}
and the homogeneous partial sums satisfy the same relation with right side zero for
$n\ge3$. For $|z|<1$,
\begin{equation}\label{eq:O_gf}
\sum_{n\ge1}A_\beta(n)\,z^{n}
=\frac{e^{-1/(1-z)}}{1-z}\int_0^{z}\frac{e^{1/(1-t)}}{t}\,
\operatorname{Li}_{\beta-1}(t)\,dt,
\end{equation}
where $\operatorname{Li}_{s}(z)=\sum_{n\ge1}z^{n}n^{-s}$.
\end{proposition}

\begin{proof}
Multiplying \eqref{eq:O_exact} by $n$ gives $n\,a_n+\sum_{j<n}A_\beta(j)=n^{1-\beta}$, and
subtracting the same relation at rank $n-1$, with $a_n=A_\beta(n)-A_\beta(n-1)$, yields
\eqref{eq:O_Arec}. In the homogeneous case the subtraction needs both ranks free of
forcing, hence $n\ge3$. For the generating function set
$F(z)=\sum_{n\ge1}A_\beta(n)z^{n}$. The recurrence \eqref{eq:O_Arec} holds for $n\ge1$
once $A_\beta(-1)=A_\beta(0)=0$ and the right side at $n=1$ is read as $1$, its telescoped
form, and summing it against $z^{n}$ gives
\[
z(1-z)^{2}F'(z)+z^{2}F(z)=(1-z)\operatorname{Li}_{\beta-1}(z),
\]
a first order linear equation whose integrating factor is $(1-z)\,e^{1/(1-z)}$. Since
$F(0)=0$, integration from the origin gives \eqref{eq:O_gf}. The right side of
\eqref{eq:O_gf} is analytic on the unit disc, so the series converges there.
\end{proof}

The factor $e^{-1/(1-z)}/(1-z)$ equals $e^{-1}\sum_{n\ge0}L_n(1)\,z^{n}$, the generating
function of the Laguerre polynomials at argument one \cite[Eq.~18.12.13]{DLMF}, in agreement with
Proposition~\ref{prop:O_laguerre}. The closed form \eqref{eq:O_gf} places on one line the
forced singularity carried by the polylogarithm, the oscillating mode carried by the
exponential factor, and the connection coefficient between the two, and it is the natural
analytic anchor for the study of the index.

\begin{proposition}\label{prop:O_harmonic}
For $n\ge1$,
\begin{equation}\label{eq:O_harmonic}
\sum_{k\le n}\bigl(1-\{n/k\}\bigr)=\gamma n+\Delta(n)-\tfrac12+\mathcal O(1/n),
\end{equation}
where $\gamma$ is the Euler constant\index[terms]{Euler constant} and
$\Delta(n)=\sum_{m\le n}\tau(m)-n\log n-(2\gamma-1)\,n$ is the remainder of the Dirichlet
divisor problem.
\end{proposition}

\begin{proof}
Since $\{n/k\}=n/k-\lfloor n/k\rfloor$ and
$\sum_{k\le n}\lfloor n/k\rfloor=\sum_{m\le n}\tau(m)$, the left side is
$n-nH_n+\sum_{m\le n}\tau(m)$. By the definition of $\Delta(n)$ it remains only to note
\[
H_n=\log n+\gamma+\frac1{2n}+\mathcal O(n^{-2}).
\]
This estimate follows directly from the defining limit for $\gamma$: expanding
$\log(1+1/k)-1/(k+1)=1/(2k^2)+\mathcal O(k^{-3})$ and summing from $k=n$ to infinity gives
$H_n-\log n-\gamma=1/(2n)+\mathcal O(n^{-2})$. Substitution proves
\eqref{eq:O_harmonic}.
\end{proof}

The index concerns the partial sums, and their recurrence \eqref{eq:O_Arec} is the ordinary
Laguerre recurrence at argument one. That observation is what closes the quarter. The
coefficients $a_n$ obey the associated Laguerre relation of
Proposition~\ref{prop:O_laguerre}, whereas the partial sums obey the ordinary one, and the
ordinary equation has two explicit solutions with known asymptotics.

\begin{lemma}\label{lem:O_basis}
Let $(u_n)$ and $(v_n)$ be defined by $u_0=1$, $u_1=0$, $v_0=0$, $v_1=1$ and, for $n\ge2$,
\begin{equation}\label{eq:O_homog}
n\,x_n-2(n-1)\,x_{n-1}+(n-1)\,x_{n-2}=0 .
\end{equation}
Then $u_n=L_n(1)$, the ordinary Laguerre polynomial evaluated at one, and with
$y=(1-z)^{-1}$,
\begin{equation}\label{eq:O_gf_uv}
\sum_{n\ge0}u_nz^{n}=\frac{e^{-z/(1-z)}}{1-z},
\qquad
\sum_{n\ge0}v_nz^{n}=y\,e^{-y}\bigl(\operatorname{Ei}(y)-\operatorname{Ei}(1)\bigr).
\end{equation}
Their Casoratian\index[terms]{Casoratian} is
\begin{equation}\label{eq:O_casoratian}
u_nv_{n-1}-u_{n-1}v_n=-\frac1n\qquad(n\ge1).
\end{equation}
With $\theta_n=2\sqrt n-\pi/4$,
\begin{align}
u_n&=\frac{e^{1/2}}{\sqrt\pi}\,n^{-1/4}\Bigl(\cos\theta_n+\mathcal O(n^{-1/2})\Bigr),
\label{eq:O_fejer_u}\\
v_n&=\frac{e^{-1/2}}{\sqrt\pi}\,n^{-1/4}
\Bigl(-\operatorname{Ei}(1)\cos\theta_n+\pi\sin\theta_n+\mathcal O(n^{-1/2})\Bigr).
\label{eq:O_fejer_v}
\end{align}
Both have complete expansions in powers of $n^{-1/2}$ whose terms are linear combinations of
$\cos\theta_n$ and $\sin\theta_n$. In particular $u_n,v_n=\mathcal O(n^{-1/4})$.
\end{lemma}
\begin{proof}
The three term recurrence of the ordinary Laguerre polynomials,
$(n+1)L_{n+1}(x)=(2n+1-x)L_n(x)-nL_{n-1}(x)$, becomes \eqref{eq:O_homog} at $x=1$ after the
shift $n\mapsto n-1$, and $L_0(1)=1$, $L_1(1)=0$ match the initial data, so $u_n=L_n(1)$. The
generating function of the Laguerre polynomials at $x=1$ gives the first identity of
\eqref{eq:O_gf_uv}. Summing \eqref{eq:O_homog} against $z^{n}$ for the initial data of $v$
gives $z(1-z)^{2}V'(z)+z^{2}V(z)=z$, whose solution with $V(0)=0$ is
$V(z)=\frac{e^{-1/(1-z)}}{1-z}\int_0^{z}\frac{e^{1/(1-t)}}{1-t}\,dt$, and the substitution
$s=(1-t)^{-1}$ gives the second identity. The Casoratian\index[terms]{Casoratian} of \eqref{eq:O_homog} satisfies
$W_n=\frac{n-1}{n}W_{n-1}$ with $W_1=-1$, which is \eqref{eq:O_casoratian}.

Formula \eqref{eq:O_fejer_u} with its complete expansion is the fixed argument expansion of
Fej\'er\index[names]{Fej\'er, L.} for $L_n(1)$
\cite[Th.~8.22.1]{Szego1975}. To determine the leading constants for $v_n$, use the Cauchy
coefficient integral for the second identity of \eqref{eq:O_gf_uv}. Its exponential phase in
the variable $y=(1-z)^{-1}$ is
\[
 -y-(n+1)\log(1-y^{-1}),
\]
whose two saddles are
$y_{\pm}=\tfrac12\pm i\sqrt{n+\tfrac34}$. Thus the local steepest descent arcs remain in fixed
subsectors about $\arg y=\pm\pi/2$. On the principal branches,
\[
 \operatorname{Ei}(y)=-E_1(-y)+i\pi\operatorname{sgn}(\Im y)
 \qquad(\Im y\ne0),
\]
and the expansion of $E_1$ is uniform in these subsectors
\cite[\S\S6.4, 6.12(i)]{DLMF}. Consequently, for every fixed $M$,
\[
 ye^{-y}\bigl(\operatorname{Ei}(y)-\operatorname{Ei}(1)\bigr)
 =\sum_{j=0}^{M-1}j!\,y^{-j}
  +\bigl(-\operatorname{Ei}(1)+i\pi\sigma\bigr)ye^{-y}
  +\mathcal O_M(|y|^{-M})
\]
on the arc where $\sigma=\operatorname{sgn}(\Im y)$. The finite sum is a polynomial in $1-z$
and has no coefficient of rank $n\ge M$, while $ye^{-y}=e^{-1}U(z)$. Combining the two
conjugate saddle contributions in
Fej\'er's calculation therefore gives \eqref{eq:O_fejer_v}, including its stated error term.

It remains to justify the complete expansions without any further saddle-point assumption.
After shifting the index, \eqref{eq:O_homog} is
\[
 x_{n+2}+p(n)x_{n+1}+q(n)x_n=0,
 \qquad
 p(n)=-\frac{2(n+1)}{n+2},\qquad q(n)=\frac{n+1}{n+2}.
\]
Writing $p(n)\sim\sum p_jn^{-j}$ and $q(n)\sim\sum q_jn^{-j}$ gives
\[
 (p_0,p_1,q_0,q_1)=(-2,2,1,-1).
\]
The characteristic equation is $(\rho-1)^2=0$, and the generic coincident-root condition is
satisfied since $2q_1=-2\ne-4=p_0p_1$. The Birkhoff--Adams theorem in the form of Wong\index[names]{Wong, R.} and
Li\index[names]{Li, H.} \cite{WongLi1992}, also recorded in \cite[\S2.9(ii)]{DLMF}, applies. Its two parameters
are
\[
 \kappa^2=\frac{2p_0p_1-4q_1}{q_0}=-4,
 \qquad
 \alpha=\frac{q_0+2q_1}{4q_0}=-\frac14.
\]
Thus the double root splits at the next order into two conjugate oscillatory modes, both with
envelope $n^{-1/4}$. More precisely, there are exact independent solutions with
\[
 w_\pm(n)\sim n^{-1/4}e^{\pm2i\sqrt n}
 \sum_{j\ge0}c_{\pm,j}n^{-j/2}.
\]
The coefficients of the recurrence are real, so the basis may be chosen conjugate. Since
\eqref{eq:O_casoratian} is nonzero, $u$ and $v$ are constant linear combinations of this basis.
Their expansions are therefore complete in powers of $n^{-1/2}$, with every coefficient a
linear combination of $\cos\theta_n$ and $\sin\theta_n$, as asserted.
\end{proof}

The threshold is the convergence abscissa of an oscillatory series, and the next lemma
isolates it.

\begin{lemma}\label{lem:O_blocks}
Let $\vartheta_m=2\sqrt m+\vartheta_0$ and, for real $p$, put
$S_p(N)=\sum_{m\le N}m^{-p}e^{i\vartheta_m}$. If $p>\tfrac12$ the series converges and its
tail is $\mathcal O_p(N^{1/2-p})$. If $p=\tfrac12$ the partial sums are bounded. The
convergence in the first case is locally uniform in a complex parameter $p$ on
$\{\Re p>\tfrac12\}$. Both assertions persist after multiplication of the summand by a
function admitting an asymptotic expansion in integral powers of $m^{-1/2}$.
\end{lemma}

\begin{proof}
Group the integers into the blocks $k^{2}\le m<(k+1)^{2}$. Writing $m=k^{2}+r$ with
$0\le r\le 2k$, expanding $\sqrt{k^{2}+r}=k+\frac{r}{2k}+\mathcal O(k^{-3}r^{2})$ and
applying Euler-Maclaurin in $r$ gives, to any fixed order,
\[\sum_{k^{2}\le m<(k+1)^{2}}m^{-p}e^{2i\sqrt m}
=e^{2ik}k^{1-2p}\Bigl(c_0+\frac{c_1(p)}{k}+\mathcal O_p(k^{-2})\Bigr),
\qquad c_0=\int_0^{2}e^{it}\,dt\neq0 .
\]
For $p>\tfrac12$ the amplitude $k^{1-2p}$ decreases to zero, the partial sums of $e^{2ik}$
are bounded because $e^{2i}\neq1$, and Dirichlet's test gives convergence, the tail from
$k\asymp\sqrt N$ being $\mathcal O(k^{1-2p})=\mathcal O(N^{1/2-p})$. At $p=\tfrac12$ the
leading block contribution is $c_0e^{2ik}$, whose partial sums are bounded, and the
subsequent terms converge by the same test. An incomplete final block obeys the same bound
directly, so the conclusion holds for every $N$ and not only at block endpoints. The
estimates are uniform on compact parameter sets, and an extra factor with an expansion in
powers of $m^{-1/2}$ produces finitely many series of the same shape with shifted exponents,
all covered by the same argument.
\end{proof}

Variation of constants against the two fundamental solutions gives the solution in closed form.

\begin{lemma}\label{lem:O_connection}
For every complex $\beta$ and every $n\ge1$,
\begin{equation}\label{eq:O_connection}
A_\beta(n)=v_n\Bigl(1+\sum_{m=2}^{n}d_\beta(m)\,u_{m-1}\Bigr)
-u_n\sum_{m=2}^{n}d_\beta(m)\,v_{m-1},
\qquad
d_\beta(m)=m^{1-\beta}-(m-1)^{1-\beta}.
\end{equation}
\end{lemma}

\begin{proof}
Normalize \eqref{eq:O_Arec} as $x_m=\frac{2(m-1)}{m}x_{m-1}-\frac{m-1}{m}x_{m-2}+r_m$ with
$r_m=d_\beta(m)/m$. By \eqref{eq:O_casoratian} an impulse $r_m$ at rank $m$ propagates to
rank $n\ge m$ as $r_m\,m\,(u_{m-1}v_n-v_{m-1}u_n)$. The initial data $A_\beta(0)=0$ and
$A_\beta(1)=1$ contribute $v_n$, and summing the impulses gives \eqref{eq:O_connection}.
\end{proof}

Since $d_\beta(m)=(1-\beta)m^{-\beta}+\mathcal O_\beta(m^{-\beta-1})$, the products in
\eqref{eq:O_connection} have leading exponent $p=\beta+\tfrac14$ in Lemma~\ref{lem:O_blocks}.
The two limits
\begin{equation}\label{eq:O_UV}
U(\beta)=1+\sum_{m\ge2}d_\beta(m)u_{m-1},
\qquad
V(\beta)=\sum_{m\ge2}d_\beta(m)v_{m-1}
\end{equation}
therefore exist for $\Re\beta>\tfrac14$, with locally uniform convergence, so $U$ and $V$ are
holomorphic on that half plane. The quarter is the abscissa at which these series cease to
converge, and nothing else.

\begin{theorem}[The quarter theorem]\label{thm:O_quarter}
The kernel $g(x)=1-\{x\}$ is a function of good variation with
\[
\alpha(g)=\tfrac14 .
\]
For every real $\beta<\tfrac14$,
\begin{equation}\label{eq:O_transparency}
A_\beta(n)=(1-\beta)\,n^{-\beta}+o(n^{-\beta})=\frac{n^{-\beta}}{g^{*}(\beta)}+o(n^{-\beta}),
\end{equation}
and for every real $\beta\ge\tfrac14$,
\begin{equation}\label{eq:O_absorption}
A_\beta(n)=\mathcal O_\beta\!\left(n^{-1/4}\right),
\end{equation}
with no logarithmic loss at the critical exponent. Every interval
$(\tfrac14,\tfrac14+\delta)$ contains an exponent at which \eqref{eq:O_transparency} fails,
and the absorbed exponent is attained, at $\beta=1$ the sequence $n^{1/4}A_1(n)$ having a
nonzero oscillatory limit superior.
\end{theorem}

\begin{proof}
Write $D$ for the operator on the left of \eqref{eq:O_Arec}. For fixed $\lambda$ and every
integer $M\ge1$, binomial expansion gives
$D(n^{-\lambda})=n^{-\lambda}\bigl(1+\sum_{j<M}q_j(\lambda)n^{-j}+\mathcal
O_{\lambda,M}(n^{-M})\bigr)$, the leading coefficient being one, and likewise
$d_\beta(n)=n^{-\beta}\bigl((1-\beta)+\sum_{j<M}b_j(\beta)n^{-j}+\mathcal
O_{\beta,M}(n^{-M})\bigr)$. Triangular matching of coefficients therefore builds, for every
$M$, a sequence $P_{\beta,M}(n)=n^{-\beta}\sum_{j<M}c_j(\beta)n^{-j}$ with $c_0(\beta)=1-\beta$
and
\begin{equation}\label{eq:O_profile}
D\,P_{\beta,M}(n)-d_\beta(n)=\mathcal O_{\beta,M}(n^{-\beta-M}).
\end{equation}
Choose $M$ with $\beta+M>\tfrac34$ and set $R(n)=A_\beta(n)-P_{\beta,M}(n)$, giving
$P_{\beta,M}(0)$ any fixed value. Applying Lemma~\ref{lem:O_connection} to the recurrence
satisfied by $R$, which differs from \eqref{eq:O_connection} only by the homogeneous term
$R(0)u_n+R(1)v_n$, and using \eqref{eq:O_profile} together with
$u_m,v_m=\mathcal O(m^{-1/4})$, both connection series converge absolutely, since
$\beta+M+\tfrac14>1$. Hence $R(n)=\mathcal O_{\beta,M}(n^{-1/4})$. For $\beta<\tfrac14$ this
error is $o(n^{-\beta})$, as is every term of the profile after the first, which proves
\eqref{eq:O_transparency}, the constant $1-\beta$ being $1/g^{*}(\beta)$ by
\eqref{eq:O_transform}.

For $\beta>\tfrac14$ the two series \eqref{eq:O_UV} converge, so the parentheses in
\eqref{eq:O_connection} are bounded. At $\beta=\tfrac14$ the leading summands of both are,
by the complete expansions of Lemma~\ref{lem:O_basis}, linear combinations of
$m^{-1/2}e^{\pm2i\sqrt m}$, whose partial sums are bounded by the critical case of
Lemma~\ref{lem:O_blocks}, and every subsequent term converges. The parentheses therefore stay
bounded at the endpoint as well, and $u_n,v_n=\mathcal O(n^{-1/4})$ gives
\eqref{eq:O_absorption} without any logarithm.

For sharpness, the tail estimate of Lemma~\ref{lem:O_blocks} applied to
\eqref{eq:O_connection} gives, for $\beta>\tfrac14$,
\begin{equation}\label{eq:O_sharp}
A_\beta(n)=U(\beta)\,v_n-V(\beta)\,u_n+\mathcal O_\beta(n^{-\beta}).
\end{equation}
At $\beta=1$ one has $d_1(m)=0$ for every $m\ge2$, hence $U(1)=1$, $V(1)=0$ and
$A_1(n)=v_n$. So $U$ is not identically zero on the connected half plane $\Re\beta>\tfrac14$,
its restriction to the real axis is real analytic, and no interval to the right of $\tfrac14$
consists entirely of its zeros. For each $\delta>0$ there is thus
$\beta_\delta\in(\tfrac14,\tfrac14+\delta)$ with $U(\beta_\delta)\neq0$, and by
\eqref{eq:O_fejer_u}, \eqref{eq:O_fejer_v} and \eqref{eq:O_sharp} the coefficient of
$n^{-1/4}\sin\theta_n$ in $A_{\beta_\delta}(n)$ equals
$\sqrt\pi\,e^{-1/2}U(\beta_\delta)\neq0$, the second solution alone carrying that component.
Hence $A_{\beta_\delta}(n)$ is not asymptotic to $(1-\beta_\delta)n^{-\beta_\delta}$, which is
$o(n^{-1/4})$, and transparency fails arbitrarily close to the right of the quarter. With the
two preceding paragraphs this gives $\alpha(g)=\tfrac14$, and
$A_1(n)=v_n$ with \eqref{eq:O_fejer_v} shows that the exponent is attained.
\end{proof}

\begin{remark}\label{rem:O_exact_zero}
At $\beta=0$ the forcing is $d_0(n)=1$ and $D(1)=1$, so the constant is a particular
solution and the initial data give the exact identity
\[A_0(n)=1-L_n(1)\qquad(n\ge0).
\]
The transparent main term and the secondary mode of size $n^{-1/4}$ are already visible at
this forcing without any asymptotic construction.
\end{remark}

\begin{numobs}\label{numobs:O_mode}
The homogeneous solutions have the envelopes predicted by \eqref{eq:O_fejer_u} and
\eqref{eq:O_fejer_v}, the octave maxima of $n^{1/4}|u_n|$ and $n^{1/4}|v_n|$ on
$[10^{6},2\cdot10^{6}]$ being $0.9301913$ and $1.2555028$ against $e^{1/2}/\sqrt\pi=0.930191$
and $\sqrt{\operatorname{Ei}(1)^{2}+\pi^{2}}\,e^{-1/2}/\sqrt\pi=1.255503$. The identity
$u_n=L_n(1)$ and the Casoratian\index[terms]{Casoratian} \eqref{eq:O_casoratian} hold with residual zero in exact
rational arithmetic up to $n=59$, the generating function of $v$ agrees with its closed form
to $10^{-18}$, and the connection formula\index[terms]{connection formula} \eqref{eq:O_connection} has residual $10^{-14}$
over six exponents, exactly zero at $\beta=1$. At $\beta=\tfrac14$ the octave maxima of
$n^{1/4}|A(n)|$ equal $1.807305$, $1.807312$, $1.807316$, $1.807317$, $1.807318$,
$1.807318$, $1.807318$, $1.807318$, $1.807319$ over the nine octaves from $10^{4}$ to
$5.12\cdot10^{6}$, showing no logarithmic growth. Below the quarter the constant $1-\beta$
emerges with a secondary term of the size $n^{-1/4}$ predicted by
\eqref{eq:O_transparency}, and the identity of Proposition~\ref{prop:O_harmonic} holds with
measured residual $-0.500000$ at $n=3\cdot10^{5}$.
\end{numobs}

\begin{remark}\label{rem:O_twins}
The same kernel read along its rows produces the divisor remainder and read along its truncated
columns produces the Laguerre mode. The two windows carry the exponent $1/4$, a theorem on the
column side and a conjectural optimum on the row side. The classical lower bound is attributed to
Hardy \cite{Hardy1916}, and Huxley\index[names]{Huxley, M. N.} \cite{Huxley2003} gives the best published upper exponent quoted
here, while the later Li\index[names]{Li, X.}--Yang\index[names]{Yang, X.} preprint record is stated in Remark~\ref{rem:huxley_exponent}.
These documentary qualifiers are not used in Theorem~\ref{thm:O_quarter}. The two quarters are
companions, not equivalents. The harmonic quarter superposes the contributions of all divisors
and is the divisor problem itself. No implication between them is claimed in either direction,
and \eqref{eq:O_harmonic} offers no route to the divisor problem, since rewriting a sum does not
improve its remainder. An operator bridge between the row observable and the triangular inverse,
if one is found, would start from a theorem rather than from a conjecture.
\end{remark}

\subsection*{Reciprocity and the direct index}

Remark~\ref{rem:O_twins} names two quarters carried by one profile and stops there. What
follows separates the two axes that produce them. The first axis is the direction of the
problem for a fixed kernel, the data being the coefficients in the direct problem and the
output of the operator in the inverse problem. The second axis is the reciprocity
$x\mapsto1/x$ of a profile defined on all of $(0,\infty)$, which produces two kernels on
$(0,1]$ out of one profile. The regularity index of the theory lives on the inverse side of
one of the two kernels, and the divisor remainder lives on the direct side of the other.
They are not two faces of one operator.

Several forcing exponents are compared in what follows, so the subscript is kept on $A_\beta$ and on the sequences attached to it. Elsewhere in the volume the exponent is fixed by its context and the plain $a_n$ and $A(n)$ of Definition~\ref{def:reg_index_fgv} are used.

\begin{definition}[Reciprocal pair]\label{def:O_reciprocal}
Let $H$ be a profile on $(0,\infty)$. Its kernel and its reciprocal kernel are the two
functions on $(0,1]$
\[
g_-^{H}(x)=H(x),
\qquad
g_+^{H}(x)=H(1/x),
\]
which agree at $x=1$. Conversely a pair of kernels agreeing at $1$ determines $H$. For a
kernel $g$ on $(0,1]$ and a sequence $a$ write
\[
\mathcal T_g\,a(n)=\sum_{k\le n}a_k\,g(k/n)
\]
for the direct sum, so that the defining equation of the theory is
$\mathcal T_{g}\,a(n)=n^{-\beta}$ read as an equation in $a$.
\end{definition}

A profile can be read in two ways, and the two readings differ by a computable amount.

\begin{proposition}\label{prop:O_identification}
For every profile $H$ and every $N\ge1$,
\begin{equation}\label{eq:O_identification}
\sum_{k\le N}H(N/k)=\sum_{k\le N}g_+^{H}(k/N)=\mathcal T_{g_+^{H}}\mathbf 1(N),
\end{equation}
while the triangular kernel of the defining equation attached to $H$ is
$G(n,k)=g_-^{H}(k/n)$. The sum on the left is therefore the direct problem of the reciprocal
kernel at the constant weight, and the index $\alpha$ of Theorem~\ref{thm:O_quarter} is the
inverse problem of the kernel.
\end{proposition}

\begin{proof}
For $k\le N$ one has $N/k\ge1$, so $H(N/k)=H\bigl(1/(k/N)\bigr)=g_+^{H}(k/N)$ with
$k/N\in(0,1]$. For $k\le n$ one has $k/n\in(0,1]$, so $G(n,k)=H(k/n)=g_-^{H}(k/n)$.
\end{proof}

The two axes give four problems, and the present appendix occupies two of them.

\medskip
\renewcommand{\arraystretch}{1.45}
{\small
\begin{tabular}{lp{5.2cm}p{5.2cm}}
\toprule
& {Direct problem, $a$ given} & {Inverse problem, $\mathcal T_ga$ given} \\
\midrule
Kernel $g_-^{H}$ & Elementary, Proposition~\ref{prop:O_direct}. No divisor remainder.
& $\alpha(g_-^{H})=\tfrac14$, Theorem~\ref{thm:O_quarter}. The Laguerre mode. \\
Reciprocal kernel $g_+^{H}$ & $\delta_{\mathbf 1}(g_+^{H})$, the divisor and circle
remainders, Corollary~\ref{cor:O_indices}. & Not studied here. \\
\bottomrule
\end{tabular}
\renewcommand{\arraystretch}{1}
}
\medskip

\noindent The direct problem of the kernel itself carries no arithmetic, which locates the
divisor remainder precisely.

\begin{proposition}\label{prop:O_direct}
For $g_-(x)=1-\{x\}$ on $(0,1]$, for every complex $\beta$ and every $n\ge1$,
\begin{equation}\label{eq:O_direct}
\mathcal T_{g_-}\bigl(k^{-\beta}\bigr)(n)
=\sum_{k\le n}k^{-\beta}\Bigl(1-\Bigl\{\frac kn\Bigr\}\Bigr)
=H_n^{(\beta)}-\frac{1}{n}\,H_n^{(\beta-1)}+n^{-\beta},
\qquad
H_n^{(\beta)}=\sum_{k\le n}k^{-\beta}.
\end{equation}
In particular $\mathcal T_{g_-}\mathbf 1(n)=\tfrac{n+1}{2}$ exactly, and for $\beta=\tfrac12$
and $\beta=1$ the right side is $\tfrac43\sqrt n+\zeta(\tfrac12)+\mathcal O(n^{-1/2})$ and
$\log n+\gamma-1+\mathcal O(n^{-1})$. No remainder of divisor type occurs at any exponent.
\end{proposition}

\begin{proof}
For $k<n$ the ratio $k/n$ lies in $(0,1)$ and $g_-(k/n)=1-k/n$, while $g_-(1)=1$. Splitting
off the rank $k=n$,
\[
\sum_{k<n}k^{-\beta}\Bigl(1-\frac kn\Bigr)+n^{-\beta}
=\bigl(H_n^{(\beta)}-n^{-\beta}\bigr)-\frac1n\bigl(H_n^{(\beta-1)}-n^{1-\beta}\bigr)+n^{-\beta},
\]
which is \eqref{eq:O_direct}. At $\beta=0$ this is $n-\tfrac{n+1}{2}+1=\tfrac{n+1}{2}$. The
two asymptotics follow from
$H_n^{(1/2)}=2\sqrt n+\zeta(\tfrac12)+\mathcal O(n^{-1/2})$,
$H_n^{(-1/2)}=\tfrac23n^{3/2}+\tfrac12\sqrt n+\zeta(-\tfrac12)+\mathcal O(n^{-1/2})$ and
$H_n^{(1)}=\log n+\gamma+\mathcal O(n^{-1})$. Every term is a partial sum of a power series
in the rank, and the divisor function\index[terms]{divisor function} never enters.
\end{proof}

The Riemann-sum discrepancy of this direct problem, its smooth second-order criterion and its
arithmetic specialization to the Ingham kernel are developed in
Section~\ref{sec:ab_riemann_discrepancy}. They belong to the general Abelian theory rather than
to this particular example.

The two readings separate, and the separation depends on only part of the data.

\begin{lemma}\label{lem:O_separation}
Let $H$ be a profile. For every $N\ge1$ the value $\mathcal T_{g_+^{H}}\mathbf 1(N)$ depends
only on the restriction of $H$ to $[1,\infty)$, and for every $n\ge1$ the row
$\bigl(G(n,k)\bigr)_{k\le n}$ depends only on the restriction of $H$ to $(0,1]$. The two
restrictions meet at the single point $x=1$. The two indices of the table are therefore
computed from disjoint data, apart from the common value $H(1)$.
\end{lemma}

\begin{proof}
If $k\le N$ then $N/k\ge1$, and if $k\le n$ then $k/n\le1$.
\end{proof}

\begin{definition}[Arithmetic completion]\label{def:O_completion}
Let $g$ be a kernel on $(0,1]$. A profile $H$ on $(0,\infty)$ is an arithmetic completion of
$g$ when $g_-^{H}=g$ and, for $x\ge1$,
\[
H(x)=w(\lfloor x\rfloor)+c_0+c_1x,
\qquad
w(q)=\sum_{j\le q}a_j,
\qquad
a_1+c_0+c_1=g(1),
\]
with $c_0$ and $c_1$ real. The sequence $a$ is the arithmetic profile of $H$ and
$A(s)=\sum_{j\ge1}a_jj^{-s}$ its Dirichlet series.
\end{definition}

The problem reduces to a Dirichlet convolution.

\begin{proposition}\label{prop:O_reduction}
For every arithmetic completion $H$ and every $N\ge1$,
\[\mathcal T_{g_+^{H}}\mathbf 1(N)=\sum_{n\le N}(a*1)(n)+c_0N+c_1NH_N,
\qquad
(a*1)(n)=\sum_{j\mid n}a_j,
\]
so the arithmetic part of the direct sum has Dirichlet series $\zeta(s)A(s)$, and the affine
part contributes to the main term alone, up to $\mathcal O(1/N)$.
\end{proposition}

\begin{proof}
By \eqref{eq:O_identification} the sum is $\sum_{k\le N}H(N/k)$, and $N/k\ge1$, so the
completion is read through $w(\lfloor N/k\rfloor)+c_0+c_1N/k$. Then
\[
\sum_{k\le N}w(\lfloor N/k\rfloor)=\sum_{k\le N}\ \sum_{j\le N/k}a_j=\sum_{jk\le N}a_j
=\sum_{n\le N}\ \sum_{j\mid n}a_j ,
\]
and the remaining terms give $c_0N+c_1N\sum_{k\le N}k^{-1}$, which by
$H_N=\log N+\gamma+\tfrac1{2N}+\mathcal O(N^{-2})$ equals
$c_0N+c_1(N\log N+\gamma N+\tfrac12)+\mathcal O(1/N)$.
\end{proof}

\begin{definition}[Direct index at weight one]\label{def:O_delta}
Let $H$ be an arithmetic completion whose series $A$ extends meromorphically to
$\{\Re s>0\}$ with finitely many poles and at most polynomial growth on vertical lines. Put
\[
M_H(N)=\sum_{\Re\varrho\ge0}\ \operatorname*{Res}_{s=\varrho}\ \frac{\zeta(s)A(s)N^{s}}{s}
\ +\ c_0N+c_1\Bigl(N\log N+\gamma N+\tfrac12\Bigr),
\]
the first sum running over the poles of $\zeta(s)A(s)/s$ in $\{\Re s\ge0\}$. The direct
index at weight one of the reciprocal kernel is
\begin{equation}\label{eq:O_delta_def}
\delta_{\mathbf 1}\bigl(g_+^{H}\bigr)=\inf\bigl\{\theta\ge0:\
\mathcal T_{g_+^{H}}\mathbf 1(N)-M_H(N)=\mathcal O_\eps\bigl(N^{\theta+\eps}\bigr)\bigr\}.
\end{equation}
The subscript records the constant weight $\mathbf 1$.
\end{definition}

The kernel of this appendix carries at least three arithmetic completions built from a
monotone or periodic staircase, and the three are the two classical lattice point problems
together with the dyadic branch of Chapter~\ref{chap:abelian}.

\begin{proposition}\label{prop:O_three}
Let $g_0=g_-$ be the kernel $1-\{x\}$ on $(0,1]$, that is $1-x$ on $(0,1)$ and $1$ at $x=1$.
The three arithmetic profiles
\[
a_j=1,
\qquad
a_j=(-1)^{j-1},
\qquad
a_j=\chi_4(j)
\]
define arithmetic completions $H_D$, $H_E$ and $H_C$ of $g_0$, with staircases
\[
w_D(q)=q,
\qquad
w_E(q)=\mathbf 1_{\{2\nmid q\}},
\qquad
w_C(q)=\mathbf 1_{\{q\equiv1,2\ (\mathrm{mod}\ 4)\}},
\]
the first taken with $c_0=1$ and $c_1=-1$, so that $H_D(x)=1-\{x\}$ on all of $(0,\infty)$,
and the other two with $c_0=c_1=0$. Their main terms are
\[
M_{H_D}(N)=\gamma N-\tfrac14,
\qquad
M_{H_E}(N)=N\log2-\tfrac14,
\qquad
M_{H_C}(N)=\tfrac{\pi}{4}N-\tfrac14,
\]
and their direct sums at weight one satisfy
\begin{align}
\mathcal T_{g_+^{H_D}}\mathbf 1(N)&=\gamma N+\Delta(N)-\tfrac12+\mathcal O(1/N),\label{eq:O_BD}\\
\mathcal T_{g_+^{H_E}}\mathbf 1(N)&=\sum_{m\le N}\tau(m)-2\sum_{m\le\lfloor N/2\rfloor}\tau(m)\notag\\
&=N\log2+\Delta(N)-2\Delta(\lfloor N/2\rfloor)
+\varrho_N\bigl(\log N+2\gamma-\log2\bigr)+\mathcal O(1/N),\label{eq:O_BE}\\
\mathcal T_{g_+^{H_C}}\mathbf 1(N)&=\frac{N_C(N)-1}{4}
=\frac{\pi}{4}N+\frac{\Delta_C(N)}{4}-\frac14,\label{eq:O_BC}
\end{align}
where $\varrho_N=N-2\lfloor N/2\rfloor\in\{0,1\}$, where $N_C(N)$ counts the lattice points
of the closed disc of radius $\sqrt N$ centered at the origin, and where
$\Delta_C(N)=N_C(N)-\pi N$. The first identity of \eqref{eq:O_BE} and the whole of
\eqref{eq:O_BC} are exact.
\end{proposition}

\begin{proof}
In the three cases $a_1=1$ and the constants match $g_0(1)=1$, so the three are completions
of the same kernel. The Dirichlet series are $A_D=\zeta$, $A_E=\eta$ and $A_C=\beta$, where
$\eta(s)=(1-2^{1-s})\zeta(s)$, so by Proposition~\ref{prop:O_reduction} the arithmetic parts
have series $\zeta^{2}$, $\zeta\eta=(1-2^{1-s})\zeta^{2}$ and $\zeta\beta$. The main terms
follow from the residues. For $\zeta^{2}$ the double pole at $s=1$ gives
$N\log N+(2\gamma-1)N$ and the simple pole of the integrand at the origin gives
$\zeta(0)^{2}=\tfrac14$, and the affine block subtracts $N\log N+\gamma N+\tfrac12$ while
adding $N$, leaving $\gamma N-\tfrac14$. For $\zeta\eta$ the factor $1-2^{1-s}$ vanishes
simply at $s=1$, so the pole there is simple with residue $\log2$, and the origin
contributes $(1-2)\zeta(0)^{2}=-\tfrac14$. For $\zeta\beta$ the pole at $s=1$ is simple with
residue $\beta(1)=\pi/4$, and the origin contributes $\zeta(0)\beta(0)=-\tfrac14$.

Identity \eqref{eq:O_BD} is Proposition~\ref{prop:O_harmonic}, through
\eqref{eq:O_identification}. For \eqref{eq:O_BE}, splitting the alternating sum by the
parity of $k$ and writing $k=2i$ in the even part,
\[
\sum_{k\le N}(-1)^{k-1}\Bigl\lfloor\frac Nk\Bigr\rfloor
=\sum_{k\le N}\Bigl\lfloor\frac Nk\Bigr\rfloor
-2\sum_{i\le\lfloor N/2\rfloor}\Bigl\lfloor\frac{\lfloor N/2\rfloor}{i}\Bigr\rfloor,
\]
by the nested floor identity $\lfloor N/(2i)\rfloor=\lfloor\lfloor N/2\rfloor/i\rfloor$, and
each of the two sums is a divisor sum. Applying Dirichlet's formula to both, with
$2\lfloor N/2\rfloor=N-\varrho_N$ and
$\lfloor N/2\rfloor\log\lfloor N/2\rfloor
=\tfrac12(N-\varrho_N)(\log N-\log2)-\tfrac{\varrho_N}{2}+\mathcal O(1/N)$, gives the second
line. For \eqref{eq:O_BC}, the two square theorem of Jacobi gives
$r_2(n)=4\sum_{d\mid n}\chi_4(d)$ for $n\ge1$, so
$\sum_{n\le N}(\chi_4*1)(n)=\tfrac14\sum_{n\le N}r_2(n)=\tfrac14(N_C(N)-1)$, the origin
being the lattice point excluded, and the last equality is the definition of $\Delta_C$.
\end{proof}

The three indices follow from the optimal exponents of the divisor problem and of the circle
problem.

\begin{corollary}\label{cor:O_indices}
Let $\theta$ and $\theta_C$ be the optimal exponents of the divisor and circle problems. Then
\[
\delta_{\mathbf 1}\bigl(g_+^{H_D}\bigr)=\theta,
\qquad
\delta_{\mathbf 1}\bigl(g_+^{H_C}\bigr)=\theta_C,
\qquad
\delta_{\mathbf 1}\bigl(g_+^{H_E}\bigr)\le\theta,
\]
while $\alpha(g_-^{H})=\tfrac14$ for all three by Theorem~\ref{thm:O_quarter}.
\end{corollary}

\begin{remark}
The classical lower bounds for $\theta$ and $\theta_C$ require two distinct primary sources.
They are not used in this corollary and are not asserted here pending exact citation reconciliation.
\end{remark}

\begin{proof}
Subtract the main terms of Proposition~\ref{prop:O_three} from \eqref{eq:O_BD},
\eqref{eq:O_BE} and \eqref{eq:O_BC}. The remainders are $\Delta(N)+\mathcal O(1)$,
$\Delta(N)-2\Delta(\lfloor N/2\rfloor)+\mathcal O(\log N)$ and $\Delta_C(N)/4$, and a term
$\mathcal O(\log N)$ is $\mathcal O_\eps(N^{\eps})$, so it does not move the infimum
\eqref{eq:O_delta_def}. The middle bound is the triangle inequality applied to the dyadic
difference. The three kernels coincide with $g_0$ by Proposition~\ref{prop:O_three}.
\end{proof}

No reverse bound is available in the middle case, since an omega theorem for $\Delta$ carries
no formal consequence for the dyadic difference $\Delta(N)-2\Delta(\lfloor N/2\rfloor)$, and
whether $\delta_{\mathbf 1}(g_+^{H_E})$ equals $\theta$ is open.

\begin{proposition}\label{prop:O_undetermined}
The profile $a_j=\mu(j)$ defines an arithmetic completion $H_\mu$ of the same kernel $g_0$,
with $\mathcal T_{g_+^{H_\mu}}\mathbf 1(N)=1$ and $M_{H_\mu}(N)=1$ for every $N\ge1$, hence
$\delta_{\mathbf 1}(g_+^{H_\mu})=0$. A kernel therefore leaves the direct index of its
reciprocal free.
\end{proposition}

\begin{proof}
Here $a_1=\mu(1)=1=g_0(1)$ and $c_0=c_1=0$, so $H_\mu$ is a completion of $g_0$. By
Proposition~\ref{prop:O_reduction} and $\mu*1=\mathbf 1_{\{n=1\}}$ the direct sum equals
$1$ for every $N$. The series is $\zeta(s)A(s)=\zeta(s)/\zeta(s)=1$, entire, so the only
pole of $\zeta(s)A(s)/s$ in $\{\Re s\ge0\}$ is the simple pole at the origin, of residue
$1$, and $M_{H_\mu}(N)=1$. The difference vanishes identically.
\end{proof}

\begin{definition}[Reciprocal defect]\label{def:O_defect}
For a profile $H$ that is an arithmetic completion of a kernel carrying a regularity index,
put
\[\mathfrak d(H)=\delta_{\mathbf 1}\bigl(g_+^{H}\bigr)-\alpha\bigl(g_-^{H}\bigr).
\]
The profile is balanced when $\mathfrak d(H)=0$. The defect compares the direct index of one
kernel with the inverse index of the other, the two being reciprocal, and it is a property
of the profile and not of either kernel alone.
\end{definition}

\begin{remark}[Two classical problems as balance statements]\label{rem:O_balance}
By Corollary~\ref{cor:O_indices} the divisor problem is the assertion $\mathfrak d(H_D)=0$
and the circle problem is the assertion $\mathfrak d(H_C)=0$. These are translations and not
reductions. The two direct indices are the two classical exponents by exact identities, the
inverse index is $\tfrac14$ by Theorem~\ref{thm:O_quarter}, and rewriting an exponent does
not estimate it. What the translation supplies is a frame in which the two classical
exponents and the Laguerre quarter are indices of one reciprocal pair, and a question that
neither classical problem asks, since Proposition~\ref{prop:O_undetermined} shows that one
kernel of the pair leaves the other free.
\end{remark}

The reciprocal kernel is not the only route to the divisor remainder. The same remainder is
the direct problem of the Ingham function itself, at the harmonic weight, and that reading
stays inside $(0,1]$ throughout.

\begin{proposition}\label{prop:O_ingham}
With $\Phi(x)=x\lfloor1/x\rfloor$ the Ingham function of Chapter~\ref{chap:ingham} and
$a_k=1/k$, for every $n\ge1$
\begin{equation}\label{eq:O_ingham}
\mathcal T_\Phi\,a(n)=\sum_{k\le n}\frac1k\,\Phi\Bigl(\frac kn\Bigr)
=\frac1n\sum_{k\le n}\Bigl\lfloor\frac nk\Bigr\rfloor
=\frac{1}{n}\sum_{m\le n}\tau(m),
\end{equation}
so that $n\,\mathcal T_\Phi a(n)=n\log n+(2\gamma-1)n+\Delta(n)$. Moreover the two readings
are tied by the pointwise identity
\begin{equation}\label{eq:O_reciprocity}
x\,H_D(1/x)=x+\Phi(x)-1,
\qquad 0<x\le1 .
\end{equation}
\end{proposition}

\begin{proof}
Since $\Phi(k/n)=(k/n)\lfloor n/k\rfloor$, the weight cancels the rank,
\[
\frac1k\,\Phi\Bigl(\frac kn\Bigr)
=\frac1k\cdot\frac kn\Bigl\lfloor\frac nk\Bigr\rfloor
=\frac1n\Bigl\lfloor\frac nk\Bigr\rfloor ,
\]
and $\sum_{k\le n}\lfloor n/k\rfloor=\sum_{m\le n}\tau(m)$. Dirichlet's formula gives the
second display. For \eqref{eq:O_reciprocity}, $H_D(1/x)=1-\{1/x\}=1-1/x+\lfloor1/x\rfloor$, and
multiplying by $x$ gives $x-1+x\lfloor1/x\rfloor$.
\end{proof}

The divisor remainder is therefore attached to the Ingham function on the direct side, at the
weight $1/k$, and the Ingham function is the one whose inverse index carries the master
equivalence of Section~\ref{sec:master_equiv}. One kernel, two weights, two exponents, and
no implication between them is claimed.

\begin{remark}[Reciprocity of the transforms]\label{rem:O_transform}
Identity \eqref{eq:O_reciprocity} has a transform shadow. With the arithmetic Mellin
transform of Chapter~\ref{chap:ingham} and $\Phi^{*}(z)=\tfrac{z}{z-1}\zeta(1-z)$ from
Proposition~\ref{prop:ingham_mellin}, the reciprocal kernel $g_+^{H_D}(x)=1-\{1/x\}$ has
\[
\bigl(g_+^{H_D}\bigr)^{*}(z)
=1+\frac{z}{z+1}\bigl(\Phi^{*}(z+1)-1\bigr)
=1+\zeta(-z)-\frac{z}{z+1},
\]
against $\bigl(g_-^{H_D}\bigr)^{*}(z)=1/(1-z)$ from \eqref{eq:O_transform}. The kernel of the
pair has a zero free transform with a single pole, the reciprocal has a transform built on
$\zeta$, and the analytic index of the reciprocal kernel is not determined here.
\end{remark}

\begin{remark}[The alternating branch and the Fibonacci sums]\label{rem:O_fibonacci}
The completion $H_E$ is the one that appears in Chapter~\ref{chap:abelian}. Its fluctuation
is the dyadic difference $\Delta(N)-2\Delta(\lfloor N/2\rfloor)$ by \eqref{eq:O_BE}, which
is the object produced by the fractional part sums of Fibonacci numbers along the even
indices, the odd indices producing the analogous dyadic difference of the circle remainder.
The three completions carry the three degree two objects $\zeta^{2}$, $\zeta\eta$ and
$\zeta\beta$, and the quarter is their common critical scale, so the parity split of
Chapter~\ref{chap:abelian} and the two lattice point problems are three completions of one
kernel rather than three unrelated coincidences.
\end{remark}

The family of Appendix~\ref{app:P} supplies a one parameter test of the balance, since each
of its members is an arithmetic completion in the sense of
Definition~\ref{def:O_completion} whose kernel has a known index.

\begin{proposition}\label{prop:O_family}
For an integer $m\ge2$ let $H_m(x)=1-\{x\}+\{x/m\}$ off the integers with $H_m(1)=1$. Then
$H_m$ is an arithmetic completion of the kernel $g_m$ of Appendix~\ref{app:P}, with
arithmetic profile $a_j=1-\mathbf 1_{\{m\mid j\}}$, affine part $c_0=1$ and $c_1=-c_m$, and
Dirichlet series $A(s)=(1-m^{-s})\zeta(s)$. Its direct sum at weight one is
\begin{equation}\label{eq:O_Rm}
\mathcal T_{g_+^{H_m}}\mathbf 1(N)
=N-c_mNH_N+\sum_{r\le N}\tau(r)-\sum_{r\le\lfloor N/m\rfloor}\tau(r)-\frac1m ,
\end{equation}
exactly at every rank, and
\begin{equation}\label{eq:O_Rm_main}
\mathcal T_{g_+^{H_m}}\mathbf 1(N)
=C_mN+\Delta(N)-\Delta\bigl(\lfloor N/m\rfloor\bigr)+\mathcal O_m(\log N),
\qquad
C_m=\gamma+\frac{1-\gamma+\log m}{m}.
\end{equation}
\end{proposition}

\begin{proof}
Off the integers $H_m(x)=\lfloor x\rfloor-\lfloor x/m\rfloor+1-c_mx$, and
$\lfloor x/m\rfloor=\lfloor\lfloor x\rfloor/m\rfloor$, so the staircase is
$w(q)=q-\lfloor q/m\rfloor$, whose increments are $a_j=1-\mathbf 1_{\{m\mid j\}}$ and whose
series is $\zeta(s)-m^{-s}\zeta(s)$. At $x=1$ the formula would give $1+1/m$ and the value is
$1$, which is the term $-1/m$ of \eqref{eq:O_Rm}, and on $(0,1]$ the profile is $1-c_mx$ off
the diagonal, so $g_-^{H_m}=g_m$. Proposition~\ref{prop:O_reduction} and
$\sum_{k\le N}a_k\lfloor N/k\rfloor=\sum_{k\le N}\lfloor N/k\rfloor
-\sum_{i\le\lfloor N/m\rfloor}\lfloor\lfloor N/m\rfloor/i\rfloor$ give \eqref{eq:O_Rm}.
Applying Dirichlet's formula to both divisor sums and
$H_N=\log N+\gamma+\mathcal O(1/N)$, the terms in $N\log N$ cancel and the coefficient of $N$
is $\gamma c_m-c_m+(\log m)/m+1$, which is $C_m$. The residue computation of
Definition~\ref{def:O_delta} returns the same constant, the double pole of
$(1-m^{-s})\zeta(s)^{2}$ at $s=1$ contributing $c_m\log N-c_m+2\gamma c_m+(\log m)/m$ and the
affine block the rest. The remaining terms are bounded by a multiple of $\log N$.
\end{proof}

Telescoping along an $m$-adic scale is what turns a pointwise bound into a global one.

\begin{lemma}\label{lem:O_telescope}
Fix an integer $m\ge2$ and a real $\theta\ge0$. Then
$\Delta(N)-\Delta\bigl(\lfloor N/m\rfloor\bigr)=\mathcal O(N^{\theta+\eps})$ for every
$\eps>0$ if and only if $\Delta(N)=\mathcal O(N^{\theta+\eps})$ for every $\eps>0$. Hence
$\delta_{\mathbf 1}\bigl(g_+^{H_m}\bigr)=\theta$, the optimal exponent of the divisor
problem, for every $m\ge2$.
\end{lemma}

\begin{proof}
One direction is immediate. For the other, $\lfloor\lfloor N/m^{j}\rfloor/m\rfloor
=\lfloor N/m^{j+1}\rfloor$ and $\Delta(0)=0$, so
$\Delta(N)=\sum_{j\ge0}\bigl[\Delta(\lfloor N/m^{j}\rfloor)
-\Delta(\lfloor N/m^{j+1}\rfloor)\bigr]$, a finite sum, and bounding each bracket by
$C(N/m^{j})^{\theta+\eps}$ gives $\Delta(N)\ll N^{\theta+\eps}$, the geometric series
converging. The index follows from \eqref{eq:O_Rm_main}, a term $\mathcal O(\log N)$ being
$\mathcal O_\eps(N^{\eps})$.
\end{proof}

The weight one in the difference is what makes the telescoping work. The dyadic difference
$\Delta(N)-2\Delta(\lfloor N/2\rfloor)$ of \eqref{eq:O_BE} carries the weight two, the same
telescoping there produces the factor $2^{J}$ against $\Delta(\lfloor N/2^{J}\rfloor)$, and
no bound survives the limit. That is why Corollary~\ref{cor:O_indices} states only one
inequality for the alternating branch while Lemma~\ref{lem:O_telescope} states an equality
here.

\begin{corollary}\label{cor:O_unbalanced}
For every integer $m\ge2$,
\[
\mathfrak d(H_m)=\theta-m ,
\]
which is negative and bounded away from zero, so no member of the family is balanced. The
kernels converge uniformly to the kernel of the present appendix, the direct index of the
reciprocal stays equal to $\theta$ throughout, and the inverse index of the kernel runs to
infinity. Balance is reached only at the limit, where the assertion $\mathfrak d=0$ is the
divisor problem.
\end{corollary}

\begin{proof}
Lemma~\ref{lem:O_telescope} gives $\delta_{\mathbf 1}(g_+^{H_m})=\theta$ and
Theorem~\ref{thm:P_index} gives $\alpha(g_-^{H_m})=\alpha(g_m)=m$, so
Definition~\ref{def:O_defect} gives the value. Since $\theta\le\tfrac12$ by the elementary
bound and $m\ge2$, the defect is at most $-\tfrac32$. The uniform convergence and the two
limiting indices are Theorem~\ref{thm:P_discontinuity}, and the limit profile is
$H_\infty(x)=1-\{x\}$, whose two indices are $\tfrac14$ and $\theta$ by
Theorem~\ref{thm:O_quarter} and Corollary~\ref{cor:O_indices}.
\end{proof}

\begin{remark}[What the family settles]\label{rem:O_family_lesson}
The two indices of a profile are independent to the fullest extent, and not only in the
sense of Proposition~\ref{prop:O_undetermined}, where the completion was free. Here the
whole profile moves inside a single family, the Abelian side is the divisor problem at every
member, and the tauberian side takes every integer value from two upward. The coincidence of
the two quarters at the limit profile is therefore a property of the limit alone. On the
tauberian side it comes from the coalescence of the two characteristic roots described in
Remark~\ref{rem:P_coalescence}, which replaces a geometric mode and an algebraic mode by the
Laguerre basis of envelope $n^{-1/4}$. On the Abelian side the quarter is the conjectural
exponent of the divisor problem, reached through Voronoi summation. Nothing links the two
mechanisms, and the equality of the two values at that one profile has no counterpart
anywhere else in the family.
\end{remark}

\begin{openproblem}\label{op:O_balance}
Characterize the profiles that are balanced in the sense of
Definition~\ref{def:O_defect}. Corollary~\ref{cor:O_unbalanced} settles a one parameter
family in the negative and locates the only balanced member at its limit, and
Proposition~\ref{prop:O_undetermined} shows that a kernel leaves the direct index of its
reciprocal free, so balance is exceptional rather than generic. What is missing is a
criterion, and no criterion is available that is not a restatement of the exponent it would
decide. The empty box of the table above, the inverse problem of a reciprocal kernel, is
untouched.
\end{openproblem}

\begin{numobs}\label{numobs:O_abelian}
Identity \eqref{eq:O_direct} is exact at the six exponents $\beta=-\tfrac12$, $0$,
$\tfrac12$, $1$, $\tfrac32$ and $2$ on every tested rank, the two asymptotics measuring
$1.328725$ against $\tfrac43$ and $-0.422769$ against $\gamma-1=-0.422784$ at $n=10^{5}$. Identity \eqref{eq:O_ingham} is exact, and so is the general form
\[
\sum_{k\le n}k^{-\beta}H_D(n/k)
=\sum_{m\le n}\sigma_{-\beta}(m)+H_n^{(\beta)}-nH_n^{(\beta+1)}
\]
for $\beta\in\{0,\tfrac12,1,2\}$. Identity \eqref{eq:O_BC} is exact up to $N=2\cdot10^{5}$
and the first identity of \eqref{eq:O_BE} up to $N=8\cdot10^{6}$. The residual of
\eqref{eq:O_BD} against $\gamma N+\Delta(N)-\tfrac12$ measures $0.0833/N$ over four decades,
its value being $1/(12N)$. Over $[10^{3},4\cdot10^{6}]$ the three fluctuations divided by
$N^{1/4}$ stay in $[0.8,3.7]$, $[-3.0,-0.3]$ and $[-0.4,0.6]$, with no visible growth, and
over $[10^{3},8\cdot10^{6}]$ the octave maxima of
$N^{-1/4}|\Delta(N)-2\Delta(\lfloor N/2\rfloor)|$ and of $N^{-1/4}|\Delta(N)|$ track each
other, $4.87$ against $5.80$ on the last octave. For the profile $a_j=\tau(j)$, whose
convolution is the Piltz function $\tau_3$, the same ratio reaches $89.7$ at
$N=4\cdot10^{6}$ and grows, the fluctuation following $N^{1/3}$ rather than $N^{1/4}$. The
transform of Remark~\ref{rem:O_transform} agrees with numerical quadrature to
$2\cdot10^{-15}$ at $z=-4+2i$.
\end{numobs}

\galleryentry{P}{The affine family with a diagonal jump}
 {$g_m(x)=1-\bigl(1-\tfrac1m\bigr)x$ on $(0,1)$ and $g_m(1)=1$, with $m\ge2$ an integer}
 {function of good variation}
 {$g_m^{*}(z)=\frac{m-z}{m(1-z)}$, single zero at $z=m$, simple pole at $z=1$}
 {$\alpha(g_m)=\eta(g_m)=m$}
 {proved, Theorem~\ref{thm:P_index}}

\rafgalleryfig{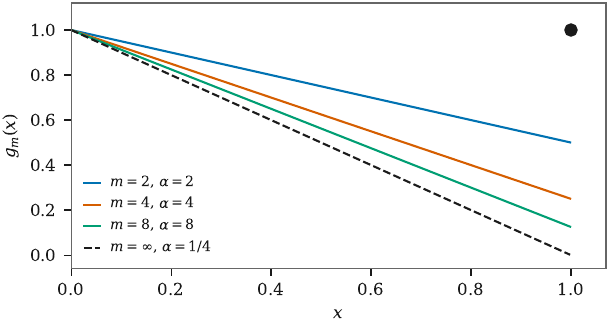}{The family $g_m$ together with its uniform limit. The profiles differ from the limit by at most $1/m$, so they converge uniformly on $(0,1]$, while their indices are $m$ and diverge. The limit has index one quarter. The diagonal jump, common to every member and to the limit, is the filled point at $x=1$.}{fig:app_P}

For an integer $m\ge2$ write $c_m=1-\tfrac1m$ and let $g_m:(0,1]\to\R$ be
\[
g_m(x)=\begin{cases}1-c_mx, & 0<x<1,\\ 1, & x=1,\end{cases}
\]
equivalently $g_m(x)=1-\{x\}+\{x/m\}$ off the diagonal, with the value at the diagonal
raised to one. The jump of height $c_m$ is what keeps the defining equation of the second
kind, since the affine branch alone takes the value $1-c_m$ at $x=1$ and would weight the
diagonal below the interior. The family is affine of decreasing slope, it converges to the
kernel of Appendix~\ref{app:O} as $m$ grows, and every statement below is uniform in $m$
except where the value $m=2$ is singled out for explicitness. The arithmetic Mellin transform
is
\begin{equation}\label{eq:P_transform}
g_m^{*}(z)=-z\int_0^1(1-c_mt)\,t^{-z-1}\,dt=\frac{m-z}{m(1-z)},
\end{equation}
with a single zero at $z=m$ and a simple pole at $z=1$, so $\eta(g_m)=m$. The pole falls
strictly inside the range of exponents below the index, and that is the feature of this
entry. The transparency constant $1/g_m^{*}$ vanishes at $\beta=1$ for every $m$, while the
equation itself stays regular there.

\begin{proposition}\label{prop:P_exact}
Put $A(r)=\sum_{k\le r}a_k$ and $B(r)=\sum_{k\le r}ka_k$, with $A(0)=B(0)=0$. The defining
equation $\sum_{k\le n}a_kg_m(k/n)=n^{-\beta}$ is equivalent to
\begin{equation}\label{eq:P_exact}
A(n)-c_m\,\frac{B(n-1)}{n}=n^{-\beta}\qquad(n\ge1),
\end{equation}
which already forces $a_1=1$, and differencing between consecutive ranks gives
\begin{equation}\label{eq:P_recurrence}
n\,A(n)-(1+c_m)(n-1)\,A(n-1)+c_m(n-1)\,A(n-2)=d_\beta(n),
\qquad d_\beta(n)=n^{1-\beta}-(n-1)^{1-\beta},
\end{equation}
valid for $n\ge2$. The indicial equation is $1-\gamma/m=0$, the numerator of
\eqref{eq:P_transform} up to a constant, so the algebraic mode is $n^{-m}$, real and simple,
and the characteristic polynomial of the recurrence at infinity is
\begin{equation}\label{eq:P_characteristic}
(X-1)(X-c_m).
\end{equation}
\end{proposition}

\begin{proof}
For $k<n$ the kernel value is $1-c_mk/n$ and at $k=n$ it is $1$, so the defining sum equals
$A(n-1)-c_mB(n-1)/n+a_n$, and $a_n=A(n)-A(n-1)$ gives \eqref{eq:P_exact}, which at $n=1$
reads $A(1)=1$. Multiplying by $n$, subtracting the same relation at rank $n-1$ and using
$B(n-1)-B(n-2)=(n-1)\bigl(A(n-1)-A(n-2)\bigr)$ yields \eqref{eq:P_recurrence}. With
$L[x](n)$ the left side and $x_n=n^{-\gamma}$, the expansions
$(n-1)^{1-\gamma}=n^{1-\gamma}\bigl(1-(1-\gamma)/n+\mathcal O(n^{-2})\bigr)$ and
$(n-1)(n-2)^{-\gamma}=n^{1-\gamma}\bigl(1+(2\gamma-1)/n+\mathcal O(n^{-2})\bigr)$ give
$L[n^{-\gamma}]=\bigl(1-\gamma(1-c_m)\bigr)n^{-\gamma}+\mathcal O(n^{-\gamma-1})$, and
$1-c_m=1/m$. Dividing \eqref{eq:P_recurrence} by $n$ and letting $n$ grow gives
\eqref{eq:P_characteristic}.
\end{proof}

The two characteristic roots are $1$ and $c_m$, separated for every finite $m$. The root
$c_m$ carries a geometric solution with a rational generating function, and the root $1$ is
the degenerate one that carries the algebraic mode.

\begin{lemma}\label{lem:P_basis}
Let $(u_n)$ and $(v_n)$ solve the homogeneous form of \eqref{eq:P_recurrence} with $u_0=0$,
$u_1=1$ and with $v$ the solution whose generating function is
\begin{equation}\label{eq:P_v_gf}
\sum_{n\ge0}v_nz^{n}=(1-z)^{m-1}(1-c_mz)^{-m},
\end{equation}
a rational function, so that $v_n=\mathcal O_m\bigl(n^{m-1}c_m^{\,n}\bigr)$. Then the
Casoratian is
\begin{equation}\label{eq:P_casoratian}
u_nv_{n-1}-u_{n-1}v_n=W_1\,\frac{c_m^{\,n-1}}{n}\qquad(n\ge1),
\end{equation}
the generating function of $u$ is
\begin{equation}\label{eq:P_gf}
\sum_{n\ge0}u_nz^{n}=(1-z)^{m-1}(1-c_mz)^{-m}\int_0^{z}(1-t)^{-m}(1-c_mt)^{m-1}\,dt,
\end{equation}
whose only non rational ingredient is the term $-c_m^{\,m-1}\log(1-z)$ of the integral, and
\begin{equation}\label{eq:P_u_asymp}
u_n=\frac{\kappa_m}{n^{m}}+\mathcal O_m\bigl(n^{-m-1}\bigr),
\qquad
\kappa_m=(-1)^{m-1}(m-1)!\;m\,(m-1)^{m-1}\neq0 .
\end{equation}
At $m=2$ everything is explicit, $v_n=(n-1)2^{-n}$ exactly, $W_1=-1$, and
\[
\sum_{n\ge0}u_nz^{n}=\frac{2z}{(2-z)^{2}}-\frac{2(1-z)\log(1-z)}{(2-z)^{2}},
\qquad
u_n=-\frac{2}{n^{2}}+\mathcal O(n^{-3}).
\]
\end{lemma}

\begin{proof}
Summing the homogeneous form of \eqref{eq:P_recurrence} against $z^{n}$ over $n\ge0$, with
$x_{-1}=x_{-2}=0$, gives $(1-z)(1-c_mz)F'(z)+c_mzF(z)=x_1$, the right side being the defect
of the relation at $n=1$. For $x_1=0$ the equation separates, and
$c_mz/\bigl((1-z)(1-c_mz)\bigr)=\tfrac{c_m}{1-c_m}\bigl(\tfrac1{1-z}-\tfrac1{1-c_mz}\bigr)$
with $c_m/(1-c_m)=m-1$ and $1/(1-c_m)=m$ integrates to \eqref{eq:P_v_gf}. Variation of
constants with $x_1=1$ gives \eqref{eq:P_gf}. For a recurrence $p_nx_n+q_nx_{n-1}+r_nx_{n-2}=0$
the Casoratian satisfies $W_n=(r_n/p_n)W_{n-1}$, here $W_n=c_m(n-1)W_{n-1}/n$, and
$\prod_{k=2}^{n}(k-1)/k=1/n$ gives \eqref{eq:P_casoratian}.

For \eqref{eq:P_u_asymp}, expand $1-c_mt=(1-c_m)+c_m(1-t)$ by the binomial theorem, so that
$(1-c_mt)^{m-1}(1-t)^{-m}=\sum_{j\le m-1}\binom{m-1}{j}(1-c_m)^{m-1-j}c_m^{\,j}(1-t)^{j-m}$,
whose only term of order $(1-t)^{-1}$ is the one with $j=m-1$, of coefficient $c_m^{\,m-1}$.
The integral in \eqref{eq:P_gf} is therefore a rational function of $z$ minus
$c_m^{\,m-1}\log(1-z)$. The rational part contributes coefficients of geometric decay, and
the singular part at $z=1$ is $c_m^{\,m-1}(1-c_m)^{-m}(1-z)^{m-1}\log\frac1{1-z}$ with
$(1-c_m)^{-m}=m^{m}$. Finally
\[
[z^{n}]\,(1-z)^{m-1}\log\frac1{1-z}
=\sum_{j\le m-1}\binom{m-1}{j}\frac{(-1)^{j}}{n-j}
=\frac{(-1)^{m-1}(m-1)!}{n(n-1)\cdots(n-m+1)},
\]
by the partial fraction identity for $\sum_j\binom{p}{j}(-1)^{j}/(x-j)$, which is
$(-1)^{m-1}(m-1)!\,n^{-m}\bigl(1+\mathcal O_m(1/n)\bigr)$. Collecting,
$\kappa_m=c_m^{\,m-1}m^{m}(-1)^{m-1}(m-1)!=(-1)^{m-1}(m-1)!\,m\,(m-1)^{m-1}$, nonzero for
every $m\ge2$. At $m=2$ the substitution $x_n=(n-1)2^{-n}$ into the homogeneous form gives
$2^{-n}(n-1)\bigl[2n-6(n-2)+4(n-3)\bigr]=0$, and \eqref{eq:P_v_gf} reduces to
$(1-z)\bigl(1-\tfrac z2\bigr)^{-2}$, whose coefficients are $(n-1)2^{-n}$.
\end{proof}

Variation of constants against the basis of the homogeneous equation gives the solution in
closed form.

\begin{lemma}\label{lem:P_connection}
For every $\beta$ and every $n\ge1$ the solution of \eqref{eq:P_recurrence} with $A(0)=0$ and
$A(1)=1$ satisfies
\begin{equation}\label{eq:P_connection}
A(n)=u_n\Bigl(1-\sum_{2\le r\le n}\frac{v_{r-1}}{W_1c_m^{\,r-1}}\,r\,d_\beta(r)\Bigr)
+v_n\sum_{2\le r\le n}\frac{u_{r-1}}{W_1c_m^{\,r-1}}\,r\,d_\beta(r).
\end{equation}
At $m=2$ the two weights are $-(r-2)$ and $-2^{r-1}u_{r-1}$ respectively.
\end{lemma}

\begin{proof}
Discrete variation of constants for $p_nx_n+q_nx_{n-1}+r_nx_{n-2}=f_n$ gives the particular
solution $\sum_{r\le n}(f_r/p_r)\bigl(u_nv_{r-1}-u_{r-1}v_n\bigr)/W_r$. Here $p_r=r$,
$f_r=d_\beta(r)$ and $W_r=W_1c_m^{\,r-1}/r$ by \eqref{eq:P_casoratian}. Adding the
homogeneous part fixed by $A(0)=0$ and $A(1)=1$, which is $u_n$, gives
\eqref{eq:P_connection}. At $m=2$ one has $W_1=-1$, $c_2^{\,r-1}=2^{1-r}$ and
$v_{r-1}=(r-2)2^{1-r}$, and $r\,d_\beta(r)/r$ has already been used.
\end{proof}

The index of the family follows, for every integer scale.

\begin{theorem}\label{thm:P_index}
For every integer $m\ge2$ the kernel $g_m$ is a function of good variation with
\[
\alpha(g_m)=\eta(g_m)=m .
\]
For every real $\beta<m$,
\begin{equation}\label{eq:P_transparency}
A_\beta(n)=\frac{n^{-\beta}}{g_m^{*}(\beta)}
+\mathcal O_{m,\beta}\bigl(n^{-\beta-1}+n^{-m}\bigr),
\qquad
\frac{1}{g_m^{*}(\beta)}=\frac{m(1-\beta)}{m-\beta},
\end{equation}
at the critical exponent $A_m(n)=\mathcal O_m(n^{-m}\log n)$, and for every $\beta>m$ one has
$A_\beta(n)=\mathcal O_{m,\beta}(n^{-m})$ with the exponent attained. At $m=2$ the constant
above the index is explicit,
\begin{equation}\label{eq:P_zeta}
A_\beta(n)=-2\,\zeta(\beta-1)\,n^{-2}+o\bigl(n^{-2}\bigr)
\qquad(\beta>2),
\end{equation}
and $\zeta(\beta-1)>1$, so the exponent is attained at every such exponent.
\end{theorem}

\begin{proof}
Start from \eqref{eq:P_connection}. Since
$d_\beta(r)=(1-\beta)r^{-\beta}+\mathcal O_\beta(r^{-\beta-1})$ and
$v_{r-1}/(W_1c_m^{\,r-1})$ is a polynomial in $r$ of degree $m-1$ with leading coefficient
fixed by \eqref{eq:P_v_gf}, the summand of the first series is
$\mathcal O_{m}(r^{m-1-\beta})$, so Euler and Maclaurin give a first series of size
$n^{m-\beta}$ for $\beta<m$, of size $\log n$ at $\beta=m$, and convergent for $\beta>m$.
Multiplied by $u_n\asymp\kappa_mn^{-m}$ this produces the three regimes of the statement, the
constant below the index matching $1/g_m^{*}(\beta)$ through \eqref{eq:P_transform}. The
second series has summand $\mathcal O_{m,\beta}(c_m^{-r}r^{-\beta-1})$ by
\eqref{eq:P_u_asymp}, hence size $\mathcal O(c_m^{-n}n^{-\beta-1})$, and
$v_n=\mathcal O_m(n^{m-1}c_m^{\,n})$ makes the whole second block
$\mathcal O_{m,\beta}(n^{m-\beta-2})$, below the first block throughout. Transparency holds
below $m$ and absorption at and above $m$, so $\alpha(g_m)=m=\eta(g_m)$.

At $m=2$ the first weight is $-(r-2)$ and the first series is $\sum_{2\le r\le n}(r-2)d_\beta(r)$.
For $\beta>2$ the identity
\[
\sum_{r\ge2}(r-2)\bigl(r^{1-\beta}-(r-1)^{1-\beta}\bigr)
=\sum_{r\ge2}\bigl[(r-2)-(r-1)\bigr]r^{1-\beta}
=1-\zeta(\beta-1),
\]
obtained by shifting the index in the second half and noting that the term $r=1$ vanishes,
turns the first block into $-u_n\zeta(\beta-1)$, which with \eqref{eq:P_u_asymp} is
\eqref{eq:P_zeta}.
\end{proof}

The family converges uniformly to the kernel of the preceding appendix while its index does not,
which is the discontinuity announced above.

\begin{theorem}\label{thm:P_discontinuity}
Let $g_\infty(x)=1-\{x\}$ be the kernel of Appendix~\ref{app:O}. Then
\[
\sup_{0<x\le1}\bigl|g_m(x)-g_\infty(x)\bigr|=\frac1m\longrightarrow0,
\qquad
\alpha(g_m)=m\longrightarrow+\infty,
\qquad
\alpha(g_\infty)=\frac14 .
\]
The regularity index is therefore neither continuous nor upper semicontinuous under uniform
convergence of kernels on $(0,1]$. The family exhibited here does not decide lower
semicontinuity, since the indices $\alpha(g_m)=m$ increase without bound while the limit index is
finite.
\end{theorem}

\begin{proof}
The two kernels differ by $x/m$ off the diagonal and agree at $x=1$, which gives the
uniform bound. The index of $g_m$ is Theorem~\ref{thm:P_index} and the index of $g_\infty$ is
Theorem~\ref{thm:O_quarter}.
\end{proof}

\begin{remark}[Where the jump comes from]\label{rem:P_coalescence}
The mechanism is visible in \eqref{eq:P_characteristic}. For every finite $m$ the two
characteristic roots $1$ and $c_m$ are separated, the root $c_m$ carrying a geometric
solution and the root $1$ carrying the algebraic mode $n^{-m}$ fixed by the indicial
equation. As $m$ grows the roots merge, $(X-1)(X-c_m)\to(X-1)^{2}$, and at the limit the
recurrence \eqref{eq:P_recurrence} becomes
$nA(n)-2(n-1)A(n-1)+(n-1)A(n-2)=d_\beta(n)$, which is \eqref{eq:O_Arec} verbatim. The
degenerate root replaces the pair of a geometric and an algebraic solution by the Laguerre
basis of Lemma~\ref{lem:O_basis}, of envelope $n^{-1/4}$, and the index collapses from
$m$ to $\tfrac14$. The family therefore joins the present appendix at $m=2$ to
Appendix~\ref{app:O} at the limit, and the quarter of Appendix~\ref{app:O} is a property of
the coalescence and not of the affine profile.
\end{remark}

\begin{remark}[The degenerate exponent]\label{rem:P_degenerate}
The pole of $g_m^{*}$ at $z=1$ makes $1/g_m^{*}(1)=0$ for every $m$, so
\eqref{eq:P_transparency} at $\beta=1$ reads $A_1(n)=\mathcal O_m(n^{-m})$ and the profile it
predicts is absent. This is the convention of Definition~\ref{def:reg_index_fgv}, met before
at $\beta=0$ in Appendix~\ref{app:B}, and it is why transparency is stated with an additive
error rather than as an asymptotic equivalence. Here the degenerate exponent is visible in
closed form, since $d_1(r)$ vanishes for every $r\ge2$ and \eqref{eq:P_connection} reduces to
$A_1(n)=u_n$ exactly, so \eqref{eq:P_u_asymp} gives
$A_1(n)=\kappa_mn^{-m}+\mathcal O_m(n^{-m-1})$. The same exponent is where the two error
terms of \eqref{eq:P_transparency} exchange dominance, the profile correction $n^{-\beta-1}$
leading below $1$ and the mode $n^{-m}$ above.
\end{remark}

\begin{remark}[Relation to Appendix~\ref{app:G}]\label{rem:P_vs_H}
The kernel of Appendix~\ref{app:G} shares the branch $1-x/2$ on $(0,\tfrac12)$ and the
analytic index $2$ with the member $m=2$ of the present family, its transform being
$P(z)/(z-1)$ with $P(z)=z-(1-z/4)2^{z}$, whose zeros form an infinite family on the line
$\Re z=2$, the endpoint $\beta=2$ being nontransparent and the two regimes on either side
of it open. Carrying the affine branch across the whole interval
replaces that family by the single real zero of \eqref{eq:P_transform}, leaving the pole at
$z=1$ and the index in place, and the equation becomes solvable. The present entry is
therefore the model in which every regime attached to the value $\eta=2$ can be read, and the
difficulty of Appendix~\ref{app:G} is isolated as the passage from one zero on the line to
the whole family.
\end{remark}

\begin{numobs}\label{numobs:P_check}
The recurrence \eqref{eq:P_recurrence} reproduces the triangular solution of the defining
equation to $8\cdot10^{-16}$ over $n\le300$ at $m=2,3,5$ and two exponents. At $m=2$ the
identity $v_n=(n-1)2^{-n}$ and the Casoratian \eqref{eq:P_casoratian} hold with residual
exactly zero in rational arithmetic up to $n=40$, and the series \eqref{eq:P_gf} matches
$u_n$ exactly to the twelfth coefficient. The transparency constants of
\eqref{eq:P_transparency} are reached at $m=2,3,4$ and $\beta\in\{0,\tfrac12,1,\tfrac32\}$,
the measured values at $n=4\cdot10^{5}$ agreeing with $m(1-\beta)/(m-\beta)$ to five decimals
away from the degenerate exponent. The constants $\kappa_m$ of \eqref{eq:P_u_asymp} measure
$-2.00008$, $24.0042$, $-648.29$, $30748$ and $-2253660$ at $m=2,\dots,6$ and $n=1.2\cdot10^{5}$,
against $-2$, $24$, $-648$, $30720$ and $-2250000$. Formula \eqref{eq:P_zeta} is reached at
$\beta=\tfrac52,3,4,6,10$, the measured values $-5.2153$, $-3.28990$, $-2.40414$, $-2.07388$
and $-2.00404$ against $-2\zeta(\beta-1)$ equal to $-5.2248$, $-3.28987$, $-2.40411$,
$-2.07386$ and $-2.00402$.
\end{numobs}

\galleryentry{Q}{The binomial harmonic kernel}
 {$G(n,k)=\tfrac12\bigl(1+\binom nk^{-1}\bigr)$, for $1\le k\le n$}
 {regular arithmetic function, $\alpha(G)=1$}
 {$G^{*}(z)=\tfrac12$, constant, no zero, so $\eta$ is not defined}
 {$\tau(G)=\alpha(G)=1$, $\eta$ not defined}
 {proved, Theorem~\ref{thm:Q_index}, transparency with $\Xi_G=1/G^{*}$ on every $\beta<1$, index attained, failure above $1$}

The kernel of this appendix is the arithmetic mean of the constant $1$ and the reciprocal of a
binomial coefficient,
\begin{equation}\label{eq:Q_kernel}
G(n,k)=\frac12\left(1+\binom nk^{-1}\right),\qquad 1\le k\le n .
\end{equation}
It is admissible, since $G(n,n)=1$. Away from the two edges the binomial coefficient is large and
the kernel is close to the constant $\tfrac12$, so the whole arithmetic content sits in the two
ranks $k=1$ and $k=n-1$, where $\binom nk^{-1}=1/n$. The appendix shows that this single edge
weight fixes the regularity index at $1$, and that it does so twice over, once as the exponent of
absorption and once as the point where two asymptotic scales cross.

\begin{proposition}\label{prop:Q_transform}
The finite probes of $G$ converge, locally uniformly on $\Re z<0$, to the constant
$G^{*}(z)=\tfrac12$. The transform has no zero, so no analytic index is available and the
regularity index of $G$ can only come from the arithmetic of the equation.
\end{proposition}

\begin{proof}
Write the probe as
$-\tfrac zn\sum_{k\le n}G(n,k)(k/n)^{-z-1}$ and split $G$ into its constant part and its
reciprocal part. The constant part contributes
$\tfrac12\cdot\bigl(-z\int_0^1t^{-z-1}\,dt\bigr)+o(1)=\tfrac12+o(1)$ by the Riemann sum of
Proposition~\ref{prop:mellin_coincidence}. For the reciprocal part, the rank $k=n$ contributes
$-z/(2n)$, the ranks $k=1$ and $k=n-1$ contribute $\mathcal{O}(|z|\,n^{\Re z-1})$ and
$\mathcal{O}(|z|/n^{2})$, and the remaining ranks are bounded by
$\tfrac{|z|}{2n}\binom n2^{-1}\sum_{k}(k/n)^{-\Re z-1}$, which is $\mathcal{O}(|z|/n^{2})$ on
$\Re z<0$. All of these tend to zero, uniformly on compact subsets of the half plane.
\end{proof}

\subsection*{The exact reduction}

The defining equation reduces to a scalar recurrence in the partial sums.

\begin{proposition}\label{prop:Q_reduction}
Let $(a_k)$ solve the defining equation at the forcing $n^{-\beta}$ and put
$A(n)=\sum_{k\le n}a_k$ and
\begin{equation}\label{eq:Q_tail}
T(n)=\sum_{k=1}^{n-1}\binom nk^{-1}a_k .
\end{equation}
Then, exactly at every rank $n\ge1$, with $A(0)=0$,
\begin{equation}\label{eq:Q_reduction}
A(n)=\tfrac12A(n-1)+n^{-\beta}-\tfrac12T(n) .
\end{equation}
\end{proposition}

\begin{proof}
Splitting the defining equation along \eqref{eq:Q_kernel} gives
$\tfrac12A(n)+\tfrac12\sum_{k\le n}\binom nk^{-1}a_k=n^{-\beta}$. The rank $k=n$ has
$\binom nn^{-1}=1$ and contributes $a_n=A(n)-A(n-1)$, the remaining ranks are $T(n)$, and
rearranging gives \eqref{eq:Q_reduction}.
\end{proof}

The recursion \eqref{eq:Q_reduction} is a contraction of ratio $\tfrac12$, so everything follows
from the size of its second member. The next lemma isolates the two edge ranks, which are the
only ones that survive.

\begin{lemma}\label{lem:Q_tail}
For a fixed $n\ge12$, suppose $|a_k|\le M$ for $1\le k\le n-1$. Then
\begin{equation}\label{eq:Q_tail_split}
T(n)=\frac{a_1+a_{n-1}}{n}+\theta(n),\qquad |\theta(n)|\le\frac{13M}{n^{2}} ,
\end{equation}
and $a_1=1$.
\end{lemma}

\begin{proof}
The ranks $k=1$ and $k=n-1$ both have $\binom nk=n$, which gives the displayed fraction, and
$\theta(n)$ is the sum of the remaining ranks $2\le k\le n-2$. On that range
$\binom nk\ge\binom n2=n(n-1)/2$, and on the shorter range $3\le k\le n-3$, of cardinality
$n-5$, the sharper bound $\binom nk\ge\binom n3=n(n-1)(n-2)/6$ holds, both by unimodality.
Keeping the two ranks $k=2$ and $k=n-2$ apart,
\[
\sum_{2\le k\le n-2}\binom nk^{-1}
\le\frac{4}{n(n-1)}+\frac{6(n-5)}{n(n-1)(n-2)}
\le\frac{4}{n(n-1)}+\frac{6}{(n-1)(n-2)} .
\]
The right side times $n^{2}$ decreases in $n$ and equals $12.24\ldots$ at $n=12$, so it is at
most $13/n^{2}$ for $n\ge12$, and multiplying by $M$ gives the bound on $\theta$. Finally the
defining equation at $n=1$ reads $a_1G(1,1)=1$ with $G(1,1)=1$, so $a_1=1$.
\end{proof}

The solution is bounded a priori, uniformly in the rank.

\begin{lemma}\label{lem:Q_bounded}
For every $\beta>0$ the solution satisfies $|A(n)|\le M_\beta$ and $|a_n|\le2M_\beta$ for every
$n$, with a constant depending on $\beta$ alone.
\end{lemma}

\begin{proof}
Let $M_n=\max_{k\le n}|A(k)|$, so that $|a_k|\le2M_{n-1}$ for every $k\le n-1$. Applying
\eqref{eq:Q_tail_split} with $M=2M_{n-1}$, and using $a_1=1$ together with
$|a_{n-1}|\le2M_{n-1}$, equation \eqref{eq:Q_reduction} gives, for $n\ge12$,
\[
|A(n)|\le\Bigl(\frac12+\frac1n+\frac{13}{n^{2}}\Bigr)M_{n-1}+n^{-\beta}+\frac{1}{2n} .
\]
For $n\ge16$ the bracket is at most $\tfrac23$, and for $\beta>0$ the two remaining terms are at
most $\tfrac32$. Hence $|A(n)|\le\tfrac23M_{n-1}+\tfrac32$, so $M_n=M_{n-1}$ as soon as
$M_{n-1}\ge\tfrac92$, and in every case $M_n\le\max\bigl(M_{15},\tfrac92\bigr)$. The finitely many
ranks below $16$ depend on $\beta$ alone.
\end{proof}

For $\beta<0$ the forcing grows, so the preceding uniform argument no longer applies. At
$\beta=0$ it still applies, but the normalized formulation below includes that endpoint and
treats the whole range $\beta\le0$ uniformly.

\begin{lemma}\label{lem:Q_bounded_nonpositive}
For every $\beta\le0$ the solution satisfies $|A(n)|\le D_\beta n^{-\beta}$ and
$|a_n|\le2D_\beta n^{-\beta}$ for every $n$, with a constant depending on $\beta$ alone.
\end{lemma}

\begin{proof}
Let $D_n=\max_{1\le k\le n}k^{\beta}|A(k)|$. Since $\beta\le0$ the map
$k\mapsto k^{-\beta}$ is nondecreasing, so $|A(k)|\le D_{n-1}(n-1)^{-\beta}$ for every
$1\le k\le n-1$. Together with $A(0)=0$, this gives
$|a_k|\le|A(k)|+|A(k-1)|\le2D_{n-1}(n-1)^{-\beta}$ throughout that range. Those are exactly
the ranks occurring in $T(n)$, so \eqref{eq:Q_tail_split} applies with
$M=2D_{n-1}(n-1)^{-\beta}$, and
\eqref{eq:Q_reduction} gives, for $n\ge12$,
\[
|A(n)|\le\tfrac12D_{n-1}(n-1)^{-\beta}+n^{-\beta}+\frac1{2n}
+\frac{D_{n-1}(n-1)^{-\beta}}{n}+\frac{13\,D_{n-1}(n-1)^{-\beta}}{n^{2}} .
\]
Multiply by $n^{\beta}$ and use the two inequalities $(n-1)^{-\beta}\le n^{-\beta}$ and
$n^{\beta}\le1$, both consequences of $\beta\le0$. The three terms carrying $D_{n-1}$ collect into
the factor $\tfrac12+\tfrac1n+\tfrac{13}{n^{2}}$, which is at most $\tfrac34$ for $n\ge12$, while
the two remaining terms are at most $\tfrac32$. Hence
\[
n^{\beta}|A(n)|\le\tfrac34D_{n-1}+\tfrac32\qquad(n\ge12).
\]
The finite number $D_{11}$ depends only on $\beta$. Set
\[
D_\beta=\max(D_{11},6).
\]
Since
\[
D_n=\max\{D_{n-1},n^{\beta}|A(n)|\}
\]
and the affine map $x\mapsto\tfrac34x+\tfrac32$ has fixed point $6$, induction gives
$D_n\le D_\beta$ at every rank. Hence $|A(n)|\le D_\beta n^{-\beta}$. Finally, for $n\ge2$,
\[
|a_n|\le |A(n)|+|A(n-1)|
\le D_\beta\bigl(n^{-\beta}+(n-1)^{-\beta}\bigr)
\le2D_\beta n^{-\beta},
\]
and the case $n=1$ follows from $a_1=1$ and $D_\beta\ge6$.
\end{proof}

\subsection*{The index}

The transparent range and the failure above it follow.

\begin{theorem}\label{thm:Q_index}
Every $\beta<1$ is transparent for the kernel \eqref{eq:Q_kernel}, with transparent
coefficient $\Xi_G(\beta)=2=1/G^{*}(\beta)$, so the identification holds throughout the
transparent range. At the endpoint $\beta=1$ the partial sums still follow an exact power,
$A(n)=n^{-1}+o(n^{-1})$, with the different constant $1$, and no $\beta>1$ is transparent, the
partial sums satisfying
$A(n)=-n^{-1}+o(n^{-1})$ there. Absorption holds at the exponent $1$ and that exponent is attained.
The kernel is therefore a regular arithmetic function\index[terms]{regular arithmetic function} with
$\tau(G)=\alpha(G)=1$.
\end{theorem}

\begin{proof}
Iterating the contraction \eqref{eq:Q_reduction} from $A(0)=0$ gives, exactly,
\begin{equation}\label{eq:Q_iterated}
A(n)=\sum_{j=0}^{n-1}2^{-j}\Bigl[(n-j)^{-\beta}-\tfrac12T(n-j)\Bigr] .
\end{equation}
Write $\beta^{-}=\min(\beta,0)$. Lemma~\ref{lem:Q_bounded} on the positive range and
Lemma~\ref{lem:Q_bounded_nonpositive} on the rest give $|a_m|\le2D_\beta m^{-\beta^{-}}$ at every
rank, so Lemma~\ref{lem:Q_tail} applies at the rank $m$ with $M=2D_\beta m^{-\beta^{-}}$ and
$T(m)=1/m+a_{m-1}/m+\mathcal{O}(m^{-\beta^{-}-2})$.

Treat the two sums in \eqref{eq:Q_iterated} separately. In the first, the ranks with
$j\ge n/2$ contribute $\mathcal{O}(2^{-n/2}n^{|\beta|})$, and on the remaining range
$(n-j)^{-\beta}=n^{-\beta}\bigl(1+\beta j/n+\mathcal{O}(j^{2}/n^{2})\bigr)$. Since
$\sum_{j\ge0}2^{-j}=2$, $\sum_{j\ge0}j2^{-j}=2$ and $\sum_{j\ge0}j^{2}2^{-j}=6$,
\[
\sum_{j}2^{-j}(n-j)^{-\beta}=2n^{-\beta}+2\beta\,n^{-\beta-1}+\mathcal{O}(n^{-\beta-2}).
\]
In the second, the same expansion applied to $1/(n-j)$ gives
$\tfrac12\sum_j2^{-j}\bigl(1/(n-j)\bigr)=1/n+\mathcal{O}(n^{-2})$, while the ranks $\theta(n-j)$
contribute $\mathcal{O}(n^{-\beta^{-}-2})$ by \eqref{eq:Q_tail_split}. Writing $\eps_n$ for
the contribution of the terms $a_{m-1}/m$,
\begin{equation}\label{eq:Q_two_towers}
A(n)=2n^{-\beta}+2\beta\,n^{-\beta-1}-\frac1n+\mathcal{O}(n^{-\beta^{-}-2})
+\mathcal{O}(n^{-\beta-2})+\eps_n ,
\qquad
\eps_n=\mathcal{O}\Bigl(\frac1n\max_{m>n/2}|a_m|\Bigr).
\end{equation}
On the nonpositive range the crude bound already settles the matter. There
Lemma~\ref{lem:Q_bounded_nonpositive} gives $\eps_n=\mathcal{O}(n^{-\beta-1})$ and every
error term of \eqref{eq:Q_two_towers} is $\mathcal{O}(n^{-\beta-1})$, so
$A(n)=2n^{-\beta}+\mathcal{O}(n^{-\beta-1})$. Differencing at the ranks $n$ and $n-1$ gives
$a_n=\mathcal{O}(n^{-\beta-1})$, hence $\eps_n=\mathcal{O}(n^{-\beta-2})$, and the display
below is reached on that range at once.
On the positive range the bound of Lemma~\ref{lem:Q_bounded} gives only
$\eps_n=\mathcal{O}(n^{-1})$, which is
too crude, so the estimate is run twice. At the first pass, dropping $\eps_n$ into the
error, \eqref{eq:Q_two_towers} already gives $A(n)=2n^{-\beta}+\mathcal{O}(n^{-\min(\beta,1)})$
uniformly on compact sets of exponents, hence, by differencing at the ranks $n$ and $n-1$,
\[
a_n=A(n)-A(n-1)=\mathcal{O}\bigl(n^{-\min(\beta+1,\,2)}\bigr),
\]
the two contributions being $-2\beta n^{-\beta-1}$ from the forced term and $n^{-2}$ from the
edge term. Feeding this back gives
$\eps_n=\mathcal{O}(n^{-\min(\beta+2,\,3)})$, which is absorbed in the two error terms
already present. At the second pass, therefore,
\begin{equation}\label{eq:Q_two_towers_final}
A(n)=2n^{-\beta}+2\beta\,n^{-\beta-1}-\frac1n
+\mathcal{O}\bigl(n^{-\min(\beta+2,\,2)}\bigr).
\end{equation}

The three clauses read off \eqref{eq:Q_two_towers_final}. For $\beta<1$ the term $2n^{-\beta}$ is
of strictly larger order than $n^{-1}$, so $A(n)=2n^{-\beta}+o(n^{-\beta})$ and $\beta$ is
transparent with coefficient $2$, which is $1/G^{*}(\beta)$ by
Proposition~\ref{prop:Q_transform}. At $\beta=1$ the first and third terms combine into
$2n^{-1}-n^{-1}=n^{-1}$ and the second is $\mathcal{O}(n^{-2})$, so $A(n)=n^{-1}+\mathcal{O}(n^{-2})$
and $\beta=1$ is transparent with coefficient $1$. For $\beta>1$ the term $-n^{-1}$ dominates
and $A(n)=-n^{-1}+o(n^{-1})$, so $n^{\beta}A(n)$ is unbounded and $\beta$ is not transparent,
while $A(n)=\mathcal{O}(n^{-1})$ gives absorption at the exponent $1$ with the exponent attained.
Every exponent below $1$ is therefore transparent, the exponent $1$ is transparent as well, and no
exponent above $1$ is. The threshold is thus $\tau(G)=1$, and every interval $[1,1+\delta)$
contains a nontransparent exponent, so that value is sharp. With the absorption just obtained the
three conditions of Remark~\ref{rem:index_three_conditions} hold, which makes $G$ a regular
arithmetic function with $\alpha(G)=1$.
\end{proof}

\begin{remark}[Where the index comes from]\label{rem:Q_source}
The transform is constant and has no zero, so the analytic index is undefined and cannot be the
source of the threshold. The threshold is the exponent of the edge weight. The rank $k=1$ carries
the weight $\tfrac12(1+1/n)$, whose deviation from the constant $\tfrac12$ is exactly
$1/(2n)$, and the coefficient it multiplies is $a_1=1$. Passing that deviation through the
geometric sum of the contraction, of total mass $\sum_j2^{-j}=2$, produces the term $-1/n$ of
\eqref{eq:Q_two_towers_final}, and the index is the exponent of that term. Nothing else in the
kernel plays a part.
\end{remark}

\begin{remark}[The coefficient at the endpoint]\label{rem:Q_coefficient}
Below the index the transparent coefficient is $2=1/G^{*}(\beta)$, in agreement with the
identification recorded after Definition~\ref{def:reg_index}, so the identification holds
throughout the transparent range of this kernel. At the endpoint the two towers of
\eqref{eq:Q_two_towers_final} come to the same order and the surviving coefficient is $1$, while
$1/G^{*}(1)=2$. The threshold is a supremum over exponents strictly below it and records
absorption at the endpoint, so nothing is asserted there and nothing is contradicted. What the
kernel exhibits is that the endpoint can carry an exact power of its own, with a constant the
transform does not supply, and that is why the identification is stated on the open range. The value $1$ is the difference of the
two masses, $2-1$, and it is proved by Theorem~\ref{thm:Q_index}.
\end{remark}

\begin{remark}[Two scales and their crossing]\label{rem:Q_towers}
Formula \eqref{eq:Q_two_towers_final} carries two independent families. The forced tower
$2n^{-\beta}+2\beta n^{-\beta-1}+\cdots$ follows the second member and moves with $\beta$. The
edge tower $-n^{-1}+\cdots$ comes from the kernel and does not move. The remainder in
\eqref{eq:Q_two_towers_final} is $\mathcal{O}(n^{-\beta-2})$ while the forced tower leads, and
$\mathcal{O}(n^{-2})$ once the edge tower leads, so the error exponent is
$\min(\beta+2,2)$ and the two regimes exchange at $\beta=1$. The index is the abscissa at which
the two towers cross, and it is the same number as the exponent of absorption. This is the
mechanism of Appendix~\ref{app:P} seen through an asymptotic expansion rather than through the
coalescence of characteristic roots.
\end{remark}

\begin{numobs}[Numerical control]\label{numobs:Q_check}
Forward substitution in the defining equation, in double precision and up to $n=6\cdot10^{4}$,
gives the following. The transparent constant is confirmed below the index, $n^{\beta}A(n)$
reaching $1.9999833$ at $\beta=-\tfrac12$, $1.9997475$ at $\beta=\tfrac14$ and $1.9959341$ at
$\beta=\tfrac12$, all at $n=6\cdot10^{4}$ and all increasing toward $2$ at the predicted rate
$n^{-(1-\beta)}$. At $\beta=0.9$ the same quantity is only $1.6672256$, against the value
$2-n^{-1/10}=1.6672\ldots$ predicted by \eqref{eq:Q_two_towers_final}, which is the same law
read near the index. The
critical value is confirmed, $nA(n)=1.0000250$ at $\beta=1$, against the predicted limit $1$ with
error $\mathcal{O}(n^{-1})$. Above the index $nA(n)$ measures $-0.3343926$ at $\beta=1.1$,
$-0.9918383$ at $\beta=\tfrac32$, $-0.9999667$ at $\beta=2$ and $-0.9999958$ at $\beta=3$,
against the predicted limit $-1$ and against the finite rank predictions
$2n^{-1/10}-1=-0.33448\ldots$ and $2n^{-1/2}-1=-0.99183\ldots$ at the first two exponents.

The coefficient $2\beta$ of \eqref{eq:Q_two_towers_final} is confirmed where the numerics can
separate it. Writing $R(n)=A(n)-2n^{-\beta}+n^{-1}$, the quantity $n^{\beta+1}R(n)$ measures
$-1.00003$ at $\beta=-\tfrac12$, $0.49972$ at $\beta=\tfrac14$ and $0.99632$ at $\beta=\tfrac12$
at $n=6\cdot10^{4}$, against $2\beta=-1$, $\tfrac12$ and $1$. Nearer the index the same
measurement does not stabilize, for the reason recorded in Remark~\ref{rem:Q_towers}. The
difficulty of measuring the coefficient is itself a measurement of the index.
\end{numobs}

\dossierchapter{sqrt2}{The broken harmonic kernel at $\sqrt2$}

This dossier carries the proofs of the statements of
Section~\ref{sec:sqrt2}, together with the exact block recurrences, the finite
certificates, and the auxiliary results they rest on. The notation is the one fixed
there, the profile and its induced kernel, the forcing equation, the block variable
$Q_n$, the orbit $(x_j)$, and the orbit trace $R_p$. A reader who wants only the
outcome may read Section~\ref{sec:sqrt2} and the closing block below, which records
what is proved, what is open, and what is not claimed.

\rafgalleryfig{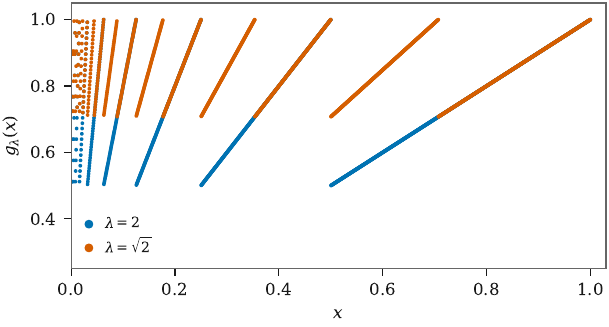}{The profile $g_\lambda$ at the two scales $\lambda=2$ and $\lambda=\sqrt2$. The teeth of the fan sit at $x=\lambda^{-i}$, so halving the scale interleaves a second family of jump points between those of the first. At $\lambda=2$ every jump point is a dyadic rational and the orbit of $n\lambda^{-i}$ modulo one is finite. At $\lambda=\sqrt2$ every second jump point is irrational, and that alternative is the Diophantine separation this dossier records.}{fig:dossier_sqrt2_profile}

\subsection*{Auxiliary statements and certificates}

The statements collected here are the ones the proofs need and the reader of
Section~\ref{sec:sqrt2} does not. They carry the exact closures, the finite
certificates, and the two negative results that delimit the remaining route.

Several forcing exponents are compared in what follows, so the subscript is kept on $A_\beta$ and on the sequences attached to it. Elsewhere in the volume the exponent is fixed by its context and the plain $a_n$ and $A(n)$ of Definition~\ref{def:reg_index_fgv} are used.

\begin{lemma}
\label{lem:w10-profile}
For every $\lambda>1$, the profile $g_\lambda$ is bounded and
Riemann-integrable on $[0,1]$. More precisely,
\[
\lambda^{-1}<g_\lambda(x)\leq1
\quad(0<x\leq1),
\qquad g_\lambda(1)=1.
\]
Its discontinuities in $(0,1]$ are exactly the points
$\lambda^{-j}$ with $j\geq1$. Hence \eqref{eq:g_lambda_def} defines an
admissible univariate profile and \eqref{eq:w10-profile} an induced
kernel. This assertion does not give the RAF property.
\end{lemma}

The mixed closure of the block variables is exact and closes after finitely many cells.

\begin{theorem}
\label{thm:w10-binary}
For $n\geq4$ and $-3\leq k\leq3$, put
\[
A_k(n)=Q_{n+k},
\qquad P_k(n)=Q_{\lfloor qn\rfloor+k},
\qquad
\mathcal S(n)=((A_k)_{k=-3}^3,(P_k)_{k=-3}^3)^{\mathsf T}.
\]
For $\epsilon\in\{0,1\}$ set
\[
n'=2n+\epsilon,
\quad P=\lfloor qn\rfloor,
\quad z=\{qn\},
\quad d=\lfloor2z+\epsilon q\rfloor,
\quad P'=2P+d,
\quad z'=\{2z+\epsilon q\}.
\]
Define
\begin{equation}
\label{eq:w10-mixed-selectors}
\alpha_{\epsilon,k}(z)=
\left\lfloor z+\frac{\epsilon+k}{q}\right\rfloor,
\quad
\gamma_k(z')=\left\lfloor\frac{k-z'}q\right\rfloor,
\end{equation}
and
\[
\xi^A_{\epsilon,k}(z)=
\mathbf1_{\{\{z+(\epsilon+k)/q\}<1/q\}},
\qquad
\xi^P_k(z')=
\mathbf1_{\{\{(k-z')/q\}<1/q\}}.
\]
The fourteen coordinates at $n'$ satisfy
\begin{align}
A'_k={}&\delta_p(n'+k)
+q\mathbf1_{2\mid\epsilon+k}A_{(\epsilon+k)/2}
-(q-1)\xi^A_{\epsilon,k}(z)P_{\alpha_{\epsilon,k}(z)},
\label{eq:w10-mixed-A}\\
P'_k={}&\delta_p(P'+k)
+q\mathbf1_{2\mid d+k}P_{(d+k)/2}
-(q-1)\xi^P_k(z')A'_{\gamma_k(z')}.
\label{eq:w10-mixed-P}
\end{align}
Every index on the right lies in $\{-3,\ldots,3\}$. If
$U_\epsilon,V_{\epsilon,z},W_{\epsilon,z},X_{z'}$ denote the coefficient
blocks displayed by \eqref{eq:w10-mixed-A} and \eqref{eq:w10-mixed-P}, then
\begin{equation}
\label{eq:w10-mixed-matrix}
\boxed{
\mathcal S(n')=
\begin{pmatrix}
U_\epsilon&V_{\epsilon,z}\\
X_{z'}U_\epsilon&W_{\epsilon,z}+X_{z'}V_{\epsilon,z}
\end{pmatrix}\mathcal S(n)
+\binom{f_A}{f_P+X_{z'}f_A}.}
\end{equation}
The four states $n=4,5,6,7$ are the finite bases. Thus every $n\geq4$
is covered by an exact affine binary system of dimension fourteen.
\end{theorem}

The cells and edges of that closure are listed once here, and the list is used throughout what
follows.

\begin{lemma}
\label{lem:w10-mixed-certificate}
The mixed closure has twenty-eight open cells in
$\{0,1\}\times[0,1)$. The ordered endpoints for $\epsilon=0$ are
\begin{equation}
\label{eq:w10-B0}
\begin{split}
\mathcal B_0=(&0,\tfrac32-q,\tfrac{3q-4}2,\tfrac{q-1}2,
1-\tfrac q2,q-1,\tfrac12,2-q,\\
&\tfrac{3(q-1)}2,\tfrac q2,\tfrac{3-q}2,
2q-2,3-\tfrac{3q}2,q-\tfrac12,1).
\end{split}
\end{equation}
The ordered endpoints for $\epsilon=1$ are
\begin{equation}
\label{eq:w10-B1}
\begin{split}
\mathcal B_1=(&0,\tfrac32-q,\tfrac{3q-4}2,3-2q,
\tfrac{q-1}2,1-\tfrac q2,\tfrac{5-3q}2,\\
&q-1,\tfrac12,2-q,\tfrac q2,\tfrac{3-q}2,
3-\tfrac{3q}2,q-\tfrac12,1).
\end{split}
\end{equation}
Let $I_{\epsilon,i}$ be the interval between consecutive endpoints.
A node $t_r=(\epsilon_r,i_r)$ records
\[
n_r=2n_{r-1}+\epsilon_r,
\qquad
\{qn_{r-1}\}\in I_{\epsilon_r,i_r}.
\]
The complete graph is
{\scriptsize
\[
\begin{array}{c|c|c@{\quad}c|c|c}
0.1&0:1,2,3&1:1,2,3&1.1&0:6,7&1:8,9\\
0.2&0:3,4&1:4,5&1.2&0:8,9&1:10\\
0.3&0:4,5&1:5,6,7&1.3&0:9,10&1:10,11\\
0.4&0:6,7&1:8,9&1.4&0:10,11&1:11,12\\
0.5&0:8,9,10,11&1:10,11,12&1.5&0:12,13,14&1:12,13,14\\
0.6&0:12,13,14&1:12,13,14&1.6&0:1,2,3&1:1,2,3\\
0.7&0:1,2,3&1:1,2,3&1.7&0:3,4&1:4,5\\
0.8&0:3,4&1:4,5&1.8&0:4,5&1:5,6,7\\
0.9&0:4,5&1:5,6,7&1.9&0:6,7&1:8,9\\
0.10&0:6,7&1:8,9&1.10&0:8,9,10,11&1:10,11,12\\
0.11&0:8,9&1:10&1.11&0:12,13,14&1:12,13,14\\
0.12&0:9,10&1:10,11&1.12&0:1,2,3&1:1,2,3\\
0.13&0:10,11&1:11,12&1.13&0:3,4&1:4,5\\
0.14&0:12,13,14&1:12,13,14&1.14&0:4,5&1:5,6,7
\end{array}
\]
}
An entry $0:1,2$ lists the targets $I_{0,1}$ and $I_{0,2}$.
The graph contains exactly $136$ edges. For the coordinate order
$A_{-3},\ldots,A_3\mid P_{-3},\ldots,P_3$, assign signs by
\[
\begin{array}{c|c}
0.1,0.2,0.13,0.14&+------\mid+++++++\\
0.3,\ldots,0.7&-------\mid+++++++\\
0.8,\ldots,0.11&------+\mid+++++++\\
0.12&+-----+\mid+++++++\\
1.1,\ldots,1.7&-------\mid+++++++\\
1.8,\ldots,1.10&-+-----\mid+++++++\\
1.11,\ldots,1.14&-------\mid+++++++
\end{array}
\]
where the space in the third pattern is typographical only. If
$\Sigma_t$ is the resulting diagonal matrix and $s\to t$ is an edge, then
\begin{equation}
\label{eq:w10-orthant}
\Sigma_t\mathcal M_t\Sigma_s\geq0
\end{equation}
entrywise. The stacked system $\mathcal M_t v=qv$ has rank fourteen.
Thus there is no common right $q$-mode. The certificate proves exact phase
coverage and compatible orthants. It does not prove the missing joint
spectral estimate.
\end{lemma}

The raw sums decompose exactly, into an orbit part and a defect part.

\begin{lemma}
\label{lem:w10-raw-decomposition}
Let $M=s(N)$, $K=s(M)$, and $\delta_m^B=B(m)-qm$. Put
\[
V_p(N)=\sum_{j\geq0}F_p(x_j),
\qquad
\mathcal J_NW=\sum_{M<n\leq N}\frac{W(n)-W(n-1)}n.
\]
Then
\begin{equation}
\label{eq:w10-raw-decomposition}
\boxed{
A_\beta(N)=\mathcal J_NV_p+G_\beta(N)
+\sum_{K<m\leq M}\frac{\delta_m^BQ_m}{mB(m)}.}
\end{equation}
For $0<p<1$,
\begin{equation}
\label{eq:w10-floorfree-main}
\mathcal J_NV_p=\frac{N^{p-1}}{g^*(1-p)}+O_p(N^{-1}).
\end{equation}
For $p=1$, it is $1/g^*(0)+O(\log N/N)$. For $p\leq0$, it is
$O_p(\log N/N)$. For every $\beta\geq0$ and $\eps>0$, the last
sum in \eqref{eq:w10-raw-decomposition} is
$O_{\beta,\eps}(N^{-1/2+\eps})$.
\end{lemma}

One Abel compression turns the defect into a centered sum.

\begin{theorem}
\label{thm:w10-Abel}
Put $M=s(N)$ and define the centered defect
\begin{equation}
\label{eq:w10-G}
G_\beta(N)=
\sum_{M<n\leq N}\Delta R_p(n)
\left(\frac1n-\frac1{qH}\right).
\end{equation}
If $K_j=B^j$, then
\begin{equation}
\label{eq:w10-G-double}
G_\beta(N)=
\sum_{j\geq0}\sum_{M<K_j(m)\leq N}
\Delta\mathcal Z_p(m)
\left(\frac1{K_j(m)}-\frac1{qH}\right).
\end{equation}
Let $J=\max\{j:x_j>0\}$, $x_{J+1}=0$, and
\begin{align}
\mathcal I_p(N)
&=\sum_{j=0}^{J}\sum_{x_{j+1}<m<x_j}
\mathcal Z_p(m)
\left(\frac1{K_j(m)}-\frac1{K_j(m+1)}\right),
\label{eq:w10-I}\\
c_0&=\frac1N-\frac1{qH},
\qquad
c_j=\frac1{K_j(x_j)}-\frac1{K_{j-1}(x_j+1)}
\quad(1\leq j\leq J).\notag\end{align}
Two exact summations by parts give
\begin{align}
G_\beta(N)&=\sum_{j=0}^{J}c_j\mathcal Z_p(x_j)+\mathcal I_p(N),
\label{eq:w10-first-Abel}\\
\sum_{j=0}^{J}c_j\mathcal Z_p(x_j)
&=c_0R_p(N)+\sum_{j=1}^{J}(c_j-c_{j-1})R_p(x_j).
\label{eq:w10-second-Abel}
\end{align}
For every $\beta\geq0$ and $\eps>0$,
\begin{equation}
\label{eq:w10-Abel-errors}
\mathcal I_p(N)=O_{\beta,\eps}(N^{-1/2+\eps}),
\qquad
\sum_{j=1}^{J}(c_j-c_{j-1})R_p(x_j)
=O_{\beta,\eps}(N^{-1+\eps}).
\end{equation}
Consequently
\begin{equation}
\label{eq:w10-compression}
\boxed{
G_\beta(N)=
\left(\frac1N-\frac1{q\lfloor N/2\rfloor}\right)R_p(N)
+O_{\beta,\eps}(N^{-1/2+\eps}).}
\end{equation}
The leading coefficient obeys
\begin{equation}
\label{eq:w10-c0-size}
\frac{q-1}{N}\leq|c_0|\leq\frac qN.
\end{equation}
\end{theorem}

The orbit trace and the pointwise quantity are equivalent term by term.

\begin{lemma}
\label{lem:w10-trace-equivalence}
For every $n\geq1$,
\begin{equation}
\label{eq:w10-new-trace-identity}
\boxed{qQ_n=R_p(2n)-R_p(\lfloor qn\rfloor).}
\end{equation}
This identity is an exact consequence of the acquired renewal. It gives
\begin{align*}
R_p(N)=o(N^p)
&\quad\Longleftrightarrow\quad Q_n=o(n^p)
&& (1/2<p\leq1),
\\
R_p(N)=O_{p,\eps}(N^{1/2+\eps})
&\quad\Longleftrightarrow\quad
Q_n=O_{p,\eps}(n^{1/2+\eps})
&& (p\leq1/2).
\end{align*}
The constants in the second equivalence may depend on $p$ and
$\eps$. Thus LOW and ABS are respectively equivalent, in their
stated domains, to the displayed pointwise estimates for $Q_n$.
\end{lemma}

The critical ray carries its own certificate, obtained by a conjugation that makes the quotient
explicit.

\begin{lemma}
\label{lem:w10-critical-certificate}
On the ray $u=1$, quotient the common coordinate of $q^{-1}M(z)$ and
conjugate by
\[
\Sigma=\operatorname{diag}(-1,-1,-1,1,1,1,1).
\]
The resulting nonnegative matrices $C_i$ use the nine phase cells
\begin{equation}
\label{eq:w10-nine-cells}
\begin{aligned}
I_4&=(0,\tau),
&I_5&=(\tau,(q-1)/2),
&I_2&=((q-1)/2,q-1),\\
I_6&=(q-1,1/2),
&I_3&=(1/2,3(q-1)/2),
&I_1&=(3(q-1)/2,q/2),\\
I_9&=(q/2,2(q-1)),
&I_8&=(2(q-1),q-1/2),
&I_7&=(q-1/2,1),
\end{aligned}
\end{equation}
where $\tau=(3q-4)/2$. Their complete graph is
\begin{equation}
\label{eq:w10-nine-graph}
\begin{aligned}
1&\to2,
&2&\to1,3,6,9,
&3&\to2,4,5,
&4&\to2,4,5,
&5&\to2,\\
6&\to7,8,
&7&\to7,8,
&8&\to1,9,
&9&\to1,3,6.
\end{aligned}
\end{equation}
It has twenty-one edges. The exact gauges $h_i^-$ and $h_i^+$ are listed
in Table~\ref{tab:w10-gauges}. Their common positive-coordinate support is
denoted by $R_i$. On every edge $i\to j$ they satisfy
\begin{equation}
\label{eq:w10-gauge-contract}
C_jh_i^-\leq a_0h_j^+,
\qquad
C_jh_i^+\leq b_0h_j^-,
\quad
a_0=\frac{9q-8}{14},
\quad
b_0=\frac{11-6q}{7},
\quad a_0b_0=\tau<1.
\end{equation}
The minimum positive gauge coordinate is
\begin{equation}
\label{eq:w10-gauge-minimum}
m_0=\frac{-106+75q}{14}>\frac{47}{10000}.
\end{equation}
It is attained only at coordinate $3$ of $h_9^+$. Every coordinate is at
most one. The actual signatures at ranks $4,5,6,7$ are $(2,1,2,3)$ and
the support check propagates $R_i$ to $R_j$ on each edge. An arbitrary
quotient image has only one possible extra component, a multiple of $e_3$
in signatures $4$ and $5$. The neutral relations are
\begin{equation}
\label{eq:w10-neutral-vector}
C_4e_3=C_5e_3=C_2e_3=e_3.
\end{equation}
Moreover
\begin{equation}
\label{eq:w10-rank-one}
C_3C_2C_1C_2=uv^{\mathsf T},
\quad
u=(0,a,a,0,\tau,\tau,0)^{\mathsf T},
\quad
v=(0,0,0,1,\tau,0,0)^{\mathsf T},
\quad a=1-q/2.
\end{equation}
The common left $q$-mode system has only the zero solution.
\end{lemma}

The weighted trace is framed above and below by the same weighted sums, which is the comparison
used below.

\begin{proposition}
\label{prop:w10-trace-Schur}
For every $\sigma>1/2$,
\begin{equation}
\label{eq:w10-trace-upper}
\sum_{N\geq1}\frac{|R_p(N)|^2}{N^{2\sigma}}
\leq C_\sigma
\sum_{m\geq1}\frac{|Q_m|^2}{m^{2\sigma}},
\end{equation}
where one admissible constant is
\begin{equation}
\label{eq:w10-Schur-constant}
C_\sigma=q^2(q+3)2^{-2\sigma}
\left(1-q^{-(2\sigma-1)/2}\right)^{-2}.
\end{equation}
Conversely,
\begin{equation}
\label{eq:w10-trace-lower}
\sum_{n\geq1}\frac{|Q_n|^2}{n^{2\sigma}}
\leq\frac2{q^2}(2^{2\sigma}+q^{2\sigma})
\sum_{N\geq1}\frac{|R_p(N)|^2}{N^{2\sigma}}.
\end{equation}
Thus the orbit trace preserves the available weighted Hilbert scale. It
does not produce the pointwise gain required by LOW or ABS.
\end{proposition}

At the density level the raw bounds improve.

\begin{proposition}
\label{prop:w10-density}
For $\beta>0$, put
\[
D_p^{\mathrm{sh}}=
\sum_{n\geq1}\frac{|\delta_p(n)|}{n}<\infty.
\]
For every $X\geq2$,
\begin{equation}
\label{eq:w10-R-L1}
\sum_{X<N\leq2X}|R_p(N)|
\leq7q(q+3)D_p^{\mathrm{sh}}
X(1+\lfloor\log_qX\rfloor).
\end{equation}
At $\beta=0$, the same argument with the truncated harmonic source gives
$O(X(\log X)^2)$.

For $\beta\geq1/2$ and every $\eps>0$, the ABS bound holds on
the block $(X,2X]$ outside at most
\begin{equation}
\label{eq:w10-ABS-exception}
O_{\beta,\eps}
\bigl(X^{1/2-\eps/2}\log X\bigr)
\end{equation}
integers. For $0<\beta<1/2$, the LOW bound holds with the stronger factor
$1/\log X$ outside at most
\begin{equation}
\label{eq:w10-LOW-exception}
O_\beta\bigl(X^\beta(\log X)^2\bigr)
\end{equation}
integers. At $\beta=0$, the exceptional count is $O((\log X)^3)$.
These density statements do not give the
uniform endpoint estimates in Definition~\ref{def:w10-low-abs}.
\end{proposition}

At the critical exponent the trace admits a lower bound.

\begin{proposition}
\label{prop:w10-critical-trace}
At $p=1/2$,
\[
\limsup_{N\to\infty}\frac{|R_{1/2}(N)|}{\sqrt N}
\geq\frac{q}{q+\sqrt q}L
>\frac{q}{q+\sqrt q}\frac{1293}{5000}.
\]
Thus the order $\sqrt N$ is attained on a subsequence of one of the two
arms $2n$ and $\lfloor qn\rfloor$. The slack $N^\eps$ in ABS is
consistent with this lower bound.
\end{proposition}

The comparison with integral scales holds only in a limited range, which locates where the
argument stops.

\begin{proposition}
\label{prop:w10-integer-comparison}
For an integral scale $\lambda\geq2$, define locally
\[
B^{(\lambda)}_0(n):=\sum_{k\leq n}k a_k^{(0)}.
\]
The established endpoint identities include
$B^{(\lambda)}_0(n)=s_\lambda(n)$ at $\beta=0$, where $s_\lambda(n)$ is the
sum of the digits of $n$ in base $\lambda$, and the exact critical formulas
at $\beta=1$, including their values at powers of $\lambda$. This notation
is distinct from the Beatty\index[terms]{Beatty sequence}\index[names]{Beatty, S.} function $B(m)=\lceil qm\rceil$. At
$\beta=0$, differences of the digital sum are subpolynomial. In contrast,
\eqref{eq:w10-HLR-asymptotic} gives polynomial coefficient growth for
$\lambda=\sqrt2$ at $\beta=1/2$. No assertion about HLR for every forcing
or every integral scale is made here.
\end{proposition}

One admissible cycle of the graph expands, and that is the exact obstruction to the full target.

\begin{proposition}
\label{prop:w10-expanding-cycle}
The $136$-edge graph contains the admissible cycle
\begin{equation}
\label{eq:w10-expanding-cycle}
(0,6)\to(0,14)\to(0,12)\to(0,9)\to(1,7)
\to(0,4)\to(0,6).
\end{equation}
For the product $P$ of the six orthant-conjugated target matrices, its
$A_0$ diagonal entry is
\begin{equation}
\label{eq:w10-expanding-entry}
P_{A_0,A_0}=-76+60q>8=q^6.
\end{equation}
The cycle is dynamically realized by the periodic phase
\[
z_0=\frac{32-2q}{63}.
\]
Since the conjugated matrices are nonnegative,
$(P^k)_{A_0,A_0}\geq(-76+60q)^k$. Therefore the proposed estimate
\[
\|\mathcal M_{t_\ell}\cdots\mathcal M_{t_1}\|
\leq C_\eps q^\ell2^{\eps\ell}
\quad\text{for every }\eps>0
\]
is false. The universal all-path estimate fails even when the class of paths is
restricted to dynamically realized paths.
More precisely, it fails for every
\[
0<\eps<\frac16\log_2\left(\frac{-76+60q}{8}\right)
=0.0243560114\ldots.
\]
This disproves a sufficient route considered earlier. It does not disprove
LOW or ABS, since the expanding direction may be absent from, or cancelled
in, the actual forced observable.
\end{proposition}

The forcing along a dyadic ray is carried by an eight-dimensional cocycle.

\begin{lemma}
\label{lem:w10-eight-cocycle}
Fix an odd integer $u$. Put
\[
N_r=u2^r,
\qquad P_r=\lfloor qN_r\rfloor,
\qquad z_r=\{qN_r\},
\qquad
r_0(u)=\min\{r\geq0:u2^r\geq4\}.
\]
For $r\geq r_0(u)$ and $0\leq c\leq3$, set
\[
A_{r,c}=Q_{N_r-c},
\qquad P^Q_{r,c}=Q_{P_r-c},
\qquad
Z_r=(A_{r,0},\ldots,A_{r,3},P^Q_{r,0},\ldots,P^Q_{r,3})^{\mathsf T}.
\]
Define
\[
\begin{aligned}
D_c(z)&=\left\lceil\frac cq-z\right\rceil,
&C_c(z)&=\left\lceil\frac{c+z}{q}\right\rceil,\\
\xi_c^A(z)&=\mathbf1_{\{\{z-c/q\}<1/q\}},
&\xi_c^P(z)&=\mathbf1_{\{\{-(c+z)/q\}<1/q\}}.
\end{aligned}
\]
Let $b=\lfloor2z\rfloor$. On indices $0,1,2,3$, the nonzero entries of
$U,V(z),W(b),X(z)$ are
\[
\begin{aligned}
U_{c,c/2}&=q&&(2\mid c),
&V(z)_{c,D_c(z)}&=-(q-1)\xi_c^A(z),\\
W(b)_{c,(c-b)/2}&=q&&(c\equiv b\pmod2),
&X(z)_{c,C_c(z)}&=-(q-1)\xi_c^P(z).
\end{aligned}
\]
With $z'=\{2z\}$, put
\begin{equation}
\label{eq:w10-eight-matrix}
M(z)=
\begin{pmatrix}
U&V(z)\\
X(z')U&W(\lfloor2z\rfloor)+X(z')V(z)
\end{pmatrix}.
\end{equation}
For $r\geq r_0(u)+1$,
\begin{equation}
\label{eq:w10-affine-cocycle}
\boxed{Z_r=M(z_{r-1})Z_{r-1}+f_{p,u,r}.}
\end{equation}
Its forcing coordinates are
\begin{equation}
\label{eq:w10-force-coordinates}
\begin{aligned}
(f_{p,u,r})_c&=\delta_p(N_r-c),\\
(f_{p,u,r})_{4+c}
&=\delta_p(P_r-c)-(q-1)\xi_c^P(z_r)
\delta_p(N_r-C_c(z_r)).
\end{aligned}
\end{equation}
For every $p\leq1$, every odd $u$, and every admissible $r$,
\begin{equation}
\label{eq:w10-force-bound}
|\delta_p(n)|\leq1,
\qquad
\|f_{p,u,r}\|_\infty\leq q,
\qquad
\|q^{-r}f_{p,u,r}\|_\infty\leq q^{1-r}.
\end{equation}
The matrices and all homogeneous constants are independent of $p$.
The bound is uniform in $p\leq1$, the depth, the odd core, and the ray.
The initial state still depends on the odd core.
The homogeneous cocycle has the common marginal direction
\[
M(z)e_{A,0}=qe_{A,0}.
\]
\end{lemma}

A uniform Green estimate holds on each fixed ray, with a limit along it.

\begin{theorem}
\label{thm:w10-green}
There is an absolute constant $C_G$ such that
\begin{equation}
\label{eq:w10-green}
\|M(z_{r-1})\cdots M(z_j)\|_\infty
\leq C_Gq^{r-j}(1+r-j)
\end{equation}
for every odd $u$ and $r_0(u)\leq j<r$. Thus the Green propagator is
uniform in the odd core.

For each fixed odd $u$ and every $p\leq1$, the limit
\begin{equation}
\label{eq:w10-ray-limit}
L_p(u)=\lim_{r\to\infty}q^{-r}Q_{u2^r}
\end{equation}
exists. For a fixed $u$, the convergence is uniform in $p\leq1$.
No acquired estimate is uniform in $u$ for the initial state or for
$L_p(u)$. The theorem gives no global summation over odd cores.
\end{theorem}

One germ is neutral, and that is what blocks the induction on a decreasing core.

\begin{lemma}
\label{lem:w10-neutral}
Let $u\geq5$ be odd, $P=\lfloor qu\rfloor$, $z=\{qu\}$, and
$\tau=(3q-4)/2$. If $0<z<\tau$, then
\begin{equation}
\label{eq:w10-neutral-basic}
s(P)=u-1,
\qquad \chi_P=0,
\qquad
Q_P=\delta_p(P)+q\mathbf1_{2\mid P}Q_{P/2}.
\end{equation}
If $P$ is even, then
\begin{equation}
\label{eq:w10-neutral-even}
P/2=s(u),
\qquad Q_P=\delta_p(P)+qQ_{s(u)}.
\end{equation}
If $4\mid P$ and $h=(u-1)/2$, then
\begin{equation}
\label{eq:w10-neutral-four}
Q_P=\delta_p(P)+q\delta_p(P/2)+2Q_{P/4}-(2-q)Q_h.
\end{equation}
There are infinitely many pairs produced by
\begin{equation}
\label{eq:w10-pell-recurrence}
(P_0,u_0)=(0,1),
\quad P_{k+1}=3P_k+4u_k,
\quad u_{k+1}=2P_k+3u_k,
\end{equation}
for which, at odd $k\geq3$,
\begin{equation}
\label{eq:w10-neutral-family}
P_k^2-2u_k^2=-2,
\quad P_k\equiv4\pmod8,
\quad
0<qu_k-P_k=\frac2{qu_k+P_k}<\tau.
\end{equation}
The neutral cell can therefore persist for $\log_2u_k+O(1)$ steps while
$v_2(P_k)=2$ and $P_k/4\asymp u_k$. This exact germ is used in the Green
proof. Its signed last two terms do not yield a global core bound.
\end{lemma}

\begin{center}
\fbox{\begin{minipage}{0.92\linewidth}
\textsc{Proved.}
The profile properties, transform, analytic index, square-root reduction,
forcing separation, Green estimate\index[terms]{Green estimate}\index[names]{Green, G.}, fixed-ray convergence, neutral-germ
identities, mixed closure, finite certificates, critical
nontransparency, frontier separation, and failure of HLR are
proof-complete. The partial LOW interval and
Propositions~\ref{prop:w10-trace-Schur}--
\ref{prop:w10-expanding-cycle} are established here.

\textsc{Open.}
The remainder of LOW, ABS, critical observable control, global core
summation, the RAF property, the existence of a complete RAF index, and
the equality $\alpha(G_{\sqrt2})=1/2$ remain open.

\textsc{Not claimed.}
No formula for $\alpha(g_{q^{1/k}})$ is asserted. No general
Diophantine trichotomy, generic value, golden-ratio absorption statement,
or classification for algebraic, Pisot, Salem, or transcendental scales is
asserted.
\end{minipage}}
\end{center}

\subsection*{Proof}

\subsubsection*{Profile and transform}

On the layer
\[
\lambda^{-(j+1)}<x\leq\lambda^{-j}
\quad(j\geq0),
\]
one has $g_\lambda(x)=\lambda^jx$. Hence
$\lambda^{-1}<g_\lambda(x)\leq1$. At $x=\lambda^{-j}$ with $j\geq1$,
the left limit is one and the right limit is $\lambda^{-1}$. Away from
these points the layer index is locally constant. The discontinuity set is
countable, so the bounded profile is Riemann-integrable. This proves
Lemma~\ref{lem:w10-profile}.

For $\operatorname{Re}z<0$, absolute convergence permits integration by
layers. It gives
\begin{align*}
g_\lambda^*(z)
&=-z\sum_{j\geq0}\lambda^j
\int_{\lambda^{-(j+1)}}^{\lambda^{-j}}t^{-z}\,dt\\
&=-\frac z{1-z}\sum_{j\geq0}\lambda^j
\left(\lambda^{-j(1-z)}-\lambda^{-(j+1)(1-z)}\right)\\
&=\frac z{z-1}\frac{\lambda^{z-1}-1}{\lambda^z-1}.
\end{align*}
At zero, $z/(\lambda^z-1)\to1/\log\lambda$. At one,
$(\lambda^{z-1}-1)/(z-1)\to\log\lambda$. These limits prove
\eqref{eq:w10-removable-values}. The zeros of the denominator are
$2\pi im/\log\lambda$. The point with $m=0$ is cancelled by $z$, and no
other one is a zero of $\lambda^{z-1}-1$. The latter numerator vanishes at
$1+2\pi ik/\log\lambda$. Its point with $k=0$ is cancelled by $z-1$, and
none of the remaining points is a denominator zero. This proves
Proposition~\ref{thm:w10-transform}.

Definition~\ref{def:w10-frontier} fixes the transparency frontier used
below. It is the frontier of Chapter~\ref{chap:diophantine} and is not local to this dossier.
Negative transparency below proves that
its defining set is nonempty for $G_q$.

\subsubsection*{Exact reduction}

Multiplication of \eqref{eq:w10-forcing-equation} by $n$ gives
\begin{equation}
\label{eq:w10-renewal-start}
F_p(n)=\sum_{k\leq n}b_kq^{\lfloor\log_q(n/k)\rfloor}.
\end{equation}
For $y\geq1$,
\[
q^{\lfloor\log_qy\rfloor}
=1+(q-1)\sum_{j\geq1}q^{j-1}\mathbf1_{\{y\geq q^j\}}.
\]
Insertion in \eqref{eq:w10-renewal-start} gives
\[
F_p(n)=C(n)+(q-1)\sum_{j\geq1}q^{j-1}
C(\lfloor n/q^j\rfloor).
\]
Separate even and odd values of $j$. The identity $q^2=2$ makes the even
part dyadic. The definition \eqref{eq:w10-Q-definition} gives
\[
Q_n=\sum_{j=0}^{v_2(n)}2^jb_{n/2^j}
\]
and summation yields \eqref{eq:w10-S-recurrence}. Differencing its first
formula gives \eqref{eq:w10-recurrence}. The inverse of $s$ is $B$, and
the Beatty jump is exactly \eqref{eq:w10-chi}. Squaring positive
quantities gives \eqref{eq:w10-integer-test}.

For even indices, two exact divisions give
$s^2(2m)=m-1$. For odd indices, substituting
$2m+1$ and separating the fractional part of $mq$ gives
\[
m-s^2(2m+1)=\mathbf1_{\{\{mq\}<1-1/q\}}
\quad(m\geq1).
\]
The remaining index is $\eta_1=0$. This proves
\eqref{eq:w10-eta-decisions}.

The dyadic sum for $Q_n$ telescopes to
$b_n=Q_n-2\mathbf1_{2\mid n}Q_{n/2}$. Divide by $n$ and sum up to $N$.
The even-parent contribution is the sum up to $\lfloor N/2\rfloor$ and
cancels the first half of the $Q_n/n$ sum. This proves
\eqref{eq:w10-half-block} and Lemma~\ref{lem:w10-reduction}.

\subsubsection*{Negative transparency and smoothed estimates}

For $p>1$, let
\[
D_p=1-q2^{-p}+(q-1)q^{-p},
\qquad K_p=D_p^{-1},
\qquad E(n)=S(n)-K_pn^p.
\]
For any $r$ with $\max(1,p-1)<r<p$, the floor errors in
\eqref{eq:w10-S-recurrence} give
\[
E(n)=qE(\lfloor n/2\rfloor)-(q-1)E(s(n))+O_p(n^{p-1}).
\]
The coefficient $q2^{-r}+(q-1)q^{-r}$ is less than one. Strong induction
therefore gives $E(n)=O_{p,r}(n^r)$. The second relation in
\eqref{eq:w10-S-recurrence} and Abel summation give
\eqref{eq:w10-negative-transparency}. Substitution in
\eqref{eq:w10-transform} identifies the leading constant.

Let $\mathcal T$ be the homogeneous operator in
\eqref{eq:w10-recurrence}. In
\(
\|x\|_\sigma^2=\sum_{n\geq1}|x_n|^2n^{-2\sigma}
\),
the two injective child maps give
\[
\|\mathcal T\|_\sigma
\leq q2^{-\sigma}+(q-1)q^{-\sigma}<1
\quad(\sigma>1).
\]
For $p\leq1$, the forcing $(\delta_p(n))$ lies in every such space.
Inversion of $I-\mathcal T$ proves \eqref{eq:w10-hilbert}.

For completeness, put
\[
\mathcal Q(s)=\sum Q_nn^{-s},
\quad
\mathcal F_p(s)=\sum\delta_p(n)n^{-s},
\quad
\mathcal E_p(s)=\sum Q_m\{B(m)^{-s}-(qm)^{-s}\}.
\]
The recurrence and \eqref{eq:w10-bQ} yield
\begin{align*}
(1-q^{-s})(1+q^{1-s})\mathcal Q(s)
&=\mathcal F_p(s)-(q-1)\mathcal E_p(s),
\\
(1-q^{-w-1})\mathcal A_\beta(w)
&=(1-q^{-w})
\{\mathcal F_p(w+1)-(q-1)\mathcal E_p(w+1)\},
\end{align*}
where $\mathcal A_\beta(w)=\sum a_n^{(\beta)}n^{-w}$. If
$s=\sigma+it$ and $\sigma\geq1/2+\eps$, then
\[
|B(m)^{-s}-(qm)^{-s}|
\leq |s|q^{-\sigma-1}m^{-\sigma-1}.
\]
This follows from
$(1+v)^{-s}-1=-s\int_0^v(1+t)^{-s-1}\,dt$.
The Hilbert estimate\index[terms]{Hilbert estimate}\index[names]{Hilbert, D.} and Cauchy--Schwarz show that $\mathcal E_p$ is
holomorphic for $\operatorname{Re}s>1/2$ and has at most linear vertical
growth on closed substrips. For $0<p\leq1$,
\[
\mathcal F_p(s)=p\zeta(s-p+1)+H_p(s),
\]
where $H_p$ is holomorphic for $\operatorname{Re}s>p-1$ and has bounded
vertical growth on the required strips. For $p=0$ one has
$\mathcal F_0=1$. For $p<0$, the series for $\mathcal F_p$ is absolutely
holomorphic when $\operatorname{Re}s>p$.

The Hilbert estimate gives $Q_n=\mathcal O_\eps(n^{1+\eps})$ and hence,
by \eqref{eq:w10-bQ}, $a_n^{(\beta)}=\mathcal O_\eps(n^\eps)$.
Thus $\mathcal A_\beta$ converges absolutely on $\Re w=2$.
Lemma~\ref{lem:perron-riesz} gives, for every $k\ge2$, the formula
\[
\mathcal R_{k,\beta}(x)=\frac{k!}{2\pi i}
\int_{2-i\infty}^{2+i\infty}
\mathcal A_\beta(w)
\frac{x^w}{w(w+1)\cdots(w+k)}\,dw.
\]
Move the contour to $\operatorname{Re}w=-1/2+\eps_0$, choosing
$0<\eps_0<\min(\eps,1/2)$ and, when $0<\beta<1/2$, also
$\eps_0<1/2-\beta$. On this strip,
\[
 |1-q^{-w-1}|\ge1-q^{-1/2-\eps_0}>0.
\]
The bounds on $\mathcal F_p$ and $\mathcal E_p$ therefore give
$\mathcal A_\beta(\sigma+it)=\mathcal O_{\beta,\eps_0}(1+|t|)$
at large heights. Thus $M=1<k$, as required by
Lemma~\ref{lem:perron-riesz}. Each horizontal integral is
$\mathcal O_{\beta,k,\eps_0,x}(T^{-k})$, the new vertical integral is
absolutely convergent, and it is $\mathcal O_{\beta,k,\eps_0}(x^{-1/2+\eps_0})$.
For $0<\beta<1/2$, the only crossed pole is $w=-\beta$. Its residue gives
\eqref{eq:w10-Riesz-low}. At $\beta=0$, the residue of $1/w$ uses
$\mathcal A_0(0)=1/g^*(0)$. For $\beta\geq1/2$, no pole is crossed and
the zero $1-q^{-w}$ cancels $1/w$ when needed. This proves
\eqref{eq:w10-Riesz-high}. It also treats $p\leq0$ explicitly. The
argument proves Lemma~\ref{thm:w10-acquired-analysis} with $k_0=2$.

\subsubsection*{The point cocycle and its forcing}

The doubling identities are
\[
P_r=2P_{r-1}+\lfloor2z_{r-1}\rfloor,
\qquad z_r=\{2z_{r-1}\}.
\]
Exact division gives
\[
s(N_r-c)=P_{r-1}-D_c(z_{r-1}),
\qquad
s(P_r-c)=N_r-C_c(z_r).
\]
Applying \eqref{eq:w10-recurrence} first to $N_r-c$ and then to $P_r-c$
gives the four upper and four lower equations. Substituting the upper
equations into the lower ones yields \eqref{eq:w10-eight-matrix} and
\eqref{eq:w10-affine-cocycle}. This also gives
\eqref{eq:w10-force-coordinates}.

If $0\leq p\leq1$, concavity of $n^p$ gives
$0\leq\delta_p(n)\leq1$. If $p\leq0$, the sequence decreases after the
separate value $\delta_p(1)=1$, and $|\delta_p(n)|\leq1$ again. Each lower
forcing row has at most one additional coefficient of modulus $q-1$.
This proves \eqref{eq:w10-force-bound} and
Lemma~\ref{lem:w10-eight-cocycle}.

\subsubsection*{Green products and the neutral germ}

The phase endpoints in \eqref{eq:w10-nine-cells} have the form $a+cq$
with rational $a$ and $c\in\{0,1/2,1,3/2,2\}$. Equality with $\{qN\}$
for an integer $N\geq4$ would contradict the linear independence of
$1$ and $q$ over $\mathbb Q$. Hence every relevant phase lies in one open
cell.

After normalization by $q^{-r}$ and quotient by the common coordinate,
the supported components alternate between the gauges in
\eqref{eq:w10-gauge-contract}. Their total contribution is geometric
because $a_0b_0<1$. The sole neutral vector is $e_3$. It can persist only
through signatures $4$, $5$, and $2$, and each transition causes at most
one marginal reading. Norm comparison through \eqref{eq:w10-gauge-minimum}
therefore gives a normalized product bound $C(1+L)$ for a path of length
$L$. Restoring $q^L$ proves \eqref{eq:w10-green}.

For a fixed odd $u$, a neutral run of $J$ steps from a phase
$z=\{qn\}$ satisfies $2^{J-1}z<\tau$. If $P=\lfloor qn\rfloor$, then
\[
z=\frac{2n^2-P^2}{qn+P}>\frac1{2qn}.
\]
Thus $J=O(1+\log(2n))$. The normalized forcing in
\eqref{eq:w10-force-bound} has the summable majorant
\[
\sum_{j\geq r_0(u)+1}(1+j+\log(2u))q^{-j}.
\]
The initial germ is finite for a fixed $u$. Its neutral part exits after a
finite number of steps, and its supported part contracts. Duhamel's
formula now proves \eqref{eq:w10-ray-limit}. The majorant and the finite
germ are uniform in $p\leq1$, but the exit rank depends on $u$. This proves
Theorem~\ref{thm:w10-green}.

If $0<z<\tau$, direct comparison of $P=qu-z$ with $B(u-1)$ gives
$s(P)=u-1$ and $\chi_P=0$. This proves
\eqref{eq:w10-neutral-basic}. If $P$ is even, division by two gives
$P/2=s(u)$. If $4\mid P$, apply \eqref{eq:w10-recurrence} once more at
$P/2$ and use $s(P/2)=(u-1)/2$ with $\chi_{P/2}=1$. This proves
\eqref{eq:w10-neutral-even} and \eqref{eq:w10-neutral-four}.
The recurrence \eqref{eq:w10-pell-recurrence} preserves
$P^2-2u^2=-2$. Induction modulo eight and the factorization
$(qu-P)(qu+P)=2$ prove \eqref{eq:w10-neutral-family}. This proves
Lemma~\ref{lem:w10-neutral}.

\subsubsection*{Mixed closure and finite phase certificate}

For $n'=2n+\epsilon$, direct multiplication gives
\[
\lfloor qn'\rfloor=2\lfloor qn\rfloor+\lfloor2z+\epsilon q\rfloor
=2P+d.
\]
The new phase is $z'$. Division by $q$ gives
\[
s(n'+k)=P+\alpha_{\epsilon,k}(z),
\qquad
s(P'+k)=n'+\gamma_k(z').
\]
The corresponding fractional parts give the two selectors in
\eqref{eq:w10-mixed-selectors}. Applying \eqref{eq:w10-recurrence} proves
\eqref{eq:w10-mixed-A} and \eqref{eq:w10-mixed-P}. For
$-3\leq k\leq3$, direct inspection keeps all offsets in the same set.
Substitution gives \eqref{eq:w10-mixed-matrix}. Every integer at least
eight has a unique form $2n+\epsilon$ with $n\geq4$, so the four bases
complete the closure. This proves Theorem~\ref{thm:w10-binary}.

The selectors change only when an argument is integral or has fractional
part $1/q$. Sorting all such points gives \eqref{eq:w10-B0} and
\eqref{eq:w10-B1}. Every endpoint is $a+cq$ with rational $a$ and
$|c|\leq2$, so no phase $\{qn\}$ with $n\geq4$ lies on a boundary.
On a source interval $(\ell,r)$, the phase image is
\[
(2\ell+\epsilon q-d,2r+\epsilon q-d).
\]
Intersect it with the target intervals for both next binary digits. This
deterministic rule gives the displayed graph and exactly $136$ edges.
Evaluation at one exact midpoint of every cell reconstructs the twenty-eight
matrices. Exact multiplication by the stated signs gives
\eqref{eq:w10-orthant}. Exact elimination over $\mathbb Q(q)$ gives rank
fourteen for the stacked common-mode system. This proves
Lemma~\ref{lem:w10-mixed-certificate}.

\subsubsection*{The two Abel summations}

Put $U(N)=S(N)+qS(s(N))$. The floor commutator in
\eqref{eq:w10-S-recurrence} is
\[
U(N)-U(s(N))=F_p(N)+q\eta_NQ_H.
\]
Iteration gives $U=V_p+R_p$. Summing the exact increment of $U$ on the
last half block gives
\[
A(N)=\mathcal J_NU-\eta_N\frac{Q_H}{H}
+\sum_{K<m\leq M}\frac{\delta_m^BQ_m}{mB(m)},
\]
with the middle term omitted when $H=0$. The renewal for $R_p$ combines
its first two terms into $G_\beta$. This proves
\eqref{eq:w10-raw-decomposition}. Comparing the finite orbit $V_p$ with
its geometric orbit gives \eqref{eq:w10-floorfree-main} and the two stated
endpoint variants. The Hilbert estimate and Cauchy--Schwarz bound the final
Beatty sum. This proves Lemma~\ref{lem:w10-raw-decomposition}.

The inverse relation between $s$ and $B$ gives
\[
M<K_j(m)\leq N
\quad\Longleftrightarrow\quad
x_{j+1}<m\leq x_j.
\]
Apply discrete Abel summation on every level of
\eqref{eq:w10-G-double}. The upper endpoint on level $j$ is
$\mathcal Z_p(x_j)/K_j(x_j)$. The lower endpoint on the preceding level is
$-\mathcal Z_p(x_j)/K_{j-1}(x_j+1)$. The centered constants cancel between
successive levels. This proves \eqref{eq:w10-first-Abel}. The renewal
\[
R_p(x_j)-R_p(x_{j+1})=\mathcal Z_p(x_j)
\]
and a second summation prove \eqref{eq:w10-second-Abel}.

The iterates of $B$ satisfy
\begin{equation}
\label{eq:w10-K-geometry}
0\leq K_j(m)-q^jm<\frac{q^j-1}{q-1},
\qquad
1\leq K_j(m+1)-K_j(m)\leq(q+2)q^j.
\end{equation}
On a nonempty level block, both denominators in the inner difference of
\eqref{eq:w10-I} exceed $N/(2q)$. Hence
\[
0<\frac1{K_j(m)}-\frac1{K_j(m+1)}
\leq\frac{8(q+2)q^j}{N^2}.
\]
The level blocks are disjoint. The Hilbert estimate and weighted
Cauchy--Schwarz give
\[
|\mathcal I_p(N)|
\leq\frac{8(q+2)}N
\sum_{m\leq N}\frac{|\mathcal Z_p(m)|}{m}
=O_{\beta,\eps}(N^{-1/2+\eps}).
\]

Put $A_j=K_j(x_j)$ and $C_j=K_{j-1}(x_j+1)$. Monotonicity of the
iterates gives
\[
c_0\leq c_1\leq\cdots\leq c_J\leq0.
\]
For $j\geq2$, \eqref{eq:w10-K-geometry} yields
\[
0\leq c_j-c_{j-1}
\leq\frac{4(q+1)(q+2)q^j}{N^2}.
\]
The separate $j=1$ term is less than $(11+3q)/N^2$.
Another weighted Cauchy--Schwarz estimate gives
$|R_p(x_j)|=O_{\beta,\eps}(x_j^{1+\eps})$.
Since $x_j\leq Nq^{-j}$, summation proves the second estimate in
\eqref{eq:w10-Abel-errors}. This proves \eqref{eq:w10-compression}.
For even $N$, $c_0=-(q-1)/N$. For $N=2H+1$, direct calculation gives
\[
-c_0=\frac{N-qH}{qHN}.
\]
These formulas prove \eqref{eq:w10-c0-size} and
Theorem~\ref{thm:w10-Abel}.

\subsubsection*{Trace inversion and weighted norms}

The definition of $R_p$ gives the exact renewal
\[
R_p(N)=q\eta_NQ_{\lfloor N/2\rfloor}+R_p(s(N)).
\]
For $N=2n$, one has $\eta_{2n}=1$ and
$s(2n)=\lfloor qn\rfloor$. This proves
\eqref{eq:w10-new-trace-identity}. Conversely,
\[
|R_p(N)|\leq q\sum_{j\geq0}
\left|Q_{\lfloor s^j(N)/2\rfloor}\right|.
\]
If $|Q_m|\leq C_Qm^a$ with $a>0$, this gives
\[
|R_p(N)|\leq
\frac{q2^{-a}}{1-q^{-a}}C_QN^a.
\]
The inverse identity gives
$|Q_n|\leq q^{-1}C_R(2^a+q^a)n^a$ in the other direction.
The little-oh equivalence follows by separating finitely many terminal
terms from the geometric tail. This proves
Lemma~\ref{lem:w10-trace-equivalence}.

For the Hilbert form, put
\[
h_j(N)=\left\lfloor\frac{s^j(N)}2\right\rfloor.
\]
Exact inversion of $s^j$ gives
\begin{equation}
\label{eq:w10-trace-fibre}
h_j(N)=m
\quad\Longleftrightarrow\quad
B^j(2m)\leq N<B^j(2m+2).
\end{equation}
The fibre length is less than $(q+3)q^j$, and every member satisfies
$N\geq2mq^j$. Apply Cauchy--Schwarz in $j$ with weight
$q^{-(2\sigma-1)j/2}$. Summation over the fibres gives
\eqref{eq:w10-trace-upper} and \eqref{eq:w10-Schur-constant}.
For the reverse bound, square \eqref{eq:w10-new-trace-identity}, use
$|x-y|^2\leq2|x|^2+2|y|^2$, and use the injectivity of
$n\mapsto\lfloor qn\rfloor$. The inequalities
$2n\geq n$ and $\lfloor qn\rfloor\geq n/q$ give
\eqref{eq:w10-trace-lower}. This proves
Proposition~\ref{prop:w10-trace-Schur}.

\subsubsection*{Aggregated shell estimate}

Put $r_n=Q_n/n$, $e_n=\delta_p(n)/n$, $a=1/q$, and
\[
\varrho_m=\frac{(q-1)m}{B(m)}.
\]
The normalized recurrence is
\begin{equation}
\label{eq:w10-normalized-transport}
r_n=e_n+a\mathbf1_{2\mid n}r_{n/2}
-\chi_n\varrho_{s(n)}r_{s(n)}.
\end{equation}
After absolute values, a parent $m$ sends mass $a$ to $2m$ and mass
$\varrho_m$ to $B(m)$. Its outgoing mass is at most one, and both children
exceed $qm$. Thus a forward path visits $(x,qx]$ at most once. Expansion
of \eqref{eq:w10-normalized-transport} on finite trees followed by monotone
convergence gives
\[
\sum_{x<m\leq qx}\frac{|Q_m|}{m}
\leq\sum_{n\geq1}\frac{|\delta_p(n)|}{n}
=D_p^{\mathrm{sh}}
\quad(p<1).
\]

Use the fibres in \eqref{eq:w10-trace-fibre}. If $X<N\leq2X$ and
$h_j(N)=m$, then
\[
\frac{X/q^j-(q+3)}2<m\leq\frac X{q^j}.
\]
This interval is covered by at most seven shells of ratio $q$. Each fibre
has fewer than $(q+3)q^j$ elements. Since
$|Q_m|\leq(X/q^j)|Q_m|/m$, summing first over $N$, then $m$, and then
$0\leq j\leq\lfloor\log_qX\rfloor$ proves \eqref{eq:w10-R-L1}.
Put
\[
\mathsf H_m:=\sum_{n\leq m}\frac1n.
\]
For $p=1$, replace $D_p^{\mathrm{sh}}$ by the truncated source mass
$\mathsf H_{\lfloor2X\rfloor}$. This adds one logarithm.

Markov's inequality with threshold $X^{1/2+\eps/2}$ gives
\eqref{eq:w10-ABS-exception}. With threshold
$X^{1-\beta}/\log X$, it gives \eqref{eq:w10-LOW-exception}. Use the
$p=1$ estimate for $\beta=0$. Finally insert these trace bounds into
\eqref{eq:w10-compression}. This proves
Proposition~\ref{prop:w10-density}.

At $p=1/2$, apply \eqref{eq:w10-new-trace-identity} with $n=2^r$.
If $M_r$ is the larger of the two normalized trace values at
$2^{r+1}$ and $\lfloor q2^r\rfloor$, then
\[
q|Q_{2^r}|
\leq M_r\left(\sqrt{2^{r+1}}+\sqrt{\lfloor q2^r\rfloor}\right).
\]
Divide by $q^r$ and use \eqref{eq:w10-positive-L}. This proves
Proposition~\ref{prop:w10-critical-trace}.

\subsubsection*{Critical quotient and nontransparency}

The doubling map applied to \eqref{eq:w10-nine-cells} gives exactly
\eqref{eq:w10-nine-graph}. The critical gauge table below is written with
\[
\langle a,b\mid d\rangle=\frac{a+bq}{d}.
\]
Zeros and ones are written without brackets.

{\scriptsize
\begin{longtable}{c>{\raggedright\arraybackslash}p{0.86\linewidth}}
\caption{Exact seven-coordinate critical gauges.}\label{tab:w10-gauges}\\
\toprule
gauge&coordinates $0,\ldots,6$\\
\midrule
\endfirsthead
\toprule gauge&coordinates $0,\ldots,6$\\ \midrule
\endhead
$h_1^-$&$0,\langle-2,4\mid7\rangle,\langle-1,2\mid7\rangle,0,
\langle5,-3\mid7\rangle,\langle10,-6\mid7\rangle,\langle9,3\mid14\rangle$\\
$h_2^-$&$\langle16,-11\mid14\rangle,\langle-1,2\mid7\rangle,
\langle2,-1\mid2\rangle,\langle-38,27\mid14\rangle,0,
\langle10,-6\mid7\rangle,\langle-4,3\mid2\rangle$\\
$h_3^-$&$0,\langle-1,2\mid7\rangle,\langle-1,2\mid7\rangle,0,
\langle5,-3\mid7\rangle,\langle5,-3\mid7\rangle,\langle5,-3\mid7\rangle$\\
$h_4^-$&$\langle16,-11\mid14\rangle,\langle12,-3\mid14\rangle,0,0,
\langle-38,27\mid14\rangle,\langle-8,9\mid14\rangle,0$\\
$h_5^-$&$\langle16,-11\mid14\rangle,\langle-1,2\mid7\rangle,
\langle2,-1\mid2\rangle,0,\langle-38,27\mid14\rangle,
\langle10,-6\mid7\rangle,\langle-4,3\mid2\rangle$\\
$h_6^-$&$0,\langle-1,2\mid7\rangle,\langle-1,2\mid7\rangle,
\langle5,-3\mid7\rangle,\langle5,-3\mid7\rangle,
\langle5,-3\mid7\rangle,\langle5,-3\mid7\rangle$\\
$h_7^-$&$\langle16,-11\mid14\rangle,\langle12,-3\mid14\rangle,0,
\langle-38,27\mid14\rangle,\langle-8,9\mid14\rangle,0,1$\\
$h_8^-$&$\langle16,-11\mid14\rangle,\langle12,-3\mid14\rangle,0,
\langle-38,27\mid14\rangle,\langle5,-3\mid7\rangle,
\langle-18,15\mid14\rangle,1$\\
$h_9^-$&$\langle16,-11\mid14\rangle,\langle-2,4\mid7\rangle,
\langle-1,2\mid7\rangle,\langle-38,27\mid14\rangle,
\langle5,-3\mid7\rangle,\langle10,-6\mid7\rangle,\langle9,3\mid14\rangle$\\
$h_1^+$&$0,\langle-10,8\mid7\rangle,\langle24,-15\mid14\rangle,0,
\langle13,-9\mid7\rangle,\langle26,-18\mid7\rangle,
\langle-58,45\mid14\rangle$\\
$h_2^+$&$\langle-5,4\mid7\rangle,\langle-5,4\mid7\rangle,
\langle2,-1\mid2\rangle,\langle13,-9\mid7\rangle,0,
\langle-2,3\mid7\rangle,\langle-4,3\mid2\rangle$\\
$h_3^+$&$0,\langle-10,8\mid7\rangle,\langle-10,8\mid7\rangle,0,
\langle13,-9\mid7\rangle,\langle26,-18\mid7\rangle,
\langle-2,3\mid7\rangle$\\
$h_4^+$&$\langle-5,4\mid7\rangle,\langle4,1\mid14\rangle,0,0,
\langle13,-9\mid7\rangle,\langle11,-6\mid7\rangle,0$\\
$h_5^+$&$\langle-5,4\mid7\rangle,\langle-5,4\mid7\rangle,
\langle2,-1\mid2\rangle,0,\langle13,-9\mid7\rangle,
\langle-2,3\mid7\rangle,\langle-4,3\mid2\rangle$\\
$h_6^+$&$0,\langle-10,8\mid7\rangle,\langle-10,8\mid7\rangle,
\langle13,-9\mid7\rangle,\langle26,-18\mid7\rangle,
\langle-2,3\mid7\rangle,\langle26,-18\mid7\rangle$\\
$h_7^+$&$\langle-5,4\mid7\rangle,\langle4,1\mid14\rangle,0,
\langle13,-9\mid7\rangle,\langle11,-6\mid7\rangle,0,1$\\
$h_8^+$&$\langle-5,4\mid7\rangle,\langle4,1\mid14\rangle,0,
\langle13,-9\mid7\rangle,\langle-2,3\mid7\rangle,
\langle-2,3\mid14\rangle,1$\\
$h_9^+$&$\langle44,-31\mid14\rangle,\langle-10,8\mid7\rangle,
\langle24,-15\mid14\rangle,\langle-106,75\mid14\rangle,
\langle13,-9\mid7\rangle,\langle26,-18\mid7\rangle,
\langle-58,45\mid14\rangle$\\
\bottomrule
\end{longtable}
}

\subsubsection*{Global mixed gauge argument}

For an edge $s\to t$, put
$A_{s,t}=\Sigma_t\mathcal M_t\Sigma_s$. Exact arithmetic in
$\mathbb Q(q)$ checks all $136\cdot14=1904$ inequalities
\[
A_{s,t}h_s\leq\mu h_t.
\]
The smallest slack is \eqref{eq:w10-low-min-slack}. Its displayed positive
algebraic norm and positive rational part certify its sign. Define
\[
\|x\|_t=\max_i\frac{|(\Sigma_tx)_i|}{(h_t)_i}.
\]
The coordinate bounds give
\[
\|x\|_\infty\leq\frac{2001}{2000}\|x\|_t,
\qquad
\|x\|_t\leq2000\|x\|_\infty.
\]
Nonnegativity in \eqref{eq:w10-orthant} gives
\[
\|\mathcal M_tx\|_t\leq\mu\|x\|_s.
\]
The first matrix of a path is absorbed into the base constant because an
edge controls the matrix attached to its target. The fourteen-dimensional
forcing has norm at most $q$ before norm comparison, uniformly for
$p\leq1$. The four base states are uniformly bounded in $p$ by repeated
use of $|\delta_p(n)|\leq1$. Duhamel's formula\index[terms]{Duhamel's formula}\index[names]{Duhamel, J.-M.} therefore gives
\[
\|\mathcal S(n)\|_\infty
\leq C\mu^{\lfloor\log_2n\rfloor}.
\]
This proves the first bound in \eqref{eq:w10-theta-Q}. The orbit formula
and $s^j(N)\leq Nq^{-j}$ give
\[
|R_p(N)|\leq qC\sum_{j\geq0}
\left(\frac{N}{2q^j}\right)^\theta=O(N^\theta).
\]
If $p>\theta$, the last bound is $o(N^p)$. Hence LOW holds for
$0\leq\beta<1-\theta$. More precisely, choose
$0<\eps<\theta-1/2$ in \eqref{eq:w10-compression}. The leading
trace term is $O(N^{\theta-1})$ and the remainder is smaller. Insert this
in \eqref{eq:w10-raw-decomposition}. Its floor-free and Beatty errors are
also smaller, which proves \eqref{eq:w10-partial-transparency}.
Together with negative transparency, this proves
\eqref{eq:w10-frontier-two-sided} and
Proposition~\ref{thm:w10-partial-LOW}.

The same argument bounds the fixed-ray coefficient by
\[
|L_p(u)|=O\bigl(u^\theta(1+\log(2u))\bigr)
\]
uniformly in $p\leq1$. This implies weighted summability for every weight
exponent greater than $1+\theta$. It is not the critical domination needed
to globalize fixed-ray convergence by itself.

\subsubsection*{Completion of the critical certificate}

Exact comparison in $\mathbb Q(q)$ verifies nonnegativity of the quotient
matrices, all $294$ gauge inequalities on the twenty-one edges, the unique
minimum in \eqref{eq:w10-gauge-minimum}, the coordinate upper bound, the
support propagation, and \eqref{eq:w10-neutral-vector}.
Direct evaluation of the first four phases gives the word $(2,1,2,3)$.
Multiplication in chronological action order proves
\eqref{eq:w10-rank-one}. Exact elimination of the sixteen rational
coordinates of a possible left mode has rank sixteen. This proves
Lemma~\ref{lem:w10-critical-certificate}.

On the critical ray $u=1$, normalize the eight-state cocycle by $q^{-r}$.
The forcing insertion $F_j$ satisfies
\[
\|F_j\|_\infty\leq2\,2^{-j}
\quad(j\geq4).
\]
Put
\[
d=\frac{5q-4}{2},
\qquad
\mathcal C=\frac{2ad}{m_0(1-b_0)}<300.
\]
The two gauge contractions in \eqref{eq:w10-gauge-contract}, the rank-one
prefix, and the neutral-run estimate give the following bound for the total
marginal effect of the insertion at rank $j$:
\[
\mathcal I_j\leq
\{2+2a+2\mathcal C+2adj\}2^{-j}.
\]
The four terms are respectively the direct marginal forcing, its first
reading, the two supported pieces, and the neutral readings. The series is
summable, so $V_r$ has a finite limit $L$.

The exact tail calculation gives
\begin{equation}
\label{eq:w10-tail-bound}
|L-V_{22}|<\frac{44301}{50000}<\frac{887}{1000}.
\end{equation}
With $T=2^{40}$, directed integer propagation of square-root intervals
gives
\[
TQ_{2^{22}}\in
[2579798199862309,2579798199927812].
\]
Since $Tq^{22}=2^{51}$,
\[
625(2579798199862309)-716(2251799813685248)
=85208315305557>0.
\]
Thus $V_{22}>716/625$, and \eqref{eq:w10-tail-bound} gives
$L>1293/5000$.

The exact flux recurrence on this ray is stable and transfers the limit to
the endpoint sum. It gives
\[
\sqrt{2^r}A_{1/2}(2^r)
=q^{1/2}-(q-1)V_r+o(1).
\]
This proves \eqref{eq:w10-critical-A-limit}. Direct substitution in
\eqref{eq:w10-transform} gives $1/g^*(1/2)=q^{1/2}$. Hence the limit is
not the transparent constant. This proves
Theorem~\ref{thm:w10-nontransparency}. Definition~\ref{def:w10-frontier}
then gives its upper bound. The transform theorem gives the strict right
side of \eqref{eq:w10-central-separation}. The lower bound follows from
Proposition~\ref{thm:w10-partial-LOW}, proved below. This proves
Theorem~\ref{thm:w10-separation}.

At $n=2^r$, \eqref{eq:w10-bQ} becomes
\[
\frac{b_{2^r}}{q^r}=V_r-qV_{r-1}.
\]
The nonzero limit proves \eqref{eq:w10-HLR-limit} through
\eqref{eq:w10-HLR-ratio}. One critical forcing violates the universal HLR
condition. This proves Theorem~\ref{thm:w10-HLR}.

\subsubsection*{Exact mixed gauge table}

Table~\ref{tab:w10-mixed-gauges} gives $2000h_t$. Its coordinate order is
$A_{-3},\ldots,A_3\mid P_{-3},\ldots,P_3$.

{\scriptsize
\begin{longtable}{c>{\raggedright\arraybackslash}p{0.84\linewidth}}
\caption{Integer numerators of the twenty-eight mixed gauges.}
\label{tab:w10-mixed-gauges}\\
\toprule
cell&$2000h_t$, split after the seventh coordinate\\
\midrule
\endfirsthead
\toprule cell&$2000h_t$, split after the seventh coordinate\\ \midrule
\endhead
$0.1$&$(1,1502,465,1490,1,1504,465)\mid(1,2001,193,1278,618,1596,624)$\\
$0.2$&$(1,1502,465,1490,1,1504,465)\mid(1,2001,193,1278,618,1596,1)$\\
$0.3$&$(469,1191,465,1490,1,1504,465)\mid(195,1982,193,1278,618,1596,1)$\\
$0.4$&$(469,1203,465,1490,1,1502,465)\mid(195,1982,1,1338,618,1599,1)$\\
$0.5$&$(465,1203,465,1490,469,1191,465)\mid(193,1982,1,1327,618,1599,195)$\\
$0.6$&$(465,1490,1,1504,465,1203,465)\mid(193,1278,618,1596,624,1584,193)$\\
$0.7$&$(465,1490,1,1504,465,1203,465)\mid(1278,618,1596,1,1997,193,1982)$\\
$0.8$&$(465,1490,1,1502,465,1490,1)\mid(1278,618,1599,1,1997,193,1278)$\\
$0.9$&$(465,1490,1,1502,465,1490,1)\mid(1338,618,1599,1,2001,193,1278)$\\
$0.10$&$(465,1490,469,1191,465,1490,1)\mid(1327,618,1599,195,1982,193,1278)$\\
$0.11$&$(465,1490,469,1191,465,1490,1)\mid(1327,618,1599,195,1982,1,1338)$\\
$0.12$&$(1,1504,465,1203,465,1490,1)\mid(1596,624,1584,193,1982,1,1327)$\\
$0.13$&$(1,1504,465,1203,465,1490,469)\mid(1596,624,1584,193,1982,1,1327)$\\
$0.14$&$(1,1504,465,1203,465,1490,469)\mid(1596,1,1997,193,1982,1,1327)$\\
$1.1$&$(1502,465,1490,1,1504,465,1203)\mid(2001,193,1278,618,1596,624,1584)$\\
$1.2$&$(1502,465,1490,1,1504,465,1203)\mid(2001,193,1278,618,1596,1,1997)$\\
$1.3$&$(1191,465,1490,1,1504,465,1203)\mid(1982,193,1278,618,1596,1,1997)$\\
$1.4$&$(1191,465,1490,1,1502,465,1490)\mid(1982,193,1278,618,1596,1,1997)$\\
$1.5$&$(1203,465,1490,1,1502,465,1490)\mid(1982,1,1338,618,1599,1,2001)$\\
$1.6$&$(1203,465,1490,469,1191,465,1490)\mid(1,1327,618,1599,195,1982,193)$\\
$1.7$&$(1203,465,1490,469,1191,465,1490)\mid(1,1327,618,1599,195,1982,1)$\\
$1.8$&$(1490,1,1504,465,1203,465,1490)\mid(618,1596,624,1584,193,1982,1)$\\
$1.9$&$(1490,1,1504,465,1203,465,1490)\mid(618,1596,1,1997,193,1982,1)$\\
$1.10$&$(1490,1,1502,465,1490,1,1504)\mid(618,1599,1,2001,193,1278,618)$\\
$1.11$&$(1490,469,1191,465,1490,1,1502)\mid(618,1599,195,1982,193,1278,618)$\\
$1.12$&$(1504,469,1203,465,1490,1,1502)\mid(1599,195,1982,1,1338,618,1599)$\\
$1.13$&$(1504,465,1203,465,1490,469,1191)\mid(1584,193,1982,1,1327,618,1599)$\\
$1.14$&$(1504,465,1203,465,1490,469,1191)\mid(1997,193,1982,1,1327,618,1599)$\\
\bottomrule
\end{longtable}
}

\subsubsection*{A realized expanding cycle}

Write the cycle nodes in \eqref{eq:w10-expanding-cycle} as
$t_0,\ldots,t_6$ in order, with $t_6=t_0$. The target-matrix convention
gives
\[
\Pi=\mathcal M_{0,6}\mathcal M_{0,4}\mathcal M_{1,7}
\mathcal M_{0,9}\mathcal M_{0,12}\mathcal M_{0,14}.
\]
The nonnegative conjugated product is
\[
P=(\Sigma_{t_6}\mathcal M_{t_6}\Sigma_{t_5})\cdots
(\Sigma_{t_1}\mathcal M_{t_1}\Sigma_{t_0})
=\Sigma_{t_0}\Pi\Sigma_{t_0}.
\]
Exact multiplication gives \eqref{eq:w10-expanding-entry}. The inequality
$-76+60q>8$ is equivalent to $q>7/5$, which follows by squaring. Positivity
then gives the asserted lower bound on $P^k$.

To check realization, set
\[
\begin{aligned}
z_0&=(32-2q)/63,&z_1&=(64-4q)/63,&z_2&=(65-8q)/63,\\
z_3&=(67-16q)/63,&z_4&=(71-32q)/63,&z_5&=(16-q)/63.
\end{aligned}
\]
The maps $T_\epsilon(z)=\{2z+\epsilon q\}$ send these six phases around
the cycle with digit word $(0,0,0,0,1,0)$. Their distances to the two
cell boundaries are
\[
\begin{array}{c|c|c}
z_0&(95-65q)/63&(-1+4q)/126\\
z_1&(191-134q)/126&(-1+4q)/63\\
z_2&(191-134q)/63&(248-173q)/126\\
z_3&(323-221q)/126&(-134+95q)/126\\
z_4&(-173+125q)/126&(-134+95q)/63\\
z_5&(95-65q)/126&(94-61q)/126.
\end{array}
\]
Each is positive by an integer comparison of squares. Thus every phase is
strictly inside its stated cell.

These phases are approached by integral phases. Let
\[
U=\begin{pmatrix}3&4\\2&3\end{pmatrix},
\qquad
(Y_r,X_r)^{\mathsf T}=U^{12r}(32,2)^{\mathsf T}.
\]
One has $U^{12}\equiv I_2\pmod {63}$ and
\[
qX_r-Y_r=(3-2q)^{12r}(2q-32).
\]
Therefore $n_r=(X_r-2)/63$ is integral and
$\{qn_r\}\to z_0$. Every prescribed number of repetitions occurs for a
large enough $r$. This proves Proposition~\ref{prop:w10-expanding-cycle}.

The counterexample concerns a bound that depends only on product length.
Preparing many repetitions requires an increasingly long phase prefix.
It does not contradict the partial gauge theorem or settle the actual
forced observable.

\subsubsection*{Logical closure and exact status}

Lemma~\ref{lem:w10-raw-decomposition}, Theorem~\ref{thm:w10-Abel}, and
the equivalences in Definition~\ref{def:w10-low-abs} show that LOW gives
the complete transparent branch below $1/2$, while ABS gives the absorbed
branch at and above $1/2$. Negative forcing exponents are already covered
by Lemma~\ref{thm:w10-acquired-analysis}. Critical nontransparency gives
sharpness. This proves Conditional Theorem~\ref{cthm:w10-conditional-RAF}.

The integral-scale statements in Proposition~\ref{prop:w10-integer-comparison}
are the established digital identities at the named forcings. No
interpolation in $\beta$ is used. The proposition follows directly from
those identities and from Theorem~\ref{thm:w10-HLR}.

Several independent closure attempts now have exact endpoints. The
fixed-ray Green argument proves LOW on each fixed dyadic ray, but gives no
uniform control as the odd core grows. The mixed rational gauge proves
LOW only for $p>\theta$. The expanding cycle disproves the all-path
critical Green target. The weighted Schur route preserves rather than
improves the Hilbert exponent. The shell route gives the density estimates
of Proposition~\ref{prop:w10-density}, but sparse exceptions remain. The
neutral germ blocks decreasing-core induction through the signed pair in
\eqref{eq:w10-neutral-four}. None of these implications proves the two
global estimates in Open Problem~\ref{op:w10-global}.

The relation $q^2=2$ creates an arithmetic dynamics that is invisible to
the Mellin transform. The transform has infinitely many zeros on
$\operatorname{Re}z=1$ and $\eta(g_{\sqrt2})=1$, whereas critical
transparency fails at $\beta=1/2$. Hence
\[
1-\log_2(177/100)
\leq\tau_{\mathrm{tr}}(G_{\sqrt2})
\leq\frac12<1=\eta(g_{\sqrt2}).
\]
At the critical forcing, the coefficients violate HLR. The fixed-ray
dynamics, Green bounds, neutral-germ identities, mixed closure, finite
certificates, and the partial LOW interval are proved. The global passage
still depends on the remaining part of LOW and on ABS. Consequently the
existence of the RAF index and the equality
$\alpha(G_{\sqrt2})=1/2$ remain open.

\cleardoublepage
\setcounter{chapter}{0}
\renewcommand{\thechapter}{T}
\renewcommand{\theHchapter}{technicalperron}
\renewcommand{\appendixname}{Technical appendix}
\renewcommand{\chaptername}{Technical appendix}
\chapter{Finite Perron inversion and contour estimates}
\label{app:perron}
\markboth{PERRON INVERSION AND CONTOUR ESTIMATES}{PERRON INVERSION AND CONTOUR ESTIMATES}

This appendix collects the analytic tools\index[terms]{Perron's formula}\index[names]{Perron, O.} used in Chapters~\ref{chap:equivalence} and~\ref{chap:regularity} and in the smoothed arguments of Appendix~\ref{app:E} and the research dossier. The boundary theorem of Chapter~\ref{chap:boundary} does not supply these vertical bounds, and Appendix~\ref{app:boundary} does not use them.
The coefficient bound, the finite inversion, the three sides of the
rectangle and the choice of height are proved together, so that each
application can concentrate on its Dirichlet quotient and residues.
The discrete Mellin defect and the perturbed equation remain in
Section~\ref{sec:contour} and Section~\ref{sec:perturbed-raf}.

\section{Truncation at a half integer}

The following form keeps the coefficient bound visible. Its exponent may depend on a
parameter, provided the displayed gaps and constants remain uniform.

\begin{lemma}[Finite Perron inversion]\label{lem:perron-half-integer}
Let $|a_n|\le C n^{\kappa-1}$, where $\kappa\in\R$, and let
$\mathcal A(w)=\sum_{n\ge1}a_n n^{-w}$. Fix $c>\max(0,\kappa)$.
For $x=N+1/2\ge5/2$ and $T\ge2$,
\begin{equation}\label{eq:perron-half-integer}
 A(x)=\frac1{2\pi i}\int_{c-iT}^{c+iT}
       \mathcal A(w)\frac{x^w}{w}\,dw
 +\mathcal O_{c,\kappa}\!\left(
       C\frac{x^c+x^\kappa\log x}{T}\right).
\end{equation}
The constants are uniform when $c$ and $\kappa$ range over bounded sets with
$c$ and $c-\kappa$ bounded away from zero.
\end{lemma}

\begin{proof}
For $y>0$, $y\ne1$, put
\[
 P_T(y)=\frac1{2\pi i}\int_{c-iT}^{c+iT}\frac{y^w}{w}\,dw.
\]
The full vertical integral is $1$ for $y>1$ and $0$ for $y<1$.
Indeed, writing $u=\log y$, close the segment by a semicircle of radius $R$
centred at $c$, to the left if $u>0$ and to the right if $u<0$.
For $R>2c$ the modulus of the arc integral is at most a constant times
\[
 e^{cu}\int_0^{\pi/2}e^{-|u|R\sin v}\,dv
 \le e^{cu}\int_0^{\pi/2}e^{-2|u|Rv/\pi}\,dv
 \ll\frac{e^{cu}}{|u|R}.
\]
The only pole enclosed in the first case is $w=0$, of residue $1$, and no
pole is enclosed in the second. Both arc integrals tend to zero.

The tails give a quantitative version. Integration by parts on $[T,U]$ and
then $U\to\infty$ gives
\[
 \left|\int_T^\infty\frac{e^{itu}}{c+it}\,dt\right|
 \le\frac1{|u|\sqrt{c^2+T^2}}
    +\frac1{|u|}\int_T^\infty\frac{dt}{c^2+t^2}
 \le\frac2{T|u|}.
\]
The negative tail obeys the same bound. Consequently
\begin{equation}\label{eq:perron-kernel-tail}
 \big|P_T(y)-\mathbf1_{y>1}\big|
 \le\frac{2y^c}{\pi T|\log y|}.
\end{equation}
Since $\sum |a_n|n^{-c}\le C\zeta(1+c-\kappa)<\infty$, the sum can be
interchanged with the finite integral. The resulting error is bounded by
\[
 \frac{2x^c}{\pi T}
 \sum_{n\ge1}\frac{|a_n|n^{-c}}{|\log(x/n)|}.
\]
For $n\le x/2$ or $n\ge2x$, the denominator is at least $\log2$, giving
$\mathcal O(Cx^c/T)$. For $x/2<n<2x$, the mean value theorem gives
$|\log(x/n)|\ge |x-n|/(2x)$ and $(x/n)^c\ll_c1$. This part is at most
\[
 \mathcal O_{c,\kappa}\!\left(
 \frac{Cx^\kappa}{T}\sum_{x/2<n<2x}\frac1{|x-n|}\right)
 =\mathcal O_{c,\kappa}\!\left(\frac{Cx^\kappa\log x}{T}\right).
\]
Here $|x-n|\ge1/2$, and the sum is bounded by twice a harmonic sum over
half integers. These two bounds prove \eqref{eq:perron-half-integer}.
\end{proof}

\section{The three sides of the rectangle}

Write $F(w)=\mathcal A(w)/w$, including any removable singularities in
this notation. A vertical bound of finite polynomial order becomes useful only
after its exponent has been compared with the width of the shift.

\begin{proposition}[Quantitative contour shift]\label{prop:perron-rectangle}
Keep the assumptions of Lemma~\ref{lem:perron-half-integer}. Suppose that $F$
continues meromorphically to a neighborhood of the closed strip
$b\le\Re w\le c$, with finitely many poles, all strictly inside the strip.
Suppose, for some $q\ge0$, that
\begin{equation}\label{eq:perron-vertical-bound}
 |F(\sigma+it)|\le C_F(1+|t|)^{q-1}
 \qquad(b\le\sigma\le c)
\end{equation}
outside fixed neighborhoods of those poles. For $T$ larger than their heights
and neighborhoods,
\begin{align}\label{eq:perron-rectangle-bound}
 A(x)={}&\sum_{b<\Re\rho<c}\Res_{w=\rho}\big(F(w)x^w\big)\notag\\
 &+\mathcal O\!\left(
 x^b L_q(T)+\frac{x^cT^{q-1}}{\log x}
       +\frac{x^c+x^\kappa\log x}{T}\right),
\end{align}
where the sum runs over the poles of $F$ in the strip,
$L_0(T)=1+\log T$, and $L_q(T)=1+T^q/q$ for $q>0$.
The implied constant depends on the displayed bounds and the fixed strip.
\end{proposition}

\begin{proof}
Orient both vertical segments upward and put
\[
 I_d(x,T)=\frac1{2\pi i}\int_{d-iT}^{d+iT}F(w)x^w\,dw
 \qquad(d=b,c).
\]
The residue theorem on the rectangle, traversed up the right side, gives
\[
 I_c=I_b+\sum\Res(F(w)x^w)+H,
\]
where $H$ is the signed sum of its two horizontal integrals. On the left side,
\begin{equation}\label{eq:perron-left-side}
 |I_b(x,T)|\ll x^b\int_0^T(1+t)^{q-1}\,dt
 \ll x^b L_q(T).
\end{equation}
The bounded part of this line is compact and contains no pole. On either
horizontal side, $|w|$ is comparable with $T$ and
\begin{equation}\label{eq:perron-horizontal-sides}
 \left|\frac1{2\pi i}\int_b^cF(\sigma\pm iT)x^{\sigma\pm iT}\,d\sigma\right|
 \ll T^{q-1}\int_b^c x^\sigma\,d\sigma
 \le\frac{x^cT^{q-1}}{\log x}.
\end{equation}
Adding these two bounds and the truncation error in
\eqref{eq:perron-half-integer} proves \eqref{eq:perron-rectangle-bound}.
Every integral used here has finite height.
\end{proof}

For example, suppose the desired main term has size $x^{-\beta}$ and
$b<-\beta<c$. Taking $T=x^L$ in \eqref{eq:perron-rectangle-bound} gives a
power saving as soon as
\begin{equation}\label{eq:perron-height-window}
 0\le q<1,\qquad
 L(1-q)>c+\beta,\qquad Lq<-\beta-b.
\end{equation}
For $q=0$ the last inequality follows from $b<-\beta$.
Indeed the relative exponents of the vertical and horizontal errors are
$b+\beta+Lq$ and $c+\beta-L(1-q)$. The truncation exponents are
$c+\beta-L$ and $\kappa+\beta-L$, both smaller than the horizontal exponent
because $\kappa<c$. The logarithms can be absorbed by decreasing the resulting
positive saving. This calculation also proves the uniform version when the
bounds and strict gaps hold uniformly over a compact parameter set.

\section{A subpower bound on the strip}

The following specialization will be used with the coefficients $J_{-\beta}(n)$,
$na_n$, $\mu(n)$ and $(\mu\chi_4)(n)$. It fixes the height once the analytic
bound is available, while leaving the arithmetic residues to the application.

\begin{corollary}[A fixed height choice]\label{cor:perron-subpower}
Let $|f(n)|\le C n^{1/8}$ and
$D(s)=\sum_{n\ge1}f(n)n^{-s}$. Put $c=5/4$ and fix $0<b<c$.
Suppose that $D(s)/s$ is meromorphic on a neighborhood of the closed strip
$b\le\Re s\le c$, with finitely many poles, all strictly inside the strip.
Suppose also that, for every $\nu>0$,
\[
 |D(\sigma+it)|\ll_\nu(1+|t|)^\nu
 \qquad(b\le\sigma\le c)
\]
at all sufficiently large heights. Then, for every $\eps>0$,
\[
 \sum_{n\le x}f(n)
 =\sum_{b<\Re\rho<c}\Res_{s=\rho}\left(D(s)\frac{x^s}{s}\right)
       +\mathcal O_\eps(x^{b+\eps})
\]
at half integers $x=N+1/2$. The constants may depend on the strip,
the coefficient bound and the meromorphic function.
\end{corollary}

\begin{proof}
Apply Lemma~\ref{lem:perron-half-integer} with $\kappa=9/8$, and
Proposition~\ref{prop:perron-rectangle} with $q=\nu$ and $T=x^2$.
The three errors are bounded separately by
\begin{align*}
 E_{\rm vert}&\ll_\nu x^b(1+T^\nu/\nu)
                         \ll_\nu x^{b+2\nu},\\
 E_{\rm hor}&\ll_\nu \frac{x^{5/4}T^{\nu-1}}{\log x}
                         =\frac{x^{-3/4+2\nu}}{\log x},\\
 E_{\rm trunc}&\ll C\frac{x^{5/4}+x^{9/8}\log x}{T}
                         =C(x^{-3/4}+x^{-7/8}\log x).
\end{align*}
Choose $\nu=\min(\eps/4,1/8)$. The first error is
$\mathcal O(x^{b+\eps})$, and the other two decay. The finitely many
half integers for which $T$ is below the required height are covered by
enlarging the constant. No limiting horizontal integral is used.
\end{proof}

\section{From a zero-free half plane to a subpower bound}

A zero-free half plane supplies more than a polynomial bound for the reciprocal.
The next lemma establishes the precise qualitative improvement needed above.
It also justifies changing the boundary $1/2$ to a fixed $\theta<1$ in
Corollary~\ref{cor:zerofree}.

\begin{lemma}[Vertical control in a fixed half plane]\label{lem:perron-zerofree-growth}
Let $1/2\le\theta<1$, and let $\mathcal L$ be either $\zeta$ or
$L(\cdot,\chi_4)$. Suppose that $\mathcal L$ has no zero in
$\Re s>\theta$. For every $b>\theta$ and every $\nu>0$,
\begin{equation}\label{eq:perron-zerofree-growth}
 |\mathcal L(\sigma+it)|+|1/\mathcal L(\sigma+it)|
       \ll_{b,\theta,\nu}(1+|t|)^\nu
       \qquad(\sigma\ge b,\ |t|\ge4).
\end{equation}
The pole of $\zeta$ at $1$ does not lie in this range of heights.
\end{lemma}

\begin{proof}
Two elementary integral formulas give the initial polynomial bound. Write
$B(u)=\sum_{n\le u}\chi_4(n)$, so that $0\le B(u)\le1$. For $\Re s>0$,
\begin{align*}
 \zeta(s)&=\frac{s}{s-1}
          -s\int_1^\infty\{u\}u^{-s-1}\,du\qquad(s\ne1),\\
 L(s,\chi_4)&=s\int_1^\infty B(u)u^{-s-1}\,du.
\end{align*}
To obtain them first take $\Re s>1$ and use
\[
 n^{-s}=s\int_n^\infty u^{-s-1}\,du.
\]
Interchanging the absolutely convergent sum with the integral gives
the counting functions $\lfloor u\rfloor$ and $B(u)$.
In the first formula use $\lfloor u\rfloor=u-\{u\}$.
The integrals on the displayed right sides converge locally uniformly for
$\Re s>0$, proving continuation to that half plane.
On each fixed strip with positive left boundary, and away from $s=1$,
they imply $|\mathcal L(s)|\ll 1+|s|$.

For $\Re s>1$ the Euler products give the absolutely convergent logarithms
\[
 \log\mathcal L(s)=\sum_p\sum_{m\ge1}
                 \frac{\xi(p)^m}{m p^{ms}},
 \qquad
 \xi(p)=1\ \hbox{or}\ \chi_4(p).
\]
In particular these logarithms are bounded on $\Re s\ge3/2$.
There is no need to consider $b\ge1$ separately, since such a half plane
is contained in one with $\theta<b<1$. Assume this latter range, and put
\[
 \delta=b-\theta,\qquad r_0=\frac12,\qquad r=2-b,\qquad
 R_1=2-\theta-\frac\delta2,\qquad
 R_2=2-\theta-\frac\delta4.
\]
Then $0<r_0<r<R_1<R_2<2$. At $s_0=2+it$, with $|t|\ge4$,
the closed disk $|s-s_0|\le R_2$ lies strictly in $\Re s>\theta$,
contains no zero and does not contain $1$. It therefore admits a
holomorphic branch
$g(z)=\log\mathcal L(s_0+z)$ on a neighborhood of the closed disk,
chosen to agree with the Euler logarithm at its center.
The polynomial bound gives
\[
 \max_{|z|\le R_2}\Re g(z)\le C+\log(1+|t|),
 \qquad |g(0)|\le C.
\]

For clarity, the two complex analysis inequalities used at this step are
the Borel--Carath\'eodory inequality and Hadamard's three-circles theorem,
in the precise forms of \cite[Theorems~6 and~7]{Simonic2023}.
If $M(r)=\max_{|z|=r}|g(z)|$, the first gives
\[
 M(R_1)\le
 \frac{2R_1}{R_2-R_1}\max_{|z|=R_2}\Re g(z)
 +\frac{R_2+R_1}{R_2-R_1}|g(0)|
 \ll_{b,\theta}\log(2+|t|).
\]
The required holomorphy and continuity on the closed disk were just verified.
On $|z|\le r_0$, the Euler logarithm gives $M(r_0)\ll1$.
The second inequality, applied to the closed annulus
$r_0\le|z|\le R_1$ within the same disk, now yields
\[
 M(r)\le M(r_0)^{1-a}M(R_1)^a
       \ll_{b,\theta}(\log(2+|t|))^a,
 \qquad
 a=\frac{\log(r/r_0)}{\log(R_1/r_0)}<1.
\]
For $b\le\sigma\le2$ the point $\sigma+it$ is in $|s-s_0|\le r$,
so the maximum principle gives the same bound for its logarithm.
For $\sigma\ge2$ absolute convergence already gives a constant bound.
Finally,
\[
 |\mathcal L(s)|+|1/\mathcal L(s)|
 \le2\exp\bigl(|\log\mathcal L(s)|\bigr)
 \le2\exp\bigl(C_{b,\theta}(\log(2+|t|))^a\bigr).
\]
Since $a<1$, the exponent is at most $\nu\log(1+|t|)$ for all
sufficiently large $|t|$. The remaining bounded range of heights is
compact for $b\le\sigma\le2$ and contains no zero or pole, so its
maximum is finite. This proves \eqref{eq:perron-zerofree-growth}.
\end{proof}

\begin{corollary}[M\"obius sums from a zero-free half plane]\label{cor:perron-mobius-halfplane}
Let $1/2\le\theta<1$. If $\zeta$ has no zero in $\Re s>\theta$, then
\[
 M(x)\ll_{\theta,\eps}x^{\theta+\eps}\qquad(\eps>0).
\]
If $L(\cdot,\chi_4)$ has no zero in that half plane, then
\[
 M_{\chi_4}(x):=\sum_{n\le x}\mu(n)\chi_4(n)
       \ll_{\theta,\eps}x^{\theta+\eps}\qquad(\eps>0).
\]
Each assertion uses the zero-free hypothesis only for its own function.
\end{corollary}

\begin{proof}
The Dirichlet series of the two coefficient sequences are $1/\zeta(s)$
and $1/L(s,\chi_4)$ for $\Re s>1$, by the Euler products. Both
coefficient sequences have modulus at most $1$.
Choose $0<\eps<1/4$ and
\[
 \delta=\min\bigl(\eps/4,(1-\theta)/2\bigr),\qquad b=\theta+\delta<1.
\]
The quotient $1/(s\mathcal L(s))$ is holomorphic on a neighborhood
of the closed strip $b\le\Re s\le5/4$. At $s=1$ in the zeta case,
the reciprocal has a zero, so the apparent singularity is removable.
Lemma~\ref{lem:perron-zerofree-growth} supplies the subpower bound.
Corollary~\ref{cor:perron-subpower}, with error exponent $\eps/2$,
has no residue to collect and gives
$\mathcal O(x^{b+\eps/2})=\mathcal O(x^{\theta+\eps})$ at half integers.
For every real $x\ge2$ set $y=\lfloor x\rfloor+1/2$.
The partial sums at $x$ and $y$ agree, and $y\asymp x$, giving the
same bound for every $x$. Larger values of $\eps$ follow by weakening
any one of the bounds already proved.
\end{proof}

Under RH the zeta assertion recovers the implication used in
Lemma~\ref{lem:littlewood}. The classical sharper logarithmic estimate
is \cite[Theorem~14.2]{Titchmarsh1986}, and the resulting M\"obius
criterion is \cite[Theorem~14.25(C)]{Titchmarsh1986}.
For the fixed primitive character $\chi_4$, an explicit logarithmic
estimate is also given in \cite[Corollary~2(a)]{Simonic2023}.
The proof above supplies the qualitative estimates needed in this volume
for every fixed zero-free boundary $\theta<1$.

\section{Perron inversion with Riesz weights}

A Riesz weight supplies enough decay for a polynomial vertical bound.\index[terms]{Mellin inversion}
The following version is used in Appendix~\ref{app:E} and the research dossier.

\begin{lemma}[Perron--Riesz inversion]\label{lem:perron-riesz}
Let $D(s)=\sum_{n\ge1}u_n n^{-s}$ converge absolutely on $\Re s=c>0$,
and let $r\ge1$ be an integer. For every real $x>0$,
\begin{equation}\label{eq:perron-riesz}
 \sum_{n\le x}u_n(1-n/x)^r
 =\frac{r!}{2\pi i}\int_{c-i\infty}^{c+i\infty}
       \frac{D(s)x^s}{s(s+1)\cdots(s+r)}\,ds.
\end{equation}
Suppose in addition that $D$ is meromorphic on a neighborhood of
$b\le\Re s\le c$, that $D(s)/[s(s+1)\cdots(s+r)]$ has finitely many
poles there, all strictly inside, and that
$|D(\sigma+it)|\ll(1+|t|)^M$ uniformly across the strip for large $|t|$.
If $r>M\ge0$, then
\begin{equation}\label{eq:perron-riesz-shift}
 \sum_{n\le x}u_n(1-n/x)^r
 =r!\sum_\rho\operatorname{Res}_{s=\rho}
       \frac{D(s)x^s}{s(s+1)\cdots(s+r)}+\mathcal O(x^b),
\end{equation}
where the sum is over the poles with $b<\Re\rho<c$.
\end{lemma}

\begin{proof}
Put $P_r(s)=s(s+1)\cdots(s+r)$. First consider
\[
 K_r(u)=\frac1{2\pi i}\int_{c-i\infty}^{c+i\infty}
                  \frac{u^s}{P_r(s)}\,ds\qquad(u>0).
\]
This integral is absolutely convergent. For $u>1$, move the line to
$\Re s=-R$, with $R>2r+1$. At a fixed $R$, the horizontal integrals
vanish as their heights tend to infinity because their integrands are
$\mathcal O_{u,R}(T^{-r-1})$. On the new vertical line,
\[
 \int_{\mathbb R}\left|\frac{u^{-R+it}}{P_r(-R+it)}\right|\,dt
 \ll_r u^{-R}\int_{\mathbb R}(R^2+t^2)^{-(r+1)/2}\,dt
 \ll_r u^{-R}R^{-r},
\]
which tends to zero as $R\to\infty$. The poles at $-j$, $0\le j\le r$,
have residues $(-1)^j u^{-j}/[j!(r-j)!]$, so the binomial identity gives
$K_r(u)=(1-u^{-1})^r/r!$. For $0<u\le1$, move instead to $\Re s=R$.
No pole is crossed, and the new integral is $\mathcal O_r(u^R R^{-r})$,
which tends to zero, including when $u=1$. Thus
$K_r(u)=(1-u^{-1})_+^r/r!$.

On the initial line, interchange the sum and integral using
\[
 x^c\sum_{n\ge1}|u_n|n^{-c}
       \int_{\mathbb R}|P_r(c+it)|^{-1}\,dt<\infty.
\]
The scalar identity with $u=x/n$ proves \eqref{eq:perron-riesz}.
It includes integer $x$, since the weight at $n=x$ is zero.

For the shift to $b$, the quotient is
$\mathcal O(|t|^{M-r-1})$ at large heights, uniformly across the strip.
For $x\ge2$, each horizontal integral is bounded by
\[
 C T^{M-r-1}\int_b^c x^\sigma\,d\sigma
 =C T^{M-r-1}\frac{x^c-x^b}{\log x}.
\]
The tails of either vertical integral are
$\mathcal O(x^\sigma T^{M-r})$ for $\sigma=b,c$.
Since $r>M$, the tails and horizontal integrals vanish as $T\to\infty$
at fixed $x$. The residue theorem gives the sum in
\eqref{eq:perron-riesz-shift}, with positive sign for a move to the left.
The full new vertical integral is $\mathcal O(x^b)$, because its compact
part has no pole and its tails are integrable. This proves the claim.
\end{proof}

\cleardoublepage
\setcounter{chapter}{0}
\renewcommand{\thechapter}{W}
\renewcommand{\theHchapter}{technicalboundary}
\renewcommand{\appendixname}{Technical appendix}
\renewcommand{\chaptername}{Technical appendix}
\chapter{Lemmas for the boundary transfer}
\label{app:boundary}
\markboth{LEMMAS FOR THE BOUNDARY TRANSFER}{LEMMAS FOR THE BOUNDARY TRANSFER}

This appendix proves the lemmas on which Chapter~\ref{chap:boundary} rests, a boundary lemma
for positive measures, the non-vanishing of Dirichlet $L$-functions on the line $\Re s=1$, and a
lemma on weighted averages of a sequence with small partial sums, together with the elementary
passage from $M(x)=o(x)$ to $\psi(x)\sim x$. None of them uses the
prime number theorem, in an arithmetic progression or otherwise. The contour estimates of
Appendix~\ref{app:perron} are not used here, and nothing here supplies them.

\section{A boundary lemma for positive measures}

\begin{proof}[Proof of Lemma~\ref{lem:positive_boundary}]
Put $B(x)=\nu([1,x])$, $c=R/\kappa$, and
\[
S(u)=e^{-\kappa u}B(e^{u})\ \ (u\ge0),\qquad S(u)=0\ \ (u<0).
\]
For $u\ge1$, positivity gives $D_\nu(\kappa+1/u)\ge e^{-(\kappa+1/u)u}B(e^{u})=e^{-1}S(u)$, and the
continuity of $D_\nu(s)-R/(s-\kappa)$ at $s=\kappa$ gives $D_\nu(\kappa+\eps)=R/\eps+\mathcal O(1)$
as $\eps\to0^{+}$. Hence $S(u)=\mathcal O(1+u)$. The monotonicity of $B$ gives the exact one-sided
comparison
\begin{equation}\label{eq:S_slow}
S(u+t)\ge e^{-\kappa t}S(u)\qquad(u,t\ge0).
\end{equation}
Let $f(u)=S(u)-c\,\mathbf 1_{[0,\infty)}(u)$. Tonelli's theorem and the substitution $x=e^{u}$ give
\[
P(s):=\int_0^{\infty}f(u)e^{-su}\,du=\frac{D_\nu(\kappa+s)}{\kappa+s}-\frac cs\qquad(\Re s>0),
\]
and with $E(w)=D_\nu(w)-R/(w-\kappa)$ and $c\kappa=R$,
\[
P(s)=\frac{E(\kappa+s)}{\kappa+s}-\frac{c}{\kappa+s},
\]
so that $P$ extends continuously to the closed half plane $\Re s\ge0$.

Let $\varphi$ be a Schwartz function whose Fourier transform is smooth with compact support. For
$\eps>0$ the function $e^{-\eps u}f(u)\mathbf 1_{[0,\infty)}(u)$ is integrable, with transform
$\xi\mapsto P(\eps+i\xi)$, and Fourier inversion for the convolution of two integrable functions
whose product of transforms is integrable gives
\[
\int_0^{\infty}f(v)e^{-\eps v}\varphi(u-v)\,dv
=\frac1{2\pi}\int_{\R}\widehat\varphi(\xi)\,P(\eps+i\xi)\,e^{i\xi u}\,d\xi .
\]
Let $\eps\to0$. On the left, dominated convergence applies, since $f(v)=\mathcal O(1+v)$ and
$\varphi(u-v)$ decays rapidly for fixed $u$. On the right, $P(\eps+i\xi)\to P(i\xi)$ uniformly on
the compact support of $\widehat\varphi$. The function $G=\widehat\varphi\,P(i\,\cdot)$ is
continuous with compact support, and the substitution $\xi\mapsto\xi+\pi/u$ gives
$2\int G(\xi)e^{i\xi u}\,d\xi=\int\big(G(\xi)-G(\xi+\pi/u)\big)e^{i\xi u}\,d\xi$, which tends to
zero as $u\to\infty$ by uniform continuity. Therefore
\begin{equation}\label{eq:S_bandlimited}
(\varphi*f)(u)\longrightarrow0,\qquad (\varphi*S)(u)\longrightarrow c\int_{\R}\varphi .
\end{equation}

At this point $S$ is only known to grow at most linearly. To bound it, take $\varphi$ non-negative
and positive near the origin, for instance the squared modulus of the inverse transform of a
non-negative smooth function with compact support, normalized to integral one. Its transform is
again smooth with compact support. For $\ell>0$ small enough that $\varphi>0$ on $[-\ell,0]$,
\eqref{eq:S_slow} gives
\[
(\varphi*S)(u)\ge\int_{-\ell}^{0}\varphi(v)\,S(u-v)\,dv\ge S(u)\int_{-\ell}^{0}\varphi(v)\,e^{\kappa v}\,dv .
\]
The last integral is positive and the left side is bounded by \eqref{eq:S_bandlimited}, so $S$ is
bounded.

Schwartz functions whose transform has compact support are dense in $L^{1}(\R)$, by the argument
closing the proof of Lemma~\ref{lem:wiener_step}. With $S$ bounded, \eqref{eq:S_bandlimited}
therefore extends to every integrable test function, and in particular the forward and backward
averages of $S$ over each interval of fixed length $\ell$ tend to $c$. By \eqref{eq:S_slow}, for
$u\ge\ell$,
\[
S(u)\le e^{\kappa\ell}\,\frac1\ell\int_0^{\ell}S(u+t)\,dt,\qquad
S(u)\ge e^{-\kappa\ell}\,\frac1\ell\int_0^{\ell}S(u-t)\,dt ,
\]
hence $e^{-\kappa\ell}c\le\liminf S(u)\le\limsup S(u)\le e^{\kappa\ell}c$. Letting $\ell\to0$ gives
$S(u)\to c$, which is the assertion. Positivity was used twice, to bound $S$ and to remove the
smoothing, and no bound on the size of the masses of $\nu$ was imposed.
\end{proof}

\section{Dirichlet \texorpdfstring{$L$}{L}-functions on the line \texorpdfstring{$\Re s=1$}{Re s=1}}

\begin{lemma}\label{lem:L_line}
Let $Q\ge2$ be fixed. For every nonprincipal character $\chi$ modulo $Q$ and every real $t$,
$L(1+it,\chi)\neq0$. In particular $L(1,\chi)\neq0$.
\end{lemma}

\begin{proof}
For nonprincipal $\chi$ the partial sums $S_\chi(u)=\sum_{m\le u}\chi(m)$ are bounded, and
\[
L(s,\chi)=s\int_1^{\infty}S_\chi(u)\,u^{-s-1}\,du
\]
is holomorphic on $\Re s>0$. For the principal
character $L(s,\chi_0)=\zeta(s)\prod_{p\mid Q}(1-p^{-s})$ has, on that half plane, a single
pole, simple, at $s=1$. Put
\[
Z_Q(s)=\prod_{\chi\bmod Q}L(s,\chi).
\]
If $f_p$ denotes the order of $p$ in $(\Z/Q\Z)^{\times}$, the characters modulo $Q$ give
$\prod_\chi(1-\chi(p)X)=(1-X^{f_p})^{\varphi(Q)/f_p}$ for $p\nmid Q$, and hence
\[
Z_Q(s)=\prod_{p\nmid Q}\big(1-p^{-f_ps}\big)^{-\varphi(Q)/f_p}\qquad(\Re s>1).
\]
The Dirichlet coefficients of $Z_Q$ are non-negative, and the coefficient of $m^{\varphi(Q)}$ is at least
$1$ whenever $m$ is prime to $Q$, since $f_p$ divides $\varphi(Q)$. The series therefore diverges at
$s=1/\varphi(Q)$, and its abscissa of convergence $\sigma_c$ satisfies $\sigma_c\ge1/\varphi(Q)>0$.

A Dirichlet series $D(s)=\sum b_nn^{-s}$ with $b_n\ge0$ and finite abscissa $\sigma_c$ cannot be
holomorphic in a neighbourhood of $\sigma_c$. This is the theorem of
Landau\index[names]{Landau, E.}. Indeed, if it were, choose a real $\sigma_0>\sigma_c$ so close to
$\sigma_c$ that the Taylor series of $D$ at $\sigma_0$ converges at a point $\sigma_0-r<\sigma_c$.
Since $(-1)^{j}D^{(j)}(\sigma_0)=\sum_nb_n(\log n)^{j}n^{-\sigma_0}\ge0$, Tonelli's theorem
applied to that Taylor series shows that $\sum_nb_nn^{-(\sigma_0-r)}$ converges, which contradicts
the definition of $\sigma_c$. The sum of the series agrees with $Z_Q$ on $\Re s>\sigma_c$, by
uniqueness of continuation from $\Re s>1$.

If some nonprincipal $L(1,\chi)$ vanished, it would cancel the only possible pole of $Z_Q$ on
$\Re s>0$, and $Z_Q$ would be holomorphic on that half plane, in particular near $\sigma_c$,
contrary to Landau's theorem. Hence $L(1,\chi)\neq0$ for every nonprincipal $\chi$, and $Z_Q$ has a
simple pole at $1$. Finally the logarithm of the Euler product of $Z_Q$ has non-negative
coefficients, and $3+4\cos\theta+\cos2\theta\ge0$ gives
\[
Z_Q(\sigma)^{3}\,|Z_Q(\sigma+it)|^{4}\,|Z_Q(\sigma+2it)|\ge1\qquad(\sigma>1).
\]
A zero of $Z_Q$ at $1+it$ with $t\neq0$ would make the left side tend to zero as
$\sigma\to1^{+}$, the pole at $1$ being simple and $Z_Q$ holomorphic at $1+2it$. So $Z_Q$ does
not vanish on the line $\Re s=1$ away from $1$, and neither does any of its factors, each of which
is finite there.
\end{proof}

The argument uses the Euler product and complex continuation, and no Tauberian theorem. The
prime number theorem in progressions deduced from it in Chapter~\ref{chap:boundary} is therefore
not obtained from itself, and it is not elementary in the sense of avoiding complex analysis.

\section{Weighted averages of a sequence with small partial sums}

\begin{lemma}\label{lem:weighted_cancellation}
Let $(b_n)$ be complex with $|b_n|\le1$ and $\sum_{n\le x}b_n=o(x)$. Let $f$ be Riemann integrable
on every interval $[\eps,1]$ with $\eps>0$, and let $|f(t)|\le C(1+\log(1/t))$ on $(0,1]$. Then
\[
\frac1x\sum_{n\le x}b_n\,f\Big(\frac nx\Big)\longrightarrow0 .
\]
\end{lemma}

\begin{proof}
Fix $\eps>0$. The terms with $n\le\eps x$ contribute at most
\[
\frac Cx\sum_{n\le\eps x}\Big(1+\log\frac xn\Big)\ll\eps\big(1+\log(1/\eps)\big)+o(1).
\]
Since $f$ is Riemann integrable on $[\eps,1]$, that interval splits into finitely many intervals
$I_1,\dots,I_r$ on which the oscillations $\omega_j$ of $f$ satisfy
$\sum_j\omega_j|I_j|<\eps$. Choose $t_j\in I_j$ and let $f_\eps=f(t_j)$ on $I_j$. On each
interval $I_j=(c_j,d_j]$ the hypothesis gives
\[
\frac1x\sum_{c_jx<n\le d_jx}b_n\longrightarrow0,
\]
so $x^{-1}\sum_{\eps x<n\le x}b_n\,f_\eps(n/x)\to0$. The rest is at most
$x^{-1}\sum_{\eps x<n\le x}|f-f_\eps|(n/x)\le\sum_j\omega_j\big(|I_j|+1/x\big)$, whose upper
limit is below $\eps$. Letting $x\to\infty$ and then $\eps\to0$ proves the lemma.
\end{proof}

\section{From the M\"obius function to the prime number theorem}

\begin{proposition}\label{prop:mobius_bridge}
If $M(x)=\sum_{n\le x}\mu(n)=o(x)$, then $\psi(x)\sim x$.
\end{proposition}

\begin{proof}
All three steps are instances of Lemma~\ref{lem:weighted_cancellation} with $b_d=\mu(d)$.
First, the identity of Meissel, $\sum_{d\le x}\mu(d)\lfloor x/d\rfloor=1$, gives
\[
M_{-1}(x)=\frac1x+\frac1x\sum_{d\le x}\mu(d)\Big\{\frac xd\Big\},
\]
and the lemma with $f(t)=\{1/t\}$ gives $M_{-1}(x)\to0$. Next, with $H(y)=\sum_{n\le y}1/n=\log
y+\gamma+R(y)$ and $R(y)=\mathcal O(1/y)$, the identity
$\sum_{d\le x}(\mu(d)/d)H(x/d)=\sum_{N\le x}N^{-1}\sum_{d\mid N}\mu(d)=1$ gives
\[
1=\sum_{d\le x}\frac{\mu(d)}{d}\log\frac xd+\gamma M_{-1}(x)+\frac1x\sum_{d\le x}\mu(d)\,f\Big(\frac dx\Big),
\qquad f(t)=\frac{R(1/t)}{t},
\]
where $f$ is bounded and Riemann integrable on every $[\eps,1]$. The lemma gives
\[
K(x):=\sum_{d\le x}\frac{\mu(d)}{d}\log\frac xd\longrightarrow1 .
\]
Finally, M\"obius inversion of $\sum_{d\mid n}\Lambda(d)=\log n$ gives
\[
\psi(x)=\sum_{d\le x}\mu(d)\log\Big(\Big\lfloor\frac xd\Big\rfloor!\Big),
\]
and the estimate $\log(\lfloor y\rfloor!)=y\log y-y+E(y)$ with $E(y)=\mathcal O(\log2y)$ gives
\[
\frac{\psi(x)}{x}=K(x)-M_{-1}(x)+\frac1x\sum_{d\le x}\mu(d)\,E\Big(\frac xd\Big).
\]
The lemma with $f(t)=E(1/t)=\mathcal O(1+\log(1/t))$ makes the last term tend to zero, so
$\psi(x)/x\to1$.
\end{proof}

\chapter*{Register of open statements}
\label{chap:register}
\addcontentsline{toc}{chapter}{Register of open statements}
\markboth{REGISTER OF OPEN STATEMENTS}{REGISTER OF OPEN STATEMENTS}
\makeatletter
\def\ps@registerhead{\let\@oddfoot\@empty\let\@evenfoot\@empty
  \def\@oddhead{{\slshape REGISTER OF OPEN STATEMENTS}\hfil\thepage}%
  \let\@evenhead\@oddhead\let\@mkboth\@gobbletwo}
\makeatother
\pagestyle{registerhead}

Open problems, conjectures and conditional theorems are marked where they occur. This register
collects them, together with passages whose proof-status notes distinguish established results,
quoted inputs and remaining questions. A proof-status note can also document a fully proved
result. The gauged index and transfer quoted from the companion volume are two such cases,
with their exact sources and uses recorded in Chapter~\ref{chap:gauge_ingham}.
The tables preserve these distinctions. Numbers and pages are those of the text.
Solutions, partial results and counterexamples are welcome, at the address given at the end of the
preface.

\section*{Open problems}

\begin{longtable}{@{}p{1.6cm}p{1.1cm}p{10.4cm}@{}}
\toprule
\textbf{Number} & \textbf{Page} & \textbf{Question} \\
\midrule
\endhead
\ref{op:xi_reciprocal} & \pageref{op:xi_reciprocal} & Which hypotheses, weaker than the weighted mass condition, force the transparent coefficient to be the reciprocal of the transform, the continuation hypotheses alone being insufficient by Proposition~\ref{prop:xi_mass_needed} \\
\ref{op:ex-green-boundary} & \pageref{op:ex-green-boundary} & A universal Schur weight, a rowwise sharpness criterion, and necessity at the boundary of the Green criterion \\
\ref{op:ex-fgv-jumps} & \pageref{op:ex-fgv-jumps} & Sharpness for finitely many signed jumps, countably many jumps beyond absolute variation, and singular continuous profiles \\
\ref{op:ex-lacunary-original} & \pageref{op:ex-lacunary-original} & Whether the real zero of the original lacunary transform is its leftmost zero, and whether that profile is a function of good variation \\
\ref{op:w10-global} & \pageref{op:w10-global} & The two global estimates LOW and ABS for the broken harmonic kernel at $\sqrt2$ \\
\ref{op:rational_gauge_general_frontier} & \pageref{op:rational_gauge_general_frontier} & Sharp frontiers of the rational kernel under negative feedback for bounded-quotient, iterated-power and irregular gauges \\
\ref{op:ab_mean_defect} & \pageref{op:ab_mean_defect} & Additional regularity ensuring a mean formula for infinite-variation profiles with zero-mass Stieltjes layers \\
\ref{op:divisor_filters} & \pageref{op:divisor_filters} & Structured realization of the quotient profile $\omega_p$ with finite exponential type or by an integer linear recurrence, for $p\ge3$ \\
\ref{op:trace-poly-general-spectrum} & \pageref{op:trace-poly-general-spectrum} & A complete Frobenius and coefficient expansion of the finite resolvent at arbitrary fixed polynomial degree \\
\ref{op:w04d-boundary-and-amplitudes} & \pageref{op:w04d-boundary-and-amplitudes} & Nonvanishing of $C_1$, first-difference control proving the pointwise spectral conjecture, and finite-rank discrete Gamma-mode expansions \\

\ref{op:G_green} & \pageref{op:G_green} & The two Green estimates for the two branch discontinuous profile \\
\ref{op:ab_riemann_discrepancy} & \pageref{op:ab_riemann_discrepancy} & Which admissible kernels outside the smooth class have a Riemann discrepancy of second order \\
\ref{op:O_balance} & \pageref{op:O_balance} & Which profiles are balanced in the sense of the reciprocal defect \\
\bottomrule
\end{longtable}

\section*{Conjectures}

\begin{longtable}{@{}p{1.6cm}p{1.1cm}p{10.4cm}@{}}
\toprule
\textbf{Number} & \textbf{Page} & \textbf{Statement} \\
\midrule
\endhead
\ref{conj:fgv_membership} & \pageref{conj:fgv_membership} & Every nonconstant positive recurrence-admissible profile that is bounded or slowly varying at the origin is a function of good variation \\
\ref{conj:ortho_spectral} & \pageref{conj:ortho_spectral} & The first spectral pair of the orthorecursive coefficients, with a nonzero amplitude \\
\ref{conj:sqrt2_collapse} & \pageref{conj:sqrt2_collapse} & The collapse $\alpha(g_{\sqrt2})=\tfrac12$, strictly below the analytic index \\
\ref{conj:anti_hlr} & \pageref{conj:anti_hlr} & The anti-HLR\index[terms]{anti-HLR statement} statement for broken harmonic functions with zeros off the convergence line \\
\ref{conj:homogeneity} & \pageref{conj:homogeneity} & The homogeneity principle\index[terms]{homogeneity principle} for kernels of constant transform \\
\ref{conj:ingham_equilibrium} & \pageref{conj:ingham_equilibrium} & Equilibrium of the Ingham function under the offset gauges $f_q(x)=q^{x}+1$ \\
\ref{conj:G_index_revised} & \pageref{conj:G_index_revised} & The index of the two branch discontinuous profile and its endpoint behavior \\
\ref{conj:J_main} & \pageref{conj:J_main} & The index and the transparent and absorbed regimes of the greatest common divisor kernel \\
\bottomrule
\end{longtable}

\noindent The statements above are the conjectures of this volume. The Riemann hypothesis
is stated once, as Conjecture~\ref{conj:rh} on page~\pageref{conj:rh}, because
Theorem~\ref{thm:tauberian_rh} makes it equivalent to $\alpha(\Phi)=\tfrac12$ and the
equivalence has to name what it is equivalent to. It is not one of them.

\section*{Conditional theorems}

\begin{longtable}{@{}p{1.6cm}p{1.1cm}p{10.4cm}@{}}
\toprule
\textbf{Number} & \textbf{Page} & \textbf{Hypothesis assumed} \\
\midrule
\endhead
\ref{cthm:ex-alpha-eta} & \pageref{cthm:ex-alpha-eta} & Subcritical transparency, absorption at and above the candidate threshold, exact inversion, integrable vertical bounds and a nonzero residue at the first zero \\
\ref{cthm:w10-conditional-RAF} & \pageref{cthm:w10-conditional-RAF} & The two global estimates LOW and ABS \\
\ref{cthm:main} & \pageref{cthm:main} & Transparency with a power saving, uniform on compacts \\
\ref{cthm:G_green_closure} & \pageref{cthm:G_green_closure} &
(G2) for the transparent main term and (G1) for its error and RAF
absorption; the sharper endpoint and above-threshold estimates of
Conjecture~\ref{conj:G_index_revised} are not included \\
\bottomrule
\end{longtable}

\section*{Statements and passages qualified by a proof-status note}

The statements and programmatic passages below are proved, or proved in part, and each carries a
note recording exactly what is established and what is not. They are listed because a reader scanning for the word
conjecture would not find them. Further notes of the same kind are attached to statements
already listed above, among the open problems, the conjectures and the conditional
theorems, and they are not repeated here.

\begin{longtable}{@{}p{1.6cm}p{1.1cm}p{10.4cm}@{}}
\toprule
\textbf{Number} & \textbf{Page} & \textbf{What the note records} \\
\midrule
\endhead
\ref{thm:transfer_law} & \pageref{thm:transfer_law} & Cases (I), (III) and (IV) proved here, and case (II), Ingham's theorem, proved in Corollary~\ref{cor:ingham_proved} and unused in later proofs \\
\ref{thm:gen_equiv} & \pageref{thm:gen_equiv} & The unconditional equivalence with the two zero hypotheses, the Perron estimates proved in Appendix~\ref{app:perron}, and the exclusion of transparency above one half \\
\ref{lem:sv_forcing} & \pageref{lem:sv_forcing} & Transfer quoted from \cite[Theorem~5.1]{CloitreFloor} under the three stated conditions on the slowly varying forcing; used only in the two following examples \\
\ref{cor:zerofree} & \pageref{cor:zerofree} & Two proved implications with different assumptions. A zero-free half plane gives a lower bound on the transparency threshold, while the reverse implication assumes membership \\
\ref{cor:H_index} & \pageref{cor:H_index} & The equivalence proved through RH; an alternative direct transfer from cumulative-resolvent decay remains open \\
\ref{thm:ortho_conditional_book} & \pageref{thm:ortho_conditional_book} & Closure proved under uniform stability; stability is open, and the stated consistency estimate restricts this route to $\sigma\le1$ \\
\ref{def:resolvent_regular} & \pageref{def:resolvent_regular} & Primitive and Volterra probes converge under (R1); the full resolvent conditions (R2)--(R3) remain programmatic and are assumed by no theorem \\
\ref{sec:contour} & \pageref{sec:contour} & The discrete defect and the application of Appendix~\ref{app:perron}. Estimate (T) is proved for Ingham under RH, while the general profile retains its quantitative hypotheses \\
\ref{numobs:gauge_critical} & \pageref{numobs:gauge_critical} & The proved decay, nonvanishing and half-index quoted from the companion volume, and their exact use in the equilibrium statement \\
\ref{prop:gauge_transfer} & \pageref{prop:gauge_transfer} & The summability and perturbation estimates quoted from the companion volume which make the gauged transfer unconditional \\
\ref{thm:fibonacci_exact} & \pageref{thm:fibonacci_exact} & The two exact Fibonacci remainders quoted from \cite[Theorem~1.1]{CloitreFibonacci}, their circle and divisor sources, and their uses in Corollary~\ref{cor:fibonacci_equivalence} and the following remarks \\
\ref{cor:fibonacci_equivalence} & \pageref{cor:fibonacci_equivalence} & The unconditional exponent equivalence: forward from Theorem~\ref{thm:fibonacci_exact}, converse from \cite[Proposition~6.1]{CloitreFibonacci}, with no later theorem depending on it \\
\ref{prop:E_below} & \pageref{prop:E_below} & The direct comparison range $\beta<-1$, completed for all $\beta<0$ by Propositions~\ref{prop:E_negative_band} and~\ref{prop:E_negative_band_complete} and assembled in Corollary~\ref{cor:E_regimes} \\
\bottomrule
\end{longtable}

\backmatter
\pagestyle{headings}
\index[terms]{FGV|see{function of good variation}}\index[terms]{RAF|see{regular arithmetic function}}\index[terms]{RH|see{Riemann hypothesis}}\index[terms]{HLR|see{Hardy--Littlewood--Ramanujan criterion}}

\printindex[names]
\printindex[terms]


\addcontentsline{toc}{chapter}{Bibliography}
\begin{thebibliography}{99}
\bibitem{Agarwal2000}\index[names]{Agarwal, R. P.} R.P. Agarwal, \emph{Difference Equations and Inequalities}, 2nd ed., Marcel Dekker, New York (2000).
\bibitem{BaezDuarte2003}\index[names]{B\'aez-Duarte, L.} L. B\'aez-Duarte, \emph{A strengthening of the Nyman-Beurling criterion for the Riemann hypothesis}, Atti Accad. Naz. Lincei Cl. Sci. Fis. Mat. Natur. Rend. Lincei (9) Mat. Appl. \textbf{14} (2003), no.~1, 5--11.

\bibitem{BalazardBB2021}\index[names]{Balazard, M.}\index[names]{Benferhat, L.}\index[names]{Bouderbala, M.} M. Balazard, L. Benferhat, M. Bouderbala, \emph{Sur la variation de certaines suites de parties fractionnaires}, Comm. Math. \textbf{29} (2021), no.~3, 407--430.
\bibitem{BenzaidLutz1987}\index[names]{Benzaid, Z.}\index[names]{Lutz, D. A.} Z. Benzaid and D. A. Lutz, \emph{Asymptotic representation of solutions of perturbed systems of linear difference equations}, Stud. Appl. Math. \textbf{77} (1987), 195--221.
\bibitem{Beurling1955}\index[names]{Beurling, A.} A. Beurling, \emph{A closure problem related to the Riemann zeta-function}, Proc. Nat. Acad. Sci. U.S.A. \textbf{41} (1955), 312--314.

\bibitem{Bingham1989}\index[names]{Bingham, N. H.}\index[names]{Goldie, C. M.}\index[names]{Teugels, J. L.} N.H. Bingham, C.M. Goldie, J.L. Teugels, \emph{Regular Variation}, Encyclopedia of Mathematics and its Applications \textbf{27}, Cambridge University Press (1989), 1st ed. 1987.
\bibitem{BinghamInoue2000a}\index[names]{Bingham, N. H.}\index[names]{Inoue, A.} N.H. Bingham and A. Inoue, \emph{Abelian, Tauberian and Mercerian theorems for arithmetic sums}, J. Math. Anal. Appl. \textbf{250} (2000), 465--493.
\bibitem{BinghamInoue2000b}\index[names]{Bingham, N. H.}\index[names]{Inoue, A.} N.H. Bingham and A. Inoue, \emph{Tauberian and Mercerian theorems for systems of kernels}, J. Math. Anal. Appl. \textbf{252} (2000), 177--197.
\bibitem{BinghamInoue2000c}\index[names]{Bingham, N. H.}\index[names]{Inoue, A.} N.H. Bingham and A. Inoue, \emph{Extension of the Drasin-Shea-Jordan theorem}, J. Math. Soc. Japan \textbf{52} (2000), no.~3, 545--559.
\bibitem{BombieriGhosh2011}\index[names]{Bombieri, E.}\index[names]{Ghosh, A.} E. Bombieri and A. Ghosh, \emph{Around the Davenport--Heilbronn function}, Uspekhi Mat. Nauk \textbf{66} (2011), no.~2, 15--66; English transl. Russian Math. Surveys \textbf{66} (2011), no.~2, 221--270.
\bibitem{Broughan2017I}\index[names]{Broughan, K.} K. Broughan, \emph{Equivalents of the Riemann Hypothesis. Volume One, Arithmetic Equivalents}, Encyclopedia of Mathematics and its Applications \textbf{164}, Cambridge University Press (2017).
\bibitem{Broughan2017II}\index[names]{Broughan, K.} K. Broughan, \emph{Equivalents of the Riemann Hypothesis. Volume Two, Analytic Equivalents}, Encyclopedia of Mathematics and its Applications \textbf{165}, Cambridge University Press (2017).

\bibitem{Borwein2009}\index[names]{Borwein, P.} P. Borwein, S. Choi, B. Rooney, A. Weirathmueller, \emph{The Riemann Hypothesis: A Resource for the Afficionado and Virtuoso Alike}, CMS Books in Mathematics, Springer, New York (2008).
\bibitem{Brunner2017}\index[names]{Brunner, H.} H. Brunner, \emph{Volterra Integral Equations: An Introduction to Theory and Applications}, Cambridge University Press, 2017.
\bibitem{Cloitre2016}\index[names]{Cloitre, B.} B. Cloitre, \emph{Good variation theory: a Tauberian approach to the Riemann hypothesis}, Int.\ J.\ Math.\ Comp.\ Sci.\ \textbf{11} (2016), no.~2, 133--149.
\bibitem{CloitreFloor}\index[names]{Cloitre, B.} B. Cloitre, \emph{A Tauberian characterisation of the Riemann hypothesis through the floor function}, \texttt{arXiv:2407.18859} (2024).
\bibitem{CloitreFibonacci}\index[names]{Cloitre, B.} B. Cloitre, \emph{Fibonacci, Dirichlet, and Gauss in a single sum}, Fibonacci Quart., to appear. Preprint \texttt{arXiv:2607.20960} (2026).
\bibitem{CloitreOrtho}\index[names]{Cloitre, B.} B. Cloitre, \emph{On the orthorecursive expansion of unity}, Int. J. Number Theory, published online (2026), DOI 10.1142/S1793042127500060; preprint \texttt{arXiv:2505.09645v2} (2025).
\bibitem{CloitreVolII}\index[names]{Cloitre, B.} B. Cloitre, \emph{Regular Arithmetic Functions, Volume II. The gauge principle applied to the Ingham function and its deformations: from theorems to routes toward the Riemann hypothesis and its generalizations}, in preparation.
\bibitem{Conrey1989}\index[names]{Conrey, J. B.} J.~B. Conrey,
\emph{More than two fifths of the zeros of the Riemann zeta function are on the critical line},
J. reine angew. Math. \textbf{399} (1989), 1--26, DOI 10.1515/crll.1989.399.1.
\bibitem{Connes2026}\index[names]{Connes, A.} A. Connes, \emph{The Riemann Hypothesis: Past, Present and a Letter Through Time}, \texttt{arXiv:2602.04022} (2026).
\bibitem{Daval2019}\index[names]{Daval, F.} F. Daval, \emph{Identit\'es int\'egrales et estimations explicites associ\'ees pour les fonctions sommatoires li\'ees \`a la fonction de M\"obius et autres fonctions arithm\'etiques}, th\`ese, Universit\'e de Lille (2019).
\bibitem{DavenportHeilbronn1936}\index[names]{Davenport, H.}\index[names]{Heilbronn, H.} H. Davenport, H. Heilbronn, \emph{On the zeros of certain Dirichlet series I, II}, J. London Math. Soc. \textbf{11} (1936), 181--185 and 307--312.
\bibitem{Deligne1974}\index[names]{Deligne, P.} P. Deligne, \emph{La conjecture de Weil. I.}, Inst. Hautes Études Sci. Publ. Math., \textbf{43} (1974), 273--307.
\bibitem{DrasinShea1976}\index[names]{Drasin, D.}\index[names]{Shea, D. F.} D. Drasin and D.F. Shea, \emph{Convolution inequalities, regular variation and exceptional sets}, J. Analyse Math. \textbf{29} (1976), 232--293.
\bibitem{Dyson1969}\index[names]{Dyson, F. J.} F.~J. Dyson, \emph{Existence of a phase-transition in a one-dimensional Ising ferromagnet}, Comm. Math. Phys. \textbf{12} (1969), 91--107.
\bibitem{DysonWoS}\index[names]{Dyson, F. J.} F.~J. Dyson, \emph{Littlewood's Tauberian theorem and Dyson's ferromagnet paper}, Web of Stories, story~26, recorded June 1998.
\bibitem{EndresSteiner}\index[names]{Egger n\'e Endres, S.}\index[names]{Steiner, F.} S. Egger n\'e Endres and F. Steiner, \emph{A new proof of the Vorono\"i summation formula}, J. Phys. A: Math. Theor. \textbf{44} (2011), no.~22, 225302; preprint \texttt{arXiv:1001.3556}.
\bibitem{ErdosSegal1978}\index[names]{Erd\H{o}s, P.}\index[names]{Segal, S. L.} P. Erdős and S.L. Segal, \emph{A note on Ingham's summation method}, J. Number Theory \textbf{10} (1978), 95--98.
\bibitem{FlajoletGourdonDumas1995}\index[names]{Flajolet, P.}\index[names]{Gourdon, X.}\index[names]{Dumas, P.} P. Flajolet, X. Gourdon, and P. Dumas, \emph{Mellin transforms and asymptotics: harmonic sums}, Theoret. Comput. Sci. \textbf{144} (1995), 3--58.
\bibitem{Folland1992}\index[names]{Folland, G. B.} G.B. Folland, \emph{Fourier Analysis and its Applications}, Wadsworth, Belmont, CA, 1992.
\bibitem{Gripenberg1990}\index[names]{Gripenberg, G.}\index[names]{Londen, S.-O.}\index[names]{Staffans, O.} G. Gripenberg, S.-O. Londen, O. Staffans, \emph{Volterra Integral and Functional Equations}, Encyclopedia of Mathematics and its Applications \textbf{34}, Cambridge University Press, 1990.
\bibitem{deHaan1970}\index[names]{Haan, L. de@de Haan, L.} L. de Haan, \emph{On Regular Variation and Its Application to the Weak Convergence of Sample Extremes}, Mathematical Centre Tracts \textbf{32}, Mathematisch Centrum, Amsterdam, 1970.
\bibitem{Hadamard1892}\index[names]{Hadamard, J.} J. Hadamard, \emph{Essai sur l'\'etude des fonctions donn\'ees par leur d\'eveloppement de Taylor}, J. Math. Pures Appl. (4) \textbf{8} (1892), 101--186.
\bibitem{HarcosMO2024}\index[names]{Harcos, G.} G. Harcos (posting as ``GH from MO''), answer to \emph{Limit involving the fractional part and the Fibonacci numbers}, MathOverflow, 30 March 2024, \url{https://mathoverflow.net/a/468048}.
\bibitem{Hardy1916}\index[names]{Hardy, G. H.} G.H. Hardy, \emph{On Dirichlet's divisor problem}, Proc. London Math. Soc. (2) \textbf{15} (1916), 1--25.
\bibitem{Hardy1949}\index[names]{Hardy, G. H.} G.~H. Hardy, \emph{Divergent Series}, Clarendon Press, Oxford (1949).
\bibitem{HardyLittlewood1914}\index[names]{Hardy, G. H.}\index[names]{Littlewood, J. E.} G. H. Hardy and J. E. Littlewood, \emph{Tauberian theorems concerning power series and Dirichlet's series whose coefficients are positive}, Proc. London Math. Soc. (2) \textbf{13} (1914), 174--191.
\bibitem{Hindry2012}\index[names]{Hindry, M.} M. Hindry, \emph{La preuve par Andr\'e Weil de l'hypoth\`ese de Riemann pour une courbe sur un corps fini}, in \emph{Henri Cartan \& Andr\'e Weil, math\'ematiciens du XX\textsuperscript{e} si\`ecle}, Journ\'ees math\'ematiques X-UPS, \'Editions de l'\'Ecole polytechnique (2012), 65--101.
\bibitem{Holder1887}\index[names]{H\"older, O.} O. H\"older, \emph{\"Ueber die Eigenschaft der Gammafunction keiner algebraischen Differentialgleichung zu gen\"ugen}, Math. Ann. \textbf{28} (1887), 1--13.
\bibitem{Huxley2003}\index[names]{Huxley, M. N.} M.N. Huxley, \emph{Exponential sums and lattice points III}, Proc. London Math. Soc. (3) \textbf{87} (2003), no.~3, 591--609.
\bibitem{Ikehara1931}\index[names]{Ikehara, S.} S. Ikehara, \emph{An extension of Landau's theorem in the analytical theory of numbers}, J. Math. Phys. \textbf{10} (1931), 1--12. DOI 10.1002/sapm19311011.
\bibitem{Ingham1942}\index[names]{Ingham, A. E.} A.~E. Ingham, \emph{On two conjectures in the theory of numbers}, Amer. J. Math. \textbf{64} (1942), 313--319.
\bibitem{Ingham1945}\index[names]{Ingham, A. E.} A.E. Ingham, \emph{Some Tauberian theorems connected with the prime number theorem}, J. London Math. Soc. \textbf{20} (1945), 171--180.
\bibitem{Ismail2005}\index[names]{Ismail, M. E. H.} M.E.H. Ismail, \emph{Classical and Quantum Orthogonal Polynomials in One Variable}, with two chapters by W. Van Assche, Encyclopedia of Mathematics and its Applications \textbf{98}, Cambridge University Press (2005).
\bibitem{Ivic2003}\index[names]{Ivi\'c, A.} A. Ivi\'c, \emph{The Riemann Zeta-Function: Theory and Applications}, Dover Publications, Mineola, NY (2003), reprint of the John Wiley \& Sons edition, New York (1985).
\bibitem{Jordan1974}\index[names]{Jordan, G. S.} G.S. Jordan, \emph{Regularly varying functions and convolutions with real kernels}, Trans. Amer. Math. Soc. \textbf{194} (1974), 177--194.
\bibitem{Jukes1971}\index[names]{Jukes, K. A.} K.~A. Jukes, \emph{On the Ingham and $(D,h(n))$ summation methods}, J. London Math. Soc. (2) \textbf{3} (1971), no.~4, 699--710.
\bibitem{Jutila2015}\index[names]{Jutila, M.} M. Jutila, \emph{Riemann's zeta-function and the divisor problem. III}, Ark. Mat. \textbf{53} (2015), 303--315.
\bibitem{KalmyninKosenko2020}\index[names]{Kalmynin, A. B.}\index[names]{Kosenko, P. R.} A.B. Kalmynin and P.R. Kosenko, \emph{Orthorecursive expansion of unity}, Int. J. Number Theory \textbf{16}(6) (2020), 1209--1226.
\bibitem{Kanemitsu2018}\index[names]{Kanemitsu, S.}\index[names]{Kuzumaki, T.}\index[names]{Tanigawa, Y.} S. Kanemitsu, T. Kuzumaki, and Y. Tanigawa, \emph{On the Wintner-Ingham-Segal summability method}, Hardy-Ramanujan J. \textbf{41} (2018), 40--47.
\bibitem{Karamata1930}\index[names]{Karamata, J.} J. Karamata, \emph{Sur un mode de croissance r\'eguli\`ere des fonctions}, Mathematica (Cluj) \textbf{4} (1930), 38--53.
\bibitem{Karamata1930b}\index[names]{Karamata, J.} J. Karamata, \emph{\"Uber die Hardy-Littlewoodschen Umkehrungen des Abelschen Stetigkeitssatzes}, Math. Z. \textbf{32} (1930), 319--320.
\bibitem{Khinchin1964}\index[names]{Khinchin, A. Ya.} A.~Ya. Khinchin, \emph{Continued Fractions}, University of Chicago Press, 1964.
\bibitem{Korevaar2004}\index[names]{Korevaar, J.} J. Korevaar, \emph{Tauberian Theory: A Century of Developments}, Grundlehren der mathematischen Wissenschaften \textbf{329}, Springer, 2004.
\bibitem{Levinson1956}\index[names]{Levinson, N.} N. Levinson, \emph{On closure problems and the zeros of the Riemann zeta function}, Proc. Amer. Math. Soc. \textbf{7} (1956), 838--845.
\bibitem{Levinson1974}\index[names]{Levinson, N.} N. Levinson, \emph{More than one third of zeros of Riemann's zeta-function are on $\sigma=1/2$}, Advances in Math. \textbf{13} (1974), 383--436.

\bibitem{LiYang2023}\index[names]{Li, X.}\index[names]{Yang, X.} X. Li and X. Yang, \emph{An improvement on Gauss's circle problem and Dirichlet's divisor problem}, arXiv:2308.14859 (2023).
\bibitem{Littlewood1911}\index[names]{Littlewood, J. E.} J.~E. Littlewood, \emph{The converse of Abel's theorem on power series}, Proc. London Math. Soc. (2) \textbf{9} (1911), 434--448.
\bibitem{Littlewood1912}\index[names]{Littlewood, J. E.} J.~E. Littlewood, \emph{Quelques cons\'equences de l'hypoth\`ese que la fonction $\zeta(s)$ de Riemann n'a pas de z\'eros dans le demi-plan $\Re s>\tfrac12$}, C.~R. Acad. Sci. Paris \textbf{154} (1912), 263--266.
\bibitem{Mordell1917}\index[names]{Mordell, L. J.} L.J. Mordell, \emph{On Mr.\ Ramanujan's Empirical Expansions of Modular Functions}, Proc. Cambridge Philos. Soc. \textbf{19} (1917), 117--124.
\bibitem{Newman1980}\index[names]{Newman, D. J.} D.~J. Newman, \emph{Simple analytic proof of the prime number theorem}, Amer. Math. Monthly \textbf{87} (1980), 693--696.
\bibitem{Nyman1950}\index[names]{Nyman, B.} B. Nyman, \emph{On some groups and semigroups of translations}, Thesis, Uppsala University (1950).

\bibitem{OEIS} OEIS Foundation Inc., \emph{The On-Line Encyclopedia of Integer Sequences}, published electronically at \texttt{https://oeis.org/}, accessed 2026.
\bibitem{DLMF}\index[names]{Olver, F. W. J.} F.~W.~J. Olver et al. (eds.), \emph{NIST Digital Library of Mathematical Functions}, \texttt{https://dlmf.nist.gov/}, Release 1.2.7 of 2026-06-15.
\bibitem{Ostrowski1920}\index[names]{Ostrowski, A.} A. Ostrowski, \emph{\"Uber Dirichletsche Reihen und algebraische Differentialgleichungen}, Math. Z. \textbf{8} (1920), 241--298.
\bibitem{Pennington1955}\index[names]{Pennington, W. B.} W.~B. Pennington, \emph{On Ingham summability and summability by Lambert series}, Proc. Cambridge Philos. Soc. \textbf{51} (1955), 65--80.
\bibitem{Pillai1933}\index[names]{Pillai, S. S.} S.~S. Pillai, \emph{On an arithmetic function}, J. Annamalai Univ. \textbf{2} (1933), 243--248.
\bibitem{Poincare1902}\index[names]{Poincar\'e, H.} H. Poincar\'e, \emph{La Science et l'Hypoth\`ese}, Flammarion, Paris, 1902, quatri\`eme partie, chapitre~IX.
\bibitem{PoincareMethode1908}\index[names]{Poincar\'e, H.} H. Poincar\'e, \emph{Science et M\'ethode}, Flammarion, Paris, 1908, livre~I, chapitre~II.
\bibitem{Polya1917}\index[names]{P\'olya, G.} G. P\'olya, \emph{\"Uber eine neue Weise, bestimmte Integrale in der analytischen Zahlentheorie zu gebrauchen}, G\"ottinger Nachrichten (1917), 149--159.
\bibitem{Rajagopal1955}\index[names]{Rajagopal, C. T.} C.~T. Rajagopal, \emph{A note on Ingham summability and summability by Lambert series}, Proc. Indian Acad. Sci. Sect. A \textbf{42} (1955), 41--50.
\bibitem{Schiff1993}\index[names]{Schiff, J. L.} J.~L. Schiff, \emph{Normal Families}, Universitext, Springer, New York, 1993.
\bibitem{Schmidt1980}\index[names]{Schmidt, W. M.} W.~M. Schmidt, \emph{Diophantine Approximation}, Lecture Notes in Mathematics \textbf{785}, Springer, 1980.
\bibitem{Segal1966}\index[names]{Segal, S. L.} S.~L. Segal, \emph{On Ingham's summation method}, Canad. J. Math. \textbf{18} (1966), 97--105.
\bibitem{Segal1975}\index[names]{Segal, S. L.} S.~L. Segal, \emph{Ingham's summability method and Riemann's hypothesis}, Proc. London Math. Soc. (3) \textbf{30} (1975), no.~2, 129--142.
\bibitem{Segal1976}\index[names]{Segal, S. L.} S.~L. Segal, \emph{Corrigendum to Ingham's summability method and Riemann's hypothesis}, Proc. London Math. Soc. (3) \textbf{34} (1977), no.~3, 438.
\bibitem{Segal1979}\index[names]{Segal, S. L.} S.~L. Segal, \emph{A second note on Ingham's summation method}, Canad. Math. Bull. \textbf{22} (1979), no.~1, 117--120.
\bibitem{Selberg1942}\index[names]{Selberg, A.} A. Selberg, \emph{On the zeros of Riemann's zeta-function}, Skr. Norske Vid.-Akad. Oslo I \textbf{10} (1942), 59 pp.
\bibitem{Selberg1992}\index[names]{Selberg, A.} A. Selberg, \emph{Old and new conjectures and results about a class of Dirichlet series}, in \emph{Proceedings of the Amalfi Conference on Analytic Number Theory}, Università di Salerno (1992), pp.\ 367--385.
\bibitem{Seneta1976}\index[names]{Seneta, E.} E. Seneta, \emph{Regularly Varying Functions}, Lecture Notes in Mathematics \textbf{508}, Springer, Berlin, 1976.
\bibitem{Simonic2023}\index[names]{Simoni\v{c}, A.} A. Simoni\v{c}, \emph{Estimates for $L$-functions in the critical strip under GRH with effective applications}, Mediterr. J. Math. \textbf{20} (2023), no.~2, article~87. DOI 10.1007/s00009-023-02289-2.
\bibitem{Szego1975}\index[names]{Szeg\H{o}, G.} G. Szeg\H{o}, \emph{Orthogonal Polynomials}, 4th ed., American Mathematical Society Colloquium Publications \textbf{XXIII}, American Mathematical Society, Providence, RI (1975).
\bibitem{Tauber1897}\index[names]{Tauber, A.} A. Tauber, \emph{Ein Satz aus der Theorie der unendlichen Reihen}, Monatsh. Math. Phys. \textbf{8} (1897), 273--277.
\bibitem{Tenenbaum2015}\index[names]{Tenenbaum, G.} G. Tenenbaum, \emph{Introduction to Analytic and Probabilistic Number Theory}, 3rd ed., Graduate Studies in Mathematics \textbf{163}, American Mathematical Society, 2015.
\bibitem{Titchmarsh1986}\index[names]{Titchmarsh, E. C.} E.~C. Titchmarsh, \emph{The Theory of the Riemann Zeta-Function}, 2nd ed., revised by D.~R. Heath-Brown, Oxford University Press, 1986.
\bibitem{Vainikko2009}\index[names]{Vainikko, G.} G. Vainikko, \emph{Cordial Volterra integral equations 1}, Numer. Funct. Anal. Optim. \textbf{30} (2009), no.~9--10, 1145--1172.
\bibitem{Vainikko2010}\index[names]{Vainikko, G.} G. Vainikko, \emph{Cordial Volterra integral equations 2}, Numer. Funct. Anal. Optim. \textbf{31} (2010), no.~2, 191--219.
\bibitem{Vainikko2014}\index[names]{Vainikko, G.} G. Vainikko, \emph{First kind cordial Volterra integral equations 2}, Numer. Funct. Anal. Optim. \textbf{35} (2014), no.~12, 1607--1637.
\bibitem{Waterhouse1969}\index[names]{Waterhouse, W. C.} W.~C. Waterhouse, \emph{Abelian varieties over finite fields}, Ann. Sci. \'Ecole Norm. Sup. (4) \textbf{2} (1969), 521--560.
\bibitem{WeilAvenir1947}\index[names]{Weil, A.} A. Weil, \emph{L'avenir des math\'ematiques}, in F. Le Lionnais (ed.), \emph{Les Grands Courants de la Pens\'ee Math\'ematique}, Cahiers du Sud, Marseille, 1948.
\bibitem{WeilMetaphysique1960}\index[names]{Weil, A.} A. Weil, \emph{De la m\'etaphysique aux math\'ematiques} (1960), in \OE uvres Scientifiques / Collected Papers, vol.~II, Springer, Berlin (1979).
\bibitem{Weyl1916}\index[names]{Weyl, H.} H. Weyl, \emph{\"Uber die Gleichverteilung von Zahlen mod. Eins}, Math. Ann. \textbf{77} (1916), 313--352.
\bibitem{Wiener1932}\index[names]{Wiener, N.} N. Wiener, \emph{Tauberian theorems}, Ann. of Math. (2) \textbf{33} (1932), 1--100.
\bibitem{Wintner1943}\index[names]{Wintner, A.} A. Wintner, \emph{Eratosthenian Averages}, Waverly Press, Baltimore (1943).
\bibitem{WongLi1992}\index[names]{Wong, R.}\index[names]{Li, H.} R. Wong and H. Li, \emph{Asymptotic expansions for second-order linear difference equations}, J. Comput. Appl. Math. \textbf{41} (1992), no.~1--2, 65--94. DOI 10.1016/0377-0427(92)90239-T.
\bibitem{YuMO2024}\index[names]{Yu, Hung-Hsun} H.-H. Yu, answer to \emph{On properties of sums involving the floor function}, MathOverflow, 7 April 2024, \url{https://mathoverflow.net/q/468499}.
\bibitem{Zacharovas2011}\index[names]{Zacharovas, V.} V. Zacharovas, \emph{A Tauberian theorem for the Ingham summation method}, Acta Arith. \textbf{148} (2011), no.~1, 31--54.
\end{thebibliography}
\end{document}